\documentclass[reqno]{amsart}
\input epsf
\address{Simons Center for Geometry and Physics,
State University of New York, Stony Brook, NY 11794-3636 U.S.A. \&
Center for Geometry and Physics, Institute for Basic Sciences (IBS), Pohang, Korea} \email{kfukaya@scgp.stonybrook.edu}
\address{Center for Geometry and Physics, Institute for Basic Sciences (IBS), Pohang, Korea \& Department of Mathematics,
POSTECH, Pohang, Korea} \email{yongoh1@postech.ac.kr}
\address{Graduate School of Mathematics,
Nagoya University, Nagoya, Japan} \email{ohta@math.nagoya-u.ac.jp}
\address{Research Institute for Mathematical Sciences, Kyoto University, Kyoto, Japan}
\email{ono@kurims.kyoto-u.ac.jp}

\usepackage{graphicx}
\usepackage{amsmath}
\usepackage{amscd}
\usepackage{amssymb}
\usepackage{amstext}
\usepackage{amsmath}
\usepackage[all]{xy}
\usepackage{yhmath}
\usepackage{mathrsfs}
\usepackage{bbding}
\usepackage[dvipdfmx]{hyperref}
\hypersetup{
 setpagesize=false,
 bookmarksnumbered=true,%
 bookmarksopen=true,%
 colorlinks=true,%
 linkcolor=blue,
 citecolor=green,
}
\makeatletter
\def\l@section{\@tocline{1}{0pt}{3mm}{8mm}{}}
\def\l@subsection{\@tocline{2}{0pt}{6mm}{10mm}{}}
\def\l@subsubsection{\@tocline{3}{0pt}{9mm}{11mm}{}}
\makeatother

\def\E{\ifmmode{\mathbb E}\else{$\mathbb E$}\fi} 
\def\N{\ifmmode{\mathbb N}\else{$\mathbb N$}\fi} 
\def\R{\ifmmode{\mathbb R}\else{$\mathbb R$}\fi} 
\def\Q{\ifmmode{\mathbb Q}\else{$\mathbb Q$}\fi} 
\def\C{\ifmmode{\mathbb C}\else{$\mathbb C$}\fi} 
\def\H{\ifmmode{\mathbb H}\else{$\mathbb H$}\fi} 
\def\Z{\ifmmode{\mathbb Z}\else{$\mathbb Z$}\fi} 
\def\P{\ifmmode{\mathbb P}\else{$\mathbb P$}\fi} 
\def\T{\ifmmode{\mathbb T}\else{$\mathbb T$}\fi} 
\def\SS{\ifmmode{\mathbb S}\else{$\mathbb S$}\fi} 
\def\DD{\ifmmode{\mathbb D}\else{$\mathbb D$}\fi} 
\def\K{\ifmmode{\mathbb K}\else{$\mathbb K$}\fi}

\newcommand{\del}{\partial}

\newcommand{\rank}{\operatorname{rank}}

\theoremstyle{theorem}
\newtheorem{thm}{Theorem}[section]
\newtheorem{cor}[thm]{Corollary}
\newtheorem{lem}[thm]{Lemma}
\newtheorem{sublem}[thm]{Sublemma}
\newtheorem{subsublemma}[thm]{Subsublemma}
\newtheorem{prop}[thm]{Proposition}

\newtheorem{lemdef}[thm]{Lemma-Definition}

\theoremstyle{definition}
\newtheorem{defn}[thm]{Definition}
\newtheorem{rem}[thm]{Remark}

\newtheorem{cons}[thm]{\rm\bfseries{Construction}}
\newtheorem{exm}[thm]{Example}
\newtheorem{conds}[thm]{Condition}
\newtheorem{conven}[thm]{Convention}
\newtheorem{proper}[thm]{Property}
\newtheorem{defnlem}[thm]{Definition-Lemma}
\newtheorem{assump}[thm]{Assumption}

\newtheorem{shitu}[thm]{Situation}
\newtheorem{choi}[thm]{Choice}
\newtheorem{notation}[thm]{Notation}

\numberwithin{equation}{section}
\makeindex
\begin{document}

\title[Kuranishi structure
and Virtual fundamental chain]{Kuranishi structure,
Pseudo-holomorphic curve,
and Virtual fundamental chain: Part 1}
\author{Kenji Fukaya, Yong-Geun Oh, Hiroshi Ohta, Kaoru Ono}

\thanks{Kenji Fukaya is supported partially by JSPS Grant-in-Aid for Scientific Research
No. 23224002 and NSF Grant No. 1406423,
Yong-Geun Oh by the IBS project IBS-R003-D1, Hiroshi Ohta by JSPS Grant-in-Aid
for Scientific Research No. 23340015 and Kaoru Ono by JSPS Grant-in-Aid for
Scientific Research, Nos. 26247006, 23224001.}

\begin{abstract}
This is the first part of the article we promised at the
end of \cite[Section 1]{foootech}.
We discuss the foundation of the virtual fundamental
chain and cycle technique, especially
its version appeared in \cite{FO} and also in
\cite[Section A1, Section 7.5]{fooobook2}, \cite[Section 12]{fooo09},
\cite{fooo091}.
In Part 1, we focus on the construction
of the virtual fundamental chain
on a single space with Kuranishi structure.
We mainly discuss the
de Rham version and so work over $\R$-coefficients,
but we also include a self-contained
account of the way how to work over $\Q$-coefficients in case the dimension of
the space with Kuranishi structure is $\le 1$.
\par
Part 1 of this document is independent of our
earlier writing \cite{foootech}.
We also do not assume the reader have any knowledge on
the pseudo-holomorphic curve, in Part 1.
\par
Part 2 (resp. Part 3), which will appear in the near future,
discusses the case of a system of Kuranishi structures
and its simultaneous perturbations
(resp. the way to implement the abstract story
in the study of moduli spaces of pseudo-holomorphic curves).
\end{abstract}
\maketitle
\par\newpage
\date{Feb 7th, 2015}


\maketitle

\tableofcontents
\newpage

\section{Introduction}

The purpose of this article is to present the virtual fundamental chain
technique and its application to the moduli space of pseudo-holomorphic curves,
in as much detail as possible.
The technique of virtual fundamental cycles and chains was invented in the year 1996 by
groups of mathematicians \cite{FO}, \cite{LiTi98}, \cite{LiuTi98}, \cite{Rua99}, \cite{Siebert}
and applied to many purposes subsequently.
Several references such as \cite{fooo06}, \cite{ChenTian}, \cite{ChenLieWang}, \cite{LuT}
provide various versions and details of this technique.

\par
In this article (besides various improvement of the detail of the proof, presentation etc.)
we add the following points to  our earlier writing \cite{foootech}
of similar nature.
\par
\begin{enumerate}
\item
We present the detail of the chain level argument,
especially in the de-Rham version.
\item
We provide a package of the statements which arise from the
virtual fundamental chain and cycle technique,
in the way one can directly quote and use without referring to their proofs.
\item
We give constructions of algebraic structures from systems of
spaces with Kuranishi structures. We separate the topological and algebraic issues from
analytic ones of the story so that this part of the construction can be
rigorously stated and proved without referring to the construction of such
system.\footnote{In our main application, such a system is constructed from
the moduli spaces of pseudo-holomorphic curves. However those results we mention here
is independent of the origin of such a system.}
\item
We discuss in detail the construction of a system of spaces with Kuranishi structures
that satisfies the axioms we mentioned in Item (3). Compared to \cite{foootech},
we put more emphasis on the consistency of  Kuranishi structures we obtain.
Namely our emphasis lies in the construction of such a system of Kuranishi structures
rather than that of a single Kuranishi structure.
\item
We also provide details of the proof of certain fundamental
properties of Gromov-Witten invariants.
In this article we clarify the fact that we can work out this proof without using continuous family of
perturbations. In other words, we need only the technique that already appeared in \cite{FO}
but not those developed in the later literature.
\item
We also explain a method of working with the $\Q$-coefficient without using
a triangulation of the perturbed moduli space.
This is possible if we only need to study moduli space of
virtual dimension 0,1 and negative. Note that these cases handle all the applications
appearing in \cite{FO}.
\end{enumerate}
\par

This article is divided into three parts.
\par
Part 1 studies space with Kuranishi structure.
Here we consider only single space with Kuranishi structure and not its system.
We also write Part 1 (and Part 2) so that it can be read without any knowledge of
pseudo-holomorphic curves. Especially readers' knowledge on analytic details on
pseudo-holomorphic curves is not required.
Occasionally we mention the moduli space of pseudo-holomorphic curves in Parts 1 and 2,
but none of the proof relies on them. This is because we do so only at the place where
we think some motivation of the construction or the definition etc. would be useful.
We believe that for the reader whose background is more on topology and/or
algebra but not on analysis, this way of writing is useful.
We also believe that our writing clarifies the fact that problem of
transversality of the moduli space of pseudo-holomorphic curves
lies in the realm of differential topology rather than
one in functional or nonlinear analysis.
\par
Item (1) mentioned above comprises Part 1.
\par
In Part 2 we study systems of Kuranishi structures and provide a systematic
way of building algebraic structures on certain chain complex arising from such a system.
We discuss the case of Floer homology (of periodic Hamiltonian system)
and of $A_{\infty}$ algebra appearing in Lagrangian Floer theory.
Item (3) above comprises Part 2.
\par
In Part 3 we present in detail the construction of the moduli spaces
of pseudo-holomorphic curves and its Kuranishi structure
which satisfies the axioms we formulate in Part 2.
We discuss them for the case of Floer homology (of periodic Hamiltonian system)
and of $A_{\infty}$ algebra appearing in Lagrangian Floer theory.
We also discuss the case of Gromov-Witten invariants.
(Item  (5) is in Part 3.)
\par
We remark Parts 1 and 2 of this article are self-contained.
It does not use the results of \cite{foootech}.
In other words, we reprove certain parts of \cite{foootech} in this article.
(We quote some general topology issue from \cite{foooshrink}.)
In Part 3, we quote \cite{foootech} in several places.
\par
Introduction to each part follows.
\par
Now we upload Part 1.
In the text below we quote results, sections etc.
of Parts 2 and 3. They will appear later on.
Note we quote Parts 2 and 3 in order to provide a
view how the results of Part 1 are used or
to provide perspective to related results.
As far as the logic of the proof concerns,
the discussion of Part 1 is complete without
referring Parts 2 and 3.

\label{sec:introduction}
\subsection{Introduction to Part 1}
\label{subsec:intro1}
In Part 1 we discuss the space with Kuranishi structure
and define its virtual fundamental chain.
\par
Kuranishi structure is a `generalization of the notion of manifold
which allows singularity'.
We entirely work in $C^{\infty}$-category.
In general, singularity in $C^{\infty}$-category can be much wilder
than one in analytic category. Our aim is not developing a general theory
of `$C^{\infty}$ analogue of scheme' but a more restricted one, which is, to define
virtual fundamental chains and cycles. Because of this, our definitions and
constructions are designed so that they are useful for that particular purpose.
\par
Now an outline of the main results proved in Part 1 is in order.
Since most basic and conceptual ideas behind the
definitions of Kuranishi structures and virtual fundamental
chains have been explained in various other previous writings of ours
(see \cite{fooo010} \cite{Fusuugaku} \cite[Section A1]{fooobook}
\cite[Section 12]{fooo09} \cite{FO} for example), here
we discuss more technical aspect of the story, that is, the way
how we rigorously realize those basic concepts.
In fact the aim of this article is {\it not} to explain the basic concepts
of the theory or to present the main ideas of the construction of virtual
fundamental chains (which, according to the opinion of the authors of this article,
is definitely more important than the details provided here)
but to present its {\it technical} details with as much details as possible
\par
Roughly speaking Kuranishi structure
of a space $X$ represents $X$ locally as
a zero set $s^{-1}(0)$ where $s$ is a section of
vector bundle $E \to U$ over an orbifold $U$.
More precisely a Kuranishi chart is a quadruple
$\mathcal U = (U,E,s,\psi)$ where
$s$ is a section of a
vector bundle $E \to U$ over an orbifold $U$ and
$\psi : s^{-1}(0) \to X$ is a homeomorphism to its image.
(See Definition \ref{kuranishineighborhooddef}.)
We call $s$ a \emph{Kuranishi map}. In fact this description in the situation of
various applications appears as a Kuranishi model of the local
description of the moduli problem.
We also need appropriate notion of coordinate changes among them.
\par
The main difference of the coordinate change between Kuranishi charts from the coordinate change
between manifold or orbifold charts lies in the fact that the coordinate change of
Kuranishi charts may not be a local `isomorphism'.
This is because the dimension of the space with Kuranishi structure is, by definition, $\dim U - {\rm rank} E$
and so the orbifolds $U_i$ of  two Kuranishi charts $(U_i,E_i,s_i,\psi_i)$ of the
same space $X$ may have different dimension. (Namely $\dim U_1 \ne \dim U_2$ in general.)
Recall that in the description of the manifold or orbifold structure one uses
a pseudo-group or a groupoid.
This description certainly is {\it not} useful for the Kuranishi structure since the most
important axiom of groupoid (the invertibility of the morphisms) is exactly the one we give up.
One may use a version of category of fractions to invert the arrows
appearing in the definition of coordinate change of Kuranishi charts
to adapt the groupoid language. This is the route taken, for example, by
D. Joyce in \cite{joyce2}, which however is not the one we take.
\par
There are two slightly different versions of the definition of
`spaces with Kuranishi type structures' in \cite{FO}.
One is what we call `a Kuranishi structure' the other `a good coordinate system'.
\par
In the manifold theory (or in the orbifold theory) we consider a maximal set of the coordinate charts
compatible with one another, and call this set a structure of $C^{\infty}$ manifold.
Equivalently for a manifold $M$ we assign a (germ of a) coordinate neighborhood to each point $p$ of $M$.
This way of defining the manifold structure is somehow more canonical.
On the other hand to perform various operations using the given manifold structure,
we usually take a locally finite cover consisting of coordinate charts in the
given manifold structure. The latter is especially necessary when we use
a partition of unity as in the case of defining integration of differential form on
a manifold $M$. Certain amount of general topology is needed to be worked out
in the manifold course to prove existence of an appropriate locally finite cover
extracted from the set of infinitely many charts defining the manifold structure.
This technicality is not entirely trivial but is never a conceptual heart of the
theory.
\footnote{When teaching a course of manifold, it is important
not to over-emphasize the importance of these issues so that the students will not
get lost in the technicalities.}
\par
The relationship between the Kuranishi structure and the good coordinate system is
similar to the one between these two ways of describing a $C^{\infty}$ structure
of manifolds. A Kuranishi structure of $X$ assigns a Kuranishi chart
$\mathcal U_p$ to each point $p \in X$.
So in particular it contains uncountably many charts.
On the other hand, a good coordinate system of $X$ consists of
locally finite (finite if $X$ is compact) Kuranishi charts $\mathcal U_{\frak p}$. ($\frak p \in \frak P$.  Here
$\frak P$ is certain index set.)
\par
Note coordinate changes between Kuranishi charts are not necessarily a local
isomorphism. In general we have coordinate change only in one direction.
In the case of Kuranishi structure the coordinate change from
$\mathcal U_{q}$ to $\mathcal U_{p}$ is defined
only when $q \in U_{p}$. (Here $U_p$ is a neighborhood of $p$.)
We remark that $q \in U_p$ does not imply $p \in U_q$.
\par
In the case of good coordinate system the set of Kuranishi charts is parameterized
by a certain partially ordered set $\frak P$ and
the coordinate change from
$\mathcal U_{\frak q}$ to $\mathcal U_{\frak p}$
is defined only when $\frak q \le \frak p$.
\begin{rem}(Chen-Tian's approach)
In the paper  \cite{ChenTian} \footnote{The paper  \cite{ChenTian} was
published in 2010 but this paper had appeared in the arXiv in 2006.} (and several other references
which follow Chen-Tian's approach)
the set of charts are parameterized by a subset of certain finite set $I$.
(Namely $\frak P \subseteq 2^I$.)
The partial order in Chen-Tian's approach is the usual inclusion of subsets.
(See \cite[Definition 2.1]{ChenTian}.)
If we examine the proof of the construction of Kuranishi structure
in  \cite{FO}, we find that the partially ordered set which appears there
is actually parameterized by a subset of certain subset of the moduli space.
(Namely $I$ is a finite subset of $X$.)
Thus our approach is not so far away from Chen-Tian's approach in
applications.
We, however, think our way of including arbitrary partially ordered set is more flexible than
Chen-Tian's. (Their approach certainly works at least for their purpose.)
\par
For example the approach using $\frak P \subseteq 2^I$ does not seem to work literally
in the case when we try to construct a system of Kuranishi structures
on moduli spaces of pseudo-holomorphic curves (without boundary) so that
it is consistent with the fiber product description of the stratum.
(See Part 3.)\footnote{One reason why it does not seem to work
is that there is an automorphism which exchange irreducible components, in general.}
One can certainly modify or use its variant to include that case.
However modifying the definition when new applications appear is not a good
way to develop general theory, especially such a lengthy one we are building.
Because then one need either to repeat all the proofs again or to say `We can prove it in the same way'.
In this article we want to avoid saying this sentence as much as possible
though it might be impossible to completely avoid it.\footnote{We had experiences
to be complained when we wrote similar sentences in our earlier papers.
We agree that it is preferable to avoid it. However it is sometimes impossible to do so,
especially in a research paper whose main focus is application and is not the detail of foundation.}
The structure of partially ordered set is
certainly the structure which remains in all the expected generalizations
and is also the very structure we need for the inductive argument we will use to work
on good coordinate system.
\end{rem}
The main reason why the notion of Kuranishi structure was introduced is to work out
the transversality issue appearing in the theory of pseudo-holomorphic curves and etc.
For this purpose it is important to perturb the Kuranishi map, denoted by $s$, so that
the resulting perturbed map becomes transversal to $0$.
The original proof of transversality theorem due to R. Thom \cite{thom}
constructed perturbation by using induction on charts. We follow
this proof. We use partial order on the index set $\frak P$
of Kuranishi charts for this purpose. Note in order to perform
this inductive construction it is essential to
use good coordinate system rather than Kuranishi structure
since it is hard to work on uncountably many coordinate charts
consisting Kuranishi structure.
Thus  construction of good coordinate system out of Kuranishi structure
is an important step in the construction of an appropriate perturbation and resolving
transversality issue.
\par
On the other hand, the fact that the Kuranishi structure is more canonical than
the good coordinate system can be seen from the following point.
We can define the product of two spaces with Kuranishi structure in a canonical way.
However the definition of product of two spaces with good coordinate system
is more complicated and is less natural.
For our application we need to construct virtual fundamental chains
in the way that is compatible with fiber product or with direct product.
For this purpose working with Kuranishi structure is more appropriate.
\par
Thus during the various constructions, we need to go from
a Kuranishi structure to a good coordinate system and back several times.
\par
Concerning this transition from one to the other, there is one difference between
the current story and the manifold theory.
In the manifold theory, one starts with a $C^{\infty}$ structure (which
consists of all the compatible coordinate charts) and then pick up an appropriate
locally finite cover. It is a fact that one can always recover original
$C^{\infty}$ structure from the charts consisting of the chosen locally finite cover.
In the current story of Kuranishi structures, we start with a Kuranishi structure $\widehat{\mathcal U}$
and construct a compatible good coordinate system $\widetriangle{\mathcal U}$.
We can use the resulting good coordinate system $\widetriangle{\mathcal U}$
to construct another Kuranishi structure $\widehat{\mathcal U^+}$.
Unfortunately $\widehat{\mathcal U^+}$ is different from
$\widehat{\mathcal U}$ in general.
In other words we lose certain information while we go from a Kuranishi structure
to its associated good coordinate system.
Some portion of the discussion in Part 1 is devoted to showing that this loss
does not affect the construction of virtual fundamental chain.
\begin{rem}
This technical trouble seems to reflect the fact
that the notion of Kuranishi structure is conceptually less canonical
than that of manifold. For example
we can construct a Kuranishi structure on (practically all the) moduli spaces
of pseudo-holomorphic curves. However the Kuranishi structure
we obtain is not unique but depends on some choices.
(We compensate this shortcoming by using an appropriate notion of cobordism.)
In case the moduli space happens to be a manifold (and the equation
defining the moduli space is Fredholm regular) the $C^{\infty}$
structure on the moduli space is certainly canonical.\footnote{except the choice
of smooth structure of gluing parameter.}
Another sign that the notion of Kuranishi structure is not
canonical enough is that we do not know how to define a morphism between two spaces
with Kuranishi structure.
\par
It seems that the route taken by Joyce \cite{joyce2}, \cite{joyce4} resolves these two issues.
See Remark \ref{joycerem} for the reason why we nevertheless use our definitions of
Kuranishi structure and good coordinate system.
\par
The history of the 20 century mathematics tells us that
in the realm of algebro-geometric or complex analytic category,
going to the ultimatum in making all the constructions,  definitions and etc.
as canonical as possible is often the correct point of view\footnote{This is the viewpoint
taken and insisted by Grotendieck.} even if it looks
cumbersome at the beginning. It also tells us that in the realm of differential geometry or in that of
$C^{\infty}$ functions, an attempt to realize a canonical construction in the ultimate level
frequently fails
\footnote{Theory of singularity of $C^{\infty}$ functions is one typical example.}
and so we are forced to find the suitable place to compromise and to be content ourselves with
being able to achieve the particular purpose we pursue.
\end{rem}
\begin{rem}
The whole exposition of Part 1 is designed in the way that most of the statements and proofs
find their analogs in the corresponding statements and proofs
of the standard theory of manifolds. Once the right statements are given, the proofs are
fairly obvious most of the times,\footnote{By this reason in many places we do not need to say much about the proof.
} although it is not entirely so because
there are also some differences between the Kuranishi structure and the manifold structure
in certain technical points.
These differences result in cumbersome and technically
heavy proofs in those cases. However we emphasize that these are rather technical points and
do not comprise the conceptional heart of the theory. Especially for the researchers who
have in their mind the construction of virtual fundamental chains in various concrete geometric
problems, thorough understanding of cumbersome details of these technicality should not be a part of
everyone's required background.
\footnote{Thorough knowledge
of such a technicality, of course, should be shared among the
people whose interest lie also on extending the technology to the extreme of its
potential border and/or  using the most delicate and difficult case of the
technology to obtain as the sharpest results as possible.}
\par
It seems to us that the whole theory now becomes nontrivial only because its presentation is heavy and lengthy.
The method of making the proof `locally trivial',
which we are taking here, has been used in various branches of mathematics
as an established method for building the foundation of a new theory, especially when the theory
is conceptually simple but meets certain complicated technicality in its rigorous details.
We choose the `Bourbaki style' way of writing
this article. Especially, we do our best effort to make explicit the assumptions
we work with and the conclusions we obtain. As its consequence, there appear so many
definitions and statements in the text, which we acknowledge is annoying but
may be inevitable.
\par
We hope that in the near future many users of virtual techniques via the Kuranishi structure
appreciate that they do not need to know anything more than a small number of basic definitions
and theorems together with the fact that the story of Kuranishi structure is mostly similar to
that of smooth manifolds. Then one can safely dispose the details such as those
we provide here and use the conclusions as a `black box'.
It is the present authors' opinion that advent of such an enlightening and agreement
in the area will be important for the development of symplectic geometry and related fields.
In fact, one main obstacle to its smooth development has been the
nuisance of working out similar, but not precisely the same, details each time when one
tries to use the moduli spaces of pseudo-holomorphic curves of various kinds in various situations.
\end{rem}
Now we provide description of the main results of Part 1 sectionwise.
In Section \ref{sec:skuraterm}, we give the definitions of Kuranishi structures (Definition \ref{kstructuredefn})
and good coordinate systems (Definition \ref{gcsystem}). These definitions are based on a
version of the definitions of an orbifold, its embedding and a vector bundle on it,
which we describe in Section \ref{sec:ofd}.
Of course the basic concepts and the mathematical contents of orbifold
were established long time ago \cite{satake}.
However there are a few different ways of defining an orbifold
in the technical point of view.
More significantly, the notion of morphisms between orbifolds is rather delicate to
define. We refer readers to the discussion of the book \cite{ofdruan} especially
its section 1 for these points.
In this article, we restrict ourselves to the world of effective orbifolds and use only embeddings as maps between them.
Then those delicate points and troubles disappear.
\par
In Section \ref{sec:fiber}, we define the fiber and the direct product of spaces with Kuranishi structure.
We remark that in the category theory the notion of fiber product is defined
in purely abstract language of objects and morphisms.
The definition of fiber product we give in this section, however, is {\it not} the one
given in the category theory.\footnote{The authors thank D. Joyce we call attention this point to us.}
In fact, we never define the notion of morphisms between spaces equipped with Kuranishi structure in general.
Therefore we define a fiber product of spaces with Kuranishi structure directly in Section \ref{sec:fiber}.
Here we consider only the fiber product over a manifold and require an appropriate transversality
to define a fiber product.
Although the definition we provide is not in the sense of category theory, they are
so natural and canonical that the basic properties expected for the
fiber product are fairly manifest. For example,  associativity of the
fiber product follows rather immediately from its definition.
We do not define a fiber or a direct product of spaces with good coordinate systems
here since its definition necessarily becomes more technical and complicated.
\par
In this article we assume all the Kuranishi structures and good coordinate systems are oriented
(Definitions \ref{kuraorient} and \ref{gcsystem}) unless otherwise mentioned.
\par
In Section \ref{sec:thick} we discuss the process going from a Kuranishi structure
to a good coordinate system and back.
As we mentioned before, we start from a Kuranishi structure
$\widehat{\mathcal U}$ and construct a good coordinate system
$\widetriangle{\mathcal U}$ as follows.
\par\medskip
\noindent
{\bf Theorem \ref{Them71restate}.}
{\it For any Kuranishi structure $\widehat{\mathcal U}$ of
$Z \subseteq X$ there exist a
good coordinate system ${\widetriangle{\mathcal U}}$
of $Z \subseteq X$
and a KG-embedding $\widehat{\mathcal U} \to {\widetriangle{\mathcal U}}$.}
\par\medskip
See Definition \ref{defn32020202} for KG-embedding.
(The good coordinate system ${\widetriangle{\mathcal U}}$ above is said to be
{\it compatible} with the Kuranishi structure $\widehat{\mathcal U}$.)
We also start from  a good coordinate system
$\widetriangle{\mathcal U}$
and construct a Kuranishi structure
$\widehat{\mathcal U^{+}}$.
(Proposition \ref{lemappgcstoKu}.)
When $\widehat{\mathcal U^{+}}$ is obtained from $\widehat{\mathcal U}$
by combining these two processes, the Kuranishi structure $\widehat{\mathcal U^{+}}$ is in general different
from $\widehat{\mathcal U}$ but $\widehat{\mathcal U^{+}}$ is related to
$\widehat{\mathcal U}$ in a way that
$\widehat{\mathcal U^{+}}$ is a thickening of $\widehat{\mathcal U}$.
We define this notion of thickening and study its properties in Section \ref{sec:thick}.
\par
To formulate the notion of compatibility between a Kuranishi structure and
a good coordinate system,
and compatibility among Kuranishi structures or among good coordinate systems,
we use the notion of embeddings.
We define 4 possible versions of such embeddings and their compositions
(there are 8 possible versions of compositions) in Section \ref{sec:thick}.
See Tables 5.1 and 5.2.
In Section \ref{sec:contgoodcoordinate}, the proof of basic results of existence of good coordinate system
compatible with given Kuranishi structure (Theorem \ref{Them71restate}) is given.
\par
In Section \ref{sec:multisection}, we define the notion of
multisections and multi-valued perturbations and their transversality.
The basic result is the following existence theorem of transversal perturbation, which is proved in Section \ref{sec:constrsec}.��
Here we leave
the precise definitions of the terminology used in the statement
to the later text.
\par\medskip\noindent
{\bf Theorem \ref{prop621}.}
{\it Let ${\widetriangle{\mathcal U}}$ be
a good coordinate system  of $Z \subseteq X$ and $\mathcal K$ its support system.
\begin{enumerate}
\item
There exists a
multivalued perturbation $\widetriangle{\frak s} = \{\frak s^n_{\frak p}\}$
of $({\widetriangle{\mathcal U}},\mathcal K)$ such that
each branch of $\frak s^n_{\frak p}$ is transversal to $0$.
\item
Suppose $\widetriangle f : (X,Z;{\widetriangle{\mathcal U}}) \to N$ is strongly
smooth and is transversal to $g : M \to N$, where $g$ is a map from
a manifold $M$.
Then we may choose $\widetriangle{\frak s}$ such that $\widetriangle f$
is strongly transversal to $g$ with respect to $\widetriangle{\frak s}$.
\end{enumerate}}
\par\medskip
A technical but nontrivial result we prove in  Section \ref{sec:multisection} is compactness of the
zero set of a multi-valued perturbation, which is a part of Corollary \ref{cor69}.
Corollary \ref{cor69} also claims the fact that the zero set of multi-valued perturbation converges to the
zero set of the Kuranishi map in Hausdorff topology as our perturbation converges to the Kuranishi map.
(We remark that the zero set of the Kuranishi map is nothing but the space $X$ itself
on which we define our good coordinate system.)
These and other related results play an important role to work out the technical detail of the proof of well-defined-ness of the
virtual fundamental chain (Propositions \ref{indepofukuracont} and \ref{prop14777}).
\par
We remark that in Section \ref{sec:multisection} we work  on a space with good coordinate system and not
on a space with Kuranishi structure.
As we mentioned before, we use an induction over the charts in the construction of
a transversal multisection. This induction works with the good coordinate system but not with
the Kuranishi structure itself.
\par
In Section \ref{sec:contfamily}, we define the notion of continuous family perturbation
(abbreviated as  CF-perturbation) and study their properties in relation to the good
coordinate systems and the integration over the fiber. We use this notion to define the integration of
differential forms on the space with good coordinate system.
(Moreover we also define the integration along the fiber.)
\par
The framework with which we define the integration along the fiber of differential form is as follows.
(See also Situation \ref{smoothcorr}.)
\begin{shitu}\label{situ13}
We consider $X$, $M$, $\widetriangle{\mathcal U}$, $\widetriangle{f}$ with the following properties.
\begin{enumerate}
\item
We are given a good coordinate system $\widetriangle{\mathcal U}
= \{(U_{\frak p},E_{\frak p},s_{\frak p},\psi_{\frak p})\}$ of a space $X$.
\item
We are given submersions $f_{\frak p} : U_{\frak p} \to M$ to a manifold $M$ such that they are compatible
with the coordinate change in an appropriate sense. (Namely we assume that there is a weakly submersive map
$\widetriangle{f}  = \{f_{\frak p} \mid \frak p \in \frak P\}
: (X,\widetriangle{\mathcal U}) \to M$.) (Definition \ref{mapkura}.)$\blacksquare$
\end{enumerate}
\end{shitu}
We next define the notion of continuous family perturbation (CF-perturbation).
We do so in 2 steps.
We first define such notion on a single chart $(U_{\frak p},E_{\frak p},s_{\frak p},\psi_{\frak p})$.
We then discuss its compatibility with coordinate change and use it to define the
notion of CF-perturbation of a good coordinate system.
We define the notion of differential form on the space with good coordinate system.
(It assigns a differential form on $U_{\frak p}$ to each $\frak p \in \frak P$ which are
compatible with coordinate changes. See Definition \ref{defndiffformgcs}.)
We use them to define integration along the fiber
\begin{equation}\label{intalongfiber}
\widetriangle{h} \mapsto \widetriangle f!(\widetriangle{h};\widetriangle{\frak S^{\epsilon}}) \in \Omega^d(M),
\end{equation}
for any sufficiently small $\epsilon > 0$.
Here the degree $d$ is  $\deg \widetriangle{h} + \dim M - \dim (X,\widetriangle{\mathcal U})$
and $\widetriangle{\frak S}$ is a CF-perturbation, which satisfies an appropriate transversality
assumption. (More precisely we assume that $\widetriangle{f}$ is strongly submersive with
respect to $\widetriangle{\frak S}$. See Definition \ref{smoothfunctiononvertK}.)
\par
In the case when $(X,\widetriangle{\mathcal U})$ is a manifold or an orbifold (that is,
when all the obstruction bundles $E_{\frak p}$ are trivial), the operation
(\ref{intalongfiber}) reduces to the standard integration along the fiber of a differential form.
Note $\widetriangle{\frak S}$ is a one-parameter family of perturbations
parameterized by $\epsilon >0$.
The integration along the fiber {\it does} depend on $\epsilon$ as well as CF-perturbation.
We  also remark that typically (\ref{intalongfiber}) diverges as $\epsilon$ goes to $0$.
We firmly believe there is {\it no} way of defining the integration along the fiber
in the way independent of the choice of CF-perturbation.
This is related to the following most basic point of the whole story of
virtual fundamental chains:
\par
In the case we need to construct a virtual fundamental \emph{chain} but \emph{not a cycle},
that is, as in the case when our (moduli) space has a boundary or a corner, the virtual fundamental chain
depends on the choice of perturbations. However we can make (many) choices in a consistent way so that
the resulting algebraic system is independent of such choices modulo certain homotopy equivalence.
\par
We will discuss this point further in Part 2.
\par
On the other hand, the integration along the fiber (\ref{intalongfiber}) is independent of various
other choices involved. Especially it is independent of the choice of partition of unity we use to define
the integration.
We define the partition of unity in Definition \ref{pounity} in the current context and prove its existence
in Proposition \ref{pounitexi}.
The above mentioned independence is proved as Proposition \ref{indepofukuracont}.
\par
In Section \ref{sec:stokes} we prove Stokes' formula for the integration along the fiber
(\ref{intalongfiber}).
We begin with discussing the boundary of a space with Kuranishi structure or with good coordinate system.
We can define the notion of a boundary $\partial M$ of a manifold with corners $M$,
which we call the normalized boundary. $\partial M$ is again a manifold with corner and
there is a map  $\partial M \to M$ which is generically one to one and is a surjection to the boundary of $M$.
Set theoretically $\partial M$ is not a subset of $M$. For example, codimension 2 corner points of $M$
appear twice in $\partial M$. These issues are not deep and are basically well-known. Because a
systematic discussion of these issues is not easy to find and also because they are needed for
a systematic study of Kuranishi structure with corner, which is important for the chain level argument (especially those appearing in Part 2), we include a systematic discussion of these issues here.
Once these points are understood for the case of manifolds, it is straightforward to generalize them
to the Kuranishi structure or to the good coordinate system with boundary and corners.
\par
Stokes' formula in this context is (\ref{stokesinintro}) given below.
In Situation \ref{situ13}, suppose $(X,\widetriangle{\mathcal U})$  has normalized boundary,
$\partial(X,\widetriangle{\mathcal U}) = (\partial X,\partial \widetriangle{\mathcal U})$.
We assume the restriction $\widetriangle{f_{\partial}}$ of $\widetriangle{f}$ to $(\partial X,\partial \widetriangle{\mathcal U})$ is
still weakly submersive.\footnote{Actually this is automatic from our definition, since the submersivity from manifolds
with corners implies that its restriction to each corner (of any fixed codimension)
is a submersion, by definition.}
Suppose also we are given a CF-perturbation $\widetriangle{\frak S}$ of
$(X,\widetriangle{\mathcal U})$. We assume $\widetriangle{f}$ is strongly submersive
with respect to  $\widetriangle{\frak S}$.
Then  $\widetriangle{\frak S}$ induces a CF-perturbation
$\widetriangle{\frak S_{\partial}}$ on the boundary with respect to which
the restriction $\widetriangle{f_{\partial}}$ of $\widetriangle{f}$ is strongly submersive.
Stokes' formula for a good coordinate system now is stated as Theorem \ref{Stokes}:
For any sufficiently small $\epsilon > 0$ we have
\begin{equation}\label{stokesinintro}
d\left(\widetriangle f!(\widetriangle h;\widetriangle{{\frak S}^{\epsilon}})\right)
=
\widetriangle f!(d\widetriangle h;\widetriangle{{\frak S}^{\epsilon}})
+
\widetriangle f_{\partial}!(\widetriangle {h_{\partial}};\widetriangle{{\frak S}_{\partial}^{\epsilon}}).
\end{equation}
This is proved in Subsection \ref{subsec:Stokesgcs}.
We note that in case $M$ (the target of $\widetriangle f$) is a point, the integration
along the fiber (\ref{intalongfiber}) is nothing but the integration of differential
form and is a real number. In that case we write
$$
\int_{(X,\widetriangle{\mathcal U},\widetriangle{\frak S^{\epsilon}})}h
= \widetriangle f!(\widetriangle h;\widetriangle{{\frak S}^{\epsilon}}).
$$
Then (\ref{stokesinintro}) becomes
\begin{equation}\label{formstakessimple}
\int_{(X,\widetriangle{\mathcal U},\widetriangle{\frak S^{\epsilon}})}dh
=
\int_{\partial(X,\widetriangle{\mathcal U},\widetriangle{\frak S^{\epsilon}})}h.
\end{equation}
If all the obstruction bundles are trivial, (\ref{formstakessimple}) is nothing but the usual Stokes' formula.
\par
In Section \ref{sec:stokes}, we also include one easy application of Stokes' formula.
Namely we prove that if $(X,\widetriangle{\mathcal U})$ is equipped with a good coordinate system {\it without boundary}
and $\widetriangle{f}$ is as in Situation \ref{situ13}, then push-forward of $1$ (which is
a differential $0$-form of $(X,\widetriangle{\mathcal U})$) gives rise to a smooth differential form on M
$$
\widetriangle f!(1;\widetriangle{{\frak S}^{\epsilon}}) \in \Omega^{\dim M - \dim (X,\widetriangle{\mathcal U})}(M).
$$
Then (\ref{stokesinintro}) implies that this form is closed. Moreover its de Rham cohomology class is
independent of the choices (of $\widetriangle{{\frak S}}$ and $\epsilon$ also).
This is a consequence of Propositions
\ref{relextendgood},
\ref{cobordisminvsmoothcor} below.
To state them we introduce the
notion of {\it smooth correspondence}
under the following situation.
Note that the assumptions on the weak submersivity as above (also see Situation \ref{situ71} below) are satisfied,
for example, for the case of Gromov-Witten invariant.
So we can use this result to prove well-defined-ness of Gromov-Witten invariant.
\begin{shitu}\label{situ71}
[Situation \ref{smoothcorr}]
Let $X$ be a compact metrizable space,
and
$\widehat{\mathcal U}$ a Kuranishi structure of $X$ (with or without boundaries or corners).
Let $M_s$ and $M_t$ be $C^{\infty}$ manifolds.
We assume $\widehat{\mathcal U}$, $M_s$ and $M_t$
are oriented\footnote{In certain situations, for example in \cite[Subsection 8.8]{fooobook2},
we discussed slightly more general case. Namely we discussed
the case when $\widehat{\mathcal U}$, $M_s$ and $M_t$ are not
necessarily orientable by introducing appropriate
$\Z_2$ local systems. See
\cite[Section A2]{fooobook2} or an appendix to Part 2.}.
\par
Let $\widehat f_s : (X;\widehat{\mathcal U}) \to M_s$ be a strongly smooth map and
$\widehat  f_t : (X;\widehat{\mathcal U}) \to M_t$ a weakly submersive strongly smooth map.
We call
$\frak X = ((X;\widehat{\mathcal U});\widehat f_s,\widehat f_t)$
a
{\it smooth correspondence}
\index{smooth correspondence ! smooth correspondence} from $M_s$ to $M_t$.
$\blacksquare$
\end{shitu}
By Theorem \ref{Them71restate} for $Z=X$, we have a good coordinate system
$(X, \widetriangle{\mathcal U})$ compatible with $(X,\widehat{\mathcal U})$.
Moreover we have $\widetriangle f_s : (X;{\widetriangle{\mathcal U}})
\to M_s$ and $\widetriangle f_t : (X;{\widetriangle{\mathcal U}})
\to M_t$ such that
$\widehat f_s$, $\widehat f_t$
are pullbacks of $\widetriangle f_s$, $\widetriangle f_t$ respectively, and
$\widetriangle f_t$ is weakly submersive
(Proposition \ref{le614} (2)).
Thus
we have the correspondence denoted by
$$
{\frak X} = ((X;\widetriangle{\mathcal U});\widetriangle f_s,
\widetriangle f_t).
$$
\par
We take a CF-perturbation
$\widetriangle{\frak S}$ of $(X;{\widetriangle{\mathcal U}})$
such that $\widetriangle f_t$ is strongly submersive
with respect to ${\widetriangle{\mathcal U}}$
(Theorem \ref{existperturbcont} (2)).
Then we define a map
\begin{equation}\label{Into:corr}
{\rm Corr}_{(\frak X,\widetriangle{{\frak S}^{\epsilon}})}
:
\Omega^k(M_s) \ni h
\to
(\widetriangle{f_t})!
(\widetriangle{f_s^*}h;\widetriangle{\frak S^{\epsilon}})
\in \Omega^{\ell+k}(M_t)
\end{equation}
(Definition \ref{defn748}), which we call the
{\it smooth correspondence map
associated to
$\frak X = ((X;\widetriangle{\mathcal U});\widetriangle f_s,\widetriangle f_t)$}.
Here
$\ell = \dim M_t -\dim (X;\widetriangle{\mathcal U})$.
Then Stokes formula \eqref{stokesinintro}
yields that for any sufficiently small $\epsilon >0$ we have
\begin{equation}
d \circ {\rm Corr}_{(\frak X,\widetriangle{{\frak S}^{\epsilon}})}
-
{\rm Corr}_{(\frak X,\widetriangle{{\frak S}^{\epsilon}})}
\circ d
=
{\rm Corr}_{(\partial\frak X ,\widetriangle{\frak S^{\partial,\epsilon}})}
\end{equation}
(Corollary \ref{Stokescorollary}).
Using this formula,
we show the following propositions.
\par\medskip\noindent
{\bf Proposition \ref{relextendgood}.}
{\it Consider Situation \ref{situ71}
and assume that our Kuranishi structure on $X$ has
no boundary.
Then the map
$
{\rm Corr}_{(\frak X,\widetriangle{{\frak S}^{\epsilon}})}
:
\Omega^k(M_s) \to \Omega^{\ell+k}(M_t)
$
defined above
is a chain map.
\par
Moreover, provided $\epsilon$ is sufficiently small, the map ${\rm Corr}_{(\frak X,\widetriangle{{\frak S}^{\epsilon}})}$
is independent of the choices of
our good coordinate system ${\widetriangle{\mathcal U}}$
and CF-perturbation $\widetriangle{\frak S}$
and of $\epsilon>0$,
up to chain homotopy.}
\par\medskip\noindent
Thus in the situation of Proposition \ref{relextendgood},
the correspondence map
${\rm Corr}_{(\frak X,\widetriangle{{\frak S}^{\epsilon}})}$
on differential forms descends to
a map on {\it cohomology} which is independent of
the choices of $\widetriangle{\mathcal U}$ and
$\widetriangle{{\frak S}^{\epsilon}}$.
We write the cohomology class as
$[{\rm Corr}_{\frak X}(h)] \in H(M_t)$ for
any closed differential form $h$ on $M_s$
by removing $\widetriangle{{\frak S}^{\epsilon}}$ from the notation.
\par\medskip\noindent
{\bf Proposition \ref{cobordisminvsmoothcor}.}
{\it Let $\frak X_i = ((X_i,\widehat{\mathcal U^i}),\widehat f_s^i,\widehat f_t^i)$ be smooth correspondence from $M_s$
to $M_t$ such that $\partial X_i = \emptyset$. Here $i=1,2$ and $M_s$, $M_t$ are independent of $i$.
We assume that there exists a smooth correspondence
$\frak Y = ((Y,\widehat{\mathcal U}),\widehat f_s,\widehat f_t)$ from $M_s$ to $M_t$ with boundary (but without corner) such that
$$
\partial \frak Y = \frak X_1 \cup -\frak X_2.
$$
Here $-\frak X_2$ is the smooth correspondence $\frak X_2$ with
opposite orientation.
Then we have
\begin{equation}\label{chomotopyrelation22}
[{\rm Corr}_{\frak X_1}(h)] = [{\rm Corr}_{\frak X_2}(h)]
\in H(M_t),
\end{equation}
where $h$ is a closed differential form on $M_s$.}
\par\medskip
Summing up these propositions,
if the Kuranishi structure of $X$ has {\it no boundary},
the virtual fundamental cycle of $X$ and the smooth correspondence map
are well-defined on {\it cohomology level}.
These are proved in Section \ref{sec:stokes}.
\par\smallskip
Besides Stokes' formula, an important property we use
for the integration along the fiber is the composition formula.
To formulate a composition formula we need to study
the fiber product of CF-perturbations.
Since a fiber product is to be better defined for the Kuranishi structure than
for the good coordinate system, we rewrite the story of CF-perturbation and integration along the
fiber with respect to the good coordinate system into the one with respect to the Kuranishi structure.
This is the content of Section \ref{sec:kuraandgood}.
The results addressed in Section \ref{sec:kuraandgood} are
really necessary for the {\it chain level argument}, in particular in the later Part of this manuscript.
\par
We can define the notion of CF-perturbation
of a Kuranishi structure in the same way as that of a good coordinate system.
Namely it assigns a CF-perturbation to each
chart $\mathcal U_p$ of $p \in X$ so that they are compatible with coordinate changes.
(See Definition \ref{defn81}.)
Note however it is very difficult to construct a CF-perturbation
with appropriate transversality property on a given
Kuranishi structure $\widehat{\mathcal U}$ of $X$,
since the proof should be by induction on charts as we mentioned several times already.
So we first construct a good coordinate system
$\widetriangle{\mathcal U}$ compatible to a given Kuranishi chart $\widehat{\mathcal U}$
and a CF-perturbation $\widetriangle{\frak S}$ on the  good coordinate system
$\widetriangle{\mathcal U}$  we obtained.
Then we construct another Kuranishi structure $\widehat{\mathcal U^{+}}$
in such a way that $\widetriangle{\frak S}$ induces a
CF-perturbation $\widehat{\frak S^+}$ of $(X,\widehat{\mathcal U^{+}})$.
(Lemma \ref{lemappgcstoKucont}.) Thus in place of constructing a CF-perturbation
of given $(X,\widehat{\mathcal U})$, we construct one
of its thickening $(X,\widehat{\mathcal U^{+}})$.
\par
In this way we have arrived in the situation where we are given a
Kuranishi structure equipped with a CF-perturbation satisfying
appropriate transversality property needed. We formulate this situation as follows.
\begin{shitu}\label{situ155}
We consider $X$, $M$, $\widehat{\mathcal U}$, $\widehat{f}$,
$\widehat{\frak S}$ with the following properties.
\begin{enumerate}
\item
We are given a Kuranishi structure $\widehat{\mathcal U} = \{(U_{p},E_{p},s_{p},\psi_{p})\}$ of a
space $X$.
\item
We are given submersions $f_{p} : U_{p} \to M$ to a manifold $M$ such that they are compatible
with coordinate change in an appropriate sense. (Namely, we assume that there is a weakly submersive map
$\widehat{f}  = \{f_{p} \mid p \in X\}
: (X,\widehat{\mathcal U}) \to M$.) (Definition \ref{mapkura}.)
\item
$\widehat{\frak S}$ is a CF-perturbation
such that $\widehat f$ is strongly transversal to $0$
with respect to $\widehat{\frak S}$.
$\blacksquare$
\end{enumerate}
\end{shitu}
We can define the notion of differential forms on the space with Kuranishi structure
in the same way as on the space equipped with a good coordinate system.
Now we define the integration along the fiber
\begin{equation}\label{intalongfiber2}
\widehat{h} \mapsto \widehat f!(\widehat{h};\widehat{\frak S^{\epsilon}}) \in \Omega^d(M)
\end{equation}
as follows.
We take a good coordinate system
$\widetriangle{\mathcal U}$ on which $\widehat{\frak S^{\epsilon}}$, $\widehat{h}$,
$\widehat{f}$ induce corresponding objects
 $\widetriangle{\frak S^{\epsilon}}$, $\widetriangle{h}$,
$\widetriangle{f}$.
Then by definition
\begin{equation}\label{formula1.5}
 \widehat f!(\widehat{h};\widehat{\frak S^{\epsilon}})
 =
\widetriangle f!(\widetriangle{h};\widetriangle{\frak S^{\epsilon}}).
\end{equation}
The main result of Section \ref{sec:kuraandgood} claims
that the right hand side of (\ref{formula1.5}) is independent of the choice of
$\widetriangle{\mathcal U}$, $\widetriangle{\frak S^{\epsilon}}$, $\widetriangle{h}$,
$\widetriangle{f}$ but depend only on $\widehat{\mathcal U}$, $\widehat{\frak S^{\epsilon}}$, $\widehat{h}$,
$\widehat{f}$. This is Theorem \ref{theorem915}.
(We note that integration along the fiber is hard to define
directly with Kuranishi structure since we use a partition of
unity to define the integration.)
\par
Stokes' formula with respect to the good coordinate system is  easily translated to one for
the Kuranishi structure. Namely we have (Proposition \ref{Stokeskura})
\begin{equation}\label{stokesinintro2}
d\left(\widehat f!(\widehat h;\widehat{{\frak S}^{\epsilon}})\right)
=
\widehat f!(d\widehat h;\widehat{{\frak S}^{\epsilon}})
+
\widehat f_{\partial}!(\widehat {h_{\partial}};\widehat{{\frak S}_{\partial}^{\epsilon}}).
\end{equation}
\par
In Section \ref{sec:composition}, we state and prove a composition formula whose
outline is in order.
We begin with $X,M_s,M_t,\widehat{\mathcal U},\widehat f_s, \widehat f_t, \widehat{\frak S}$
satisfying the following properties.
(Here the indices $s$ and $t$ stand for the source and the target, respectively.)
\begin{shitu}\label{situ1665}
\begin{enumerate}
\item
$X,M_t,\widehat{\mathcal U},\widehat{f_t}$ play the role of
$X$, $M$, $\widehat{\mathcal U}$, $\widehat{f}$
in Situation \ref{situ155}.
\item
$X,M_s,\widehat{\mathcal U},\widehat{f_s}$ play the role of
$X$, $M$, $\widehat{\mathcal U}$, $\widehat{f}$
of Situation \ref{situ155} except we do not assume weak submersivity of
$\widehat f_s$ but assume only strong smoothness.
\item
$\widehat{\frak S}$ is a CF-perturbation
such that $X,M_t,\widehat{\mathcal U},\widehat{f_t}$ together with
$\widehat{\frak S}$ satisfy the transversality
assumption required to define
integration along the fiber by $\widehat{f_t}$.
(Namely we assume $\widehat f_t$ is strongly submersive with respect to
$\widehat{\frak S}$. See Definition \ref{defn929292}.)$\blacksquare$
\end{enumerate}
\end{shitu}
We call $(X,M_s,M_t,\widehat{\mathcal U},\widehat{f_s}, \widehat{f_t})$ a \emph{smooth
correspondence} (See Situation \ref{smoothcorr}) and
$(X,M_s,M_t,\widehat{\mathcal U},\widehat{\frak S},\widehat{f_s}, \widehat{f_t})$
a \emph{perturbed smooth correspondence}
(Definition \ref{defn839}).
To each such correspondence $\frak X = (X,M_s,M_t,\widehat{\mathcal U},\widehat{f_s}, \widehat{f_t},
\widehat{\frak S})$ and sufficiently small $\epsilon > 0$,
we associate a linear map
$$
{\rm Corr}^{\epsilon}_{\frak X} : \Omega^*(M_s) \to \Omega^{*+d}(M_t)
$$
by
\begin{equation}\label{correspondintro}
{\rm Corr}^{\epsilon}_{\frak X}(h)
=  \widehat{f_t}!(\widehat f_s^*\widehat{h};\widehat{\frak S^{\epsilon}}).
\end{equation}
(See Definition \ref{def92111}.)
Here $\widehat{f_s^*}$ is the pull-back operation which assigns a differential form on
$(X,\widehat{\mathcal U})$ to a differential form on $M_s$.
The pull-back is defined for an arbitrary strongly smooth map $\widehat{f_{s}}$.
(See Definition \ref{defn75555}. Here we do not need to assume weak or strong
submersivity.) We define an integer $d$ to be
$$
d = \dim M_t -  \dim (X,\widehat{\mathcal U}).
$$
Our composition formula claims that the assignment $\frak X \mapsto {\rm Corr}^{\epsilon}_{\frak X}$
is compatible with compositions.
\par
Suppose that both
\begin{eqnarray*}
\frak X_{21} & = & (X_{21},M_1,M_2,\widehat{\mathcal U_{21}},
\widehat{\frak S_{21}},\widehat{f_{1,21}}, \widehat{f_{2,21}})\\
\frak X_{32} & = & (X_{32},M_2,M_3,\widehat{\mathcal U_{32}},
\widehat{\frak S_{32}},
\widehat{f_{2,32}}, \widehat{f_{3,32}})
\end{eqnarray*}
satisfy Situation \ref{situ1665}. We define their composition
$\frak X_{31} = \frak X_{32} \circ \frak X_{21}$ as follows.
\par
First we consider the fiber product of Kuranishi structures
$$
(X_{31},\widehat{\mathcal U_{31}}) = (X_{21},\widehat{\mathcal U_{21}}) \times_{M_2} (X_{32},\widehat{\mathcal U_{32}})
$$
to define the space with Kuranishi structure $(X_{31},\widehat{\mathcal U_{31}})$ on which
the CF-perturbations $\widehat{\frak S_{21}}$ and $\widehat{\frak S_{32}}$
induce a CF-perturbation. (Definition \ref{defn837}.)
$\widehat{f_{1,21}}$ induces $\widehat f_{1,31} : (X_{31},\widehat{\mathcal U_{31}}) \to M_1$
and
$\widehat{f_{3,32}}$ induces $\widehat f_{3,31} : (X_{31},\widehat{\mathcal U_{31}}) \to M_3$.
We then put
$$
\frak X_{31} = \frak X_{32} \circ \frak X_{21}
=
(X_{31},M_1,M_3,\widehat{\mathcal U_{31}},
\widehat{\frak S_{31}},
\widehat{f_{1,31}}, \widehat{f_{3,31}}).
$$
It satisfies Situation \ref{situ1665}.
(Lemma \ref{lem838} etc.)
Now the composition formula is stated as
\par\medskip\noindent
{\bf Theorem \ref{compformulaprof}.}
{\it Suppose that $\tilde{\frak X}_{i+1 i} =
(X_{i+1 i},\widehat{\mathcal U_{i+1 i}},\widehat{\frak S_{i+1 i}},
\widehat{f_{i,i+1 i}},
\widehat{f_{i+1,i+1 i}})$
are perturbed smooth correspondences for $i=1,2$.
Then
\begin{equation}\label{formula814}
{\rm Corr}^{\epsilon}_{\tilde{\frak X}_{32}\circ\tilde{\frak X}_{21}}
=
{\rm Corr}^{\epsilon}_{\tilde{\frak X}_{32}} \circ {\rm Corr}^{\epsilon}_{\tilde{\frak X}_{21}}
\end{equation}
for each sufficiently small $\epsilon >0$.}
\par\medskip
\begin{rem}
Integration along the fiber of differential form on
spaces with Kuranishi structure is written in \cite[Section 12]{fooo09}.
Especially Stokes' formula and the Composition formula was given in
\cite[Lemma 12.13]{fooo09} and \cite[Lemma 12.15]{fooo09} respectively.
Here we present them in greater details.
In \cite{fooo09}, the process going from a Kuranishi structure
to a good coordinate system and back was not written explicitly.
(One reason is because the main focus of \cite{fooo09} lies in its application
to the Lagrangian Floer theory of torus orbits of toric manifolds
but not in the foundation of the general theory.)
Here we provide thorough detail.
This theory is actually very similar to the manifold theory.
\end{rem}
\par
In Section \ref{sec:contgoodcoordinate},
we prove the existence of a good coordinate system that is compatible
with the given Kuranishi structure (Theorem \ref{Them71restate}).
We also prove its several variations.
The proof we give there is basically the same as those presented in
\cite{FO} which itself is more detailed in \cite{foootech}.
We simplify the proof in several places as well as provide details of several
points whose proofs were rather sketchy in the previous writings.
We separate the discussion on general topology issue from other parts
and put it in a separate paper \cite{foooshrink}.
\par
Section \ref{sec:contfamilyconstr} is devoted to the proof of existence of a CF-perturbation
satisfying appropriate transversality properties.
This proof is split into 3 parts.
\par
In the first part (Subsection \ref{subsec:confapersingle})
we prove such an existence result for a single Kuranishi chart.
For this, we use the language of sheaf. In particular we prove that the assignment
$$
U \mapsto \{\text{all CF-perturbations on $U$}\}
$$
for each open subset $U \subset U_{\frak p}$ defines a sheaf.
We also consider the sub-sheaf consisting of CF-perturbations
satisfying appropriate transversality properties. The main result is stated as a softness of these sheaves.
(Proposition \ref{prop123123}.)
\par
The second part of the proof (Subsection \ref{subsec:extembandcfp})
discusses the case where we have a coordinate change of Kuranishi chart from $\mathcal U_1$
to $\mathcal U_2$. Assuming we are given a CF-perturbation on $\mathcal U_1$ with
various transversality properties, we show that we can find a
CF-perturbation of $\mathcal U_2$ with certain transversality properties such that
these two CF-perturbations are compatible with the coordinate change
(Proposition \ref{prop1221}).
(Actually Proposition \ref{prop1221} includes a relative version of this statement,
 which we also need for our exposition.)
\par
In the third part (Subsection \ref{subsec:cfpgoodcsys}) we complete the proof of existence
theorem of a CF-perturbation (Theorem \ref{existperturbcont} and its variants)
combining the results of the earlier two subsections.
\par
 As far as the de Rham version is concerned, the results up to Section \ref{sec:contfamilyconstr}
 provide a package we need for the case we work on a single space with Kuranishi structure.
 The contents of Sections \ref{sec:constrsec} and \ref{sec:onezerodim} will be used in
 Parts 2 and 3 to verify that most part of the story works when the ground field is $\Q$.
 The results of Sections \ref{sec:constrsec} and \ref{sec:onezerodim}
 are not necessary over the ground field $\R$ or $\C$.
 Actually de Rham version of the story is in various sense easier to work out than proving
 the corresponding results in the singular homology version.
 (Maybe the reason is similar to the reason why teaching an under-graduate homology theory of manifold
 based on de Rham theory is easier than teaching on based on the singular homology theory.
 This is especially so when one teaches cup product.)
 \par
 In Section \ref{sec:constrsec}, we prove an existence theorem of multi-valued perturbation
(See Theorem  \ref{prop621}). The proof is similar to the proof of Theorem \ref{existperturbcont}.
There are two differences:
\par
The first difference is as follows. In the case we work with the de Rham complex,
we can construct a CF-perturbation of $(X,\widetriangle{\mathcal U})$
so that a given weakly submersive map $\widetriangle f : (X,\widetriangle{\mathcal U}) \to M$
becomes strongly submersive and use it. (This means that the restriction of $f_{\frak p}$ to the zero
set of perturbed section is a submersion.)
This makes it possible to define the integration along the fiber.
On the other hand we cannot expect to find multi-valued perturbation
such that restriction of $f_{\frak p}$ to its zero set is submersive.
(We can find a multi-valued perturbation that is transversal to $0$.)
This difference makes it a bit harder to use multi-valued perturbation
to work out chain level argument.
We can still do it but we need to break symmetry more.
This point is explained, for example, in
\cite[Subsection 7.2.2]{fooobook2},
\cite{fooo091} and etc..
In this article we use multi-valued perturbation only in the case when the dimension of our
good coordinate system is 1,0 or negative because of this issue.
\par
The second difference is rather technical and is explained in Subsection \ref{subsec:nastyreason}.
\par
In Sections \ref{sec:onezerodim}, we discuss the virtual fundamental chain (over $\Q$)
through a multi-valued perturbation. The de Rham version of this section is
Sections \ref{sec:contfamily} and \ref{sec:stokes}.
In this section, we study only the case where the (virtual) dimension of our space
$(X,\widetriangle{\mathcal U})$
with good coordinate system is negative, $0$ or $1$ and prove the following.
\begin{enumerate}
\item
In case $\dim (X,\widetriangle{\mathcal U}) < 0$ the zero set of multisection which is
transversal to $0$ is an empty set.
\item
In case $\dim (X,\widetriangle{\mathcal U}) = 0$ the zero set of multisection which is
transversal to $0$ consists of
finitely many points and is away from the boundary of $X$.
By defining appropriate weight to each point of this zero set and
taking weighted sum, we can define a rational number
which is the virtual fundamental chain.
(This number in general depends on the choice of multisection.)
\item
Suppose $\dim (X,\widetriangle{\mathcal U}) = 1$.
Then  $\dim \partial(X,\widetriangle{\mathcal U}) = 0$.
By Item (2) we can define a virtual fundamental chain
$[\partial(X,\widetriangle{\mathcal U},\widetriangle{\frak s^{n}})]$
of
$\partial(X,\widetriangle{\mathcal U})$.
(Here $\widetriangle{\frak s^{n}}$ is the multi-valued perturbation we use to define it.)
Then we have:
\begin{equation}
[\partial(X,\widetriangle{\mathcal U},\widetriangle{\frak s^{n}})] = 0.
\end{equation}
\end{enumerate}
Item (1) is Lemma \ref{lem1311111} (1).
Item (2) is Lemma \ref{lem1311111} (2),
Lemma \ref{lem134} and Definition \ref{defn1355}.
Item (3) is Theorem \ref{prop13777}, which corresponds to
Stokes' formula in our situation and is the main result
of Section \ref{sec:onezerodim}.
\par\medskip\noindent
{\bf Theorem \ref{prop13777}.}
{\it Let $({\widetriangle{\mathcal U}},\widetriangle{\frak s})$ be a good coordinate system with multivalued perturbation
\index{multivalued perturbation ! good coordinate system with multivalued perturbation}
\index{good coordinate system ! with multivalued perturbation}
of $X$
(see Definition \ref{gcswithperturb}) and assume
$\dim (X,{\widetriangle{\mathcal U}}) = 1$.
We consider its normalized boundary
$\partial(X,{\widetriangle{\mathcal U}}) =
(\partial X,\partial{\widetriangle{\mathcal U}})$
where $\widetriangle{\frak s}$ induces a multivalued perturbation
$\widetriangle{\frak s_{\partial}}$ thereof and
$(\partial{\widetriangle{\mathcal U}},\widetriangle{\frak s_{\partial}^{n}})$
is a good coordinate system with multivalued perturbation of
$\partial X$ with $\dim (\partial X, \partial \widetriangle{\mathcal U})=0$.
Then the following formula holds.
$$
[(\partial X,\partial{\widetriangle{\mathcal U}},\widetriangle{\frak s_{\partial}^{n}})]
=0.
$$}
\par\medskip
Furthermore we also show
\par\medskip\noindent
{\bf Corollary \ref{cobordisminvsmoothcormulti}.}
{\it Let $\frak X_i = (X_i,\widehat{\mathcal U^i})$, $i=1,2$ be spaces with
Kuranishi structure without boundary
of dimension $0$.
Suppose that there exists a space with Kuranishi structure
$\frak Y = (Y,\widehat{\mathcal U})$  (but without corner) such that
$$
\partial \frak Y = \frak X_1 \cup -\frak X_2.
$$
Here $-\frak X_2$ is the smooth correspondence $\frak X_2$ with
opposite orientation.
Then we have
\begin{equation}\label{chomotopyrelation2}
[(X_1,\widehat{\mathcal U^1})]
=
[(X_2,\widehat{\mathcal U^2})].
\end{equation}}
\par\medskip
The proof of Theorem \ref{prop13777}, which was given in \cite{FO}, goes
as follows.
We take a transversal multisection of $(X,\widetriangle{\mathcal U})$
extending the given one, $\widetriangle{\frak s^{n}}$, on the boundary.
We may take it so that its zero set has a triangulation
and it, when equipped with an appropriate weight, defines a chain.
(In our situation it is a singular chain of a space which consists of a single point.)
The boundary of this chain is the degree $0$ chain which is the
virtual fundamental chain $[\partial(X,\widetriangle{\mathcal U},
\widetriangle{\frak s^{n}})]$.
Therefore it is zero.
(Degree zero singular chain of a point is zero if it is homologous to zero.)
\par
This proof is correct as it is.  The only nontrivial issue (besides the existence of
transversal multisection) to be clarified is the triangulability of the zero set.
Nontrivial point of the argument is the following.
If we take each branch of multisection $\frak s_i^{n}$,
its zero set is a one dimensional manifold.
We need to check however union of various $(\frak s_i^{n})^{-1}(0)$
(for different branches) has triangulation.
In case the intersection of  $(\frak s_i^{n})^{-1}(0)$ with $(\frak s_j^{n})^{-1}(0)$
is wild this triangulability fails.
Other point is a discussion of the behavior of the zero set at the locus where the number of branches
changes.
As we mentioned in \cite{FO} we can resolve these points by taking
an appropriate generic choice of perturbations.
(See \cite[page 946]{FO}.)
In our situation where our space (with good coordinate system) is 1 dimensional it is not
so difficult to work it out.\footnote{In \cite{FO} the case when the dimension is higher is
also discussed. Those cases are more nontrivial to handle. However none of the application in
\cite{FO} uses the case when the dimension of $(X,\widetriangle{\mathcal U})$ is higher
than $1$.}
\par
We will discuss triangulation of the zero set of multisection elsewhere in more
detail. In this article, we provide an alternative proof of Theorem \ref{prop13777}
which may be simpler and more transparent.
In this proof we take a function $f : (X,\widetriangle{\mathcal U})  \to \R_{\ge 0}$
such that $f^{-1}(0) = \partial X$ and its gradient vector field of $f$ along $\partial X$
points inward.  For generic $s > 0$, the level set $f^{-1}(s) = X_s$ carries a good coordinate
system induced from $\widetriangle{\mathcal U}$ which we write
 $\widetriangle{\mathcal U_s}$. The dimension of $(X_s,\widetriangle{\mathcal U_s})$ is zero.
Using compactness of $X$ we find that $X_s$ is an empty set for sufficiently large $s$.
Since $(X_0,\widetriangle{\mathcal U_0}) = \partial(X,\widetriangle{\mathcal U})$,
to prove Theorem \ref{prop13777}, it suffices to prove that the virtual fundamental chain
$[(X_s,\widetriangle{\mathcal U_s},\widetriangle{\frak s^{n}}\vert_{\widetriangle{\mathcal U_s}})] \in \Q$ is independent of $s$.
Note the zero set of $(\widetriangle{\frak s^{n}})^{-1}(0) \cap f^{-1}(0)$ is zero dimensional, which is
a finite set. So the required independence follows by locally studying the zero set of $\widetriangle{\frak s^{n}}$.
\begin{rem}\label{joycerem}
We would like to again mention a relationship with
\cite{joyce2}, \cite{joyce4}.
In \cite{joyce4}, Joyce gave an alternative version of `space with Kuranishi structure'.
In his version he relaxes the condition of compatibility of coordinate changes so that
it is required only at the zero set of Kuranishi map. He requires compatibility including the derivative
up to the first order. In that way Joyce succeeds in inverting the arrows of the coordinate
change of Kuranishi charts so that the trouble coming from noninvertibility of
coordinate change disappears. His way has an advantage that one can define the notion of morphisms
between `spaces with Kuranishi structures'.
\par
We however use our version of Kuranishi structure and good coordinate system in this article.
The reason is as follows. Our goal is to define a system of operators from a system of
smooth correspondences (which is the object such as Situation \ref{situ1665} (1)(2)).
We need to choose a chain model on which we realize our operations forming the algebraic structure.
Our choice in this article is the de Rham complex.
(In some other occasion we also use the singular chain complex.)
Note the space $X$ which has Kuranishi structure
may have pathological topology in general. So singular homology does not behave
nicely for $X$. In the case of de Rham model, the
situation gets even worse. Namely it seems impossible to define the
notion of differential forms on $X$.
\footnote{Since \v Cech cohomology behaves nicely with respect to the projective limit, it might be the best choice
if we want to define chain model directly on the topological space $X$.}
\par
Thus, we need to take a union of charts of $X$ which has a positive size
to work with de Rham or singular homology.
Namely we need a system of spaces $U_{\frak p}$ which is
a manifold or orbifold and containing $X$.
Both singular homology and de Rham cohomology of such spaces behave nicely.
We remark that to define the notion of differential forms of
a good coordinate system or of Kuranishi structure,
we need compatibility of coordinate change in our sense, that is stronger than
Joyce assumed in \cite{joyce4}.
Namely, we also need to assume compatibility at some points
outside $X$ (i.e., outside the zero set of Kuranishi map $s$).
By this reason, it seems that we can use neither de Rham cohomology
nor singular homology directly if we use the definition of Kuranishi structure
in the sense of \cite{joyce4}.
\par
As far as we understand, Joyce's plan is to use a version of Kuranishi homology
(\cite{joyce})
as the cohomology theory which makes sense in his version of Kuranishi structure.
It seems likely that this approach works.
One potential trouble however is the Poincar\'e duality.
Joyce in \cite{joyce} provides a chain level intersection paring.
However the intersection `number' in his chain level intersection paring
is not a number but is an element of some huge complex
(whose cohomology group is $\Q$).
Though this construction provides the same amount of information in the
homology level there is a trouble using it for the chain level
argument.
While we work on the chain level argument,
sometimes we need to convert some input variables (of algebraic
operation we will obtain) to output variables.
We use the (chain level) Poincar\'e duality for this purpose.
Note the pairing
$$
(u,v) \mapsto \int_M u \wedge v \in \R
$$
works in the chain level in de Rham theory. This identifies an element of de Rham complex with an element in its dual.
Although the dual space of the space of differential forms
is the set of distributions and is different from the set of differential forms,
the difference between the spaces of differential forms and of distributions is relatively small so that we can
still use the chain level Poincar\'e duality
to convert certain input variables to output variables.
It seems that to realize this Poinca\'e duality in the situation of Kuranishi
homology, one needs to work more on the side of homological algebra.
The amount of algebraic work to be done for this purpose might be very heavy,
although it is plausible.
\end{rem}
\begin{rem}
As we mentioned in the beginning, we previously wrote an article  \cite{foootech} which
provides a detailed explanation of similar nature as that of the present article.
We refer readers to \cite[Part 6]{foootech} for some documentation of certain
activities in which we were involved  concerning the foundation of the virtual
fundamental chain techniques, around the time when \cite{foootech} was written.
After we posted \cite{foootech} in arXiv in 2013 September,
we have continued our effort of accommodating the demand for more details of this technique
which came from some part of symplectic geometry community.

The first named author, together with other mathematicians,
organized a semester-long program in the Simons Center for Geometry and Physics
to discuss the foundation of the virtual fundamental chain techniques.
Two one-week long conferences were held on the subject as well as a series of eleven
lectures are presented by the first named author which
are closely related to the content of this article.
(The video of the conferences and the first named author's lectures are available
in the web page of the Simons Center for Geometry and Physics.)
In addition, in our attempt to clarify the `Hausdorffness issue' raised by
D. McDuff and K. Werheim  in their joint lectures given in the Institute for Advanced Study
in early 2012, we separately wrote a paper \cite{foooshrink}.

While we have been writing this article and during these activities occurred
(that is, between September 2013 and March 2015), the present authors have not been aware of
any explicit mathematical questions unanswered on the foundation of virtual
fundamental chain or cycle technique or on its application to the moduli space of pseudo-holomorphic curves.
\footnote{We disregard `objections' directed to the points whose answers we had already provided
before they were presented.}
\end{rem}

\newpage
\section{Notations}
\label{sec:notations}
\begin{enumerate}
\item[$\circ$]
${\rm Int}\, A$, $\ring A$: Interior of a subset $A$ of a topological space.
\item[$\circ$]
$\overline A$: Closure of a subset $A$ of a topological space.
\item[$\circ$]
${\rm Perm}(k)$: The permutation group of order $k!$.
\item[$\circ$]
${\rm Supp}(h)$, ${\rm Supp}(f)$: The support of a differential form $h$,
a function $f$, etc..
\item[$\circ$]
$\varphi^{\star}\mathscr F$: Pullback of a sheaf $\mathscr F$ by
a map $\varphi$.
\item[$\circ$]
$X$: A paracompact metrizable space. (Part I).
\item[$\circ$]
$Z$: A compact subspace of $X$. (Part I).
\item[$\circ$]
$\mathcal U =(U,\mathcal E,\psi,s)$:
A Kuranishi chart, Definition \ref{defnKchart}.
\item[$\circ$]
$\mathcal U\vert_{U_0} =(U_0,\mathcal E\vert_{U_0},\psi\vert_{U_0\cap s^{-1}(0)},s\vert_{U_0})$:
open subchart of $\mathcal U =(U,\mathcal E,\psi ,s)$, Definition \ref{defnKchart}.
\item[$\circ$]
$\Phi = (\varphi,\widehat\varphi)$:
Embedding of Kuranishi charts, Definition \ref{defKchart}.
\item[$\circ$]
$o_p, o_p(q)$: Points in a Kuranishi neighborhood $U_p$ of $p$. Definition \ref{kuranishineighborhooddef}.
\item[$\circ$]
$\Phi_{21} = (U_{21},\varphi_{21},\widehat\varphi_{21})$:
Coordinate change of Kuranishi charts from $\mathcal U_1$ to $\mathcal U_2$, Definition \ref{coordinatechangedef}.
\item[$\circ$]
$\widehat{\mathcal U} = (\{\mathcal U_p\},\{\Phi_{pq}\})$:
Kuranishi structure, Definition \ref{kstructuredefn}.
\item[$\circ$]
$(X,\widehat{\mathcal U})$, $(X,Z;\widehat{\mathcal U})$:
K-space, relative K-space, Definition \ref{Kspacedef}.
\item[$\circ$]
${\widetriangle{\mathcal U}}
= (({\frak P},\le), \{\mathcal U_{\frak p}\},
\{\Phi_{\frak p\frak q}\})$:
Good coordinate system, Definition \ref{gcsystem}.
\item[$\circ$]
$\vert{\widetriangle{\mathcal U}}\vert$:
Definition \ref{defofveruver}.
\item[$\circ$]
$\widehat\Phi : \widehat{\mathcal U} \to \widehat{\mathcal U'}$:
KK-embedding. An embedding of Kuranishi structures, Definition \ref{defn311}.
\item[$\circ$]
$\widetriangle\Phi : \widetriangle{\mathcal U} \to \widetriangle{\mathcal U'}$
: GG-embedding.
An embedding of good coordinate systems, Definition \ref{defn31222}.
\item[$\circ$]
$\widehat\Phi : \widehat{\mathcal U} \to \widetriangle{\mathcal U}$:
KG-embedding,
An embedding of a Kuranishi structure to a good coordinate system,
Definition \ref{defn32020202}.
\item[$\circ$]
$\widehat\Phi : \widetriangle{\mathcal U} \to\widehat{\mathcal U}$:
GK-embedding.
An embedding of good coordinate system to a Kuranishi structure,
Definition \ref{embgoodtokura}.
\item[$\circ$]
$\widehat f : (X,Z;\widehat{\mathcal U}) \to Y$ and
$\widetriangle f : (X,Z;\widetriangle{\mathcal U}) \to Y$
:
Strongly continuous map, Definitions \ref{mapkura} and
\ref{definition32727}.
\item[$\circ$]
$(X,Z;\widehat{\mathcal U})  \times_{N} M$,
$(X_1,Z_1;\widehat{\mathcal U}_1)  \times_{M}
(X_2,Z_2;\widehat{\mathcal U}_2)$:
Fiber product of Kuranishi structures,
Definition \ref{firberproddukuda}.
\item[$\circ$]
$S_k(X,Z;\widehat{\mathcal U})$, $S_k(X,Z;\widehat{\mathcal U})$:
Corner structure stratification, Definition \ref{dimstratifidef}.
\item[$\circ$]
$\mathcal S_{\frak d}(X,Z;\widehat{\mathcal U})$,
$\mathcal S_{\frak d}(X,Z;\widetriangle{\mathcal U})$:
Dimension stratification,
Definition \ref{stratadim}
\item[$\circ$]
$\widehat{\mathcal U} < \widehat{\mathcal U^+}$:
$\widehat{\mathcal U^+}$ is a thickening of $\widehat{\mathcal U}$.
Definition \ref{thickening}.
\item[$\circ$]
$\mathcal S_{\frak p}(X,Z;{\widetriangle {\mathcal U}};\mathcal K)$:
Definition \ref{situ61} (4).
\item[$\circ$]
$
\mathcal K = \{\mathcal K_{\frak p}\mid {\frak p
\in \frak P}\}$:
A support system. Definition \ref{situ61} (2).
\item[$\circ$]
$(\mathcal K^1,\mathcal K^2)$ or
$(\mathcal K^-,\mathcal K^+)$:
A support pair, Definition \ref{situ61} (3).
\item[$\circ$]
$\mathcal K^1 < \mathcal K^2$: Definition \ref{situ61}.
\item[$\circ$]
$\vert\mathcal K\vert$:  Definition \ref{situ61}.
\item[$\circ$]
$B_{\delta}(A)$: Metric open ball,
(\ref{defmetricball}).
\item[$\circ$]
$\mathcal S_x = (W_x,\omega_x,\{{\frak s}_x^{\epsilon}\})$:
CF-perturbation (=continuous family perturbation) on one orbifold chart. Definition
\ref{defn73ss}.
\item[$\circ$]
$\mathcal S_x^{\epsilon} = (W_x,\omega_x,{\frak s}_x^{\epsilon})$ for
each $\epsilon >0$:
Definition \ref{C0convconti}.
\item[$\circ$]
$\frak S = \{(\frak V_{\frak r},\mathcal S_{\frak r})\mid{\frak r\in
\frak R}\}$:
Representative of a CF-perturbation on Kuranishi chart
$\mathcal U$.
Definition \ref{semiglobalocntpert}.
\par\noindent
Here $\frak V_{\frak r}=(V_{\frak r},E_{\frak r},\Gamma_{\frak r},\phi_{\frak r},\widehat{\phi}_{\frak r})$ is an orbifold chart of $(U,\mathcal E)$ and
$\mathcal S_{\frak r}  = (W_{\frak r} ,\omega_{\frak r}, \{{\frak s}_{\frak r} ^{\epsilon}\})$ is a
CF-perturbation of $\mathcal U$
on $\frak V_{\frak r}$.
\item[$\circ$]
$\frak S^{\epsilon} = \{(\frak V_{\frak r},\mathcal S_{\frak r}^{\epsilon})\mid{\frak r\in
\frak R}\}$
for each $\epsilon >0$. Definition \ref{semiglobalocntpert}.
\item[$\circ$]
$\widetriangle{\frak S} = \{\frak S_{\frak p} \mid \frak p \in \frak P\}$:
CF-perturbation of good coordinate system.
Definition \ref{defn7732}.
\item[$\circ$]
$\widehat{\frak S}$:
CF-perturbation of Kuranishi structure.
Definition \ref{defn81}.
\item[$\circ$]
$\mathscr{S}$: Sheaf of CF-perturbations.
Proposition \ref{prop721}.
\item[$\circ$]
$
\mathscr S_{\pitchfork 0}
$,
$
\mathscr S_{f \pitchfork}
$, $
\mathscr S_{f \pitchfork g}
$:
Subsheaves of $\mathscr S$. Definition \ref{strosubsemiloc}.
\item[$\circ$]
${\widetriangle f}!(\widetriangle h;\widetriangle{{\frak S}^{\epsilon}})$:
Pushout or integration along the fiber of $\widetriangle{h}$ by
$(\widetriangle{f},\widetriangle{{\frak S}^{\epsilon}})$
on good coordinate system.
Definition \ref{pushforwardKuranishi}

\item[$\circ$]
${\widehat f}!(\widehat h;\widehat{{\frak S}^{\epsilon}})$:
Pushout or integration along the fiber of $\widehat{h}$ by
$(\widehat{f},\widehat{{\frak S}^{\epsilon}})$
on Kuranishi structure.
Definition \ref{deflemgg}.
\item[$\circ$]
${\rm Corr}_{(\frak X,\widetriangle{{\frak S}^{\epsilon}})}$:
Smooth correspondence associated to good coordinate system.
Definition \ref{defn748}.
\item[$\circ$]
${\rm Corr}_{(\frak X,\widehat{\frak S^{\epsilon}})}$:
Smooth correspondence of Kuranishi structure
Definition \ref{def92111}.
\item[$\circ$]
$(\frak s_{\frak p}^{n})^{-1}(0)$:
The zero set of multisection.
\item[$\circ$]
$\Pi((\mathfrak S^{\epsilon})^{-1}(0))$:
Support set of a CF-perturbation $\mathfrak S^{\epsilon}$.
Definition \ref{defn767}.
\item[$\circ$]
$(V,\Gamma,\phi)$: Orbifold chart, Definitions \ref{2661}, \ref{defn26550}.
\item[$\circ$]
$(V,E,\Gamma,\phi,\widehat\phi)$: Orbifold chart of
a vector bundle, Definitions \ref{defn2613}, \ref{defn2655}.
\item[$\circ$]
$(X,\mathcal E)$: Orbibundle, Definition \ref{defn2820}.
\end{enumerate}
{\bf Convention on the way to use several
notations.}
\begin{enumerate}
\item
[$\widehat{}$ \, and \,$\widetriangle{}$]
\enskip We use `hat' such as $\widehat{\mathcal U}$,
$\widehat{f}$, $\widehat{\frak S}$, $\widehat{h}$
of an object defined on a Kuranishi structure $\widehat{\mathcal U}$.
We use `triangle' such as $\widetriangle{\mathcal U}$,
$\widetriangle{f}$, $\widetriangle{\frak S}$, $\widetriangle{h}$
of an object defined on a good coordinate system $\widetriangle{\mathcal U}$.
\item[$p$ and $\frak p$]
For a Kuranishi structure $\widehat{\mathcal U}$ on $Z\subseteq X$
we write $\mathcal U_p$ for its Kuranishi chart, where
$p \in Z$. (We use an italic letter $p$.)
For a good coordinate system $\widetriangle{\mathcal U}$ on $Z\subseteq X$
we write $\mathcal U_{\frak p}$ for its Kuranishi chart, where
$\frak p \in \frak P$. (We use a German character $\frak p$.) Here $\frak P$ is
a partial ordered set.
\item
[$\blacksquare$]
\enskip
The mark $\blacksquare$ indicates the end of Situation.
See Situation \ref{opensuborbifoldchart}, for example.
\end{enumerate}
\newpage

\section{Kuranishi structure and good coordinate system}
\label{sec:skuraterm}

\subsection{Kuranish structure}
\label{subsec:kuranishi}

The notion of Kuranishi structure in this document is the same as one
in \cite[Section A1]{fooo09} and \cite{foootech},
except we include the existence of tangent bundle in the definition
of Kuranishi structure.
The notion of good coordinate system in this document
is the same as one in \cite{foootech}.
We introduce some more notations which are useful to
shorten the account of this article.
We refer Section \ref{sec:ofd},
for the definition of (effective) orbifold,
vector bundle on an oribifold, and their embeddings.
Our orbifold is always assumed to be {\it effective} unless otherwise mentioned explicitly.
\par
Throughout Part 1, $X$ is always a separable metrizable space.
\begin{defn}\label{defnKchart}
A {\it Kuranishi chart}
\index{Kuranishi chart ! Kuranishi chart} of $X$ is $\mathcal U =(U,\mathcal E,\psi,s)$ with
the following properties.
\begin{enumerate}
\item
$U$ is an orbifold.
\item
$\mathcal E$ is a vector bundle on $U$.
\item
$s$ is a smooth section of $\mathcal E$.
\item
$\psi : s^{-1}(0) \to X$ is a homeomorphism onto an open set.
\end{enumerate}
We call $U$ a {\it Kuranishi neighborhood}, $\mathcal E$ an
{\it obstruction bundle},
$s$ a {\it Kuranishi map} and $\psi$ a {\it parametrization}.
\index{Kuranishi chart ! Kuranishi neighborhood}
\index{Kuranishi chart ! Kuranishi map}
\index{Kuranishi chart ! obstruction bundle}
\index{Kuranishi chart ! parametrization}

If $U'$ is an open subset of $U$, then by restricting $\mathcal E$, $\psi$ and $s$
to $U'$, we obtain a
Kuranishi chart which we write $\mathcal U\vert_{U'}$
and call it an {\it open subchart}.\index{Kuranishi chart ! open subchart}
\par
The {\it dimension} \index{dimension ! of Kuranishi chart} $\mathcal U =(U,\mathcal E,\psi,s)$
is by definition
$$
\dim \mathcal U = \dim U  - {\rm rank} \mathcal E.
$$
Here ${\rm rank} \mathcal E$ is the dimension  of the fiber $\mathcal E \to U$.
\par
We say that $\mathcal U =(U,\mathcal E,\psi,s)$ is {\it orientable}
\index{orientable !  Kuranishi chart} if $U$ and $E$ are orientable.
An {\it orientation} \index{orientation !  Kuranishi chart}
of $\mathcal U =(U,\mathcal E,\psi,s)$ is a pair of orientations of
$U$ and of $\mathcal E$.
An open subchart of an oriented Kuranishi chart is oriented.
\end{defn}
\begin{defn}\label{defKchart}
Let $\mathcal U = (U,\mathcal E,\psi,s)$, $\mathcal U'
= (U',\mathcal E',\psi',s')$ be Kuranishi charts of $X$.
An {\it embedding} of Kuranishi charts $: \mathcal U\to \mathcal U'$
\index{Kuranishi chart ! embedding of Kuranishi charts}
\index{embedding ! of Kuranishi charts}
is  a pair $\Phi = (\varphi,\widehat\varphi)$ with the following properties.
\begin{enumerate}
\item
$\varphi : U \to U'$ is an embedding of orbifolds.
(See Definition \ref{def262220}.)
\item
$\widehat\varphi : \mathcal E \to \mathcal E'$ is an embedding of vector bundles
over $\varphi$.
(See Definition \ref{defn2820}.)
\item
$\widehat\varphi \circ s = s' \circ \varphi$.
\item
$\psi' \circ \varphi = \psi$ holds on $s^{-1}(0)$.
\item
For each $x \in U$ with $s(x) = 0$, the (covariant) derivative
$
D_{\varphi(x)}s'
$
induces an isomorphism
\begin{equation}\label{form3.1111}
\frac{T_{\varphi(x)}U'}{(D_x\varphi)(T_xU)}
\cong
\frac{\mathcal E'_{\varphi(x)}}{\widehat\varphi(\mathcal E_x)}.
\end{equation}
\end{enumerate}
In other words, the map (\ref{form3.1111}) is the right vertical arrow
of the next commutative diagram.
\begin{equation}\label{diagrampart1333}
\begin{CD}
T_xU
@ > {D_x\varphi} >>
T_{\varphi(x)}U'
@ >>>
\frac{T_{\varphi(x)}U'}{(D_x\varphi)(T_xU)}
\\
@ V{D_xs}VV @ VV{D_{\varphi(x)}s'}V @VVV\\
\widehat\varphi(\mathcal E_x) @>{\widehat{\varphi}}>>\mathcal E'_{\varphi(x)}
@>>>\frac{\mathcal E'_{\varphi(x)}}{\widehat\varphi(\mathcal E_x)}
\end{CD}
\end{equation}
\par
If $\dim U = \dim U'$ in addition, we call $(\varphi,\widehat\varphi)$ an {\it open embedding}.
\index{embedding ! open embedding of Kuranishi chart}
\index{Kuranishi chart ! open embedding of Kuranishi chart}
\end{defn}
\begin{defn}
In the situation of
Definition \ref{defKchart}, suppose $\mathcal U$ and $\mathcal U'$
are oriented. Then the orientations induce trivialization
of
${\rm Det} TU' \otimes {\rm Det}{\mathcal E}'$ and of
${\rm Det} TU \otimes {\rm Det} {\mathcal E}$.
(Here ${\rm Det} {\mathcal E}$ is a real line bundle which is the determinant
line bundle of  $\mathcal E$. ${\rm Det} TU'$ etc. are defined in the same way.)
\par
We say $\Phi = (\varphi,\widehat\varphi)$ is {\it orientation preserving}
\index{orientation preserving ! embedding of Kuranishi chart}
if the isomorphism
$$
{\rm Det}T_{\varphi(x)}U'  \otimes ({\rm Det}T_{x}U)^*
\cong
{\rm Det}\mathcal E'_{\varphi(x)}  \otimes ({\rm Det}\mathcal E_x)^*
$$
induced by (\ref{form3.1111}) is compatible with
these trivializations.
\end{defn}
The composition of embeddings of Kuranishi charts is again an embedding of
Kuranishi charts.
There is an obvious embedding of Kuranishi charts from $\mathcal U$ to itself,
that is, the identity.
We can define the notion of {\it isomorphism} of Kuranishi charts
by using the above two facts in an obvious way.
\index{Kuranishi chart ! isomorphism of Kuranishi chart}
\par
\begin{defn}\label{kuranishineighborhooddef}
For $A \subseteq X$, a {\it Kuranishi neighborhood} of $A$ is a Kuranishi chart
such that ${\rm Im}(\psi)$ contains $A$.
In case $A=\{p\}$ we call it a Kuranishi neighborhood of $p$
or a Kuranishi chart at $p$.
\index{Kuranishi chart ! Kuranishi neighborhood of $A$}
\par
When $\mathcal U_p = (U_p,\mathcal E_p,\psi_p,s_p)$ is a Kuranishi neighborhood of $p$,
we denote by $o_p \in U_p$ the point such that
$s_p(o_p)= 0$ and $\psi_p(o_p) = p$.
If $q \in {\rm Im}(\psi_p)$ we denote by $o_p(q) \in U_p$ the point such that
$s_p(o_p(q))= 0$ and $\psi_p(o_p(q)) = q$.
Note such $o_p$ and $o_p(q)$ are unique.
\end{defn}

\begin{defn}\label{coordinatechangedef}
Let $\mathcal U_1 = (U_1,\mathcal E_1,\psi_1,s_1)$, $\mathcal U_2
= (U_2,\mathcal E_2,\psi_2,s_2)$ be Kuranishi charts of $X$.
A {\it coordinate change in weak sense} (resp. {\it in strong sense}) from $\mathcal U_1$
\index{Kuranishi chart ! coordinate change in weak sense}
\index{Kuranishi chart ! coordinate change in strong sense}
to
$\mathcal U_2$ is $\Phi_{21} = (U_{21},\varphi_{21},\widehat\varphi_{21})$
with the following properties (1) and (2) (resp. (1), (2) and (3)):
\begin{enumerate}
\item
$U_{21}$ is an open subset of $U_1$.
\item
$(\varphi_{21},\widehat\varphi_{21})$ is an
embedding of Kuranishi charts $: \mathcal U_1\vert_{U_{21}}
\to \mathcal U_2$.
\item
$\psi_1(s_1^{-1}(0) \cap U_{21}) = {\rm Im}(\psi_1) \cap   {\rm Im}(\psi_2)$.
\end{enumerate}
In case $\mathcal U_1$ and $\mathcal U_2$ are oriented $\Phi_{21}$
is said to be {\it orientation preserving} \index{orientation preserving ! of
coordinate change} if it is so as an embedding.
\end{defn}

\begin{rem}
We use coordinate changes in weak sense for Kuranishi structures (Definition \ref{kstructuredefn}),
while we use coordinate changes in strong sense for good coordinate systems (Definition \ref{gcsystem}).
From now on, {\it coordinate changes appearing in Kuranishi structures are in weak sense}
and {\it coordinate changes appearing in good coordinate systems are in strong sense.}
\end{rem}
\begin{conven}
Hereafter in Part 1,  $Z$ is assumed to be a compact subset of $X$, unless otherwise specified.
\end{conven}
\begin{defn}\label{kstructuredefn}
A {\it Kuranishi structure} $\widehat{\mathcal U}$ of $Z \subseteq X$
\index{Kuranishi structure ! Kuranishi structure} assigns a Kuranishi neighborhood $\mathcal U_p
= (U_p,\mathcal E_p,\psi_p,s_p)$
of $p$ to each $p \in Z$ and a coordinate change in weak sense
$\Phi_{pq} = (U_{pq},\varphi_{pq},\widehat\varphi_{pq}) : \mathcal U_q \to \mathcal U_p$
to each $p$, $q \in {\rm Im}(\psi_p) \cap Z$ such that
$q \in \psi_q(s_q^{-1}(0) \cap U_{pq})$ and the following holds
for each $p$, $q \in {\rm Im}(\psi_p) \cap Z$, $r \in \psi_q(s_q^{-1}(0) \cap U_{pq}) \cap Z$.
\par
We put
$U_{pqr} = \varphi_{qr}^{-1}(U_{pq}) \cap U_{pr}$. Then we have
\begin{equation}\label{form3333}
\Phi_{pr}\vert_{U_{pqr}} = \Phi_{pq}\circ
\Phi_{qr}\vert_{U_{pqr}}.
\end{equation}
\par
We call $Z$ the {\it support set }
\index{Kuranishi structure ! support set} of our Kuranishi structure.
\par
We also require that the dimension of $\mathcal U_p$ is independent of
$p$ and call it the {\it dimension} of $\widehat{\mathcal U}$.
\index{dimension ! of Kuranishi structure}
\end{defn}
\begin{rem}\label{rem3737}
We require that the equality (\ref{form3333}) holds on
the domain where both sides are defined.
This is always the case of this kinds of equality when we require such an
equality between the maps whose domain is a Kuranishi chart
that is a member of a Kuranishi structure.
\par
On the other hand, in case when we study  maps
whose domain is a Kuranishi chart that is a member of a good coordinate system
defined in Definition \ref{gcsystem}, we sometimes require  other conditions such as the condition that
the domains of the two maps coincide.
(See Definition \ref{defn31222} (1) for example.) We mention explicitly
those conditions when we require it.
\end{rem}
\begin{defn}\label{kuraorient}
We say the Kuranishi structure $(\{\mathcal U_p\},\{\Phi_{pq}\})$
is {\it orientable} \index{orientable ! Kuranishi structure}
if we can choose orientation of $\mathcal U_p$
such that all $\Phi_{pq}$ are orientation preserving.
\par
The notion of orientation of an orientable Kuranishi structure
and of oriented Kuranishi structure  is defined in an obvious way.
\end{defn}
\begin{defn}\label{Kspacedef}
A {\it K-space}  is a pair $(X,\widehat{\mathcal U})$
of a paracompact metrizable space $X$ and a Kuranishi structure
$\widehat{\mathcal U}$ of $X$.
\index{K-space ! K-space}
\par
A {\it relative K-space} is a triple $(X,Z;\widehat{\mathcal U})$,
where $Z \subseteq X$ is a compact subspace and $\widehat{\mathcal U}$
is a Kuranishi structure of $Z \subseteq X$.
\end{defn}
\begin{rem}
In \cite{FO,fooobook2,foootech} we assumed that the orbifold appearing in Kuranishi
structure is a global quotient. Namely we assumed $U = V/\Gamma$ where $V$ is a
manifold and $\Gamma$ is a finite group acting on $V$ effectively and smoothly.
There is no practical difference of the definition since we can always
replace $U_p$ by a smaller open subset so that
it becomes of the form $U_p = V_p/\Gamma_p$.
\end{rem}
\begin{rem}
In \cite{joyce} etc. a space with Kuranishi structure is called
Kuranishi space. However
the name `Kuranishi space' has been used for a long time
for the deformation space of complex structure,
which Kuranishi discovered in his celebrated work. The Kuranishi structure
in our sense is much
inspired by Kuranishi's work, but a space with Kuranishi structure
is different from the deformation space of
complex structure (Kuranishi space). So we call it
K-space in this document.
\par
We also insist to call $s$ a Kuranishi map. This is the
main notion discovered by Kuranishi.
\end{rem}
From now on when we write Kuranishi neighborhood of $p$ as $\mathcal U_p$, $\mathcal U'_p$
etc. we use the notation $V_p,U_p$ etc. by
$
\mathcal U_p =  (U_p,\mathcal E_p,\psi_p,s_p)
$.

\subsection{Good coordinate system}
\label{subsec:goodcoordinate}
\begin{defn}\label{gcsystem}
A {\it good coordinate system} of  $Z \subseteq X$
\index{good coordinate system ! good coordinate system}  is
$${\widetriangle{\mathcal U}}
= (({\frak P},\le), \{\mathcal U_{\frak p}\mid \frak p \in \frak P\},
\{\Phi_{\frak p\frak q}
\mid \frak q \le \frak p\})$$
such that:
\begin{enumerate}
\item
$({\frak P},\le)$ is a partially ordered set.
We assume $\# \frak P$ is finite.
\item
$\mathcal U_{\frak p}$ is a Kuranishi chart of $X$.
\item
$
\bigcup_{\frak p \in \frak P} U_{\frak p} \supseteq Z.
$
\item
Let $\frak q \le \frak p$. Then
$\Phi_{\frak p \frak q} = (U_{\frak p \frak q},\varphi_{\frak p \frak q},\widehat\varphi_{\frak p \frak q})$ is a coordinate change
in strong sense:
$\mathcal U_{\frak q} \to \mathcal U_{\frak p}$
in the sense of Definition \ref{coordinatechangedef}.
\item
If $\frak r \le \frak q \le \frak p$, then by putting
$U_{\frak p\frak q\frak r} = \varphi_{\frak q\frak r}^{-1}(U_{\frak p\frak q}) \cap U_{\frak p\frak r}$ we have
\begin{equation}
\Phi_{\frak p \frak r}\vert_{U_{\frak p \frak q \frak r}} = \Phi_{\frak p \frak q}\circ
\Phi_{\frak q \frak r}\vert_{U_{\frak p \frak q \frak r}}.
\end{equation}
\item
If ${\rm Im}(\psi_{\frak p}) \cap {\rm Im}(\psi_{\frak q})
\ne \emptyset$, then either
$\frak p \le \frak q$ or $\frak q \le \frak p$ holds.
\item
We define a relation $\sim $ on the disjoint union
$\coprod_{\frak p \in \frak P}U_{\frak p}$ as follows.
Let $x \in U_{\frak p}, y \in U_{\frak q}$. We define $x\sim y$ if and
only if one of the following holds:
\begin{enumerate}
\item $\frak p = \frak q$ and $x=y$.
\item $\frak p \le \frak q$ and $y = \varphi_{\frak q\frak p}(x)$.
\item $\frak q\le \frak p$ and $x = \varphi_{\frak p\frak q}(y)$.
\end{enumerate}
Then $\sim$ is an equivalence relation.
\item
The quotient of $\coprod_{\frak p \in \frak P}U_{\frak p}/\sim$
by this equivalence relation is Hausdorff with respect to the
quotient topology.
\end{enumerate}
\par
In case $Z = X$ we call it a good coordinate system of $X$.
\par
In case ${\widetriangle{\mathcal U}}$ satisfies only
(1)-(6), we call it {\it a good coordinate system in the weak sense}.
\index{good coordinate system ! in weak sense}
\par
We call $Z$ the {\it support set}
\index{good coordinate system ! support set} of our good coordinate system.
\par
We also require that the dimension of $\mathcal U_{\frak p}$ is independent of
$\frak p$ and call it the {\it dimension} of $\widetriangle{\mathcal U}$.
\index{dimension ! of good coordinate system}
\par
We say a good coordinate system structure ${\widetriangle{\mathcal U}}
= (({\frak P},\le), \{\mathcal U_{\frak p}\mid \frak p \in \frak P\},
\{\Phi_{\frak p\frak q}
\mid \frak q \le \frak p\}) $
is {\it orientable} \index{orientable ! good coordinate system}
if we can choose orientation of $\mathcal U_{\frak p}$
such that all $\Phi_{\frak p\frak q}$ are orientation preserving.
The notion of orientation of orientable good coordinate system
and of oriented good coordinate system is defined in an obvious way.
\end{defn}
\begin{defn}\label{defofveruver}
We write $\vert{\widetriangle{\mathcal U}}\vert$ the quotient
set of the equivalence relation in Definition \ref{gcsystem} (7).
(See  \cite[Remark 5.20]{foootech} for its topology.)
\end{defn}
\begin{rem}
\begin{enumerate}
\item
Condition (7) above was not included in the definition of
good coordinate system in \cite{FO}, \cite{fooobook2}. This condition is due to Joyce \cite{joyce}.
\item
The fact that Condition (8) makes the argument used in our construction of
the perturbations more transparent became
clearer during the discussion at the google group
Kuranishi.
This condition was not included in the definition of
good coordinate system in \cite{FO}, \cite{fooobook2}.
The authors thank the members of the google group
Kuranishi who helped us much to polish the discussion here.
\item
However, we note that, for each `good coordinate system'
in the sense of \cite{FO}, \cite{fooobook2}, for which
(7), (8) are not necessarily satisfied, we can
{\it always} shrink $U_{\frak p}$ so that (7), (8)
are satisfied.
The detailed proof of this fact is given in \cite{foooshrink}.
Therefore, the statements in \cite{fooobook2} based on the definition of
`good coordinate system'
in the sense of \cite{fooobook2}, is correct  {\it
as stated there without change}\footnote{
Note however that there is an error related to the notion of `germ of Kuranishi
neighborhood' in \cite{FO}, which was explained in \cite[Subsection 34.1]{foootech}.
This error had been corrected in \cite[Section A1]{fooobook2}.}.
\item
Throughout this document, we denote by $\mathcal U_{\frak p}$ etc.
(where the index $\frak p$ is a German character) a Kuranishi chart
which is a member of good coordinate system,
and by $\mathcal U_{p}$ etc.
(where the index $p$ is an italic letter) a Kuranishi chart
which is a member of Kuranishi structure.
\end{enumerate}
\end{rem}
From now on, when we write coordinate change $\mathcal U_{\frak q} \to \mathcal U_{\frak p}$
as $\Phi_{\frak p\frak q}$ we use the notations
$U_{\frak p\frak q}$,  $\varphi_{\frak p \frak q}$ etc. where
$\Phi_{\frak p\frak q} = (U_{\frak p \frak q},\varphi_{\frak p \frak q},\widehat\varphi_{\frak p \frak q})$.
The same remark applies to $\Phi_{pq}$.
\begin{lem}\label{sumchart}
Let $\mathcal U_i
= (U_i,\mathcal E_i,\psi_i,s_i)$,
$(i=1,2)$ be Kuranishi charts
and $\Phi_{21} : \mathcal U_1 \to \mathcal U_2$
a coordinate change.
We assume that $\dim U_1 = \dim U_2$ and the map
$$
U_{21} \to U_1 \times U_2
$$
defined by $x \mapsto (\varphi_{21}(x),x)$ is proper.
Then
there exists a Kuranishi chart $\mathcal U_3
= (U_3,\mathcal E_3,\psi_3,s_3)$ and open KK-embeddings
$\Phi_{3i} = (\varphi_{3i},\widehat\varphi_{3i})
: \mathcal U_i \to \mathcal U_3$ $(i=1,2)$ such that
\begin{enumerate}
\item
$
\Phi_{32}\circ \Phi_{21}
= \Phi_{31}\vert_{U_{21}}.
$
\item
$U_3 = {\rm Im}(\varphi_{31}) \cup {\rm Im}(\varphi_{32})$.
\end{enumerate}
\end{lem}
\begin{proof}
We can glue two orbifolds (of the same dimension) $U_1$ and $U_2$
by the diffeomorphism $U_{21} \to U_2$ to its image (that is an open set).
Here $U_{21} \subset U_1$ is also an open set.
By the properness we assumed the glued space is Hausdorff.
Therefore, by gluing we obtain an orbifold, which we denote by $U_3$.
We can glue $\mathcal E_1$ and $\mathcal E_2$ to obtain $\mathcal E_3$.
The rest of the proof is obvious.
\end{proof}
We note that $\mathcal U_3$ is unique up to isomorphism that is
compatible with $\Phi_{3i}$. (The proof of this
uniqueness is easy and is left to the reader.)
We call the chart $\mathcal U_3$ the {\it sum chart}.\index{Kuranishi chart ! sum chart}
\begin{lem}\label{lem310}
Let  $\mathcal U_i
= (U_i,\mathcal E_i,\psi_i,s_i)$
$(i=1,2,3)$ and $\Phi_{ij}$ ($1\le i\le j\le 3$) be as in
Lemma \ref{sumchart}.
Let $\mathcal U_0= (U_0,\mathcal E_0,\psi_0,s_0)$
be another
Kuranishi chart and
$
\Phi_{0i} : 
\mathcal U_i \to \mathcal U_0
$
(resp.
$
\Phi_{i0} : 
\mathcal U_0 \to \mathcal U_i
$) embeddings of Kuranishi charts for $i=1,2$.
We assume
$$
\aligned
&\Phi_{02} \circ \Phi_{21}
= \Phi_{01}\vert_{U_{21}},\\
&\text{(resp.
$
U_{10} \cap U_{20} = \varphi_{10}^{-1}(U_{21})
$
and
$ \Phi_{21}  \circ \Phi_{10}\vert_{U_{10} \cap U_{20} }
=  \Phi_{20}\vert_{U_{10} \cap U_{20} }$.)}
\endaligned$$
Then
there exists a unique embedding of Kuranishi charts
$
\Phi_{03} 
: \mathcal U_3 \to \mathcal U_0
$
(resp.
$
\Phi_{30} 
: \mathcal U_0 \to \mathcal U_3
$)
such that
$$
\Phi_{03} \circ \Phi_{3i}
=
\Phi_{0i}
\qquad \text{(resp. $\Phi_{3i} \circ \Phi_{i0}
=
\Phi_{30}$).}
$$
\end{lem}
\begin{proof}
We note that our orbifold is always assumed to be effective
and we only consider embeddings as maps between them.
As its consequence, two such maps coincide if they coincide set-theoretically.
Moreover, if we are given  smooth maps (embedding) on  open
subsets of our orbifolds so that they coincide on the intersection
of the domain, then we can glue them to obtain a smooth map
(embedding).
Lemma \ref{lem310} is an immediate consequence of these facts.
\end{proof}
\begin{rem}
It seems that Lemma \ref{lem310} will become false
if we include noneffective orbifold.
\end{rem}
\subsection{Embedding of Kuranishi structures I}
\label{subsec:embedding1}
\begin{defn}\label{defn311}
Let $\widehat{\mathcal U} =(\{\mathcal U_p\},\{\Phi_{pq}\})$,
$\widehat{\mathcal U'} =(\{\mathcal U'_p\},\{\Phi'_{pq}\})$ be Kuranishi structures of $Z \subseteq X$.
A {\it strict KK-embedding}\index{embedding ! strict KK-embedding} $\widehat\Phi  = \{\Phi_{p}\}$
from $\widehat{\mathcal U}$ to $\widehat{\mathcal U'}$
assigns, to each $p \in Z$,  an embedding of Kuranishi charts
$\Phi_p = (\varphi_p,\widehat{\varphi}_p) : \mathcal U_p \to \mathcal U'_p$
such that for each $q \in {\rm Im}(\psi_p) \cap Z$ we have
the following:
\begin{enumerate}
\item[$\circledast$]
$
\Phi_p \circ \Phi_{pq}\vert_{U_{pq} \cap \varphi_q^{-1}(U'_{pq})}
=
\Phi'_{pq}\circ \Phi_q\vert_{U_{pq} \cap \varphi_q^{-1}(U'_{pq})}.
$
\end{enumerate}
(See Remark \ref{rem3737}.)
\par
We say that it is an {\it open KK-embedding}
\index{embedding ! open KK-embedding} if
$\dim U_p = \dim U'_p$ for each $p$.
\par
We say that
$\widehat{\mathcal U}$ is an {\it open substructure}
\index{Kuranishi structure ! open substructure} of $\widehat{\mathcal U'}$ if
there exists an open KK-embedding
$\widehat{\mathcal U} \to \widehat{\mathcal U'}$.
\par
A {\it KK-embedding} from $\widehat{\mathcal U}$ to $\widehat{\mathcal U^+}$\index{embedding ! KK-embedding}
is a strict KK-embedding $\widehat{\mathcal U_0} \to \widehat{\mathcal U^+}$
from an open substructure $\widehat{\mathcal U_0}$ of  $\widehat{\mathcal U}$.
\end{defn}
\begin{defn}
In the situation of Definition \ref{defn311} we assume that
$\widehat{\mathcal U}$ and $\widehat{\mathcal U'}$ are oriented.
We say that $\widehat\Phi  = \{\Phi_{p}\}$ is {\it orientation
preserving} \index{orientation
preserving ! embedding} if each of $\Phi_{p}$ is orientation preserving.
\end{defn}
\begin{rem}
The notion orientation preserving embedding can be defined
for other types of embeddings (there are 4 types of them see Table 5.1.)
in an obvious way.
\end{rem}
\begin{conven}\label{conv323}
Hereafter in this article we assume all the Kuranishi charts, Kuranish structures
and good coordinate systems are oriented unless otherwise mentioned explicitly.
We also assume all the coordinate change and embedding among Kuranishi charts,
Kuranish structures and good coordinate systems
are orientation preserving  unless otherwise mentioned explicitly.
\end{conven}

\begin{defn}\label{defn31222}
Let ${\widetriangle{\mathcal U}} =(\frak P,\{\mathcal U_{\frak p}\},\{\Phi_{\frak p\frak q}\})$,
${\widetriangle{\mathcal U'}} =(\frak P',\{\mathcal U'_{\frak p'}\},\{\Phi'_{\frak p'\frak q'}\})$ be good coordinate systems of $Z \subseteq X$.
A {\it GG-embedding}
\index{embedding ! GG-embedding} $\widetriangle{\Phi}  = \{\Phi_{\frak p}\}$
from ${\widetriangle{\mathcal U}}$ to ${\widetriangle{\mathcal U'}}$
assigns an order preserving map $\frak i : \frak P \to \frak P'$ and, to each ${\frak p}\in \frak P$,  an embedding of Kuranishi charts
$\Phi_{\frak p} = (\varphi_{\frak p},
\widehat{\varphi}_{\frak p}) : \mathcal U_{\frak p} \to \mathcal U'_{\frak i({\frak p})}$
such that for each $\frak q \le {\frak p}$ we have
the following:
\begin{enumerate}
\item
$
U_{\frak p\frak q} = \varphi_{
\frak q}^{-1}(U'_{\frak i(\frak p)\frak i(\frak q)}).
$
\item
$
\Phi_{\frak p} \circ \Phi_{\frak p\frak q}
=
\Phi'_{\frak i(\frak p)\frak i(\frak q)}\circ \Phi_{\frak q}\vert_{U_{\frak p\frak q}}.
$
\end{enumerate}
\begin{equation}\label{diag33--}
\begin{CD}
\mathcal U_{\frak q}\vert_{U_{\frak p\frak q}
} @ > {\Phi_{\frak q}} >>
{\mathcal U}'_{\frak i(\frak q)}\vert_{{U}'_{\frak i(\frak p)\frak i(\frak q)}}  \\
@ V{\Phi_{\frak p\frak q}}VV @ VV{\Phi'_{\frak i(\frak p)\frak i(\frak q)}}V\\
\mathcal U_{\frak p} @ > {\Phi_{\frak p}} >>{\mathcal U}'_{\frak i(\frak p)}
\end{CD}
\end{equation}We say that $\widetriangle\Phi$ is a {\it weakly open GG-embedding}
\index{embedding ! weakly open GG-embedding} if
$\dim U_{\frak p} = \dim U'_{\frak i(\frak p)}$ for each $\frak p$.
We say it is an {\it open GG-embedding}
\index{embedding ! open GG-embedding} if
$\frak P = \frak P'$ in addition. We say it is a {\it strongly open GG-embedding}
\index{embedding ! strongly open GG-embedding} if
$$
\psi_{\frak p}(U_{\frak p} \cap s_{\frak p}^{-1}(0)) = \psi_{\frak p}(U'_{\frak p} \cap (s'_{\frak p})^{-1}(0))
$$
holds in addition.
\par
We say that
${\widetriangle{\mathcal U}}$ is an {\it open substructure}
\index{good coordinate system ! open substructure}
(resp. {\it weakly open substructure}, {\it strongly open substructure})
\index{good coordinate system ! weakly open substructure}
\index{good coordinate system ! strongly open substructure}
of ${\widetriangle{\mathcal U'}}$ if
there exists an open  (resp. weakly open, strongly open) GG-embedding
${\widetriangle{\mathcal U}} \to {\widetriangle{\mathcal U'}}$.
\par
We say a GG-embedding $\widetriangle{\Phi}$ is an {\it isomorphism}
\index{good coordinate system ! isomorphism}
\index{embedding ! isomorphism of GG-embeddings}
if the map $\frak i$ is a bijection and $\widehat{\varphi}_{\frak p}$ is an
isomorphism for each $\frak p$.
\end{defn}
\begin{rem}
Definition \ref{defn31222}, especially Item (1), implies that a GG-embedding
${\widetriangle{\mathcal U}} \to  {\widetriangle{\mathcal U'}}$
induces an injective continuous map
$\vert {\widetriangle{\mathcal U}}\vert
\to \vert {\widetriangle{\mathcal U'}}\vert$.
\end{rem}
\begin{lem}\label{lem320}
Let $\widetriangle{\mathcal U} = (\frak P,\{\mathcal U_{\frak p}\},\{\Phi_{\frak p\frak q}\})$
be a good coordinate system of $Z \subseteq X$ and let
$U^0_{\frak p} \subseteq U_{\frak p}$ be given open subsets such that
$
Z \subset \bigcup_{\frak p \in \frak P} \psi_{\frak p}(s_{\frak p}^{-1}(0) \cap U^0_{\frak p}).
$
Then there exists a unique coordinate change $\Phi_{\frak p\frak q}^0$ such that
$(\frak P,\{\mathcal U_{\frak p}\vert_{U^0_{\frak p}}\},\{\Phi^0_{\frak p\frak q}\})$
is an open substructure of $\widetriangle{\mathcal U}$.
\end{lem}
\begin{proof}
Let $\Phi_{\frak p\frak q} = (U_{\frak p\frak q},\varphi_{\frak p\frak q},
\widehat\varphi_{\frak p\frak q})$.
We put
\begin{equation}\label{eqform3636}
U^0_{\frak p\frak q}
=
U_{\frak p\frak q} \cap U^0_{\frak q} \cap \varphi_{\frak p\frak q}^{-1}(U^0_{\frak p})
\end{equation}
and
$\Phi^0_{\frak p\frak q} = (U^0_{\frak p\frak q},\varphi_{\frak p\frak q}\vert_{U^0_{\frak p\frak q}},
\widehat\varphi_{\frak p\frak q}\vert_{U^0_{\frak p\frak q}})$.
It is easy to see that $(\frak P,\{\mathcal U_{\frak p}\vert_{U^0_{\frak p}}\},\{\Phi^0_{\frak p\frak q}\})$
is an open substructure of $\widetriangle{\mathcal U}$.
\par
On the other hand, if $(\frak P,\{\mathcal U_{\frak p}\vert_{U^0_{\frak p}}\},\{\Phi^0_{\frak p\frak q}\})$
is an open substructure of $\widetriangle{\mathcal U}$,
then Definition \ref{defn31222} (1) implies that the domain
$U^0_{\frak p\frak q}$ of $\Phi_{\frak p\frak q}$ must be as in (\ref{eqform3636}).
\end{proof}
\begin{lem}\label{lem321321}
Let $\widehat{\mathcal U} = (\{\mathcal U_{p}\},\{\Phi_{pq}\})$
be a Kuranishi structure of $Z \subseteq X$ and
$U^0_{p} \subseteq U_{p}$  open subsets containing $p$.
Then there exists $\Phi_{pq}^0$ such that
$(\{\mathcal U_{p}\vert_{U^0_{p}}\},\{\Phi^0_{pq}\})$
is an open substructure of $\widehat{\mathcal U}$.
\end{lem}
\begin{rem}
The uniqueness dose not hold in Lemma \ref{lem321321}
since there is no condition similar to
Definition \ref{defn31222} (1) for Kuranishi structure.
They have, however, a common open substructure.
\end{rem}
\begin{proof}
We put
\begin{equation}\label{eqform3636revrev}
U^0_{pq}
=
U_{pq} \cap U^0_{q} \cap \varphi_{pq}^{-1}(U^0_{p})
\end{equation}
and
$\Phi^0_{pq} = (U^0_{pq},\varphi_{pq}\vert_{U^0_{pq}},
\widehat\varphi_{pq}\vert_{U^0_{pq}})$.
It is easy to see that $(\{\mathcal U_{p}\vert_{U^0_{p}}\},\{\Phi^0_{pq}\})$
is a Kuranishi structure.
\end{proof}
\begin{defn}\label{defn32020202}
Let $\widehat{\mathcal U}$ be a Kuranishi structure and
${\widetriangle{\mathcal U}}$  a good coordinate system
of $Z \subseteq X$.
A {\it strict KG-embedding}
\index{embedding ! strict KG-embedding} of $\widehat{\mathcal U}$ to ${\widetriangle{\mathcal U}}$
assigns, for each $p\in Z$, $\frak p\in \frak P$ with
$p \in {\rm Im}(\psi_{\frak p})$,
an embedding of Kuranishi charts
$
\Phi_{\frak p p} = (\varphi_{\frak p p},\widehat\varphi_{\frak p p}) :
\mathcal U_p \to \mathcal U_{\frak p}$ with the following
properties.
\par
If $\frak p,\frak q\in \frak P$,  $\frak q \le \frak p$,
$p \in {\rm Im}(\psi_{\frak p}) \cap Z$,
$q \in {\rm Im}(\psi_{p}) \cap \psi_{\frak q}(U_{\frak p\frak q}) \cap Z$,
then
the following diagram commutes.
\begin{equation}\label{diag33}
\begin{CD}
\mathcal U_{q}\vert_{U_{pq}
\cap \varphi_{\frak q q}^{-1}(U_{\frak p \frak q})} @ > {\Phi_{\frak q q}} >>
{\mathcal U}_{\frak q}\vert_{{U}_{\frak p\frak q}}  \\
@ V{\Phi_{pq}}VV @ VV{\Phi_{\frak p \frak q}}V\\
\mathcal U_{p} @ > {\Phi_{\frak p p}} >>{\mathcal U}_{\frak p}
\end{CD}
\end{equation}
\par
A {\it KG-embedding} of $\widehat{\mathcal U}$ to ${\widetriangle{\mathcal U}}$
is by definition a strict KG-embedding of an open substructure of $\widehat{\mathcal U}$
to ${\widetriangle{\mathcal U}}$.
\index{embedding ! KG-embedding}
\par
If $\widehat{\mathcal U_0} \to \mathcal{\widehat{\mathcal U}}$
is a
KK-embedding,
${\widetriangle{\mathcal U}}
\to \widetriangle{{\mathcal U}^+}$ is a GG-embedding,
and $\widehat{\mathcal U} \to {\widetriangle{\mathcal U}}$ is a strict
KG-embedding
(resp. KG-embedding), then
the {\it composition}
$$
\widehat{\mathcal U}_0 \to \widehat{\mathcal U}
    \to  {\widetriangle{\mathcal U}} \to \widetriangle{{\mathcal U}^+}
$$
is defined as a strict KG-embedding  (resp.  a KG-embedding).
\index{embedding ! composition }
(See Definition \ref{definition516161} (3).)
\end{defn}

The next result is  the same as
\cite[Lemma 6.3]{FO}, \cite[Theorem 7.1]{foootech}.
We will reproduce its proof together with various
addenda in Section \ref{sec:contgoodcoordinate}.
In particular, Theorem \ref{Them71restate} is proved in Subsection \ref{subsec:constgcsabs}.

\begin{thm}\label{Them71restate}
For any Kuranishi structure $\widehat{\mathcal U}$ of
$Z \subseteq X$ there exist a
good coordinate system ${\widetriangle{\mathcal U}}$
of $Z \subseteq X$
and a KG-embedding $\widehat{\mathcal U} \to {\widetriangle{\mathcal U}}$.
\end{thm}
\begin{rem}
According to Convention \ref{conv323},
Theorem \ref{Them71restate} contains the statement that
${\widetriangle{\mathcal U}}$ and the KG-embedding $\widehat{\mathcal U}\to {\widetriangle{\mathcal U}}$
are oriented (provided $\widehat{\mathcal U}$ is oriented).
We do not repeat this kinds of remarks later on.
\end{rem}

\begin{defn}\label{gcsystemcompa}
A good coordinate system $\widetriangle{\mathcal U}$ is said to be
{\it compatible} with  a Kuranishi structure
${\widehat{\mathcal U}}$ if there exists a KG-embedding
$\widehat{\mathcal U} \to {\widetriangle{\mathcal U}}$.
\index{compatibility ! between good coordinate system and Kuranishi structure}
\index{good coordinate system ! compatibility between good coordinate system and Kuranishi structure}
\index{Kuranishi structure ! compatibility between good coordinate system and Kuranishi structure}
\end{defn}
\begin{rem}
\begin{enumerate}
\item
The next terminology is due to Joyce \cite{joyce}.
\begin{defn}
A good coordinate system is said to be
{\it excellent} if $\frak P \subset \Z_{\ge 0}$,
$\le$ is the standard inequality on  $\Z_{\ge 0}$
and $\dim U_{\frak p} = \frak p$.
\end{defn}
Starting with
an arbitrary good coordinate system
$\widetriangle{\mathcal U} = (\{U_{\frak p}\},\{\Phi_{\frak p\frak q}\})$,
we can construct an excellent good coordinate system
$\widetriangle{\mathcal U'}$ as follows.
Suppose $\dim U_{\frak p_1} = \dim U_{\frak p_2}$.
If ${\rm Im}(\psi_{\frak p_1})$ is disjoint from ${\rm Im}(\psi_{\frak p_2})$,
we take its disjoint union as a new chart and remove these two charts.
Suppose ${\rm Im}(\psi_{\frak p_1}) \cap {\rm Im}(\psi_{\frak p_2}) \ne \emptyset$.
Then we may assume $\frak p_1 < \frak p_2$.
Since an embedding between two orbifolds of the same
dimension is necessarily a diffeomorphism,
the coordinate change $\Phi_{\frak p_2\frak p_1}$ is an isomorphism.
We can use this observation and Lemma \ref{sumchart} to construct the sum chart
of $\mathcal U_{\frak p_1}$ and $\mathcal U_{\frak p_2}$.
We take it as a new chart and remove $\mathcal U_{\frak p_1}$ and $\mathcal U_{\frak p_2}$.
The coordinate change between sum charts and other charts can be defined by using
Lemma \ref{lem310}.
\par
We can continue this process finitely many times until we get an excellent coordinate system.
Note there is a weakly open GG-embedding
$\widetriangle{\mathcal U} \to \widetriangle{\mathcal U'}$.
\item
We note that in \cite[Section 7]{foootech} we introduced
the notion of mixed neighborhood.
It is basically equivalent to  the notion of
a pair of an excellent
good coordinate system ${\widetriangle{\mathcal U}}$
and a KG-embedding $\widehat{\mathcal U} \to {\widetriangle{\mathcal U}}$.
(The only difference is that the conclusion of \cite[Lemma 7.32]{foootech}
is not assumed in \cite[Definition 7.15]{foootech}.
This difference is not essential at all because of \cite[Lemma 7.32]{foootech}.)
Therefore Theorem \ref{Them71restate} is
equivalent to \cite[Theorem 7.1]{foootech}.
\end{enumerate}
\end{rem}
\begin{defn}\label{mapkura}
Let $\widehat{\mathcal U}$ be a Kuranishi structure  of $Z \subseteq X$.
\begin{enumerate}
\item
A {\it strongly continuous map $\widehat f$}
\index{map from Kuranishi structure ! strongly continuous map}
from $(X,Z;\widehat{\mathcal U})$ to
a topological space $Y$ assigns a continuous map $f_{p}$ from
$U_{p}$ to $Y$ for each $p\in X$ such that
$f_p \circ \varphi_{pq} = f_q$ holds on $U_{pq}$.
\item
In the situation of (1), the map $f:Z \to Y$ defined by
$f(p) = f_p(p)$ is a continuous map from $Z$ to $Y$.
We call $f : Z \to Y$ the {\it underlying continuous map} of $\widehat f$.
\item
We require that the underlying continuous map $f : Z \to Y$ is extended
to a continuous map $f : X \to Y$ and include it to the data
defining a strong continuous map.
\item
When $Y$ is a smooth manifold, we say $\widehat f$ is {\it strongly smooth}
\index{map from Kuranishi structure ! strongly smooth}
if each of $f_p$ is smooth.
\item
A strongly smooth map is said to be {\it weakly submersive}
\index{map from Kuranishi structure ! weakly submersive} if each of
$f_p$ is a submersion.
\end{enumerate}
We sometimes say $f$ is a strongly continuous map
(resp. a strongly smooth map, a weakly submersive map) in place of
$\widehat f$ is a strongly continuous map
(resp. a strongly smooth map, weakly submersive),
by an abuse of notation.
\end{defn}
\begin{rem}\label{remrem327}
The continuity claimed in (2) follows from the next diagram
\begin{equation}
\begin{CD}
s_p^{-1}(0) @ > {\psi_{p}} >>
Z  \\
@ V{f_p}VV @ VV{f}V\\
Y @ = Y
\end{CD}
\end{equation}
whose comutativity is a consequence of
$f_p \circ \varphi_{pq} = f_q$.
\end{rem}
\begin{rem}
In \cite{joyce}, Joyce used the terminology `strong submersion'
instead of `weak submersion'
which we have been used since \cite{FO}.
We insist on using the terminology `weakly submersive'
by the following two reasons.
\begin{enumerate}
\item
Let $V_1 \subset V_2$ be a submanifold, $f_2$  a map from $V_2$
and $f_1$  a restriction of $f_2$ to $V_1$.
The condition that $f_2$ is smooth is stronger than the condition
that $f_1$ is smooth. So we used the word strongly smooth.
On the other hand, in case $f_2$ is smooth, the condition that
$f_2$ is a submersion at each point of $V_1$
is weaker than the condition that $f_1$ is
a submersion. So we used the word weak submersion.
\item
Later we will use the terminology `strongly submersive' for a different notion.
See Definitions \ref{transofdvect} and \ref{transkurakuravect}.
\end{enumerate}
\end{rem}

\begin{defn}\label{definition32727}
Let ${\widetriangle{\mathcal U}}$ be a good coordinate system of
$Z \subseteq X$.
\begin{enumerate}
\item
A {\it strongly continuous map $\widetriangle f$} from $(X,Z;{\widetriangle{\mathcal U}})$
to
a topological space $Y$ assigns a continuous map $f_{\frak p}$ from
$U_{\frak p}$ to $Y$ to each $\frak p\in \frak P$ such that
$f_{\frak p} \circ \varphi_{{\frak p}{\frak q}} = f_{\frak q}$ holds on
$U_{{\frak p}{\frak q}}$.
\item
In the situation of (1),
the map $f: Z \to Y$ defined by
$f(p) = f_{\frak p}(o_{\frak p}(p))$ (for $p \in {\rm Im}(\psi_{\frak p})\cap Z$) is a continuous map from $Z$ to $Y$.
\footnote{The proof of continuity is the same as Remark \ref{remrem327}.}
We call $f : Z \to Y$ the {\it underlying continuous map} of $\widetriangle f$.
\item
We require that the underlying continuous map $f : Z \to Y$ is extended
to a continuous map $f : X \to Y$ and include it to the data
defining strong continuous map.
\item
When $Y$ is a smooth manifold, we say $\widetriangle f$ is {\it strongly smooth}
if each of $f_{\frak p}$ is smooth.
\item
A strongly smooth map is said to be {\it weakly submersive} if each of
$f_{\frak p}$ is a submersion.
\end{enumerate}
We sometimes say $f$ is a strongly continuous map
(resp. a strongly smooth map, weakly submersive) in place of
$\widetriangle f$ is a strongly continuous map
(resp. a strongly smooth map, weakly submersive),
by an abuse of notation.
\end{defn}
\begin{defn}\label{defn320}
\begin{enumerate}
\item
If $\widehat f : (X,Z;\widehat{\mathcal U'}) \to Y$
is a strongly continuous map and
$\widehat\Phi = \{\Phi_{p}\} : \widehat{\mathcal U} \to \widehat{\mathcal U'}$
is a KK-embedding,
then $f_p \circ \varphi_p : U_p \to Y$  defines a
strongly continuous map, which we call
the {\it pullback}
\index{map from Kuranishi structure ! pullback of strongly continuous map} and write $\widehat f \circ \widehat\Phi$.
\item
Let $\widehat{\Phi}$ be a KK-embedding.
If $\widehat f$ is strongly smooth, then so is $\widehat f \circ \widehat\Phi$.
If $\widehat{\Phi}$ is an open embedding and
$\widehat f$ is weakly submersive,
then $\widehat f \circ \widehat\Phi$ is also weakly submersive.
\item
The good coordinate system version of pullback of maps
can be defined in the same way as (1). A similar statement as (2) holds as well.
\item
If $\widetriangle f : (X,Z;{\widetriangle{\mathcal U}}) \to Y$
is a strongly continuous map from a good coordinate system and
$\widehat\Phi : \widehat{\mathcal U} \to {\widetriangle{\mathcal U}}$
is a KG-embedding, then the pullback  $\widetriangle f\circ \widehat\Phi$
can be defined in the same way as (1).
A similar statement as (2) holds as well.
\end{enumerate}
\end{defn}


\section{Fiber product of Kuranishi structures}
\label{sec:fiber}

\subsection{Fiber product}
\label{subsec:fibprod}

Before studying fiber product we consider direct product.
Let $X_i$, $i=1,2$ be separable metrizable spaces, $Z_i \subseteq X_i$ compact  subsets, and $\widehat{\mathcal U}_i$ Kuranishi structures
of $Z_i \subseteq X_i$.
We will define a Kuranishi structure  of
the direct product $Z_1 \times Z_2 \subseteq X_1 \times X_2$.
\begin{defn}\label{directproduct}
For $p_i \in Z_i$ let
$\mathcal U_{p_i} = (U_{p_i},\mathcal E_{p_i},\psi_{p_i},s_{p_i})$ be their Kuranishi neighborhoods.
Then the Kuranishi neighborhood of $p = (p_1,p_2) \in Z_1 \times Z_2$ is
$\mathcal U_{p} = \mathcal U_{p_1} \times \mathcal U_{p_2} = (U_p,\mathcal E_p,\psi_p,s_p)$ where
$$
 (U_p,\mathcal E_p,\psi_p,s_p) =
(U_{p_1}\times U_{p_2},\mathcal E_{p_1}\times \mathcal E_{p_2},\psi_{p_1}\times \psi_{p_2},s_{p_1}\times
s_{p_2}) .
$$
This system satisfies the condition of
Kuranishi neighborhood (Definition \ref{kuranishineighborhooddef}.)
\par
Suppose $q_i \in Z_i$ and $q = (q_1,q_2) \in Z_1 \times Z_2$.
If
$q \in \psi_p(s^{-1}_p(0))$, then it is easy to see that
$q_i \in \psi_{p_i}(s^{-1}_{p_i}(0))$ for $i=1,2$.
Therefore there exist coordinate changes
$\Phi_{p_iq_i} = ({\varphi_{p_iq_i}},\widehat{\varphi}_{p_iq_i},h_{p_iq_i})$ from
$\mathcal U_{q_i}$ to
$\mathcal U_{p_i}$.
We define
$$
\aligned
\Phi_{pq} &=
\Phi_{p_1q_1} \times  \Phi_{p_2q_2}=
(U_{pq},{\varphi}_{pq},\widehat{\varphi}_{pq})
 \\
&=
(U_{p_1q_1} \times U_{p_2q_2},{\varphi}_{p_1q_1}\times {\varphi}_{p_2q_2},\widehat{\varphi}_{p_1q_1}\times \widehat{\varphi}_{p_2q_2}).
\endaligned$$
This satisfies the condition of coordinate change
of
Kuranishi charts (Definition \ref{coordinatechangedef}).
\par
Then it is also easy to show that $(\{\mathcal U_{p_1} \times \mathcal U_{p_2}\},\{\Phi_{p_1q_1} \times  \Phi_{p_2q_2}\})$
defines a Kuranishi structure of $Z_1 \times Z_2 \subseteq X_1 \times X_2$ in the sense of Definition
\ref{kstructuredefn}.
(We note that effectivity of an orbifold is preserved by the
direct product.)
We call this Kuranishi structure the {\it direct product Kuranishi structure}.
\index{Kuranishi structure ! direct product}
\index{fiber product ! direct product of Kuranishi structures}
\end{defn}
We can easily prove  that
the direct product of oriented Kuranishi structures  (\cite[Definition 4.5]{foootech}) is also
oriented.
\par\medskip
Next we study fiber product.
Let $(X,Z;\widehat{\mathcal U})$ be a relative K-space
and $\widehat f = \{f_p\}: (X,Z;\widehat{\mathcal U}) \to N$  a strongly smooth map,
where $N$ is a smooth manifold of finite dimension.
Let $f' : M \to N$ be a smooth map between smooth manifolds.
We assume $M$ is compact.
In this section  we define a Kuranishi structure
on the pair of topological spaces
\begin{equation}\label{formula4141}
\aligned
Z \times_N M
&= \{(p,q) \in Z \times M \mid f(p) = f'(q)\},
\\
X \times_N M
&= \{(p,q) \in X \times M \mid f(p) = f'(q)\},
\endaligned
\end{equation}
that is the fiber product.
\index{fiber product ! of topological spaces}
The assumption we need to require is certain transversality,
which we define below.
\begin{defn}\label{transverse1}
We say $\widehat f$ is {\it weakly transversal} to $f'$
\index{weakly transversal to a smooth map}
on $Z \subseteq X$ if the following holds.
Let  $(p,q) \in Z \times_N M$ and
$\mathcal U_p = (U_p,E_p,s_p,\psi_p)$ be a Kuranishi neighborhood
of $p$.
We then require that
for each $(x,y) \in U_p\times M$ with
$f_p(x) = f'(y)$ we have
\begin{equation}\label{transformula}
(d_xf_p)(T_xU_p) + (d_yf')(T_yM) = T_{f(x)}N.
\end{equation}
Let us explain the meaning of (\ref{transformula}).
We take an orbifold chart $(V_p(x),\Gamma_p(x),\phi_p(x))$ of
$U_p$ at $x$.
(Definition \ref{defn26550}.)
The composition
$$
V_p(x) \overset{\phi_p(x)}\longrightarrow  U_p \overset{f_p}\longrightarrow N
$$
is by assumption a smooth map which we write $f_{p,x}$.
(\ref{transformula}) means
$$
(d_{o(x)}f_{p.x})(T_{o_p(x)}V_p(x)) + (d_yf')(T_yM) = T_{f(x)}N
$$
where $o(x) \in V_p(x)$ is the base point which satisfies
$\phi_p(x)(o(x)) = x$.
\end{defn}
\begin{exm}\label{fiberexa}
\begin{enumerate}
\item
If $\widehat f : (X,Z;\widehat{\mathcal U}) \to N$ is weakly submersive in the sense of Definition \ref{mapkura},
then for any $f' : M \to N$, $f$ is weakly transversal to $f'$.
\item
If $f' : M \to N$ is a submersion, then any strongly smooth map
$\widehat f :  (X,Z;\widehat{\mathcal U}) \to N$ is weakly transversal to $f'$.
\item
The pullback of map in Definition \ref{defn320}
by an open embedding preserves the weak transversality.
\item
Suppose
$(X_i,Z_i)$ $(i=1,2)$ have Kuranishi structures
$\widehat{\mathcal U_i}$ and
the maps
$\widehat f_i : (X_i,Z_i;\widehat{\mathcal U_i}) \to N$
are strongly smooth.
We put:
$$
Z_1 \times_N Z_2
=
\{(p,q) \in Z_1 \times Z_2 \mid f_1(p) = f_2(q) \}.
$$
Let $(p,q)  \in Z_1 \times_N Z_2$.
We denote by $\mathcal U_p$,
$\mathcal U_q$ the Kuranishi neighborhoods of $p, \, q$ respectively
and assume
the condition
\begin{equation}\label{2transverse}
(d_x(f_1)_p)(T_xU_p) + (d_y(f_2)_q)(T_yU_q) = T_{(f_1)_p(x)}N
\end{equation}
for each $(x,y) \in U_p\times U_q$ with $(f_1)_p(x) = (f_2)_q(y)$.
(The precise meaning of (\ref{2transverse}) can be
defined in a similar way as the case of (\ref{transformula}).)
It is easy to see that  (\ref{2transverse}) is equivalent to the next
condition.
We consider the map
$$
f = (f_1,f_2) : X_1 \times X_2 \to N\times N.
$$
We use the direct product Kuranishi structure (Definition \ref{directproduct})
on $Z_1\times Z_2 \subseteq X_1\times X_2$.
Then
(\ref{2transverse}) holds if and only if
$f$ is transversal to the diagonal embedding
$N \to N \times N$ in the sense of Definition \ref{transverse1}.
\item
We can generalize the situation of (4) to the case when three or more
factors are involved. In fact, in the study of the moduli space of
pseudo-holomorphic curves, we will encounter the situation
where we consider the fiber product of various factors
which are organized by a tree or a graph.
See Parts 2 and 3 or \cite[Subsection 7.1.1]{fooobook2}.
\end{enumerate}
\end{exm}
\begin{defn}
In the situation of Example \ref{fiberexa} (4) we say
$\widehat f_1$ is {\it weakly transversal} to $\widehat f_2$ if (\ref{2transverse})
is satisfied.
\index{transversality ! between two strongly smooth maps}
\end{defn}
Now we assume that $\widehat f : (X,Z;\widehat{\mathcal U}) \to N$ is weakly transversal to $f' : M \to N$
in the sense of Definition \ref{transverse1} and define a
Kuranishi structure on the fiber product (\ref{formula4141}).
Recalling that a Kuranishi neighborhood $U_p$ of $p \in Z$ is assumed to be
an effective orbifold, we find the following.
\begin{lem}\label{fibereffec}
\index{fiber product ! effectivity of fiber product}
For each $(p,\frak x) \in Z \times_N M$ the fiber product $U_p\times_N M$  is
again an effective orbifold.
\end{lem}
\begin{proof}
Let $(p, {\mathfrak x}) \in Z  \times_N M$.  Pick an orbifold chart $(V_p,
\Gamma_p)$ at $p$.
Denote by ${\tilde o}_p \in V_p$ the point, which is mapped to $o_p$ under
the quotient map
$V_p \to U_p$ by $\Gamma_p$ and by $\widetilde{f}_p$ the lift of $f_p:U_p
\to N$.
Since $\widetilde{f}_p$ is $\Gamma_p$-invariant, we find that
$K_{{\tilde o}_p}={\rm Ker~} d\widetilde{f}_p $ at ${\tilde o}_p$ is
transversal to the tangent space
$T_{{\tilde o}_p} V_p^{\Gamma_p}$ of the fixed point set $V_p^{\Gamma_p}$ by
$\Gamma_p$-action.
Hence the tangent space $T_{[{\tilde o}_p, x]} V_p \times_N M$ contains
$K_{{\tilde o}_p}$.
Since $\Gamma_p$ acts trivially on $V_p^{\Gamma_p}$ and $\Gamma_p$ acts
effectively on $V_p$,
$\Gamma_p$ acts effectively on $K_{{\tilde o}_p}$.
It implies that the fiber product $U_p \times_N M$ is an effective orbifold.
\end{proof}
\par
Let $\mathcal U_p$ be the given Kuranishi neighborhood of
$p$ and $(p,\frak x) \in Z \times_N M$.
We define
$
U_{(p,\frak x)} = U_p \times_N M.
$
Note $U_{(p,\frak x)}$ is a smooth orbifold by  Definition \ref{transverse1}.
The bundle $\mathcal E_{(p,\frak x)}$ is the pullback of $\mathcal E_p$ by the
map $U_{(p,\frak x)} \to U_p$ that is the projection to the first factor.
The section $s_p$ induces $s_{(p,\frak x)}$ of $\mathcal E_{(p,\frak x)}$
in an obvious way.
Note
$
s_{(p,\frak x)}^{-1}(0) = s_p^{-1}(0) \times_N M.
$
Therefore $\psi_p :  s_p^{-1}(0) \to X$
induces
$$
\psi_{(p,\frak x)} : s_{(p,\frak x)}^{-1}(0)
= s_p^{-1}(0) \times_N M
\to X \times_N M.
$$
It is easy to see that $\psi_{(p,\frak x)}$ is a homeomorphism
onto a neighborhood of $(p,\frak x)$.
\par
In sum we have the following:
\begin{lem}\label{lem25}
$\mathcal U_{(p,\frak x)} = (U_{(p,\frak x)},\mathcal E_{(p,\frak x)},s_{(p,\frak x)},\psi_{(p,\frak x)})$
is a Kuranishi neighborhood of $(p,\frak x) \in X\times_N M$.
\end{lem}
We next consider the coordinate change.
Let  $(p,\frak x), (q,\frak y) \in Z\times_N M$. We assume
$
(q,\frak y) = \psi_{(p,\frak x)}(x,y)
$
where $(x,y) \in  V_{(p,\frak x)}$.
By definition we have
$
q = \psi_p(x).
$
Therefore by  Definition \ref{kstructuredefn} there exists a coordinate change
$\Phi_{pq} = (U_{pq},\varphi_{pq},\widehat{\varphi}_{pq})$
from $\mathcal U_q$
to
$\mathcal U_p$
in the sense of Definition \ref{coordinatechangedef}.
\par
Now we put
\begin{enumerate}
\item
$
U_{(p,\frak x),(q,\frak y)} = U_{pq} \times_N M
$.
\item
$\varphi_{(p,\frak x),(q,\frak y)} = \varphi_{pq} \times_N {\rm id}
: U_{pq} \times_N M \to U_p \times_N M.$
\item
$\widehat\varphi_{(p,\frak x),(q,\frak y)} = \widehat\varphi_{pq} \times_N {\rm id} :
\mathcal E_q\vert_{U_{pq}} \times_N M \to \mathcal E_p\times_N M $.
\end{enumerate}
\begin{lem}\label{lem27}
$\Phi_{(p,\frak x),(q,\frak y)} = (U_{(p,\frak x),(q,\frak y)}
,\varphi_{(p,\frak x),(q,\frak y)},\widehat\varphi_{(p,\frak x),(q,\frak y)})$ defines a coordinate change
from $\mathcal U_{(q,\frak y)}$
to $\mathcal U_{(p,\frak x)}$.
\end{lem}
The proof is immediate from the definition.
\begin{lem}\label{fiberkura22}
Suppose $Z \subseteq X$ has a Kuranishi structure
$\widehat{\mathcal U}$ and $\widehat f : (X,Z;\widehat{\mathcal U}) \to N$ is weakly transversal to $f' : M \to N$.
Then the Kuranishi neighborhoods in Lemma \ref{lem25}
together with coordinate changes in Lemma \ref{lem27}
define a Kuranishi structure of the fiber product
of $Z \times_N M \subseteq X \times_N M$.
\end{lem}
The proof is again immediate from the definition.
\begin{defn}\label{firberproddukuda}
\begin{enumerate}
\item
Suppose $\widehat f : (X,Z;\widehat{\mathcal U}) \to N$ is weakly transversal to $f' : M \to N$.
We call the Kuranishi structure obtained in Lemma \ref{fiberkura22},
the {\it fiber product Kuranishi structure}
and write the resulting relative K-space by
$$
(X,Z;\widehat{\mathcal U}) \,{}_{f}\times_{f'} M
\quad
\text{or}
\quad
 (X,Z;\widehat{\mathcal U})  \times_{N} M.
$$
\item
Suppose $\widehat f_i : (X_i,Z_i;\widehat{\mathcal U}_i) \to N$ are strongly
smooth maps. We assume $\widehat f_1$ and $\widehat f_2$ are weakly transversal
in the sense of Example \ref{fiberexa} (4).
Then we define the {\it fiber product}
\index{Kuranishi structure ! fiber product of Kuranishi structures}
\index{fiber product ! of Kuranishi structures}
$$
(X_1,Z_1;\widehat{\mathcal U}_1) \,{}_{f_1}\times_{f_2}
(X_2,Z_2;\widehat{\mathcal U}_2)
\quad
\text{or}
\quad
(X_1,Z_1;\widehat{\mathcal U}_1)  \times_{M}
(X_2,Z_2;\widehat{\mathcal U}_2)
$$
as the fiber product
$$
\left(
(X_1,Z_1;\widehat{\mathcal U}_1) \times
(X_2,Z_2;\widehat{\mathcal U}_2)\right)
\,\,{}_{f_1 \times f_2}\times_{i} \Delta_M.
$$
Here $i : \Delta_M \to M\times M$ is the embedding of the diagonal.
\end{enumerate}
\end{defn}
For the purpose of reference we also include another obvious statements.
\begin{lem}
We consider the situation of Lemma \ref{fiberkura22}.
\begin{enumerate}
\item
If $\widehat g : (X,Z;\widehat{\mathcal U}) \to M'$ is also a strongly continuous map, then
it induces a strongly continuous map
$Z \times_N M \to M'$.
It is weakly submersive if $(\widehat f,\widehat g) : (X,Z;\widehat{\mathcal U}) \to N \times M'$ is weakly
submersive.
\item
If $\widehat f$ is weakly submersive, then
the projection $Z \times_N M \to M$ is weakly submersive .
\end{enumerate}
\end{lem}
\begin{lem}
Let $\widehat{\mathcal U}$, $\widehat{\mathcal U^+}$ be Kuranishi structures of
$Z \subseteq X$ and  $\widehat{\Phi} : \widehat{\mathcal U} \to \widehat{\mathcal U^+}$ a KK-embedding.
Let $\widehat f : (X,Z;\widehat{\mathcal U^+}) \to N$ be a
strongly smooth map and $g : M \to N$ a smooth map between manifolds.
\begin{enumerate}
\item
If $\widehat f\circ \widehat{\Phi} : (X,Z;\widehat{\mathcal U}) \to N$
is weakly transversal to $g$, then
$\widehat f : (X,Z;\widehat{\mathcal U^+}) \to N$ is
weakly transversal to $g$.
\item
In the situation of (1),
$\widehat{\Phi}$ induces a KK-embedding
$$
\widehat{\Phi} \times_N M :
\widehat{\mathcal U} \times_N M
\to
\widehat{\mathcal U^+} \times_N M.
$$
\end{enumerate}
The same conclusions hold if we replace $g : M \to N$
by a strongly smooth map from a relative K-space
$g : (X',Z';\widehat{\mathcal U'}) \to N$.
\end{lem}
\subsection{Boundary and corner I}
\label{subsec:bdrycorn1}
So far we study the case when our Kuranishi structures do not
have boundary or corner.
Its generalization to the case when our Kuranishi structure and/or
the manifold $M$ has boundary or corner is straightforward.
We however state them for the completeness' sake.
Later we need to and will study boundary and corner more systematically.
(See Subsections \ref{subsec:kuranishiemb2}, \ref{subsection:normbdry2},
and Part 2.)
\begin{defn}
An orbifold with corner is a space locally
homeomorphic to $V/\Gamma$ where $V$ is a smooth
manifold with corner and $\Gamma$ is a finite group acting
smoothly and effectively on $V$.
We assume the smoothness of the coordinate change
as usual.
\end{defn}
See Definition \ref{orbifolddefn} for more precise and detailed definition.
\begin{defn}\label{defn4111}
Let $M$ be an orbifold with corner.
It has the following canonical stratification $\{S_k(M)\}$.
The stratum $S_k(M)$ is the closure of the set of the points whose neighborhoods are
diffeomorphic to  open neighborhoods of $0$
of the space $([0,1)^k \times \R^{n-k})/\Gamma$.
We call this stratification the {\it corner structure stratification}
\index{corner ! corner structure stratification of an orbifold}
\index{orbifold ! corner structure stratification}
\index{stratification ! corner structure stratification of an orbifold} of $M$.
\par
It is easy to see that $\overset{\circ}S_k(M) = S_k(M) \setminus S_{k+1}(M)$ carries a structure of a smooth
orbifold of dimension $n-k$ without boundary or corner.
However this orbifold may not be effective.
In this document, we {\it assume} the next condition in addition
as a part of the definition of orbifold with corners.
\end{defn}
\begin{conven}\label{effectivitycorner}(Corner effectivity hypothesis)
\index{corner ! corner effectivity assumption}
When we say $M$ is an orbifold with corners,
we assume the orbifold $\overset{\circ}S_k(M)$ is always an {\it effective} orbifold
in this document.
\end{conven}
\begin{defn}\label{dimstratifidef}
For a relative K-space $(X,Z;\widehat{\mathcal U})$, we put
\begin{equation}
\aligned
S_k(X,Z;\widehat{\mathcal U})
&=
\{p \in Z \mid o_p \in S_k(U_p)\},
\\
\overset{\circ}S_k(X,Z;\widehat{\mathcal U})
&=
S_k(X,Z;\widehat{\mathcal U})
\setminus
\bigcup_{k' > k}S_{k'}(X,Z;\widehat{\mathcal U}),
\endaligned
\end{equation}
where $\mathcal U_p = (U_p,E_p,s_p,\psi_p)$ is the
Kuranishi neighborhood of $p$.
We call this stratification the
{\it corner structure stratification}
\index{corner ! corner structure stratification of Kuranishi structure}
\index{stratification ! corner structure stratification of Kuranishi structure}
\index{Kuranishi structure ! corner structure stratification}
of $\widehat{\mathcal U}$.
We can define corner structure stratification $\{S_k(X,Z;\widetriangle{\mathcal U})\}$
of a good coordinate system $\widetriangle{\mathcal U}$ in the same way.
\end{defn}
\begin{lem}
For any compact subset $K$ of $\overset{\circ}S_k(X,Z;\widehat{\mathcal U}) \subseteq X$, the Kuranishi structure $\widehat{\mathcal U}$ induces
a Kuranishi structure without boundary of
dimension $\dim (X,Z;\widehat{\mathcal U}) - k$ on
$K \subseteq \overset{\circ}S_k(X,Z;\widehat{\mathcal U})$.
\par
The same conclusion holds for good coordinate system.
\end{lem}
\begin{proof}
We put
$$
\overset{\circ}{\mathcal S}_k (\mathcal U_p) =
\left(\overset{\circ}S_k (\mathcal U_p), E_p\vert_{\overset{\circ}S_k (\mathcal U_p)},
s_p\vert_{\overset{\circ}S_k (\mathcal U_p)}, \psi_p\vert_{\overset{\circ}S_k (\mathcal U_p)}\right).
$$
Then
we define a Kuranishi neighborhood of
$K \subseteq \overset{\circ}S_k(X,Z;\widehat{\mathcal U})$
at $p$ by
$
\overset{\circ}{\mathcal S}_k (\mathcal U_p).
$
Suppose $q = o_p(q) \in \psi_p(s_p^{-1}(0)) \cap Z$. Then
$q \in S_{k'}(X,Z;\widehat{\mathcal U})$ if and only
if $o_p(q) \in S_{k'}(U_p)$.
Using this fact, we can restrict coordinate changes to
$
\overset{\circ}{\mathcal S}_k (\mathcal U_p)
$
to obtain desired coordinate changes. The compatibility conditions
follow from ones of $\widehat{\mathcal U}$.
\end{proof}
\begin{rem}
\begin{enumerate}
\item
In general the above Kuranishi structure on $K \subseteq \overset{\circ}S_k(X,Z;\widehat{\mathcal U})$
may not be orientable even if $\widehat{\mathcal U}$ is orientable.
\item
In case $k=1$  the above Kuranishi structure of $K \subseteq \overset{\circ}S_1(X,Z;\widehat{\mathcal U})$
is orientable if $\widehat{\mathcal U}$ is orientable.
\item
The Kuranishi structure induced to the normalized corner
of  $(X,Z;\widehat{\mathcal U})$
(see Part 2) is orientable if  $\widehat{\mathcal U}$ is orientable.
\end{enumerate}
\end{rem}
\begin{defn}\label{defn417}
Let $M_1$ and $M_2$ be smooth orbifolds with corner,
$N$ a smooth orbifold without boundary or corner and
$f_i : M_i \to N$ smooth maps.
We say that $f_1$ is {\it transversal}
\index{transversality ! of maps from orbifolds with corners}
to $f_2$ if for each $k_1, k_2$ the
restriction $f_1 : \overset{\circ}S_{k_1}(M_1) \to N$  is
transversal to $f_2 : \overset{\circ}S_{k_2}(M_2) \to N$.
\par
We can define the case of strongly continuous maps from relative K-spaces
with corners to a manifold in the same way.
\index{transversality ! of maps from Kuranishi structure with corners}
The case of good coordinate system is the same.
\end{defn}

\begin{lem}\label{fiberkurabdry}
\begin{enumerate}
\item
Suppose that $Z \subseteq X$ has a Kuranishi structure with boundary and/or corner
and $\widehat f : (X,Z;\widehat{\mathcal U}) \to N$ is weakly transversal to $f' : M \to N$.
Then the fiber product
$Z {}\times_{N} M \subseteq X {}\times_{N} M$ has a Kuranishi structure with
corner.
\item
If $Z_i \subseteq X_i$ has a Kuranishi structure with boundary and/or corner
and $\widehat f_i : (X_i,Z_i;\widehat{\mathcal U_i}) \to N$ a
strongly smooth map to a manifold.
Suppose they are weakly transversal to each other.
Then the fiber product
$Z_1 \times_N Z_2 \subseteq X_1 \times_N X_2$ has a Kuranishi structure with corners.
\end{enumerate}
\end{lem}
The proof is immediate from definition.
\begin{defn}
We call the Kuranishi structure obtained in Lemma \ref{fiberkurabdry},
the {\it fiber product Kuranishi structure}.
\index{Kuranishi structure ! fiber product of Kuranishi structures with corners}
\index{fiber product ! of Kuranishi structures with corners}
\end{defn}
\subsection{Basic property of fiber product}
\label{subsec:fibbasic}
One important property of fiber product is its
associativity, which we state below
\footnote{The fiber product in the sense of category theory is
always associative if it exists. Since we do {\it not} study morphism between K-spaces, the fiber product we defined is {\it not} a fiber
product in the sense of category theory. Therefore we need to prove
its associativity. However it is obvious in our case.}.
We consider the following situation.
Suppose $(X_i,Z_i)$ have Kuranishi structures for $i=1,2,3$ and let
$\widehat f_1 : (X_1,Z_1;\widehat{\mathcal U_1}) \to M_1$,
$\widehat f_2 = (\widehat f_{2,1},\widehat f_{2,2}) : (X_2,Z_2;\widehat{\mathcal U_2}) \to M_1 \times M_2$,
$\widehat f_3 : (X_3,Z_3;\widehat{\mathcal U_3}) \to M_2$  be maps  which are weakly smooth.
We assume $\widehat f_1$ is transversal to $\widehat f_{2,1}$
and $\widehat f_{2,2}$ is transversal to $\widehat f_{3}$.

\begin{lem}\label{lemassoc}
In the above situation, the
following three conditions are equivalent.
\begin{enumerate}
\item
The map $\widehat f_{3} : (X_3,Z_3;\widehat{\mathcal U_3}) \to M_2$ is transversal
to the map $\widehat f'_{2,2} : (X_1,Z_1;\widehat{\mathcal U_1})\times_{M_1} (X_2,Z_2;\widehat{\mathcal U_2}) \to M_2$,
which is induced by $\widehat f_{2,2}$.
\item
The map
$\widehat f_{1} : (X_1,Z_1;\widehat{\mathcal U_1}) \to M_1$ is weakly transversal
to the map $\widehat f'_{2,1} : (X_2,Z_2;\widehat{\mathcal U_2})\times_{M_2} (X_3,Z_3;\widehat{\mathcal U_3}) \to M_1$,
which is induced by $\widehat f_{2,1}$.
\item
The map
\begin{equation}\label{form4545}
(\widehat f_1,\widehat f_2,\widehat f_3) :
(X_1,Z_1;\widehat{\mathcal U_1}) \times (X_2,Z_2;\widehat{\mathcal U_2}) \times (X_3,Z_3;\widehat{\mathcal U_3}) \to M_1^2 \times M_2^2
\end{equation}
is weakly transversal to
$$
\Delta = \{(x_1,x_2,y_1,y_2) \in M_1 \times M_1 \times M_2 \times M_2 \mid x_1=x_2,\,\, y_1=y_2\}.
$$
Here we use the direct product Kuranishi structure in the left hand side of  (\ref{form4545}).
\end{enumerate}
In case those three equivalent conditions are satisfied, we have
\begin{equation}\label{associa}
\aligned
&\left((X_1,Z_1;\widehat{\mathcal U_1}) \times_{M_1} (X_2,Z_2;\widehat{\mathcal U_2})\right) \times_{M_2} (X_3,Z_3;\widehat{\mathcal U_3})
\\
&\cong
(X_1,Z_1;\widehat{\mathcal U_1}) \times_{M_1} \left((X_2,Z_2;\widehat{\mathcal U_2}) \times_{M_2} (X_3,Z_3;\widehat{\mathcal U_3})\right).
\endaligned
\end{equation}
\end{lem}
Here the isomorphism $\cong$ in (\ref{associa}) is defined as follows.
\begin{defn}\label{defniso}
Suppose $(X_1,Z_1;\widehat{\mathcal U_1})$ and $(X_2,Z_2;\widehat{\mathcal U_2})$ are relative K-spaces.
Let $f : (X_1,Z_1) \to (X_2,Z_2)$ be a homeomorphism.
An {\it isomorphism} of relative K-spaces between $(X_1,Z_1;\widehat{\mathcal U_1})$ and $(X_2,Z_2;\widehat{\mathcal U_2})$
assigns the maps
$f_p, \hat f_p$ to each $p \in X_1$ such that
the following holds.
Let $\mathcal U_p^1, \mathcal U_{f(p)}^2$ be the Kuranishi charts of $p, f(p)$ in $X_1, X_2$, respectively.
\index{K-space ! isomorphism of K-spaces}
\begin{enumerate}
\item
$
f_p : U^1_p \to U^2_{f(p)}
$
is a diffeomorphism of orbifolds.
\item
$\hat f_p : \mathcal E^1_p \to \mathcal E^2_{f(p)}$ is a
bundle isomorphism over $f_p$.
\item
$s^p_2 \circ f_p = \hat f_p \circ s^p_1$.
\item
$\psi_{f(p)}^2 \circ f_p = \psi_{p}^1$ on $s_{p}^{-1}(0)$.
\item
$f_p(o^1_p) = o^2_p$.
\end{enumerate}
\end{defn}
\begin{rem}
This definition of isomorphism is too restrictive to be a natural notion
of isomorphism between Kuranishi structures.
To find a correct notion of morphisms between K-spaces
and of isomorphism between them is interesting and is a highly nontrivial
problem. We do not study it here since it is not necessary for our purpose.
A slightly better notion is an equivalence as germs of Kuranishi structures.
See \cite{Fu1}.
\end{rem}
The proof of Lemma \ref{lemassoc} is easy and is omitted.
\par\medskip
In the previous literature such as \cite[Section A1.2]{fooobook2}
we defined a fiber product using the notion of
good coordinate system.
There is one difficulty in defining the fiber product with a space equipped
with good coordinate system, which we explain below.
\par
For $i=1,2$, suppose that $X_i$ have good coordinate systems that are defined by
$\frak P_i$,
$\mathcal U_{p_i}
= (U^i_{\frak p},\mathcal E^i_{\frak p},\psi^i_{\frak p},s^i_{\frak p})$,
and
$\Phi_{\frak p_i\frak q_i}
= (U^i_{\frak p\frak q},\widehat\varphi^i_{\frak p\frak q},\varphi^i_{\frak p\frak q})$.
Let
$\widehat{f}_i = \{(f_i)_{\frak p}\} : (X_i,\mathcal U_{p_i}) \to Y$ be strongly smooth maps.
We assume that $f_{1,\frak p_1} : U^1_{\frak p_1} \to Y$
is transversal to $f_{2,\frak p_2} : U^2_{\frak p_2} \to Y$
for each $\frak p_1 \in \frak P_1$ and $\frak p_2 \in \frak P_2$.
Then we define
$$
U_{(\frak p_1,\frak p_2)} = U^1_{\frak p_1}\times_Y U^1_{\frak p_2}
$$
and define other objects $\mathcal E_{(\frak p_1,\frak p_2)},\psi_{(\frak p_1,\frak p_2)},s_{(\frak p_1,\frak p_2)}$
by taking fiber product in a similar way and to define a good coordinate system.
\par
This is written in \cite[Section A1.2]{fooobook2}.
A point to take care of in this construction is as follows.
(This point is mentioned in \cite[Remark 10 page 165]{fooo010}
and is discussed in detail by Joyce in \cite{joyce}.)
\par
Let $\frak p_i, \frak q_i \in \frak P_i$ such that $\frak q_i < \frak p_i$.
We assume that the fiber product
$$
(U^1_{\frak p_1\frak q_1} \cap (\frak s^1_{\frak q_1})^{-1}(0)) \times_Y (U^2_{\frak p_2\frak q_2} \cap (\frak s^2_{\frak q_2})^{-1}(0))
$$
is nonempty. Then we have
$$
\psi_{(\frak q_1,\frak p_2)}(\frak s^{-1}_{(\frak q_1,\frak p_2)}(0)) \cap
\psi_{(\frak p_1,\frak q_2)}(\frak s^{-1}_{(\frak p_1,\frak q_2)}(0))
\ne \emptyset.
$$
On the other hand, neither $(\frak q_1,\frak p_2) \le (\frak q_2,\frak p_1)$ nor
$(\frak q_2,\frak p_1) \le (\frak q_1,\frak p_2)$.
In fact it may happen that
$
\dim U_{\frak q_i} < \dim U_{\frak p_i}
$. In such a case there is no way to define
$\varphi_{(\frak q_1,\frak p_2),(\frak p_1,\frak q_2)}$ or
$\varphi_{(\frak p_1,\frak q_2),(\frak q_1,\frak p_2)}$.
\begin{rem}
Note the same problem already occurs while we study the direct product.
\end{rem}
We can resolve this problem by shrinking $U_{(\frak p_1,\frak p_2)}$ appropriately.
Joyce \cite{joyce} gave a beautiful canonical way to perform this shrinking
process so that the resulting fiber product is associative.
(See also  \cite[Figure 14]{fooo010}.)
\par
In case we have a multisection (multivalued perturbation)
on the Kuranishi structures on $X_i$ so that the fiber product over $Y$ is
transversal on its zero set,
we use the fiber product of this multisection.
This is especially important when we work in the chain level.
Joyce \cite{joyce}  did not seem to discuss this point since, for his
purpose in \cite{joyce}, it is unnecessary.
We have no doubt that we can incorporate the construction of multisection
to Joyce's fiber product so that we can perturb the Kuranishi structure in a way consistent with the
fiber product and is also associativity of fiber product holds together with
perturbation.
\par
However, in this article we take a slightly different way.
We define the fiber product among the spaces with Kuranishi structures themselves not  those
with good coordinate systems.
Then the above mentioned problem does not occur.
In other words, Kuranishi chart (of the fiber product) is defined
as the fiber product of Kuranishi charts
{\it without shrinking it}.
(Lemmata \ref{lem25}, \ref{lem27}, \ref{fiberkura22}.)
Moreover associativity holds obviously. (Lemma \ref{fiberkurabdry}).
\par
On the other hand, the compatibility of the multisection with
fiber product still needs to be taken care of.
In fact, to find a multisection with appropriate transversality properties,
we used a good coordinate system.
So we need to perform certain process to move from a good coordinate system
back to a Kuranishi structure together with multisections on it.
We will discuss this point in Sections \ref{sec:thick} - \ref{sec:contfamilyconstr}.
\par
\section{Thickening of a Kuranishi structure}
\label{sec:thick}

\subsection{Background of introducing the notion of thickening}
\label{subsec:thickback}
Let $X$ be a paracompact metrizable space, and let
$\widehat{\mathcal U}
= (\{\mathcal U_p\},\{\Phi_{pq}\})$ be a Kuranishi structure on it.
We consider a system of multisections\footnote{We
will discuss multisection in Section \ref{sec:multisection}.
Here we just mention it to motivate the definition
we give in this section. The readers who
do not know the definition of mutisection can safely skip
the part before Definition \ref{stratadim}.}
$\{ \frak s_p \}$ of the vector bundle $\mathcal E_p \to U_p$
for each $p$ with the following property:
\begin{enumerate}
\item[($\bigstar$)]
For each $p$ and $q \in {\rm Im} (\psi_{p})$
the pullback of $\frak s_p$ to $U_{pq}$ that is a
multisection of $\varphi_{pq}^*\mathcal E_p$ is the image of
the multisection $\frak s_q$ by the bundle
embedding $\widehat{\varphi}_{pq}$.
\end{enumerate}
This is a kind of  obvious condition of multisection (multivalued perturbation)
that is compatible with the Kuranishi structure $\widehat{\mathcal U}$.
(We define such a notion precisely later in Definition \ref{compapertKuranishi}.)
\par
However, we need to note the following:
Let us take a good coordinate system ${\widetriangle{\mathcal U}}=
((\frak P,\le) ,
\{\mathcal U_{\frak p} \mid \frak p\in \frak P\},
\{\Phi_{\frak p\frak q} \mid \frak p,\frak q\in \frak P, \frak q \le \frak p\})$ such that
${\widetriangle{\mathcal U}}$ is compatible
with $\widehat{\mathcal U}$ in the sense of
Definition \ref{gcsystemcompa} and use it to define
a system of multisections $\frak s_{\frak p}$
on $U_{\frak p}$.
Then it is usually impossible to use $\frak s_{\frak p}$
to obtain a system of multisections $\frak s_p$ of  $\widehat{\mathcal U}$
that has property ($\bigstar$).
\par
The reason is as follows.
Let $p\in X$. We take $\frak p \in \frak P$ such that
$p \in {\rm Im}(\psi_{\frak p})$.
By definition of the compatibility of good coordinate system
and Kuranishi structure, there exists
an embedding $\Phi_{\frak p p} : \mathcal U_p \to
\mathcal U_{\frak p}$. We put
$\Phi_{\frak p p} = (\varphi_{\frak p p},\widehat{\varphi}_{\frak p p})$.
Then we consider
$\varphi_{\frak p p}(o_p) \in U_{\frak p}$.
(Here $\psi_p(o_p) = p$.)
By definition (See Definition \ref{defn62}.)
\begin{equation}
\frak s_{\frak p}(o_p) \in (\mathcal E_{\frak p}\vert_{o_p})^l.
\end{equation}
On the other hand, $\widehat{\varphi}_{\frak p p}$ restricts to a linear embedding
$
\mathcal E_p\vert_{o_p} \to \mathcal E_{\frak p}\vert_{\psi_{\frak p p}(o_p)}.
$
It induces
\begin{equation}\label{32formula}
(\mathcal E_p\vert_{o_p})^l \to (\mathcal E_{\frak p}\vert_{\psi_{\frak p p}(o_p)})^l.
\end{equation}
By inspecting the construction of
the multisection $\frak s_{\frak p}$ given in \cite[p 955]{FO},  we find that
\begin{equation}\label{33formula}
\frak s_{\frak p}(\varphi_{\frak p p}(o_p)) \notin
{\rm Im} (\ref{32formula})
\end{equation}
in general.
So $\frak s_{\frak p}$ cannot be pulled back to a multisection of $\mathcal E_p$
on $U_{\frak p p}$.
\par
To explain the
reason why (\ref{33formula}) occurs, we introduce some notations.
\begin{defn}\label{stratadim}
For a Kuranishi structure $\widehat{\mathcal U}$ of $Z \subseteq X$, we define
the {\it dimension stratification}
\index{stratification ! dimension stratification}
\index{dimension stratification} of $Z$ by
\begin{equation}\label{defstratifi}
\mathcal S_{\frak d}(X,Z;\widehat{\mathcal U})
=
\{ p \in Z \mid \dim U_p \ge \frak d\}.
\end{equation}
Here $\frak d \in \Z_{\ge 0}$.
\par
When ${\widetriangle{\mathcal U}}$ is a good coordinate
system of $Z \subseteq X$, we define the dimension
stratification of $Z$
by
\begin{equation}\label{defstratifi2}
\mathcal S_{\frak d}(X,Z;{\widetriangle{\mathcal U}})
=
\{ p \in Z \mid \exists \,\frak p, \, \dim U_{\frak p} \ge \frak d,\,
p \in {\rm Im} \psi_{\frak p}\}.
\end{equation}
\end{defn}
\begin{lem}
\begin{enumerate}
\item
$\mathcal S_{\frak d}(X,Z;\widehat{\mathcal U})$,
is a closed subset of $Z$.
$\mathcal S_{\frak d}(X,Z;{\widetriangle{\mathcal U}})$ is an open subset of $Z$.
Moreover if $\frak d' < \frak d$, then
$$
\mathcal S_{\frak d}(X,Z;\widehat{\mathcal U}) \subseteq
\mathcal S_{\frak d'}(X,Z;\widehat{\mathcal U}),
\mathcal \qquad
\mathcal S_{\frak d}(X,Z;\widetriangle{\mathcal U}) \subseteq
\mathcal S_{\frak d'}(X,Z;\widetriangle{\mathcal U}).
$$
\item
If $\widehat{\mathcal U}$ is embedded into $\widehat{\mathcal U'}$ then
$$
\mathcal S_{\frak d}(X,Z;\widehat{\mathcal U}) \subseteq
\mathcal S_{\frak d}(X,Z;\widehat{\mathcal U'}).
$$
The equality holds if and only if the embedding from $\widehat{\mathcal U}$ to $\widehat{\mathcal U'}$
is an open embedding.
\par
The same holds for GG-embeddings.
(The equality in
$
\mathcal S_{\frak d}(X,Z;\widetriangle{\mathcal U}) \subseteq
\mathcal S_{\frak d}(X,Z;\widetriangle{\mathcal U'})
$ holds if the embedding is a strongly open embedding.)
\par
If the good coordinate system ${\widetriangle{\mathcal U}}$
is compatible with $\widehat{\mathcal U}$, then
\begin{equation}\label{dimandUUU}
\mathcal S_{\frak d}(X,Z;\widehat{\mathcal U}) \subseteq
\mathcal S_{\frak d}(X,Z;{\widetriangle{\mathcal U}}).
\end{equation}
\end{enumerate}
\end{lem}
The proof is immediate from definition.
However,
we note that the equality almost never holds in
(\ref{dimandUUU}).
Namely
$\mathcal S_{\frak d}(X,Z;\widetriangle{\mathcal U})$
contains an open neighborhood of
$\mathcal S_{\frak d}(X,Z;\widehat{\mathcal U})$.
This is the reason why (\ref{33formula}) occurs.
\subsection{Definition of thickening}
\label{subsec:thickening}
To go around this problem we introduce the notion of
thickening.
\begin{defn}\label{thickening}
Let $\widehat{\mathcal U}$
be a  Kuranishi structures of   $Z \subseteq X$.
We say that $(\widehat{\mathcal U^+},\widehat\Phi)$
is a {\it thickening}
\index{thickening of Kuranishi structure} of $\widehat{\mathcal U}$
if the following condition is satisfied.
\begin{enumerate}
\item $\widehat{\mathcal U^+}$
is a Kuranishi structure of  $Z \subseteq X$ and
$\widehat\Phi : \widehat{\mathcal U} \to \widehat{\mathcal U^+}$
is a KK-embedding.
\item
For each  $p \in Z$ there exists
a neighborhood $O_p$ of $p$ in
$\psi_p((s_p)^{-1}(0))
\cap \psi^+_p((s^+_p)^{-1}(0)) \subset X$ with the following properties.
\par
For each $q \in O_p \cap Z$ there exists a neighborhood $W_p(q)$
of $o_{p}(q)$ in $U_{p}$
such that:
\begin{enumerate}
\item
$\varphi_p(W_p(q)) \subseteq \varphi^+_{pq}(U_{pq}^+)$.
\item
For any $x \in W_p(q)$, $y\in U_{pq}^+$ with
$ \varphi_p(x) = \varphi^+_{pq}(y)$, we have
$$
\widehat\varphi_p(E_p\vert_x)
\subseteq
\widehat\varphi^+_{pq}(E^+_q\vert_y).
$$
\end{enumerate}
We sometimes say $\widehat{\mathcal U^+}$ is a
thickening of $\widehat{\mathcal U}$ by an abuse of notation.
\par
We write $\widehat{\mathcal U} < \widehat{\mathcal U^+}$
if  $\widehat{\mathcal U^+}$ is a
thickening of $\widehat{\mathcal U}$.
\end{enumerate}
\end{defn}
\begin{equation}
\xymatrix{
U^+_{pq}\ar[rr]^{\varphi^+_{pq}}   &&    U^+_p \\
W_p(q)    \ar@{.>}[u] \ar@{^{(}-{>}}[rr]  &&
U_p \ar[u]_{\varphi_p}
\\
q  \ar@{|->}[d]  &  &s_p^{-1}(0) \ar@{^{(}-{>}}[u] \ar[d]_{\psi_p} \\
O_p\ar@{^{(}-{>}}[rr] && X
}
\end{equation}
\begin{rem}
Condition (2) above implies that
\begin{enumerate}
\item[(\FiveStarOpen)]
$S_{\frak d}(X,Z;\widehat{\mathcal U^+})$
is a neighborhood of $S_{\frak d}(X,Z;\widehat{\mathcal U})$.
\end{enumerate}
In fact if $q \in O_p \cap Z$, then
$$
\dim U^{+}_q \ge \dim W_q(p) = \dim U_p
$$
by Condition (2)(a).
In particular, $\widehat{\mathcal U}$ is almost never a thickening of itself.
\par
In the case $\dim U^{+}_p$ is strictly greater than
$\dim U_p$ and $\dim U^+_q$, Condition (\FiveStarOpen) may not imply
Condition (2)(a).
\end{rem}

\begin{lem}
If $(\widehat{\mathcal U^+},\widehat\Phi)$ is a thickening of
$\widehat{\mathcal U}$ and
$(\widehat{\mathcal U^{++}},\widehat{\Phi^+})$
is a thickening of
$\widehat{\mathcal U^+}$,
then
$(\widehat{\mathcal U^{++}},\widehat{\Phi^+}\circ \widehat\Phi)$
is a thickening of $\widehat{\mathcal U}$.
\end{lem}
\begin{proof}
Let $O_p$, $W_p(q)$ be as in Definition \ref{thickening} (2) (a)
for $\widehat\Phi : \widehat{\mathcal U} \to \widehat{\mathcal U^+}$
and let $O^+_p$, $W^+_p(q)$ be one for
$\widehat{\Phi^+} : \widehat{\mathcal U^+} \to \widehat{\mathcal U^{++}}$.
We put $O_p^{++} = O_p \cap O^+_p$
and
$W_p^{++}(q) = W_p(q) \cap \varphi_p^{-1}(W^+_p(q))$.
Suppose $q \in O_p^{++}$. Then we have
$$
(\varphi_p^+\circ\varphi_p)(W_p^{++}(q))=
\varphi_p^+(\varphi_p(W_p^{++}(q)))
\subseteq \varphi_p^+(W_p^+(q))
\subset \varphi_{p}^{++}(U_{pq}^{++}).
$$
Thus we have checked Definition \ref{thickening} (2) (a).
The proof of Definition \ref{thickening} (2) (b) is similar.
\end{proof}
\subsection{Existence of thickening}
\label{subsec:thickexi}
We next prove the existence of thickening.
We need some notation.
\begin{defn}\label{situ61}
Let
${\widetriangle {\mathcal U}}$
be a good coordinate system of $Z \subseteq X$.
\begin{enumerate}
\item
A {\it support system}
\index{support ! support system}
\index{good coordinate system ! support system}
of ${\widetriangle {\mathcal U}}$ is $
\mathcal K = \{\mathcal K_{\frak p}\mid {\frak p
\in \frak P}\}$
where $\mathcal K_{\frak p} \subset U_{\frak p}$ is a compact subset  for each
$\frak p \in \frak P$ such that
it is a closure of  an open subset $\overset{\circ}{{\mathcal K}_{\frak p}}$
of $ U_{\frak p}$,
and
\begin{equation}\label{form5.65.6}
\bigcup_{\frak p \in \frak P}
\psi_{\frak p}(\overset{\circ}{\mathcal K}_{\frak p} \cap s_{\frak p}^{-1}(0))
\supseteq Z.
\end{equation}
\item
A {\it support pair} $(\mathcal K^1,\mathcal K^2)$ is a pair of support systems $(\mathcal K^i_{\frak p})_{\frak p
\in \frak P}$ $i=1,2$,  such that
\index{support ! support pair}
\index{good coordinate system ! support pair}
\begin{equation}
{\mathcal K}_{\frak p}^1 \subset \overset{\circ}{{\mathcal K}_{\frak p}^2}.
\end{equation}
\par
We write  $\mathcal K^1 < \mathcal K^2$ if $(\mathcal K^1,\mathcal K^2)$
is a support pair.
\item
When $\mathcal K$ is a support system,
we define
$$
\vert \mathcal K\vert = \left(\coprod_{\frak p \in \frak P} \mathcal K_{\frak p}\right) / \sim.
$$
Here, for $x \in \mathcal K_{\frak p}$, $y \in \mathcal K_{\frak q}$,
the relation $x \sim y$ is defined by:
$x = \varphi_{\frak p \frak q}(y)$ or $y = \varphi_{\frak q \frak p}(x)$.
On $\vert \mathcal K\vert$, we put the induced topology from $\vert \widehat{\mathcal U}\vert$.
Then it follows from the definition and \cite[Proposition 5.17]{foootech}
or \cite[Proposition 5.1]{foooshrink},
that the space $\vert \mathcal K\vert$
is metrizable.
\item
When $\mathcal K$ is a support system we define
\begin{equation}
\aligned
\mathcal S_{\frak p}(X,Z;{\widetriangle {\mathcal U}};\mathcal K)
&=
\bigcup_{\frak q \ge \frak p}\psi_{\frak q}(s_{\frak q}^{-1}(0)\cap \mathcal K_{\frak q})\cap Z, \\
\overset{\circ}{\mathcal S}_{\frak p}(X,Z;{\widetriangle {\mathcal U}};\mathcal K)
&=
\mathcal S_{\frak p}(X,Z;{\widetriangle {\mathcal U}};\mathcal K)
\setminus
\bigcup_{\frak q > \frak p}\mathcal S_{\frak q}(X;{\widetriangle{\mathcal U}};\mathcal K).
\endaligned
\end{equation}
\end{enumerate}
\end{defn}

\begin{prop}\label{exithicken}
For any Kuranishi structure $\widehat{\mathcal U}$ there exists its thickening.
\end{prop}
\begin{proof}
By Theorem \ref{Them71restate}
we have a good coordinate system ${\widetriangle{\mathcal U}}
= (\frak P,\{\mathcal U_{\frak p}\},\{\Phi_{\frak p \frak q}\})$
compatible with $\widehat{\mathcal U}$
and its support pair $(\mathcal K^-,\mathcal K^+)$.
By compatibility, there exists a KG-embedding
$\widehat{\mathcal U} \to {\widetriangle{\mathcal U}}$
which we denote by
$\widehat\Phi = \{\Phi_{\frak p p}\mid p\in X,\frak p \in \frak P, p\in \psi_{\frak p}(s^{-1}_{\frak p}(0))\}$.
\par
The first step is to define
$\widehat{\mathcal U^+}$.
\begin{lem}
Let $p \in Z$. There exists unique $\frak p \in \frak P$
such that
$
p \in
\overset{\circ}{\mathcal S}_{\frak p}(X,Z;{\widetriangle {\mathcal U}};\mathcal K^-).
$
\end{lem}
\begin{proof}
By definition
$\overset{\circ}{\mathcal S}_{\frak p}(X,Z;{\widetriangle {\mathcal U}};\mathcal K^-)$ are
disjoint from one another for different $\frak p$'s. By (\ref{form5.65.6}) they cover $Z$.
This finishes the proof.
\end{proof}
We take an open neighborhood $U^+_p$ of $o_{\frak p}(p)
\in \mathcal K^-_{\frak p}$ in
$\mathcal K^+_{\frak p}$ such that
\begin{equation}\label{eq5959}
\psi_{\frak p}(s^{-1}_{\frak p}(0) \cap U^+_{p})
\cap \psi_{\frak q}(s_{\frak q}^{-1}(0) \cap \mathcal K^-_{\frak q}) \ne \emptyset
\,\,\Rightarrow\,\,
\frak q \le \frak p.
\end{equation}
Such a neighborhood exists by Condition (6) of Definition \ref{gcsystem}.
We define
$$
\mathcal U^+_p =
\mathcal U_{\frak p}\vert_{U^+_p}.
$$
We next define a coordinate change.
Let $q \in \psi_{\frak p}(s^{-1}_{\frak p}(0) \cap U^+_p)$.
Take the unique $\frak q \in \frak P$ such that
$
q \in
\overset{\circ}{\mathcal S}_{\frak q}(X,Z;{\widetriangle {\mathcal U}};\mathcal K^-)
$.
Since $q \in \psi_{\frak p}(s^{-1}_{\frak p}(0) \cap U^+_{p})
\cap \psi_{\frak q}(s_{\frak q}^{-1}(0) \cap \mathcal K^-_{\frak q})$,
(\ref{eq5959}) implies $\frak q \le \frak p$.
We put
\begin{equation}\label{U+pq}
U^+_{pq} = U_{q}^+ \cap U_{\frak p\frak q}
\cap \varphi_{\frak p\frak q}^{-1}(U_p^+).
\end{equation}
This is a subset of $U_q^+$ $\subset U_{\frak q}$ and contains $o^+_q =
o_{\frak q}(q)$.
We define
$$
\Phi^+_{pq} = \Phi_{\frak p\frak q}\vert_{U^+_{pq}}.
$$
Clearly $\Phi^+_{pq}$ is a coordinate change from $U_{q}^+$
to $U_{p}^+$.
Using Definition \ref{gcsystem} applied to ${\widetriangle{\mathcal U}}$,
we can easily show that $\mathcal U^+_p$ and $\Phi^+_{pq}$
define a Kuranishi structure.
We denote it by $\widehat{\mathcal U^+}$.
\par
We next define an open substructure $\widehat{\mathcal U_0}$ of
$\widehat{\mathcal U}$ and a strict embedding
$\widehat{\mathcal U_0} \to \widehat{\mathcal U^+}$.
Let $p \in Z$ and we take $\frak p$ such that
$
p \in
\overset{\circ}{\mathcal S}_{\frak p}(X,Z;{\widetriangle{\mathcal U}};\mathcal K^-)
$. We put
$$
U^0_p = \varphi_{\frak p p}^{-1}(U^+_p), \qquad
\mathcal U^0_p = \mathcal U_p\vert_{U^0_p}.
$$
By restricting the coordinate change
$\Phi_{pq}$ of  $\widehat{\mathcal U}$ to
$\varphi_{pq}^{-1}(U^0_p) \cap U_{qp} \cap U_q^0$,
we obtain $\widehat{\mathcal U_0}$ that is an
open substructure of $\widehat{\mathcal U}$.
Then
$
\Phi_{p} = \Phi_{\frak p p}\vert_{U^0_p}
$
is defined.
\par
Definition \ref{defn311} $\circledast$ follows from the fact that
$\widehat{\Phi} : \widehat{\mathcal U} \to {\widetriangle{\mathcal U}}$
is a KG-embedding.
\par
We finally prove that $\widehat{\mathcal U^{+}}$ is a thickening.
Let $
p \in
\overset{\circ}{\mathcal S}_{\frak p}(X,Z;{\widetriangle {\mathcal U}};\mathcal K^-)
$.
We choose $O_p$, a neighborhood of $p$ in $X$ so that the following
condition ($\divideontimes$) is satisfied.
\begin{enumerate}
\item[($\divideontimes$)]
If $O_p \cap \psi_{\frak q}(s_{\frak q}^{-1}(0) \cap \mathcal K_{\frak q}^-) \ne \emptyset$, then
$p \in \psi_{\frak q}(s_{\frak q}^{-1}(0) \cap
{\rm Int}\,\mathcal K_{\frak q}^+)$.
\end{enumerate}
Let $q \in O_p \cap Z$. We take $\frak q$ such that
$
q \in
\overset{\circ}{\mathcal S}_{\frak q}(X,Z;{\widetriangle {\mathcal U}};\mathcal K^-)
$.
By Condition ($\divideontimes$) we have
$p \in \psi_{\frak q}(s_{\frak q}^{-1}(0) \cap
{\rm Int}\,\mathcal K_{\frak q}^+)$.
Therefore there exists an embedding
$
\Phi_{\frak q p} : \mathcal U_p \to \mathcal U_{\frak q}.
$
\par
Recall
from \eqref{U+pq} that $q \in U^+_{\frak p \frak q} \subset U_{\frak p \frak q}$.
We put
$$
W_{p}(q) = \varphi_{\frak q p}^{-1}(U_q^+) \cap U^0_p \cap \varphi_{\frak p \frak q}^{-1}(U_{\frak p \frak q}).
$$
This is an open subset of $U_p^0$ and contains $o_p(q)$.
Now we have
$$\varphi_{\frak q p}(W_p(q)) \subset U^+_q \cap \varphi_{\frak q p}( \varphi_{\frak p p}^{-1}(U^+_p)) \cap U_{\frak p \frak q}.$$
Since
$$
\varphi_{\frak p \frak q} \bigl( \varphi_{\frak q p}( \varphi_{\frak p p}^{-1}(U^+_p)) \cap U_{\frak p \frak q} \bigr) \subset U^+_p,$$
we have
$$
\varphi_{\frak p p}(W_p(q)) = \varphi_{\frak p \frak q}(\varphi_{\frak q p}(W_p(q))
\subset  \varphi_{\frak p \frak q}(U_{p q}^+)
= \varphi_{pq}(U_{p q}^+).
$$
Thus we have proved Definition \ref{thickening} (2)(a).
\par
Using $\widehat \varphi_{\frak q p}$ we can prove
Definition \ref{thickening} (2)(b) in the same way.
The proof of Proposition \ref{exithicken} is complete.
\end{proof}

\subsection{Embedding of Kuranishi structures II}
\label{subsec:kuranishiemb2}

\begin{defn}\label{embgoodtokura}
Let ${\widetriangle{\mathcal U}} = (\frak P,\{\mathcal U_{\frak p}\},\{\Phi_{\frak p\frak q}\})$ be a
good coordinate system of $Z \subseteq X$ and $\widehat{\mathcal U^+}$  a Kuranishi structure of $Z \subseteq X$.
A {\it GK-embedding} $\widehat\Phi : {\widetriangle{\mathcal U}} \to \widehat{\mathcal U^+}$ is
a collection $\{(U_{\frak p}(p),\Phi_{p\frak p})\}$ with the following properties.
\index{embedding ! GK-embedding}
\begin{enumerate}
\item $(U_{\frak p}(p),\Phi_{p \frak p})$ is defined when $\frak p \in \frak P$ and $p \in \psi_{\frak p}(s_{\frak p}^{-1}(0)) \cap Z$.
\item
$U_{\frak p}(p)$ is an open neighborhood of $o_{\frak p}(p)$ in $U_{\frak p}$ where $\psi_{\frak p}(o_{\frak p}(p)) = p$.
\item
$\Phi_{p \frak p} : \mathcal U_{\frak p}\vert_{U_{\frak p}(p)} \to {\mathcal U}^+_p$ is an embedding of Kuranishi charts.
\item
If $\frak q \in \frak P$, $\frak q \le \frak p$, $q \in \psi_{\frak p}(U_{\frak p}(p) \cap s_{\frak p}^{-1}(0))$ and
$q \in \psi_{\frak q}(s_{\frak q}^{-1}(0)) \cap Z$,
then $q \in \psi_p^+(U_p^+ \cap (s_{p}^+)^{-1}(0))$ and the following diagram commutes.
\end{enumerate}
\begin{equation}\label{diagram58}
\begin{CD}
\mathcal U_{\frak q}\vert_{
\varphi_{\frak p\frak q}^{-1}(U_{\frak p}(p))\cap (\varphi_{q\frak q}^+)^{-1}({U}_{pq}^+)} @ > {\Phi_{q\frak q}} >>
{\mathcal U}^+_{q}\vert_{{U}_{pq}^+}  \\
@ V{\Phi{\frak p\frak q}}VV @ VV{\Phi_{pq}^+}V\\
\mathcal U_{\frak p}\vert_{U_{\frak p}(p)} @ > {\Phi_{p\frak p}} >>{\mathcal U}_{p}^+
\end{CD}
\end{equation}
\end{defn}
\begin{rem}
\begin{enumerate}
\item
The case $p=q$ (but $\frak p \ne \frak q$) is included in Definition \ref{embgoodtokura} (4).
The case $\frak p =  \frak q$ (but $p \ne q$) is also included.
\item
Note $q \in \psi_p^+(U_p^+ \cap (s_{\frak p}^+)^{-1}(0))$ (in Definition \ref{embgoodtokura} (4)) follows
from the assumptions ($\frak q \in \frak P$ and $q \in \psi_{\frak q}(s_{\frak q}^{-1}(0))$, $\frak q \le \frak p$)
and Definition \ref{embgoodtokura} (3).
\item
We include $U_{\frak p}(p)$ as a part of the data to define an embedding.
We sometimes need to shrink it.
Such a process is included in the discussion below.
\end{enumerate}
\end{rem}
\begin{defnlem}
We can pullback a strongly continuous (resp. strongly smooth) map
from Kuranishi structure to one from a good coordinate system
by a GK-embedding.
Weak submersivity is preserved by an open embedding.
\end{defnlem}
The proof is immediate from the definition.
\par
Below we define compositions of embeddings
of various types. See the table below.
In the tables below, KS = Kuranishi structure, GCS = good coordinate system.
\par\medskip
\begin{center}
\begin{tabular}{|c|c|c|c|c|c|}
\hline
source & target & symbol
& Definition & name & comment \\
\hline
{KS} & {KS} &
$\widehat{\mathcal U} \to \widehat{\mathcal U^+}$
& Definition \ref{defn311} & KK-embedding & (1)\\
{KS} & {GCS} &
$\widehat{\mathcal U} \to \widetriangle{\mathcal U}$ & Definition
\ref{defn32020202} & KG-embedding & (2)\\
{GCS} & {KS} &
$\widetriangle{\mathcal U} \to \widehat{\mathcal U^+}$ & Definition \ref{embgoodtokura}
& GK-embedding & (3)\\
{GCS} & {GCS} &  $\widetriangle{\mathcal U} \to \widetriangle{\mathcal U^+}$ & Definition \ref{defn31222} & GG-embedding & (4)\\
\hline
\end{tabular}
\end{center}
$$
{\mbox {Table 5.1}: \text{Definition of embedding}}
$$
\noindent
Comments:
(1) Strict version and non-strict version exist. (2)
Strict version and non-strict version exist. (3) None. (4) None.
\par\medskip
\begin{center}
\begin{tabular}{|c|c|c|c|c|}
\hline
1st structure & 2nd structure & 3rd structure
& definition & comment \\
\hline
{KS} & {KS} & {KS}
& Definition \ref{definition516161} & (1)\\
{KS} & {KS} & {GCS}  & Definition \ref{defn32020202} & (2)\\
{KS} & {GCS} & {KS}  & Lemma \ref{lem513} & (3)\\
{KS} & {GCS} & {GCS}  & Definition \ref{defn32020202} & (4)\\
{GCS} & {KS} & {KS}
& Definition-Lemma \ref{kuragoodkuracomp} & (5)\\
{GCS} & {KS} & {GCS}  & Definition-Lemma \ref{henacomp} & (6)\\
{GCS} & {GCS} & {KS}  & Definition-Lemma \ref{kuragoodkuracomp} & (7)\\
{GCS} & {GCS} & {GCS}  & Definition \ref{defn31222} & (8)\\
\hline
\end{tabular}
\end{center}
$$
{\mbox {Table 5.2}: \text{Composition of embeddings}
\atop
\text{1st structure} \to \text{2nd structure} \to  \text{3rd structure} }
$$
\par\smallskip
\noindent
Comments:
(1) Strict version and non-strict version exist. Composition
of non-strict version is well-defined only up to equivalence. (2)
Strict version and non-strict version exist.
Composition
of non-strict version is well-defined only up to equivalence.
(Definition \ref{definition516161} (3))
(3) Need to restrict to an open substructure. The composition becomes
a thickening.
(4) Strict version and non-strict version exist.
(5) None.
(6) Need to restrict to a weakly open substructure.
(7)  None. (8) None.
\begin{defn}
\begin{enumerate}
\item
Let $\widehat\Phi = \{(U_{\frak p}(p),\Phi_{\frak p p})\} : {\widetriangle{\mathcal U}} \to \widehat{\mathcal U^+}$ be a
GK-embedding.
Suppose for each $\frak p \in \frak P$ and $p \in \psi_{\frak p}(s_{\frak p}^{-1}(0)) \cap Z$ we are given
an open subset $U'_{\frak p}(p)$ of $U_{\frak p}(p)$ such that $o_{\frak p}(p) \in U'_{\frak p}(p)$.
Then $\{(U'_{\frak p}(p),\Phi_{\frak p p}\vert_{U'_{\frak p}(p)})\}$ is also an embedding ${\widetriangle{\mathcal U}} \to \widehat{\mathcal U^+}$.
We call it an {\it open restriction of the embedding} $\widehat\Phi$.
\index{embedding ! open restriction of the embedding}
\item
Two embeddings ${\widetriangle{\mathcal U}} \to \widehat{\mathcal U^+}$ are said to be {\it equivalent}
if they have a common open restriction. This is obviously an equivalence relation.
\end{enumerate}
\end{defn}
\begin{defnlem}\label{kuragoodkuracomp}
Let $\widehat\Phi : {\widetriangle{\mathcal U}} \to \widehat{\mathcal U^+}$ be a  GK-embedding,
$\widetriangle\Phi_1 : {\widetriangle{\mathcal U'}} \to {\widetriangle{\mathcal U}}$
a GG-embedding and
$\widehat\Phi_2 : \widehat{\mathcal U^+} \to \widehat{\mathcal U^{++}}$  a strict KK-embedding.
Then we can define the composition
$$
{\widetriangle{\mathcal U'}} \longrightarrow {\widetriangle{\mathcal U}} \longrightarrow
\widehat{\mathcal U^+} \longrightarrow \widehat{\mathcal U^{++}}
$$
which is a GK-embedding.
\end{defnlem}
The proof is easy and left to the reader.
\par

\begin{lem}\label{lem513}
If
$\widehat{\Phi^1} : {\widehat{\mathcal U}} \to {\widetriangle{\mathcal U}}$ is a KG-embedding
and
$\widehat{\Phi^2} : {\widetriangle{\mathcal U}} \to \widehat{\mathcal U^+}$
is  a GK-embedding, then we can find an
open sub-structure  ${\widehat{\mathcal U_0}}$ of ${\widehat{\mathcal U}}$ such that
the composition of
$
{\widehat{\mathcal U_0}} \to {\widehat{\mathcal U}}
\to {\widetriangle{\mathcal U}}$
and
${\widetriangle{\mathcal U}} \to \widehat{\mathcal U^+}
$
is defined as a strict KK-embedding
$:
{\widehat{\mathcal U_0}}  \to \widehat{\mathcal U^+}
$.
\par
$\widehat{\mathcal U^+}$ is a thickening of ${\widehat{\mathcal U_0}}$.
\end{lem}
\begin{proof}
Replacing $\widehat{\mathcal U}$ by its open substructure,
we may assume that $\widehat{\Phi^1}$ is a strict KG-embedding.
Let $p \in Z$.
We define $U^0_{p} \subset U_p$
by
$$
U^0_{p}  = \bigcap_{\frak p : p \in \psi_{\frak p}(s_{\frak p}^{-1}(0))}
(\varphi_{\frak p p}^1)^{-1}(U_{\frak p}(p)).
$$
Here we define $U_{\frak p}(p)$ by
$\Phi^2 = \{(U_{\frak p}(p),\Phi^2_{p \frak p})\}$.
We use them to define our open substructure ${\widehat{\mathcal U_0}} = {\widehat{\mathcal U}}\vert_{\{U_p^0\}}$.
Then
$$
\Phi^2_{p\frak p}\circ\Phi^1_{\frak p p}\vert_{U^0_{p}}  : \mathcal U_p\vert_{U^0_{p}} \to \mathcal U^+_p
$$
is well-defined and defines the required embedding.
\par
We can prove that
$\widehat{\mathcal U^+}$ is a thickening of ${\widehat{\mathcal U_0}}$
in the same way as the proof of Proposition \ref{exithicken}.
\end{proof}
Lemma \ref{lem513} implies that there exists a KK-embedding
${\widehat{\mathcal U}}
\to {\widehat{\mathcal U^+}}$ in the situation of Lemma \ref{lem513}.
This embedding is well-defined in the following sense.
\begin{defn}\label{embedkuraequiv}
Let ${\widehat{\mathcal U}}$ and ${\widehat{\mathcal U^+}}$ be Kuranishi structures.
Suppose ${\widehat{\mathcal U_{0,i}}}$ $(i=1,2)$ are open substructures of ${\widehat{\mathcal U}}$
and $\widehat\Phi^{0,i} : {\widehat{\mathcal U_{0,i}}} \to {\widehat{\mathcal U^+}}$ are strict KK-embeddings.
We say they are {\it equivalent} if there exists an open neighborhood $(U^{00})_{p}$ of $o_p$ in $U_p$
such that
$$
U^{00}_{p} \subset U^{0,1}_{p} \cap U^{0,2}_{p},
\qquad
\Phi^{0,1}_p\vert_{U^{00}_{p}} = \Phi^{0,2}_p\vert_{U^{00}_{p}}.
$$
Here $U^{0,i}_{p}$ is the Kuranishi neighborhood of $p$
assigned by ${\widehat{\mathcal U_{0,i}}}$.
\par
We can define an equivalence between two KG-embeddings $\widehat{\mathcal U}
\to \widetriangle{\mathcal U}$ in the same way.
\end{defn}
\begin{defn}\label{definition516161}
\begin{enumerate}
\item The {\it composition} of two strict  KK-embeddings
 $\widehat{\mathcal U} \to \widehat{\mathcal U'}$,
$\widehat{\mathcal U'} \to \widehat{\mathcal U''}$ is defined
in an obvious way and it is a strict KK-embedding
$\widehat{\mathcal U} \to \widehat{\mathcal U''}$.
\item
We compose
two KK-embeddings
$\widehat{\mathcal U} \to \widehat{\mathcal U'}$,
$\widehat{\mathcal U'} \to \widehat{\mathcal U''}$
and obtain a  KK-embedding $\widehat{\mathcal U} \to \widehat{\mathcal U''}$.
The composition is well-defined up to equivalence defined
in Definition \ref{embedkuraequiv}.
\item
The composition of a KK-embedding $\widehat{\mathcal U} \to \widehat{\mathcal U'}$
and a KG-embedding $\widehat{\mathcal U'} \to \widetriangle{\mathcal U}$
in Definition \ref{defn32020202} is well-defined up to equivalence.
\end{enumerate}
\end{defn}
\begin{lem}\label{deflem517}
Let $\widehat\Phi : {\widetriangle{\mathcal U}} \to \widehat{\mathcal U^+}$ be a  GK-embedding and $\widehat{\mathcal U^+_0}$ an open substructure
of $\widehat{\mathcal U^+}$.
Then there exists a GK-embedding $\widehat{\Phi^0} : {\widetriangle{\mathcal U}} \to \widehat{\mathcal U^+_0}$
such that the composition
${\widetriangle{\mathcal U}} \to \widehat{\mathcal U^+_0} \to \widehat{\mathcal U^+}$
is an open restriction of $\widehat\Phi$.
\end{lem}
\begin{proof}
Let $\widehat\Phi = \{(U_{\frak p}(p),\Phi_{p \frak p})\}$.
We put
$
U^0_{\frak p}(p) = \varphi_{p \frak p}^{-1}(U_p^0) \subset U_{\frak p}(p).
$
We define $ \{(U^0_{\frak p}(p),\Phi_{p \frak p}\vert_{U^0_{\frak p}(p)})\}$ and obtain the required embedding.
\end{proof}
The composition of embeddings in another case is slightly nontrivial.
\begin{defnlem}\label{henacomp}
{\it Let $\widehat{\Phi} : {\widetriangle{\mathcal U}} \to \widehat{\mathcal U}$
be a GK-embedding, and
$\widehat{\Phi^+} :{\widehat{\mathcal U}} \to  {\widetriangle{\mathcal U^+}}$
a KG-embedding.
Then there exists a weakly open substructure ${\widetriangle{\mathcal U_0}}$ of
${\widetriangle{\mathcal U}}$ such that
$
{\widetriangle{\mathcal U_0}} \longrightarrow {\widetriangle{\mathcal U}} \overset{\widehat{\Phi}}{\longrightarrow}
\widehat{\mathcal U}
$
and
$
\widehat{\mathcal U} \overset{\widehat{\Phi^+}}{\longrightarrow}  {\widetriangle{\mathcal U^+}}
$
can be composed to a GG-embedding
$
{\widetriangle{\mathcal U_0}} \longrightarrow {\widetriangle{\mathcal U^+}}.
$
}
\end{defnlem}
We will prove Definition-Lemma \ref{henacomp} in Subsection
\ref{subsec:proofofwdofcomp}, where we use it.
\begin{rem}
We may introduce the notion of germ of Kuranishi structures and use it to describe these situations.
See Section \cite{Fu1}.
\end{rem}
\begin{defn}\label{defn517}
Let $\widehat{\mathcal U}$ be a Kuranishi structure on $X$ and $\widehat{\mathcal U^+}$
be its thickening.
We say a good coordinate system ${\widetriangle{\mathcal U}}$ is {\it in between $\widehat{\mathcal U}$
and $\widehat{\mathcal U^+}$} and write
$$
\widehat{\mathcal U} < {\widetriangle{\mathcal U}} < \widehat{\mathcal U^+},
$$
if the following holds.
\begin{enumerate}
\item
There exist embeddings
$\widehat{\mathcal U} \to {\widetriangle{\mathcal U}}$
and  ${\widetriangle{\mathcal U}} \to \widehat{\mathcal U^+}$.
\item
The composition
$\widehat{\mathcal U} \to {\widetriangle{\mathcal U}}
 \to \widehat{\mathcal U^+}$
is equivalent to the given
embedding $\widehat{\mathcal U} \to \widehat{\mathcal U^+}$ in the sense of
Definition \ref{embedkuraequiv}.
\end{enumerate}
\end{defn}
\begin{prop}\label{prop518}
Let $\widehat{\mathcal U}$ be a Kuranishi structure of $Z \subseteq X$ and $\widehat{\mathcal U^+}$ its thickening.
Then there exists a good coordinate system ${\widetriangle{\mathcal U}}$
in between $\widehat{\mathcal U}$
and $\widehat{\mathcal U^+}$.
\end{prop}
We also use  the following version thereof.
\begin{prop}\label{prop519}
Let $\widehat{\mathcal U}$ be a Kuranishi structure of $Z \subseteq X$
and $\widehat{\mathcal U^+_a}$ $(a=1,2)$  thickenings
of $\widehat{\mathcal U}$. Then
there exists a good coordinate system ${\widetriangle{\mathcal U}}$
and embeddings
$$
\widehat{\mathcal U} \to {\widetriangle{\mathcal U}},
\qquad
{\widetriangle{\mathcal U}} \to \widehat{\mathcal U^+_a}
\quad
 (a=1,2)
$$
such that their compositions
$\widehat{\mathcal U} \to {\widetriangle{\mathcal U}} \to
\widehat{\mathcal U^+_a}
$ are equivalent to the
given embedding $\widehat{\mathcal U} \to \widehat{\mathcal U^+_a}$.
\end{prop}
\begin{equation}
\xymatrix{
&&&  \widehat{\mathcal U^+_1}   \\
\widehat{\mathcal U} \ar[rrru] \ar@{.>}[rr] \ar@{>}[rrrd] &&  \widetriangle{\mathcal U} \ar@{.>}[ru] \ar@{.>}[rd] \\
&&& \widehat{\mathcal U^+_2}
}
\end{equation}
The proofs are almost the same as the proof of Theorem \ref{Them71restate}.
We will review and prove them later in Subsection \ref{subsec:constgcsrel}.

\section{Multivalued perturbation}
\label{sec:multisection}

\subsection{Multisection}
\label{subsec multi}

\subsubsection{Multisection on an orbifold}
\label{subsubsec multiofd}

We next define the notion of multivalued perturbations
associated to a given good coordinate system.
We will slightly  modify the previously given definition to
make explicit certain properties which we used
to study its zero set (in  \cite[Section 2.6]{foootech} for example.)\footnote
{Note the definition of multisection we use here
exactly the same  as one of the smooth multisection in \cite{FO}.}
\par
We begin with a review of multisections.
We first introduce
certain notations on vector bundles on orbifolds.
See Section \ref{sec:ofd} for detail.
\begin{defn}\label{defn61}
Let $U$ be an orbifold and $\mathcal E$  a vector bundle on it.
\begin{enumerate}
\item
(Definitions \ref{2661} (1) and \ref{defn26550} (1)(3))
Let $x \in U$.
We call $(V_x,\Gamma_x,\phi_x)$ an {\it orbifold chart}
of $U$ at $x$ if the following holds.
\begin{enumerate}
\item
$V_x$ is a smooth manifold on which
a finite group $\Gamma_x$ acts effectively and smoothly.
\item $\phi_x : V_x
\to U$ is a $\Gamma_x$-invariant map which induces a
diffeomorphism  $\overline\phi_x :  V_x/\Gamma_x
\to U$
\footnote{See Definition \ref{defn285555} (1) for the definition of diffeomorphism here.} onto an open neighborhood $U_x$ of $x$.
\item
We require that there exists a unique point $o_x \in V_x$ such that
$o_x$ is a fixed point of all elements of $\Gamma_x$ and
$\phi_x([o_x]) = x$.
\end{enumerate}
\item
(Definition \ref{defn2655} (3))
A {\it trivialization of our obstruction bundle $\mathcal E = \widetilde{\mathcal E}/\Gamma_x$ on $U_x$}
is by definition $(E_x,\widehat\phi_x)$
such that
\begin{enumerate}
\item
$E_x$ is a vector space on which $\Gamma_x$ acts linearly.
\item
$\widehat\phi_x
 : V_x \times E_x \to \widetilde{\mathcal E}$
is a $\Gamma_x$-invariant smooth map which induces an isomorphism of vector bundles
$(V_x \times E_x)/\Gamma_x \to \mathcal E\vert_{V_x/\Gamma_x}$.
\end{enumerate}
\item
(Definition \ref{defn2655} (2)(4))
We call  $\frak V_x = (V_x,\Gamma_x,E_x,\phi_x,\widehat\phi_x)$
an {\it orbifold chart} of $(U,\mathcal E)$.
\end{enumerate}
\end{defn}
\begin{defn}\label{defn62}
Let $\frak V_x = (V_x,\Gamma_x,E_x,\phi_x,\widehat\phi_x)$ be
an orbifold chart of $(U,\mathcal E)$.
\begin{enumerate}
\item
A {\it smooth $\ell$-multisection of $\mathcal E$} on an orbifold chart $\frak V_x$
is $\frak s = (\frak s_1,\dots,\frak s_{\ell})$ with the following properties.
\index{multisection ! multisection, smooth $\ell$-multisection on orbifold chart}
\index{orbifold ! multisection, smooth $\ell$-multisection on orbifold chart}
\begin{enumerate}
\item
$\frak s$ is a smooth map $V_x \to E_x^{\ell}$.
\item
For each $y \in V_x$ and $\gamma\in \Gamma_x$ there exists
$\sigma \in {\rm Perm}(\ell)$ such that
$$
\frak s_{\sigma(k)}(y) = \gamma \frak s_{k}(y).
$$
\end{enumerate}
Hereafter we simply say {\it $\ell$-multisection} in place of smooth $\ell$-multisection.
\item
Two  $\ell$-multisections $(\frak s_1,\dots,\frak s_{\ell})$ and $(\frak s'_1,\dots,\frak s'_{\ell})$
are said to be {\it equivalent as $\ell$-multisections on $\frak V_x$}
if  for each $y$ there exists a permutation $\sigma \in {\rm Perm}(\ell)$
such that
$
\frak s'_i(y) = \frak s_{\sigma(i)}(y)
$.
\index{multisection ! equivalence as $\ell$-multisections}
\item
The {\it $\ell'$-iteration} of $\ell$-multisection $\frak s$ is the
$\ell'\ell$-multisection $\frak s'$ such that
$\frak s'_k = \frak s_m$ for $k\equiv m \mod \ell$.
We denote the $\ell'$-iteration of $\frak s$ by $\frak s^{\times \ell'}$.
\index{multisection ! iteration of $\ell$-multisection}
\item
Let $\frak s_{(1)}$ be an $\ell_1$-multisection and $\frak s_{(2)}$ an  $\ell_2$-multisection.
We say $\frak s_{(1)}$ is {\it equivalent to $\frak s_{(2)}$ as multisections} if
$\frak s_{(1)}^{\times \ell_{2}}$ is equivalent to $\frak s_{(2)}^{\times \ell_{1}}$
as $\ell_{1}\ell_{2}$-multisections.
\index{multisection ! equivalence as multisections}
\item
It is easy to see that Item (4) defines an equivalence relation.
We say its equivalence class a {\it multisection} on our orbifold chart.
\index{multisection ! on an orbifold chart}
\index{orbifold ! multisection on an orbifold chart}
\end{enumerate}
\end{defn}

\begin{shitu}\label{opensuborbifoldchart}
Let $\frak V_x = (V_x,\Gamma_x,E_x,\phi_x,\widehat\phi_x)$
and $\frak V'_{x'} = (V'_{x'},\Gamma'_{x'},E'_{x'},\phi'_{x'},\widehat\phi'_{x'})$
be two orbifold charts of a vector bundle $(U,\mathcal E)$.
We assume that $\phi'_{x'}(V'_{x'}) \subset
\phi_{x}(V_{x})$ and that
there exist
$\tilde\varphi_{xx'} : V'_{x'} \to V_{x}$,
$h_{xx'} : \Gamma'_{x'} \to \Gamma_{x}$ such that
$h_{xx'}$ is an injective group homomorphism and
$\tilde\varphi_{xx'}$ is an $h_{xx'}$ equivariant smooth open embedding such that
they induce the composition map
$$
(\overline{\phi_x})^{-1}\circ \overline{\phi'_{x'}} : V'_{x'}/\Gamma'_{x'} \to V_{x}/\Gamma_{x},
$$
where $\overline{\phi_x}$ (resp. $\overline{\phi'_{x'}}$) is induced by $\phi_x$ (resp. $\phi'_{x'}$).
In other words,
$$
\phi'_{x'}(y) \equiv \phi_{x}(\tilde\varphi_{xx'}(y)) \mod \Gamma_x.
$$
Moreover we assume that the composition
$$
\left(\overline{\widehat\phi_x}\right)^{-1}\circ \overline{\widehat\phi'_{x'}} : (V'_{x'}\times E'_{x'})/\Gamma'_{x'} \to (V_{x}\times E_{x})/\Gamma_{x}
$$
is induced by a smooth map
$\breve\varphi_{xx'} : V'_{x'} \times E'_{x'} \to E_x$
that is linear in $E'_{x'}$ factor.
In other words
$$
\widehat\phi'_{x'}(y,v) \equiv \widehat\phi_x(\tilde\varphi_{xx'}(y),\breve\varphi_{xx'}(y,v))
\mod \Gamma_x.
$$
We put
$
\Phi_{xx'} = (h_{xx'},\tilde\varphi_{xx'},\breve\varphi_{xx'}).
$$\blacksquare$
\end{shitu}
\begin{defn}
We call $\Phi_{xx'}$ a {\it coordinate change}
\index{orbifold ! coordinate change of orbifold chart of vector bundles}
from an orbifold chart
$\frak V_{x'}$ to $\frak V_{x}$.
\end{defn}
\begin{rem}
We put
$
\tilde{\hat{\varphi}}_{xx'}(y,v) = (y,\breve\varphi_{xx'}(y,v))
$.
Then
$$
(h_{xx'},\tilde\varphi_{xx'},\tilde{\hat{\varphi}}_{xx'})
$$
is a {\it local representative of an embedding of vector bundles},
\index{embedding ! local representative of embedding of vector bundle on orbifold}
    ${\rm id} : \mathcal E\vert_{{\rm Im}\overline{\widehat{\phi^{\prime}}}_{x'}} \to \mathcal E\vert_{{\rm Im}\overline{\widehat\phi}_{x}}$
in the sense of Definition \ref{def28262826}.
Moreover it is a fiberwise isomorphism.
\end{rem}
\begin{defn}\label{multilocalrest}
Let $\frak s_x$
be an $\ell$-multisection on $\frak V_x$ and
$\Phi_{xx'} : \frak V_{x'} \to \frak V_x$  a coordinate change.
We define the {\it restriction} $\Phi_{xx'}^*\frak s$ by
$$
(\Phi_{xx'}^*\frak s)_k(y)
=
g_y^{-1}\frak s_k(\tilde\varphi_{xx'}(y))
$$
where
$
g_y : E_{x'} \to E_x
$
is defined by
$
g_y(v) = \breve\varphi_{xx'}(y,v).
$
\index{multisection ! restriction by coordinate change}
\end{defn}
\begin{lem}
If $\frak s_x$ is equivalent to $\frak s'_x$ as multisections
then $\Phi_{xx'}^*\frak s$  is equivalent to
$\Phi_{xx'}^*\frak s'$ as multisections.
\end{lem}
We omit the proof.
(See the proof of a similar lemma,
Lemma \ref{lem77}.)
\par
Here is a notational remark.
\begin{rem}\label{xkararhe}
So far we have written $\frak V_x = (V_x,\Gamma_x,E_x,\phi_x,\widehat\phi_x)$.
The point $x$ plays no particular role
except we assume the existence of $o_x \in V_x$
such that $\phi_x(o_x) = x$ and $o_x$ is
fixed by all the elements of $\Gamma_x$.
If we change the choice of such $x$, the constructions so far
do not change at all.
So, from now on, we do not specify $x$ in our notation of local
orbifold chart but only assume an existence of such $x$.
We will write for example
$\frak V_{\frak r} = (V_{\frak r},\Gamma_{\frak r},E_{\frak r},\phi_{\frak r},\widehat\phi_{\frak r})$
instead of $\frak V_{x} = (V_{x},\Gamma_{x},E_{x},\phi_{x},\widehat\phi_{x})$.
\end{rem}\label{lem606969}
\begin{defn}\label{defn6868}
Let $U$ be an orbifold and $\mathcal E$ a vector bundle on it.
\begin{enumerate}
\item
A {\it representative of a multisection of $\mathcal E$ on $U$}
is $(\{\frak V_{\frak r} \mid \frak r \in \frak R\},\{\frak s_{\frak r}
 \mid \frak r \in \frak R\})$
with the following properties.
\begin{enumerate}
\item
$\frak V_{\frak r}$ is an orbifold chart of a vector bundle $(U,\mathcal E)$
such that
$\bigcup_{\frak r \in \frak R}{U_{\frak r}} = U$.
\item
$\frak s_{\frak r}$ is a multisection of $\frak V_{\frak r}$.
\item
For any $y \in \frak V_{\frak r_1} \cap \frak V_{\frak r_2}$,
there exist an orbifold chart $\frak V_y$ and coordinate changes
$\Phi_{\frak r_i y} : \frak V_y \to \frak V_{\frak r_i}$ such that
$\Phi_{\frak r_1 y}^*\frak s_{\frak r_1}$ is equivalent to
$\Phi_{\frak r_2 y}^*\frak s_{\frak r_2}$.
\end{enumerate}
\index{multisection ! representative of multisection of a vector bundle on an orbifold}
\item
Let $(\{\frak V^{(i)}_{\frak r^{(i)}} \mid \frak r^{(i)} \in \frak R^{(i)}\},\{\frak s^{(i)}_{\frak r^{(i)}}
 \mid \frak r^{(i)} \in \frak R^{(i)}\})$
be representatives of multisections of $\mathcal E$ on $U$ for $i=1,2$.
We say they are {\it equivalent} if the following holds.
\par
For any $x \in \frak V^{(1)}_{\frak r^{(1)}_1} \cap \frak V^{(2)}_{\frak r^{(2)}_2}$,
there exist an orbifold chart $\frak V_x$ and coordinate changes
$\Phi_{\frak r_i x} : \frak V_x \to \frak V^{(i)}_{\frak r_i^{(i)}}$
($i=1,2$) such that
$\Phi_{\frak r^{(1)}_1 x}^*\frak s_{\frak r^{(1)}_1}^{(1)}$ is equivalent to
$\Phi_{\frak r^{(2)}_2 x}^*\frak s_{\frak r^{(2)}_2}^{(2)}$.
\par
An equivalence class of this equivalence relation is called
a {\it multisection of $(U,\mathcal E)$}.
\index{multisection ! of a vector bundle on an orbifold}
\index{orbifold ! multisection of a vector bundle on an orbifold}
\item(See \cite[Definition 3.10]{FO}.)
Let $\frak s^{n}$ be a sequence of multisections of $(U,\mathcal E)$.
We say that {\it it converges to a multisection $\frak s$ in $C^{k}$-topology}
($k$ is any of $0,1,\dots,\infty$)
if there exists a representative
$(\{\frak V_{\frak r} \mid \frak r \in \frak R\},\{\frak s^{n}_{\frak r}
\mid \frak r \in \frak R\})$ of $\frak s^{n}$ for sufficiently large $n$ and
$(\{\frak V_{\frak r} \mid \frak r \in \frak R\},\{\frak s_{\frak r}
\mid \frak r \in \frak R\})$ of $\frak s$ such that
$\frak s^{n}_{\frak r}$ converges to $\frak s_{\frak r}$ in compact
$C^k$-topology for each $\frak r$.
\index{multisection ! $C^k$-convergence}
We note that we assume $\frak V_{\frak r}$ and $\frak R$ are
independent of $n$.
\end{enumerate}
\end{defn}
\begin{defn}\label{definition610}
Let $\frak s$ be a multisection of  a vector bundle $(U,\mathcal E)$ on orbifold $U$
and $x \in U$.
We put
$\frak s = [(\{\frak V_{\frak r} \mid \frak r \in \frak R\},\{\frak s_{\frak r}
 \mid \frak r \in \frak R\})]$.
We take an orbifold chart $\frak V_x$ at $x$.
A map germ $[s]$, where $s : O_{x} \to E_x$,  is said to be a
{\it branch of $\frak s$ at $x$} if the following holds.
\begin{enumerate}
\item
$O_x$ is a neighborhood of $o_x$ in $V_x$.
\item
Let $\frak r \in \frak R$  such that
$x \in U_{\frak r}$.
Then there exists $k$ such that
$$
\widehat\phi_{\frak r}(\tilde\varphi_{\frak r x}(y),\frak s_{\frak r,k}(\tilde\varphi_{\frak r x}(y))) = \widehat\phi_{x}(y,s(y))
$$
if $y$ is
on a neighborhood of $x$ in $O_x$. Here
$\tilde\varphi_{\frak r x}$ is a part of a coordinate change
$\frak V_x \to \frak V_{\frak r}$.
\end{enumerate}
\index{multisection ! branch of}
\end{defn}

\subsubsection{Multisection on a good coordinate system}
\label{subsubsec multigcs}
Let
${\widetriangle{\mathcal U}}=(\frak P, \{ \mathcal U_{\frak p}\}, \{\Phi_{\frak p \frak q}\})$
be a good coordinate system of $Z \subseteq X$.

\begin{defn}\label{defn612}
Let
$\mathcal K =  \{\mathcal K_{\frak p}\}$
be a support system of a good coordinate system ${\widetriangle {\mathcal U}}$ and $\frak s_{\frak p}^{n}$  multisections of $E_{\frak p}$ on a neighborhood of $\mathcal K_{\frak p}$ for each $n \in \Z_{\ge 0}$ and $\frak p$.
We say that $\widetriangle{\frak s}
= \{\widetriangle{\frak s^n} \mid n \in \Z_{\ge 0}\} = \{\frak s_{\frak p}^{n} \mid n \in \Z_{\ge 0},
\frak p \in \frak P\}$ is a {\it multivalued perturbation of $({\widetriangle U},\mathcal K)$}
if the following conditions are satisfied.
\begin{enumerate}
\item
$
\frak s_{\frak p}^{n} \circ \varphi_{\frak p\frak q}
=
\widehat\varphi_{\frak p\frak q}\circ \frak s_{\frak q}^{n}
$
on a neighborhood of $\mathcal K_{\frak q}
\cap \varphi_{\frak p,\frak q}^{-1}(\mathcal K_{\frak p})$.
\item
$
\lim_{n\to \infty} \frak s_{\frak p}^{n} = s_{\frak p}
$
in $C^1$-topology on a neighborhood of $\mathcal K_{\frak p}$.
\end{enumerate}
\par
A {\it multivalued perturbation of $\widetriangle{\mathcal U}$}
\index{multivalued perturbation ! of good coordinate system} is
a collection $\{\frak s_{\frak p}^{n}\}$
such that (1)(2) hold for some support system $\mathcal K$.
\end{defn}
Note that the Kuranishi map $s_{\frak p}$, which is a single valued
section of $\mathcal E_{\frak p}$, can be regarded as a multisection
by Lemma \ref{lem2626}.
The $C^1$-convergence in Definition \ref{defn612} (2)
therefore is defined in Definition \ref{defn6868} (3).

Below we will elaborate on the equality in Definition \ref{defn612} (1)
further.
Let $x \in  \mathcal K_{\frak q}
\cap \varphi_{\frak p\frak q}^{-1}(\mathcal K_{\frak p})$
and $x' = \varphi_{\frak p\frak q}(x) \in \mathcal K_{\frak p}$.
We can take orbifold charts $\frak V_{x}$ of $(U_{\frak q},\mathcal E_{\frak q})$,
$\frak V_{x'}$ of $(U_{\frak p},\mathcal E_{\frak p})$
such that $(\varphi_{\frak p\frak q},\widehat\varphi_{\frak p\frak q})$
has a local representative $(h_{\frak p\frak q;x'x},\tilde\varphi_{\frak p\frak q;x'x},
\tilde{\hat{\varphi}}_{\frak p\frak q;x'x})$ with respect to these orbifold charts.
(Lemma \ref{lem2622}.)
We define
$\breve\varphi_{\frak p\frak q;x'x} :
V_x \times E_x \to E_{x'}$ by the relation
$$
\tilde{\hat{\varphi}}_{\frak p\frak q;x'x}(y,v)
=
(\tilde\varphi_{\frak p\frak q;x'x}(y),\breve\varphi_{\frak p\frak q;x'x}(y,v)).
$$
\par
We may choose $\frak V_{x}$ and $\frak V_{x'}$ so small  that
$\frak s_{\frak p}^{n}$ and $\frak s_{\frak q}^{n}$
have representatives on the charts.
Let $\frak s^n_{\frak p,x'}$ and $\frak s^n_{\frak q,x}$
be the representatives, which are $\ell_1$ and $\ell_2$ multisections,
respectively. By taking an appropriate iteration we may assume
$\ell_1 = \ell_2 = \ell$. We then require
\begin{equation}\label{form6161}
\frak s^n_{\frak p,x';k}(\tilde\varphi_{\frak p\frak q;x'x}(y)) =
\breve\varphi_{\frak p\frak q;x'x}(y,\frak s^n_{\frak q,x;\rho_y(k)}(y))
\end{equation}
for $y \in V_x$, $k=1,\dots,\ell$,
where $\rho_y \in {\rm Perm}(\ell)$.
(\ref{form6161}) is the precise form of
Definition \ref{defn612} (1).
\par\medskip
We will use the $C^1$-convergence to prove certain
important properties which we use in the next subsection.
To state this property we need to prepare some notations.
\par
We denote the normal bundle of our embedding
$\varphi_{\frak p\frak q}: U_{\frak p\frak q} \to U_{\frak p}$ by
$$
N(\varphi_{\frak p\frak q};U_{\frak p})
: = \frac{\varphi_{\frak p\frak q}^*TU_{\frak p}}{TU_{\frak q}\vert_{U_{\frak p\frak q}}}
$$
It defines a vector bundle over $U_{\frak p\frak q}$.
For a compact subset $K \subset U_{\frak q}$ we denote
by $N_K(\varphi_{\frak p\frak q};U_{\frak p})$ the restriction of this vector bundle to $K \cap U_{\frak p\frak q}$.
\begin{shitu}\label{shit612712}
We fix a Riamannian metric on $U_{\frak p}$.
It induces a metric on $N(\varphi_{\frak p\frak q};U_{\frak p})$.
We denote by $N^\delta(\varphi_{\frak p\frak q};U_{\frak q})$
the $\delta$-disc bundle thereof for $\delta > 0$.
Using the normal exponential map of the embedding $\varphi_{\frak p\frak q}$,
we have a diffeomorphism:
\begin{equation}
{\rm Exp} : N_{K}^\delta(\varphi_{\frak p\frak q};U_{\frak p}) \to U_{\frak p\frak q} \times U_{\frak p}
\end{equation}
which is given by ${\rm Exp}(x,v) = (x, \exp_x v)$ where $\exp_x$ is the exponential map of
the metric given on $U_{\frak p}$.
(See \cite[Lemma 6.5]{fooooverZ}.)
Here $\delta$ is a positive number depending on a compact
subset $K$ of $U_{\frak q}$ and the embedding $\varphi_{\frak p\frak q}$.
We put
\begin{equation}
BN_{\delta'}(K;U_{\frak p})
= \bigcup_{x \in U_{\frak p\frak q} \cap K}
{\rm exp}_x(N_K^{\delta'}(\varphi_{\frak p\frak q};U_{\frak p})) \subset  U_{\frak p}
\end{equation}
for $\delta' \le \delta$. We denote by $\pi_{\delta'}: BN_{\delta'}(K;U_{\frak p}) \to U_{\frak p\frak q} \cap K$
the composition
of ${\rm Ext}^{-1}$ with
the projection $N(\varphi_{\frak p\frak q};U_{\frak p}) \to U_{\frak q}$
of the vector bundle.
Note that on the image of $\varphi_{\frak p\frak q}$,
the obstruction bundle $\mathcal E_{\frak p}$ has a
subbundle $\widetilde\varphi_{\frak p\frak q}(\mathcal E_{\frak q})$.
We consider a sub-bundle
$$
\mathcal E_{\frak q;\frak p} \subset \mathcal E_{\frak p}
= \pi_{\delta}^*(\widetilde\varphi_{\frak p\frak q}(\mathcal E_{\frak q}))
$$
on $B_{\delta}(K;U_{\frak p})$ which restricts to the bundle
$\widetilde\varphi_{\frak p\frak q}(\mathcal E_{\frak q})$
on $\varphi_{\frak p\frak q}(K) \subset U_{\frak p}$.
We take the quotient bundle
$
\mathcal E_{\frak p}/\mathcal E_{\frak q;\frak p}
$
and consider
$$
\overline{s_{\frak p}} \equiv s_{\frak p}
\mod \mathcal E_{\frak q;\frak p}
$$
that is a section of
$\mathcal E_{\frak p}/\mathcal E_{\frak q;\frak p}$.
In a similar way for
each branch $\frak s^{n}_{\frak p;k}$ of $\frak s^{n}_{\frak p}$
we obtain
\begin{equation}\label{coderivativeoffras}
\overline{\frak s^{n}_{\frak p;k}}(y) \in (\mathcal E_{\frak p}/\mathcal E_{\frak q;\frak p})_y.
\end{equation}
$\blacksquare$
\end{shitu}
We denote by
$E: BN_{\delta}(K;U_{\frak p}) \to N_{K}^\delta(\varphi_{\frak p\frak q};U_{\frak q})$
the inverse of ${\rm Exp}$ on its image. Namely
\begin{equation}\label{invexpmap}
E(y) = (x,v)
\qquad
\Leftrightarrow
\qquad
x = \pi_\delta(y), \, v = \exp_x^{-1}(y).
\end{equation}

\begin{lem}\label{lem611}
Suppose we are in Situation \ref{shit612712}.
There exist $c>0$, $\delta_0>0$ and $n_0 \in \Z_{\ge 0}$ such that
for $y \in BN_{\delta_0}(K;U_{\frak p})$
\begin{equation}\label{normailityestimate}
\vert \overline{\frak s^{n}_{\frak p;k}}(y) \vert \ge c \vert E(y)\vert
\end{equation}
and
\begin{equation}\label{normailityestimate22}
\vert\overline{s_{\frak p}}(y) \vert\ge c \vert E(y)\vert
\end{equation}
hold for any branch $\overline{\frak s^{n}_{\frak p;k}}$, if $n > n_0$ and
$d(\pi_{\delta}(y),s_{\frak q}^{-1}(0)) < \delta_0$.
\end{lem}
\begin{proof}
We choose and fix a connection of the
quotient bundle $\mathcal E_{\frak p}/\mathcal E_{\frak q;\frak p}$.
For $y = {\rm Exp}(x,0) \in \varphi_{\frak p,\frak q}(K)$
we consider the covariant derivative
\begin{equation}\label{normalonimage}
V \mapsto D_V \overline{s_{\frak p}} :  (N_K(\varphi_{\frak p\frak q};U_{\frak p}))_x
\to (\mathcal E_{\frak p}/\mathcal E_{\frak q;\frak p})_x.
\end{equation}
By Definition \ref{defKchart} (5),  the map (\ref{normalonimage}) is an isomorphism
if $x \in \frak s_q^{-1}(0)$ in addition.
Therefore we may choose $\delta_0$ so that
the map (\ref{normalonimage}) is an isomorphism
if $d(x,s_{\frak q}^{-1}(0)) < \delta_0$.
Then the existence of $\delta_0$, $c$ satisfying (\ref{normailityestimate22})
is an immediate consequence of the fact $s_{\frak p}$
is smooth.
The inequality
(\ref{normailityestimate}) then is a consequence of $C^1$-convergence.
\end{proof}
\begin{rem}\label{rem612612}
Note the set theoretical fiber of a vector bundle over an orbifold
is a quotient of a vector space by a finite group.
The value $\overline{s_{\frak p}}(y)$
is well-defined as an element of vector space
if we fix a local trivialization.
When we do not specify the local trivialization,
the value $\overline{s_{\frak p}}(y)$ is defined as an element of a quotient of a
vector space by a finite group.
The left hand side of (\ref{normailityestimate22})
so makes sense.
\par
In case of multisection $\overline{\frak s^{n}_{\frak p;k}}(y)$,
this is well-defined as an element of vector space if we fix
a local trivialization.
When we change the local trivialization it changes by the
permutation of $k$ and a finite group action.
Therefore the validity of
(\ref{normailityestimate}) for all branches $k$ is independent of the
choice of trivialization.
\end{rem}
\begin{rem}\label{rem614}
Note that in Definition \ref{defn62} (2) we allow the permutation $\sigma$ to
 {\it depend} on $y$ which lies in a neighborhood of $x$.
By this reason the notion of branch of
multisection should be studied rather carefully.
Here is an example:
We define
$$
e_{\epsilon_1\epsilon_2}(t)
=
\begin{cases}
\epsilon_1 e^{ -1/\vert t\vert} & t\le 0 \\
\epsilon_2 e^{ -1/\vert t\vert} & t\ge 0.
\end{cases}
$$
Here $\epsilon_i$ is either plus or minus.
Four functions $e_{++}$, $e_{--}$, $e_{-+}$, $e_{+-}$
are all smooth functions on $\R$.
We define 2-multisections $\frak s$ and $\frak s'$ as follows.
$$
\frak s = (e_{++},e_{--}),
\qquad
\frak s' = (e_{+-},e_{-+}).
$$
They are 2-multisections on $\R$ of a trivial line bundle.
It is easy to see that $\frak s$ is equivalent to
$\frak s'$ in the sense of Definition \ref{defn62} (2).
(This definition coincides with \cite{FO}.)
However, it is impossible to choose $\sigma$ appeared in
 Definition \ref{defn62} (2) in a way
independent of $y$.
See Subsection \ref{subsec:nastyreason}
for more discussion about this point.
\end{rem}
It is convenient to remove the assumption $d(\pi_{\delta}(y),s_{\frak q}^{-1}(0)) < \delta_0$
in Lemma \ref{lem611}.
Lemmata \ref{lem6767} and \ref{lem6768} below say that we can always do so by
replacing our good coordinate system
by a strongly open GG-embedding.
\begin{conds}\label{conds6.17}
Let $\widetriangle{\mathcal U}$ be a good coordinate system of $(X,Z)$
and $\mathcal K$ its support system.
We consider the following condition for them.
\par
For each $\frak p > \frak q$, there exist Riemannian metrics on $U_{\frak p}$ and
the sub-bundle $\mathcal E_{\frak q;\frak p}$ as in Situation \ref{shit612712}.
(Here the compact set $K$ appearing in  Situation \ref{shit612712} is  $\varphi_{\frak p\frak q}(\mathcal K_{\frak q})$.)
and we have a connection on $\mathcal E_{\frak p}/\mathcal E_{\frak q;\frak p}$.
They satisfy the following.
\par
If $x \in \varphi_{\frak p\frak q}(\mathcal K_{\frak q})$,
$V \in (N_K(\varphi_{\frak p\frak q};U_{\frak p}))_x$, $V \ne 0$ then
\begin{equation}\label{nmlderi0}
D_V \overline{s_{\frak p}} \ne 0.
\end{equation}
\end{conds}
\begin{rem}
Note the left hand side of (\ref{nmlderi0})  is a covariant derivative
which depends on the choice of the connection in general.
It is independent of the choice of the connection when $x \in s_{\frak p}^{-1}(0))$.
\end{rem}
\begin{lem}\label{lem6767}
Suppose $\widetriangle{\mathcal U}$ and $\mathcal K$
satisfy Condition \ref{conds6.17}.
Then
there exist $c>0$ and $n_0 \in \Z_{\ge 0}$
\begin{equation}\label{normailityestimate3}
\vert \overline{s_{\frak p}}(y) \vert \ge c \vert E(y)\vert
\end{equation}
and
\begin{equation}\label{normailityestimate42}
\vert \overline{\frak s^{n}_{\frak p;k}}(y) \vert \ge c \vert E(y) \vert
\end{equation}
hold for all branch $\overline{\frak s^{n}_{\frak p;k}}$, if $n > n_0$,
$x \in \mathcal K_q$ and $y \in BN_{\delta_0}(K;U_{\frak p})$.
\end{lem}
Here the map $E$ is as in (\ref{invexpmap}).
The proof is the same as the (second half of the) proof of Lemma \ref{lem611}.
\par
We take a metric $d$ on $\vert \mathcal K^+\vert$ in the next lemma and
$B_{\delta}(A) = \{x \in \vert \mathcal K^+\vert  \mid d(x,A) < \delta\}$.
We also regard $Z \subset \vert\mathcal K^+\vert$ by $\psi$.
\begin{lem}\label{lem6768}
Let $\widetriangle{\mathcal U}$ be a good coordinate system
of $(X,Z)$ and $(\mathcal K^-,\mathcal K^+)$ be
its support pair. Then there exists
a support system $\mathcal K^{-\prime}$ and $\delta_0 > 0$ with the following properties.
\begin{enumerate}
\item
$(\mathcal K^{-\prime},\mathcal K^+)$ is a support pair.
\item
$\widetriangle{\mathcal U}$, $\mathcal K^{-\prime}$
satisfy Condition \ref{conds6.17}.
\item
\begin{equation}\label{neweq612}
B_{\delta}(\mathcal K^-_{\frak q} \cap Z)  \cap  Z \subseteq
B_{\delta}(\mathcal K^{-\prime}_{\frak q} \cap Z),
\end{equation}
for any $\frak q \in \frak P$.
\end{enumerate}
\end{lem}
\begin{proof}
Take a support system $\mathcal K^{-+}$ such that
$\mathcal K^- < \mathcal K^{-+} < \mathcal K^{+}$.
Apply (the first half of the proof of) Lemma \ref{lem611}
to $K = \mathcal K_{\frak q}^{-+}$ and obtain $\delta_{\frak q}$ such that
\begin{equation}\label{newform6132}
D_V \overline{s_{\frak p}} \ne 0
\end{equation}
for $V \in (N_K(\varphi_{\frak p\frak q};U_{\frak p}))_x$,
$x \in B_{\delta_{\frak q}}(\mathcal K_{\frak q}^{-+} \cap s_{\frak q}^{-1}(0)) \cap U_{\frak q}$.
\par
Put $\delta_1 = \min\{\delta_{\frak q} \mid \frak q \in \frak P\}$.
Take $\delta_2>0$ such that
\begin{equation}\label{newform613}
B_{2\delta_2}(\mathcal K^-_{\frak q}) \cap U_{\frak q} \subset \overset{\circ}{{\mathcal K}_{\frak q}^{-+}},
\qquad
B_{2\delta_2}(\mathcal K^{-+}_{\frak q}) \cap U_{\frak q} \subset \overset{\circ}{{\mathcal K}_{\frak q}^{+}}
\end{equation}
for all $\frak q$.
We take
\begin{equation}\label{newform614}
\mathcal K^{-\prime}_{\frak q} = {\rm Close}(B_{\delta_2}(\mathcal K^-_{\frak q} \cap Z) \cap U_{\frak q}).
\end{equation}
We take $\delta_0$ smaller than $\min\{\delta_1,\delta_2\}$.
\par
Item (1) follows from (\ref{newform613}).
Item (2) follows from (\ref{newform6132}) and (\ref{newform613}).
\par
We prove hat
Item (3) holds if we replace $\delta_0$ by a smaller positive number if necessary.
We remark that for any $\delta>0$ there exists $\delta' > 0$ such that
\begin{equation}\label{form616}
B_{\delta'}(\mathcal K^-_{\frak q} \cap Z) \cap U_{\frak p}
\subset
N_{{\rm Close}(B_{\delta}(\mathcal K^-_{\frak q} \cap Z) \cap U_{\frak q})}^{\delta}(\varphi_{\frak p\frak q};U_{\frak p}).
\end{equation}
(\ref{form616}) and Lemma \ref{lem611}  implies
\begin{equation}\label{newform617}
B_{\delta'}(\mathcal K^-_{\frak q} \cap Z) \cap U_{\frak p}
\subset U_{\frak q}
\end{equation}
for sufficiently small $\delta'$ and $\frak q < \frak p$. (\ref{neweq612})
follows from (\ref{newform613}), (\ref{newform614}) and (\ref{newform617}).
\end{proof}
\begin{rem}\label{rem61261222}
Note the way taken in \cite{foootech} or in the earlier literatures
such as \cite{FO},\cite{fooobook2}
is different from that in this document and
proceed as follows.
We fix the extension of the subbundle
$\mathcal E_{\frak q;\frak p}$ and fix the
choice of the
splitting
$
\mathcal E_{\frak q}
\equiv
\mathcal E_{\frak q;\frak p} \oplus \frac{\mathcal E_{\frak q}}{\mathcal E_{\frak q;\frak p}}.
$
We then assumed  the equality
\begin{equation}\label{formula611}
\Pi_{\frac{\mathcal E_{\frak q}}{\mathcal E_{\frak q;\frak p}}}(\frak s_{\frak p,k}^n) =
\Pi_{\frac{\mathcal E_{\frak q}}{\mathcal E_{\frak q;\frak p}}}(s_{\frak p})
\end{equation}
for any branch $\frak s_{\frak p,k}^n$ of our multisection $\frak s_{\frak p}^n$.
(See \cite[(6.4.4)]{FO} for example.)\footnote{
We remark that actually we can make sense the equality (\ref{formula611})
without taking and fixing the
splitting
$
\mathcal E_{\frak p}
\equiv
\mathcal E_{\frak q;\frak p} \oplus \frac{\mathcal E_{\frak q}}{\mathcal E_{\frak q;\frak p}}
$, since the projection
$\mathcal E_{\frak p} \to \frac{\mathcal E_{\frak q}}{\mathcal E_{\frak q;\frak p}}$
is well-defined without fixing this splitting.}
We did {\it not} assume $\frak s_{\frak p,k}^n$ converge to
$s_{\frak p}$ in $C^1$-topology but assumed only
$C^0$-convergence.
\par
However (\ref{formula611}) together with the fact
$$
\vert\Pi_{\frac{\mathcal E_{\frak q}}{\mathcal E_{\frak q;\frak p}}}(s_{\frak p})\vert > c\vert v\vert
$$
(which follows from the proof of Lemma \ref{lem6767} (\ref{normailityestimate3}))
is enough to prove (\ref{normailityestimate42}).
\par
We have slightly modified the definition here,
since by assuming $C^1$-convergence as in Definition \ref{defn612} (2)
we can prove Lemma \ref{lem6767}, which we will use instead of
the assumption (\ref{formula611}) made in the earlier literature.
In fact, (\ref{normailityestimate42}) is the property we need.
(See the proof of Sublemma \ref{subsub67} Case 4.)
\par
We however emphasize that the results using
the definition in  the earlier literature is {\it literally correct without change}
by the proof given there. We here are improving the presentation
of the proof of the earlier literatures
but are {\it not} correcting the proof therein.
\end{rem}

\par
In Part 2 of this document, we need to study a family of
multivalued perturbations.
The following notion is useful for the study of a family of
multivalued perturbations.
\begin{defn}\label{uniformmulivalupert}
A $\sigma$ parameterized family of multivalued perturbations
$\{ \{{\frak s}^n_{\sigma}\} \mid \sigma \in \mathscr A\}$
of $(\widetriangle{\mathcal U},\mathcal K)$
is said to be a {\it uniform family}
\index{uniform family ! of multivalued perturbations}
\index{multivalued perturbation ! uniform family} if the convergence in Definition \ref{defn612}
is uniform.
More precisely, we require the following.
\par
For each $\epsilon$ there exists $n(\epsilon)$
such that if $n > n(\epsilon)$ then
\begin{equation}
\vert s(y) - s_{\frak p}(y) \vert < \epsilon,
\qquad
\vert (Ds)(y) - (Ds_{\frak p})(y) \vert < \epsilon
\end{equation}
hold for any branch $s$ of ${\frak s}^n_{\sigma}$
at any point $y \in \mathcal K_{\frak p}$
for any $\frak p \in \frak P$, $\sigma \in \mathscr A$.
\end{defn}
\begin{rem}
Note that we assume that $(\widetriangle{\mathcal U},\mathcal K)$ is
independent of $\sigma$.
\end{rem}
Inspecting the proof of Lemma \ref{lem6768}, we
have the following.
\begin{lem}\label{lem618}
If $\{ \{{\frak s}^n_{\sigma}\} \mid \sigma \in \mathscr A\}$
is a uniform family, then the constants $n_0$, $c$ and $\delta_0$
in Lemma  \ref{lem6768} can be taken independent of $\sigma$.
\end{lem}
\begin{rem}\label{rmCkconverge}
Suppose $\mathscr A$ consists of one point.
Then the condition assumed in Lemma \ref{lem618}
is slightly weaker than $C^1$ convergence in the sense of
Definition \ref{defn6868} (3).
In fact, in Definition \ref{defn6868} (3), we assumed,
for example, that we may choose that the number of
branches of $\frak s^n$ is independent of $n$.
\par
If we define $C^k$ convergence of multisections to
a {\it multisection} using branch in the same way as above,
it seems rather complicated.
Note we discuss here the case of $C^1$ convergence of multisections to
a {\it single-valued} section, the Kuranishi map.
\end{rem}

\subsection{Support system and the zero set of multisection}
\label{subset:supportzeromulti}

We use Lemma \ref{lem6767} to prove Propositions \ref{splem2} and \ref{lem715} below.
In this subsection we use a metric on subsets of
$
\vert \widetriangle{\mathcal U} \vert = \bigcup_{\frak p \in \frak P} U_{\frak p} / \sim
$, which we choose as follows.
We start with support systems $\mathcal K^i$, $i=1,2,3$
with $\mathcal K^1 <\mathcal K^2 <\mathcal K^3$.
(See Definition \ref{situ61} (2) for this notation.)
The union of the images of $\mathcal K_{\frak p}^3$ in
$\vert \widetriangle{\mathcal U}\vert$
is denoted by $\vert \mathcal K^3 \vert$. The quotient topology on
$\vert \mathcal K^3 \vert$ is metrizable.
(See \cite[Proposition 5.1]{foooshrink}.)
We use this topology or its induced topology.
The space $X$ can be regarded as a subspace of
$\vert \mathcal K^3 \vert$ and of $\vert \mathcal K^2 \vert$ or
$\vert \mathcal K^1 \vert$.
We take a metric on a compact neighborhood of $\vert \mathcal K^3\vert$ in
$\vert \widetriangle{\mathcal U}\vert$ and use the induced
metric on various spaces appearing below.
Note that all the spaces $\mathcal K^i_{\frak p}$ etc. are contained in
the compact neighborhood of $\vert \mathcal K^3\vert$ so have this metric.
\par
For a subset $A \subset \vert \mathcal K^3 \vert$ we put
\begin{equation}\label{defmetricball}
B_{\delta}(A) = \{ x \in \vert \mathcal K^3\vert \mid d(x,A) < \delta\}.
\end{equation}
For a point $p \in \vert \mathcal K^3 \vert$,
we define
$B_{\delta}(p) := B_{\delta}(\{p\})$.
\par
Sometimes we identify a subset $\mathcal K^i_{\frak p}$
with its image in  $\vert \mathcal K^3\vert$.
Then for example,
for a subset $A \subset \mathcal K^3_{\frak q} \cap U_{\frak p\frak q}$,
we identify $A$ with $\varphi_{\frak p \frak q}(A)$.
\begin{prop}\label{splem2}
Let $\mathcal K^- < \mathcal K^+ < \mathcal K^2 < \mathcal K^{3}$
and let $\{\frak s_{\frak p}^{n}\}$ be a multivalued perturbation of $(\widetriangle{\mathcal U},\mathcal K^3)$.
Then there exist $\delta >0$, $n_0 \in \Z_{\ge 0}$ such that
for any $\frak q \in \frak P$ and $n >n_0$
\begin{equation}
B_{\delta}(\mathcal K_{\frak q}^- \cap Z) \cap
\bigcup_{\frak p} ((\frak s_{\frak p}^{n})^{-1}(0)
\cap \mathcal K^2_{\frak p})
\subset
\mathcal K_{\frak q}^+.
\end{equation}
\end{prop}
Here $(\frak s_{\frak p}^{n})^{-1}(0)$ is the set of the points in $\mathcal K_{\frak p}^3$ such that at least
one of the branches of $\frak s_{\frak p}^{n}$ vanishes.
\begin{proof}
We first remark that
in view of Lemma \ref{lem6768} it suffices to prove the proposition
when $\mathcal K^-$
satisfies Condition \ref{conds6.17} in addition.
In fact we apply Lemma \ref{lem6768} to obtain $\mathcal K^{-\prime}$
satisfying Condition \ref{conds6.17} in addition.
By using  Lemma \ref{lem6768}, we have
$$
B_{\delta}(\mathcal K_{\frak q}^- \cap Z) \cap
\bigcup_{\frak p} ((\frak s_{\frak p}^{n})^{-1}(0)
\cap \mathcal K^2_{\frak p})
\subset
B_{\delta}(\mathcal K_{\frak q}^{-\prime} \cap Z) \cap
\bigcup_{\frak p} ((\frak s_{\frak p}^{n})^{-1}(0)
\cap \mathcal K^2_{\frak p})
$$
Therefore (\ref{splem2}) with $\mathcal K^-$ replaced by $\mathcal K^{-\prime}$
implies (\ref{splem2}).
\par
Hereafter we assume the condition.
\par
Let $x \in \mathcal K_{\frak q}^- \cap Z$. We first show the following lemma.
\begin{lem}\label{sublem66}
There exist $\delta_{x,\frak q} > 0$ and $n_{x,\frak q} > 0$
such that for $n > n_{x,\frak q}$
$$
B_{\delta_{x,\frak q}}(x) \cap Z \subset \mathcal K_{\frak q}^+,
\quad B_{\delta_{x,\frak q}}(x) \cap
\bigcup_{\frak p} (\frak s_{\frak p}^{n})^{-1}(0)
\subset \mathcal K_{\frak q}^+.
$$
\end{lem}
\begin{proof}
During the proof of Lemma \ref{sublem66} we fix
$x$ and $\frak q$. The constants in Sublemmata \ref{sublem618},
\ref{subsub67} depend on $x$ and $\frak q$.
\begin{sublem}\label{sublem618}
There exists $\delta_1>0$  with the following properties.
\begin{enumerate}
\item
If $x \notin \mathcal K_{\frak p}^2$, then
$B_{\delta_1}(x) \cap  \mathcal K_{\frak p}^2 = \emptyset$.
\item
If $x \in \mathcal K_{\frak p}^-$, then
$B_{\delta_1}(x) \cap
{\mathcal K_{\frak p}^2} = B_{\delta_1}(x) \cap {\mathcal K_{\frak p}^+}$.
\item
If $x \in \mathcal K_{\frak p}^2$, $\frak q \le \frak p$, then
$
B_{\delta_1}(x) \cap \mathcal K_{\frak q}^2 \subset \mathcal K_{\frak p}^{3} \cap \mathcal K_{\frak q}^2.
$
\item
If $x \in \mathcal K_{\frak p}^2$, $\frak q \ge \frak p$, then
$
B_{\delta_1}(x) \cap \mathcal K_{\frak p}^2 \subset \mathcal K_{\frak q}^+.
$
\end{enumerate}
\end{sublem}
\begin{proof}
Statement (1) follows from compactness of $\mathcal K_{\frak p}^2$,
(2) from $\mathcal K_{\frak p}^- \subset {\rm Int}\,
\mathcal K_{\frak p}^+$,
(3) from $\mathcal K^2_{\frak p} \subset {\rm Int}\, {\mathcal K}_{\frak p}^{3}$, and
(4) from $x \in \mathcal K_{\frak q}^- \cap Z$
and $\mathcal K_{\frak q}^- \subset {\rm Int}\,
\mathcal K_{\frak q}^+$, respectively.
\end{proof}
\begin{sublem}\label{subsub67}
There exists $\delta_{2,\frak p}$ for each $\frak p \in \frak P$
and $n_{1,\frak p} \in \Z_{\ge 0}$ such that
$$
\aligned
B_{\delta_{2,\frak p}}(x) \cap (s_{\frak p})^{-1}(0)
\cap  \mathcal K_{\frak p}^2 \subset \mathcal K_{\frak q}^+,
\quad
B_{\delta_{2,\frak p}}(x) \cap  (\frak s_{\frak p}^{n})^{-1}(0)
\cap  \mathcal K_{\frak p}^2 \subset \mathcal K_{\frak q}^+
\endaligned
$$
hold for $n > n_{1,\frak p}$.
\end{sublem}
\begin{proof}
We discuss 4 cases separately.
In the first 3 cases we will prove
\begin{equation}\label{form6969}
B_{\delta_{2,\frak p}}(x) \cap \mathcal K^2_{\frak p}
\subset \mathcal K_{\frak q}^+.
\end{equation}
(\ref{form6969}) obviously  implies the required inclusion in those cases.
\par\smallskip
\noindent(Case 1)
Neither $\frak p \le \frak q$ nor $\frak q \le \frak p$.
In this case we may choose $\delta_{2,\frak p} = \delta_1$ since the left hand side of (\ref{form6969}) is an empty set by Sublemma \ref{sublem618} (1).
\par\smallskip
\noindent(Case 2)
$\frak p = \frak q$. We take $\delta_{2,\frak p} = \delta_1$. Then (\ref{form6969}) follows from Sublemma \ref{sublem618} (2).
\par\smallskip
\noindent(Case 3)
$\frak p < \frak q$.
We take $\delta_{2,\frak p} = \delta_1$. Then (\ref{form6969}) follows from Sublemma \ref{sublem618} (4).
\par\smallskip
\noindent(Case 4)
$\frak q < \frak p$.
This is the most important case.
Let $d_{\frak p}$ be a metric function induced by a Riemannian
metric $g_{\frak p}$ of $U_{\frak p}$.
Note the metric $d$ which we used to define the metric ball
in Sublemma \ref{sublem618} is different from $d_{\frak p}$.
However they define the same topology.
We next prove that there exist, $\delta_3>0$ and $c>0$ such that
\begin{equation}\label{form67}
\vert s_{\frak p}(y)\vert \ge c d_{\frak p}(y,\mathcal K^2_{\frak q})
\end{equation}
holds for $y \in B_{\delta_{3}}(x) \cap \mathcal K^2_{\frak p}$.
To prove this we show the next subsublemma.
\begin{subsublemma}\label{sublem628}
There exists $\delta_3 > 0$ such that
if $y \in B_{\delta_3}(x)\cap \mathcal K_{\frak p}^2$
then there exists a minimal  $g_{\frak p}$-geodesic $\ell  : [0,d] \to U_{\frak p}$ of length $d$ such that
$\ell(0) \in \varphi_{\frak p\frak q}(U_{\frak p\frak q})\cap \mathcal K_{\frak p}^{3} \cap \mathcal K_{\frak q}^2$ and
$d = d_{\frak p}(y,\mathcal K^2_{\frak q})$.
\end{subsublemma}
\begin{proof}
Since  ${\mathcal K}_{\frak q}^{2}$ is a relatively compact subset of
 ${\mathcal K}_{\frak q}^{3}$ and $x \in {\mathcal K}_{\frak q}^{2}$,
there exists $\delta_3 > 0$  such that $d_{\frak p}(x,y) < \delta_3$ implies
that there exists $z \in \mathcal K^3_{\frak q}$ with the property that
the minimal geodesic joining $z$ and $y$ is perpendicular to $\mathcal K^3_{\frak q}$
and that  $d_{\frak p}(z,y) \le d_{\frak p}(x,y)$.
We can now use Sublemma \ref{sublem618} (3) to show that we may choose $\delta_3$ such that
$z \in \varphi_{\frak p\frak q}(U_{\frak p\frak q})\cap \mathcal K_{\frak p}^{3} \cap \mathcal K_{\frak q}^2$.
The subsublemma follows.
\end{proof}
The inequality (\ref{form67}) follows from Subsublemma \ref{sublem628} and
Lemma \ref{lem6767} (\ref{normailityestimate3}).
\par
 Now
(\ref{form67}) implies that
$$
B_{\delta_{3}}(x) \cap (s_{\frak p})^{-1}(0) \subset \mathcal K_{\frak q}^2.
$$
It implies
$B_{\delta_{3}}(x) \cap (s_{\frak p})^{-1}(0) \subset \mathcal K_{\frak q}^+$
by Sublemma \ref{sublem618} (2) applied to
$\frak p = \frak q$.
\par
We use Lemma \ref{lem6767} (\ref{normailityestimate42})
and Subsublemma \ref{sublem628} in the same way as
the proof of (\ref{form67}) to show:
\begin{equation}\label{estimateofsepsilon}
\vert \frak s_{\frak p}^{n}(y)\vert \ge c d_{\frak p}(y,U_{\frak q})
\end{equation}
for all $y \in B_{\delta_{4}}(x) \cap \mathcal K^2_{\frak p}$
and sufficiently large $n$.
(Note that this inequality holds for any branch of $\frak s^{n}_{\frak q}$.)
Using (\ref{estimateofsepsilon}) in place of
(\ref{form67}) we prove
$B_{\delta_{4}}(x) \cap  (\frak s_{\frak p}^{n})^{-1}(0) \subset \mathcal K_{\frak q}^2$
in the same way as above.
Then
$B_{\delta_{4}}(x) \cap  (\frak s_{\frak p}^{n})^{-1}(0) \subset \mathcal K_{\frak q}^+$
follows from Sublemma \ref{sublem618} (2).
Thus $\delta_{2,\frak p} = \min\{ \delta_1,\delta_3,\delta_4\}$ has the required properties.
The proof of Sublemma \ref{subsub67} is complete.
\end{proof}
We put $\delta_{x,\frak q} = \min\{\delta_{2,\frak p}\mid \frak p \in \frak P\}$ and
$n_{x,\frak q} = \max \{n_{1,\frak p}  \mid \frak p \in \frak P\}$.
Then Lemma \ref{sublem66} follows from Sublemma \ref{subsub67}.
\end{proof}
Now we take finitely many points $x_i \in \mathcal K_{\frak q}^- \cap Z$, $i =1,\dots, N_{\frak q}$ such that
$$
\bigcup_{i=1}^{N_{\frak q}} B_{\delta_{x,\frak q} }(x_i) \supset \mathcal K_{\frak q}^- \cap Z.
$$
We put $\frak U_{\frak q} = \bigcup_{i=1}^{N_{\frak q}} B_{\delta_{x_i,\frak q}}(x_i)$. Then for any $n \ge \max \{ n_{x_i, \frak q} \mid i = 1, \dots ,N_{\frak q} \}$
we have
$$
\frak U_{\frak q} \supset \mathcal K_{\frak q}^- \cap Z, \qquad
\frak U_{\frak q} \cap  \bigcup_{\frak p} ((\frak s_{\frak p}^{n})^{-1}(0)) \subset \mathcal K_{\frak q}^+.
$$
Since $\frak U_{\frak q}$ is open and $\mathcal K_{\frak q}^- \cap Z$ is compact, there exists
$\delta_{\frak q} > 0$ such that
$
B_{\delta_{\frak q}}(\mathcal K_{\frak q}^- \cap Z) \subset \frak U_{\frak q}.
$
It is easy to see that
$
\delta = \min \{ \delta_{\frak q} \mid \frak q \in \frak P\}
$
and
$
n_0 = \max \{ n_{x_i,\frak q} \mid i=1,\dots, N_{\frak q},\,
\frak q \in \frak P\}
$
have the required properties.
The proof of Proposition \ref{splem2} is complete.
\end{proof}
\begin{prop}\label{lem715}
Let
$\mathcal K^1,\mathcal K^2,\mathcal K^3$
be a triple of support systems of a good coordinate system ${\widetriangle {\mathcal U}}$ of $Z \subseteq X$ with $\mathcal K^1<\mathcal K^2<\mathcal K^3$ and
$\widetriangle{\frak s} = \{\frak s_{\frak p}^{n}\}$  a multivalued perturbation of $({\widetriangle {\mathcal U}},\mathcal K^3)$.
Then there exists a neighborhood $\frak U(Z)$ of
$Z$ in $\vert \mathcal K^2 \vert$ and $n_0\in \Z_{\ge 0}$ such that
the following holds for any $n > n_0$.
\begin{equation}
\left(\bigcup_{\frak p}((\frak s_{\frak p}^{n})^{-1}(0)
\cap \mathcal K^1_{\frak p})\right)
\cap \frak U(Z)
=
\left(\bigcup_{\frak p}((\frak s_{\frak p}^{n})^{-1}(0)
\cap \mathcal K^{2}_{\frak p})\right)
\cap \frak U(Z).
\end{equation}
\end{prop}
\begin{proof}
The inclusion $\subseteq$ is obvious for any $\frak U(Z)$.
We will prove the  inclusion of the opposite direction.
\begin{lem}\label{splem1}
For each $\delta >0$, there exists $\delta' > 0$ such that
for every $\frak p\in \frak P$
\begin{equation}
B_{\delta'}(Z) \cap \mathcal K_{\frak p}^2  \subset B_{\delta}(Z \cap \mathcal K_{\frak p}^2).
\end{equation}
Here we put
$
B_{\delta}(A) = \{ x\in \vert \mathcal K^2 \vert \mid d(x,A) < \delta\}.
$
\end{lem}
\begin{proof}
The proof is by contradiction.
If the lemma does not hold, there exist $\frak p \in \frak P$, $\delta > 0$,
a sequence $\delta_i \to 0$, and points $x_i \in B_{\delta_i}(Z) \cap \mathcal K_{\frak p}^2$
such that $x_i \notin B_{\delta}(Z \cap \mathcal K_{\frak p}^2)$.
\par
Since $\mathcal K_{\frak p}^2$ is compact, we may assume that the sequence $x_i$ converges to
a point $x \in \mathcal K_{\frak p}^2$. Then $x \in \mathcal K_{\frak p}^2 \cap Z$.
Therefore $x _i \in B_{\delta}(Z \cap \mathcal K_{\frak p}^2)$ for sufficiently large $i$.
This is a contradiction.
\end{proof}

By Lemma \ref{splem1}
\begin{equation}\label{form615454}
\aligned
\left(\bigcup_{\frak q}((\frak s_{\frak q}^{n})^{-1}(0)
\cap \mathcal K^{2}_{\frak q})\right)
\cap B_{\delta'}(Z)
&=
\bigcup_{\frak q}\left((\frak s_{\frak q}^{n})^{-1}(0)
\cap \mathcal K^{2}_{\frak q}
\cap B_{\delta'}(Z)\right) \\
&\subseteq
\bigcup_{\frak q}
\left((\frak s_{\frak q}^{n})^{-1}(0)
\cap B_{\delta}(Z \cap \mathcal K^{2}_{\frak q}) \cap \mathcal K^2_{\frak q}\right).
\endaligned
\end{equation}
for sufficiently large $n$.
We take  a support system $\mathcal K^0 = (\mathcal K_{\frak p}^0)_{\frak p \in \frak P}$ such that
$\mathcal K^0 < \mathcal K^1$.
Since $\bigcup_{\frak q} Z \cap \mathcal K^{0}_{\frak q} = Z$, we have
\begin{equation}\label{+next66611}
\bigcup_{\frak q}
\left((\frak s_{\frak q}^{n})^{-1}(0)
\cap B_{\delta}(Z \cap \mathcal K^{2}_{\frak q}) \cap \mathcal K^2_{\frak q}\right)
\subseteq
\bigcup_{\frak p,\frak q}
\left((\frak s_{\frak p}^{n})^{-1}(0)
\cap B_{\delta}(Z \cap \mathcal K^{0}_{\frak q}) \cap \mathcal K^2_{\frak p}\right).
\end{equation}
In fact
$$
\aligned
\bigcup_{\frak q}
\left((\frak s_{\frak p}^{n})^{-1}(0)
\cap B_{\delta}(Z \cap \mathcal K^{0}_{\frak q}) \cap \mathcal K^2_{\frak p}\right)
&=
(\frak s_{\frak p}^{n})^{-1}(0)
\cap B_{\delta}(Z)\cap \mathcal K^2_{\frak p}
\\
&\supseteq
(\frak s_{\frak p}^{n})^{-1}(0)
\cap B_{\delta}(Z \cap \mathcal K^{2}_{\frak p})
\cap \mathcal K^2_{\frak p}.
\endaligned
$$
We apply
Proposition  \ref{splem2} to $(\mathcal K^-,\mathcal K^+) = (\mathcal K^0,\mathcal K^1)$ and
obtain
\begin{equation}\label{form617617}
\bigcup_{\frak p}
\left((\frak s_{\frak p}^{n})^{-1}(0)
\cap B_{\delta}(Z \cap \mathcal K^{0}_{\frak q})
\cap \mathcal K^{2}_{\frak p}\right)
\subset
\mathcal K^{1}_{\frak q}
\end{equation}
for sufficiently large $n$.
Note
$$
\bigcup_{\frak p}
\left((\frak s_{\frak p}^{n})^{-1}(0)
\cap  \mathcal K^{1}_{\frak q}\right)
=
(\frak s_{\frak q}^{n})^{-1}(0)
\cap  \mathcal K^{1}_{\frak q}.
$$
Hence (\ref{form617617}) implies
$$
\bigcup_{\frak p,\frak q}
\left((\frak s_{\frak p}^{n})^{-1}(0)
\cap B_{\delta}(Z \cap \mathcal K^{0}_{\frak q})
\cap \mathcal K^{2}_{\frak p}\right)
\subseteq
\bigcup_{\frak q}
\left((\frak s_{\frak q}^{n})^{-1}(0)
\cap \mathcal K^{1}_{\frak q}\right).
$$
Combined with (\ref{form615454}) and (\ref{+next66611}), we have
$$
\left(\bigcup_{\frak q}(\frak s_{\frak q}^{n})^{-1}(0)
\cap \mathcal K^{2}_{\frak q})\right)
\cap B_{\delta'}(Z)
\subseteq
\bigcup_{\frak q}\left(
(\frak s_{\frak q}^{n})^{-1}(0)
\cap  \mathcal K^{1}_{\frak q}\right).
$$
We take $\frak U(Z) = B_{\delta'}(Z)$.
The proof of Proposition \ref{lem715} is then complete.
\end{proof}

Proposition \ref{lem715} corresponds to \cite[Lemma 6.6]{foootech}. The proof
we gave above is based on the same idea.
\begin{cor}\label{cor69}
There exist a neighborhood $\frak U(Z)$ of $Z$
in $\vert\mathcal K^2\vert$ and $n_0 \in \Z_{\ge 0}$ such that the space
$\left(\bigcup_{\frak p}((\frak s_{\frak p}^{n})^{-1}(0)
\cap \overset{\circ}{\mathcal K^2_{\frak p}})\right)
\cap \frak U(Z)$ is compact if $n > n_0$.
Moreover,
\begin{equation}\label{hdfconv}
\lim_{n\to \infty}
\left(\bigcup_{\frak p}((\frak s_{\frak p}^{n})^{-1}(0)
\cap  \overset{\circ}{\mathcal K^2_{\frak p}})\right)
\cap \frak U(Z)
\subseteq X.
\end{equation}
Here the limit is taken
in Hausdorff topology.
\end{cor}
The first claim corresponds to \cite[Lemma 6.11]{foootech} and the second
claim corresponds to \cite[Lemma 6.12]{foootech}.
Using Proposition \ref{lem715} the proof is also the same
as those lemmata.
We reproduce them here for reader's convenience.
\begin{proof}
Proposition \ref{lem715} implies that
\begin{equation}\label{form619619}
\left(\bigcup_{\frak p}((\frak s_{\frak p}^{n})^{-1}(0)
\cap \overset{\circ}{\mathcal K^2_{\frak p}})\right)
\cap \frak U(Z)
=
\left(\bigcup_{\frak p}((\frak s_{\frak p}^{n})^{-1}(0)
\cap {\mathcal K^1_{\frak p}})\right)
\cap \frak U(Z).
\end{equation}
We may assume that $\frak U(Z)$ is compact.
Then
$
\left(\bigcup_{\frak p}((\frak s_{\frak p}^{n})^{-1}(0)
\cap {\mathcal K^1_{\frak p}})\right)
\cap \frak U(Z)
$ is compact.
The compactness of
$\left(\bigcup_{\frak p}((\frak s_{\frak p}^{n})^{-1}(0)
\cap \overset{\circ}{\mathcal K^2_{\frak p}})\right)
\cap \frak U(Z)$ follows from (\ref{form619619}).
\par
We next prove (\ref{hdfconv}).
Suppose (\ref{hdfconv}) does not hold for any
$\frak U(Z)$.
Then there exist $\frak p \in \frak P$, $\delta > 0$, $n_i \to \infty$, and $x_i$ such that
$x_i\in (\frak s_{\frak p}^{n_i})^{-1}(0)
\cap {\mathcal K^1_{\frak p}}
\cap \frak U(Z)$
and $d(x_i,X) \ge \delta$ for all $i$.
We may assume that $x_i$ converges to $x$.
(Note we may assume that $\frak U(Z)$ is compact.)
Then
$x\in (s_{\frak p})^{-1}(0)
\cap {\mathcal K^1_{\frak p}}
\cap \frak U(Z)$.
Therefore $x \in X$. This contradicts to $d(x,X) \ge \delta > 0$.
\end{proof}
\begin{rem}
To derive the estimate \eqref{estimateofsepsilon} we used
Definition \ref{defn6868} (3) for the notion of $C^1$-convergence of the
multivalued perturbations $\{ \frak s_{\frak p}^n \}$. However,
we note that \eqref{estimateofsepsilon} can be also obtained by using
a slightly weaker notion of $C^1$-convergence defined by Definition \ref{uniformmulivalupert}. See also Remark \ref{rmCkconverge}.
Therefore Propositions \ref{splem2}, \ref{lem715} also hold
even if we use this weaker notion of $C^1$-convergence of multivalued perturbations.
Indeed, the next proposition, which concerns uniformity of the constants appearing
in Propositions \ref{splem2}, \ref{lem715} and Corollary \ref{cor69},
is also obtained under the weaker notion of $C^1$-convergence.
\end{rem}
\begin{prop}\label{lem627}
Let $\{ \{{\frak s}^n_{\sigma}\mid n \in \Z_{\ge 0}\} \mid \sigma \in \mathscr A\}$
be a uniform family of multivalued perturbations
of $(\widetriangle{\mathcal U},\mathcal K^3)$.
\begin{enumerate}
\item
In Proposition \ref{splem2} the constants $\delta$ and $n_0$
can be taken independent of $\sigma$.
\item
In Proposition \ref{lem715} the set $\frak U(Z)$ and the
constant $n_0$ can be taken independent of $\sigma$.
\item
In Corollary \ref{cor69} the set $\frak U(Z)$ and the
constant $n_0$ can be taken independent of $\sigma$.
Moreover
\begin{equation}
\lim_{n\to \infty}\sup\left\{ d_H\left(X,\left(\bigcup_{\frak p}((\frak s_{\sigma,\frak p}^{n})^{-1}(0)
\cap  \overset{\circ}{\mathcal K^2_{\frak p}})\right)
\cap \frak U(Z)\right) \mid \sigma \in \mathscr A\right\}
= 0.
\end{equation}
\end{enumerate}
\end{prop}
\begin{proof}
This is a consequence of Lemma \ref{lem618}
and the proofs of Propositions \ref{splem2}, \ref{lem715} and Corollary \ref{cor69}.
\end{proof}

\begin{defn}\label{transofdvect}
(Orbifold case)
\begin{enumerate}
\item
Let $\frak s$ be a multisection of a vector bundle $\mathcal E$ on an orbifold $U$.
We say it is {\it transversal to $0$} on $K\subset U$
\index{transversality ! of multisection on orbifold} if for each $x \in K$
and any branch $\frak s_k$ of $\frak s$ at $x$ such that $\frak s_k(x) = 0$, $\frak s_k$ is transversal to $0$.
(Note $\frak s_k : V_x \to E_x$ is a smooth map and $(V_x,\Gamma_x,E_x,\phi_x,\widehat\phi_x)$
is an orbifold chart of $(U,\mathcal E)$.)
\item
In the situation of (1),
let $f : U \to N$ be a smooth map to a manifold.
We say $f$ is {\it strongly submersive} with respect to $\frak s$,
\index{strongly submersive (w.r.t. multisection) ! map on orbifold}
if for any branch $\frak s_k$ of $\frak s$ at $x \in K$ such that $\frak s_k(x) = 0$,
the composition
\begin{equation}\label{form618618}
\frak s_k^{-1}(0) \hookrightarrow V \overset{\!\!\!\!\!\phi_x}{\longrightarrow U}
\overset{\!\! f}\longrightarrow N
\end{equation}
is a submersion.
\item
In the situation of (2), let $g : M \to N$ be a smooth map between manifolds.
Suppose the multisection $\frak s$ is transversal to $0$ on $K$.
We say $(\frak s,f)$ is {\it strongly transversal to $g$}
if $\frak s$ is transversal to $0$ and,
for any branch $\frak s_k$ of $\frak s$ at $x \in K$ such that $\frak s_k(x) = 0$,
the composition
(\ref{form618618})
is transversal to $g$.
\index{strongly transversal (w.r.t. multisection) ! to a map on orbifold}
\end{enumerate}
\end{defn}
\begin{defn}\label{transkurakuravect}
(Good coordinate system case)
\begin{enumerate}
\item
Let $\widetriangle{\frak s}$ be a multisection of $(\widetriangle{\mathcal U},\mathcal K)$,
where $\widetriangle{\mathcal U}$ is
a good coordinate system
of $Z\subseteq X$ and $\mathcal K$ its support system.
We say it is {\it transversal to $0$}
\index{transversality ! of multisection on good coordinate system} if for each
$\frak p \in \frak P$,
$\frak s_{\frak p}$ is transversal to $0$ on $\mathcal K_{\frak p}$.
\item
In the situation of (1),
let $\widetriangle f : (X,Z;\widehat{\mathcal U}) \to N$ be a smooth map to a manifold.
We say $f$ is {\it strongly submersive} with respect to $\widetriangle{\frak s}$,
\index{strongly submersive (w.r.t. multisection) ! map on good coordinate system}
if for each $\frak p \in \frak P$,
$f_{\frak p}$ is
strongly submersive  with respect to $\frak s_{\frak p}$.
\item
In the situation of (2), let $g : M \to N$ be a smooth map between manifolds.
Suppose the multisection $\widetriangle{\frak s}$ is transversal to $0$ on $K$.
We say $(\widetriangle{\frak s},f)$ is {\it strongly transversal to $g$}
\index{strongly transversal (w.r.t. multisection) ! to a map on good coordinate system}
\index{good coordinate system ! strongly transversal (w.r.t. multisection) to a map on good coordinate system}
if for each $\frak p \in \frak P$,
$f_{\frak p}$ is
strongly transversal to $g$  with respect to $\frak s_{\frak p}$.
\item
A multivalued perturbation of a good coordinate system
$\widetriangle{\frak s} = \{\widetriangle{\frak s^n}\}$ is said to be
{\it transversal to $0$} if  $\widetriangle{\frak s^n}$ is transversal
to $0$ for sufficiently large $n$.
\item
Strong submersivity of maps on good coordinate system with
respect to a multivalued perturbation is defined in the same way.
The definition of {\it strong submersivity of a map
$\widetriangle f : (X,Z;\widehat{\mathcal U}) \to N$ on good coordinate system
to $g : M \to N$ with respect to a multivalued perturbation}
is defined in the same way.
\end{enumerate}
\end{defn}
\begin{thm}\label{prop621}
Let ${\widetriangle{\mathcal U}}$ be
a good coordinate system  of $Z \subseteq X$ and $\mathcal K$ its support system.
\begin{enumerate}
\item
There exists a
multivalued perturbation $\widetriangle{\frak s} = \{\frak s^n_{\frak p}\}$
of $({\widetriangle{\mathcal U}},\mathcal K)$ such that
each branch of $\frak s^n_{\frak p}$ is transversal to $0$.
\item
Suppose $\widetriangle f : (X,Z;{\widetriangle{\mathcal U}}) \to N$ is strongly
smooth and is transversal to $g : M \to N$, where $g$ is a map from
a manifold $M$.
Then we may choose $\widetriangle{\frak s}$ such that $\widetriangle f$
is strongly transversal to $g$ with respect to $\widetriangle{\frak s}$.
\end{enumerate}
\end{thm}
This is actually proved during the proof of \cite[Proposition 6.3]{foootech}.
We will prove it in Section \ref{sec:constrsec}.

\subsection{Embedding of Kuranishi structure and multisection}
\label{subset:embmultikura}

\begin{defn}\label{compapertKuranishi}
Let $\widehat{\mathcal U}$ be a Kuranishi structure of $Z \subseteq X$.
A {\it strictly compatible multivalued perturbation
of $\widehat{\mathcal U}$} is a collection
\index{multivalued perturbation ! strictly compatible multivalued perturbation of
Kuranishi structure}
\index{compatibility ! strictly compatible multivalued perturbation of
Kuranishi structure}
\index{multivalued perturbation ! of
Kuranishi structure}
$\widehat{\frak s} =
\{\widehat{\frak s^n}\} = \{ \frak s^{n}_p\}_{p \in Z}$ such that
$\frak s^{n}_p$ is a multisection of $E_{p}$ on $U_p$ for each $p
\in X$
and $n \in \Z_{\ge 0}$, which have the following properties.
\begin{enumerate}
\item
$
\frak s_{p}^{n} \circ \varphi_{pq}
=
\widehat\varphi_{pq}\circ \frak s_{q}^{n}
$
on $U_{pq} $.
\item
$
\lim_{n\to \infty} \frak s_{p}^{n} = s_{p}
$
in $C^1$ sense on $U_{p}$.
\end{enumerate}
The precise meaning of (1), (2) above is the same as in
the case of
Definition \ref{defn612}.
\end{defn}
\begin{rem}
We use the terminology, strictly compatible multivalued perturbations, in
Definition \ref{compapertKuranishi}.
The phrase `strictly compatible' indicates that this is rather a strong condition and is hard to realize.
For example, we may not expect such a perturbation exists for a given Kuranishi structure.
Namely we need to replace the given Kuranishi structure to its appropriate thickening
to obtain strictly compatible multivalued perturbation.
(See Proposition \ref{lemappgcstoKu}.)
Nevertheless we usually omit  the phrase `strictly compatible' and simply say multivalued perturbation.
\end{rem}
\begin{defn}\label{defn692}
\begin{enumerate}
\item
Let $\widehat{\Phi}
= \{\Phi_{p}\}: \widehat{\mathcal U} \to {\widehat{\mathcal U'}}$
be a strict KK-embedding of
Kuranishi structures.
Let $\{\frak s^{n}_{p}\}$ and  $\{\frak s^{\prime n}_{p}\}$
be multivalued perturbations
of  ${\widehat{\mathcal U}}$ and $\widehat{\mathcal U'}$,
respectively.
We say $\{\frak s^{n}_{p}\}$ and  $\{\frak s^{\prime n}_{p}\}$ are
{\it  compatible} with $\widehat{\Phi}$ if
\index{embedding ! compatibility of multivalued perturbation with KK-embedding}
\index{multivalued perturbation ! compatibility with KK-embedding}
\index{compatibility ! of multivalued perturbation with KK-embedding}
$
\frak s_{p}^{\prime n} \circ \varphi_{p}
=
\widehat\varphi_{p}\circ \frak s_{p}^{n}.
$
\item
Let $\widehat{\mathcal U_0}$ be an open substructure of a Kuranishi structure
$\widehat{\mathcal U}$.
Let $\{\frak s^{n}_{p}\}$ be a
multivalued perturbation of  ${\widehat{\mathcal U}}$.
Then $\{\frak s^{n}_{p}\vert_{U_p^0}\}$
is a
multivalued perturbation of  $\widehat{\mathcal U_0}$.
We call it the {\it restriction} of $\{\frak s^{n}_{p}\}$
and write $\{\frak s^{n}_{p}\}\vert_{\widehat{\mathcal U}^0}$.
\index{multisection ! restriction
to open substructure of Kuranishi structure}
\item
Let $\widehat{\Phi}
= \{\Phi_{p}\}: \widehat{\mathcal U} \to {\widehat{\mathcal U'}}$
be a (not necessary strict) KK-embedding of
Kuranishi structures.
Let $\{\frak s^{n}_{p}\}$ and  $\{\frak s^{\prime n}_{p}\}$
be multivalued perturbations of  ${\widehat{\mathcal U}}$ and $\widehat{\mathcal U'}$,
respectively.
We say $\{\frak s^{n}_{p}\}$ and  $\{\frak s^{\prime n}_{p}\}$ are
{\it  compatible} with $\widehat{\Phi}$
if a restriction  $\{\frak s^{n}_{p}\}\vert_{\widehat{\mathcal U_0}}$
is compatible to $\{\frak s^{\prime n}_{p}\}$
with respect to a strict embedding
$\widehat{\mathcal U_0} \to \widehat{\mathcal U'}$.
Here $\widehat{\mathcal U_0}$ is an open substructure of
$\widehat{\mathcal U}$
\end{enumerate}
\end{defn}
\begin{defn}\label{defn69233}
Let ${\widetriangle{\Phi} }
= \{\Phi_{p}\}: {\widetriangle{\mathcal U}} \to {{\widetriangle{\mathcal U'}}}$
be a GG-embedding, and
$\widetriangle{\frak s} = \{\frak s^{n}_{\frak p}\}$ and  $
\widetriangle{\frak s'} = \{\frak s^{\prime n}_{\frak p}\}$
multivalued perturbations
of ${{\widetriangle{\mathcal U}}}$ and ${\widetriangle{\mathcal U'}}$,
respectively.
We say $\{\frak s^{n}_{\frak p}\}$ and  $\{\frak s^{\prime n}_{\frak p}\}$ are
{\it compatible} with ${\widetriangle{\Phi}}$ if
\index{multivalued perturbation ! compatibility with GG-embedding}
\index{embedding ! compatibility of multivalued perturbation
with GG-embedding}
\index{compatibility ! of multivalued perturbation
with GG-embedding}
$
\frak s_{\frak p}^{\prime n} \circ \varphi_{\frak p}
=
\widehat\varphi_{\frak p}\circ \frak s_{\frak p}^{n}.
$
\end{defn}
\begin{defn}\label{defn69}
Let $\widehat{\Phi}
= (\{U_{\frak p}(p)\},\{\Phi_{p\frak p}\}): {\widetriangle{\mathcal U}} \to \widehat{\mathcal U^+}$
be a GK-embedding.
Let $\widetriangle{\frak s} = \{\frak s^{n}_{\frak p}\}$ and  $
\widehat{\frak s} = \{\frak s^{n}_{p}\}$
be multivalued perturbations of  ${\widetriangle{\mathcal U}}$ and $\widehat{\mathcal U^+}$,
respectively.
We say $\{\frak s^{n}_{\frak p}\}$ and  $\{\frak s^{n}_{p}\}$ are
{\it  compatible} with $\widehat{\Phi}$ if
\index{compatibility ! of multivalued perturbations
with GK-embedding}
\index{multivalued perturbation ! compatibility with GK-embedding}
\index{embedding ! compatibility of multivalued perturbations with GK-embedding}
$
\frak s_{p}^{n} \circ \varphi_{p\frak p}
=
\widehat\varphi_{p\frak p}\circ \frak s_{\frak p}^{n}
$
holds on $U_{\frak p}(p)$.
\end{defn}

There are various obvious statements about the composition of embeddings
and its compatibilities with the multivalued perturbations.
We leave to the interested readers to state and prove them.

\begin{defn}\label{situ610}
Let ${\widetriangle{\mathcal U}}$, ${\widetriangle{\mathcal U_0}}$
be good coordinate systems of $Z \subseteq X$.
An open GG-embedding $\widetriangle\Phi : {\widetriangle{\mathcal U_0}}
\to {\widetriangle{\mathcal U}}$
is said to be {\it relatively compact}
\index{embedding ! relatively compact open embedding of good coordinate systems}
if, for each $\frak p$, the subset $\varphi_{\frak p}(U^0_{\frak p})$
is relatively compact in $U_{\frak p}$.
\end{defn}
\begin{prop}\label{lemappgcstoKu}
Let ${\widetriangle{\mathcal U_0}}
\to {\widetriangle{\mathcal U}}$
be a relatively compact open GG-embedding of
good coordinate systems of
$Z \subseteq X$.
Then there exist a Kuranishi structure $\widehat{\mathcal U}$
and a GK-embedding
${\widetriangle{\mathcal U_0}} \to \widehat{\mathcal U}$
with the following properties.
\begin{enumerate}
\item
Let $\mathcal K$ be a support system of
${\widetriangle{\mathcal U}}$ and
${\widetriangle{\frak s}} = \{\frak s^{n}_{\frak p}\}$  a multivalued perturbation of
$({\widetriangle{\mathcal U}},\mathcal K)$.
We assume
$
\overline{\varphi_{\frak p}(U^0_{\frak p})} \subset
{\rm Int} \,\, \mathcal K_{\frak p}
$.
\par
Then there exists a multivalued perturbation $\widehat{\frak s} = \{\frak s^{n}_p\}$
of ${\widehat{\mathcal U}}$ such that
$\widetriangle{\frak s}\vert_{{\widetriangle{\mathcal U}}_0}$ and ${\widehat{\frak s}}$
are compatible with the embedding
${\widetriangle{\mathcal U_0}}
\to {\widehat{\mathcal U}}$.
\item
If ${\widetriangle f} : (X,Z;{\widetriangle{\mathcal U}}) \to Y$
is a strongly continuous (resp. strongly smooth) map,
then there exists ${\widehat f} : (X,Z;{\widehat{\mathcal U}}) \to Y$
such that ${\widetriangle f}\vert_{{\widetriangle{\mathcal U_0}}}$
is a pullback of ${\widehat f}$ by the embedding
${\widetriangle{\mathcal U_0}}
\to {\widehat{\mathcal U}}$.
\end{enumerate}
\end{prop}
We will prove Proposition \ref{lemappgcstoKu}
together with the following relative version.
\begin{prop}\label{prop634}
In the situation of Proposition \ref{lemappgcstoKu} (1)
we assume the following in addition.
\begin{enumerate}
\item[(\EightStarTaper)]
$\widehat{\mathcal U^+}$ is a Kuranishi structure of $Z
\subseteq X$ and
${\widetriangle{\mathcal U}} \to \widehat{\mathcal U^+}$
is a GK-embedding.
\end{enumerate}
Then we may choose
${\widehat{\mathcal U}}$ and ${\widetriangle{\mathcal U_0}}
\to {\widehat{\mathcal U}}$
in Proposition \ref{lemappgcstoKu} such that
the following holds in addition.
\par
There exists a KK-embedding
${\widehat{\mathcal U}} \to \widehat{\mathcal U^+}$
such that:
\begin{enumerate}
\item[(i)]
The composition of
${\widetriangle{\mathcal U_0}}
\to {\widehat{\mathcal U}}$ and
${\widehat{\mathcal U}} \to \widehat{\mathcal U^+}$
is equivalent to the composition of the embeddings
${\widetriangle{\mathcal U_0}} \to {\widetriangle{\mathcal U}}$,
${\widetriangle{\mathcal U}} \to \widehat{\mathcal U^+}$.
\item[(ii)]
$\widehat{\mathcal U^+}$ is a thickening of ${\widehat{\mathcal U}}$.
\item[(iii)]
In the situation of Proposition \ref{lemappgcstoKu} (2), we assume that
$\widetriangle f : (X,Z; {\widetriangle{\mathcal U}}) \to Y$ is a pullback of
$\widehat{f^+} : (X,Z; {\widehat{\mathcal U^+}}) \to Y$, in addition.
Then we may take $\widehat f$ such that it is a pullback of $\widehat{f^+}$.
\end{enumerate}
\end{prop}
\begin{equation}
\xymatrix{
&  &  & Y \\
\widetriangle {\mathcal U} \ar[urrr]^{\widetriangle f}
\ar[r]  & \widehat{\mathcal U^+} \ar[urr]_{\widehat f^+} && \\
\widetriangle {\mathcal U^0}  \ar@{^{(}-{>}}[u]\ar@{.>}[r] &\widehat {\mathcal U}
\ar@{.>}[u]\ar@{.>}[uurr]_{\widehat f} &&
}
\end{equation}
\begin{proof}
[Proof of Propositions \ref{lemappgcstoKu}
and  \ref{prop634}]
We put $\mathcal K^0_{\frak p} = \overline{\varphi_{\frak p}(U^0_{\frak p})}$
and $\mathcal K^0 = \{\mathcal K^0_{\frak p}\}$.
Then $(\mathcal K^0,\mathcal K)$ is a support pair.
Let $p \in Z$. We take $\frak p_p$ with the following properties.
\begin{proper}\label{property619}
\begin{enumerate}
\item
$p \in \psi_{\frak p_p}(s_{\frak p_p}^{-1}(0) \cap \mathcal K^0_{\frak p_p})$.
\item
If $p \in \psi_{\frak q}(s_{\frak q}^{-1}(0) \cap \mathcal K^0_{\frak q})$ then $\frak q \le \frak p_p$.
\end{enumerate}
\end{proper}
Existence of such $\frak p_p$ follows from Definition \ref{gcsystem} (6).
We take an open neighborhood $U_p$ of $p$ in $\mathcal K_{\frak p_p}$
with the following properties.
\begin{proper}\label{property620}
\begin{enumerate}
\item
$U_p$ is relatively compact in ${\rm Int}\,\mathcal K_{\frak p_p}$.
\item If $\psi_{p}(\overline{U_p}\cap s_p^{-1}(0)) \cap
\psi_{\frak q}(\mathcal K^0_{\frak q} \cap s_{\frak q}^{-1}(0)) \cap Z \ne \emptyset$,
then $\frak q \le \frak p_p$.
\item
In the situation of Proposition \ref{prop634}
we require $U_p \subset U_{\frak p_p}(p)$ in addition.
Here $U_{\frak p_p}(p)$
is an open neighborhood of $o_p$ in $U_{\frak p}$ that
appears in the definition of the GK-embedding ${\widetriangle{\mathcal U}}
\to \widehat{\mathcal U^+}$. (Definition \ref{defn69}.)
\end{enumerate}
\end{proper}
We define a Kuranishi chart $\mathcal U_p$ by
$\mathcal U_{\frak p_p}\vert_{U_p}$.
\par
We next define a coordinate change among them.
Let $q \in \psi_{\frak p_p}(s_{\frak p_p}^{-1}(0) \cap U_p)$.
Property \ref{property619}
and Property \ref{property620} (2) imply that $\frak p_q \le \frak p_p$.
Therefore there exists a coordinate change $\Phi_{\frak p_p\frak p_q}$
of the good coordinate system ${\widetriangle{\mathcal U}}$.
The coordinate change from $\mathcal U_q$ to $\mathcal U_p$
is by definition the restriction of $\Phi_{\frak p_p\frak p_q}$
to $U_q \cap \varphi_{\frak p_p\frak p_q}^{-1}(U_q)$.
Compatibility of the coordinate changes follows from the
compatibility of the coordinate changes of ${\widetriangle{\mathcal U}}$
and the commutativity of Diagram (\ref{diag33}).
We thus obtain the required Kuranishi structure $\widehat{\mathcal U}$.
\par
Firstly we prove Proposition \ref{lemappgcstoKu}  (1).
We define the GK-embedding
${\widetriangle{\mathcal U_0}}
\to \widehat{\mathcal U}$.
Let $p \in U_{\frak p}^0 \cap X$.
Since $\mathcal K_{\frak p}^0$ is the closure of $U_{\frak p}^0$,
Property \ref{property619} (2) implies $\frak p \le \frak p_p$.
Therefore there exists a coordinate change
$\Phi_{\frak p_p\frak p} : \mathcal U_{\frak p} \to \mathcal U_{\frak p_p}$.
We put
$$
U^0_{\frak p}(p) = \varphi_{\frak  p_p\frak p}^{-1}(U_p)\cap U_{\frak p}^0,
\qquad
\Phi^0_{p\frak p} = \Phi_{\frak  p_p\frak p}\vert_{U^0_{\frak p}(p)}.
$$
It is easy to see that they define the required GK-embedding ${\widetriangle{\mathcal U_0}}
\to \widehat{\mathcal U}$.
\par
Secondly we prove Proposition \ref{lemappgcstoKu}  (2) (3).
We define $\frak s^{n}_p = \frak s^{n}_{\frak p_p}\vert_{U_p}$.
Its compatibility with coordinate change follows from one
of $\frak s^{n}_{\frak p}$.
Thus we obtain a multivalued perturbation $\{\frak s_{p}^{n}\}$.
The strong compatibility of it with the GK-embedding
$ {\widetriangle{\mathcal U_0}}
\to \widehat{\mathcal U}$ follows from the strong
compatibility of $\frak s^{n}_{\frak p}$ with coordinate change.
\par
If ${\widetriangle f} = \{f_{\frak p}\}$, then we define
$
f_{p} = f_{\frak p_p}\vert_{U_p}$.
It is easy to see that it has required properties.
\par
Thirdly we prove Proposition \ref{prop634} (i).
We define a KK-embedding $: \widehat{\mathcal U} \to \widehat{\mathcal U^+}$.
Let $p \in X$.
We consider $\varphi_{p\frak p_p} : \mathcal U_{\frak p_p}\vert_{U_{\frak p}(p)} \to \mathcal U^+_{p}$ that is a part of the data defining the
GK-embedding ${\widetriangle{\mathcal U}}
\to \widehat{\mathcal U^+}$ (Definition \ref{defn69}.)
By Property \ref{property620} (3) we can
restrict it to $U_p$. It is easy to see that they are compatible with coordinate changes and
define the required KK-embedding $: \widehat{\mathcal U} \to \widehat{\mathcal U^+}$.
\par
Commutativity of Diagram \ref{diagram58} implies that the composition
${\widetriangle{\mathcal U_0}}
\to \widehat{\mathcal U} \to \widehat{\mathcal U}^{+}$
is the given embedding
${\widetriangle{\mathcal U_0}} \to \widehat{\mathcal U}^{+}$.
\par
We finally prove
Proposition \ref{prop634} (ii)(iii), that is, $\widehat{\mathcal U^+}$
is a thickening of $\widehat{\mathcal U}$.
Let $p \in X$. We put
$$
O_p
=
\psi_{\frak p_p}
(
s_{\frak p_p}^{-1}(0) \cap U_{p} \cap {\rm Int}\, \mathcal K_{\frak p_p}
).
$$
Note $U_{p} \subset U_{\frak p_p}(p) \subset U_{\frak p_p}$.
\par
Let $q \in O_p$. Then
$$
q \in \psi_{\frak p_p}(s_{\frak p_p}^{-1}(0) \cap
{\rm Int}\, \mathcal K_{\frak p_p})
\cap \psi_p^+((s_{p}^+)^{-1}(0)).
$$
Therefore there exist $O_{\frak p_p}(q) \subset U_{\frak p_p}$
and $\varphi_{q \frak p_p} : O_{\frak p_p}(q)
\to U_q^+$.
We put
$$
W_{p}(q) = U_p \cap \varphi^{-1}_{q \frak p_p}(U_{pq}^+).
$$
Then
$$
\varphi_p(W_p(q))
=
\varphi_{p \frak p_p}(W_p(q))
\subset
\varphi^+_{pq}(\varphi_{q \frak p_p}(W_p(q)))
\subset
\varphi^+_{pq}(U_{pq}^+).
$$
We have thus checked Definition \ref{thickening} (2)(a).
Definition \ref{thickening} (2)(b) can be checked
in the same way by using
$\widehat\varphi_{q \frak p_p}$.
We have thus proved (ii).
\par
(iii) is a consequence of the fact that $U^+_{p}$ is an open subset of
$U_{\frak p_p}$ and $f^+_{p}$ is a restriction of $f_{\frak p_p}$.
\par
The proof of Propositions \ref{lemappgcstoKu}
and  \ref{prop634} is now complete.
\end{proof}
Propositions \ref{lemappgcstoKu}
and  \ref{prop634} provide a way to transfer  a
multivalued perturbation of a good coordinate system to that
of a Kuranishi structure.
The next results describe the way of transferring them in the opposite direction.
For this, we need one more definition.

\begin{defn}\label{defn6928}
Let $\widehat{\Phi}
= \{\Phi_{\frak p p}\} : \widehat{\mathcal U} \to {\widetriangle{\mathcal U}}$
be a KG-embedding.
Let $\widehat{\frak s} = \{\frak s^{n}_{p}\}$ and $
\widetriangle{\frak s} = \{\frak s^{n}_{\frak p}\}$
be multivalued perturbations of  ${\widehat{\mathcal U}}$ and ${\widetriangle{\mathcal U}}$,
respectively.
We say $\{\frak s^{n}_{p}\}$ and $\{\frak s^{n}_{\frak p}\}$  are
{\it compatible} with $\widehat{\Phi}$ if they satisfy
$$
  \varphi_{\frak p p} \circ \frak s_{p}^{n}
=
\frak s_{\frak p}^{n}\circ \widehat\varphi_{\frak p p}
$$
on $U_{p}$.
\index{compatibility ! between good coordinate system and Kuranishi structure}
\index{good coordinate system ! compatibility between good coordinate system and Kuranishi structure}
\index{Kuranishi structure ! compatibility between good coordinate system and Kuranishi structure}
\end{defn}
\begin{prop}\label{le614}
Let $\widehat{\mathcal U}$ be a Kuranishi structure on $Z \subseteq X$
and $\widehat{\frak s} = \{\frak s^{n}_{p}\}$
a multivalued perturbation of ${\widehat{\mathcal U}}$.
Then we can take a good coordinate system $\widetriangle{\mathcal U}$
and the strict KG-embedding $\widehat{\Phi} : \widehat{\mathcal U_0} \to \widetriangle{\mathcal U}$
in Theorem \ref{Them71restate} so that the following holds in addition.
\begin{enumerate}
\item
There exists a multivalued perturbation ${\widetriangle{\frak s}} = \{\frak s^{n}_{\frak p}\}$ of ${\widetriangle{\mathcal U}}$
such that
${\widehat{\frak s}}\vert_{\widehat{\mathcal U_0}}$ and
${\widetriangle{\frak s}}$
are compatible with the embedding $\widehat{\Phi}$.
\item
If ${\widehat f} : (X,Z;{\widehat{\mathcal U}})\to Y$
is a strongly continuous map, then there
exists ${\widetriangle f} : (X,Z;{\widetriangle{\mathcal U}})
\to Y$ such that  ${\widetriangle f}\circ \widehat\Phi$
is a pullback of ${\widehat f}$.
If ${\widehat f}$ is strongly smooth
(resp. weakly submersive) then so is ${\widetriangle f}$.
The transversality to $M \to Y$ is also preserved.
\end{enumerate}
\end{prop}
\begin{prop}\label{pro616}
Suppose we are in the situation of Propositions \ref{prop518}
(resp. \ref{prop519})
and \ref{le614}.
Then we can take the GK-embedding
$\widehat{\Phi^+} :  \widetriangle{\mathcal U} \to
\widehat{\mathcal U^+}$ in Proposition \ref{prop518}
(resp. the GK-embeddings $\widehat{\Phi^+_a} :  \widetriangle{\mathcal U} \to
\widehat{\mathcal U^+_a}$ in Proposition \ref{prop519}
($a=1,2$))
so that the following holds.
\begin{enumerate}
\item
If ${\widehat{\frak s^+}}$ is a multivalued perturbation of ${\widehat{\mathcal U^+}}$
such that ${\widehat{\frak s^+}}$, ${\widehat{\frak s}}$ are
strongly compatible
with the KK-embedding ${\widehat{\mathcal U}} \to \widehat{\mathcal U^+}$,
then we may choose $\widetriangle{\frak s}$
such that  $\widetriangle{\frak s}$, $\widehat{\frak s^+}$ are
compatible
with the embedding $\widehat{\Phi^+}$.
(resp.
If  ${\widehat{\frak s^+_a}}$ ($a=1,2$) is a multivalued perturbation of ${\widehat{\mathcal U^+_a}}$
such that ${\widehat{\frak s^+_a}}$, ${\widehat{\frak s}}$ are
strongly compatible
with the embedding ${\widehat{\mathcal U}}
\to \widehat{\mathcal U^+}$,
then we may choose $\widetriangle{\frak s}$
such that  $\widetriangle{\frak s}$, $\widehat{\frak s^+_a}$ are both
compatible
with the embedding $\widehat{\Phi^+_a}$.)
\item
If ${\widehat{f^+}} : (X,Z;\widehat{\mathcal U^+}) \to Y$ is a strongly
continuous map
so that the pull back of ${\widehat{f^+}}$ is
${\widehat f} : (X,Z;{\widehat{\mathcal U}})\to Y$,
then we may choose $\widetriangle{f} : (X,Z;\widetriangle{\mathcal U}) \to Y$
such that  ${\widehat{f^+}} \circ \widehat{\Phi^+} = \widetriangle{f}$.
(resp.
If ${\widehat{f^+_a}} : (X,Z;\widehat{\mathcal U_a^+}) \to Y$ are strongly
continuous maps
such that the pull back of both  ${\widehat{f^+_a}}$ ($a=1,2$) are
${\widehat f} : (X,Z;{\widehat{\mathcal U}})\to Y$,
then we may choose $\widetriangle{f} : (X,Z;\widetriangle{\mathcal U}) \to Y$
such that  ${\widehat{f^+_a}} \circ \widehat{\Phi^+_a} = \widetriangle{f}$.)
\end{enumerate}
\end{prop}
We will prove Propositions \ref{le614} and \ref{pro616} in
Subsection \ref{subsec:movingmulsectionetc}.

\subsection{General strategy of construction of virtual fundamental chain}
\label{bigremarkinsec6}
In this subsection we summarize a general strategy we will take
and
show how the results of this section will be used in the strategy.
\par
Let us start with a Kuranishi structure $\widehat{\mathcal U}$ of $X$.
\par\medskip
\begin{enumerate}
\item[{\bf Step 1}:]
We find a good coordinate system ${\widetriangle{\mathcal U}}$ such that
$\widehat{\mathcal U} < {\widetriangle{\mathcal U}}$,
which means that $\widehat{\mathcal U}$ is embedded in
${\widetriangle{\mathcal U}}$.
(Theorem \ref{Them71restate}.)
\item[{\bf Step 2}:]
We find a multivalued perturbation ${\widetriangle{\frak s}}$ of ${\widetriangle{\mathcal U}}$
that has various transversality properties.
(Theorem \ref{prop621}.)
\item [{\bf Step 3}:]
We obtain a virtual fundamental chain
associated to the perturbations of
${\widetriangle{\mathcal U}}$.
\item[{\bf Step 4}:]
We next apply Proposition \ref{lemappgcstoKu}
to obtain Kuranishi structure $\widehat{\mathcal U^+}$
such that
$$
\widehat{\mathcal U} < {\widetriangle{\mathcal U}}< {\widehat{\mathcal U^+}}
$$
and a multivalued perturbation ${\widehat{\frak s^+}}$ of it.
\item[ {\bf Step 5}:]
We apply Proposition \ref{prop518}
to $\widehat{\mathcal U^{+}}$
and obtain a good coordinate system
${\widetriangle{\mathcal U^{+}}}$ such that
$$
\widehat{\mathcal U} < {\widetriangle{\mathcal U}}< {\widehat{\mathcal U^+}} <{\widetriangle{\mathcal U^{+}}}.
$$
Moreover, the multivalued perturbation ${\widehat{\frak s^+}}$
induces a multivalued perturbation ${\widetriangle{\frak s^+}}$ of ${\widetriangle{\mathcal U^{+}}}$.
\item[{\bf Step 6}:]
The transversality of ${\widetriangle{\frak s}}$ implies one of ${\widehat{\frak s^+}}$
and then one of ${\widetriangle{\frak s^+}}$.
\item[{\bf Step 7}:]
We obtain a virtual fundamental chain
associated to ${\widetriangle{\frak s^+}}$.
\end{enumerate}
\par\smallskip
Now an important statement
is that the virtual fundamental chain
obtained in Step 3 is the same as the
virtual fundamental {\it chain} obtained
in Step 7 for any sufficiently large $n$.
(Roughly speaking, this is a consequence of
Proposition \ref{lem715}.)
Its de Rham version is Proposition \ref{integralinvembprop}.
\par
This statement can be used as follows.
Note the constructions in Step 5 is not unique.
Namely for each given ${\widehat{\mathcal U^+}}$
there are many possible
choices of ${\widetriangle{\mathcal U^{+}}}$.
However the virtual fundamental chain associated to it
is independent of such choices.
Moreover it coincides with the virtual fundamental chain
obtained in Step 3.
\par
In other words, we can recover the virtual fundamental chain of
${\widetriangle{\frak s}}$ (that is defined on
${\widetriangle{\mathcal U}}$ in Step 3)
from ${\widehat{\frak s^+}}$ that is defined
on ${\widehat{\mathcal U^+}}$.
\par
By this reason, we can forget the good coordinate system ${\widetriangle{\mathcal U}}$
and remember only the Kuranishi structure ${\widehat{\mathcal U^+}}$
and ${\widehat{\frak s^+}}$ on it.
Since Kuranishi structure behaves better with fiber product than
good coordinate system, we can use this fact to make the whole construction
compatible with the fiber product description of the boundaries.
\par
In Sections \ref{sec:contfamily}-\ref{sec:contfamilyconstr}, we will work out this process in the
de Rham model in great detail,
where we will use {\it CF-perturbations}, which is an abbreviation of {\it continuous family perturbations}, instead of multivalued perturbations.
\par
In Part 2 of this document, we will discuss in detail the way how we make the whole construction compatible when we start with the
system of Kuranishi structures.

\section{CF-perturbation
and integration along the fiber (pushout)}
\label{sec:contfamily}

\subsection{Introduction to Sections \ref{sec:contfamily}-\ref{sec:composition}}
\label{subsec:seccontintro}

As we mentioned in Introduction, we study systems of
Kuranishi structures so that the boundary of each of its
member is described by the fiber product of the other members.
We will obtain an algebraic structure on certain chain complexes
which realize the homology groups of certain spaces.
They are the spaces over which we take fiber product
between members of the system of the Kuranishi structures.
To work out this process we need to make a choice
of the homology theory we use.
The choices are de Rham cohomology, singular homology,
\v Cech cohomology,
Morse homology or Kuranishi homology (see \cite{joyce}), and etc..
In \cite{fooobook} we took the most standard choice, that is, the singular
homology.
In this article, we mainly use  de Rham cohomology.
There are three advantages in using de Rham cohomology.
One is that it seems shortest to write a detailed and rigorous proof.
The second is that it is easiest to keep as much symmetry as possible.
The third is, by using de Rham cohomology,
we might clarify some direct relation to quantum field
theory, (especially in the case of perturbation of the constant maps).
There are certain disadvantage in using de Rham cohomology.
The most serious disadvantage is that we can work only
over real or complex numbers as a ground field.
Certain technical points which appear when we use singular homology
will be explained elsewhere.
The way to use Morse homology is discussed in \cite{FOOO08III}.
Actually Morse homology is one the authors of the present
document had used around 20 years ago in
\cite{Fuk97III},\cite{Oh96II}, etc. See \cite[Remark 1.32]{fooobook}.
\par
The situation we work with is as follows.

\begin{shitu}\label{smoothcorr}(See \cite[Section 12]{fooo09})
Let $X$ be a compact metrizable space,
and
$\widehat{\mathcal U}$ a Kuranishi structure of $X$ (with or without boundaries or corners).
Let $M_s$ and $M_t$ be $C^{\infty}$ manifolds.
We assume $\widehat{\mathcal U}$, $M_s$ and $M_t$
are oriented\footnote{In certain situations, for example in \cite[Subsection 8.8]{fooobook2},
we discussed slightly more general case. Namely we discussed
the case when $\widehat{\mathcal U}$, $M_s$ and $M_t$ are not
necessarily orientable by introducing appropriate
$\Z_2$ local systems. See \cite[Section A2]{fooobook2}.}.
\par
Let $\widehat f_s : (X;\widehat{\mathcal U}) \to M_s$ be a strongly smooth map and
$\widehat  f_t : (X;\widehat{\mathcal U}) \to M_t$ a weakly submersive strongly smooth map.
We call $\frak X = ((X;\widehat{\mathcal U});\widehat f_s,\widehat f_t)$ a
{\it smooth correspondence}
\index{smooth correspondence ! smooth correspondence} from $M_s$ to $M_t$.$\blacksquare$
\end{shitu}
Our goal in Sections \ref{sec:contfamily}-\ref{sec:composition} is to associate a linear map
\begin{equation}\label{smoothcorr}
{\rm Corr}_{\frak X} : \Omega^k(M_s) \to \Omega^{k+ \ell}(M_t)
\end{equation}
to a smooth correspondence $\frak X$ and study its properties.
Here $\Omega^k(M_s)$ is the set of smooth $k$ forms on $M_s$ and
$$
\ell = \dim M_t  - \dim (X,\widehat{\mathcal U}).
$$
\par
If $(X;\widehat{\mathcal U})$ is a smooth orbifold, namely if
the obstruction bundles are all $0$, then we can define (\ref{smoothcorr}) by
$$
{\rm Corr}_{\frak X}(h) = f_t! (f_s^*(h)).
$$
Here $f^*_s : \Omega^{k}(M_s) \to \Omega^{k}(X)$ is the pullback of the differential form and
$f_t! : \Omega^{k}(X) \to \Omega^{k+\ell}(M_t)$ is {\it the integration along the fiber},
or {\it pushout}, which is characterized by
$$
\int_{M_t} f_t!(v)\wedge \rho  = \int_X v \wedge f_t^*\rho.
$$
(Note $\ell \le 0$ in this case.)
Existence of such $f_t!$ is a consequence of the fact that $f_t$ is a proper submersion.
(In our situation where $(X;\widehat{\mathcal U})$ is an orbifold this is a consequence of the weak submersivity of $\widehat f_t$.)
\par
When the obstruction bundle is nontrivial, we need to perturb the
space $X$ so that integration along the fiber is well-defined.
However, taking a {\it multivalued perturbation} of $\widehat{\mathcal U}$ discussed in
Section \ref{sec:multisection}
is not
good enough for our purpose unless $M_t$ is a point.
Let us elaborate on this point below.
\par
Suppose $\widetriangle{\frak s} = \{\frak s_{\frak p}^{n}\}$ is a multivalued perturbation
of ${\widetriangle{\mathcal U}}$ where ${\widetriangle{\mathcal U}}$ is a good
coordinate system compatible with ${\widehat{\mathcal U}}$.
If we assume that $\frak s_{\frak p}^{n}$ is transversal to
$0$, then in case
$\dim{\widehat{\mathcal U}} = \deg h$, we can define the number
$$
\int_{\bigcup_{\frak p}(\frak s_{\frak p}^{n})^{-1}(0)} f_s^*(h) \in \R.
$$
However we can not expect the map
$$
f_t\vert_{\bigcup_{\frak p}(\frak s_{\frak p}^{n})^{-1}(0)} :
\bigcup_{\frak p}(\frak s_{\frak p}^{n})^{-1}(0) \to M_t
$$
is a `submersion' in any reasonable sense.
In fact, there may happen the case when
$\dim{\widehat{\mathcal U}}$ is strictly smaller than $\dim M_t$.
Therefore the integration along the fiber
$f_t\vert_{\bigcup_{\frak p}(\frak s_{\frak p}^{\epsilon})^{-1}(0)}!$
sends a differential form to a distributional form which may not be
smooth.
Thus we need to find an appropriate way to smooth it
to define (\ref{smoothcorr}).
The way we take here is to use {\it CF-perturbation}, which is an abbreviation of {\it continuous family perturbations}.
We had discussed this construction in \cite[Section 7.5]{fooobook2}, \cite[Section 12]{fooo010},
\cite[Section 4]{fooo091}, \cite[Section 12]{fooo09}.

Now the outline of Sections \ref{sec:contfamily}-\ref{sec:composition} is as follows.
We review and describe the CF-perturbation and the
integration along the fiber in greater detail and
then combine them with the process to transfer
various objects from a Kuranishi structure
to a good coordinate system and back.
More specifically, we first introduce the notion of a CF-perturbation of
a single Kuranishi chart $\mathcal U$ in Subsection \ref{subsec:conper1chart}, where we find in Proposition \ref{prop721} that the set of CF-perturbations of $\mathcal U$ turns out to be
a sheaf $\mathscr S$.
We introduce several subsheaves of $\mathscr S$ which satisfy various  transversality conditions.
Using these subsheaves we define the pushout of differential forms.
Next in
Subsection \ref{subsec:conpergcs},
we generalize these results to the case of
CF-perturbations of a good coordinate system.
Then we can formulate the pushout of differential forms and smooth correspondences in a good coordinate system.
We also formulate and prove Stokes' formula for good coordinate system in Section
\ref{sec:stokes}.
So far, everything here are discussed based on good coordinate system.
However, as mentioned at the end of Section \ref{sec:fiber},
it is more convenient and natural to use Kuranishi structure itself rather than
good coordinate system, when we study
the fiber product of $K$-spaces.
For this purpose, we start from a Kuranishi structure and
use certain embeddings into a good coordinate system and/or
another Kuranishi structure
introduced in Section \ref{sec:thick} to
translate the above results based on the good coordinate system into ones based on the Kuranishi structure and study their relationship.
As a result, we show in Theorem
\ref{theorem915} that the pushout of differential forms for Kuranishi structure
is indeed independent of choice of good coordinate system.
After these foundational results on the pushout of differential forms
are prepared,
we prove a basic result about smooth correspondence, which is called {\it composition formula of smooth correspondence}, in Theorem \ref{compformulaprof}.
The proof of the existence of a CF-perturbation for any $K$-space
is postponed till Section \ref{sec:contfamilyconstr}, where we also prove
 in Proposition \ref{prop123123} that the sheaf $\mathscr S$ of CF-perturbations, together with the several subsheaves mentioned above, is soft.

\subsection{CF-perturbation on a single
Kuranishi chart}
\label{subsec:conper1chart}

We first consider the
situation where we have only one Kuranishi chart as follows.
After that we will introduce a CF-perturbation of good coordinate system
in Subsection \ref{subsec:conpergcs}.
A CF-perturbation of Kuranishi structure will be defined in
Subsection \ref{subsec:contfamiKura}.

\begin{shitu}\label{smoothcorrsingle}
Let $\mathcal U = (U,\mathcal E,s,\psi)$ be a Kuranishi chart of $X$, and
$f : U \to M$ a smooth submersion to a smooth manifold $M$, and
$h$ a differential form on $U$ which has compact support. Assume that
$U, \mathcal E$ and $M$ are oriented.$\blacksquare$
\end{shitu}
\subsubsection{CF-perturbation on one
orbifold chart}
\label{cftoneorbifoldchart}
Under Situation \ref{smoothcorrsingle}
let
$\frak V_x = (V_x,\Gamma_x,E_x,\phi_x,\widehat\phi_x)$
be an orbifold chart of $(U,\mathcal E)$.
(Definition \ref{defn2613}.)
We assume $(V_x,\Gamma_x,\phi_x)$ is an oriented orbifold chart.
(Definition \ref{defn281010} (5).)
Since $f$ is a submersion,
the composition
$
f\circ \phi_x : V_x  \to U \to M
$
is a smooth submersion, which is denoted by $f_x$.

\begin{defn}\label{defn73ss}
A {\it CF-perturbation (=continuous family perturbation)}
\index{CF-perturbation ! CF-perturbation} of $\mathcal U$ on our orbifold chart
$\frak V_x = (V_x,\Gamma_x,E_x,\phi_x,\widehat\phi_x)$ consists of
$\mathcal S_x = (W_x,\omega_x,\{{\frak s}_x^{\epsilon}\})$,
$0 < \epsilon  \le 1$,  with the following properties:
\begin{enumerate}
\item
$W_x$ is an open neighborhood of $0$ of a finite dimensional vector space
$\widehat W_x$ on which $\Gamma_x$ acts linearly.
$W_x$ is $\Gamma_x$ invariant.
\item
$
{\frak s}_x^{\epsilon} : V_x \times W_x \to E_x
$
is a $\Gamma_x$-equivariant smooth map for each $0 < \epsilon  \le 1$.
\item
For $y \in V_x$ , $\xi \in W_x$ we have
\begin{equation}\label{C0convconti}
\lim_{\epsilon\to 0} {\frak s}_x^{\epsilon}(y,\xi)
= s_x(y)
\end{equation}
in compact $C^1$-topology
on $V_x \times W_x$.
\item
$\omega_x$ is a smooth differential form on $W_x$ of
degree $\dim W_x$ that is
$\Gamma_x$ invariant, of compact support and
$$
\int_{W_x} \omega_x = 1.
$$
We assume $\omega_x = \vert\omega_x\vert{\rm vol}_x$ here ${\rm vol}_x$ is a volume form of the oriented manifold $W_x$
and $\vert\omega_x\vert$ is a non-negative function.
\end{enumerate}
For each $0 < \epsilon \le 1$, we denote the
restriction of $\mathcal S_x$ at $\epsilon$, by $\mathcal S_x^{\epsilon} = (W_x,\omega_x,{\frak s}_x^{\epsilon})$.
\end{defn}
\begin{rem}
\begin{enumerate}
\item
In Definition \ref{defn73ss} (3) we regard $s_x : V_x \to  E_x$
as a $\Gamma_{x}$ equivariant map that is a local representative of
the Kuranishi map in the sense of Definition \ref{defnlocex}.
\item
In our earlier writings, we used a family of {\it multi}-sections parameterized by $W_x$.
Here we use a family of sections parameterized by $W_x$ on $V_x$ such that
it is $\Gamma_x$ equivariant as a map from $V_x \times W_x$.
We also allow $W_x$ to have a nontrivial $\Gamma_x$ action.
This formulation seems simpler.
\item
We may regard $ {\frak s}_x^{\epsilon}$ as a local representative
of a section of the vector bundle
$(V_x \times W_x \times E_x)/\Gamma_x \to (V_x \times W_x)/\Gamma_x$.
(Lemma \ref{lem2627}.)
\end{enumerate}
\end{rem}
\begin{defn}\label{conmultiequiv11}
Let $\mathcal S^i_x = (W^i_x,\omega^i_x,\{{\frak s}_{x}^{\epsilon,i}\})$
$(i=1,2)$ be two CF-perturbations
of $\mathcal U$ on $\frak V_x$.
\begin{enumerate}
\item
We say
$\mathcal S^1_x$ is a {\it projection} of
\index{CF-perturbation ! projection of}
$\mathcal S^2_x$,
if there exist a map $\Pi :
\widehat W^2_x \to \widehat W^1_x$
with the following properties.
\begin{enumerate}
\item
$\Pi$ is a $\Gamma_x$ equivalent linear projection
which sends $W^2_x$ to $W^1_x$ and satisfies
$
\Pi!(\omega^2_x) = \omega^1_x.
$
\item
For each $y \in V_x$ and $\xi \in W_x^1$
we have
$$
{\frak s}_{x}^{\epsilon,1}(y,\Pi(\xi)) = {\frak s}_{x}^{\epsilon,2}(y,\xi).
$$
\end{enumerate}
\item
We say $\mathcal S^1_x$ is {\it equivalent} to $\mathcal S^2_x$
\index{CF-perturbation ! equivalence on an orbifold chart}
on $\frak V_x$ if there exist $N$ and $\mathcal S^{(i)}_x$
for $i=0,\dots,2N$ with the following properties.
\begin{enumerate}
\item
$\mathcal S^{(i)}_x$ is a CF-perturbation
of $\mathcal U$ on $\frak V_x$.
\item
$\mathcal S^{(0)}_x = \mathcal S^1_x$,
$\mathcal S^{(2N)}_x = \mathcal S^2_x$.
\item
$\mathcal S^{(2k-1)}_x$ and $\mathcal S^{(2k+1)}_x$
are both projections of $\mathcal S^{(2k)}_x$.
\end{enumerate}
\end{enumerate}
It is easy to see that the relations defined in Definition \ref{conmultiequiv11}
(2) is an equivalence relation.
\end{defn}
\begin{equation}
\xymatrix{
& \mathcal S_x^{(1)} \ar[ld]\ar[rd] && \cdots\ar[ld]\ar[rd] &&
\mathcal S_x^{(2N-1)} \ar[ld]\ar[rd]  \\
\mathcal S_x^{(0)} && \mathcal S_x^{(2)} & \cdots&
\mathcal S_x^{(2N-2)} && \mathcal S_x^{(2N)}
}
\nonumber\end{equation}

\begin{defn}\label{contipertlocalrest}
Let $
\Phi_{xx'} = (h_{xx'},\widetilde\varphi_{xx'},\breve\varphi_{xx'})
$ be a coordinate change from $\frak V_{x'}$ to $\frak V_{x}$.
(See  Situation \ref{opensuborbifoldchart}.)
Let $\mathcal S_x = (W_x,\omega_x,{\frak s}_{x}^{\epsilon})$
be a CF-perturbation
of $\mathcal U$ on $\frak V_x$.
We define its {\it restriction} $\Phi_{xx'}^*\mathcal S_x$ by
\index{CF-perturbation ! restriction by
coordinate change of orbifold charts}
$$
\Phi_{xx'}^*\mathcal S_x
=
(W_x,\omega_x,\{{\frak s}^{\epsilon \prime}_{x'}\})
$$
where ${\frak s}^{\epsilon\prime}_{x'}$ is defined as follows.
We associate a linear isomorphism
$
g_y : E_{x'} \to E_x
$
to each $y \in V'_{x'}$ by
$
\breve\varphi_{xx'}(y,v) = g_y(v).
$
Then we put
\begin{equation}\label{pullbacklocaldefcont}
{\frak s}^{\epsilon \prime}_{x'}(y,\xi)
= g_y^{-1}({\frak s}^{\epsilon}_{x}(\tilde\varphi_{xx'}(y),\xi)).
\end{equation}
\end{defn}
\begin{lem}\label{lem77}
\begin{enumerate}
\item
If $\mathcal S^1_x$ is equivalent to $\mathcal S^2_x$, then
$
\Phi_{xx'}^*\mathcal S^1_x
$
is equivalent to
$
\Phi_{xx'}^*\mathcal S^2_x
$.
\item
The restriction $\Phi_{xx'}^*\mathcal S_x$ may
depend on the choice of $\Phi_{xx'} = (h_{xx'},\tilde\varphi_{xx'},\breve\varphi_{xx'})$.
However the equivalence class of $\Phi_{xx'}^*\mathcal S_x$
is independent of such a choice.
\end{enumerate}
\end{lem}
\begin{proof}
(1) is obvious as far as we use the same
$(h_{xx'},\tilde\varphi_{xx'},\breve\varphi_{xx'})$.
We will prove (2).
Let
$(h^{i}_{xx'},\tilde\varphi^{i}_{xx'},\breve\varphi^{i}_{xx'})$
($i=1,2$) be two choices and
${\frak s}^{\epsilon i \prime}_{x'}(y,\xi)$  the restrictions
obtained by these two choices for $i=1,2$, respectively.
Then by Lemma \ref{lem2715} there exists $\gamma \in \Gamma_x$ such that
$$
h^{2}_{xx'}(\mu) = \gamma h^{1}_{xx'}(\mu) \gamma^{-1},
\quad
\tilde\varphi^{2}_{xx'} = \gamma \tilde\varphi^{1}_{xx'},
\quad
\breve\varphi^{2}_{xx'} = \gamma \breve\varphi^{1}_{xx'}.
$$
We put
$
\breve\varphi^{i}_{xx'}(y,v) = g^{i}_y(v).
$
Then
$
g^2_y(v) = \gamma g^1_y(v).
$
Therefore (\ref{pullbacklocaldefcont}) implies
\begin{equation}\label{lem78profcal}
\aligned
{\frak s}^{\epsilon 2 \prime}_{x'}(y,\xi)
&= (g_y^2)^{-1}({\frak s}^{\epsilon \prime}_{x}(\tilde\varphi^2_{xx'}(y),\xi))
\\
&=(g_y^1)^{-1}\gamma^{-1}({\frak s}^{\epsilon\prime}_{x}(\gamma\tilde\varphi^1_{xx'}(y),\xi))\\
&=(g_y^1)^{-1}({\frak s}^{\epsilon 1\prime}_{x}(\tilde\varphi^1_{xx'}(y),\gamma^{-1}\xi))
\\
&={\frak s}^{\epsilon 1 \prime}_{x'}(y,\gamma^{-1}\xi).
\endaligned
\end{equation}
We note that $\gamma^{-1}$ induces a $\Gamma_{x'}$
linear isomorphism from $(W_x,h^2_{xx'})$ to  $(W_x,h^1_{xx'})$.
(Here the $\Gamma_{x'}$ action on $W_x$ is induced
from the $\Gamma_{x}$ action by $h^{i}_{xx'}$ in case of $(W_x,h^i_{xx'})$,
$i=1,2$.
Then the map $\xi \mapsto \gamma^{-1}\xi$ is $\Gamma_x$ equivariant as a map
$W_x \to W_x$, where $\Gamma_x$ acts in a different way on the source and the target.)
Therefore ${\frak s}^{\epsilon 2 \prime}_{x'}(y,\xi)$ is equivalent to
${\frak s}^{\epsilon 1 \prime}_{x'}(y,\xi)$.
\end{proof}
\par
We next define the pushout of a differential form by
using a CF-perturbation.

\begin{defn}\label{submersivepertconlocloc}
In Situation \ref{smoothcorrsingle},
let $\mathcal S_x = (W_x,\omega_x,\{{\frak s}_{x}^{\epsilon}\})$
be a CF-perturbation
of $\mathcal U$ on $\frak V_x$.
\begin{enumerate}
\item
We say $\mathcal S_x$ is {\it transversal to $0$}
\index{CF-perturbation ! transversal to $0$ !
on one orbifold chart}
if there exists $\epsilon_0 > 0$ such that the map ${\frak s}_{x}^{\epsilon}$ is transversal to $0$ on
a neighborhood of the support of $\omega_x$ for all $0 < \epsilon < \epsilon_0$.
In particular
$$
({\frak s}_{x}^{\epsilon})^{-1}(0)
=
\{(y,\xi) \in V_x \times W_x\mid {\frak s}_{x}^{\epsilon}(y,\xi) = 0\}
$$
is a smooth submanifold of $V_x \times W_x$ on
a neighborhood of the support of $\omega_x$.
\item
We say $f_x = f \circ \phi_x$ is {\it strongly submersive} with respect to $(\frak V_x,\mathcal S_x)$ if $\mathcal S_x$ is transversal to $0$ and
there exists $\epsilon_0 > 0$ such that the map
\begin{equation}\label{form7575}
f_x \circ \pi_1\vert_{({\frak s}_{x}^{\epsilon})^{-1}(0)} :
({\frak s}_{x,k}^{\epsilon})^{-1}(0) \to M
\end{equation}
is a submersion on
a neighborhood of the support of $\omega_x$,
for all $0 < \epsilon < \epsilon_0$. Here $\pi_1 : V_x \times W_x \to V_x$ is the
projection.
\index{strongly submersive (w.r.t. CF-perturbation) !
on one chart}
\item
Let $g : N \to M$ be a smooth map between manifolds.
We say {\it $f_x$ is strongly transversal to $g$} with respect to $(\frak V_x,\mathcal S_x)$
if \index{transversal (w.r.t. CF-perturbation) ! to a map on one chart}
$\mathcal S_x$ is transversal to $0$ and there exists $\epsilon_0 > 0$ such that  the map
(\ref{form7575}) is transversal to $g$, for all $0 < \epsilon < \epsilon_0$. Here $\pi_1 : V_x \times W_x \to V_x$ is the
projection.
\end{enumerate}
\end{defn}
\begin{lem}\label{newlem79}
Suppose $\mathcal S^1_x$ is equivalent to $\mathcal S^2_x$.
\begin{enumerate}
\item
$\mathcal S^1_x$ is transversal to $0$ if and only if  $\mathcal S^2_x$
is transversal to $0$.
\item
$f_x$ is strongly submersive with respect to $(\frak V_x,\mathcal S^1_x)$ if and only if
$f_x$ is strongly submersive with respect to $(\frak V_x,\mathcal S^2_x)$.
\item
$f_x$ is strongly transversal  to $g : N \to M$ with respect to $(\frak V_x,\mathcal S^1_x)$ if and only if
$f_x$ is strongly transversal  to $g$ with respect to $(\frak V_x,\mathcal S^2_x)$.
\end{enumerate}
\end{lem}
\begin{proof}
It suffices to prove the lemma for the case when $\mathcal S^2_x$
is a projection of $\mathcal S^1_x$.
This case follows
from the fact that $\omega_x = \vert\omega_x\vert{\rm vol}_x$ where ${\rm vol}_x$ is a volume form of $W_x$
and $\vert\omega_x\vert$ is a non negative function, which is a part of Definition \ref{defn73ss} (4).
\end{proof}
We recall that
a smooth differential form $h$ on $V_x/\Gamma_x$
is identified with a $\Gamma_x$ invariant
smooth differential form $\tilde h$ on $V_x$.
(Definition \ref{defn281010} (2).)
The {\it support} \index{support ! of differential form
on orbifold}\index{differential form ! support of differential form
on orbifold} of $h$ is the quotient of the support of
$\tilde h$ by $\Gamma_x$ and is a closed subset of $V_x/\Gamma_x
\cong U_x$.
We denote it by ${\rm Supp}(h)$.

\begin{defn}
In Situation \ref{smoothcorrsingle},
let $\mathcal S_x = (W_x,\omega_x,\{{\frak s}_{x}^{\epsilon}\})$
be a CF-perturbation
of $\mathcal U$ on $\frak V_x$.
Let $h$ be a smooth differential form on $U_x$ that has compact support.
Then we define a smooth differential form $f_x!(h;\mathcal S^{\epsilon}_x)$
on $M$
for each $\epsilon > 0$ by the equation (\ref{form72}) below.
We call it the {\it pushout} of $h$ with respect to $f_x$, $\mathcal S_x$.
\index{integration along the fiber (pushout) ! of differential form - one chart}
\par
Let $\rho$ be a smooth differential form on $M$. Then
we require
\begin{equation}\label{form72}
\# \Gamma_x
\int_{M} f_x!(h;\mathcal S^{\epsilon}_x) \wedge \rho
=
\int_{({\frak s}_{x}^{\epsilon})^{-1}(0)}
\pi_1^*(\tilde h) \wedge \pi_1^*(f_x^*\rho) \wedge \pi_2^*\omega_x.
\end{equation}
Here $\pi_1$ (resp. $\pi_2$) is the projection of $V_x \times W_x$ to the
first (resp. second) factor.
\end{defn}
Unique existence of such $f_x!(h;\mathcal S^{\epsilon}_x)$
is an immediate consequence of the existence of pushforward
of a smooth form by a proper submersion.
\begin{rem}
In the left hand side of (\ref{form72}) we crucially use the fact
that $\Gamma_x$ is a finite group. It seems that
this is the {\it only} place we use the finiteness of $\Gamma_x$
when we use de Rham theory to realize virtual fundamental chain.
We might try to use a CF-perturbation
and de Rham version together with an appropriate model
of equivariant cohomology to study
virtual fundamental chain in case the isotropy group
can be a continuous compact group of positive dimension,
such as the case of gauge theory or pseudo-holomorphic curves
in a symplectic manifold acted by a compact Lie group.
\end{rem}
\begin{lem}\label{pushequivlocal}
If $\mathcal S_x^1$ is equivalent to $\mathcal S_x^2$, then
$$
f_x!(h;\mathcal S^{1,\epsilon}_x)
=
f_x!(h;\mathcal S^{2,\epsilon}_x).
$$
\end{lem}
\begin{proof}
It suffices to prove the equality in the case when
$\mathcal S_x^1$ is a projection of $\mathcal S_x^2$.
This is immediate from  definition.
\end{proof}
\begin{lem}\label{resandtrans}
Suppose we are in the situation of Definition \ref{contipertlocalrest}
and Situation \ref{smoothcorrsingle}.
\begin{enumerate}
\item
If $(\frak V_x,\mathcal S_x)$ is transversal to $0$, then
 $(\frak V_{x'},\Phi_{xx'}^*\mathcal S_x)$ is transversal to $0$.
\item
If
$f_x$ is strongly submersive with respect to
$(\frak V_x,\mathcal S_x)$,
then $f_{x'}$ is strongly submersive with
respect to $(\frak V_{x'},\Phi_{xx'}^*\mathcal S_x)$.
\item
If
$f_x$ is strongly transversal to $g : N \to M$ with respect to
$(\frak V_x,\mathcal S_x)$,
then
$f_{x'}$ is strongly transversal to $g : N \to M$ with respect to
 $(\frak V_{x'},\Phi_{xx'}^*\mathcal S_x)$.
\end{enumerate}
\end{lem}
The proof is immediate from definition.
\subsubsection{CF-perturbation on a single Kuranishi chart}
\label{cftonekurachart}
In Subsubsection \ref{cftoneorbifoldchart}, we studied locally on a single
chart $U_x = V_x/\Gamma_x$.
We next work globally on an orbifold $U$.
We apply Remark \ref{xkararhe} hereafter.
\begin{rem}
In this subsubsection we consider a Kuranish chart $\mathcal U = (U,\mathcal E,s,\psi)$.
However the parametrization $\psi$ does not play any role in this subsubsection.
\end{rem}
\begin{defn}\label{semiglobalocntpert}
Let $\mathcal U =(U,\mathcal E,s,\psi)$ be a Kuranishi chart.
A {\it representative of a CF-perturbation of $\mathcal U$}
 \index{CF-perturbation ! representative on Kuranishi chart} is the following object $\frak S = \{(\frak V_{\frak r},\mathcal S_{\frak r})\mid{\frak r\in
\frak R}\}$.
\begin{enumerate}
\item
$\{U_{\frak r} \mid \frak r \in \frak R\}$ is a
family of open subsets of $U$
such that
$
U = \bigcup_{\frak r \in \frak R}U_{\frak r}.
$
\item
$\frak V_{\frak r} = (V_{\frak r},\Gamma_{\frak r},E_{\frak r},\phi_{\frak r},\widehat\phi_{\frak r})$ is an orbifold chart of $(U,\mathcal E)$ such that
$\phi_{\frak r}(V_{\frak r}) = U_{\frak r}$.
\item
$\mathcal S_{\frak r}  = (W_{\frak r} ,\omega_{\frak r}, \{{\frak s}_{\frak r} ^{\epsilon}\})$ is a
CF-perturbation of $\mathcal U$
on $\frak V_{\frak r}$.
\item
For each $x \in U_{\frak r_1} \cap U_{\frak r_2}$, there exists
an orbifold chart $\frak V_{\frak r}$ with the
following properties:
\begin{enumerate}
\item
$
x \in U_{\frak r} \subset U_{\frak r_1} \cap U_{\frak r_2}.
$
\item
The restriction of $\mathcal S_{\frak r_1}$ to $\frak V_{\frak r}$
is equivalent to the restriction of $\mathcal S_{\frak r_2}$ to $\frak V_{\frak r}$.
\end{enumerate}
\end{enumerate}
For each $\epsilon >0$, we write
$$
\frak S^{\epsilon}
= \{(\frak V_{\frak r},\mathcal S_{\frak r}^{\epsilon})\mid{\frak r\in
\frak R}\}.
$$
See Definition \ref{defn73ss} for the notation $\mathcal S^{\epsilon}_{\frak r}.$
\end{defn}
\begin{defn}\label{projectcontfamilocal}
Let $\mathcal U = (U,\mathcal E,s,\psi)$ be as in Definition \ref{semiglobalocntpert}
and $\frak S^i = \{(\frak V_{\frak r}^i,\mathcal S_{\frak r}^i)\mid{\frak r\in
\frak R^i}\}$ ($i=1,2$)  representatives of
CF-perturbations of $\mathcal U$.
We say that $\frak S^2$ is {\it equivalent} to $\frak S^1$
if, for each $x \in U_{\frak r_1} \cap U_{\frak r_2}$, there exists
an orbifold chart $\frak V_{\frak r}$ with the
following properties:
\begin{enumerate}
\item
$x \in U_{\frak r} \subset U_{\frak r_1} \cap U_{\frak r_2}.$
\item
The restriction of $\mathcal S^1_{\frak r_1}$ to $\frak V_{\frak r}$
is equivalent to the restriction of $\mathcal S^2_{\frak r_2}$ to $\frak V_{\frak r}$.
\end{enumerate}
\end{defn}
\begin{defn}\label{defn71717}
Suppose we are in Situation \ref{smoothcorrsingle}. Let
$\frak S = \{\frak S_{\frak r}\}$
be a representative of a CF-perturbation of $\mathcal U$.
Let $U' \subseteq U$ be an open subset.
\par
Let $\frak S_{\frak r} = (\frak V_{\frak r},\mathcal S_{\frak r})$,
$\frak V_{\frak r} = (V_{\frak r},\Gamma_{\frak r},E_{\frak r},\phi_{\frak r},\widehat\phi_{\frak r})$.
If ${\rm Im}(\phi_{\frak r}) \cap U' = \emptyset$, then we remove $\frak r$ from $\frak R$.
Let $\frak R_0$ be obtained by removing all such $\frak r$ from $\frak R$.
If ${\rm Im}(\phi_{\frak r}) \cap U' \ne \emptyset$, then $\frak V_{\frak r}\vert_{{\rm Im}\phi_{\frak r} \cap U'}$ is an orbifold chart
of $(U,\mathcal E)$, which we write $\frak V_{\frak r}\vert_{U'} = (V'_{\frak r},\Gamma'_{\frak r},E_{\frak r},\phi'_{\frak r},\widehat\phi'_{\frak r})$.
Let $\mathcal S_{\frak r} = (W_{\frak r} ,\omega_{\frak r} ,\{{\frak s}^{\epsilon}_{\frak r}\})$.
We define
\begin{equation}\label{formform7676}
\mathcal S_{\frak r}\vert_{U'}
=
(W_{\frak r} ,\omega_{\frak r} ,\{{\frak s}^{\epsilon}_{\frak r}\vert_{V'_{\frak r} \times W_{\frak r}}\}).
\end{equation}
Now we put:
$$
\frak S\vert_{U'}
=
\{(\frak V_{\frak r}\vert_{U'},\mathcal S_{\frak r}\vert_{U'}) \mid \frak r \in \frak R_0 \}.
$$
It is a representative of a CF-perturbation of $\mathcal U\vert_{U'}$,
which we call the {\it restriction}
\index{CF-perturbation ! restriction to an open subset of orbifold} of $\frak S$ to $U'$.
\end{defn}
\begin{lem}\label{lem718718}
If $\frak S^1$ is equivalent to $\frak S^2$, then
$\frak S^1\vert_{\Omega}$ is equivalent to $\frak S^2\vert_{\Omega}$.
\end{lem}
The proof is immediate from definition.
\begin{defn}\label{rem718}
Suppose we are in Situation \ref{smoothcorrsingle}.
\begin{enumerate}
\item
A {\it CF-perturbation}\index{CF-perturbation ! CF-perturbation} on $\mathcal U$ is an equivalence class of
a representative of a CF-perturbation
\index{CF-perturbation ! on orbifold}
with respect to the equivalence relation in Definition \ref{projectcontfamilocal}.
\item
Let $\Omega$ be an open subset of $U$.
We denote by $\mathscr{S}(\Omega)$ the set of all CF-perturbations on  $\mathcal U\vert_{\Omega}$.
\item
Let $\Omega_1 \subset \Omega_2 \subset U$. Then using Lemma \ref{lem718718}, the restriction defined in
Definition \ref{defn71717} induces a map
\begin{equation}
\frak i_{\Omega_1\Omega_2} : \mathscr{S}(\Omega_2) \to \mathscr{S}(\Omega_1).
\end{equation}
We call this map {\it restriction map}\index{sheaf $\mathscr{S}$ ! restriction map}.
\item
$\Omega \mapsto \mathscr{S}(\Omega)$ together with $\frak i_{\Omega_1 \Omega_2}$ defines a presheaf.
(In fact, $\frak i_{\Omega_1\Omega_2}\circ \frak i_{\Omega_2\Omega_3} = \frak i_{\Omega_1\Omega_3}$ holds obviously.)
The next proposition says that it is indeed a sheaf.
We call it the {\it  sheaf of CF-perturbations} on $\mathcal U$
\index{CF-perturbation ! sheaf of}
\index{sheaf $\mathscr{S}$ ! of CF-perturbations}.
\end{enumerate}
\end{defn}
\begin{rem}
Hereafter, by a slight abuse of notation, we also use the symbol $\frak S$ for a CF-perturbation.
Namely we use this symbol not only for a representative of continuous family perturbation but also
for its equivalence class, by an abuse of notation.
\end{rem}
\begin{prop}\label{prop721}
The presheaf $\mathscr{S}$ is a sheaf.
\end{prop}
\begin{proof}
This is mostly a tautology. We provide a proof for completeness' sake.
\par
Let
$
\bigcup_{a\in A} \Omega_{a} = \Omega
$
be an open cover of $\Omega$.
Suppose $\frak S_a \in \mathscr{S}(\Omega_a)$ and
$\{(\frak V_{a,\frak r},\mathcal S_{a,\frak r})\mid{\frak r\in
\frak R_a}\}$ is a representative of $\frak S_a$.
We put $\Omega_{ab} = \Omega_a \cap \Omega_{b}$.
We assume
\begin{equation}\label{form7878}
\frak i_{\Omega_{ab}\Omega_a}(\frak S_{a}) = \frak i_{\Omega_{ab}\Omega_b}(\frak S_{b}).
\end{equation}
To prove Proposition \ref{prop721} it suffices to show that there exists
a unique $\frak S \in \mathscr{S}(\Omega)$ such that
\begin{equation}\label{form7979}
\frak i_{\Omega_{a}\Omega}(\frak S) = \frak S_{a}.
\end{equation}
\begin{proof}[Proof of uniqueness]
Suppose that $\frak S, \frak S' \in \mathscr{S}(\Omega)$ both satisfy (\ref{form7979}).
Let
$\{(\frak V_{\frak r},\mathcal S_{\frak r})\mid{\frak r\in
\frak R}\}$ and
$\{(\frak V'_{\frak r'},\mathcal S'_{\frak r'})\mid{\frak r'\in
\frak R'}\}$
be representatives of $\frak S$ and $\frak S'$, respectively.
We will prove that they are equivalent.
\par
Let $x \in \Omega$.
There exist $\frak r \in \frak R$, $\frak r' \in \frak R'$
such that $x \in U_{\frak r} \cap U'_{\frak r'}$, where
$U_{\frak r} = {\rm Im}(\phi_{\frak r})$, $U'_{\frak r'} = {\rm Im}(\phi'_{\frak r'})$.
We take $a$ such that $x \in \Omega_a$.
By (\ref{form7979}), $\mathcal S_{\frak r}\vert_{U_{\frak r} \cap \Omega_a}$ is
equivalent to $\mathcal S'_{\frak r'}\vert_{U'_{\frak r'} \cap \Omega_a}$.
Therefore there exists an orbifold chart $\frak V_x$
such that
$U_x \subset U_{\frak r} \cap U'_{\frak r'}$ and
the restriction of $\mathcal S_{\frak r}\vert_{U_{\frak r} \cap \Omega_a}$
to $\frak V_x$ is equivalent to the restriction of
$\mathcal S'_{\frak r'}\vert_{U'_{\frak r'} \cap \Omega_a}$
to $\frak V_x$.
Thus the restriction of $\mathcal S_{\frak r}$
to $\frak V_x$ is equivalent to the restriction of
$\mathcal S'_{\frak r'}$
to $\frak V_x$.
\par
Since this holds for any $x \in \Omega$,
$\{(\frak V_{\frak r},\mathcal S_{\frak r})\mid{\frak r\in
\frak R}\}$ is equivalent to
$\{(\frak V'_{\frak r'},\mathcal S'_{\frak r'})\mid{\frak r'\in
\frak R'}\}$ by definition.
\end{proof}
\begin{proof}[Proof of existence]
Let
$\frak V_{\frak r,a} = (V_{\frak r,a},\Gamma_{\frak r,a},E_{\frak r,a},\phi_{\frak r,a},\widehat\phi_{\frak r,a})$
and $U_{\frak r,a} = \phi_{\frak r,a}(V_{\frak r,a})$.
Then
$\{U_{\frak r,a} \mid a \in A,\,\, \frak r \in \frak R_a\}$ is an open cover of $\Omega$.
We put
$$
\frak S = \coprod_{a \in A} \{(\frak V_{a,\frak r},\mathcal S_{a,\frak r})\mid{\frak r\in
\frak R_a}\}.
$$
To show $\frak S \in \mathscr{S}(\Omega)$ it suffices to check Definition \ref{semiglobalocntpert} (4).
This is a consequence of (\ref{form7878}) and the definitions of
the equivalence
of representatives of CF-perturbations and of the restriction.
We can check (\ref{form7979}) also by (\ref{form7878}) in the same way as
in the proof of uniqueness.
\end{proof}
The proof of Proposition \ref{prop721} is now complete.
\end{proof}
We recall that the {\it stalk}
\index{sheaf $\mathscr{S}$ ! stalk of the sheaf of CF-perturbations}
\index{CF-perturbation ! stalk of the sheaf of} $\mathscr S_x$ of the sheaf $\mathscr S$ at $x \in U$ is by definition
\begin{equation}\label{form72000}
\mathscr S_x = \varinjlim_{\Omega \ni x}\mathscr S(\Omega).
\end{equation}
\begin{defnlem}
The stalk $\mathscr S_x$ is identified with the set of the equivalence classes
of the equivalence relation defined in Item (2) on the set which is defined in
Item (1).
\begin{enumerate}
\item
We consider the set $\widetilde{\mathscr S}_x$ of pairs $(\frak V_{\frak r},\frak S_{\frak r})$ where
$\frak V_{\frak r}$ is an orbifold chart of $(U,\mathcal E)$ at $x$ and
$\frak S_{\frak r}$ is a CF-perturbation on $\frak V_{\frak r}$.
\item
Let $(\frak V_{\frak r},\frak S_{\frak r}), (\frak V_{\frak r'},\frak S_{\frak r'}) \in \widetilde{\mathscr S}_x$.
We say that they are equivalent if there exists an orbifold chart $\frak V_{x}$ at $x$ such that $U_x \subset U_{\frak r} \cap U_{\frak r'}$ and
the restriction of $\frak S_{\frak r}$ to $\frak V_{x}$ equivalent to
the restriction of $\frak S_{\frak r'}$ to $\frak V_{x}$.
\end{enumerate}
\end{defnlem}
The proof of obvious.
The next lemma is standard in sheaf theory.
\begin{lem}\label{lem723}
The set $\mathscr S(\Omega)$ is identified with
the following object:
\begin{enumerate}
\item
For each $x \in \Omega$ it associates $\frak S_x \in \mathscr S_x$.
\item
For each $x \in \Omega$, there exists a representative  $(\frak V_{\frak r},\frak S_{\frak r})$ of $\frak S_x$,
such that for each $y \in \phi_{\frak r}(V_{\frak r})$ the germ $\frak S_y$ is represented by
$(\frak V_{\frak r},\frak S_{\frak r})$.
\end{enumerate}
\end{lem}

\begin{defn}
Let $K \subseteq U$ be a closed subset. A {\it CF-perturbation
of $K \subseteq U$} \index{CF-perturbation ! on a closed set of an orbifold} is an element of the inductive limit:
$
\varinjlim_{\Omega \supset K}\mathscr S(\Omega)
$.
We denote the set of all CF-perturbations
of $K \subseteq U$ by $\mathscr S(K)$.
Namely
\begin{equation}\label{sheafoverclosedset}
\mathscr S(K) = \varinjlim_{U\supset \Omega \supset K;~ \Omega \text{ open}}\mathscr S(\Omega).
\end{equation}
\end{defn}

\subsection{Integration along the fiber (pushout) on a single Kuranishi chart}
\label{subsec:intonechart}

\begin{defnlem}\label{strosubsemiloc}
Suppose we are in  Situation \ref{smoothcorrsingle}.
Let $\Omega \subset U$ be an open subset and
let $\frak S \in \mathscr S(\Omega)$ be a
CF-perturbation.
We consider its representative $\{(\frak V_{\frak r},\mathcal S_{\frak r})\mid{\frak r\in
\frak R}\}$.
\begin{enumerate}
\item
We say that  $(\mathcal U,\frak S)$ is {\it transversal to $0$}
if, for each $\frak r$,
$(\frak V_{\frak r},\mathcal S_{\frak r})$ is transversal to $0$.
This is independent of the choice of representative.
\index{CF-perturbation ! transversal to $0$ ! on orbifold}
\item
We say that $f$ is {\it strongly submersive with respect to $(\mathcal U,\frak S)$}
if, for each $\frak r$, the map
$f$ is strongly submersive with respect to $(\frak V_{\frak r},\mathcal S_{\frak r})$.
This is independent of the choice of representative.
\index{strongly submersive (w.r.t. CF-perturbation) ! on orbifold}
\item
Let $g : N \to M$ be a smooth map between manifolds.
We say that the map $f$ is {\it strongly transversal to $g$ with respect to $(\mathcal U,\frak S)$}
if, for each $\frak r$,
$f$ is strongly transversal to $g$ with respect to $(\frak V_{\frak r},\mathcal S_{\frak r})$.
This is independent of the choice of representative.
\index{transversal (w.r.t. CF-perturbation) ! to a map on orbifold}
\item
We denote by
$
\mathscr S_{\pitchfork 0}(\Omega)
$
the set of all $\frak S \in \mathscr S(\Omega)$  transversal to $0$,
$
\mathscr S_{f \pitchfork}(\Omega)
$
the set of all $\frak S \in \mathscr S(\Omega)$ such that $f$ is  strongly submersive with respect to $(\mathcal U,\frak S)$,
and by
$
\mathscr S_{f \pitchfork g}(\Omega)
$
the set of all $\frak S \in \mathscr S(\Omega)$ such that $f$ is  strongly transversal to $g$ with respect to $(\mathcal U,\frak S)$.
They are subsheaves of $\mathscr S$.
\item
For a closed subset $K \subseteq U$
we define $
\mathscr S_{\pitchfork 0}(K)
$,
$
\mathscr S_{f \pitchfork}(K)
$
and $
\mathscr S_{f \pitchfork g}(K)
$
in the same way as (\ref{sheafoverclosedset}).
\end{enumerate}
\end{defnlem}
\begin{proof}
The statements
(1), (2), and (3) follow from
Lemma \ref{newlem79}. (4) is a consequence of the definition.
\end{proof}
Next we introduce the notion of
{\it member of a CF-perturbation}, which is an analogue of the notion of branch of
\index{CF-perturbation ! member of}
\index{member !{\it see CF-perturbation}} multisection.
\begin{defn}\label{memberonechart}
Let $\frak V_{\frak r}$ be an orbifold chart of $(U,\mathcal E)$
and $\mathcal S_{\frak r}$  a CF-perturbation of
$(\mathcal U,\frak V_{\frak r})$.
Let $x \in V_{\frak r}$. We take $\epsilon \in (0,1]$.
\par
Consider the germ of a map $y \mapsto \frak s(y)$,
$O_x \to E_{\frak r}$
(where $O_x$ is a neighborhood of $x \in V_{\frak r}$).
We say $\frak s$ is a {\it member of $\mathcal S^{\epsilon}_{\frak r}$ at $x$}
\index{CF-perturbation ! member of} if there exists
$\xi \in W_x$
such that the germ of $y \mapsto \frak s_{\frak r}^{\epsilon}(y,\xi)$
at $x$ is $\frak s$ and $\xi \in {\rm supp}\,\omega_x$.
(See Definition \ref{defn73ss} for the notation $\mathcal S^{\epsilon}_{\frak r}$.)
\end{defn}
\begin{rem}
The member of $\mathcal S^{\epsilon}_{\frak r}$ at $x$ depends on $\epsilon$.
In other words, we define the notion of the member of $\mathcal S^{\epsilon}$
for each $\epsilon  \in (0,1]$.
\end{rem}
\begin{lem}\label{lem721}
In the situation of Definition \ref{memberonechart},
let $\mathcal S'_{\frak r}$ be a
CF-perturbation of $(\mathcal U,\frak V_{\frak r})$
that it is equivalent to $\mathcal S_{\frak r}$.
Then $\frak s$ is a member of $\mathcal S^{\epsilon}_{\frak r}$
if and only if $\frak s$ is a member of $\mathcal S^{\prime\epsilon}_{\frak r}$.
\end{lem}
\begin{proof}
It suffices to consider only the case when $\mathcal S'_{\frak r}$
is a projection of $\mathcal S_{\frak r}$. This also follows
from the fact that $\omega_x = \vert\omega_x\vert{\rm vol}_x$ where ${\rm vol}_x$ is a volume form of $W_x$
and $\vert\omega_x\vert$ is a non negative function, which is a part of Definition \ref{defn73ss} (4).
\end{proof}
Therefore we can define an element of the stalk $\mathscr S_x$ at $x$
of the sheaf $\mathscr S$ of CF-perturbations.
\par
We now recall from Definition \ref{defn73ss} that we denoted by
$
\mathcal S_x^\epsilon = (W_x, \omega_x,\frak s_x^\epsilon)
$
the restriction of the CF-perturbation
$\mathcal S_x = (W_x, \omega_x,\{\frak s_x^\epsilon\})$ at $\epsilon$. Recall by definition that
$\{\frak s_x^\epsilon\}$ is an $\epsilon$-dependent family of parameterized section, i.e.,
a map $\frak s_x^\epsilon: V_x \times W_x \to E_x$.
\begin{defn}
Let $\frak S \in \mathscr S(\Omega)$ and $x \in \Omega$.
A {\it member of $\frak S^{\epsilon}$ at $x$}
\index{CF-perturbation ! member of} is a member of the germ of $\frak S^{\epsilon}$ at $x$.
\end{defn}
\begin{defn}
Let $U$ be an orbifold and $K\subset U$ a compact subset
and $\{U_{\frak r}\}$ a set of finitely many open subsets such that $\cup_{\frak r} U_{\frak r} \supset K$.
A {\it partition of unity} \index{orbifold ! partition of unity} subordinate to $\{U_{\frak r}\}$ on $K$
assigns a smooth function
$\chi_{\frak r}$ on $U$ to each $\frak r$ such that:
\begin{enumerate}
\item
${\rm supp} \chi_{\frak r}\subset U_{\frak r}$.
\item
$
\sum_{\frak r\in\frak R} \chi_{\frak r}
 \equiv 1$ on a neighborhood of $K$.
\end{enumerate}
\end{defn}
It is standard and easy to prove that a partition of unity always exists on an orbifold.
\begin{defn}
Suppose we are in Situation
\ref{smoothcorrsingle}.
Let $h$ be a smooth differential form of compact support in $U$.
Let $\frak S \in \mathscr S_{f \pitchfork}({\rm Supp}(h))$.
Let $\{\chi_{\frak r}\}$ be a smooth partition of unity subordinate to
the covering $\{U_{\frak r}\}$ on ${\rm Supp}(h)$.
We define the {\it pushout of $h$
\index{integration along the fiber (pushout) ! of differential form on single Kuranishi chart} by $f$ with respect to  $\frak S^{\epsilon}$}
by
$$
f!(h;\frak S^{\epsilon})
=
\sum_{\frak r\in\frak R}
f!(\chi_{\frak r}h;\mathcal S^{\epsilon}_{\frak r})
$$
for each $\epsilon >0$.
It is a smooth form on $M$ of degree
$$
\deg f!(h,\mathcal S^{\epsilon})
= \deg h + \dim M - \dim \mathcal U,
$$
where
$
\dim \mathcal U = \dim U - \rank \mathcal E.
$
The well-defined-ness follows from Lemma \ref{lem721} (3).
\par
We also call pushout {\it integration along the fiber.}
\index{integration along the fiber (pushout) ! integration along the fiber}
\end{defn}
\begin{rem}
In general
$f!(h;\frak S^{\epsilon})$ {\it depends} on $\epsilon$.
Moreover $\lim_{\epsilon\to 0}f!(h;\frak S^{\epsilon})$ typically diverges.
\end{rem}
\begin{lem}\label{lem721}
\begin{enumerate}
\item
$
f!(h_1+h_2;\frak S^{\epsilon})
=
f!(h_1;\frak S^{\epsilon})
+
f!(h_2;\frak S^{\epsilon})
$
and
$
f!(ch;\frak S^{\epsilon}) = cf!(h;\frak S^{\epsilon})
$
for $c \in \R$.
\item
Pushout of $h$ is independent of the choice of partition of
unity.
\item
If $\frak S^1$ is equivalent to $\frak S^2$ then
\begin{equation}\label{icchiprojectint}
f!(h;\frak S^{1,\epsilon})
=
f!(h;\frak S^{2,\epsilon}).
\end{equation}
\end{enumerate}
\end{lem}
\begin{proof}
(1) is obvious as far as we use the same partition of unity on both
sides.
We will prove (2) and (3) at the same time.
We take a partition of unity $\{\chi^i_{\frak r} \mid \frak r \in \frak R_i\}$
subordinate to $\frak S^i$
for $i=1,2$ and will prove (\ref{icchiprojectint}). Here
we use those partitions of unity to define the left and right hand sides of
(\ref{icchiprojectint}), respectively.
The case $\frak S^1 = \frak S^2$ will be (2).
\par
We put $h_0 = \chi_{\frak r_0}\widehat h$.
In view of (1) it suffices to prove
\begin{equation}\label{foru735}
f!(h_0;\mathcal S^{1,\epsilon}_{\frak r_0})
=
\sum_{\frak r \in \frak R_2}
f!(\chi^2_{\frak r}h_0;\mathcal S^{2,\epsilon}_{\frak r}).
\end{equation}
To prove (\ref{foru735}), it suffices to show
the next equality for each $\frak r \in \frak R_2$.
\begin{equation}\label{foru74}
f!(\chi^2_{\frak r}h_0;\mathcal S^{1,\epsilon}_{\frak r_0})
=
f!(\chi^2_{\frak r}h_0;\mathcal S^{2,\epsilon}_{\frak r}).
\end{equation}
We put $h_1 = \chi^2_{\frak r}h_0$.
Let $K = {\rm Supp}(h_1)$.
For each $x \in K$ there exists $\frak V_x$ such that
the restriction of $\mathcal S^{1,\epsilon}_{\frak r_0}$
to $U_{x}$ is equivalent to the restriction of
$\mathcal S^{2,\epsilon}_{\frak r}$ to $U_x$.
We cover $K$ by a finitely many such $U_{x_i}$,
$i=1,\dots,N$.
Let $\{\chi'_i \mid i=1,\dots,N\}$ be a partition of unity on $K$
subordinate to
the covering $\{U_{x_i}\}$.
Then we obtain:
$$
\aligned
f!(\chi^2_{\frak r}h_0;\mathcal S^{1,\epsilon}_{\frak r_0})
&=
\sum_{i=1}^N
f!(\chi'_{i}h_1;\mathcal S^{1,\epsilon}_{\frak r_0}\vert_{U_{x_i}})\\
&=
\sum_{i=1}^N
f\vert_{U_{x_i}}!(\chi'_{i}h_1;\mathcal S^{2,\epsilon}_{\frak r}\vert_{U_{x_i}})
=
f!(\chi^2_{\frak r}h_0 ;\mathcal S^{2,\epsilon}_{\frak r}).
\endaligned$$
The proof of Lemma \ref{lem721} is now complete.
\end{proof}

\subsection{CF-perturbations of good coordinate system}
\label{subsec:conpergcs}

To consider a CF-perturbation of good coordinate system, we
first study the pullback of a CF-perturbation by an
embedding of Kuranishi charts in this subsubsection.
Using the pullback, we then define the notion of a CF-perturbation of
good coordinate system.

\subsubsection{Embedding of Kuranishi charts and CF-perturbations}
\label{subsub:embcfp}
\begin{notation}
To highlight the dependence on the Kuranishi structure, we write $\mathscr S^{\mathcal U^1}$, $\mathscr S^{\mathcal U^2}$ etc. in place of  $\mathscr S$.
Namely $\mathscr S^{\mathcal U^2}(\Omega)$ is the set of all continuous family perturbations of $\mathcal U^2\vert_{\Omega}$.
\end{notation}
\begin{shitu}\label{contfamipullbacksitu}
Let $\mathcal U^i = (U^i,\mathcal E^i,s^i,\psi^i)$ $(i=1,2)$ be Kuranishi charts and
$\Phi_{21} = (\varphi_{21},\widehat\varphi_{21}): \mathcal U^1 \to \mathcal U^2$
an embedding of Kuranishi charts.
Let
$\frak S^2 \in \mathscr{S}^{\mathcal U^2}(U^2)$ and let $\{\frak S^2_{\frak r}\} = \{ (\frak V^2_{\frak r},\mathcal S^2_{\frak r})\}$
be its representative.$\blacksquare$
\end{shitu}

We will define the pullback $\Phi_{21}^*\frak S^2  \in \mathscr{S}
^{\mathcal U^1}(U_1)$.
We need certain conditions for this pullback to be defined.
\begin{conds}\label{exitpullbackcont1}
Suppose we are in Situation \ref{contfamipullbacksitu}.
We require that there exists
an orbifold chart $\frak V^1_{\frak r} = (V^1_{\frak r},\Gamma^1_{\frak r},E^1_{\frak r},\phi^1_{\frak r},\widehat\phi^1_{\frak r})$
such that $\phi^1_{\frak r}(V^1_{\frak r}) = \varphi_{21}^{-1}(U_{\frak r}^2)$.
 \end{conds}
(Recall $U_{\frak r}^2 = \phi_{\frak r}^2(V_{\frak r}^2)$ and
$\frak V^2_{\frak r} = (V^2_{\frak r},\Gamma^2_{\frak r},E^2_{\frak r},\phi^2_{\frak r},\widehat\phi^2_{\frak r})$.)
\begin{lem}
For any given $\{ (\frak V^2_{\frak r},\mathcal S^2_{\frak r})\}$ there exists
$\{ (\frak V^{2 \prime}_{\frak r},\mathcal S^{2 \prime}_{\frak r})\}$ which is  equivalent to $\{ (\frak V^2_{\frak r},\mathcal S^2_{\frak r})\}$ and satisfies Condition \ref{exitpullbackcont1}.
\end{lem}
\begin{proof}
For each $x \in U_2$ there exists $\frak V_x$
such that
$\varphi_{21}^{-1}(U_{x}^2)$ has an orbifold chart
and $\frak V_x \subset \frak V_{\frak r}$ for some $\frak r$.
We cover $U_2$ by such $\frak V_x$ to obtain the
required $\{ (\frak V^{2 \prime}_{\frak r},\mathcal S^{2 \prime}_{\frak r})\}$.
\end{proof}
Therefore we may assume Condition \ref{exitpullbackcont1}.
Then we can represent the orbifold embedding
$(\varphi_{21},\widehat\varphi_{21}) : (U^1,\mathcal E^1) \to (U^2,\mathcal E^2)$  in terms of the  orbifold charts
$\frak V^1_{\frak r}$, $\frak V^2_{\frak r}$ by
$(h^{\frak r}_{21},\tilde\varphi^{\frak r}_{21},\breve\varphi^{\frak r}_{21})$
that have the following properties.
\begin{proper}
\begin{enumerate}\label{proper728}
\item
$h^{\frak r}_{21} : \Gamma_1^{\frak r} \to   \Gamma_2^{\frak r}$ is an injective group homomorphism.
\item
$\tilde\varphi^{\frak r}_{21} : V_1^{\frak r} \to  V_2^{\frak r}$ is an $h^{\frak r}_{21}$-equivariant smooth
embedding of manifolds.
\item
$h^{\frak r}_{21}$ and $\tilde\varphi^{\frak r}_{21}$ induce an orbifold embedding
$$
\left(\overline{\phi_2^{\frak r}}\right)^{-1} \circ \varphi_{21} \circ \overline{\phi^{\frak r}_{1}} :  V^1_{\frak r}/\Gamma^1_{\frak r}
\to  V^2_{\frak r}/\Gamma^2_{\frak r}.
$$
\item $\breve\varphi^{\frak r}_{21} : V_1^{\frak r} \times  E_1^{\frak r} \to  E_2^{\frak r}$ is an $h^{\frak r}_{21}$-equivariant smooth
map such that for each $y\in V_1^{\frak r}$ the map $v \mapsto \breve\varphi^{\frak r}_{21}(y,v)$
is a linear embedding $E_1^{\frak r} \to  E_2^{\frak r}$.
\item
$\breve\varphi^{\frak r}_{21}$ induces a smooth embedding of vector bundles:
$$
\left(\overline{\widehat\phi^{\frak r}_{2}}\right)^{-1} \circ \widehat\varphi_{21} \circ \overline{\widehat\phi^{\frak r}_{1}} : (V^1_{\frak r} \times E^1_{\frak r})/\Gamma^1_{\frak r}
\to(V^2_{\frak r} \times E^2_{\frak r})/\Gamma^2_{\frak r}.
$$
In other words for each $(y,v) \in V^1_{\frak r} \times E^1_{\frak r}$ we have
$$
\widehat\varphi_{21}(\widehat\phi^{\frak r}_{1}(y,v))
= \widehat\phi^{\frak r}_{2}(\tilde\varphi^{\frak r}_{21}(y),\breve\varphi^{\frak r}_{21}(y,v)).
$$
\end{enumerate}
\end{proper}
See Lemma \ref{lem2622}.
\begin{rem} The map
$(h^{\frak r}_{21},\tilde\varphi^{\frak r}_{21},\breve\varphi^{\frak r}_{21})$
satisfying (1)-(5) above is not unique.
\end{rem}
\begin{conds}\label{exitpullbackcont2}
We consider Situation \ref{contfamipullbacksitu} and assume Condition \ref{exitpullbackcont1}.
We take $(h^{\frak r}_{21},\tilde\varphi^{\frak r}_{21},\breve\varphi^{\frak r}_{21})$ which satisfies Property \ref{proper728}.
Let $\frak S^2_{\frak r} = (W^2_{\frak r},\omega^2_{\frak r},\{\frak s^{2,\epsilon}_{\frak r}\})$.
Then for each $y \in V^1_{\frak r}$, $\xi \in W^2_{\frak r}$,
we require
\begin{equation}\label{s2haimage}
{\frak s}^{2,\epsilon}_{\frak r}(y,\xi) \in {\rm Im}(g_y)
\end{equation}
where
$g_y : E^1_{\frak r} \to E^2_{\frak r}$ is defined by
\begin{equation}\label{defofgy}
\breve\varphi^{\frak r}_{21}(y,v) = g_y(v).
\end{equation}
\end{conds}
\begin{defn}
In Situation \ref{contfamipullbacksitu} we say that $\{ (\frak V^2_{\frak r},\mathcal S^2_{\frak r})\}$ {\it can be
pulled back to $\mathcal U^1$} by $\Phi_{21}$ if and only if Conditions \ref{exitpullbackcont1} and \ref{exitpullbackcont2}
are satisfied.
\par
The {\it pullback} $\Phi_{21}^*\{ (\frak V^2_{\frak r},\mathcal S^2_{\frak r})\}$ of $\{ (\frak V^2_{\frak r},\mathcal S^2_{\frak r})\}$
is by definition $\{ (\frak V^1_{\frak r},\mathcal S^1_{\frak r}) \mid \frak r \in \frak R_0\}$
which we define below.
\begin{enumerate}
\item
$\frak R_0$ is the set of all $\frak r$
such that $U_{\frak r}^2 \cap \varphi_{21}(U_{21}) \ne \emptyset$.
\item
$\frak V^1_{\frak r}$ then is given by Condition \ref{exitpullbackcont1}.
\item
$\mathcal S^1_{\frak r} = (W^2_{\frak r},\omega^2_{\frak r},\{s^{1,\epsilon}_{\frak r}\})$,
where $\frak s^{1,\epsilon}_{\frak r}$ is defined by
\begin{equation}\label{formula718}
\frak s^{1,\epsilon}_{\frak r}(y,\xi) = g_y^{-1}(\frak s^{2,\epsilon}_{\frak r}(y,\xi)).
\end{equation}
Here $g_y$ is as in (\ref{defofgy}).
The right hand side exists because of (\ref{s2haimage}).
\end{enumerate}
\end{defn}
It is easy to see that  $\{ (\frak V^1_{\frak r},\mathcal S^1_{\frak r}) \mid
\frak r \in \frak R_0\}$ is a representative of a continuous
family perturbation of $\mathcal U^1$.
\begin{lem}\label{lem7414}
\begin{enumerate}
\item
If $\{ (\frak V^2_{\frak r},\mathcal S^2_{\frak r})\}$ can be
pulled back to $\mathcal U^1$
 by $\Phi_{21}$ and
$\{ (\frak V^{\prime 2}_{\frak r'},\mathcal S^{\prime 2}_{\frak r'})\}$
is equivalent to $\{ (\frak V^2_{\frak r},\mathcal S^2_{\frak r})\}$, then
$\{ (\frak V^{\prime 2}_{\frak r'},\mathcal S^{\prime 2}_{\frak r'})\}$ can be
pulled back to $\mathcal U^1$ by $\Phi_{21}$.
Moreover $\Phi^*_{21}\{ (\frak V^2_{\frak r},\mathcal S^2_{\frak r})\}$ is equivalent
to $\Phi^*_{21}\{ (\frak V^{\prime 2}_{\frak r'},\mathcal S^{\prime 2}_{\frak r'})\}$.
\item
The pullback $\Phi^*_{21}\{ (\frak V^2_{\frak r},\mathcal S^2_{\frak r})\}$ is
independent of the choice of $(h^{\frak r}_{21},\varphi^{\frak r}_{21},\widehat\varphi^{\frak r}_{21})$
up to equivalence.
\end{enumerate}
\end{lem}
\begin{proof}
To prove (1)
it suffices to consider only the case when $\mathcal S'_{\frak r}$
is a projection of $\mathcal S_{\frak r}$, which follows again
from the fact that $\omega_x = \vert\omega_x\vert{\rm vol}_x$ where ${\rm vol}_x$ is a volume form of $W_x$
and $\vert\omega_x\vert$ is a non-negative function, which is a part of Definition \ref{defn73ss} (4).
The assertion (2) follows from Lemma \ref{lem2715}.
\end{proof}
\begin{defnlem}\label{deflem743}
Suppose we are in Situation \ref{contfamipullbacksitu}.
\begin{enumerate}
\item
We denote by $\mathscr S^{\mathcal U^1\triangleright \mathcal U^2}(U^2)$ the set of all
elements of $\mathscr S^{\mathcal U^2}(U^2)$ whose representative
can be pulled back to $\mathcal U^1$.
This is well-defined by Lemma \ref{lem7414} (1).
\item
The pullback
$$
\Phi^*_{21} : \mathscr S^{\mathcal U^1\triangleright \mathcal U^2}(U^2) \to
\mathscr S^{\mathcal U^1}(U^1)
$$
is defined by Lemma \ref{lem7414} (2).
\item
$\Omega \mapsto \mathscr S^{\mathcal U^1\triangleright \mathcal U^2}(\Omega)$
is a subsheaf of $\mathscr S^{\mathcal U^2}$.
\item
The restriction map $\Phi^*_{21}$ is induced by a sheaf morphism:
\begin{equation}
\Phi^{*}_{21} : \varphi^{\star}_{21}\mathscr S^{\mathcal U^1\triangleright \mathcal U^2} \to \mathscr S^{\mathcal U^1}.
\end{equation}
Here the left hand side is the pullback sheaf.
(In this document, we use $\star$ to denote the pullback sheaf to distinguish
it from pullback map.)
\par
In other words, the following diagram commutes for any open sets $\Omega, \Omega'$.
\begin{equation}
\begin{CD}
\mathscr S^{\mathcal U^1\triangleright \mathcal U^2}(\Omega) @ > {\frak i_{\Omega'\Omega}} >>
\mathscr S^{\mathcal U^1\triangleright \mathcal U^2}(\Omega')  \\
@ V{\Phi^*_{21}}VV
@ VV{\Phi^*_{21}}V\\
\mathscr S^{\mathcal U^1}(U^1 \cap \Omega) @ >>{\frak i_{\Omega'\cap U^1 \Omega\cap U^1 }}>
\mathscr S^{\mathcal U^1}(U^1 \cap \Omega')
\end{CD}
\end{equation}
\end{enumerate}
\end{defnlem}
\begin{proof}
We can check the assertion directly. So we omit the proof.
\end{proof}
\begin{lem}
Let $\Phi_{i+1 i} : \mathcal U^i \to \mathcal U^{i+1}$
($i=1,2$) be
embeddings of Kuranishi charts.
We put $\Phi_{31} = \Phi_{32}\circ \Phi_{21}$.
\begin{enumerate}
\item
We have
\begin{equation}
(\Phi^{*}_{32})^{-1}(\mathscr S^{\mathcal U^1\triangleright \mathcal U^2})
\cap \mathscr S^{\mathcal U^2\triangleright \mathcal U^3}
\subseteq
\mathscr S^{\mathcal U^1\triangleright \mathcal U^3}
\end{equation}
as subsheaves of $\mathscr S^{\mathcal U^3}$.
\item
The next diagram commutes.
\begin{equation}
\begin{CD}
 \Phi^{\star}_{31}((\Phi^{*}_{32})^{-1}(\mathscr S^{\mathcal U^1\triangleright \mathcal U^2})
\cap \mathscr S^{\mathcal U^2\triangleright \mathcal U^3}) @ > {\Phi^*_{32}} >>
 \Phi^{\star}_{21}(\mathscr S^{\mathcal U^1\triangleright \mathcal U^2})  \\
@ V{}VV
@ VV{\Phi^*_{21}}V\\
 \Phi^{\star}_{31}\mathscr S^{\mathcal U^1\triangleright \mathcal U^3} @ >>{\Phi^*_{31}}>
\mathscr S^{\mathcal U^1}
\end{CD}
\end{equation}
In other words, the next diagram commutes for each $\Omega \subset U^3$.
\begin{equation}
\begin{CD}
(\Phi^{*}_{32})^{-1}\mathscr S^{\mathcal U^1\triangleright \mathcal U^2}(\Omega \cap U^2) \cap
\mathscr S^{\mathcal U^2\triangleright \mathcal U^3}(\Omega)
@ > {\Phi^*_{32}} >>
\mathscr S^{\mathcal U^1\triangleright \mathcal U^2}(\Omega \cap U^2)  \\
@ V{}VV
@ VV{\Phi^*_{21}}V\\
\mathscr S^{\mathcal U^1\triangleright \mathcal U^3}
(\Omega) @ >>{\Phi^*_{31}}>
\mathscr S^{\mathcal U^1}(\Omega \cap U^1)
\end{CD}
\end{equation}
\end{enumerate}
\end{lem}
\begin{proof}
This is a consequence of (\ref{formula718}).
\end{proof}
Next we generalize Lemma \ref{lem6767} to the case of
continuous families.
\begin{lem}\label{lem6767cont}
Suppose $\widetriangle{\mathcal U}$ and $\mathcal K$
satisfy Condition \ref{conds6.17}.
Then
there exists $c>0$, $\delta_0>0$ and $\epsilon_0 >0$ such that
for each $\epsilon < \epsilon_0$ and
member $\frak s^{\epsilon}$ of $\frak S^{\epsilon}_{\frak p}$ at $y = {\rm Exp}(x,v) \in BN_{\delta_0}(K;U_{\frak p})$,
$x \in \mathcal K_q$,
we have
\begin{equation}\label{normailityestimate4}
\vert \frak s^{\epsilon}(y) \vert \ge c \vert v\vert
\end{equation}
\end{lem}
The proof is the same as the proof of Lemma \ref{lem6767}.

\begin{rem}
The constants $\epsilon_0$, $c$, $\delta_0$ can be taken to be independent of the choice of representative
of CF-perturbations.
In fact, the notion of member is independent of the choice of representatives
of CF-perturbation.
\end{rem}

\subsubsection{CF-perturbations on good coordinate system}
\label{subsub:goodcsyscfp}

\begin{defn}\label{defn7732}
Let ${\widetriangle{\mathcal U}} = \{\mathcal U_{\frak p} \mid \frak p \in \frak P\}$ be a good coordinate system of $Z \subseteq X$.
A {\it CF-perturbation of $({\widetriangle{\mathcal U}},\mathcal K)$}
\index{CF-perturbation ! of good coordinate system}
\index{good coordinate system ! CF-perturbation of} is
by definition $\widetriangle{\frak S} = \{\frak S_{\frak p} \mid \frak p \in \frak P\}$ with the following properties.
\begin{enumerate}
\item
$\frak S_{\frak p}
\in \mathscr S^{\mathcal U_{\frak p}}(\mathcal K_{\frak p})$.
\item
If $\frak q \le \frak p$ then $\frak S_{\frak p} \in
\mathscr S^{\mathcal U^{\frak q}\triangleright \mathcal U^{\frak p}}(\mathcal K_{\frak p})$.
\item
The pull back $\Phi_{\frak p\frak q}^*(\frak S_{\frak p})$ is equivalent
to $\frak S_{\frak q}$ as an element of
$\mathscr S^{\mathcal U_{\frak q}}(\varphi_{\frak p\frak q}^{-1}(\mathcal K_{\frak p}) \cap \mathcal K_{\frak q})$.
\end{enumerate}
\end{defn}
\begin{defn}\label{smoothfunctiononvertK}
Let $\widetriangle{\frak S}= \{\frak S_{\frak p} \mid \frak p \in \frak P\}$ be a CF-perturbation
of a good coordinate system ${\widetriangle{\mathcal U}}$
and $\mathcal K$ its support system.
\begin{enumerate}
\item
We say $\widetriangle{\frak S}$ is {\it transversal to 0} if
\index{CF-perturbation ! transversal to $0$ ! on good coordinate system}
each of $\frak S_{\frak p}$ is transversal to $0$.
\item
Let $\widetriangle f : (X,Z;{\widetriangle{\mathcal U}}) \to M$ be
a strongly smooth map that is weakly submersive.
We say that $\widetriangle f$ is {\it strongly submersive with respect to $\widetriangle{\frak S}$
\index{strongly submersive (w.r.t. CF-perturbation) ! on good coordinate system}  on
$\mathcal K$} if for each $\frak p \in \frak P$ the map
$f_{\frak p}$ is strongly submersive with respect to $\frak S_{\frak p}$
on $\mathcal K_{\frak p}$ in the sense of Definition \ref{submersivepertconlocloc}.
\item
Let $g : N \to M$ be a smooth map between smooth manifolds.
We say that $\widetriangle f$ is {\it strongly transversal to $g$ with respect to $\widetriangle{\frak S}$ on
$\mathcal K$} \index{strongly transversal (w.r.t. CF-perturbation) ! to a map on good coordinate system}  if for each $\frak p \in \frak P$ the map
$f_{\frak p}$ is strongly transversal to $g$ with respect to $\frak S_{\frak p}$
on $\mathcal K_{\frak p}$.
\end{enumerate}
\end{defn}
\begin{thm}\label{existperturbcont}
Let ${\widetriangle{\mathcal U}}$ be a good coordinate system
of $Z \subseteq X$
and $\mathcal K$ its support system.
\begin{enumerate}
\item
There exists a CF-perturbation $\widetriangle{\frak  S}$
of $({\widetriangle{\mathcal U}},\mathcal K)$ transversal to $0$.
\item
If $\widetriangle f : (X,Z;{\widetriangle{\mathcal U}}) \to M$ is a
weakly submersive strongly smooth map, then
we may take $\widetriangle{\frak  S}$
with respect to which
$\widetriangle f$ is strongly submersive.
\item
If $\widetriangle f : (X,Z;{\widetriangle{\mathcal U}}) \to M$ is a
strongly smooth map which is weakly transversal to $g : N \to M$, then
we may take $\widetriangle{\frak  S}$
with respect to which $\widetriangle f$ is strongly transversal to $g$.
\end{enumerate}
\end{thm}
The proof of Theorem \ref{existperturbcont} is given in Subsection \ref{subsec:cfpgoodcsys}.

For a later use we include its relative version, Proposition \ref{prop7582752}
and Lemma \ref{lem753753}.
To state this relative version we need some digression.

\begin{defn}\label{defn735f}
Let $X$ be a separable metrizable space and $Z_1, Z_2 \subseteq X$
compact subsets. We assume $Z_1 \subset \ring{Z_2}$.
\begin{enumerate}
\item
For each $i=1,2$, let ${\widetriangle{\mathcal U^i}}
= (\frak P_i,\{\mathcal U^i_{\frak p}\},\{\Phi^i_{\frak p\frak q}\})$
be a good coordinate system of $Z_i \subseteq X$.
We say $\widetriangle{\mathcal U^2}$ {\it strictly extends}
\index{extend ! strictly extend, two good coordinate systems}
\index{good coordinate system ! strictly extend, two good coordinate systems}
$\widetriangle{\mathcal U^1}$ if the following holds.
\begin{enumerate}
\item
$
\frak P_1 = \{ \frak p \in \frak P_2 \mid
\psi_{\frak p}((s_{\frak p})^{-1}(0)) \cap Z_1 \ne \emptyset
\}.
$
The  partial order of $\frak P_1$ is the restriction of one of $\frak P_2$.
\item
If $\frak p \in \frak P_1$, then
$\mathcal U^1_{\frak p}$ is an open subchart of  $\mathcal U^2_{\frak p}$.
Moreover ${\rm Im} (\psi^1_{\frak p}) \cap Z_1 = {\rm Im} (\psi^2_{\frak p}) \cap Z_1$.
\item If $\frak p, \frak q \in \frak P_1$, then
$\Phi_{\frak p\frak q}^1$ is a restriction of $\Phi^2_{\frak p\frak q}$.
\end{enumerate}
Note the case $Z_2 = X$ is included in this definition.
\item
In the situation of (a), we say
${\widetriangle{\mathcal U}^2}$  {\it extends}
${\widetriangle{\mathcal U^1}}$
(resp. {\it weakly extends})
\index{extend ! two good coordinate systems}
\index{extend ! weakly extend, two good coordinate systems}
if it strictly extends
an open substructure
(resp. weakly open substructure) of
${\widetriangle{\mathcal U^1}}$.
\item
Let $\widehat{\mathcal U^2}
= (\{\mathcal U_p^2 \mid p \in Z_2\},\{\Phi_{pq}^2 \mid q \in {\rm Im}(\psi_{p}), \, p,q \in Z_2\})$ be a Kuranish structure
of $Z_2 \subset X$.
Then $(\{\mathcal U_p^2 \mid p \in Z_1\},\{\Phi_{pq}^2 \mid q \in {\rm Im}(\psi_{p}), \, p,q \in Z_1\})$ is
a Kuranishi structure of $Z_1 \subseteq X$. We call it the {\it restriction}
\index{Kuranishi structure ! restriction to compact subset}
of $\widehat{\mathcal U^2}$ and write $\widehat{\mathcal U^2}\vert_{Z_1}$.
\end{enumerate}
\end{defn}
\begin{lem}
For each $i=1,2$ let ${\widetriangle{\mathcal U^i}}
= (\frak P_i,\{\mathcal U^i_{\frak p}\},\{\Phi^i_{\frak p\frak q}\})$
be a good coordinate system of $Z_i \subseteq X$,
and $\widehat{\mathcal U^2}$  a Kuranishi structure of
$Z_2 \subseteq X$ such that $\widetriangle{\mathcal U^2}$
is compatible with $\widehat{\mathcal U^2}$.
Suppose that $Z_1 \subset \ring{Z_2}$ and
$\widetriangle{\mathcal U^2}$ strictly extends
$\widetriangle{\mathcal U^1}$.
\par
Then there exists a KG-embedding $\widehat{\mathcal U^2}\vert_{Z_1} \to \widetriangle{\mathcal U^1}$
with the following property.
Let $p \in Z_1$ and $p \in {\rm Im}(\psi^1_{\frak p})$.
Let $\Phi_{\frak p p}^1 : \mathcal U^2_p\vert_{U^1_{p}}
\to \mathcal U^1_{\frak p}$ be a part of the KG-embedding
$\widehat{\mathcal U^2}\vert_{Z_1} \to \widetriangle{\mathcal U^1}$.
(Here $\mathcal U^2_p\vert_{U^1_{p}}$ is an open
subchart of $\mathcal U^2_p$.)
Let $\Phi_{\frak p p}^2 : \mathcal U^2_p \to \mathcal U^2_{\frak p}$ be
a part of the KG-embedding $: \widehat{\mathcal U^2} \to \widetriangle{\mathcal U^2}$.
The next diagram commutes:
\begin{equation}\label{form726}
\begin{CD}
\mathcal U^2_p\vert_{U^1_{p}}
@ >>>
\mathcal U^2_p   \\
@ V{\Phi_{\frak p p}^1}VV
@ VV{\Phi_{\frak p p}^2}V\\
\mathcal U^1_{\frak p}
 @ >>>
\mathcal U^2_{\frak p}
\end{CD}
\end{equation}
Note that the horizontal arrows are embeddings as open substructures.
\par
Moreover, the KG embedding
$\widehat{\mathcal U^2}\vert_{Z_1} \to \widetriangle{\mathcal U^1}$
satisfying this property is unique up to equivalence in the sense of Definition \ref{embedkuraequiv}.
\par
The same conclusion holds for extension or weak extension.
\end{lem}
\begin{proof}
Let $p \in Z_1$ and $p \in {\rm Im}(\psi^1_{\frak p})$.
By Definition \ref{defn735f} (1)(b),
$\mathcal U^1_{\frak p}$ is an open subchart of $\mathcal U^2_{\frak p}$.
In particular, $U^1_{\frak p}$ is an open subset of $U^2_{\frak p}$.
Let $\Phi_p^2 : \mathcal U^2_p \to \mathcal U^2_{\frak p}$ be
a part of the KG embedding$: \widehat{\mathcal U^2} \to \widetriangle{\mathcal U^2}$.
We put
$$
U^1_{p} = (\varphi^2_p)^{-1}(U^1_{\frak p}) \subseteq U^2_{p}.
$$
By Lemma \ref{lem321321}, there exists a Kuranishi structure of $Z_1 \subseteq X$
whose Kuranishi chart is $\mathcal U^2_p\vert_{U^1_{p}}$.
We can define $\Phi_p^1$ by restricting $\Phi_p^2$.
The commutativity of Diagram (\ref{form726}) is obvious.
\end{proof}
We recall Definition \ref{gcsystemcompa} for the compatibility used in the next proposition.
\begin{prop}\label{prop7582752}
Let $X$ be a separable metrizable space, $Z_1, Z_2 \subseteq X$  be
compact subsets. We assume $Z_1 \subset \ring{Z_2}$.
Let $\widehat{\mathcal U^2}$ be a Kuranishi structure of $Z_2 \subseteq X$
and $\widetriangle{\mathcal U^1}$ a good coordinate system of $Z_1 \subseteq X$
which is compatible with $\widehat{\mathcal U^2}\vert_{Z_1}$.
Then there exists a good coordinate system $\widetriangle{\mathcal U^2}$ of $Z_2 \subseteq  X$
such that
\begin{enumerate}
\item
$\widetriangle{\mathcal U^2}$ extends $\widetriangle{\mathcal U^1}$.
\item
$\widetriangle{\mathcal U^2}$ is compatible with $\widehat{\mathcal U^2}$.
\item
Diagram (\ref{form726}) commutes.
\end{enumerate}
\end{prop}
The proof is given in Subsection \ref{subsec:moreversionegcs2}.
\begin{lem}\label{lem753753}
Suppose we are in the situation of Proposition \ref{prop7582752}.
We may choose $\widetriangle{\mathcal U^2}$ such that the following
holds.
\par
Let $\widetriangle{\mathcal U^1_0}$  be an open substructure of
$\widetriangle{\mathcal U^1}$ strictly extended to $\widetriangle{\mathcal U^2}$.
Let $\widetriangle {f^1} = \{f^1_{\frak p}\} : (X,Z_1;\widetriangle{\mathcal U^1_0}) \to M$
and $\widehat {f^2} = \{f^2_{p}\} : (X,Z_2;\widehat{\mathcal U^2}) \to M$
be strongly continuous maps.
Assume that the equality
$
f^1_{\frak p} \circ \varphi^1_{\frak p p} = f^2_{p}
$
holds on $U^1_{p} \subset U^2_p$  for each $p \in {\rm Im}\,(\psi_{\frak p}) \cap Z_1$. Then
the following hold:
\begin{enumerate}
\item
There exists a strongly continuous map
$\widetriangle {f^2} = \{f^2_{\frak p}\} : (X,Z_2;\widetriangle{\mathcal U^2}) \to M$
such that:
\begin{enumerate}
\item
$
f^2_{\frak p} \circ \varphi^2_{\frak p p} = f^2_{p}
$
for $p \in {\rm Im}\,(\psi_{\frak p}) \cap Z_2$.
\item
$f^2_{\frak p} = f^1_{0,\frak p}$ on $U^1_{\frak p} \subset U^2_{\frak p}$.
\end{enumerate}
\begin{equation}
\xymatrix{
&  &  & M \\
\widetriangle {\mathcal U^1_0} \ar[urrr]^{\widetriangle f^1}
\ar[r]  & \widetriangle {\mathcal U^2} \ar@{.>}[urr]|{\widetriangle f^2} && \\
\widehat {\mathcal U^2}\vert_{Z_1}  \ar[u]\ar[r] &\widehat {\mathcal U^2}
\ar[u]^{\{f^2_{\frak p p}\}} \ar[uurr]_{\widehat f^2} &&
}
\end{equation}
\item
If $\widetriangle {f^1}$ and $\widehat {f^2}$ are strongly smooth (resp. weakly submersive)
then so is $\widetriangle {f^2}$.
\item
If $\widetriangle {f^1}$ and $\widehat {f^2}$ are strongly transversal to $g : N \to M$
then so is $\widetriangle {f^2}$.
\end{enumerate}
\end{lem}
The proof is given in Subsection \ref{subsec:moreversionegcs2}.
\begin{defn}\label{defn738}
Let ${\widetriangle{\mathcal U^i}}
= (\frak P_i,\{\mathcal U^i_{\frak p}\},\{\Phi^i_{\frak p\frak q}\})$
be good coordinate systems of $Z_i \subseteq X$, $i=1,2$
and $\mathcal K^i$ support systems of ${\widetriangle{\mathcal U^i}}$ for $i=1,2$.
\begin{enumerate}
\item
Suppose $\widetriangle{\mathcal U^2}$ strictly extends
$\widetriangle{\mathcal U^1}$.
We say that $\mathcal K^1$ is {\it compatible with} $\mathcal K^2$ if
for each $\frak p \in \frak P_1$ we have
$$
\mathcal K^1_{\frak p} \subset \ring{\mathcal K^2_{\frak p}}.
$$
(We note that $U^1_{\frak p} \subset U^2_{\frak p}$ by Definition \ref{defn735f} (1)(c).)
In this situation we say $(\widetriangle{\mathcal U^2},\mathcal K^2)$ {\it strictly extends}
\index{extend ! strictly extend, pairs of good coordinate systems and support systems}
\index{good coordinate system ! strictly extend, pairs of good coordinate systems and support systems}
$(\widetriangle{\mathcal U^1},\mathcal K^1)$.
\item
Suppose $\widetriangle{\mathcal U^2}$  extends
$\widetriangle{\mathcal U^1}$.
We say that $\mathcal K^1$ is {\it compatible with} $\mathcal K^2$ if
the following holds.
\par
By definition there exists an open substructure $\widetriangle{\mathcal U^1_0} = \{\mathcal U^1_{\frak p,0}\}$
of $\widetriangle{\mathcal U^1}$ such that $\widetriangle{\mathcal U^2}$ strictly extends
$\widetriangle{\mathcal U^1_0}$.
For each $\frak p \in \frak P_1$ we require
$$
\mathcal K^1_{\frak p} \subset U^1_{\frak p,0} \cap \ring{\mathcal K^2_{\frak p}}.
$$
In this situation we say $(\widetriangle{\mathcal U^2},\mathcal K^2)$ {\it extends}
\index{extend ! pairs of good coordinate systems and support systems}
\index{good coordinate system ! extend, pairs of good coordinate systems and support systems}
$(\widetriangle{\mathcal U^1},\mathcal K^1)$.
\end{enumerate}
\end{defn}
\begin{defn}\label{defn738pert}
For each $i=1,2$ let ${\widetriangle{\mathcal U^i}}
= (\frak P_i,\{\mathcal U^i_{\frak p}\},\{\Phi^i_{\frak p\frak q}\})$
be a good coordinate system of $Z_i \subseteq X$,
$\mathcal K^i$ a support system of ${\widetriangle{\mathcal U^i}}$, and
$\widetriangle{\frak S^i}$ a CF-perturbation
of $({\widetriangle{\mathcal U^i}},\mathcal K^i)$.
\begin{enumerate}
\item
Suppose $(\widetriangle{\mathcal U^2},\mathcal K^2)$ strictly extends
$(\widetriangle{\mathcal U^1},\mathcal K^1)$.
We say $\widetriangle{\frak S^2}$ {\it strictly extends}
\index{CF-perturbation ! strictly extend}
\index{extend ! strictly extend, CF-perturbation} $\widetriangle{\frak S^1}$
if the restriction of ${\frak S^2_{\frak p}}$ to $\mathcal K^1_{\frak p}$ is ${\frak S^1_{\frak p}}$
for each $\frak p \in \frak P_1$.
\item
Suppose $(\widetriangle{\mathcal U^2},\mathcal K^2)$ extends
$(\widetriangle{\mathcal U^1},\mathcal K^1)$.
We say $\widetriangle{\frak S^2}$ {\it extends}
\index{extend ! CF-perturbation}
\index{CF-perturbation ! extend} $\widetriangle{\frak S^1}$
if the restriction of ${\frak S^2_{\frak p}}$ to $\mathcal K^1_{\frak p}$ is $\frak S^1_{\frak p}$
for each $\frak p \in \frak P_1$.
\end{enumerate}
\end{defn}
\begin{rem}
In Definition \ref{defn738pert} (2)
we note that $\mathcal K_{\frak p}^1$ can be
regarded as a support system of an open substructure of $\widetriangle{\mathcal U^1}$ by Definition \ref{defn738} (2).
\end{rem}
\begin{prop}\label{existperturbcontrel}
Suppose we are in the situation of Proposition \ref{prop7582752}.
We may choose $\widetriangle{\mathcal U^2}$ such that the following
holds.
\par
Let $\widetriangle{\mathcal U^1_0}$ be an open substructure of
$\widetriangle{\mathcal U^1}$  strictly extended to $\widetriangle{\mathcal U^2}$.
Let $\mathcal K^1$, $\mathcal K^2$ be support systems of ${\widetriangle{\mathcal U^1}}$, ${\widetriangle{\mathcal U^2}}$ respectively such that
$(\widetriangle{\mathcal U^2},\mathcal K^2)$
extends $(\widetriangle{\mathcal U^1},\mathcal K^1)$.
Let
$\widetriangle{\frak S^1}$ be a CF-perturbation
of $({\widetriangle{\mathcal U^1}},\mathcal K^1)$.
Then there exists a CF-perturbation $\widetriangle{\frak S^2}$
of $({\widetriangle{\mathcal U^2}},\mathcal K^2)$ which
extends $\widetriangle{\frak S^1}$.
Moreover the following holds.
\begin{enumerate}
\item
If $\widetriangle{\frak S^1}$ is transversal to $0$,
so is $\widetriangle{\frak S^2}$.
\item
Suppose we are in the situation of Lemma \ref{lem753753} (1) in addition.
We assume $\widetriangle {f^1}$ is strongly submersive
with respect to $\widetriangle{\frak S^1}$
and $\widetriangle {f^2}$ is weakly submersive.
Then $\widetriangle {f^2}$ is strongly submersive
with respect to $\widetriangle{\frak S^2}$.
\item
Suppose we are in the situation of Lemma \ref{lem753753} (1) in addition.
We assume $\widetriangle {f^1}$ is strongly transversal to $g : N \to M$
with respect to $\widetriangle{\frak S^1}$ and $\widetriangle {f^2}$ is
weakly transversal to $g$.
Then $\widetriangle {f^2}$ is strongly transversal to $g : N \to M$ with
respect to $\widetriangle{\frak S^2}$.
\end{enumerate}
\end{prop}
The proof of Proposition \ref{existperturbcontrel} is given in Subsection \ref{subsec:cfpgoodcsys}.
\par
Before we end this subsection, we consider a family of CF-perturbations.
The next definition is a CF-perturbation
version of Definition \ref{uniformmulivalupert}.
\begin{defn}\label{uniformcongpert}
A $\sigma$-parameterized family of CF-perturbations
$\{ \widetriangle{\frak S_{\sigma}} \mid \sigma \in \mathscr A\}$
($\widetriangle{\frak S_{\sigma}}
= \{\mathcal S^{\epsilon}_{\sigma,\frak p}\}$)
of $(\widetriangle{\mathcal U},\mathcal K)$
is called a {\it uniform family}
\index{uniform family ! of CF-perturbations}
\index{CF-perturbation ! uniform family} if the convergence in Definition \ref{defn73ss} (3)
is uniform.
More precisely, we require the following.
\par
For each $\frak o>0$ there exists $\epsilon_0(\frak o)>0$
such that if $0< \epsilon < \epsilon_0(\frak o)$, then
\begin{equation}
\vert \frak s(y) - s_p(y) \vert < \frak o,
\qquad
\vert (D\frak s)(y) - (Ds_p)(y) \vert < \frak o
\end{equation}
hold for any $\frak  s$ which is a member of $\mathcal S^{\epsilon}_{\sigma,\frak p}$
at any point $y \in \mathcal K_{\frak p}$
for any $\frak p \in \frak P$ and $\sigma \in \mathscr A$.
\end{defn}
\begin{lem}\label{lemma748}
If $\{ \widetriangle{\frak S_{\sigma}} \mid \sigma \in \mathscr A\}$
is a uniform family of CF-perturbations of
$(\widetriangle{\mathcal U},\mathcal K)$, then the
constant $\epsilon_0$ in Lemma \ref{lem6767cont}
can be taken independent of $\sigma$.
\end{lem}
\begin{proof}
In the same way as Lemma \ref{lem618} follows from
the proof of Lemma \ref{lem611},
this follows from the proof of Lemma \ref{lem6767cont},
which is the same as the proof of Lemma \ref{lem6767}.
\end{proof}

\subsection{Partition of unity associated to a good coordinate system}
\label{subsec:partitionunigcs}

We next define the notion of partition of unity on spaces with good coordinate system.
\begin{shitu}\label{situ5757}
Let $\widetriangle{\mathcal U}=(\frak P, \{\mathcal U_{\frak p} \},\{ \Phi_{\frak p \frak q}\})$ be a good coordinate
system of $Z \subseteq X$ and $\mathcal K$  its
support system.$\blacksquare$
\end{shitu}
\begin{defn}\label{def758}
In Situation \ref{situ5757}, let $\Omega$ be an open subset of $\vert{\mathcal K}\vert$
and $f : \Omega \to \R$  a continuous function on it.
We say $f$ is {\it strongly smooth}
\index{strongly smooth function}
\index{map from Kuranishi structure ! strongly smooth function} if the restriction of
$f$ to $\mathcal K_{\frak p} \cap \Omega$ is smooth for any
$\frak p \in \frak P$.
\end{defn}
We remark that the strongly smooth function  $f : \Omega \to \R$ is nothing but a
strongly smooth map to $\R$ regarded as a manifold in the sense of
Definition \ref{definition32727}.
\par
We take a support system $\mathcal K^+$ such that $(\mathcal K ,\mathcal K^+)$ is
a support pair.
We take a metric $d$ on $\vert{\mathcal K}^+ \vert$ and use it in the next definition.
\begin{defn}\label{def761}
In Situation \ref{situ5757},
let $\delta > 0$ be a  positive number.
We put
\begin{equation}\label{formula7120}
\mathcal K_{\frak p}(2\delta)
=
\{
x \in U_{\frak p} \mid d(x,\mathcal K_{\frak p})
\le 2\delta
\}.
\end{equation}
$2\delta$-neighborhood of $\mathcal K_{\frak p}$.
We assume $\mathcal K_{\frak p}(2\delta)$ is compact.

We put $\mathcal K(2\delta) = \{\mathcal K_{\frak p}(2\delta)\}_{\frak p \in \frak P}$,
which is a support system of  $\widetriangle{\mathcal U}$.
We put
\begin{equation}\label{formula712}
\Omega_{\frak p}(\mathcal K,\delta)
=
B_{\delta}(\mathcal K_{\frak p})
=
\{x \in \vert \mathcal K(2\delta)\vert
\mid d(x,\mathcal K_{\frak p}) < \delta\}.
\end{equation}
\end{defn}
\par
\centerline{
\epsfbox{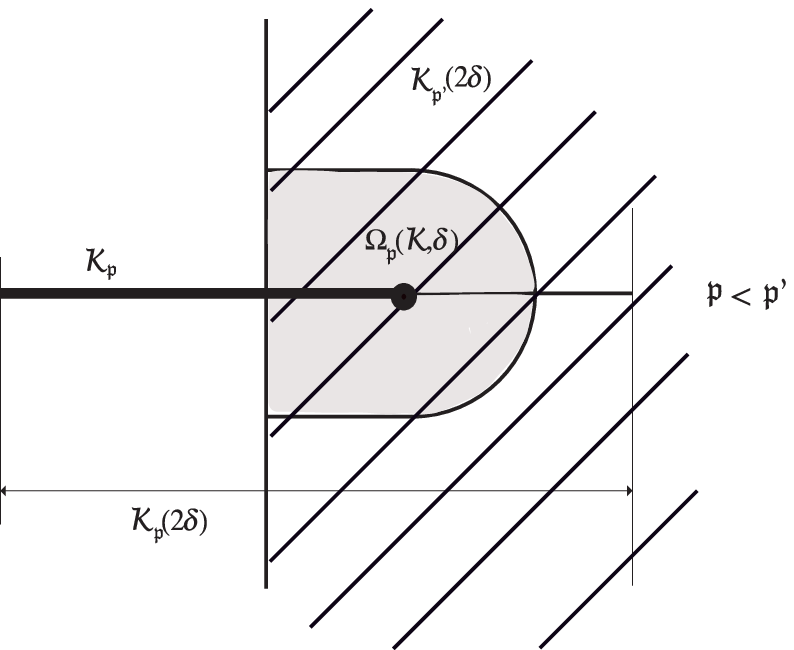}}
\par
\centerline{\bf Figure 7.1}

\begin{lem}\label{rem74444}
If $\frak q < \frak p$ and
$\delta>0$ is sufficiently small, then
$$
\Omega_{\frak p}(\mathcal K,\delta) \cap \mathcal K_{\frak q}(2\delta)
\subset \mathcal K_{\frak p}(2\delta).
$$
Moreover,
$\Omega_{\frak p}(\mathcal K,\delta) \cap \mathcal K_{\frak p}(2\delta)
\subset \rm{Int}\,\mathcal K_{\frak p}(2\delta)$ and
$\Omega_{\frak p}(\mathcal K,\delta) \cap \mathcal K_{\frak p}(2\delta)$
is an orbifold.
\end{lem}
\begin{proof}
The first claim is a consequence of
(\ref{formula712}).
The second claim follows from the fact that ${\rm Int}\,\mathcal K_{\frak p}(2\delta)$ is open
in $\bigcup_{\frak q \le \frak p} \mathcal K_{\frak q}(2\delta)$.
\end{proof}
\begin{defn}\label{pounity}
We say $\{\chi_{\frak p}\}$ is a {\it strongly smooth
partition of unity }
of the quintuple $(X,Z,{\widetriangle{\mathcal U}},\mathcal K,\delta)$
if the following holds.
\begin{enumerate}
\item
$\chi_{\frak p} : \vert \mathcal K(2\delta)\vert \to [0,1]$ is a strongly
smooth function.
\item
${\rm supp} \, \chi_{\frak p} \subseteq \Omega_{\frak p}(\mathcal K,\delta)$.
\item
There exists an open neighborhood $\frak U$ of $Z$ in $\vert \mathcal K(2\delta)\vert$
such that for each point $x \in \frak U$, we have
\begin{equation}\label{pounitymainequ}
\sum_{\frak p} \chi_{\frak p}(x) = 1.
\end{equation}
\end{enumerate}
\end{defn}

\begin{rem}
In our earlier writings such as \cite[Section 12]{fooo09},
we defined a partition of unity in a slightly different way. Namely
we require $\chi_{\frak p}$ to be defined on $U_{\frak p}$ and
 we required
\begin{equation}\label{pounitymainequ2}
\chi_{\frak p}(x)
+
\sum_{\frak q > \frak p, x \in U_{\frak q\frak p}}\chi_{\frak q}
(\varphi_{\frak q\frak p}(x))
+
\sum_{\frak p > \frak q, x \in {\rm Im}\varphi_{\frak p\frak q}}
\chi^{\delta}(\rho(x;U_{\frak p\frak q}))
\chi_{\frak q}(\pi(x))= 1
\end{equation}
instead of  (\ref{pounitymainequ}).
Here $\pi$ is the projection of a tubular neighborhood
of $\varphi_{\frak p\frak q}(U_{\frak p\frak q})$
in $U_{\frak p}$, $\rho$ is a tubular distance function
of this tubular neighborhood
\footnote{The tubular distance function is, roughly speaking, the distance from
the image $\varphi_{\frak p\frak q}(U_{\frak p\frak q})$.
We do not give the precise definition here since we do not use this notion.
See \cite{Math73}.} and $\chi^{\delta} : [0,\infty) \to [0,1]$
is a smooth function such that it is 1 in a neighborhood of $0$ and is $0$ on
$(\delta,\infty)$.
We required (\ref{pounitymainequ2}) for all $x \in \mathcal K^1_{\frak p}$.
Formula (\ref{pounitymainequ2}) depends on the choice of
the tubular neighborhood and we need to take
certain compatible system of tubular neighborhoods
(see \cite{Math73})
to define it.
\par
Note the covering
$
\vert {\widetriangle{\mathcal U}} \vert
= \bigcup_{\frak p \in \frak P} U_{\frak p}
$
is {\it not} an open covering. In fact $U_{\frak p}$
is not an open subset of $\vert {\widetriangle{\mathcal U}} \vert$
unless $\frak p$ is maximal. So a partition of unity in the above sense
is different from one in the usual sense.
Definition \ref{pounity} seems to be simpler than the one
in the previous literature.
Also it is more elementary in the sense that we do not use
any compatible system of tubular neighborhoods.
However, by using the third term of (\ref{pounitymainequ2})
we can extend the function $\chi_{\frak q}$
to its neighborhood in $\vert\mathcal K\vert$ and
the compatibility of tubular neighborhoods implies that
it becomes a strongly smooth function.
So the present definition is not very different from the earlier one.
\end{rem}
In the rest of this subsection we prove the existence of strongly smooth
partition of unity.
The proof is very  similar to the corresponding one in manifold theory.
However we prove it for completeness' sake.
We begin with the following lemma.
\begin{lem}\label{bumpfunctionlemma}
For any open set $W$ of $\vert \mathcal K(2\delta)\vert$ containing
a compact subset $K$ of $Z$ there exists a strongly smooth function $g : \vert \mathcal K(2\delta)\vert
\to \R$ that has a
compact support in $W$ and is $1$ on a neighborhood of $K$.
\end{lem}
\begin{proof}
Let $\mathcal K^+$ be a support system such that
$\mathcal K(2\delta) < \mathcal K^+$.
We take a metric $d$ on $\vert \mathcal K^+\vert$.
(See \cite[Proposition 5.1]{foooshrink}.)
For each $x \in K$ we consider
$$
\epsilon_x = \inf\{ d(x,\mathcal K_\frak q(2\delta)) \mid
\frak q \in \frak P, \,\,
x \notin \mathcal K_\frak q(2\delta)\}.
$$
Let $\frak p_x \in \frak P$ be the element
which is maximum among elements
$\frak p \in \frak P$ such that
$x \in \mathcal K_\frak p(2\delta)$.
Let $W^+_x$ be an open neighborhood of $x$ in
$
\mathcal K_{\frak p_x}^+ \cap B_{\epsilon_x/2}(x,\vert\mathcal K^+\vert).
$
We may choose it so small  that
if $x \in \mathcal K_{\frak q}(2\delta)$
then $W^+_x \cap \mathcal K_{\frak q}(2\delta)$ is open.
(Here we use the fact that
$x \in \mathcal K_{\frak q}(2\delta)$ implies $\frak q \le \frak p_x$.)
Therefore  $W_x = W^+_x \cap \ring{\mathcal K}_{\frak p_x}(2\delta)$
is open in $\vert\mathcal K(2\delta)\vert$.
We may choose $W^+_x$ small so that $W_x \subset W$.
\par
We can take a smooth function $f_{x} : W^+_x \to [0,1]$ that has a compact support
and is 1 in a neighborhood $Q^+_x$ of $x$.
(This is because $W^+_x$ is an orbifold.)
We restrict $f_x$ to $W_x$ and extend it by $0$ to $\vert\mathcal K\vert$,
which we denote by the same symbol.
Then $f_x$ is a strongly smooth function with a compact support in $W_x$ and equal to $1$ on an open neighborhood
$Q_x= Q^+_x \cap \vert\mathcal K(2\delta)\vert$ of $x$.
\par
We find finitely many points $x_1,\dots,x_N$ of $K$ so that
$$
\bigcup_{i=1}^N Q_{x_i} \supset K.
$$
We put
$$
f = \sum_{i=1}^N f_{x_i}.
$$
Then the function $f$ is a strongly smooth function on $\vert\mathcal K(2\delta)\vert$
with a compact support in $W$ and satisfies
$
f(x) \ge 1
$
if $x \in Q = \bigcup_{i=1}^N Q_{x_i}$. Let $\rho : \R \to \R$ be a nondecreasing smooth map
such that $\rho(s) = 0$ if $s$ is in a neighborhood of $0$
and $\rho(s) = 1$ if $s \ge 1$.
It is easy to see that $g = \rho\circ f$ has the required property.
\end{proof}
\begin{prop}\label{pounitexi}\index{orbifold ! partition of unity}
If $\delta>0$ is sufficiently small, then
a strongly smooth
partition of unity
of $(X,Z,{\widetriangle{\mathcal U}},\mathcal K,\delta)$
exists.
\end{prop}
\begin{proof}
We put
$$
\mathcal K_{\frak p}(-\delta)
=
\{x \in \mathcal K_{\frak p} \mid d(x,U_{\frak p} \setminus \mathcal K_{\frak p}) \ge \delta\}.
$$
It is easy to see that
$$
\bigcup_{\frak p} \mathcal K_{\frak p}(-\delta) \supset Z.
$$
for sufficiently small $\delta>0$.
We apply Lemma \ref{bumpfunctionlemma} to
$(K,W) = ({\mathcal K_{\frak p}(-\delta)},\Omega_{\frak p}(\mathcal K,\delta))$
to obtain $f_{\frak p}$.
Then there exists a neighborhood $W'$ of $Z$ such that
$$
\sum_{\frak p} f_{\frak p}(x) \ge 1/2
$$
for $x\in W'$.
We apply Lemma \ref{bumpfunctionlemma} to
$(K,W) = (Z,W')$ to obtain $f : \vert\mathcal K(2\delta)\vert \to [0,1]$.
Now we define
\begin{equation}
\chi_{\frak p}(x) =
\begin{cases}
\displaystyle
\frac{f(x)f_{\frak p}(x)}{\sum_{\frak p} f_{\frak p}(x)}
&\text{if $x \in W'$}\\
 0
&\text{if $x \notin W'$}.
\end{cases}
\end{equation}
Then it is easy to see that this is a strongly smooth partition of unity.
\end{proof}

\subsection{Differential form on good coordinate system and Kuranishi structure}
\label{subsec:differentialforms}

In this subsection we define differential forms of good coordinate system and
Kuranishi structure, and
give several basic definitions for differential forms.

\begin{defn}\index{differential form ! of good coordinate system}
\label{defndiffformgcs}
Let ${\widetriangle{\mathcal U}}=(\frak P, \{\mathcal U_{\frak p}\}, \{ \Psi_{\frak p \frak q}\})$ be a good coordinate system
and $\mathcal K$  its support system.
A {\it smooth differential $k$ form $\widetriangle h$ of
$({\widetriangle{\mathcal U}},\mathcal K)$}
by definition assigns a smooth $k$ form $h_{\frak p}$ on a
neighborhood of $\mathcal K_{\frak p}$ for each $\frak p \in \frak P$ such that
the next equality holds on $\varphi_{\frak p\frak q}^{-1}(\mathcal K_{\frak p})\cap
\mathcal K_{\frak q}$.
\begin{equation}\label{compatidefform}
\varphi_{\frak p\frak q}^* h_{\frak p} = h_{\frak q}.
\end{equation}
A {\it smooth differential $k$ form $\widetriangle h$ of ${\widetriangle{\mathcal U}}$}
by definition assigns a smooth differential $k$ form $h_{\frak p}$ on $U_{\frak p}$
for each $\frak p$ such that (\ref{compatidefform})
is satisfied on $U_{\frak p\frak q}$.
\end{defn}
\begin{defn}
A {\it differential $k$ form $\widehat h$ of a Kuranishi structure}
\index{differential form ! of Kuranishi structure}
$\widehat{\mathcal U}$ of $Z\subseteq X$ assigns a differential $k$-form
$h_p$ on $U_p$ for each $p \in Z$ such that
$\varphi_{pq}^*h_p = h_q$.
\end{defn}
\begin{defn}\label{defn75555}
\begin{enumerate}
\item
Let $\widetriangle f = \{f_{\frak p}\} :
(X,Z;{\widetriangle{\mathcal U}}) \to M$ be a strongly smooth
map  and $h$ a smooth differential $k$ form
on $M_s$.
Then
$\widetriangle f^*h =  \{f_{\frak p}^*h\}$ is a smooth differential $k$ form on
${\widetriangle{\mathcal U}}$,
which we call  the {\it pullback of $h$ by $\widetriangle f = \{f_{\frak p}\}$}
\index{differential form ! pullback by a strongly smooth map}
and denote by $\widetriangle f^* h$.
\item
A smooth differential $0$ form is nothing but a strongly smooth function
\index{strongly smooth function}
\index{differential form ! strongly smooth function}
in the sense of Definition \ref{def758}.
\item
If $\widetriangle {h^i} = \{h^i_{\frak p}\}$ are smooth differential $k_i$ forms on ${\widetriangle{\mathcal U}}$
for $i = 1,2$, then $\{h^1_{\frak p} \wedge h^2_{\frak p}\}$ is a
smooth differential $k_1 + k_2$ form on ${\widetriangle{\mathcal U}}$.
We call it the {\it wedge product}
\index{differential form ! wedge product} and denote it by $\widetriangle {h^1} \wedge \widetriangle {h^2}$.
\item
In particular, we can define a product of a smooth differential form and a strongly smooth function.
\item
The sum of smooth differential forms of the same degree is defined
by taking the sum for each $\frak p \in \frak P$.
\item
If $\widetriangle h = \{h_{\frak p}\}$ is a smooth differential $k$ form on ${\widetriangle{\mathcal U}}$, then
$\{dh_{\frak p}\}$ is a smooth differential
$k+1$ form on ${\widetriangle{\mathcal U}}$,
which we call the {\it exterior derivative} of $\widetriangle h$
\index{differential form !exterior derivative} and denote by $d\widetriangle h$.
\item The {\it support} ${\rm Supp}(\widetriangle h)$ of $\widetriangle h$ is the union of \index{differential form ! support}
the supports of $h_{\frak p}$, $\frak p \in \frak P$, which is a subset of $\vert \mathcal K\vert$.
\end{enumerate}
\par
(1)-(6) have  obvious versions in the case of differential forms of a
Kuranishi structures. (7) is modified as follows.
\par
\begin{enumerate}
\item[(7)']
If $\widehat h = \{h_{p} \mid p \in Z\}$ is a smooth differential $k$ form on
a Kuranishi structure ${\widehat{\mathcal U}}$ of $Z \subseteq X$,
its {\it support} ${\rm Supp}(\widehat h)$ is the set of the points $p \in Z$ such that
$\widehat h_p$ is nonzero on any neighborhood of $o_p$ in $U_p$.
\end{enumerate}
Note that ${\rm Supp}(\widehat h)$ is a subset of $Z$ in this case.
\end{defn}
\begin{defn}
Let ${\widetriangle{\mathcal U}}$ be a good coordinate system
of $Z \subseteq X$ and $\widetriangle h = \{h_{\frak p}\}$
a differential form on it.
We say that $\widetriangle h$ has a {\it compact support in
$\ring{Z}$} if
\begin{equation}
{\rm Supp}(\widetriangle h) \cap X
\subset \ring Z.
\end{equation}
Here the intersection in the left hand side is taken on
$\vert{\widetriangle{\mathcal U}}\vert$.
\end{defn}

\subsection{Integration along the fiber (pushout) on a
good coordinate system}
\label{subsec:integrationgcs}

To define the pushout (integration along the fiber) of a differential form
using a CF-perturbation,
we need a CF-perturbation version of
Propositions \ref{splem2}, \ref{lem715}.
To state them
we introduce the notation of {\it support set} of a CF-perturbation.
\begin{defn}\label{defn767}
\begin{enumerate}
\item
Let $\mathcal U$ be a Kuranishi chart,
$\mathfrak S$ a CF-perturbation of $\mathcal U$
and $\{(\frak V_{\frak r},\mathcal S_{\frak r}) \mid \frak r\in
\frak R\}$ its representative.
For each $\epsilon >0$ we define the {\it support set}
$\Pi((\mathfrak S^{\epsilon})^{-1}(0)) $ of $\frak S$
as the set of all $x \in U$ with the following property:
\par
There exist $\frak r \in \frak R$ and $y \in V_{\frak r}$, $\xi \in W_{\frak r}$
such that
$$
\phi_{\frak r}([y]) = x, \qquad
s^{\epsilon}_{\frak r}(y,\xi) = 0,
\qquad \xi \in {\rm Supp}(\omega_x).
$$
This definition is independent of the choice of representative
because of Definition \ref{defn73ss} (4).
\item
Let ${\widetriangle{\mathcal U}}$ be a good coordinate system,
$\mathcal K$ its support system and $\widetriangle{\mathfrak S}
= \{\mathfrak S_{\frak p}\}$ a CF-perturbation
of $({\widetriangle{\mathcal U}},\mathcal K)$.
The {\it support set $\Pi((\widetriangle{\mathfrak S^{\epsilon}})^{-1}(0))$ of
$\widetriangle{\mathfrak S}$}
\index{CF-perturbation ! support set}
is defined by
$$
\Pi((\widetriangle {{\mathfrak S}^{\epsilon}})^{-1}(0))
= \bigcup_{\frak p \in \frak P}
\left(\mathcal K_{\frak p} \cap \Pi(({\mathfrak S}_{\frak p}^{\epsilon})^{-1}(0))\right)
$$
which is a subset of $\vert \mathcal K\vert$.
\end{enumerate}
\end{defn}
The CF-perturbation version of Propositions \ref{splem2}
is as follows
\begin{lem}\label{lem739}
Let ${\widetriangle{\mathcal U}}$ a good coordinate system of $Z \subseteq X$,
$\mathcal K$ its support system and $\widetriangle{\mathfrak S}
= \{\mathfrak S_{\frak p} \mid \frak p \in \frak P\}$ a CF-perturbation
of $({\widetriangle{\mathcal U}},\mathcal K)$.
Let $\mathcal K^- < \mathcal K^+ < \mathcal K' < \mathcal K$.
Then there exist positive numbers $\delta_0$  and $\epsilon_0$
such that
$$
B_{\delta}(\mathcal K_{\frak q}^{-} \cap Z) \cap
\bigcup_{\frak p \in \frak P}
(\mathcal K'_{\frak p} \cap \Pi((\widetriangle{{\mathfrak S}_{\frak p}^{\epsilon}})^{-1}(0)))
\subset
\mathcal K_{\frak q}^{+}
$$
 for $\delta < \delta_0, 0 < \epsilon < \epsilon_0$.
\end{lem}
\begin{proof}
Using Lemma \ref{lem6767cont}, the proof is the same as the proof of Proposition  \ref{splem2}.
\end{proof}
The next lemma is the CF-perturbation version of
Proposition \ref{lem715}.
\begin{lem}\label{lem740}
Let $\mathcal K^1,\mathcal K^2, \mathcal K^3$ be support systems of
a good coordinate system ${\widetriangle{\mathcal U}}$  of $Z \subseteq X$
with $\mathcal K^1 < \mathcal K^2 < \mathcal K^3$
and $\widetriangle{\mathfrak S}$
a CF-perturbation of $({\widetriangle{\mathcal U}},\mathcal K^3)$.
Then
there exists a neighborhood $\frak U(Z)$ of $Z$
in $\vert \mathcal K^2\vert$ such that  for $0 < \epsilon < \epsilon_0$
\begin{equation}
 \frak U(Z) \cap
\bigcup_{\frak p}(\Pi((\widetriangle{{\mathfrak S}_{\frak p}^{\epsilon}})^{-1}(0)
)\cap \mathcal K^1_{\frak p})
=
\frak U(Z)\cap
\bigcup_{\frak p}(\Pi((\widetriangle{{\mathfrak S}_{\frak p}^{\epsilon}})^{-1}(0)
)\cap \mathcal K^2_{\frak p}).
\end{equation}
\end{lem}
\begin{proof}
The proof is the same as the proof of Proposition \ref{lem715}.
\end{proof}
The next lemma is the CF-perturbation version of
Corollary \ref{cor69}.

\begin{lem}\label{lem743}
In the situation of Lemma \ref{lem740},
$\left(\bigcup_{\frak p}((\frak s_{\frak p}^{\epsilon})^{-1}(0))
\cap \ring{\mathcal K}^1_{\frak p})\right)
\cap \frak U(Z)$ is compact for a sufficiently small neighborhood $\frak U(Z)$
of $Z$ in $\vert\mathcal K^2\vert$.
Moreover, we have
$$
\lim_{\epsilon\to 0}
\left(\bigcup_{\frak p}(\Pi((\widetriangle{{\mathfrak S}_{\frak p}^{\epsilon}})^{-1}(0)
\cap \ring{\mathcal K}^1_{\frak p})\right)
\cap \frak U(Z)
\subseteq X
$$
in Hausdorff topology.
\end{lem}
\begin{proof}
The proof is the same as Corollary \ref{cor69}.
\end{proof}
Now to define the pushout of a differential form we consider the following situation.
\begin{shitu}\label{situ774}
Let
$\widetriangle{\mathcal U}=(\frak P, \{\mathcal U_{\frak p} \},\{ \Phi_{\frak p \frak q}\})$ be a good coordinate system,
$\mathcal K$ its support system,
$\widetriangle h = \{h_{\frak p}\}$ a differential form on  ${\widetriangle{\mathcal U}}$,
$\widetriangle f
: (X,Z;\widetriangle{\mathcal U}) \to M$  a
strongly smooth map, and $\widetriangle{\frak S}$ a CF-perturbation
of $({\widetriangle{\mathcal U}}, \mathcal K)$.
We assume that
\begin{enumerate}
\item
${\widetriangle f}$ is strongly submersive with
respect to $\widetriangle{\frak S}$.
\item
$\widetriangle h$ has a compact support in $\ring Z$.$\blacksquare$
\end{enumerate}
\end{shitu}
\begin{choi}\label{cho776}
In Situation \ref{situ774} we make the following
choices.
\begin{enumerate}
\item
A triple of support systems $\mathcal K^1,\mathcal K^2, \mathcal K^3$
with
$\mathcal K^1 < \mathcal K^2< \mathcal K^3 < \mathcal K$.
\item
A constant $\delta>0$ such that:
\begin{enumerate}
\item
$\mathcal K^1(2\delta)$ is
compact. (Definition \ref{def761}.)
\item
There exists a strongly smooth partition of unity  $\{\chi_{\frak p}\}$
of $(X,Z,{\widetriangle{\mathcal U}},\mathcal K^1,\delta)$.
(Proposition \ref{pounitexi}.)
\item
$\mathcal K^1(2\delta) < \mathcal K^2$.
\item
$\delta$ satisfies the conclusion of Lemma \ref{lem739}
for $\mathcal K^- = \mathcal K^1$,
 $\mathcal K^+ = \mathcal K^2$,  $\mathcal K' = \mathcal K^3$.
\item
We put
$$
\delta_0 = \inf \{d(\mathcal K^2_{\frak p},\mathcal K^2_{\frak q}) \mid \text{neither $\frak p\le\frak q$ nor
$\frak q \le \frak p$}\},
$$
where we use the metric $d$ of $\vert\mathcal K\vert$. Then $\delta < \delta_0/2$.
\end{enumerate}
\item
We take a strongly smooth partition of unity $\{\chi_{\frak p}\}$
of $(X,Z,{\widetriangle{\mathcal U}},\mathcal K^1,\delta)$.
\item
We take an open neighborhood $\frak U(X)$ of $Z$ in
$\vert\mathcal K\vert$ such that the conclusion of Lemma \ref{lem740}
holds.
\end{enumerate}
\end{choi}
\begin{defn}\label{pushforwardKuranishi}
In Situation \ref{situ774},
we make Choice \ref{cho776}.
We define a smooth differential
form
$
{\widetriangle f}!(\widetriangle h;\widetriangle{{\frak S}^{\epsilon}})$
on the manifold $M$
by (\ref{formula714}). We call it
the {\it pushout} or {\it integration along the fiber}
\index{integration along the fiber (pushout) ! of differential form on good coordinate system}
of $\widetriangle h$ by $(\widetriangle f,\widetriangle{{\frak S}^{\epsilon}})$.
\par
\begin{equation}\label{formula714}
{\widetriangle f}!(\widetriangle h;\widetriangle{{\frak S}^{\epsilon}})
=
\sum_{\frak p\in \frak P}
f_{\frak p}!(\chi_{\frak p}h_{\frak p};{\frak S}^{\epsilon}_{\frak p}
\vert_{\frak U(Z) \cap \mathcal K^1_{\frak p}(2\delta)}).
\end{equation}
We note that the restriction of $\chi_{\frak p}h_{\frak p}$ to
$\mathcal K^1_{\frak p}(2\delta)$ has compact support in
${\rm Int}\,\mathcal K^1_{\frak p}$.
Therefore the right hand side of (\ref{formula714}) makes sense.
\par
The degree is given by
\begin{equation}
\deg {\widetriangle f}!(\widetriangle h;\widetriangle{{\frak S}^{\epsilon}})
=
\deg \widehat h + \dim M - \dim \widetriangle{\mathcal U}.
\end{equation}
\end{defn}
\begin{defn}\label{defnspadesuit}
Let $F_a : (0,\epsilon_a) \to \mathscr X$ be a family
of maps parameterized by $a \in \mathscr B$.
We say that $F_a$ is {\it independent of the choice
of $a$ in the sense of $\spadesuit$}
\index{independent of the choice in the sense of $\spadesuit$} if the following holds.
\begin{enumerate}
\item[$\spadesuit$]
For $a_1,a_2  \in \mathscr B$ there exists $0 < \epsilon_0
< \min\{\epsilon_{a_1},\epsilon_{a_2}\}$
such that $F_{a_1}(\epsilon) = F_{a_2}(\epsilon)$
for all $\epsilon < \epsilon_0$.
\end{enumerate}
\end{defn}
\begin{prop}\label{indepofukuracont}
In Situation \ref{situ774}
the pushout
${\widetriangle f}!(\widetriangle h;\widetriangle{{\frak S}^{\epsilon}})$
is independent of Choice \ref{cho776}
in the sense of $\spadesuit$.
\end{prop}
\begin{rem}
However ${\widetriangle f}!(\widetriangle h;\widetriangle{{\frak S}^{\epsilon}})$
depends on the choices of $\epsilon$ and $\widetriangle{\frak S}$.
\end{rem}

\begin{proof}[Proof of Proposition \ref{indepofukuracont}]
We first show the independence of $\frak U(Z)$.
Lemma \ref{lem743} implies that if $\frak U'(Z)  \subset \frak U(Z)$
is another open neighborhood of $Z$ then,
for sufficiently small $\epsilon$, the value of the
right hand side  of (\ref{formula714}) does not change if we replace $\frak U(Z)$
by $\frak U'(Z)$.
(We use Situation \ref{situ774} (2) also here.)
This implies independence of $\frak U(Z)$.
Moreover it implies that we can always replace $\frak U(Z)$ by
a smaller open neighborhood of $Z$.
\par
We next show independence of $\mathcal K^2$, $\mathcal K^3$.
Let $\mathcal K^{2 \prime}$, $\mathcal K^{3 \prime}$ be alternative choices of $\mathcal K^2$, $\mathcal K^3$.
We take $\mathcal K^{2 \prime\prime}_{\frak p} =\mathcal K^2_{\frak p}  \cup \mathcal K^{2 \prime}_{\frak p}$, $\mathcal K^{3 \prime\prime}_{\frak p} =\mathcal K^3_{\frak p}  \cup \mathcal K^{3 \prime}_{\frak p}$.
Then $\mathcal K^{2 \prime\prime}$,
$\mathcal K^{3 \prime\prime}$ are support systems.
Note in Definition \ref{pushforwardKuranishi}, the support system $\mathcal K^2$, $\mathcal K^3$
appears only when we apply Lemma \ref{lem740} to obtain
$\frak U(Z)$ and Lemma \ref{lem739} to obtain $\delta$. Since we can always replace $\frak U(Z)$ by a
smaller neighborhood of $Z$ (as far as $\epsilon>0$ is sufficiently small)
and
since we do not need to change $\delta$ in Lemma \ref{lem739} when we replace $\mathcal K^2$, $\mathcal K^3$ by
$\mathcal K^{2 \prime\prime} \supset \mathcal K^2$,
$\mathcal K^{3 \prime\prime} \supset \mathcal K^3$,
it follows that we obtain the same number in (\ref{formula714})
if we replace $\mathcal K^2$ or $\mathcal K^{2 \prime}$ by
$\mathcal K^{2 \prime\prime}$, as far as $\epsilon>0$ is sufficiently small.
This implies independence of $\mathcal K^2$, $\mathcal K^3$.
\par
It remains to prove the independence of $\mathcal K^1$ and of $\{\chi_{\frak p}\}$, $\delta$.
We will prove the independence for this case below.
\par
Let $\mathcal K^1_{\frak p}, \chi_{\frak p}, \delta$ and $\mathcal K^{1 \prime}_{\frak p},
\chi'_{\frak p}, \delta'$ be
two such choices. We
take $\mathcal K^{1 \prime\prime}_{\frak p}
= \mathcal K^1_{\frak p} \cup \mathcal K^{1 \prime}_{\frak p}$.
Then  $\mathcal K^{1 \prime\prime} < \mathcal K^2 <\mathcal K^3$. We can also take $\delta'' > 0$ and a
strongly smooth partition of unity $\{\chi''_{\frak p}\}$ of
$(X,Z,{\widetriangle{\mathcal U}},\mathcal K^{1 \prime\prime},\delta'')$.
So it suffices to prove that
the pushout defined by $\{\mathcal K^1_{\frak p}\}, \{\chi_{\frak p}\}, \delta$ coincides
with one defined by  $\{\mathcal K^{1 \prime\prime}_{\frak p}\},\{\chi''_{\frak p}\}, \delta''$
and that
the pushout defined by $\{\mathcal K^{1 \prime}_{\frak p}\}, \{\chi'_{\frak p}\}
, \delta'$
 coincides
with one defined by  $\{\mathcal K^{1 \prime\prime}_{\frak p}\},\{\chi''_{\frak p}\},
\delta''$.
In other words, we may assume $\mathcal K^1_{\frak p}
\subset \mathcal K^{1 \prime}_{\frak p}$.
We will prove the independence in this case.
\par
We  observe that
\begin{equation}\label{form715}
\widetriangle f!(\widetriangle h_1+\widetriangle h_2;\widetriangle{{\frak S}^{\epsilon}})
=
\widetriangle f!(\widetriangle h_1;\widetriangle{{\frak S}^{\epsilon}})
+
\widetriangle f!(\widetriangle h_2;\widetriangle{{\frak S}^{\epsilon}})
\end{equation}
as far as we use the same strongly smooth partition of unity in
all these three terms.
(This is a consequence of Lemma \ref{lem721} (2).)
We take $\frak p_0 \in \frak P$ and put
\begin{equation}\label{form738373}
\widetriangle h_0 = \chi'_{\frak p_0} \widetriangle h.
\end{equation}
In view of (\ref{form715}) we find that to prove Proposition \ref{indepofukuracont} it suffices
to show the next formula.
\begin{equation}\label{form716}
f_{\frak p_0}!((\widetriangle h_0)_{\frak p_0};{\frak S}^{\epsilon}_{\frak p_0}
\vert_{\frak U(Z) \cap \mathcal K^{1\prime}_{\frak p_0}(2\delta')})
=
\sum_{\frak p}f_{\frak p}!((\chi_{\frak p}\widetriangle h_0)_{\frak p};{\frak S}^{\epsilon}_{\frak p}
\vert_{\frak U(Z) \cap \mathcal K^1_{\frak p}(2\delta)}).
\end{equation}
By taking $\epsilon>0$ small, we may assume $\sum \chi_{\frak p} = 1$ on the intersection of
$\frak U(Z)$ and
the support set $\Pi((\widetriangle{{\mathfrak S}_{\frak p_0}^{\epsilon}})^{-1}(0))$.
(This is a consequence of Lemma \ref{lem743}, Definition \ref{pounity} (3)
and Situation \ref{situ774} (2).)
Therefore, to prove (\ref{form716}) it suffices to prove the next formula
for each $\frak p$.
\begin{equation}\label{form717}
f_{\frak p_0}!((\chi_{\frak p}\widetriangle h_0)_{\frak p_0};{\frak S}^{\epsilon}_{\frak p_0}
\vert_{\frak U(X) \cap \mathcal K^{1\prime}_{\frak p_0}(2\delta')})
=
f_{\frak p}!((\chi_{\frak p}\widetriangle h_0)_{\frak p};{\frak S}^{\epsilon}_{\frak p}
\vert_{\frak U(X) \cap \mathcal K^1_{\frak p}(2\delta)})
\end{equation}
whose proof is now in order.
\par
In case $\frak p_0 = \frak p$, (\ref{form717}) follows from
$$
{\rm Supp} (\chi_{\frak p} h_{0,\frak p_0})
= {\rm Supp} (\chi_{\frak p_0} h_{0,\frak p_0})
\subseteq
\mathcal K^1_{\frak p_0}(2\delta) \cap \mathcal K^{1\prime}_{\frak p_0}(2\delta').
$$
If neither $\frak p \le \frak p_0$ nor $\frak p_0 \le \frak p$, then
both sides of (\ref{form717}) are zero because
\begin{equation}\label{form743}
{\rm Supp}\,(\chi_{\frak p} h_{0,\frak p_0}) \subseteq
\Omega_{\frak p_0}(\mathcal K^{1\prime},\delta') \cap \Omega_{\frak p}(\mathcal K^{1},\delta)
= \emptyset.
\end{equation}
Note the second equality of (\ref{form743}) is a consequence of Choice \ref{cho776} (2)(e).
\par
We will discuss the other two cases below.
\par\smallskip
\noindent(Case 1): $\frak p > \frak p_0$.
\par
By definition $\Omega_{\frak p_0}(\mathcal K^{1 \prime},\delta')
\subset B_{\delta'}(\mathcal K^{1\prime}_{\frak p_0})$.
Therefore by (\ref{form738373}) the support of $\widetriangle h_0$ is in $B_{\delta'}(\mathcal K^{1\prime}_{\frak p_0})$.
By taking $\epsilon>0$ small, Lemma \ref{lem739} implies
$$
{\rm Supp}(\widetriangle h_0) \cap \Pi((\widetriangle{{\frak S}^{\epsilon}})^{-1}(0))
\subset
\mathcal K^{2\prime}_{\frak p_0} \cap B_{\delta'}(\mathcal K^{1\prime}_{\frak p_0})
\subset
\mathcal K^{1\prime}_{\frak p_0}(2\delta').
$$
Therefore
\begin{equation}\label{form741741}
{\rm Supp}(\chi_{\frak p}\widetriangle h_0) \cap \Pi((\widetriangle{{\frak S}^{\epsilon}})^{-1}(0))
\cap \frak U(Z)
\subseteq
\mathcal K^{1\prime}_{\frak p_0}(2\delta')
\cap
\mathcal K^{1}_{\frak p}(2\delta)
\cap \frak U(Z).
\end{equation}
Then (\ref{form717}) follows from Definition \ref{defn7732} (2)(3).
\par\smallskip
\noindent(Case 2): $\frak p_0 > \frak p$.
\par
By definition $\Omega_{\frak p}(\mathcal K^1,\delta)
\subset B_{\delta}(\mathcal K^1_{\frak p})$.
Therefore the support of $\chi_{\frak p}\widetriangle h_0$ is in $B_{\delta}(\mathcal K^1_{\frak p})$.
Therefore by taking $\epsilon>0$ small, Lemma \ref{lem739} implies
$$
{\rm Supp}(\chi_{\frak p}\widetriangle h_0) \cap \Pi((\widehat{{\frak S}^{\epsilon}})^{-1}(0))
\subset \mathcal K^2_{\frak p} \cap B_{\delta}(\mathcal K^1_{\frak p_0})
\subset \mathcal K^1_{\frak p}(2\delta).
$$
It implies (\ref{form741741}).
Then (\ref{form717}) follows from Definition \ref{defn7732} (2)(3).
\par
The proof of Proposition \ref{indepofukuracont} is complete.
\end{proof}
\begin{rem}
The pushforward (\ref{formula714}) is independent of the choice of
the support system $\mathcal K$ appearing in
Situation \ref{situ774}, as far as $\widetriangle{\frak S}$ and $\widetriangle h$ are defined on it.
In fact $\mathcal K$ does not appear in the definition.
\end{rem}
\begin{lem}\label{lem782}
Let $h,h_1,h_2$ be differential forms on $(X,Z;\widetriangle{\mathcal U})$
and $c_1, c_2 \in \R$.
\begin{enumerate}
\item
$
\widetriangle f!(c_1\widetriangle h_1 + c_1\widetriangle h_2;\widetriangle{{\frak S}^{\epsilon}})
=
c_1\widetriangle f!(\widetriangle h_1;\widetriangle{{\frak S}^{\epsilon}})
+
c_2\widetriangle f!(\widetriangle h_2;\widetriangle{{\frak S}^{\epsilon}}).
$
\item
If $\rho \in \Omega^*(M)$, then
$
\widetriangle f!(\widetriangle h \wedge \widetriangle f^*(\rho);\widetriangle{{\frak S}^{\epsilon}})
=
\widetriangle f!(\widetriangle h;\widetriangle{{\frak S}^{\epsilon}}) \wedge \rho.
$
\end{enumerate}
\end{lem}
\begin{proof}
Formula (1) follows from Lemma \ref{lem721} (1).
Formula (2) is immediate from definition.
\end{proof}
We now define the smooth correspondence.
\begin{cons}\label{cordefjyunbi}
Suppose we are in Situation \ref{smoothcorr}.
We construct objects as in Situation \ref{situ774} as follows.
\par
We put $Z=X$.
We take a good coordinate system
${\widetriangle{\mathcal U}}$ compatible with
${\widehat{\mathcal U}}$
such that $\widehat f_s$ and $\widehat f_t$
are pullbacks of $\widetriangle f_s : (X;{\widetriangle{\mathcal U}})
\to M_s$ and $\widetriangle f_t : (X;{\widetriangle{\mathcal U}})
\to M_t$, respectively.
Moreover we may take
$\widetriangle f_t$ so that it is weakly submersive.
(Proposition \ref{le614} (2).)
\par
We take a CF-perturbation
$\widetriangle{\frak S}$ of $(X;{\widetriangle{\mathcal U}})$
such that $\widetriangle f_t$ is strongly submersive
with respect to ${\widetriangle{\mathcal U}}$.
(Theorem \ref{existperturbcont} (2).)
\par
Let ${\mathcal K}$ be a support system of  ${\widetriangle{\mathcal U}}$.
We consider the differential form
$\widetriangle f_s^*h$ on $(X,{\widetriangle{\mathcal U}})$.
\end{cons}
We denote the correspondence by
$$
{\frak X} = ((X;\widetriangle{\mathcal U});\widetriangle f_s,
\widetriangle f_t).
$$
\begin{defn}\label{defn748}
Using Construction \ref{cordefjyunbi}, we  define
\begin{equation}\label{scordef}
{\rm Corr}_{(\frak X,\widetriangle{{\frak S}^{\epsilon}})}(h)
=
({\widetriangle f}_t)!(\widetriangle f_s^*h;\widetriangle{{\frak S}^{\epsilon}}).
\end{equation}
We call the linear map
$$
{\rm Corr}_{(\frak X,\widetriangle{{\frak S}^{\epsilon}})} :
\Omega^*(M_s)
\to
\Omega^{*+ \dim M_t - \dim \widetriangle{\mathcal U}}(M_t)
$$
\index{smooth correspondence ! associated to good coordinate system}
\index{good coordinate system ! smooth correspondence associated to good coordinate system}
the {\it smooth correspondence map
associated to
$\widetriangle{\frak X} = ((X;\widetriangle{\mathcal U});\widetriangle f_s,
\widetriangle f_t)$}.
\end{defn}
\begin{rem}\label{rem:785}
\begin{enumerate}
\item
Proposition \ref{indepofukuracont} implies that the right hand side of (\ref{scordef})
is independent of various choices appearing in
Definition \ref{pushforwardKuranishi}
 if $\epsilon>0$ is
sufficiently small.
However, it {\it does} depend on $\widetriangle{\frak S}$ and
$\epsilon > 0$.
So we keep the symbol $\epsilon$ in the notation
$\widetriangle{{\frak S}^{\epsilon}}$
of the left hand side of (\ref{scordef}).
\item
There seems to be no way to define a
smooth correspondence in the way that it becomes independent of the
choices {\it in the chain level}.
This is related to an important point of the
story, that is,  the construction of various structures
from system of moduli spaces are well-defined only
up to homotopy equivalence and only
as a whole.
(We however note that the method of \cite{joyce}
seems to be a way to minimize this dependence.)
This is the fundamental issue which
appears in {\it any} approach.
For example, it should remain to exist in the infinite dimensional
approach to construct virtual fundamental chain,
such as those by Li-Tian \cite{LiTi98}, Liu-Tian \cite{LiuTi98},
Siebert \cite{Siebert},
Chen-Tian \cite{ChenTian}, Chen-Li-Wang \cite{ChenLieWang}
or Hofer-Wyscoski-Zehnder \cite{hwze}.
\item
Dependence of the good coordinate system
${\widetriangle{\mathcal U}}$
and the other choices made in Construction \ref{cordefjyunbi} will be discussed in Section \ref{sec:kuraandgood}.
\end{enumerate}
\end{rem}
In Proposition \ref{indepofukuracont}, we have proved independence of the
pushout of various choices for sufficiently small $\epsilon>0$.
On the other hand, how $\epsilon$ must be small depends on our good coordinate system
and CF-perturbation on it.
In certain situations appearing in applications,
we need to estimate this required smallness of $\epsilon$
uniformly from below when our CF-perturbations
vary in a certain family.
The next proposition can be used for such a purpose.
\begin{defn}\label{intheclubsuit}
Let $F_{\sigma,a} : (0,\epsilon_a) \to \mathscr X$ be a family
of maps parameterized by $a \in \mathscr B$
and $\sigma \in \mathscr A$.
We say that $F_{\sigma,a}$ is {\it
uniformly independent of the choice
of $a$ in the sense of $\clubsuit$}
\index{uniformly independent of the choice in the sense of $\clubsuit$} if the following holds.
\begin{enumerate}
\item[$\clubsuit$]
For $a_1,a_2  \in \mathscr B$ there exists $0 < \epsilon_0
< \min\{\epsilon_{a_1},\epsilon_{a_2}\}$ independent of $\sigma$
such that $F_{\sigma,a_1}(\epsilon) = F_{\sigma,a_2}(\epsilon)$
for $0<\epsilon < \epsilon_0$ and any $\sigma \in \mathscr A$.
\end{enumerate}
\end{defn}
\begin{prop}\label{lem761}
We assume $\{ \widetriangle{\frak S_{\sigma}} \mid \sigma \in \mathscr A\}$  is a
uniform family of CF-perturbations parameterized by $\sigma \in \mathscr A$
in the sense of Definition \ref{uniformcongpert}.
Then we can make Choice \ref{cho776} in a way independent of $\sigma$.
\par
Moreover the pushout
${\widetriangle f}!(\widetriangle h;\widetriangle{{\frak S}_{\sigma}^{\epsilon}})$
of Proposition \ref{indepofukuracont}  is uniformly
 independent of  Choice \ref{cho776}
in the sense of $\clubsuit$.
\end{prop}
\begin{proof}
From the proof of Proposition \ref{indepofukuracont}, the proof of
Proposition \ref{lem761} follows from the next lemma.
\end{proof}

\begin{lem}\label{lem6878787}
Let $\{ \widetriangle{\frak S_{\sigma}} \mid \sigma \in \mathscr A\}$
be a uniform family of CF-perturbations.
Then the following holds.
\begin{enumerate}
\item
In Lemma \ref{lem739} the constants $\delta_0$ and $\epsilon_0$
can be taken in dependent of $\sigma$.
\item
In Lemma \ref{lem740} the set $\frak U(Z)$ and the
constant $\epsilon_0$ can be taken independent of $\sigma$.
\item
In Lemma \ref{lem743} the set $\frak U(Z)$ and the
constant $\epsilon_0$ can be taken independent of $\sigma$.
Moreover
\begin{equation}
\lim_{n\to \infty}\sup\left\{ d_H\left(X, \bigcup_{\frak p}(\Pi((\widetriangle{{\mathfrak S}_{\sigma,\frak p}^{\epsilon}})^{-1}(0)
\cap  \overset{\circ}{\mathcal K^2_{\frak p}}
\cap \frak U(Z))\right) \mid \sigma \in \mathscr A\right\}
= 0.
\end{equation}
\end{enumerate}
\end{lem}
\begin{proof}
Using Lemma \ref{lemma748},
the proof of Lemma \ref{lem6878787}
is the same as that of Proposition \ref{lem627}.
\end{proof}

\section{Stokes' formula}
\label{sec:stokes}
\subsection{Boundary and corner II}
\label{subsection:normbdry2}
In this section, we state and prove Stokes' formula.
We first discuss the notion of boundary or corner
of an orbifold and of a Kuranishi structure in more detail.
(The discussion below is a detailed version of \cite[the last paragraph of page 762]{fooobook2}.
See also \cite[page 11]{joyce}.
\cite{joyce3} gives a systematic account on this issue.
\par
Let $U$ be an orbifold with boundary and corner.
We defined its corner structure stratification $S_k(U)$ and $\overset{\circ}S_k(U)$ in Definition \ref{defn4111}.
Note $\overset{\circ}S_k(U)$ is an orbifold of dimension $\dim U -k$ without boundary.
However we also note that we can {\it not} find a structure of orbifold with corner on $S_k(U)$
such that  $\overset{\circ}S_0(S_k(U)) = \overset{\circ}S_k(U)$.
\begin{exm}
Let $U = \R_{\ge 0}^2$. Then
$S_1(U)$ is homeomorphic to $\R$ and
$S_2(U)$ is one point identified with $0 \in \R = S_1(U)$.
\end{exm}
To obtain an orbifold with corner from $S_1(U)$ we need to modify it at its
boundaries and corners. Let us first consider the case of manifolds.
\begin{lem}\label{lemma749}
Suppose $U$ is a manifold with corner.
Then there exists a manifold with corner, denoted by $\partial U$, and a map
$\pi : \partial U \to S_1(U)$ with the following properties.
\begin{enumerate}
\item
For each $k$, $\pi$ induces a map
\begin{equation}\label{maponkstratum}
\overset{\circ}S_k(\partial U) \to \overset{\circ}S_{k+1}(U).
\end{equation}
\item
The map (\ref{maponkstratum}) is a $(k+1)$-fold covering map.
\item
$\pi$ is a smooth map
$
\partial U \to U.
$
\end{enumerate}
\end{lem}
\begin{rem}
The smoothness claimed in Lemma \ref{lemma749} (3) is defined as follows.
Let $U$ be any smooth manifold with corner. We
can find a smooth manifold without boundary or
corner $U^+$ of the same dimension
as $U$ and an embedding $U \to U^+$, such that
for each point $p\in U$ there exists a coordinate of $U^+$ at $p$
by which $U$ is identified with an open subset of $[0,1)^{\dim U}$
by a diffeomorphism from $U^+$ onto an open subset of $\R^{\dim U}$.
Then a map $F : U_1 \to U_2$ between two manifolds with corners is
said to be smooth if $F$ extends to $F^+$ that is a smooth map
from a neighborhood of $U_1$ in $U_1^+$ to $U_2^+$.
\end{rem}
\begin{proof}
We fix a Riemannian metric on $U$ so that each $\epsilon$-ball $B_{\epsilon}(p)$ is convex.
Let $p \in \overset{\circ}S_k(U)$.
We consider $\alpha \in \pi_0(B_{\epsilon}(p) \cap \overset{\circ}S_{1}(U))$.
The set of all such pairs $(p,\alpha)$ with $k\ge 1$ is our $\partial U$.
The map $(p,\alpha) \to p$ is the map $\pi$.
\par
By identifying $U$ locally with $[0,\infty)^n$,
it is easy to construct the structure of manifold with corner on $\partial U$ and prove that they have the required properties.
\end{proof}
\begin{defn}\label{defbdrofd}
Let $U$ be an orbifold.
We call $\partial U$ a {\it normalized boundary}
\index{normalized boundary ! of orbifold}
\index{orbifold ! normalized boundary}
\index{boundary ! {\it see normalized boundary}}
 of $U$ and $\pi : \partial U \to S_1(U)$
the \index{normalized boundary ! normalization map from normalized boundary of
orbifold} {\it normalization map}.
\end{defn}
\begin{lem}\label{lem750}
\begin{enumerate}
\item
Let $U$ and $U'$ be as in Lemma \ref{lemma749} and
$F : U \to U'$ a diffeomorphism. Then
$F$ uniquely induces a diffeomorphism
$$
 F^{\partial} : \partial U \to \partial U'
$$
such that $\pi \circ  F^{\partial} = F \circ \pi$.
\item
Suppose a finite group $\Gamma$ acts on $U$, where $U$ is  as in Lemma \ref{lemma749}.
Suppose also that each connected component of $\overset{\circ}S_k(U)/\Gamma$
is an effective orbifold for each $k$.
Then $\Gamma$ acts on $\partial U$ so that $\pi$ is $\Gamma$ equivariant
and each connected component of $\overset{\circ}S_k(\partial U)/\Gamma$
is an effective orbifold for each $k$.
\item
Let $U$ be as in Lemma \ref{lemma749} and $U'$ its open subset.
Then there exits an open embedding $\partial U' \to \partial U$ which
commutes with $\pi$.
\end{enumerate}
\end{lem}
\begin{proof}
(1) is immediate from the construction.
Then the uniqueness implies (2) and (3).
\end{proof}
Now we consider the case of an orbifold $U$.
We cover $U$ by orbifold charts $\{ (V_i,\Gamma_i,\phi_i)\}$.
We apply Lemma  \ref{lemma749} to $V_i$ and obtain $\partial V_i$.
Then $\Gamma_i$ action on $V_i$ induces one on $\partial V_i$ by Lemma \ref{lem750} (2).
We thus obtain orbifolds $\partial V_i/\Gamma_i$.
Using Lemma \ref{lem750} (1) and (3) we can glue $\partial V_i/\Gamma_i$ for various $i$
and obtain $\partial U$. We obtain also a map $\pi : \partial U \to S_1(U)$.
It induces a map
$\overset{\circ}S_k(\partial U) \to \overset{\circ}S_{k+1}(U)$.

\begin{rem}\label{rem8ten6}
\begin{enumerate}
\item
We note that the map $\overset{\circ}S_k(\partial U)
\to \overset{\circ}S_{k+1}(U)$ is a
$(k+1)$-fold orbifold covering  of orbifolds in the sense we will define in Part 2.
\item
In particular,
$\overset{\circ}S_0(\partial U) \to \overset{\circ}S_{1}(U)$
is a diffeomorphism of orbifolds.
\item
We also note that
$\overset{\circ}S_k(\partial U) \to \overset{\circ}S_{k+1}(U)$
is not necessarily a $k+1$ to $1$ map set-theoretically.
The following is a counter example.
Let
$$
U = (\R_{\ge 0}^2 \times \R)/\Z_2
$$
where the action is $(a,b,c) \mapsto (b,a,-c)$.
Then $\partial U \cong \R_{\ge 0} \times \R$,
$S_1(\partial U) \cong \R$, $S_2(U) \cong \R/\Z_2$ and the map $\pi$ is
canonical projection
$\R \to \R/\Z_2$ on $S_1(\partial U)$.
So it is generically 2 to 1 but is 1 to 1 at $0$.
\end{enumerate}
\end{rem}
Next we consider the case of Kuranishi structure.
We recall the following notation from Definition \ref{dimstratifidef}.
$$
S_k(X,Z;\widehat{\mathcal U}) = \{ p \in Z \mid o_p \in S_k(U_p)\},
\qquad
\overset{\circ}{S_k}(X,Z;\widehat{\mathcal U}) = \{ p \in Z \mid o_p \in \overset{\circ}{S_k}(U_p)\},
$$
where $\widehat{\mathcal U}$ is a Kuranishi structure of $Z \subseteq X$
and
$$
\aligned
S_k(X,Z;{\widetriangle{\mathcal U}}) &= \{ p \in Z \mid \exists \frak p \exists x \in S_k(U_{\frak p}), {\rm st}\,\,
s_{\frak p}(x) = 0, \psi_{\frak p}(x) = p\},
\\
\overset{\circ}{S_k}(X,Z;{\widetriangle{\mathcal U}}) &= \{ p \in Z \mid \exists \frak p \exists x \in \overset{\circ}{S_k}(U_{\frak p}), {\rm st}\,\,
s_{\frak p}(x) = 0, \psi_{\frak p}(x) = p\},
\endaligned
$$
where $\widetriangle{\mathcal U}$ is a good coordinate system of $Z \subseteq X$.
\begin{rem}
We can rewrite the set $S_k(X,Z;{\widetriangle{\mathcal U}})$ as
$$
S_k(X,Z;{\widetriangle{\mathcal U}}) = \{ p \in X \mid \forall \frak p \forall x \in U_{\frak p},
s_{\frak p}(x) = 0, \psi_{\frak p}(x) = p \,\,\Rightarrow\,\, x \in {S_k}(U_{\frak p})
\}.
$$
A similar remark applies to $\overset{\circ}{S_k}(X,Z;{\widetriangle{\mathcal U}})$.
\end{rem}
\begin{lemdef}\label{lem754}
\begin{enumerate}
\item
Any compact subset of the space
$\overset{\circ}{S_k}(X,Z;\widehat{\mathcal U})$
(resp. $\overset{\circ}{S_k}(X,Z;{\widetriangle{\mathcal U}})$) has Kuranishi structure without boundary
(resp. good coordinate system without boundary) and of
dimension $\dim (X,Z;\widehat{\mathcal U}) - k$
(resp. $\dim (X,Z;{\widetriangle{\mathcal U}}) - k$).
\item
There exist a relative K-space with corner
$\partial(X,Z;\widehat{\mathcal U})$
(resp. $\partial(X,Z;{\widetriangle{\mathcal U}})$)
whose underlying topological spaces are
$(\partial X,\partial Z)$
and a continuous map between their underlying topological spaces
$\pi : \partial Z \to S_1(X,Z;\widehat{\mathcal U})$
(resp. $\pi : \partial Z \to S_1(X,Z;\widetriangle{\mathcal U})$)
such that the following holds.
We call $\partial(X,Z;\widehat{\mathcal U})$,
$\partial(X,Z;{\widetriangle{\mathcal U}})$
the {\rm normalized boundary}
\index{normalized boundary ! of Kuranishi structure}
\index{Kuranishi structure ! normalized boundary} of
$(X,Z;\widehat{\mathcal U})$,
$(X,Z;{\widetriangle{\mathcal U}})$, respectively.
\begin{enumerate}
\item
If $\pi(\tilde p) = p$, $\tilde p \in \partial Z$, $p \in Z$, then
the Kuranishi neighborhood of $\tilde p$ is obtained
by restricting $\mathcal U_p$ to
$\partial U_p$, which is as in Definition  \ref{defbdrofd}.
\item
The coordinates of $\partial(X,Z;{\widetriangle{\mathcal U}})$
are obtained by restricting $\mathcal U_{\frak p}$ to $\partial U_{\frak p}$.
\item
The coordinate change of $\partial(X,Z;\widehat{\mathcal U})$
(resp. $\partial(X,Z;{\widetriangle{\mathcal U}})$)
is obtained by restricting one of $\partial \mathcal U_{p}$
(resp. $\partial \mathcal U_{\frak p}$).
\item
The restriction of $\pi$ induces a map
$$
\overset{\circ}{S_0}(\partial(X,Z;\widehat{\mathcal U}))
\to \overset{\circ}{S_1}(X,Z;\widehat{\mathcal U})
$$
that is an isomorphism of Kuranishi structures.
\item
The restriction of $\pi$ induces a map
$$
\overset{\circ}{S_0}(\partial(X,Z;{\widetriangle{\mathcal U}}))
\to \overset{\circ}{S_1}(X,Z;{\widetriangle{\mathcal U}})
$$
that is an isomorphism of good coordinate systems.
\item
In the case of Kuranishi structure and $Z \ne X$,
we need to replace ${\widehat{\mathcal U}}$ by
its open substructure.
\end{enumerate}
\item
Various kinds of embeddings among Kuranishi structures
and/or good coordinate systems induce
embeddings
of their normalized boundaries.
\end{enumerate}
\end{lemdef}
\begin{proof}
We first prove (2).
Let $\mathcal U = (U,\mathcal E,s,\psi)$ be a Kuranishi chart.
We restrict $\mathcal E$ and $s$ to $\partial U$ and obtain
$\partial U,\mathcal E^{\partial}, s^{\partial}$.
We will define underlying topological spaces $\partial X$,
$\partial Z$,
parametrization $\psi^{\partial}$ and the coordinate change.
\par
Let $\Phi_{21}=(U_{21}, \varphi_{21},\widehat{\varphi}_{21}) : \mathcal U_1 \to \mathcal U_2$ be
a coordinate change of Kuranishi charts.
We note that we required the following condition for an embedding of orbifolds
$\varphi_{21} : U_1 \to U_2$.
\begin{equation}\label{cornerorbemb}
\varphi_{21}(S_k(U_1)) \subset S_k(U_2),
\qquad
\overset{\circ}{S_k}(U_1) = \varphi_{21}^{-1}(\overset{\circ}{S_k}(U_2)).
\end{equation}
We can generalize Lemma \ref{lem750} (3) to the case
when $U_1, U_2$ are orbifolds.
Moreover by the condition (\ref{cornerorbemb})
we can generalize Lemma \ref{lem750} (1)
to the embedding of orbifolds with corners.
Thus $\varphi_{21}$ induces
$\varphi^{\partial}_{21} : \partial U_1 \to \partial U_2$.
In the same way
$\widehat\varphi_{21} : \mathcal E_1 \to \mathcal E_2$
induces $\widehat\varphi_{21}^{\partial} : \mathcal E_1^{\partial}
\to \mathcal E_2^{\partial}$.
\par
Thus the data consisting of the coordinate change of the
Kuranishi charts given as in (2) (a),(b) are defined by (2) (c),
except the underlying topological space
 $\partial X$, $\partial Z$ and
parametrization $\psi$.
\par
Below we will construct $\partial X$, $\partial Z$ and
$\psi$.
\par
We first consider the case
of good coordinate system and $X = Z$.
Let $\widetriangle{\mathcal U} =(\frak P,\{\mathcal U_{\frak p}\},\{\Phi_{\frak p\frak q}\})$.
We consider $\partial U_{\frak p}$ and
$\varphi_{\frak p\frak q}^{\partial}$, defined as above.
We glue the spaces $\partial U_{\frak p}$ by  $\varphi_{\frak p\frak q}^{\partial}$
and obtain a topological space
$\vert\partial \widetriangle{\mathcal U}\vert$.
The zero sets of the Kuranishi maps $s_{\frak p}^{\partial}$ on  $\partial U_{\frak p}$
are glued to define a subspace of $\vert\partial \widetriangle{\mathcal U}\vert$,
which we define to be $\partial X$.
$\partial X$ is Hausdorff and metrizable.\footnote{To find a topology which s metrizable,
we consider a support system $\mathcal K$ and
use the fact that $\partial X
= \bigcup \mathcal K_{\frak p} \cap (s_{\frak p}^{\partial})^{-1}(0)$
and \cite[Proposition 5.1]{foooshrink}.}
We define $\psi_{\frak p}^{\partial} : (s_{\frak p}^{\partial})^{-1}(0)
\to \partial X$ by mapping a point of
$(s_{\frak p}^{\partial})^{-1}(0)$  to its equivalence class.
Then $\mathcal U^{\partial}_{\frak p} = (\partial U_{\frak p},\mathcal E_{\frak p}^{\partial},s_{\frak p}^{\partial},
\psi_{\frak p}^{\partial})$
is a Kuranishi chart of $\partial X$.
We put $\partial U_{\frak p\frak q} = U_{\frak p\frak q} \cap
\partial U_{\frak q}$.
Then
$\Phi_{\frak p\frak q}^{\partial} =
(\partial U_{\frak p\frak q},\varphi_{21}^{\partial},\widehat\varphi_{21}^{\partial})
$
is a coordinate change $\mathcal U^{\partial}_{\frak p}
\to \mathcal U^{\partial}_{\frak q}$.
Thus we obtain a good coordinate system $\partial(X,Z;{\widetriangle{\mathcal U}})$ in case $Z = X$.
\par
Next we consider the case
of good coordinate system but $X \ne Z$.
We glue the zero sets of Kuranishi map on  $\partial U_{\frak p}$
in the same way as above
to obtain a topological space $\partial X$.
(See Remark \ref{rem8989}.)
We define  the subset  $\partial Z \subset \partial X$
by
$$
\partial Z =
\bigcup_{\frak p \in \frak P}\{ x \in \partial U_{\frak p} \mid s_{\frak p}(x) = 0,
\psi_{\frak p}(\pi(x)) \in Z\}.
$$
Here we identify $\partial U_{\frak p}$ with  its image in $\vert\partial \widetriangle{\mathcal U}\vert$
and the union is taken in $\vert\partial \widetriangle{\mathcal U}\vert$.
The rest of the proof is the same as the case of $X = Z$.
\par
Finally we consider the case of Kuranishi structure.
We take a good coordinate system $\widetriangle{\mathcal U}$
compatible to $\widehat{\mathcal U}$
and use the case of good coordinate system
to define $\partial X$ and $\partial Z$.
\par
It now remains to define the parameterization $\psi^{\partial}_p :
(s^{\partial}_p)^{-1}(0) \to \partial X$.
Let $\mathcal U_p$ be a Kuranishi neighborhood of $p$ which is
a part of the data of
$\widehat{\mathcal U}$.
In case the embedding $\widehat{\mathcal U} \to \widetriangle{\mathcal U}$
is strict, $(s^{\partial}_p)^{-1}(0) \subset U_{\frak p}$
for some $\frak p \in \frak P$.
Therefore we obtain $\psi^{\partial}_p$  by restricting the parametrization
map $\psi_{\frak p}$ of the good coordinate system $\partial\widetriangle{\mathcal U}$.
\par
In the case $Z \ne X$ we replace $\widehat{\mathcal U}$
by its open substructure $\widehat{\mathcal U_0}$
such that there exists a strict embedding
$\widehat{\mathcal U_0} \to \widetriangle{\mathcal U}$.
\par
Suppose $Z = X$. We define $\psi^{\partial}_p : (s^{\partial}_p)^{-1}(0) \to \partial X$ (without taking open substructure) as follows.
Let $x \in (s^{\partial}_p)^{-1}(0) \subset s_p^{-1}(0)$
and $q = \psi_p(x) \in X$.
We have $o_q \in U_q$ such that $\varphi_{pq}(o_q) = x$.
(\ref{cornerorbemb}) implies $o_q \in \partial U_{q}$.
Moreover $o_q \in \partial U_{0,q}$.
(Here $U_{0,q}$ is the Kuranishi neighborhood of
the open substructure $\widehat{\mathcal U_0}$.)
Therefore $o_q$ may be regarded as an element of $\partial U_{\frak p}$
for some $\frak p \in \frak P$.
We define $\psi_p(x)$ to be the equivalence class of
$o_q \in \partial U_{\frak p}$, which is an element of $\partial X$.
\par
Therefore the proof of the statement (2) is complete.
\par
The statement (1) can be proved in the same way and the statement
(3) is obvious from definition.
\end{proof}
\begin{rem}\label{rem8989}
Here is a technical remark about the way to define
underlying topological space $\partial X$ in Lemma-Definition \ref{lem754}.
\begin{enumerate}
\item
Let $\widehat{\mathcal U}$ be a Kuranishi structure
of $Z \subseteq X$. Then the parametrization
$\psi_p : s_p^{-1}(0) \to X$ has $X$ as a target space.
So $\widehat{\mathcal U}$ is {\it not} a Kuranishi structure of $Z$ itself.
\item
The data consisting of $\widehat{\mathcal U}$
contain enough information to determine
which points of $Z$ lie in the boundary.
However the data do not contain such information
for the points of $X \setminus Z$ which are far away from $Z$.
\par
The space $\partial X$ that we defined in the proof of  Lemma-Definition
\ref{lem754}
consists of points which correspond to the  `boundary points' of $X$
that is sufficiently close to $Z$.
\par
Since the image of $\psi_p$, $p\in Z$ lies in a
neighborhood of $Z$, we need only a
neighborhood of $Z$ in $X$ to define
Kuranishi structure of $Z \subseteq X$.
This is the reason why it suffices to define $\partial X$
in a neighborhood of $Z$.
\item
On the other hand, as a consequence of (2), the topological space
$\partial X$ is not canonically determined from $(X,Z;\widehat{\mathcal U})$.
For example, the following phenomenon happens.
Let $\widehat {\mathcal U}$ be a Kuranishi structure of
$Z_2 \subseteq X$ and $Z_1 \subset \ring Z_2$.
We restrict ${\mathcal U}$ to $Z_1$ to obtain
${\mathcal U}\vert_{Z_1}$.
We consider
$
\partial(X,Z_1;\widehat{\mathcal U}\vert_{Z_1})
$
and
$
\partial(X,Z_2;\widehat{\mathcal U})
$.
Let $(\partial_1X,\partial Z_1)$
and $(\partial_2X,\partial Z_2)$
be their underlying topological spaces.
Then $\partial_1X$ may not be the same as $\partial_2X$.
\end{enumerate}
There is no such an issue in the absolute case $Z=X$.
In the applications the case $Z \ne X$ appears only
together with a means of defining $\del X$ given.
\end{rem}
We note that all the arguments of Section \ref{sec:contfamily}
work
for Kuranishi structure
or good coordinate system
with boundaries or corners.
\par
The next lemma describes the way how to restrict
CF-perturbations and etc.  to the normalized
boundary.

\begin{lem}\label{lemma755}
Let $\widetriangle{\mathcal U}$ be a good coordinate system
of $Z \subseteq X$,
${\mathcal K}^1,{\mathcal K}^2,{\mathcal K}^3$ a triple of
support systems of $\widetriangle{\mathcal U}$
 with ${\mathcal K}^1 <{\mathcal K}^2<{\mathcal K}^3$, and
$\widetriangle{\frak S^{\epsilon}}$ a CF-perturbation of
$({\widetriangle{\mathcal U}},\mathcal K^3)$.
\begin{enumerate}
\item
$\{\partial U_{\frak p} \cap \mathcal K^i_{\frak p}\}$ is a
support system of $\partial (X,Z;{\widetriangle{\mathcal U}})$,
which we denote by $\mathcal K^i_{\partial}$.
${\mathcal K}_{\partial}^1,{\mathcal K}_{\partial}^2,{\mathcal K}_{\partial}^3$ are support systems with
${\mathcal K}_{\partial}^1 < {\mathcal K}_{\partial}^2 < {\mathcal K}_{\partial}^3$.
Here $\partial U_{\frak p} \cap \mathcal K^i_{\frak p}
= \pi^{-1}(\mathcal K^i_{\frak p}) \subset \partial U_{\frak p}$.
\item
\begin{enumerate}
\item
For each $\frak p$ there exists ${\frak S}^{\partial}_{\frak p}$
that is a CF-perturbation of $\mathcal K^3_{\partial,\frak p} \subset \partial U_{\frak p}$.
\item
The  restriction of ${\frak S}^{\partial}_{\frak p}$
to $\overset{\circ}{S_0}(\partial U_{\frak p})$
is identified with the restriction of ${\frak S}_{\frak p}$
to $\overset{\circ}{S_1}(U_{\frak p})$ by the diffeomorphism
in Lemma \ref{lem754} (2)(e).
\item
The collection $\{{\frak S}^{\partial}_{\frak p}\}$ is a
CF-perturbation of
$(\partial(X,Z;{\widetriangle{\mathcal U}}),\mathcal K^3_{\partial})$,
which we denote by $\widetriangle{\frak S^{\partial}}$.
\item
If $\widetriangle{\frak S}$ varies in a uniform family
(in the sense of Definition \ref{uniformcongpert}) then $\widetriangle{\frak S^{\partial}}$
varies in a uniform family.
\end{enumerate}
\item
\begin{enumerate}
\item
A strongly continuous map $\widetriangle f :
(X,{\widetriangle{\mathcal U}})\to M$ induces a strongly continuous map $\widetriangle{f_{\partial}}
: \partial(X,Z;{\widetriangle{\mathcal U}})\to M$.
\item
The restriction of $\widetriangle{f_{\partial}}$ to
$\overset{\circ}{S_0}(\partial(X,{\widetriangle{\mathcal U}}))$
coincides with the restriction of $\widetriangle f$
to $\overset{\circ}{S_1}(X,{\widetriangle{\mathcal U}})$.
\item
If $\widetriangle f$ is strongly smooth (resp. weakly submersive),
so is $\widetriangle{f_{\partial}}$.
\end{enumerate}
\item
\begin{enumerate}
\item
If $\widetriangle{\frak S}$ is transversal to $0$,
so is $\widetriangle{\frak S^{\partial}}$.
\item
If $\widetriangle f$ is strongly submersive with respect to
$\widetriangle{\frak S}$  then
$\widetriangle{f_{\partial}}$ is strongly submersive with respect to
$\widetriangle{\frak S^{\partial}}$.
\item
Let $g : N \to M$ be a smooth map between smooth manifolds.
If $\widetriangle f$ is strongly smooth and weakly transversal to $g$ then
so is $\widetriangle{f_{\partial}}$.
\end{enumerate}
\item
\begin{enumerate}
\item
A differential form $\widetriangle h$  on $(X,Z;{\widetriangle{\mathcal U}})$
induces a differential form on $\partial(X,Z;{\widetriangle{\mathcal U}})$,
which we write $\widetriangle{h_{\partial}}$.
\item
The restriction of $\widetriangle{h_{\partial}}$ to
$\overset{\circ}{S_0}(\partial(X,Z;{\widetriangle{\mathcal U}}))$
coincides with the restriction of $\widetriangle h$
to $\overset{\circ}{S_1}(X,Z;{\widetriangle{\mathcal U}})$.
\item
In particular,
a strongly continuous function on $(X,Z,{\widetriangle{\mathcal U}})$
induces one on  $\partial(X,Z;{\widetriangle{\mathcal U}})$,
such that (b) above applies.
\end{enumerate}
\item
If $\{\chi_{\frak p}\}$  is a strongly smooth
partition of unity of $(X,Z,{\widetriangle{\mathcal U}},\mathcal K^2,\delta)$
then
$\{(\chi_{\frak p})_{\partial}\}$  is a strongly smooth
partition of unity of $(\partial X,\partial Z,\partial{\widetriangle{\mathcal U}},\mathcal K_{\partial}^2,\delta)$.
Here $(\chi_{\frak p})_{\partial}$ is one induced from
$\{\chi_{\frak p}\}$ as in (5) (c) above.
\end{enumerate}
\end{lem}
\begin{proof}
In the case when $U$ is an orbifold with corners, various transversality
or submersivity are defined by requiring the conditions
not only to $ \overset{\circ}{S_0}(U)$
(the interior point) but also
to all $\overset{\circ}{S_k}(U)$.
Once we observe this point all the statements are obvious from the definition.
\end{proof}
\subsection{Stokes' formula for a good coordinate system}
\label{subsec:Stokesgcs}
Now we are ready to state and prove Stokes' formula.
\begin{thm}{\rm (Stokes' formula, \cite[Lemma 12.13]{fooo09})}\label{Stokes}
Assume that we are in the situation of Lemma \ref{lemma755}
(1), (2), (3) (a)(b)(c), (4) (a)(b), (5) (a)(b).
Then, for each sufficiently small $\epsilon>0$, we have
\begin{equation}
d\left(\widetriangle f!(\widetriangle h;\widetriangle{{\frak S}^{\epsilon}})\right)
=
\widetriangle f!(d\widetriangle h;\widetriangle{{\frak S}^{\epsilon}})
+
\widetriangle f_{\partial}!(\widetriangle {h_{\partial}};\widetriangle{{\frak S}_{\partial}^{\epsilon}}).
\end{equation}
\end{thm}
\begin{proof}
Let $\{\chi_{\frak p}\}$ be a strongly smooth partition of unity
of $(X,Z,{\widetriangle{\mathcal U}},\mathcal K^2,\delta)$. Let
 $(\chi_{\frak p})_{\partial}$ and ${\mathcal K}_{\partial}^1$
be defined by Lemma \ref{lemma755}.
We put $h_0 = \chi_{\frak p} h_{\frak p}$.
It suffices to show
\begin{equation}\label{eq722}
\aligned
&d\left(f_{\frak p}!(h_0;{\frak S}_{\frak p}^{\epsilon})\vert_{\frak U(Z) \cap \mathcal K_{\frak p}^1(2\delta)}\right)
\\
&=
f_{\frak p}!(dh_0;{\frak S}_{\frak p}^{\epsilon}
\vert_{\frak U(Z) \cap \mathcal K_{\frak p}^1(2\delta)})
+
f^{\partial}_{\frak p}!(h_0;{\frak S}^{\partial,\epsilon}_{\frak p}
\vert_{\frak U(Z) \cap \mathcal K_{\partial,\frak p}^1(2\delta)}),
\endaligned
\end{equation}
where
$
\widetriangle{{\frak S}^{\partial,\epsilon}}
= \{{\frak S}^{\partial,\epsilon}_{\frak p} \mid \frak p \in \frak P\}.
$
Let ${\frak S}_{\frak p} = \{(\frak V_{\frak r},\mathcal S_{\frak r}^{\frak p}) \mid \frak r \in \frak R\}$ and $\{\chi_{\frak r}\}$ a
partition of unity subordinate to $\{U_{\frak r}\}$.
We put $h_1 = \chi_{\frak r}h_0$ and $f_{\frak r}
= f_{\frak p}\vert_{U_{\frak r}}$.
To prove (\ref{eq722}), it suffices to prove:
\begin{equation}\label{eq723}
\aligned
&d\left(f_{\frak r}!(h_1;{\mathcal S}_{\frak r}^{\epsilon})\vert_{\frak U(Z) \cap U_{\frak r}}\right)
\\
&=
f_{t,\frak r}!(dh_1;{\mathcal S}_{\frak r}^{\epsilon})\vert_{\frak U(Z) \cap U_{\frak r}})
+
f_{\frak r}^{\partial}!(h_1;{\mathcal S}^{\partial,\epsilon}_{\frak r}\vert_{\frak U(Z) \cap \partial U_{\frak r}}),
\endaligned
\end{equation}
where
$
\widetriangle{{\frak S}^{\epsilon}_{\partial,\frak p}}
=
\{(\partial \frak V_{\frak r},\mathcal S^{\partial,\epsilon}_{\frak r})\}
$.
(\ref{eq723}) follows from the next lemma.
\begin{lem}\label{lem759sss}
Let $\Omega$ be an open neighborhood of $0$ in
$[0,1)^m \times \R^{n-m}$ and $f : \Omega \to M$  a
smooth map. Let $h$ be a smooth differential $k$ form
on $\Omega$ with compact support and $\rho$  a
differential $(n-k-1)$-form on $M$.
Then we have
$$
(-1)^k\int_{\Omega} h \wedge f^*d\rho
=
\int_{\Omega \cap \partial ([0,1)^m \times \R^{n-m})}h
\wedge f^*\rho
+
\int_{\Omega} dh \wedge \rho.
$$
\end{lem}
\par
Lemma \ref{lem759sss} is an immediate consequence of the usual Stokes' formula.
Thus the proof of Theorem \ref{Stokes} is complete.
\end{proof}
Using Stokes' formula we can immediately
prove the following basic properties of
smooth correspondence.
\begin{cor}\label{Stokescorollary}
In Situation \ref{smoothcorr} we apply Construction \ref{cordefjyunbi}.
Let
$
{\rm Corr}_{((X,\widetriangle{\mathcal U}),\widetriangle{{\frak S}^{\epsilon}})}
:
\Omega^k(M_s) \to \Omega^{\ell+k}(M_t)
$
be the map obtained by Definition \ref{defn748}.
(Here $\ell = \dim M_t - \dim (X,\widehat{\mathcal U})$.)
We define the boundary by
$$
\partial (X,\widetriangle{\mathcal U}) = (\partial(X,\widetriangle{\mathcal U}),
\widetriangle f_s\vert_{\partial(X,\widetriangle{\mathcal U})},
\widetriangle f_t\vert_{\partial(X,\widetriangle{\mathcal U})}).
$$
$\widetriangle{{\frak S}^{\epsilon}}$
induces a CF-perturbation
$\widetriangle{{\frak S}^{{\partial},\epsilon}}$ of it
as  in
Lemma \ref{lemma755} (4).
$\partial (X,\widetriangle{\mathcal U})$ and $\widetriangle{{\frak S}^{{\partial},\epsilon}}$
define a map
$
{\rm Corr}_{(\partial(X,\widetriangle{\mathcal U}), \widetriangle{\frak S^{{\partial},\epsilon}})}
:
\Omega^k(M_s) \to \Omega^{\ell+k+1}(M_t).
$
Then for any sufficiently small $\epsilon >0$, we have
\begin{equation}
d \circ {\rm Corr}_{((X,\widetriangle{\mathcal U}),\widetriangle{{\frak S}^{\epsilon}})}
-
{\rm Corr}_{((X,\widetriangle{\mathcal U}),\widetriangle{{\frak S}^{\epsilon}})}
\circ d
=
{\rm Corr}_{(\partial(X,\widetriangle{\mathcal U}) ,\widetriangle{\frak S^{\partial,\epsilon}})}.
\end{equation}
In particular, ${\rm Corr}_{((X,\widetriangle{\mathcal U}),\widetriangle{{\frak S}^{\epsilon}})}$
is a chain map if $\widehat{\mathcal U}$ is a Kuranishi structure
without boundary.
\end{cor}
\begin{proof}
This is immediate from Theorem \ref{Stokes}.
\end{proof}
\begin{lem}\label{lem761rev}
We assume $\widetriangle{\frak S_{\sigma}}$  is a uniform
family in the sense of Definition \ref{uniformcongpert}.
Then the positive number $\epsilon$ in
Theorem \ref{Stokes} and Corollary \ref{Stokescorollary}
can be taken independent of $\sigma$.
\end{lem}
The proof is the same as that  of Lemma \ref{lem761}.

\subsection{Well-defined-ness of virtual fundamental cycle}
\label{subsec:JFCstokes}

We  use Corollary \ref{Stokescorollary} to prove well-defined-ness
of virtual cohomology class,
and well-defined-ness of the smooth correspondence {\it in the
cohomology level},
when Kuranishi structure has no boundary.

\begin{prop}\label{relextendgood}
Consider Situation \ref{smoothcorr}
and assume that our Kuranishi structure on $X$ has
no boundary.
Then the map
$
{\rm Corr}_{(\frak X,\widetriangle{{\frak S}^{\epsilon}})}
:
\Omega^k(M_s) \to \Omega^{\ell+k}(M_t)
$
defined in Definition \ref{defn748}
is a chain map.
\par
Moreover, provided $\epsilon$ is sufficiently small, the map ${\rm Corr}_{(\frak X,\widetriangle{{\frak S}^{\epsilon}})}$
is independent of the choices of
our good coordinate system ${\widetriangle{\mathcal U}}$
and CF-perturbation $\widetriangle{\frak S}$
and of $\epsilon>0$,
up to chain homotopy.
\end{prop}
\begin{proof}
The first half is repetition of Corollary \ref{Stokescorollary}.
We will prove the independence of the definition up to chain homotopy below.
Let ${\widetriangle{\mathcal U}}$,
${\widetriangle{\mathcal U'}}$ be two choices
of good coordinate system and
$\widetriangle{\frak S}$,
$\widetriangle{{\frak S}^{\prime}}$ CF-perturbations of
$(X;{\widetriangle{\mathcal U}})$,
$(X;{\widetriangle{\mathcal U'}})$ respectively.
(During the proof of Proposition \ref{relextendgood},
we do not need to discuss the choice of
support system, since the correspondence map is independent of it.)
\par
We put direct product Kuranishi structure on $X\times [0,1]$.
During the proof of Proposition \ref{relextendgood}, we do not need to make
a specific choice of
support system because Proposition \ref{indepofukuracont} (see also
Remark \ref{rem:785}) shows the map
${\rm Corr}_{(\frak X,\widetriangle{{\frak S}^{\epsilon}})}$
is independent thereof.
We identify $X = X \times \{0\}$.
Then the good coordinate system
${\widetriangle{\mathcal U}}$
induces ${\widetriangle{\mathcal U}}\times [0,1/3)$ on $X \times [0,1/3)$
such that
$\partial
(X \times [0,1/3),{\widetriangle{\mathcal U}}\times [0,1/3))
$
is isomorphic to
$
(X;{\widetriangle{\mathcal U}})
$.
Similarly we have a good coordinate system
${\widetriangle{\mathcal U'}}\times (2/3,1]$ on
$X \times (2/3,1]$ such that
$\partial(X \times (2/3,1];{\widetriangle{\mathcal U'}}\times (2/3,1])$
is isomorphic to
$
(X;{\widetriangle{\mathcal U'}})
$
with opposite orientation.
Here the notion of isomorphism of good coordinate
system is defined in Definition \ref{defn31222}.
Then, by Proposition  \ref{prop7582752},
there exists a good coordinate system
${\widetriangle{\mathcal U'}}\times [0,1]$
such that
\begin{equation}\label{boundarycobor77}
\partial
(X\times [0,1];{\widetriangle{\mathcal U'}}\times [0,1])
=
(X;{\widetriangle{\mathcal U}})
\cup
-(X,{\widetriangle{\mathcal U'}}).
\end{equation}
\par
We next consider two choices of CF-perturbations,
which we denote by $\widetriangle{\frak S}$
and $\widetriangle{{\frak S}^{\prime}}$.
We assume that $\widetriangle{f_t}$ is
strongly submersive with respect to both of them.
We define
$\widetriangle{\frak S}\times [0,1/3)$
and
$\widetriangle{{\frak S}'} \times (2/3,1]$
as follows.
\par
We consider Situation \ref{smoothcorrsingle}.
Let $\frak V_x = (V_x,\Gamma_x,E_x,\phi_x,\widehat\phi_x)$
be an orbifold chart of $(U,\mathcal E)$ and
$\mathcal S_x = (W_x,\omega_x,{\frak s}_{x}^{\epsilon})$
a CF-perturbation
on it. (Definition \ref{defn73ss}.)
Suppose $(f_t)_x$ is strongly submersive with respect to $\mathcal S_x$.
We take
$\frak V_x
\times [0,1/3) = (V_x\times [0,1/3) ,\Gamma_x,E_x\times [0,1/3) ,\phi_x\times
{\rm id},\widehat\phi_x\times {\rm id})$
that is an orbifold chart of $(U\times [0,1/3),\mathcal E\times [0,1/3))$.
In an obvious way $\mathcal S_x$ induces a CF-perturbation of it,
with respect to which $(f_t)_x \circ \pi$  is strongly submersive.
Here $\pi : X \times [0,1/3) \to X$ is the projection.
We denote it by
$\mathcal S_x\times [0,1/3)$.
(See Definition \ref{defn1022} for detail.)
\par
We perform this construction of multiplying $[0,1/3)$
for each chart
(once for each orbifold chart and once for each
Kuranishi chart) then it is fairly obvious that they are compatible
with various coordinate changes. Thus we obtain
$\widetriangle{\frak S} \times [0,1/3)$ that is
a CF-perturbation of $X \times [0,1/3)$.
We obtain $\widetriangle{{\frak S}'} \times (2/3,1]$
in the same way.
\par
Now we use Proposition \ref{existperturbcontrel}
with $Z_1 = X \times \{0,1\}$, $Z_2 = X \times [0,1]$.
Then we obtain a CF-perturbation
$\widetriangle{{\frak S}^{[0,1]}}$
of $X \times [0,1]$ such that its restriction to $X \times \{0\}$ and $X \times \{1\}$
are $\widetriangle{\frak S}$ and
$\widetriangle{{\frak S}'}$, respectively.
\par
Now we use Corollary \ref{Stokescorollary} and (\ref{boundarycobor77})
to show:
\begin{equation}\label{chomotopyrelation}
\aligned
&d\circ {\rm Corr}_{(\frak X\times [0,1],
\widetriangle{{{\frak S}^{[0,1]}}^{\epsilon}})}
+
{\rm Corr}_{(\frak X\times [0,1],\widetriangle{{{\frak S}^{[0,1]}}^{\epsilon}})}\circ d \\
&=
{\rm Corr}_{(\frak X,\widetriangle{{\frak S}^{\epsilon}})}
-
{\rm Corr}_{(\frak X,\widetriangle{{\frak S}^{ \prime \epsilon}})}.
\endaligned
\end{equation}
The independence of sufficiently small $\epsilon > 0$
follows from the following facts: For each $c > 0$ the family
$\epsilon \mapsto \widetriangle{{\frak S}^{c\epsilon}}$
is also a CF-perturbation.
\par
The proof of Proposition \ref{relextendgood} is complete.
\end{proof}
Therefore in the situation of Proposition \ref{relextendgood},
the correspondence map
${\rm Corr}_{(\frak X,\widetriangle{{\frak S}^{\epsilon}})}$
on differential forms descends to
a map on cohomology which is independent of
the choices of $\widetriangle{\mathcal U}$ and
$\widetriangle{{\frak S}^{\epsilon}}$.
We write the cohomology class as
$[{\rm Corr}_{\frak X}(h)] \in H(M_t)$ for
any closed differential form $h$ on $M_s$
by removing $\widetriangle{{\frak S}^{\epsilon}}$ from the notation.
\par
In Proposition \ref{relextendgood} we fixed our Kuranishi structure
$\widehat{\mathcal U}$ on $X$.
In fact, we can prove the same conclusion
under  milder assumption.
\begin{prop}\label{cobordisminvsmoothcor}
Let $\frak X_i = ((X_i,\widehat{\mathcal U^i}),\widehat f_s^i,\widehat f_t^i)$ be smooth correspondence from $M_s$
to $M_t$ such that
$\partial X_i = \emptyset$. Here $i=1,2$ and $M_s$, $M_t$ are independent of $i$.
We assume that there exists a smooth correspondence
$\frak Y = ((Y,\widehat{\mathcal U}),\widehat f_s,\widehat f_t)$ from $M_s$ to $M_t$ with boundary (but without corner) such that
$$
\partial \frak Y = \frak X_1 \cup -\frak X_2.
$$
Here $-\frak X_2$ is the smooth correspondence $\frak X_2$ with
opposite orientation.
Then we have
\begin{equation}\label{chomotopyrelation22}
[{\rm Corr}_{\frak X_1}(h)] = [{\rm Corr}_{\frak X_2}(h)]
\in H(M_t),
\end{equation}
where $h$ is a closed differential form on $M_s$.
\end{prop}
\begin{proof}
We take a good coordinate system $\widetriangle{\mathcal U}$ of $Y$
and a KG embedding $(Y,\widehat{\mathcal U}) \to (Y,\widetriangle{\mathcal U})$.
(Theorem \ref{Them71restate}.)
$\widehat f_t$ is pulled back from
$\widetriangle f_t : (Y,\widetriangle{\mathcal U}) \to M_t$.
$\widehat f_s$ is also pulled back from $\widetriangle f_s$
(Proposition \ref{le614} (2).)
We also obtain a CF perturbation $\widetriangle{\frak S}$ of
$\widetriangle{\mathcal U}$ with respect to which
$\widetriangle f$ is strongly submersive
(Theorem \ref{existperturbcont} (2)).
They restrict to
$(X_i,\widetriangle{\mathcal U_i})$,
$\widetriangle{\frak S_i}$
and
$\widetriangle f_t^i : (Y,\widetriangle{\mathcal U}) \to M_t$,
$\widetriangle f_s^i$.
\par
We remark that $\widetriangle f_t^i$ is strongly transversal to
$\widetriangle{\frak S_i}$.
(This is the consequence of the definition of strong transversality.
Namely we required the transversality on each of the strata
of corner structure stratification (Definition \ref{defn417}).)
\par
By Proposition \ref{relextendgood}
we can use
$(X_i,\widetriangle{\mathcal U_i})$,
$\widetriangle{\frak S_i}$,
$\widetriangle f_t^i$,
$\widetriangle f_s^i$
to define smooth correspondence
${\rm Corr}_{\frak X_i}$ (in the cohomology level.)
\par
The proposition now follows from
Corollary \ref{Stokescorollary} and (\ref{boundarycobor77})
applied to $(Y,\widetriangle{\mathcal U})$,
$\widetriangle{\frak S}$, $\widetriangle f_t$, $\widetriangle f_t$.
Namely we can calculate in the same way as
(\ref{chomotopyrelation}).
\end{proof}

\begin{rem}
The proofs of Propositions \ref{relextendgood}, \ref{cobordisminvsmoothcor}
(Formulae (\ref{chomotopyrelation}), (\ref{chomotopyrelation22})) are
a prototype of the proofs of various similar equalities
which appear in our construction of structures and
proof of its independence.
We will apply a similar method in a slightly complicated
situation in Part 2 systematically.
\end{rem}
\section{From good coordinate system to Kuranishi structure and back
with CF-perturbations}
\label{sec:kuraandgood}

As we explained at the end of Section \ref{sec:fiber}, it is more
canonical to define the notion of fiber product
of spaces with Kuranishi structure  than
to define that of  fiber product
of spaces with good coordinate system.
On the other hand, in Section \ref{sec:contfamily}, we gave the definition of
CF-perturbation and of the
pushout of differential forms by using good coordinate system.
In this section, we describe the way how we go from a good coordinate
system to a Kuranishi structure and back together with CF-perturbations on them,
and prove Theorem \ref{theorem915} that we can define the pushout
by using Kuranishi structure itself in the way that the outcome is independent of
auxiliary choice of good coordinate system.

\subsection{CF-perturbation and embedding of Kuranish structure}
\label{subsec:contfamiKura}

\begin{defn}\label{defn81}
Let $\widehat{\mathcal U}$ be a Kuranishi structure on $Z \subseteq X$.
A {\it CF-perturbation $\widehat{\frak S}$
of $\widehat{\mathcal U}$}
\index{CF-perturbation ! of Kuranishi structure}
\index{Kuranishi structure ! CF-perturbation of} assigns
$\frak S_p$ for each $p \in Z$ with the following properties.
\begin{enumerate}
\item
$\frak S_p$ is a CF-perturbation
of $\mathcal U_p$.
\item
If $q \in {\rm Im}(\psi_{p}) \cap Z$, then $\frak S_p$
can be pulled back by $\Phi_{pq}$.
Namely
\begin{equation}
\frak S_p \in \mathscr S^{\mathcal U_q \triangleright \mathcal U_p}(U_p).
\end{equation}
\item
If $q \in {\rm Im}(\psi_{p}) \cap Z$, then $\frak S_p$, $\frak S_{q}$
are  compatible with  $\Phi_{pq}$.
Namely
\begin{equation}\label{form91}
\Phi_{pq}^*(\frak S_p) = \frak S_q\vert_{U_{pq}} \in \mathscr S^{\mathcal U_q}(U_{pq}).
\end{equation}
\end{enumerate}
\end{defn}

\begin{defn}\label{defn929292}
Suppose we are in the situation of Definition \ref{defn81}. Let $\widehat f :
(X,Z;\widehat{\mathcal U}) \to M$ be a strongly smooth map.
Here $M$ is a smooth manifold.
\begin{enumerate}
\item
We say $\widehat{\frak S}$ is strictly {\it transversal to $0$} if each $\frak S_p$ is
transversal to $0$. \index{CF-perturbation ! transversal to $0$ ! on Kuranishi structure}
We say $\widehat{\frak S}$ is  {\it transversal to $0$} if its restriction
to an open substructure is strictly so.
\item
We say $\widehat f$ is strictly {\it strongly submersive with respect to $\widehat{\frak S}$}
\index{strongly submersive (w.r.t. CF-perturbation) ! map on Kuranishi structure}
if each of  $f_p$ is strongly submersive with respect to $\frak S_p$.
We say $\widehat f$ is  {\it strongly submersive with respect to $\widehat{\frak S}$}
if its restriction
to an open substructure is strictly so.
\item
We say $\widehat f$ is strictly {\it strongly transversal to $g: N \to M$ with respect to $\widehat{\frak S}$}
if each of  $f_p$ is strongly transversal to $g: N \to M$ with respect to $\frak S_p$.
We say $\widehat f$ is {\it strongly transversal to $g: N \to M$ with respect to $\widehat{\frak S}$}
\index{strongly transversal (w.r.t. CF-perturbation) ! to a map on Kuranishi structure}
if its restriction
to an open substructure is strictly so.
\end{enumerate}
\end{defn}

We next define compatibility of CF-perturbations
with various embeddings of Kuranishi structures and/or
good coordinate systems and
prove versions of several lemmata in Section \ref{sec:multisection}
corresponding to the current context of CF-perturbations.

\begin{defn}\label{defn83k}
Let $\widehat{\mathcal U}$ and $\widehat{\mathcal U^+}$ be
Kuranishi structures of $Z \subseteq X$,
${\widetriangle{\mathcal U}}$ and ${\widetriangle{\mathcal U^+}}$
good coordinate systems of $Z \subseteq X$.
Let $\mathcal K$ and $\mathcal K^+$ be support systems of
${\widetriangle{\mathcal U}}$ and ${\widetriangle{\mathcal U^+}}$,
respectively.
Let $\widehat{\frak S}$, $\widehat{\frak S^+}$,
${\widetriangle{\frak S}}$,   ${\widetriangle{\frak S^+}}$
be CF-perturbations of
$\widehat{\mathcal U}$, $\widehat{\mathcal U^+}$,
$({\widetriangle{\mathcal U}},\mathcal K)$,
$({\widetriangle{\mathcal U^+}},\mathcal K^+)$,
respectively.
\begin{enumerate}
\item
Let $\widehat\Phi : \widehat{\mathcal U} \to \widehat{\mathcal U^+}$
be a strict KK-embedding.
We say $\widehat{\frak S}$, $\widehat{\frak S^+}$ are
{\it compatible}
\index{compatibility ! of support systems with GG-embedding}
\index{embedding !compatibility of support systems with GG-embedding} with $\widehat\Phi$ if the following holds for each $p$.
\begin{enumerate}
\item
$\frak S^+_p \in \mathscr S^{\mathcal U_p \triangleright \mathcal U^+_p}(U_{p}).$
Here we use the embedding $\Phi_{p}$ to define the subsheaf $\mathscr S^{\mathcal U_p \triangleright \mathcal U^+_p}$.
\item
$
\Phi_{p}^*(\frak S^+_p) = \frak S_p \in \mathscr S^{\mathcal U_p}(U_{p})$.
\end{enumerate}
\item
Let $\widehat\Phi : \widehat{\mathcal U} \to \widehat{\mathcal U^+}$
be a KK-embedding.
We say $\widehat{\frak S}$, $\widehat{\frak S^+}$ are
{\it compatible}
\index{compatibility ! of CF-perturbation with KK-embedding}
\index{CF-perturbation ! compatibility with KK-embedding}
\index{embedding ! compatibility of CF-perturbation with KK-embedding}
with $\widehat\Phi$
if there exist an open substructure $\widehat{\mathcal U_0}$, a
CF-perturbation $\widehat{\frak S_0}$ of $\widehat{\mathcal U_0}$ and
a strict KK-embedding $\widehat\Phi_0 : \widehat{\mathcal U_0} \to \widehat{\mathcal U^+}$
such that $\widehat{\frak S_0}$, $\widehat{\frak S^+}$ are compatible with $\widehat\Phi_0$ and
$\widehat{\frak S_0}$, $\widehat{\frak S}$ are compatible with the open embedding
$\widehat{\mathcal U_0} \to \widehat{\mathcal U}$.
\item
Let $\widetriangle\Phi = (\{\Phi_{\frak p}\},\frak i) : {\widetriangle{\mathcal U}} \to {\widetriangle{\mathcal U^+}}$
be a GG-embedding.
We say that $\mathcal K, \mathcal K^+$ is {\it compatible}
\index{compatibility ! of CF-perturbation with GK-embedding}
\index{CF-perturbation ! compatibility with GK-embedding}
\index{embedding ! compatibility of CF-perturbation with GK-embedding}  with $\widehat\Phi$ if
$
\varphi_{\frak p}(\mathcal K_{\frak p}) \subset \ring\mathcal K^+_{\frak i(\frak p)}
$
for each $\frak p \in \frak P$.
\item
In the situation of (3),
we say ${\widetriangle{\frak S}}$, ${\widetriangle{\frak S^+}}$ are
{\it compatible} \index{compatibility ! of CF-perturbation with GG-embedding}
 \index{CF-perturbation ! compatibility with GG-embedding}
  \index{embedding ! compatibility of CF-perturbation with GG-embedding}  with $\widehat\Phi$ if
the following holds for each $\frak p \in \frak P$.
\begin{enumerate}
\item
${\frak S}^+_{\frak i(\frak p)} \in \mathscr S^{\mathcal U_{\frak p} \triangleright \mathcal U^+_{\frak i(\frak p)}}(\mathcal K_{{\frak i(\frak p)}}).$
Here we use the embedding $\Phi_{\frak p} : \mathcal U_{\frak p} \to \mathcal U^+_{\frak i(\frak p)}$ to define
the subsheaf $\mathscr S^{\mathcal U_{\frak p} \triangleright \mathcal U^+_{\frak i(\frak p)}}$.
\item
$
\Phi_{\frak p}^*(\frak S^+_{\frak i(\frak p)}) = \frak S_{\frak p} \in \mathscr S^{\mathcal U_{\frak p}}(\mathcal K_{\frak p})$.
\end{enumerate}
\item
Let $\widehat\Phi : {\widehat{\mathcal U}} \to {\widetriangle{\mathcal U}}$
be a strict KG-embedding.
We say ${\widehat{\frak S}}$, ${\widetriangle{\frak S}}$ are
{\it compatible}
\index{compatibility ! of CF-perturbation with KG-embedding}
\index{CF-perturbation ! compatibility with KG-embedding}
\index{embedding ! compatibility of CF-perturbation with KG-embedding} with $\widehat\Phi$ if  the following holds for each $\frak p$
and $p \in \psi_{\frak p}(\mathcal K_{\frak p}\cap s_{\frak p}^{-1}(0)) \cap Z$.
\begin{enumerate}
\item
${\frak S}_{\frak p} \in \mathscr S^{\mathcal U_{p} \triangleright \mathcal U_{\frak p}}(\mathcal K_{\frak p}).$
Here we use the embedding $\Phi_{\frak p p} : \mathcal U_p \to \mathcal U_{\frak p}$
to define the subsheaf $ \mathscr S^{\mathcal U_{p} \triangleright \mathcal U_{\frak p}}$.
\item
$\Phi_{p}^*({\frak S}_{\frak p}) = {\frak S}_{p} \in \mathscr S^{\mathcal U_{p}}(U_{p})$.
\end{enumerate}
\item
In case $\widehat\Phi : {\widehat{\mathcal U}} \to {\widetriangle{\mathcal U}}$
is a KG-embedding, we can define compatibility of ${\widehat{\frak S}}$, ${\widetriangle{\frak S}}$
with $\widehat\Phi$ in the same way as Item (2) (using Items (1) and (5)).
\item
Let $\widehat\Phi
= (\{U_{\frak p}(p)\},\{
\Phi_{p \frak p}\}) : {\widetriangle{\mathcal U}} \to \widehat{\mathcal U}$
be a GK-embedding.
We say ${\widetriangle{\frak S}}$,
$\widehat{\frak S}$ are
{\it compatible} \index{compatibility ! of CF-perturbation with GK-embedding}
\index{CF-perturbation ! compatibility with GK-embedding}
\index{embedding ! compatibility of CF-perturbation with GK-embedding}  with $\widehat\Phi$ if the following holds for each
$\frak p$ and $p \in \psi_{\frak p}(\mathcal K_{\frak p} \cap s_{\frak p}^{-1}(0)) \cap Z$.
\begin{enumerate}
\item
${\frak S}_{p} \in \mathscr S^{\mathcal U_{\frak p} \triangleright \mathcal U_{p}}(U_{p}).$
Here we use the embedding $\Phi_{p \frak p} : \mathcal U_{\frak p}\vert_{U_{\frak p}(p)}
\to \mathcal U_p$ to define the subsheaf $\mathcal U_{\frak p} \triangleright \mathcal U_{p}$.
\item
$\Phi_{p\frak p}^*{\frak S}_p = \widehat{\frak S}_{\frak p}\vert_{U_{\frak p}(p)}
\in \mathscr S^{\mathcal U_{p}}(U_{\frak p}(p))$.
\end{enumerate}
\end{enumerate}
\end{defn}
With these definitions of compatibility, we now prove the
compatibilities relevant to various embeddings.
\begin{lem}\label{lem9494}
Let $\widetriangle\Phi = (\{\Phi_{\frak p}\},\frak i) : \widetriangle{\mathcal U} \to \widetriangle{\mathcal U^+}$
be a GG-embedding.
\begin{enumerate}
\item
If $\mathcal K$ is a support system of $\widetriangle{\mathcal U}$,
then there exist a support system $\mathcal K^+$ of $\widetriangle{\mathcal U^+}$
such that $\mathcal K$, $\mathcal K^+$ are compatible with $\widetriangle\Phi$.
\item
If $\mathcal K_i$ ($i=1,\dots,m$) are support systems of $\widetriangle{\mathcal U}$ with $\mathcal K_i < \mathcal K_{i+1}$
then there exist support systems $\mathcal K^+_i$ ($i=1,\dots,m$) of $\widetriangle{\mathcal U^+}$
such that $\mathcal K_i$, $\mathcal K_i^+$ are compatible with $\widetriangle\Phi$ and $\mathcal K^+_i < \mathcal K^+_{i+1}$.
\end{enumerate}
\end{lem}
\begin{proof}
(1)  Let $\mathcal K = (\mathcal K_{\frak p})$.
Let
$
\mathcal K_{0,\frak p_+}$ be a
closure of a sufficiently small neighborhood of $\bigcup_{\frak p \in \frak P \atop
\frak i(\frak p) = \frak p_+} \varphi_{\frak p}(\mathcal K_{\frak p})
$
for $\frak p_+ \in \frak P_+$.
It is easy to see that $\mathcal K_{0}^+ = (\mathcal K_{0,\frak p_+})$ is a support system of $\widetriangle{\mathcal U^+}$.
Any $\mathcal K^+ > \mathcal K_{0}^+$ has required properties.
The proof of (2) is similar by using upward induction on $i$.
\end{proof}
\begin{rem}
It seems possible to prove the following.
For each $\mathcal K^+$ there exists $\mathcal K$ such that
$\mathcal K$, $\mathcal K^+$ are compatible with $\widetriangle\Phi$.
Its proof seems to be more complicated than that of Lemma \ref{lem9494}.
We do not try to prove it here since we do not use it.
\end{rem}
\begin{lem}\label{compatiwithwkemb}
Let $\widetriangle{\Phi} : {\widetriangle{\mathcal U}} \to {\widetriangle{\mathcal U^+}}$
be a weakly open GG-embedding
and $\mathcal K$, $\mathcal K^+$ support systems of
${\widetriangle{\mathcal U}}$, ${\widetriangle{\mathcal U^+}}$, respectively, which are compatible with $\widetriangle\Phi$.
Then for any CF-perturbation
$\widetriangle{\frak S^+}$ of $({\widetriangle{\mathcal U^+}},\mathcal K^+)$,
there exists a unique CF-perturbation $\widetriangle{\frak S}$
of $({\widetriangle{\mathcal U}},\mathcal K)$ such that
$\widetriangle{\frak S^+}$ and $\widetriangle{\frak S}$
are compatible with $\widetriangle{\Phi}$.
\end{lem}
\begin{proof}
For any $\frak p \in \frak P$ we restrict ${\frak S_{\frak i(\frak p)}^+}$
to ${\mathcal U}_{\frak p}$ to obtain ${\frak S}_{\frak p}$.
We thus obtain $\widetriangle{\frak S}$ .
Since normal bundles are trivial in the case of weakly open embedding, the compatibility is
automatic.
\end{proof}
\begin{lem}\label{lem9797}
Various transversality or submersivity of the target of
an open KK-embedding imply those of the source.
The same holds for a weakly open GG-embedding.
\end{lem}
\begin{proof}
This is an easy consequence of the definition.
\end{proof}

\begin{lem}
The notion of compatibility of CF-perturbations to embeddings is
preserved by the composition of embeddings
of various kinds.
\end{lem}
The proof is obvious.
\par
The next lemma is a CF-perturbation version of
Proposition \ref{lemappgcstoKu}.

\begin{lem}\label{lemappgcstoKucont}
In the situation of Proposition \ref{lemappgcstoKu},
let $\mathcal K_0$, $\mathcal K$ be support systems of
$\widetriangle{\mathcal U_0}$, ${\widetriangle{\mathcal U}}$, respectively,
which are compatible with the open embedding $\widetriangle{\mathcal U_0}\to {\widetriangle{\mathcal U}}$.
Let $\widetriangle{\frak S}$ be a CF-perturbation of $({\widetriangle{\mathcal U}},\mathcal K)$,
which restricts  to a CF-perturbation
${\widetriangle{\frak S_0}}$
of $({\widetriangle{\mathcal U_0}},\mathcal K_0)$.
Then the following holds.
\begin{enumerate}
\item
There exists a CF-perturbation
$\widehat{\frak S}$ of $\widehat{\mathcal U}$ such that
$\widetriangle{\frak S_0}$ and $\widehat{\frak S}$ are compatible with
the GK-embedding ${\widetriangle{\mathcal U_0}} \to \widehat{\mathcal U}$.
\item
In the situation of Proposition \ref{lemappgcstoKu} (2),
if ${\widetriangle f}$ is strongly submersive with
respect to ${\widetriangle{\frak S}}$,
then ${\widehat f}$ is strongly submersive
with respect to $\widehat{\frak S}$.
The transversality to $M \to Y$ is also preserved.
\end{enumerate}
\end{lem}
\begin{proof}
The proof of Lemma \ref{lemappgcstoKucont} is the same as that of
Proposition \ref{lemappgcstoKu}.
In fact, the Kuranishi chart of $\widehat{\mathcal U}$ is a restriction of a Kuranishi chart of
${\widetriangle{\mathcal U_0}}$. Since $\mathcal U_{0,\frak p} \subset \mathcal K_{\frak p}$
(see Proposition \ref{lemappgcstoKu} (1)), we can restrict ${\widetriangle{\frak S}}$
to the Kuranishi charts of  $\widehat{\mathcal U}$.
\end{proof}
We next state CF-perturbation versions of Propositions \ref{le614} and \ref{pro616}.
(In Lemmas \ref{le714} and \ref{le7155} we do not specify support system
for the CF-perturbations of good coordinate system.
We take one but do not mention them.)

\begin{lem}\label{le714}
Let $\widehat{\mathcal U}$ be a Kuranishi structure on $Z \subseteq X$.
Then we can take a good coordinate system $\widetriangle{\mathcal U}$
and the KG-embedding $\widehat{\Phi} : \widehat{\mathcal U_0} \to \widetriangle{\mathcal U}$
in Theorem \ref{Them71restate} so that the following holds in addition.
\begin{enumerate}
\item
If $\widehat h$ is a differential form of $\widehat{\mathcal U}$, then
there exists a differential form ${\widetriangle h}$ on
${{\widetriangle {\mathcal U}}}$ such that
$\widehat{\Phi}^*({\widetriangle h}) = \widehat h\vert_{\widehat{\mathcal U_0}}$.
If $\widehat h$ has a compact support in $\ring Z$, then  ${\widetriangle h}$
has a compact support in $\vert\widetriangle {\mathcal U}\vert$
and ${\rm Supp}({\widetriangle h}) \cap Z \subset \overset{\circ} Z$.
\item
If $\widehat{\frak S}$
a CF-perturbation of ${\widehat{\mathcal U}}$, then
there exists a CF-perturbation ${\widetriangle{\frak S}}$
of ${\widetriangle{\mathcal U}}$
such that
${\widehat{\frak S}}\vert_{\widehat{\mathcal U_0}}$ and
${\widetriangle{\frak S}}$
are compatible with the KG-embedding $\widehat{\Phi}$.
\item
In the situation of (2) the following holds.
\begin{enumerate}
\item
If $\widehat{\frak S}$
is transversal to $0$ then
so is ${\widetriangle{\frak S}}$.
\item
If ${\widehat f}$ is strongly submersive
with respect to $\widehat{\frak S}$, then
${\widetriangle f}$ is strongly submersive with
respect to ${\widetriangle{\frak S}}$.
\item
If
${\widehat f}$ is strongly transversal to $g : M \to Y$
with respect to ${\widehat{\frak S}}$
then ${\widetriangle f}$ is strongly transversal to $g$ with respect to ${\widetriangle{\frak S}}$.
\end{enumerate}
\end{enumerate}
\end{lem}
\begin{lem}\label{le7155}
Suppose we are in the situation of Propositions \ref{prop518}
(resp. Proposition \ref{prop519})
and \ref{le614}.
Then we can take the GK-embedding
$\widehat{\Phi^+} :  \widetriangle{\mathcal U} \to
\widehat{\mathcal U^+}$ in Proposition \ref{prop518}
(resp. the GK-embeddings $\widehat{\Phi^+_a} :  \widetriangle{\mathcal U} \to
\widehat{\mathcal U^+_a}$ in Proposition \ref{prop519}
($a=1,2$))
so that the following holds.
\begin{enumerate}
\item
If  ${\widehat{\frak S^+}}$ is a CF-perturbation of ${\widehat{\mathcal U^+}}$
such that ${\widehat{\frak S^+}}$, ${\widehat{\frak S}}$ are
strongly compatible
with the embedding ${\widehat{\mathcal U}} \to \widehat{\mathcal U^+}$,
then we may choose $\widetriangle{\frak S}$
such that  $\widetriangle{\frak S}$, $\widehat{\frak S^+}$ are
compatible
with the embedding $\widehat{\Phi^+}$.
(resp.
If  ${\widehat{\frak S^+_a}}$ ($a=1,2$) is a CF-perturbation of ${\widehat{\mathcal U^+_a}}$
such that ${\widehat{\frak S^+_a}}$, ${\widehat{\frak S}}$ are
strongly compatible
with the embedding ${\widehat{\mathcal U}} \to \widehat{\mathcal U^+}$,
then we may choose $\widetriangle{\frak S}$
such that  $\widetriangle{\frak S}$, $\widehat{\frak S^+_a}$ are both
compatible
with the embedding $\widehat{\Phi^+_a}$.)
\item
If $\widehat{\frak S}$
is transversal to $0$, so is $\widetriangle{\frak S}$.
\item\
In the situation of Proposition \ref{pro616} suppose $Y$ is a manifold $M$. Then
if ${\widehat f}$  is strongly submersive with respect to ${\widehat{\frak S}}$,
then ${\widetriangle{f}}$  is strongly submersive with respect to $\widetriangle{\frak S}$.
\item
In the situation of Proposition \ref{pro616} suppose $M$ is a manifold. Then
if ${\widehat f}$ is strongly transversal to $g : N \to M$ with respect to ${\widehat{\frak S}}$, then
${\widetriangle{f}}$  is strongly transversal to $g : N \to M$ with respect to $\widetriangle{\frak S}$.
\end{enumerate}
\end{lem}

The proofs of Lemmata \ref{le714} and \ref{le7155} are given in Subsection \ref{subsec:movingmulsectionetc}.
\subsection{Integration along the fiber (pushout) for Kuranishi structure}
\label{subsec:intKurast}

\begin{shitu}\label{sitsu8main}
Let $\widehat{\mathcal U}$ be a Kuranishi structure on $X$ and
${\widehat{\frak S}}$ a CF-perturbation of $(X,Z;\widehat{\mathcal U})$.
Let $\widehat f : (X,Z;\widehat{\mathcal U}) \to M$
be a strongly smooth map that
is strongly submersive with respect to ${\widehat{\frak S}}$.
Let $\widehat h$ be a differential form on $\widehat{\mathcal U}$.
By Lemma
\ref{le714}, we obtain
${\widetriangle{\mathcal U}}$,
$\widehat{\Phi}$,
${\widetriangle{\frak S}}$,
${\widetriangle f}$,
${\widetriangle{ h}}$.$\blacksquare$
\end{shitu}
\begin{defn}\label{deflemgg}\index{pushout ! {\it see integration along the fiber}}
In Situation \ref{sitsu8main}, we define the {\it pushout}, or
the {\it integration along the fiber}
$\widehat f !(\widehat{ h};\widehat{{\frak S}^{\epsilon}})$ by
\begin{equation}\label{form9393}
\widehat f !(\widehat{ h};\widehat{{\frak S}^{\epsilon}})
=
{\widetriangle f} !\left({\widetriangle{h}};{\widetriangle{\frak S^{\epsilon}}}\right).
\end{equation}
Here the right hand side is defined in Definition \ref{pushforwardKuranishi}.
Hereafter we mostly use the terminology `pushout' in this document.
\end{defn}
\begin{thm}\label{theorem915}
The right hand side of (\ref{form9393}) is independent of choices of
${\widetriangle{\mathcal U}}$,
$\widehat{\Phi}$,
${\widetriangle{\frak S}}$,
${\widetriangle f}$,
${\widetriangle{ h}}$
in the sense of $\spadesuit$ of Definition \ref{defnspadesuit},
but depends only on ${\widehat{\mathcal U}}$,
${\widehat{\frak S}}$,
${\widehat f}$,
${\widehat{ h}}$ and $\epsilon$.
\end{thm}
The proof uses Proposition \ref{integralinvembprop}.
To state it we consider the following situation.
\begin{shitu}\label{situ9016}
Let ${\widetriangle{\mathcal U}}$, ${\widetriangle{\mathcal U^+}}$
be good coordinate systems of $Z \subseteq X$,
$\widetriangle\Phi :
{\widetriangle{\mathcal U}} \to {\widetriangle{\mathcal U^+}}$
a GG-embedding, and
$\mathcal K$, $\mathcal K^+$  the  respective
support systems of ${\widetriangle{\mathcal U}}$, ${\widetriangle{\mathcal U^+}}$  compatible with $\widetriangle\Phi$.
Let
${\widetriangle{\frak S}}$,   ${\widetriangle{\frak S^+}}$
be CF-perturbations of
$(\widetriangle{\mathcal U},\mathcal K)$,
$(\widetriangle{\mathcal U^+},\mathcal K^+)$,
respectively.
Let ${\widetriangle{h^+}}$ be a differential form on ${\widetriangle{\frak U^+}}$ which has
a compact support in $\ring Z$
and  $\widetriangle {f^+} : (X,Z;{\widetriangle{\mathcal U^+}}) \to M$
a strongly smooth map.
We put $\widetriangle h = \widetriangle\Phi^* \widetriangle{h^+}$
and $\widetriangle f = \widetriangle {f^+} \circ \widetriangle\Phi : (X,Z;{\widetriangle{\mathcal U}}) \to M$.
\par
We assume that  $\widetriangle f$ is strongly submersive with respect to $\widetriangle{\frak S}$ and
 $\widetriangle{f^+}$ is strongly submersive with respect to $\widetriangle{{\frak S}^+}$.
$\blacksquare$
\end{shitu}
\begin{prop}\label{integralinvembprop}
In Situation \ref{situ9016} we have
\begin{equation}
{\widetriangle f} !\left( {\widetriangle{h}};{\widetriangle{\frak S^{\epsilon}}}
\right)
=
{\widetriangle{f^+}} !\left({\widetriangle{h^+}};{\widetriangle{\frak S^{+\epsilon}}}
\right)
\end{equation}
for each sufficiently small $\epsilon >0$.
\end{prop}
\begin{proof}[Proof of Proposition \ref{integralinvembprop} $\Rightarrow$
Theorem \ref{theorem915}]
We use Definition-Lemma \ref{henacomp} also in this proof.
\begin{equation}
\xymatrix{
& \widetriangle{\mathcal U_1} \ar[rr]  && \widehat{\mathcal U_1^+} \ar[rr]  &&\widetriangle{\mathcal U_1^{+}} \\
\widehat{\mathcal U}     \ar[ru] \ar[rr]\ar[rd] && \widetriangle{\mathcal U_3} \ar[ru]\ar[rd] \\
& \widetriangle{\mathcal U_2} \ar[rr]  && \widehat{\mathcal U_2^+} \ar[rr]  &&\widetriangle{\mathcal U_2^{+}}
}
\end{equation}
Let
${\widetriangle{\mathcal U_i}}$,
$\widehat{\Phi_i}$,
${\widetriangle{\frak S_i}}$,
${\widetriangle f_i}$,
${\widetriangle{h_i}}$, $i=1,2$ be two choices.
By Lemma \ref{le7155} we have
${\widehat{\mathcal U_i^+}}$,
${\widehat{\frak S_i^+}}$,
${\widehat f_i^+}$,
${\widehat{h_i^+}}$, $i=1,2$
and GK-embeddings $\widehat{\Phi^+_i} : {\widetriangle{\mathcal U_i}} \to {\widehat{\mathcal U_i^+}}$
to which various objects are compatible.
\par
By Lemma \ref{le7155} we obtain
${\widetriangle{\mathcal U_3}}$,
$\widehat{\Phi_3}$,
${\widetriangle{\frak S_3}}$,
${\widetriangle f_3}$,
${\widetriangle{h_3}}$
and GK-embeddings
$\widehat{\Phi^{+-}_i} : {\widetriangle{\mathcal U_3}} \to {\widehat{\mathcal U_i^+}}$
to which various objects are compatible.
\par
By Lemma \ref{le714}, we obtain
${\widetriangle{\mathcal U_i^{+}}}$,
$\widehat{\Phi_i^{+}}$,
${\widetriangle{\frak S_i^{+}}}$,
${\widetriangle f_i^{+}}$,
${\widetriangle{ h_i^{+}}}$
and KG-embeddings
$\widehat{\Phi^{+}_i} : {\widehat{\mathcal U^+_i}} \to {\widetriangle{\mathcal U_i^{+}}}$
to which various objects are compatible.
\par
Now we claim
\begin{equation}\label{formula9696}
{\widetriangle{f_1}} !\left({\widetriangle{h_1}};{\widetriangle{\frak S_1^{\epsilon}}}
\right)
=
{\widetriangle{f_1^{+}}} !\left({\widetriangle{h_1^{+}}};{\widetriangle{\frak S_1^{+\epsilon}}}
\right).
\end{equation}
In fact by Definition-Lemma \ref{henacomp}
there exists a weakly open substructure ${\widetriangle{\mathcal U_{0,1}}}$ of
${\widetriangle{\mathcal U_{1}}}$ and a GG-embedding
${\widetriangle{\mathcal U_{0,1}}} \to {\widetriangle{\mathcal U_1^{+}}}$.
By Lemma \ref{compatiwithwkemb} we can restrict ${\widetriangle{\frak S_1}}$
to ${\widetriangle{\frak S_{0,1}}}$, as well as other objects.
Strong submersivity is preserved by Lemma \ref{lem9797}.
Therefore by Proposition \ref{integralinvembprop} we find
$$
{\widetriangle{f_1}} !\left({\widetriangle{h_1}};{\widetriangle{\frak S_1^{\epsilon}}}
\right)
=
{\widetriangle{f_{0,1}}} !\left({\widetriangle{h_{0,1}}};{\widetriangle{\frak S_{0,1}^{\epsilon}}}
\right)
=
{\widetriangle{f_1^{+}}} !\left({\widetriangle{h_1^{+}}};{\widetriangle{\frak S_1^{+\epsilon}}}
\right).
$$
Here $\widetriangle{h_{0,1}}$ is the pull back of $\widetriangle{h_{1}}$
to ${\widetriangle{\mathcal U_{0,1}}}$.

We have thus proved (\ref{formula9696}).
Using the same argument three more times, we obtain
$$
\aligned
{\widetriangle{f_1}} !\left({\widetriangle{h_1}};{\widetriangle{\frak S_1^{\epsilon}}}
\right)
&=
{\widetriangle{f_1^{+}}} !\left({\widetriangle{h_1^{+}}};{\widetriangle{\frak S_1^{+\epsilon}}}
\right)=
{\widetriangle{f_3}} !\left({\widetriangle{h_3}};{\widetriangle{\frak S_3^{\epsilon}}}\right)  \\
&=
{\widetriangle{f_2^{+}}} !\left({\widetriangle{h_2^{+}}};{\widetriangle{\frak S_2^{+\epsilon}}}
\right) =
{\widetriangle{f_2}} !\left({\widetriangle{h_2}};{\widetriangle{\frak S_2^{\epsilon}}}\right).
\endaligned
$$
We have thus proved the required independence of $\widehat{\Phi}$,
${\widetriangle{\frak S}}$,
${\widetriangle f}$,
${\widetriangle{ h}}$.
\end{proof}

\subsection{Proof of Definition-Lemma \ref{henacomp}}
\label{subsec:proofofwdofcomp}

\begin{proof}[Proof of Definition-Lemma \ref{henacomp}]
Recalling the notation for GK-embedding in Definition \ref{embgoodtokura},
we put
$$
\widehat{\Phi} = \{(U_{\frak p}(p),\Phi_{p\frak p})\} ~:~
\widetriangle{\mathcal U} \to \widehat{\mathcal U}
$$
and a KG-embedding $\widehat{\Phi^+}  :
\widehat{\mathcal U} \to {\widetriangle{\mathcal U^+}}$.
We take a support system $\mathcal K$ of
${\widetriangle{\mathcal U}}$ and
$\mathcal K^{+}$ of
${\widetriangle{\mathcal U^+}}$, respectively.
For $\frak p \in \frak P$, $\frak q \in \frak P^+$ we define
\begin{equation}\label{defXpq}
Z_{\frak p\frak q} = (\mathcal K_{\frak p} \cap Z)
\cap  (\mathcal K^{+}_{\frak q} \cap Z).
\end{equation}
Here and hereafter the set theoretical symbols such as equality and the intersection in (\ref{defXpq}) etc.. are regarded as those among the
subsets of
$\vert{\widetriangle {\mathcal U^+}}\vert$.
\par
We will use Lemma \ref{lem816} to obtain the partial ordered set $\frak P_0$
which is a part of weakly open substructure $\widetriangle{\mathcal U^0}$
of $\widetriangle{\mathcal U}$.
(In this subsection we write  $\widetriangle{\mathcal U^0}$ in place of
$\widetriangle{\mathcal U_0}$.)

\begin{lem}\label{lem816}
There exist a finite subset $A_{\frak p\frak q}$ of $Z_{\frak p\frak q}$
for each $\frak p \in \frak P$, $\frak q \in \frak P^+$
and a subset $U_{(\frak p,p)}$ of $U_{\frak p}$
for each $\frak p$ and $p \in A_{\frak p\frak q}$ such that
they have the following properties.
\begin{enumerate}
\item
$p \in U_{(\frak p,p)}$. $U_{(\frak p,p)}$ is an open subset of $U_{\frak p}$.
\item
$U_{(\frak p,p)} \subset U_{\frak p}(p)$.
\item
If $p \in A_{\frak p\frak q}$, $p' \in A_{\frak p'\frak q'}$,
$\frak p \le \frak p'$ and
$$
\varphi_{\frak p'\frak p}^{-1}(U_{(\frak p',p')}) \cap U_{(\frak p,p)}
\ne \emptyset,
$$
then $\frak q\le \frak q'$.
\item
For each $\frak p_0 \in \frak P$, $\frak q_0 \in \frak P^+$ we have
$$
\bigcup_{\frak p, \frak q :
\frak p_0\le \frak p, \frak q_0 \le \frak q,
\atop p\in A_{\frak p\frak q}}
\left(
U_{(\frak p,p)}\cap Z
\right)
\supseteq Z_{\frak p_0\frak q_0}.
$$
\item
If $(\frak p,\frak q) \ne (\frak p',\frak q')$ then
$A_{\frak p\frak q} \cap A_{\frak p'\frak q'} = \emptyset$.
\end{enumerate}
\end{lem}
We recall:
\begin{defn}
Let $(\frak P,\le)$ be a partially ordered set.
A subset $\frak I \subseteq \frak P$ is said to be an {\it ideal}
\index{ideal} if $\frak p \in \frak I$,
$\frak p' \ge \frak p$ implies $\frak p' \in \frak I$.
\end{defn}
\begin{proof}
We define a partial order on $\frak P \times \frak P^+$ such that
$(\frak p,\frak q) \le (\frak p',\frak q') $ if and only if
`$\frak p\le \frak p'$'  $\wedge$ `$\frak q\le \frak q'$'.
(Note if $\frak p < \frak p'$ and $\frak q >\frak q'$, neither
$(\frak p,\frak q) \le (\frak p',\frak q') $ nor
$(\frak p,\frak q) \ge (\frak p',\frak q') $ hold.)
Let $\frak I \subset \frak P \times \frak P^+$ be an ideal.
We will prove the following by induction on
$\#\frak I$.
\begin{sublem}\label{setsekigomasubme}
For each $(\frak p,\frak q) \in \frak I$
there exist  a finite subset $A_{\frak p\frak q}$ of $Z_{\frak p\frak q}$
and a subset $U_{(\frak p,p)}$ of $U_{\frak p}$
for each $p \in A_{\frak p\frak q}$ such that
they satisfy (1)(2)(5) of Lemma \ref{lem816} and the following
conditions (3)' and (4)'.
\begin{enumerate}
\item[(3)']
\begin{enumerate}
\item
If $(\frak p,\frak q), (\frak p',\frak q') \in \frak I$, then
 Lemma \ref{lem816} (3) holds.
\item
If $(\frak p,\frak q) \in \frak I$, $p \in A_{\frak p\frak q}$
and $(\frak p',\frak q') \in \frak P\times \frak P^+$ satisfies
$$
\overline{U_{(\frak p,p)}} \cap Z_{\frak p'\frak q'} \ne \emptyset,
$$
then $(\frak p,\frak q) \ge (\frak p',\frak q')$.
\end{enumerate}
\item[(4)']
For each $(\frak p_0,\frak q_0) \in \frak I$ we have
$$
\bigcup_{\frak p, \frak q :
\frak p_0\le \frak p, \frak q_0 \le \frak q,
\atop p\in A_{\frak p\frak q}}
\left(
U_{(\frak p,p)}\cap Z
\right)
\supseteq Z_{\frak p_0\frak q_0}.
$$
\end{enumerate}
\end{sublem}
Note that we do not assume $(\frak p',\frak q') \in \frak I$ in Sublemma
\ref{setsekigomasubme} (3)' (b).
\begin{proof}
The case $\frak I = \emptyset$ is trivial.
Suppose Sublemma \ref{setsekigomasubme} is proved for all
$\frak I'$ with
$\#\frak I' < \#\frak I$. We will prove the
case of $\frak I$.
\par
Let $(\frak p_1,\frak q_1)$ be a minimal element of
$\frak I$. Then
$\frak I_- =\frak I \setminus \{(\frak p_1,\frak q_1)\}$
is an ideal of $\frak P \times \frak P^+$.
By induction hypothesis, we obtain $A_{\frak p\frak q}$
for $(\frak p,\frak q) \in \frak I_-$
and $U_{(\frak p,p)}$
for each $p \in A_{\frak p\frak q}$, $(\frak p,\frak q) \in \frak I_-$.
By induction hypothesis, Sublemma \ref{setsekigomasubme} (4)',
the set
$$
O
= \bigcup_{(\frak p, \frak q) \in \frak I :
(\frak p_1,\frak q_1) < (\frak p,\frak q)
\atop p\in A_{\frak p\frak q}}
\left(
U_{(\frak p,p)}\cap Z_{\frak p_1\frak q_1}
\right)
$$
is an open neighborhood of
$$
L = \left(\bigcup_{(\frak p, \frak q) \in \frak I :
(\frak p_1,\frak q_1) < (\frak p,\frak q)
}
Z_{\frak p\frak q}\right) \cap Z_{\frak p_1\frak q_1}
$$
in $Z_{\frak p_1\frak q_1}$.
\begin{subsublemma}\label{sublem8199}
If $x \in Z_{\frak p_1\frak q_1} \setminus O$
and $x \in Z_{\frak p \frak q}$, then
$(\frak p,\frak q) \le (\frak p_1,\frak q_1)$.
\end{subsublemma}
Note that we do {\it not} assume $(\frak p,\frak q) \in \frak I$.
\begin{proof}
Since
$$
x \in \mathcal K_{\frak p} \cap \mathcal K_{\frak p_1}
\cap
\mathcal K^+_{\frak q} \cap \mathcal K^+_{\frak q_1} \cap Z,
$$
Definition \ref{gcsystem} (5) implies that
`$\frak p \le \frak p_1$ or $\frak p \ge \frak p_1$'
holds and
`$\frak q \le \frak q_1$ or $\frak q \ge \frak q_1$'
holds also.
\par
Suppose $\frak p >\frak p_1$. Then we claim
$(\frak p,\frak q)
> (\frak p_1,\frak q_1)$ can not occur.
In fact, if $(\frak p,\frak q)
> (\frak p_1,\frak q_1)$, then
$(\frak p,\frak q) \in \frak I_-$ because $\frak I$ is an ideal.
This contradicts to $x \in Z_{\frak p_1\frak q_1} \setminus O$.
(We use the induction hypothesis  Sublemma \ref{setsekigomasubme} (4)'
here.)
Therefore $\frak q < \frak q_1$
must hold.
Then
$x \in \mathcal K_{\frak p}  \cap \mathcal K^+_{\frak q_1} \cap Z$
and $(\frak p,\frak q_1) > (\frak p_1,\frak q_1) $.
This contradicts $x \notin O$.
We can find a contradiction from $\frak q > \frak q_1$
in a similar way.
Therefore
we obtain $(\frak p,\frak q) \le (\frak p_1,\frak q_1)$.
\end{proof}
\begin{subsublemma}\label{propety818}
For each $x \in Z_{\frak p_1\frak q_1} \setminus O$, there exists
its neighborhood $W_x$ in $U_{\frak p_1}$ with the following
properties.
\begin{enumerate}
\item
$x \in W_x$ and $W_x$ is open in $U_{\frak p_1}$.
\item
$W_x \subset U_{\frak p_1}(x)$.
\item
If $W_x \cap Z_{\frak p\frak q} \ne \emptyset$ then $(\frak p,\frak q) \le (\frak p_1,\frak q_1)$.
\item
If $\frak p \ge \frak p_1$, $p \in A_{\frak p\frak q}$, $(\frak p,\frak q) \in \frak I_-$
and
$W_x \cap \varphi_{\frak p \frak p_1}^{-1}({U_{(\frak p,p)}}) \ne \emptyset$, then $\frak q \ge \frak q_1$.
\end{enumerate}
\end{subsublemma}
\begin{proof}
Since $Z_{\frak p \frak q}$ is a closed set, Subsublemma \ref{sublem8199} implies
that (3) holds for a sufficiently small neighborhood $W_x$ of $x$.
\par
We next prove that (4)
holds for a sufficiently small neighborhood $W_x$ of $x$.
Suppose $\frak p \ge \frak p_1$, $p \in A_{\frak p\frak q}$, $(\frak p,\frak q) \in \frak I_-$
and $x \in \overline{U_{(\frak p,p)}}$.
Then
$
x \in \overline{U_{(\frak p,p)}} \cap  Z_{\frak p_1\frak q_1}.
$
We apply the induction hypothesis   Sublemma \ref{setsekigomasubme} (3)' (b)
to $\frak I_-$ and  find $(\frak p,\frak q) \ge (\frak p_1,\frak q_1)$.
In particular, $\frak q \ge \frak q_1$.
Then we can take a sufficiently small neighborhood $W_x$ so that
$$
\frak p \ge \frak p_1, p \in A_{\frak p\frak q}, (\frak p,\frak q) \in \frak I_-,
W_x \cap \overline{\varphi_{\frak p p}^{-1}(U_{(\frak p,p)})}
\ne \emptyset
\,\,\Rightarrow\,\,
\frak q \ge \frak q_1.
$$
Since $W_x$ is open,
the condition  $W_x \cap \varphi_{\frak p \frak p_1}^{-1}(\overline{U_{(\frak p,q)}}) \ne \emptyset$
is equivalent to the condition
$W_x \cap \varphi_{\frak p \frak p_1}^{-1}({U_{(\frak p,q)}}) \ne \emptyset$.
Thus we have proved (4).
\end{proof}
We take an open neighborhood $W^0_x$ of $x$ such that
$\overline{W^0_x} \subset W_x$.
We take a finite subset $A_{\frak p_1\frak q_1} \subset Z_{\frak p_1\frak q_1} \setminus O$
such that
\begin{equation}\label{eq8333}
Z_{\frak p_1\frak q_1} \setminus O \subset
\bigcup_{x \in A_{\frak p_1\frak q_1}} W^0_x.
\end{equation}
Lemma \ref{lem816} (5) is obvious from definition.
\par
For $x \in A_{\frak p_1\frak q_1}$, we put
\begin{equation}\label{eq09999}
U_{(\frak p_1,x)} = W^0_x.
\end{equation}
\begin{subsublemma}\label{cond819}
There exists an open neighborhood
$U'_{(\frak p,p)}$ of $p$ for $(\frak p,\frak q) \in \frak I_-$ and
$p \in A_{\frak p\frak q}$ such that the following holds.
\begin{enumerate}
\item
$U'_{(\frak p,p)} \subset U_{(\frak p,p)}$.
\item
Sublemma \ref{setsekigomasubme} (1)(2)(4)' hold for $U'_{(\frak p,p)}$.
\item
If $p \in A_{\frak p\frak q}$, $(\frak p,\frak q) \in \frak I_-$, $\frak p_1\ge \frak p$,
$x \in A_{\frak p_1\frak q_1}$,
then
$$
\varphi_{\frak p_1\frak p}^{-1}(U_{(\frak p_1,x)}) \cap U'_{(\frak p,p)}
= \emptyset.
$$
\end{enumerate}
\end{subsublemma}
\begin{proof}
We take
\begin{equation}\label{form99}
U'_{(\frak p,p)} = U_{(\frak p,p)} \setminus \bigcup_{x \in A_{\frak p_1\frak q_1}}\overline{U_{(\frak p_1,x)}}.
\end{equation}
Here we regard $U_{(\frak p,p)}$ and $\overline{U_{(\frak p_1,x)}}$ as subsets of
$\vert {\widetriangle {\mathcal U^+}}\vert$.
(1) (3) are immediate.
We will prove (2).
By Subsublemma \ref{propety818} (3) and (\ref{eq09999}), we have
\begin{equation}\label{form1000}
Z_{\frak p\frak q} \cap \overline{U_{(\frak p_1,x)}}
= \emptyset,
\end{equation}
for each
$(\frak p,\frak q) \in \frak I_-$,
$x \in A_{\frak p_1\frak q_1}$.
Therefore $p \in A_{\frak p\frak q}$, $(\frak p,\frak q) \in \frak I_-$ imply
$p \notin \overline{U_{(\frak p_1,x)}}$.
Hence $p \in U'_{(\frak p,p)} \subset U_{(\frak p,p)}$.
This implies that Sublemma \ref{setsekigomasubme} (1)(2) hold
for $U'_{(\frak p,p)}$.
Sublemma \ref{setsekigomasubme} (4)' is a consequence of (\ref{form1000})
and (\ref{form99}).
\end{proof}
Hereafter we write $U_{(\frak p,p)}$ in place of $U'_{(\frak p,p)}$.
\par
We will prove that they have the properties claimed in
Sublemma \ref{setsekigomasubme}.
Sublemma \ref{setsekigomasubme}  (1),(2)
follow from Subsublemma \ref{propety818} (1),(2) and the
induction hypothesis.
Sublemma \ref{setsekigomasubme}  (4)' follows from
(\ref{eq8333}) and induction hypothesis
(which is claimed as Subsublemma \ref{cond819} (2)).
\begin{proof}[Proof of Sublemma \ref{setsekigomasubme} (3)' (a)]
Suppose
$(\frak p,\frak q), (\frak p',\frak q') \in \frak I$,
$p \in A_{\frak p\frak q}$, $p' \in A_{\frak p'\frak q'}$,
$\frak p \le \frak p'$ and
$
\varphi_{\frak p'\frak p}^{-1}(U_{(\frak p',p'))}) \cap U_{(\frak p,p)}
\ne \emptyset.
$
We will prove $\frak q \le \frak q'$.
\par
The case $(\frak p,\frak q), (\frak p',\frak q') \in \frak I_-$
follows from the induction hypothesis.
\par
Suppose $(\frak p',\frak q') = (\frak p_1,\frak q_1)$.
Then
$
\varphi_{\frak p_1\frak p}^{-1}(U_{(\frak p_1,p'))}) \cap U_{(\frak p,p)}
\ne \emptyset.
$
Subsublemma \ref{cond819}  (3)
implies that $(\frak p,\frak q) \notin \frak I_-$.
Therefore $(\frak p,\frak q) = (\frak p_1,\frak q_1)$.
Hence $\frak q \le \frak q'$ as required.
\par
We next assume $(\frak p_1,\frak q_1) = (\frak p,\frak q)$.
Then Subsublemma \ref{propety818} (4) implies $\frak q_1 \le \frak q'$
as required.
\end{proof}
\begin{proof}[Proof of Sublemma \ref{setsekigomasubme} (3)' (b)]
Suppose
$(\frak p,\frak q) \in \frak I$, $p \in A_{\frak p\frak q}$
and
$
\overline{U_{(\frak p,p)}} \cap Z_{\frak p'\frak q'} \ne \emptyset
$.
We will prove $(\frak p,\frak q) \ge (\frak p',\frak q')$.
\par
The case $(\frak p,\frak q) \in \frak I_-$
follows from the induction hypothesis.
Suppose $(\frak p,\frak q) = (\frak p_1,\frak q_1)$ .
Then $
\overline{U_{(\frak p_1,p)}} \cap Z_{\frak p'\frak q'} \ne \emptyset
$.
Note
$
\overline U_{(\frak p_1,x)} \subset \overline W^0_x \subset W_x.
$
Therefore Subsublemma \ref{propety818} (3) implies
$(\frak p',\frak q') \le (\frak p_1,\frak q_1)$, as required.
\end{proof}
\par
Therefore the proof of Sublemma \ref{setsekigomasubme} is now complete.
\end{proof}
Lemma \ref{lem816} is the case $\frak I = \frak P \times \frak P^+$
of Sublemma \ref{setsekigomasubme}.
\end{proof}
Now we put
$$
\frak P_0 = \bigcup_{(\frak p,\frak q) \in \frak P \times
\frak P^+} A_{\frak p\frak q} \times \{(\frak p,\frak q)\}.
$$
We choose any linear order on $A_{\frak p\frak q}$ and
define a partial order on $\frak P_0$ by the following:
$$
(x,(\frak p,\frak q))
\le (x',(\frak p',\frak q'))
\quad
\text{if and only if}
\quad
\begin{cases}
(\frak p,\frak q) < (\frak p',\frak q'))\\
\text{or $(\frak p,\frak q) = (\frak p',\frak q')$, $x\le x'$}.
\end{cases}
$$
We define
$$
 U^0_{(x,(\frak p,\frak q))} = U_{(\frak p,x)},
\quad
\mathcal U^0_{(x,(\frak p,\frak q))}
= \mathcal U_{\frak p}\vert_{ U^0_{(x,(\frak p,\frak q))} }.
$$
We define coordinate changes among them by restricting those
of ${\widetriangle{\mathcal U}}$.
We thus obtain a good coordinate system ${\widetriangle{\mathcal U^0}}$.
(Note we use Lemma \ref{lem816} (3) to check
Definition \ref{gcsystem} (5).)
\par
We will define a weakly open embedding
${\widetriangle{\mathcal U^0}} \to {\widetriangle{\mathcal U}}$.
We first define a map
$\frak P_0 \to \frak P$ by sending $(x,(\frak p,\frak q)) \mapsto \frak p$.
This is order preserving.
We also have an open embedding of Kuranishi charts
$\mathcal U^0_{(x,(\frak p,\frak q))} = \mathcal U_{\frak p}\vert_{U^0_{(x,(\frak p,\frak q))}}
\to \mathcal U_{\frak p}$.
They obviously commute with coordinate change.
\par
We next define the embedding
$
{\widetriangle{\mathcal U^0}} \to {\widetriangle{\mathcal U^+}}
$
that will be the composition of
${\widetriangle{\mathcal U^0}}
\to {\widetriangle{\mathcal U}}
\to {\widehat{\mathcal U}}
$
and
$ {\widehat{\mathcal U}} \to {\widetriangle{\mathcal U^+}}$.
We define a map
$\frak P_0 \to \frak P^+$ by sending
$(x,(\frak p,\frak q)) \mapsto \frak q$.
This is an order preserving map.
We next define a map
$\mathcal U^0_{(x,(\frak p,\frak q))}
\to \mathcal U^+_{\frak q}$
as the composition of
$$
\mathcal U^0_{(x,(\frak p,\frak q))}
\to
\mathcal U_{\frak p}\vert_{U_{\frak p}(x)}
\to \mathcal U_x
\to \mathcal U^+_{\frak q}.
$$
Here the first map
$\mathcal U^0_{(x,(\frak p,\frak q))}
\to
\mathcal U_{\frak p}\vert_{U_{\frak p}(x)}$
is an open embedding that exists by
Lemma \ref{lem816} (2).
The second map
$\mathcal U_{\frak p}\vert_{U_{\frak p}(x)}
\to \mathcal U_x$
is a part of the GK-embedding
${\widetriangle{\mathcal U}}
\to {\widehat{\mathcal U}}$.
The third map
$ \mathcal U_x
\to \mathcal U^+_{\frak q}$
is a part of the KG-embedding
${\widehat{\mathcal U}}
\to {\widetriangle{\mathcal U^+}}$.
The proof of Definition-Lemma \ref{henacomp}
is now complete.
\end{proof}

\subsection{Proof of Proposition \ref{integralinvembprop}}
\label{subsec:intwdonKura}

\begin{proof}[Proof of Proposition \ref{integralinvembprop}]
Let $(\mathcal K^+_1,\mathcal K^+_2)$ (resp. $(\mathcal K_1,\mathcal K_2)$)
be a support pair of ${\widetriangle{\mathcal U^+}}$ (resp. ${\widetriangle{\mathcal U}}$).
We may choose them so that
$
\varphi_{\frak p}(\mathcal K^i_{\frak p}) \subseteq \mathcal K^{i+}_{\frak i(\frak p)}
$.
Let $\mathcal K_2^+ < \mathcal K_3^+$ and
$\mathcal K_2 < \mathcal K_3$.
\par
We will choose $\delta_+$, $\delta$ and $\frak U(Z)$, later.
Let $\{\chi^+_{\frak p^+}\}$ (resp. $\{\chi_{\frak p}\}$) be a strongly smooth partition of unity
of $(X,Z,{\widetriangle{\mathcal U^+}},{\widetriangle{\frak S^+}},\delta_+)$
(resp. $(X,Z,{\widetriangle{\mathcal U}},{\widetriangle{\frak S}},\delta)$).
By inspecting the proof of Proposition \ref{pounitexi},
we can take $\chi_{\frak p}$ so that it is not only a strongly smooth function on
$\vert \mathcal K_2\vert$ but also one on $\vert \mathcal K^+_2\vert$.
\par
We take $\frak p^+_0 \in \frak P^+$ and set $h_0 = \chi^+_{\frak p^+_0}h_{\frak p^+_0}$.
To prove Proposition \ref{integralinvembprop} it suffices to show
\begin{equation}\label{formula84}
f_{\frak p^+_0}!\left( h_0; \frak S^{+\epsilon}_{\frak p^+_0} \vert_{\frak U(Z)\cap
\mathcal K^{+1}_{\frak p_0^+}(2\delta_+)}\right)
=
\sum_{\frak p \in \frak P}
f_{\frak p}!( (\chi_{\frak p}h_0)_{\frak p}; \frak S^{\epsilon}_{\frak p} \vert_{\frak U(Z)\cap
\mathcal K^{1}_{\frak p}(2\delta)}).
\end{equation}
By taking $\epsilon >0$ sufficiently small, we may assume $\sum \chi_{\frak p}= 1$ on
$\frak U(Z) \cap \Pi((\frak S^{+ \epsilon}_{\frak p^+_0})^{-1}(0))$.
(This is a consequence of Lemma \ref{lem739} and Definition \ref{pounity} (3).
Note the differential form $(\chi_{\frak p}h_0)_{\frak p_0}$
is defined since the function $\chi_{\frak p}$ is strongly smooth
on $\vert\mathcal K_2^+\vert$.)
Therefore to prove (\ref{formula84}) it suffices to show
\begin{equation}\label{formula85}
f_{\frak p^+_0}!\left( (\chi_{\frak p}h_0)_{\frak p_0^+}; \frak S^{+\epsilon}_{\frak p^+_0}
\vert_{\frak U(Z)\cap
\mathcal K^{+1}_{\frak p_0^+}(2\delta_+)}\right)
=
f_{\frak p}!( (\chi_{\frak p}h_0)_{\frak p}; \frak S^{\epsilon}_{\frak p} \vert_{\frak U(Z)\cap
\mathcal K^{1}_{\frak p}(2\delta)})
\end{equation}
for each $\frak p$. We will prove it below. There are three cases.
\par\medskip
\noindent(Case 1)  Neither $\frak i(\frak p) \le \frak p^+_0$ nor $\frak i(\frak p) \ge \frak p^+_0$:
In this case we have
$$
\mathcal K^{1+}_{\frak p^+_0} \cap \mathcal K^{1}_{\frak p}
\subset
\mathcal K^{1+}_{\frak p^+_0} \cap \mathcal K^{1 +}_{\frak i(\frak p)}
= \emptyset.
$$
Therefore
in the same way as the proof of (\ref{form743}), we can choose $\delta$,$\delta_+$  small so that
$$
\Omega^{1+}_{\frak p^+_0}(\mathcal K^{1+},\delta_+) \cap
\Omega^{1}_{\frak p}(\mathcal K^{1},\delta) = \emptyset.
$$
Then both sides of (\ref{formula85}) are zero.
\par\medskip
\noindent(Case 2)  $\frak i(\frak p) \le \frak p^+_0$:
We consider the embedding
$$
\mathcal U_{\frak p}   \,\,\overset{\Phi_{\frak p}}\longrightarrow \,\, \mathcal U^+_{\frak i(\frak p)}
 \,\,\overset{\Phi_{\frak p^+_0 \frak i(\frak p)}}\longrightarrow \,\, \mathcal U^+_{\frak p_0}.
$$
In the same way as the proof of (\ref{form741741}) we can choose $\delta$, $\delta_+$,
$\frak U(Z)$ small so that
\begin{equation}\label{form914914}
{\rm Supp}(\chi_{\frak p}\widetriangle h_0) \cap \Pi((\widetriangle{\frak S^{+\epsilon}})^{-1}(0))
\cap \frak U(Z)
\subset \mathcal K^{1 +}_{\frak p^+_0}(2\delta^+)
\cap
\mathcal K^{1}_{\frak p}(2\delta)
\cap \frak U(Z).
\end{equation}
Then (\ref{formula85}) follows.
\par\medskip
\noindent(Case 3)  $\frak i(\frak p) \ge \frak p^+_0$:
In the same way as the proof of Proposition \ref{indepofukuracont} Case 2,
we can choose $\delta$, $\delta_+$,
$\frak U(Z)$ small so that
(\ref{form914914}) holds. (\ref{formula85}) is its consequence.
\par
Thus the proof of Proposition \ref{integralinvembprop} is complete.
\end{proof}
\subsection{CF-perturbations of correspondences}
\label{subsec:confamicor}
\begin{defn}\label{def92111}
We consider Situation \ref{smoothcorr}.
Let $\widehat{\frak S}$ be a CF-perturbation of $\widehat{\mathcal U}$
such that $\widehat {f_t}$ is strongly submersive
with respect to $\widehat{\frak S}$.
\par
We call such $\widehat{\frak S}$ a {\it CF-perturbation
of Kuranishi correspondence}
\index{CF-perturbation ! of Kuranishi correspondence} $\frak X$.
We then define
\begin{equation}\label{defn924}
{\rm Corr}_{(\frak X,\widehat{\frak S^{\epsilon}})}
: \Omega^k(M_s) \to   \Omega^{k+\ell}(M_t)
\end{equation}
by
\begin{equation}
{\rm Corr}_{(\frak X,\widehat{\frak S^{\epsilon}})}(h)
=
\widehat{f_t} !((\widehat{f_s})^*h;{\widehat{\frak S}}^{\epsilon}).
\end{equation}
This is well-defined by Definition-Lemma \ref{deflemgg}.
We call the linear map ${\rm Corr}_{(\frak X,\widehat{\frak S^{\epsilon}})}$
a {\it smooth correspondence map}
\index{smooth correspondence ! smooth correspondence} of Kuranishi structure and
$\ell$ the {\it degree}
\index{smooth correspondence ! degree} of smooth correspondence $\frak X$
and write it $\deg \frak X$.
\end{defn}
The next lemma says that in Situation \ref{smoothcorr}
we can always thicken our Kuranishi structure so that
the assumptions of Definition \ref{def92111} is satisfied.
\begin{lem}
For each  smooth correspondence $(X,\widehat{\mathcal U},\widehat{f_s},\widehat{f_t})$
as in Situation \ref{smoothcorr} there exist
$\widehat{\mathcal U^+},\widehat{\frak S^+},\widehat{f^+_s},\widehat{f^+_t}$
with the following properties.
\begin{enumerate}
\item
$(X,\widehat{\mathcal U^+},\widehat{f^+_s},\widehat{f^+_t})$
is a Kuranishi correspondence and
$\widehat{\frak S^+}$ a CF-perturbation of Kuranishi correspondence.
\item
$\widehat{\mathcal U^+}$ is a thickening of $\widehat{\mathcal U}$.
\item
Let $\widehat\Phi : \widehat{\mathcal U} \to \widehat{\mathcal U^+}$
be the KK-embedding. Then $\widehat{f^+_s}$ and $\widehat{f^+_t}$
induce $\widehat{f_s}$ and $\widehat{f_t}$
by $
\widehat\Phi
$.
\end{enumerate}
\end{lem}
\begin{proof}
This is an immediate consequence of
Lemmata \ref{lemappgcstoKucont} and \ref{le714}.
\end{proof}
\subsection{Stokes' formula for Kuranishi structure}
\label{subsec:StokesKura}
We have Stokes' formula in Theorem \ref{Stokes}, which
is the formula for good coordinate system.
In this subsection we translate it to one for Kuranishi structure.

\begin{shitu}\label{stokeskurashitu}
Let $\widehat{\mathcal U}$ be a Kuranishi structure of $Z \subseteq X$,
$\widehat{\frak S}$  its CF-perturbation,
$\widehat f :  (X,Z;\widehat{\mathcal U}) \to M$ a strongly
submersive map with respect to $\widehat{\frak S}$,
and $\widehat h$ a differential form on $ (X,Z;\widehat{\mathcal U})$.
\par
Let $\partial(X,Z,\widehat{\mathcal U},\widehat{\frak S})
= (\partial X,\partial Z,\widehat{\mathcal U_{\partial}},\widehat{\frak S_{\partial}})$,
where $(\partial X,\partial Z;\widehat{\mathcal U_{\partial}})$
is the normalized boundary of $(X,Z;\widehat{\mathcal U})$ on which
$\widehat{\frak S}$ induces a CF-perturbation
$\widehat{\frak S_{\partial}}$ by Lemma \ref{lemma755} (2).
Since $\widehat f$  induces a map $\widehat{f_{\partial}} : (\partial X,
\partial Z;\widehat{\mathcal U_{\partial}})
\to M$, which is strongly submersive with respect to $\widehat{\frak S_{\partial}}$
if $\widehat f$ is strongly submersive with respect to $\widehat{\frak S}$
(Lemma \ref{lemma755} (4)).
\par
Let $\widehat{h_{\partial}}$ be the restriction of $\widehat h$ to
$(\partial X,\partial Z;\widehat{\mathcal U_{\partial}})$.$\blacksquare$
\end{shitu}
\begin{prop}{\rm (Stokes' formula for Kuranishi structure.)}\label{Stokeskura}
In Situation \ref{stokeskurashitu} we have the next formula for each sufficiently small $\epsilon>0$:
\begin{equation}
d\left(\widehat f !(\widehat h;\widehat{\frak S^{\epsilon}})\right)
=
\widehat f !(d\widehat h;\widehat{\frak S^{\epsilon}})
+
\widehat f_{\partial}!(\widehat {h_{\partial}};\widehat{\frak S_{\partial}^{\epsilon}}).
\end{equation}
\end{prop}
\begin{proof}
By Lemmata \ref{lemappgcstoKucont} and \ref{le714},
there exist a good coordinate system ${\widetriangle{\mathcal U}}$
and an KG-embedding
${\widehat{\mathcal U}} \to {\widetriangle{\mathcal U}}$.
Moreover there exists a CF-perturbation
${\widetriangle{\frak S}}$ of ${\widetriangle{\mathcal U}}$
such that ${\widehat{\frak S}}$, ${\widetriangle{\frak S}}$
are compatible with the KG-embedding
${\widehat{\mathcal U}} \to {\widetriangle{\mathcal U}}$.
Furthermore there exist a strongly smooth map
${\widetriangle f} : (X,Z;{\widetriangle{\mathcal U}})
\to M$ and a differential form $\widetriangle h$, which are pulled back to $\widehat f$
and $\widehat h$, by the KG-embedding
${\widehat{\mathcal U}} \to {\widetriangle{\mathcal U}}$.
Then ${\widetriangle{\mathcal U}}$, ${\widetriangle{\frak S}}$, ${\widetriangle f}$
and $\widetriangle h$, induce
${\widetriangle{\mathcal U_{\partial}}}$, ${\widetriangle{\frak S_{\partial}}}$,
$\widetriangle{f_{\partial}}$ and $\widetriangle{h_{\partial}}$
on the boundary, respectively, which are compatible with corresponding objects on
$(\partial X,\partial Z;\widehat{\mathcal U_{\partial}})$.
Thus Proposition \ref{Stokeskura} follows
by applying Theorem \ref{Stokes}  to
${\widetriangle{\mathcal U}}$, ${\widetriangle{\frak S}}$, ${\widetriangle f}$,
$\widetriangle h$,
and
${\widetriangle{\mathcal U_{\partial}}}$, ${\widetriangle{\frak S_{\partial}}}$, $\widetriangle{f_{\partial}}$,  $\widetriangle{h_{\partial}}$.
\end{proof}
The next corollary is an immediate consequence of
Proposition \ref{Stokeskura}.
\begin{cor}
In the situation of Definition \ref{def92111} we have
 the next formula for each sufficiently small $\epsilon >0$:
$$
d \circ {\rm Corr}_{(\frak X,\widehat{\frak S^{\epsilon}})}
= {\rm Corr}_{(\frak X,\widehat{\frak S^{\epsilon}})} \circ d
+  {\rm Corr}_{\partial(\frak X,\widehat{\frak S^{\epsilon}})}.
$$
\end{cor}

\subsection{Uniformity of CF-perturbations on
Kuranishi structure}
\label{subsec:unfcfpKura}
In this subsection,
we collect various facts which we use to show the existence of
uniform bound of the constants $\epsilon$ that appear
in Theorem \ref{theorem915} etc..

\begin{defn}
Let $\widehat{\mathcal U}$ be a Kuranishi structure
on $Z \subseteq X$ and
$\widehat{\frak S_{\sigma}}$ be a $\sigma \in \mathscr A$
parameterized family of CF-perturbations.
We say that $\widehat{\frak S_{\sigma}}$ is a {\it uniform family}
\index{uniform family ! of CF-perturbations}
\index{CF-perturbation ! uniform family} if the convergence in Definition \ref{defn73ss}
is uniform.
More precisely, we require the following.
\par
For each $\frak o$ there exists $\epsilon_0 >0$
such that if $0<\epsilon < \epsilon_0$,
$p \in Z$ then
\begin{equation}
\vert s(y) - s_p(y) \vert < \frak o,
\qquad
\vert (Ds)(y) - (Ds_p)(y) \vert < \frak o,
\end{equation}
hold for any $s$ which is any member of ${\frak S}^{\epsilon,p}_{\sigma}$
at any point $y \in U_{p}$ and $\sigma \in \mathscr A$.
\end{defn}

\begin{lem}\label{lem92929}
\begin{enumerate}
\item
In the situation of Lemma \ref{lemappgcstoKucont}
if $\widetriangle{\frak S}$ varies in a uniform family
then $\widehat{\frak S}$ varies in a uniform family.
\item
In the situation of Lemma \ref{le714} (2), if
${\widehat{\frak S}}$ varies in a uniform family
then $\widetriangle{\frak S}$ varies in a uniform family.
\item
In the situation of Lemma \ref{le7155},
if ${\widehat{\frak S^+}}$, ${\widehat{\frak S}}$
vary in a uniform family
(resp. ${\widehat{\frak S^+_a}}$ ($a=1,2$)
, ${\widehat{\frak S}}$ vary in a uniform family)
then $\widetriangle{\frak S}$
varies in a uniform family.
\end{enumerate}
\end{lem}
The proof will be given at the end of  Subsection \ref{subsec:movingmulsectionetc}.

\begin{prop}
In the situation of Theorem \ref{theorem915}
suppose ${\widetriangle{\frak S_{\sigma}}}$
varies in a uniform family.
(We require that ${\widetriangle{\mathcal U}}$,
$\widetriangle{\Phi}$
are independent of the parameter $\sigma$.)
Then the pushout $\widehat f !(\widehat{ h};{\widehat{\frak S_{\sigma}}})$
is uniformly  independent of the choices in the sense of $\clubsuit$
in Definition \ref{intheclubsuit}.
We may choose the constant $\epsilon$ in Proposition \ref{Stokeskura}
independent of $\sigma$ also.
\end{prop}
\begin{proof}
Using Lemma \ref{lem92929} the proof goes in the same way as
the proof of Proposition \ref{lem761}.
\end{proof}
\begin{rem}
We can choose $\epsilon_0$ independent of $\widehat f$ and $\widehat h$.
\end{rem}

\section{Composition formula of smooth correspondences}
\label{sec:composition}

The purpose of this section is to provide thorough technical
detail of the proof of \cite[Lemma 12.15]{fooo09}
= Theorem \ref{compformulaprof},
where fiber product of Kuranishi structures
is used as a way to define composition of
smooth correspondences.
For this purpose we work out the plan described in
Subsection \ref{bigremarkinsec6}
in the de Rham model.

\subsection{Direct product and CF-perturbation}
\label{subsec:dprocontper}

Firstly,
we begin with defining direct product
of CF-perturbations.

\begin{shitu}\label{sit824}
For each $i=1,2$,
$\mathcal U_i = (U_i,\mathcal E_i,s_i,\psi_i)$ is a Kuranishi chart of $X$,
$x_i \in U_i$, $\frak V_{x_i}^i = (V^i_{x_i},\Gamma^i_{x_i},E^i_{x_i},\psi^i_{x_i},\widehat\psi^i_{x_i})$
is an orbifold chart of $(U_i,\mathcal E_i)$ as in Definition \ref{defn73ss}.
Let $\mathcal S_{x_i}^i = (W_{x_i}^i,\omega_{x_i}^i,\frak s_{x_i}^{i\epsilon})$
be a CF-perturbation $\mathcal U_i$ on $\frak V_{x_i}^i$.$\blacksquare$
\end{shitu}

\begin{defn}\label{defn1022}
In Situation \ref{sit824},
we define the {\it direct product of $\mathcal S_{x_1}^1$ and $\mathcal S_{x_2}^2$}  \index{CF-perturbation ! direct product of two CF-perturbations
on one chart}
\index{fiber product ! direct product of two CF-perturbations
on one chart}
by
$$
\mathcal S_{x_1}^1 \times \mathcal S_{x_2}^2
=
(W_{x_1}^1 \times W_{x_2}^2, \omega_{x_1}^1 \times \omega_{x_2}^2, \frak s_{x_1}^{1\epsilon} \times
\frak s_{x_2}^{2\epsilon}),
$$
where
$$
( \frak s_{x_1}^{1\epsilon} \times
\frak s_{x_2}^{2\epsilon})(y_1,y_2,\xi_1,\xi_2)
=
( \frak s_{x_1}^{1\epsilon}(y_1,\xi_1),\frak s_{x_2}^{2\epsilon}(y_2,\xi_2))
$$
for $y_i \in V^i_{x_i}$ $\xi_i \in W_{x_i}^i$.
\end{defn}
\begin{lem}\label{lem826}
\begin{enumerate}
\item $\mathcal S_{x_1}^1 \times \mathcal S_{x_2}^2$
is a CF-perturbation of $\mathcal U_1 \times \mathcal U_2$.
\item
If $\mathcal S_{x_i}^i$ are equivalent to $\mathcal S_{x_i}^{i\prime}$ for $i=1,2$,
then
$\mathcal S_{x_1}^1 \times \mathcal S_{x_2}^2$  is equivalent to
$\mathcal S_{x_1}^{1\prime} \times \mathcal S_{x_2}^{2\prime}$.
\item
Let $\Phi^i : \frak V_{x'_i}^{i\prime} \to \frak V_{x_i}^i$ be an embedding
of orbifold chart. Then
$$
(\Phi^1)^*\mathcal S_{x_1}^1 \times (\Phi^2)^*\mathcal S_{x_2}^2
$$
is equivalent to
$$
(\Phi^1 \times \Phi^2)^*(\mathcal S_{x_1}^1 \times \mathcal S_{x_2}^2).
$$
\end{enumerate}
\end{lem}
This is a direct consequence of definitions.
\begin{lemdef}
Suppose we are in Situation \ref{sit824}.
\begin{enumerate}
\item
Let
$
\frak S^i = \{(\frak V^i_{\frak r_i},\mathcal S^i_{\frak r_i}) \mid \frak r_i \in \frak R_i\}
$
be representatives of CF-perturbations of $\mathcal U^i$ for $i=1,2$.
Then
$$
\{ (\frak V^1_{\frak r_1} \times \frak V^2_{\frak r_2},\mathcal S^1_{\frak r_1} \times \mathcal S^2_{\frak r_2})
\mid (\frak r_1,\frak r_2) \in \frak R_1 \times \frak R_2\}
$$
is a representative of a CF-perturbation of $\mathcal U^1 \times \mathcal U^2$.
We call it the {\rm direct product}
and write $\frak S^1 \times \frak S^2$.
\item
If $\frak S^i$ is equivalent to $\frak S^{i\prime}$, then
$\frak S^1 \times \frak S^2$ is equivalent to
$\frak S^{1\prime} \times \frak S^{2\prime}$.
\item
Therefore we can define direct product
of CF-perturbations.
\item
Direct product defines a sheaf morphism
\begin{equation}\label{form101}
\pi_1^{\star} \mathscr S^{\mathcal U_1} \times
\pi_2^{\star} \mathscr S^{\mathcal U_2}
\to
\mathscr S^{\mathcal U_1 \times \mathcal U_2},
\end{equation}
where $\pi_i : U_1 \times U_2 \to U_i$ are projections.
\end{enumerate}
\end{lemdef}
\begin{proof}
This is an immediate consequence of Lemma \ref{lem826}.
\end{proof}
\begin{lem}\label{lem828}
Let $\Phi^i : \mathcal U^{i} \to \mathcal U^{i +}$ be embeddings of Kuranishi charts and
$\frak S^i$, $\frak S^{i + }$
CF-perturbations of $\mathcal U^{i}$, $ \mathcal U^{i +}$, for $i=1,2$, respectively.
\begin{enumerate}
\item
If $\frak S^{i +}$ can be pulled back by $\Phi^i$ for $i=1,2$, then
$\frak S^{1 +} \times \frak S^{2 +}$ can be pulled back by $\Phi^1 \times \Phi^2$.
\item
If $\frak S^{i +}$, $\frak S^{i}$ are  compatible with $\Phi^i$  for $i=1,2$,
then $\frak S^{1 +} \times \frak S^{2 +}$
and $\frak S^{1} \times \frak S^{2}$ are  compatible with $\Phi^1 \times \Phi^2$.
\item The next diagram commutes:
\begin{equation}\label{diagrampart1010}
\begin{CD}
(\Phi^1)^{\star}\pi_1^{\star} \mathscr S^{\mathcal U^1\triangleright \mathcal U^{1+}} \times
(\Phi^2)^{\star}\pi_2^{\star} \mathscr S^{\mathcal U^2\triangleright \mathcal U^{2+}}
@ > {(\ref{form101})} >>
(\Phi^1\times \Phi^2)^{\star} \mathscr S^{(\mathcal U_1
\times \mathcal U^2) \triangleright (\mathcal U^{1+} \times
\mathcal U^{2+})}
\\
 @ VV{(\Phi^1)^* \times (\Phi^2)^*}V @VV{(\Phi^1\times \Phi^2)^*}V\\
\pi_1^{\star} \mathscr S^{\mathcal U^1} \times
\pi_2^{\star} \mathscr S^{\mathcal U^2}
@>> {(\ref{form101})} > \mathscr S^{\mathcal U^1 \times \mathcal U^2}
\end{CD}
\nonumber
\end{equation}
\end{enumerate}
\end{lem}
This is a direct consequence of the definitions.
\begin{lemdef}
Let $\widehat{\mathcal U^i}= (\{\mathcal U^i_{p_i}\},\{\Phi^i_{p_iq_i}\})$ be Kuranishi structures of $Z_i \subseteq X_i$
for $i=1,2$,
and $\widehat{\mathcal U^1} \times \widehat{\mathcal U^2}$
the direct product Kuranishi structure on $Z_1 \times Z_2 \subseteq X_1 \times X_2$.
Let $\widehat{\frak S^i} = \{\frak S^i_{p_i}\}$ be CF-perturbations of $\mathcal U^i_{p_i}$.
Then $ \{\frak S^1_{p_1} \times \frak S^2_{p_2}\}$ defines a CF-perturbation of
$\widehat{\mathcal U^1} \times \widehat{\mathcal U^2}$. We call it the {\rm direct product} of CF
perturbations and denote it by
$\widehat{\frak S^1} \times \widehat{\frak S^2}$.
\index{CF-perturbation ! direct product of CF-perturbations
of Kuranishi structures}
\index{fiber product ! direct product of CF-perturbations
of Kuranishi structures}
\index{Kuranishi chart ! direct product of CF-perturbations
of Kuranishi structures}
\end{lemdef}
\begin{proof}
This is an immediate consequence of Lemma \ref{lem828}.
\end{proof}
We have thus defined the direct product of CF-perturbations.
\begin{rem}
We have defined the notion of direct product of CF-perturbations
of Kuranishi structures, but {\it not} one of good coordinate systems. The reason is explained
at the end of Section \ref{sec:fiber}.
\end{rem}
\subsection{Fiber product and CF-perturbation}
\label{subsec:fprocontper}
We next discuss the case of fiber product.
\begin{defn}
Let $\mathcal U = (U,\mathcal E,s,\psi)$ be a Kuranishi chart of $X$.
For
$x \in U$ let $\frak V_x$ be an orbifold chart of
$(U,\mathcal E)$ and $\mathcal S_x = (W_x,\omega_x,\{\frak s^{\epsilon}_x\})$  a CF-perturbation of $\mathcal U$ on $\frak V_x$.
Let $f : U \to M$ be a smooth map to a manifold $M$ and $g : N \to M$  a
smooth map from a manifold $N$.
Suppose  $f$ is strongly transversal to $g$ with respect to $\mathcal S_x$
in the sense of Definition \ref{submersivepertconlocloc} (3).
Then we take the fiber product $X {}_f\times_g N$,
fiber product Kuranishi chart $(\frak V_x) {}_f\times_g N$ and
$(\mathcal  S_x){}_f\times_g N=  (W_x,\omega_x,\{(\frak s^{\epsilon}_x)\,\,{}_f\times_g N\})$.
Here
$$
((\frak s^{\epsilon}_x) {}_f\times_g N) : ((V_x) {}_{f}\times_g N) \times W_x
\to E_x
$$
is defined by
$$
((\frak s^{\epsilon}_x) {}_f\times_g N)((y,z),\xi) = \frak s^{\epsilon}_x(y,\xi).
$$
We call $(\mathcal  S_x){}_f\times_g N$ the  {\it fiber product CF-perturbation}. \index{CF-perturbation ! fiber product CF-perturbation on
one chart}
\index{fiber product ! fiber product CF-perturbation on
one chart}
\index{Kuranishi chart ! fiber product CF-perturbation on
one chart}
It is a CF-perturbation of $\mathcal U {}_f\times_g N$.
\end{defn}
\begin{lem}\label{lem826222}
\begin{enumerate}
\item
If $\mathcal S_{x}$ is equivalent to $\mathcal S'_{x}$ and
 $f$ is strongly transversal to $g$ with respect to $\mathcal S_x$,
 $f$ is strongly transversal to $g$ with respect to $\mathcal S'_x$.
Moreover
$(\mathcal  S_x){}_f\times_g N$  is equivalent to
$(\mathcal  S'_x){}_f\times_g N$.
\item
Let $\Phi : \frak V_{x'}^{\prime} \to \frak V_{x}$ be an embedding
of orbifold charts.
If  $f \circ \Phi$ is strongly transversal to $g$ with respect to
$\Phi^*\frak S_x$, $f$ is strongly transversal to $g$ with respect to
$\frak S_x$. Moreover
$$
\Phi^*((\mathcal S_{x}) {}_f \times_g N)
$$
is equivalent to
$$
(\Phi^*(\mathcal S_{x})) {}_f \times_g N.
$$
\end{enumerate}
\end{lem}
This is a direct consequence of the definitions.
\begin{lemdef}\label{lemdef834}
Let
$\mathcal U = (U,\mathcal E,s,\psi)$ be a Kuranishi chart of $X$,
and $\frak S = \{(\frak V_{\frak r},\mathcal S_{\frak r})
\mid \frak r \in \frak R\}$
a representative of a CF-perturbation of  $\mathcal U$.
Let $f : U \to M$ be a smooth map to a manifold $M$ and $g : N \to M$  a
smooth map from a manifold $N$.
\begin{enumerate}
\item
If $f$ is strongly transversal to $g$ with respect to $\frak S$
in the sense of Definition-Lemma \ref{strosubsemiloc} (3),
then
$$
\frak S {}_f\times_g N= \{((\frak V_{\frak r}) {}_f\times_g N,(\mathcal S_{\frak r}) {}_f\times_g N) \mid \frak r \in \frak R\}
$$
is a CF-perturbation of $\mathcal U{}_f\times_g N$.
\item
If $\frak S$ is equivalent to $\frak S'$, then
$\frak S {}_f\times_g N$ is equivalent to $\frak S' {}_f\times_g N$.
\item
Therefore we can define a {\rm fiber product of
CF-perturbations}
\index{fiber product ! of CF-perturbations with a map on a single
Kuranishi chart}
\index{CF-perturbation ! fiber product of
CF-perturbations with a map on a single
Kuranishi chart} with a map $g : N \to M$
when $\widehat f$ is strongly transversal to $g$.
\end{enumerate}
\end{lemdef}
\begin{proof}
This is an immediate consequence of Lemma \ref{lem826222}.
\end{proof}
\begin{lem}\label{lem828222}
Let $\Phi : \mathcal U \to \mathcal U^+$ be an embedding of Kuranishi charts
and
$\frak S$, $\frak S^{+ }$
CF-perturbations of $\mathcal U$, $ \mathcal U^+$, respectively.
Suppose $\frak S$, $\frak S^{+ }$ are compatible with $\Phi$.
Let $f^+ : \mathcal U^+ \to M$ be a strongly smooth map to a manifold $M$,
$f  = f^+ \circ \varphi : \mathcal U \to M$, and $g : N \to M$  a smooth map from a manifold $N$.
We assume $f, f^+$ are strongly transversal to $g$ with respect to $\frak S$, $\frak S^{+ }$, respectively.
\par
Then
$\frak S^+ \,{}_{f^+}\times_g N$,  $\frak S \,{}_f\times_g N$
are compatible to $\Phi \times {\rm id} :
\mathcal U \,{}_f\times_g N \to \mathcal U^+ \,{}_{f^+}\times_g N$.
\end{lem}
This is a direct consequence of the definitions.
\begin{lemdef}\label{lemdef835}
Let $\widehat{\mathcal U}= (\{\mathcal U_{p}\},\{\Phi_{pq}\})$ be a Kuranishi structure of $Z \subseteq X$
and $\widehat{\frak S} = \{\frak S_{p}\}$ a CF-perturbation of $\widehat{\mathcal U}$.
Suppose that a strongly smooth map $\widehat f : (X,\widehat{\mathcal U}) \to M$ is strongly transversal
to a smooth map $g : N \to M$
with respect to $\widehat{\frak S}$ in the sense of
Definition \ref{smoothfunctiononvertK} (3).
\par
Then
$\{(\frak S_{p}) {}_f\times_g N\}$ is a CF-perturbation.
We call it {\rm a fiber product CF-perturbations}
\index{CF-perturbation ! fiber product CF-perturbation of Kuranishi structure}
\index{fiber product ! fiber product CF-perturbation of Kuranishi structure}
 and write
$(\frak S_{p}) {}_f\times_g N$.
\end{lemdef}
Lemma-Definition \ref{lemdef835} is a consequence of Lemma \ref{lem828222}.
\begin{defn}\label{defn837}
Let  $\widehat{\mathcal U^i}= (\{\mathcal U^i_{p_i}\},\{\Phi^i_{p_iq_i}\})$ be Kuranishi structures of $Z_i \subseteq X_i$
and  $\widehat{\mathcal U^1} \times \widehat{\mathcal U^2}$  the direct product Kuranishi structure on $Z_1 \times Z_2 \subseteq  X_1 \times X_2$.
Let $\widehat{\frak S^i} = \{\frak S^i_{p_i}\}$ be  CF-perturbations of $\mathcal U^i_{p_i}$
and $\widehat{\frak S^1} \times \widehat{\frak S^2}$ their direct product.
Let $\widehat{f^i} : (X_i,\widehat{\mathcal U^i}) \to M$  be strongly smooth maps to a manifold $M$.
\begin{enumerate}
\item
We say that $\widehat{f^1}$ is {\it strongly transversal to $\widehat{f^2}$ with respect to
$\widehat{\frak S^1}$, $\widehat{\frak S^2}$},
\index{strongly transversal (w.r.t. CF-perturbation) ! two maps from Kuranishi structures}
if and only if
$$
(\widehat{f^1},\widehat{f^2}) : (X_1\times X_2,\widehat{\mathcal U^1} \times \widehat{\mathcal U^2})
\to M \times M
$$
is strongly transversal to the diagonal $\Delta_M =\{(x,x) \mid x \in M\}$,
with respect to the direct product $\widehat{\frak S^1}\times \widehat{\frak S^2}$
in the sense of Definition \ref{smoothfunctiononvertK} (3).
\item
In the situation of (1), we define the {\it fiber product of CF-perturbations} \index{fiber product ! of two
CF-perturbations of Kuranishi structures}
\index{CF-perturbation ! fiber product of two
CF-perturbations of Kuranishi structures}
by
$$
(\widehat{\frak S^1}) {}_{\widehat{f^1}} \times_{\widehat{f^2}}  (\widehat{\frak S^2})
=
(\widehat{\frak S^1}\times\widehat{\frak S^2}) \,{}_{(\widehat{f^1},\widehat{f^2})}\times_{M\times M} \Delta_M.
$$
Here the right hand side is defined by Lemma-Definition \ref{lemdef835}.
\end{enumerate}
\end{defn}
\begin{lem}\label{lem838}
\begin{enumerate}
\item
Suppose we are in the situation of Definition \ref{defn837}.
If $\widehat{f^1}$ is strongly submersive
with respect to $\widehat{\frak S^1}$,
then $\widehat{f^1}$ is strongly transversal to any $\widehat{f^2}$
with respect to $\widehat{\frak S^1}$ and $\widehat{\frak S^2}$, provided $\widehat{\frak S^2}$ is transversal to $0$.
\item
In the situation of Definition \ref{defn837},
we assume $\widehat{f^1}$ is strongly submersive
with respect to $\widehat{\frak S^1}$.
Let $\widehat{f^3} : (X_2,\widehat{\mathcal U^2}) \to N$
be another strongly smooth map such that
$\widehat{f^3}$ is strongly submersive with respect to
$\widehat{\frak S^2}$.
\par
Then
$$
\widetilde{\widehat{f^3}}
: (X_1\times X_2,\widehat{\mathcal U^1} \times \widehat{\mathcal U^2})
\to N
$$
is strongly submersive with respect to
$(\widehat{\frak S}^1) {}_{\widehat{f^1}} \times_{\widehat{f^2}}  (\widehat{\frak S}^2)$.
Here
$
\widetilde{\widehat {f^3}}$
is the map induced from $\widehat{f^3}$.
\item
The fiber product of uniform family of CF-perturbations
is uniform.
\end{enumerate}
\end{lem}
\begin{proof}
It suffices to prove the corresponding statement on a single orbifold chart.
Namely for each $i=1,2$ we consider $\frak V_{{\frak r}_i}^i$ an orbifold chart of a Kuranishi
neighborhood of $\widehat{\mathcal U}^i$,
a CF-perturbation $\mathcal S^i_{{\frak r}_i}$ of it,
and maps
$f^i_{{\frak r}_i} : U^i_{{\frak r}_i} \to M$,
$f^3_{{\frak r}_2} : U^2_{{\frak r}_2} \to N$.
We will prove this case below.
\par\medskip
\noindent
{\it Proof of (1)} :
By assumption
$$
f^1_{{\frak r}_1}\vert_{(\mathcal S^{1\epsilon}_{{\frak r}_1})^{-1}(0)}
: (\mathcal S^{1\epsilon}_{{\frak r}_1})^{-1}(0) \to M
$$
is a submersion.
Therefore it is transversal to
$$
f^2_{{\frak r}_2}\vert_{(\mathcal S^{2\epsilon}_{{\frak r}_2})^{-1}(0)}
: (\mathcal S^{2\epsilon}_{{\frak r}_2})^{-1}(0) \to M
$$
as required.
\par\smallskip
\noindent
{\it Proof of (2)} :
Let
$
(y_i,\xi_i) \in (\mathcal S^{i\epsilon}_{{\frak r}_i})^{-1}(0).
$
Here $y_i \in V^i_{{\frak r}_i}$, $\xi_i \in W^i_{{\frak r}_i}$
where
$\frak V_{x_i}^i = (V^i_{{\frak r}_i},\Gamma^i_{x_i},E^i_{{\frak r}_i},\psi^i_{{\frak r}_i},\widehat{\psi}^i_{{\frak r}_i})$,
$\mathcal S^{i\epsilon}_{{\frak r}_i} = (W^i_{{\frak r}_i},\omega^i_{{\frak r}_i},\frak s^{i\epsilon}_{{\frak r}_i})$.
\par
Suppose
$f^1_{{\frak r}_1}(y_1) = f^2_{{\frak r}_2}(y_2) = z$
and $f^3_{{\frak r}_2}(y_2) = w$.
We consider
$$
\aligned
(d_{y_1}f^1_{{\frak r}_1}
\oplus
d_{y_2}f^2_{{\frak r}_2})
\oplus
d_{y_2}f^3_{{\frak r}_2}
:
T_{(y_1,\xi_1)}(\mathcal S^{1\epsilon}_{{\frak r}_1})^{-1}(0)
&\oplus
T_{(y_2,\xi_2)}(\mathcal S^{2\epsilon}_{{\frak r}_2})^{-1}(0)
\\
&\to
T_zM \oplus T_z M \oplus T_wN.
\endaligned
$$
Let $\frak v \in T_wN$.
Then there exists $\tilde{\frak v_2}\in
T_{(y_2,\xi_2)}(\mathcal S^{2\epsilon}_{{\frak r}_2})^{-1}(0)$
such that
\begin{equation}\label{formula810}
(d_{y_2}f^3_{{\frak r}_2})(\tilde{\frak v_2}) = \frak v.
\end{equation}
Then there exists
$\tilde{\frak v_1}\in
T_{(y_1,\xi_1)}(\mathcal S^{1\epsilon}_{{\frak r}_1})^{-1}(0)$
such that
\begin{equation}\label{formula811}
(d_{y_1}f^1_{{\frak r}_1})(\tilde{\frak v_1})
=
(d_{y_2}f^2_{{\frak r}_2})(\tilde{\frak v_2}).
\end{equation}
(\ref{formula811}) implies that
$$(\tilde{\frak v_1},\tilde{\frak v_2})
\in
T_{((y_1),(y_2))}
(
(\mathcal S^{1\epsilon}_{{\frak r}_1})^{-1}(0) \,{}_{f^1}\times_{f^2}
(\mathcal S^{2\epsilon}_{{\frak r}_1})^{-1}(0)
)
$$
and (\ref{formula810}) implies that
$$
(d_{((y_1),(y_2))} \overline{f^3})
(\tilde{\frak v_1},\tilde{\frak v_2})
= \frak v.
$$
Here
$\overline{f^3} :
(\mathcal S^{1\epsilon}_{{\frak r}_1})^{-1}(0) \,{}_{f^1}\times_{f^2}
(\mathcal S^{2\epsilon}_{{\frak r}_1})^{-1}(0)\to N$
is a local representative of $\widetilde{\widehat{f^3}}$.
We have thus proved the required submersivity.
\par
The proof of (3) is obvious from the definition.
\end{proof}
\subsection{Composition of smooth correspondences}
\label{subsec:compsmcor}
In this subsection we define composition of smooth correspondences
and its perturbation. Let us consider the following situation.

\begin{shitu}\label{compositu}
Let $(X_{21},\widehat{\mathcal U_{21}})$,
$(X_{32},\widehat{\mathcal U_{32}})$ be K-spaces
and $M_{i}$ ($i=1,2,3$) smooth manifolds.
Let
$$
\widehat{f_{i,ji}} : (X_{ji},\widehat{\mathcal U_{ji}}) \to M_i,
\qquad
\widehat{f_{j,ji}} : (X_{ji},\widehat{\mathcal U_{ji}}) \to M_j
$$
be strongly smooth maps for $(i,j) = (1,2)$ or $(2,3)$.
We assume
$\widehat{f_{2,21}}$ and $\widehat{f_{3,32}}$ are weakly submersive.
These facts imply that
$$
\frak X_{i+1 i} = ((X_{i+1 i},\widehat{\mathcal U_{i+1 i}}),\widehat{f_{i,i+1 i}},
\widehat{f_{i+1,i+1 i}})
$$
is a smooth correspondence from $M_i$ to $M_{i+1}$
for $i=1,2$.
(Lemma \ref{lem838} (2).)
Let $\widehat{\frak S_{i+1 i}}$ be a CF-
perturbation of $(X_{i+1 i},\widehat{\mathcal U_{i+1 i}})$
for each $i=1,2$.
We assume that $\widehat{f_{i+1,i+1 i}}$ is  strongly submersive
with respect to $\widehat{\frak S_{i+1 i}}$ for each $i=1,2$.$\blacksquare$
\end{shitu}
\begin{defn}
In Situation \ref{compositu}, we put
\begin{equation}
X_{31} = X_{21} \times_{M_2} X_{32}
=
\{(x_{21},x_{32}) \in X_{21} \times X_{32}\mid f_{2,21}(x_{21}) = f_{2,32}(x_{32})\}.
\end{equation}
We put the fiber product Kuranishi structure
\begin{equation}\label{compkurafiber}
\widehat{\mathcal U_{31}}
=
\widehat{\mathcal U_{21}} \times_{M_2} \widehat{\mathcal U_{32}}
\end{equation}
on $X_{31}$ and
define
\begin{equation}\label{mapf131}
\widehat{f_{1,31}} : (X_{31},\widehat{\mathcal U_{31}}) \to M_1,
\qquad
\widehat{f_{3,31}} : (X_{31},\widehat{\mathcal U_{31}}) \to M_3
\end{equation}
as the compositions
$$
(X_{31},\widehat{\mathcal U_{31}}) \to (X_{21},\widehat{\mathcal U_{21}}) \to M_1,
\qquad
(X_{31},\widehat{\mathcal U_{31}}) \to (X_{32},\widehat{\mathcal U_{32}}) \to M_3,
$$
where the first arrows are obvious projections.
We write
$$
\frak X_{21} \times_{M_2} \frak X_{32}
=
((X_{31},\widehat{\mathcal U_{31}}),\widehat{f_{1,31}},\widehat{f_{3,31}})
$$
and call it the {\it composition of smooth correspondences}
\index{smooth correspondence ! composition}
$\frak X_{21}$ and $ \frak X_{32}$.
We also denote it by $\frak X_{32} \circ \frak X_{21}$.
\begin{equation}\label{diagram1010}
\xymatrix{
&& \frak X_{31} \ar[ld]\ar[rd] \\
& \frak X_{21} \ar[ld]\ar[rd] && \frak X_{32}\ar[ld]\ar[rd] \\
M_1 && M_2 && M_3}
\end{equation}
\end{defn}
\begin{rem}
Note that we did not define the `maps' $\frak X_{31} \to \frak X_{21}$,
$\frak X_{31} \to \frak X_{32}$ in Diagram (\ref{diagram1010}).
This is because we never defined the notion of morphism between K-spaces
in this document.
However, the maps $\frak X_{31} \to M_1$ and $\frak X_{31} \to M_3$
are defined by composing the map $f_{21,p}$ or $f_{32,q}$ and the
projection on each chart.
\end{rem}
\begin{lem}
\begin{enumerate}
\item
The fiber product (\ref{compkurafiber}) is well-defined.
\item
The map
$(X_{31},\widehat{\mathcal U_{31}}) \to M_3$
is weakly submersive.
\end{enumerate}
\end{lem}
\begin{proof}
(1) By assumption, $\widehat{f_{2,21}}$ is weakly submersive.
This implies well-defined-ness of  (\ref{compkurafiber}).
\par
(2) Let $(p,q) \in X_{31}$, i.e., $p\in X_{21}$, $q \in X_{32}$
$f_{2,21}(p) = f_{2,32}(q)$.
We put $x = f_{2,21}(p) = f_{3;32}(q)$ and $y \in f_{3,32}(q)$.
By assumption
$$
d_{o_p}(f_{2,21})_p : T_{o_p}U_p \to T_xM_2,
\quad
d_{o_q}(f_{3,32})_q : T_{o_q}U_q \to T_yM_3
$$
are surjective. Let $v_3 \in T_yM_3$.
There exists $\tilde v_3 \in T_{o_q}U_q$ such that
$(d_{o_q}(f_{3,32})_q)(\tilde v_3)  = v_3$.
There exists $\tilde v_2 \in T_{o_p}U_p$ such that
$(d_{o_p}(f_{2,21})_p)(\tilde v_2) = (d_{o_q}(f_{2,32})_q)(\tilde v_3)$.
Then $(\tilde v_2,\tilde v_3) \in T_{o_{(p,q)}}U_{(p,q)}$ and
$(d_{o_{(p,q)}}(f_{3,31})_{(p,q)})(\tilde v_2,\tilde v_3) = v_3$ as required.
\end{proof}
\begin{defn}\label{defn839}
\begin{enumerate}
\item
Let
$\frak X = ((X,\widehat{\mathcal U}),\widehat{f_s},\widehat{f_t})$ be a smooth correspondence
and $\widehat{\frak S}$ a CF-perturbation of
$(X,\widehat{\mathcal U})$.
We say that $(X,\widehat{\mathcal U},\widehat{\frak S},\widehat{f_s},
\widehat{f_t})$
is a {\it perturbed smooth correspondence}
\index{smooth correspondence ! perturbed smooth correspondence}
if $\widehat{f_t}$ is strongly submersive with respect to
$\widehat{\frak S}$.
\item
A perturbed smooth correspondence
$\tilde{\frak X} = (\frak X,\widehat{\frak S}) =
(X,\widehat{\mathcal U},\widehat{\frak S},\widehat{f_s},\widehat{f_t})$
from $M_s$ to $M_t$ induces a linear map
$
\Omega^*(M_s) \to \Omega^{* + \deg\frak X} (M_t)
$
by (\ref{defn924}).
We write it as ${\rm Corr}^{\epsilon}_{\tilde{\frak X}}$.
\item
In Situation \ref{compositu}, let
$\widehat{{\frak S}_{i+1 i}}$ be a CF-perturbation of
$(X_{i+1 i},\widehat{{\mathcal U}_{i+1 i}})$ for each $i=1,2$.
Suppose
$\tilde{\frak X}_{i+1 i} = (X_{i+1 i},\widehat{\mathcal U_{i+1 i}},\widehat{\frak S_{i+1 i}},\widehat{f_{i,i+1 i}},\widehat{f_{i+1, i+1 i}})$ is a perturbed smooth correspondence for each $i=1,2$.
Then by Lemma \ref{lem838}
\begin{equation}\label{fppersmcorr}
\aligned
(X_{21}\,{}_{f_{2,21}}\times_{ f_{2,32}}  X_{32},\,
&\widehat{\mathcal U_{21}}\,{}_{\widehat{f_{2,21}}}\times _{\widehat{f_{2,32}}},\widehat{\mathcal U_{32}},\\
&(\widehat{{\frak S}_{21}}) {}_{\widehat {f_{2,21}}} \times_{\widehat{f_{2,32}}}  (\widehat{\frak S_{32}}),\widehat{f_{1,31}},\widehat{f_{3,31}})
\endaligned
\end{equation}
is a perturbed smooth correspondence from $M_1$ to $M_3$.
We call (\ref{fppersmcorr}) the {\it composition}
\index{smooth correspondence ! composition of perturbed smooth correspondences} of
$\tilde{\frak X}_{21}$ and $\tilde{\frak X}_{32}$
and write
$\tilde{\frak X}_{32} \circ \tilde{\frak X}_{21}$.
Here
\begin{equation}
\aligned
& \widehat{f_{1,31}} : (X_{21}\,{}_{f_{2,21}}\times_{ f_{2,32}}  X_{32},\,
\widehat{\mathcal U_{21}}\,{}_{\widehat{f_{2,21}}}\times _{\widehat{f_{2,32}}},\widehat{\mathcal U_{32}}) \to M_1, \\
& \widehat{f_{3,31}} : (X_{21}\,{}_{f_{2,21}}\times_{f_{2,32}}  X_{32},\,
\widehat{\mathcal U_{21}}\,{}_{\widehat{f_{2,21}}}\times _{\widehat{f_{2,32}}},\widehat{\mathcal U_{32}}) \to M_3
\endaligned
\end{equation}
are maps as in \eqref{mapf131}.
\end{enumerate}
\end{defn}

\subsection{Composition formula}
\label{subsec:compformulaub}
The main result of this section is the following.

\begin{thm}\label{compformulaprof}{\rm(Composition formula, \cite[Lemma 12.15]{fooo09})}
Suppose that $\tilde{\frak X}_{i+1 i} =
(X_{i+1 i},\widehat{\mathcal U_{i+1 i}},\widehat{\frak S_{i+1 i}},
\widehat{f_{i,i+1 i}},
\widehat{f_{i+1,i+1 i}})$
are perturbed smooth correspondences for $i=1,2$.
Then
\begin{equation}\label{formula814}
{\rm Corr}^{\epsilon}_{\tilde{\frak X}_{32}\circ\tilde{\frak X}_{21}}
=
{\rm Corr}^{\epsilon}_{\tilde{\frak X}_{32}} \circ {\rm Corr}^{\epsilon}_{\tilde{\frak X}_{21}}
\end{equation}
for each sufficiently small $\epsilon >0$.
\end{thm}
\begin{rem}
Note that ${\rm Corr}^{\epsilon}_{\tilde{\frak X}_{**}}$
depends on the positive number $\epsilon$.
\end{rem}
\begin{proof}
Let $h_1$ and $h_3$ be differential forms on $M_1$ and $M_3$,
respectively.
It suffices to show the next formula.
\begin{equation}\label{fpm109109}
\int_{M_3}
{\rm Corr}^{\epsilon}_{\tilde{\frak X}_{32}\circ\tilde{\frak X}_{21}}(h_1) \wedge h_3
=
\int_{M_3}
{\rm Corr}^{\epsilon}_{\tilde{\frak X}_{32}}({\rm Corr}^{\epsilon}_{\tilde{\frak X}_{21}}(h_1)) \wedge h_3.
\end{equation}
We use the following notation.
\begin{defn}
In Situation \ref{sitsu8main}, we consider the case when $M$ is a point
and put:
$$
\int_{(X,Z,\widehat{\mathcal U},{\widehat{\frak S^{\epsilon}}})}\widehat{ h}
=
\widehat f !(\widehat{h};{\widehat{\frak S^{\epsilon}}}).
$$
We call it the {\it integration of $\widehat h$ over
$(X,Z,\widehat{\mathcal U},{\widehat{\frak S^{\epsilon}}})$.}
\index{integration ! of differential form by
CF-perturbation on Kuranishi structure}
It is a real number depending on
$(X,Z,\widehat{\mathcal U},\widehat{\frak S}), \epsilon$ and $\widehat h$.
We also define
$$
\int_{(X,{\mathcal U},{\frak S}^{\epsilon})} h
=
f !({h};{{\frak S}}^{\epsilon}).
$$
Here ${\mathcal U}$ is a Kuranishi chart of $X$, $h$ is a differential form of compact support on $U$,
${\frak S}^{\epsilon}$ is a CF-perturbation of
$\mathcal U$ on the support of $h$ and $f : U \to $ a point
is a trivial map, such that $f$ is strongly submersive with
respect to ${\frak S}^{\epsilon}$.
(Note the strong submersivity in the case when the map is trivial
is nothing but transversality of the CF-perturbation to $0$.)
We call it the {\it integration of $h$
\index{integration ! of differential form on K-space} over
$(\widehat{\mathcal U},\widehat{\frak S^{\epsilon}})$.}
\par
In the notation above we omit $Z$ if $Z=X$.
\end{defn}
Using this notation the right hand side of (\ref{fpm109109}) is
\begin{equation}
\aligned
&\int_{(X_{32},\widehat{\mathcal U_{32}},\widehat{\frak S_{32}^{\epsilon}})} (\widehat{f_{2,32}})^*({\rm Corr}_{\tilde{\frak X}_{21}}(\widehat{h_1})) \wedge
(\widehat{f_{3,32}})^* \widehat{h_3}\\
&=
\int_{(X_{32},\widehat{\mathcal U_{32}},\widehat{\frak S_{32}^{\epsilon}})} (\widehat{f_{2,32}})^*(\widehat f_{1,21}!(\widehat{f_{1,21}})^*(\widehat{h_1});
\widehat{{\frak S}_{21}^{\epsilon}})) \wedge
(\widehat{f_{3,32}})^* \widehat{h_3}.
\endaligned
\end{equation}
On the other hand the left hand side of (\ref{fpm109109}) is
\begin{equation}
\int_{(X_{31},\widehat{\mathcal U_{31}},{\widehat{\frak S^{\epsilon}_{31}}})}
(\widehat{f_{1,31}})^*(\widehat{h_1}) \wedge (\widehat{f_{3,31}})^*(\widehat{h_3}).
\end{equation}
Therefore (\ref{fpm109109}) follows from the next proposition.
\begin{prop}\label{fubiniprop}
For $i=1,2$, let $\widehat{\mathcal U_i}$ be Kuranishi structures of
$Z_i \subseteq X_i$,
$\widehat{\frak S_i}$ their CF-perturbations,
$\widehat h_i$  smooth differential forms on $(X_i,\widehat{\mathcal U_i})$
which have compact support in $\ring Z_i$, and
$\widehat f_i : (X_i,Z_i;\widehat{\mathcal U_i}) \to M$
strongly smooth maps.
Supposed that $\widehat f_1$ is strongly submersive with respect to
$\widehat{\frak S_1}$ and
$\widehat{\frak S_2}$ is transversal to $0$
and denote by $(X,Z,\widehat{\mathcal U},\widehat{\frak S})$
the fiber product
$$
(X_1,Z_1,\widehat{\mathcal U_1},\widehat{\frak S_1})
\,{}_{\widehat{f_1}}
\times_{\widehat{f_2}}\,\,
(X_2,Z_2,\widehat{\mathcal U_2},\widehat{\frak S_2})
$$
over $M$. Then
\begin{equation}\label{fubini}
\int_{(X,Z,\widehat{\mathcal U},\widehat{\frak S^{\epsilon}})} \widehat{h_1}\wedge
\widehat{h_2}
=
\int_{(X_2,Z_2,\widehat{\mathcal U}_2,\widehat{\frak S_2^{\epsilon}})}
(\widehat{f_2})^*(\widehat{f_1}!(\widehat{h_1};\widehat{\frak S_1^{\epsilon}}))
\wedge \widehat{h_2}.
\end{equation}
\end{prop}
\begin{rem}
We may regard Formula (\ref{fubini}) as a version of Fubini formula.
\end{rem}
$$
\xymatrix{
& (X,Z,\widehat{\mathcal U},\widehat{\frak S}) \ar[ld]\ar[rd]\\
(X_1,Z_1,\widehat{\mathcal U_1},\widehat{\frak S_1}) \ar[rd]^{\widehat{f_1}} && (X_2,Z_2,\widehat{\mathcal U_2},\widehat{\frak S_2}) \ar[ld]_{\widehat{f_2}} \\
& M
}
$$
\begin{proof}[Proof of Proposition \ref{fubiniprop}]
For $i=1,2$
let ${\widetriangle{\mathcal U_i}}$ be good coordinate systems of
$X_i$ and $\widehat{\mathcal U_i} \to {\widetriangle{\mathcal U_i}}$
KG-embeddings.
We may assume that there exist CF-perturbations ${\widetriangle{\frak S_i}}$ of
$(X_i,Z_i;{\widetriangle{\mathcal U_i}})$ such that
${\widetriangle{\frak S_i}}$, ${\widehat{\frak S_i}}$
are compatible with the KG-embeddings
and $\widehat{h_i}$, $\widehat{f_i}$ are pull back of
differential forms on ${\widetriangle{\mathcal U_i}}$
and of strongly smooth maps on ${\widetriangle{\mathcal U_i}}$,
which we also denote by $\widetriangle{h_i}$, $\widetriangle{f_i}$, respectively.
(Theorem \ref{Them71restate}, Proposition \ref{le614} (2), Theorem \ref{existperturbcont}.)
\par
Let $\{\chi^i_{\frak p_i}\}$ be strongly smooth partitions of unity
of $(X_i,{\widetriangle{\mathcal U_i}})$.
The functions $\chi^i_{\frak p_i}$ can be regarded as strongly smooth maps
$(X_i,Z_i;{\widetriangle{\mathcal U_i}}) \to \R$.
Therefore they induce strongly smooth maps
$(X_i,Z_i;{\widehat{\mathcal U_i}}) \to \R$, which we also denote by
$\chi^i_{\frak p_i}$.
Then they induce strongly smooth functions
on the fiber product
$(X,Z;\widehat{\mathcal U})$.
\par
To prove (\ref{fubini}) it suffices to prove the next formula for each
$\frak p_1$, $\frak p_2$ and $\epsilon$.
\begin{equation}\label{fubinionechart1}
\aligned
&\int_{(X,Z,\widehat{\mathcal U},\widehat{\frak S^{\epsilon}})}
\chi^1_{\frak p_1}\widehat{h_1}\wedge
\chi^2_{\frak p_2}\widehat{h_2} \\
&=
\int_{({\mathcal U}_{2,\frak p_2},{\frak S}^{\epsilon}_{2,\frak p_2})}
f_{2,\frak p_2}^*(f_{1,\frak p_1}!
(\chi^1_{\frak p_1}h_{1,\frak p_1};{\frak S}_{1,\frak p_1}^{\epsilon}))
\wedge \chi^2_{\frak p_2} h_{2,\frak p_2}.
\endaligned
\end{equation}
We will use the following result.

\begin{prop}\label{lemma8845}
Let $Z_{(1)}$ and $Z_{(2)}$ be compact subsets of $X$ such that
$Z_{(1)} \subset \ring Z_{(2)}$.
Let ${\widehat{\mathcal U}}$ be a
Kuranishi structure of $Z_{(2)} \subseteq X$
and $\widehat h$ a differential form
on ${\widehat{\mathcal U}}$.
Let $\widehat{f_{(2)}} : (X,Z_{(2)};{\widehat{\mathcal U}}) \to M$ be
a strongly smooth map and
$\widehat{\frak S_{(2)}}$ a CF-perturbation of ${\widehat{\mathcal U}}$.
Denote by $\widehat{\frak S_{(1)}}$
the restriction of $\widehat{\frak S_{(2)}}$ to ${\widehat{\mathcal U}\vert_{Z_{(1)}}}$ (Definition \ref{defn71717}).
Suppose
\begin{enumerate}
\item
$\widehat{f_{(2)}}$ is strongly submersive with respect to $\widehat{\frak S_{(2)}}$.
\item
$\widehat h$ has compact support in $\ring Z_{(1)}$.
\end{enumerate}
We denote by $\widehat h\vert_{Z_{(1)}}$ the
differential form on ${\widehat{\mathcal U}\vert_{Z_{(1)}}}$ induced by $\widehat h$
via the condition (2). Then
\begin{equation}
\widehat{f_{(2)}}!(\widehat h;\widehat{\frak S_{(2)}^{\epsilon}})
=
\widehat{f_{(1)}}!(\widehat h\vert_{Z_{(1)}};\widehat{\frak S_{(1)}^{\epsilon}}).
\end{equation}
where $\widehat{f_{(1)}}$ is the restriction of $\widehat{f_{(2)}}$ to
$(X,Z_{(1)};{\widehat{\mathcal U}\vert_{Z_{(1)}}})$.
\end{prop}
\begin{proof}
By using differential forms $\rho$ on $M$ it suffices
to consider the case $\deg h = \dim \widehat{\mathcal U}$ and
prove the next formula (see Lemma \ref{lem782} (2))
\begin{equation}\label{form1014}
\int_{(X,Z_{(2)},\widehat{\mathcal U},\widehat{\frak S_{(2)}^{\epsilon}})} \widehat h
=
\int_{(X,Z_{(1)},\widehat{\mathcal U}\vert_{Z_{(1)}},\widehat{\frak S_{(1)}^{\epsilon}})} \widehat h.
\end{equation}
We take good coordinate systems $\widetriangle{\mathcal U_{(i)}}$
of $Z_{(i)} \subseteq X$ such that:
\begin{enumerate}
\item
$\widetriangle{\mathcal U_{(1)}}$ is compatible with
$\widehat{\mathcal U}\vert_{Z_{(1)}}$
and $\widetriangle{\mathcal U_{(2)}}$ is compatible with
$\widehat{\mathcal U}$.
\item
$\widetriangle{\mathcal U_{(2)}}$ strictly extends
$\widetriangle{\mathcal U_{(1)}}$.
\item
There are CF-perturbations $\widetriangle{\frak S_{(i)}}$,
differential forms $\widetriangle{h_{(i)}}$ on $(X,Z_{(i)};\widetriangle{\mathcal U_{(i)}})$
which are compatible each other and are
compatible with corresponding objects on $\widetriangle{\mathcal U_{(1)}}$
and on $\widehat{\mathcal U}\vert_{Z_{(1)}}$.
\item
$\widetriangle{\frak S_i}$ are transversal to $0$.
\end{enumerate}
Existence of such objects is a consequence of Theorem \ref{Them71restate},
Propositions \ref{prop7582752} and \ref{existperturbcontrel}.
\par
Let $\widetriangle{\mathcal U_{(i)}} = (\frak P_{(i)},\{\mathcal U_{{(i)},\frak p}\},
\{\Phi_{{(i)},\frak p\frak q}\})$.
By Definition \ref{defn735f} (1)(a),
$
\frak P_{(1)} = \{\frak p \in \frak P_{(2)} \mid {\rm Im}(\psi_{{(2)},\frak p})
\cap Z_{(1)} \ne \emptyset\}
$
and $\mathcal U_{{(1)},\frak p}$ is an open subchart
of $\mathcal U_{{(2)},\frak p}$ for $\frak p \in \frak P_{(1)} \subset \frak P_{(2)}$.
We choose support systems $\mathcal K^{(i)}$
and a neighborhood $\frak U_{1}(Z_{(1)})$ of $Z_{(1)}$ in
$\vert\widetriangle{\mathcal U_{(2)}} \vert$ with the following properties.
\begin{enumerate}
\item[(a)]
If $\frak p \in \frak P_{(1)}\subset \frak P_{(2)}$
then $\mathcal K^{(1)}_{\frak p} \cap \frak U_{(1)}(Z_{(1)})= \mathcal K^{(2)}_{\frak p}
\cap \frak U_{1}(Z_{(1)})$.
\item[(b)]
If $\frak p \in \frak P_{(2)} \setminus \frak P_{(1)}$
then $\frak U_{1}(Z_{(1)}) \cap \mathcal K_{\frak p}^{(2)}= \emptyset$.
\end{enumerate}
We take support pairs $(\mathcal K^{1,{(i)}},\mathcal K^{2,{(i)}})$
of $\widetriangle{\mathcal U^{(i)}}$
such that $\mathcal K^{2,{(i)}} < \mathcal K^{{(i)}}$.
By (a), we may assume that there exists a neighborhood $\frak U_2(Z_{(1)})$ of $Z_{(1)}$ in
$\vert\widetriangle{\mathcal U_{(2)}} \vert$ such that
\begin{enumerate}
\item[(c)]
If $\frak p \in \frak P_{(1)} \subset \frak P_{(2)}$
then $\mathcal K^{j,(1)}_{\frak p} \cap \frak U_2(Z_{(1)})= \mathcal K^{j,(2)}_{\frak p}
\cap \frak U_2(Z_{(1)})$ for $j=1,2$.
\end{enumerate}
We can take $\delta_{(i)}$ and partition of unities $\{\chi^{(i)}_{\frak p}\}$
of $(X,Z_{(i)},{\widetriangle{\mathcal U_{(i)}}},\mathcal K^{1,{(i)}},\delta_{(i)})$
and a neighborhood $\frak U_3(Z_{(1)})$ of $Z_{(1)}$ in
$\vert\widetriangle{\mathcal U_{(2)}} \vert$
for $i=1,2$,
such that:
\begin{enumerate}
\item[(d)]
If $\frak p \in \frak P_{(1)} \subset \frak P_{(2)}$
then $\chi^{(1)}_{\frak p} = \chi^{(2)}_{\frak p}$ on $\frak U_3(Z_{(1)})$.
\end{enumerate}
By definition we have
\begin{equation}\label{form1015}
\int_{(X,Z_{(i)},\widehat{\mathcal U_{(i)}},\widehat{\frak S_{(i)}^{\epsilon}})}
\widehat h
=
\sum_{\frak p \in \frak P_{(i)}}
\int_{\mathcal K^{1,{(i)}}_{\frak p}(2\delta_{(i)}) \cap \frak U(Z_{(i)})} \chi^{(i)}_{\frak p} h_{\frak p}
\end{equation}
for sufficiently small neighborhoods $\frak U(Z_{(i)})$ of $Z_{(i)}$ in $\vert\widetriangle{\mathcal U_{(i)}}\vert$.
By (a)(b)(c)(d) above and using the fact that $\widehat h$ has
a compact support in $\ring Z_{(1)}$, (\ref{form1015})
implies (\ref{form1014}).
In fact, if $\frak p \in \frak P_{(2)} \setminus \frak P_{(1)}$
then (b) implies that
$$
\int_{\mathcal K^{1,{(2)}}_{\frak p}(2\delta_{(2)}) \cap \frak U(Z_{(2)})} \chi^{(2)}_{\frak p} h_{\frak p} = 0
$$
and if $\frak p \in \frak P_{(1)}$ then (c)(d) imply
$$
\int_{\mathcal K^{1,{(1)}}_{\frak p}(2\delta_{(1)}) \cap \frak U(Z_{(1)})} \chi^{(1)}_{\frak p} h_{\frak p} =
\int_{\mathcal K^{1,{(2)}}_{\frak p}(2\delta_{(2)}) \cap \frak U(Z_{(2)})} \chi^{(2)}_{\frak p} h_{\frak p}.
$$
Thus the proof of Proposition \ref{lemma8845} is complete.
\end{proof}
We continue the proof of Proposition \ref{fubiniprop}.
We consider the fiber product Kuranishi chart
$\mathcal U_{1,\frak p_1} \,{}_{f_{1,\frak p_1}}\times_{f_{2,\frak p_2}}
\, \mathcal U_{2,\frak p_2}$
and the fiber product CF-perturbation
$
\frak S_{1,\frak p_1}\,{}_{f_{1,\frak p_1}}\times_{f_{2,\frak p_2}}
\, \frak S_{2,\frak p_2}
$
of it.
\par
We apply Lemma \ref{lemma8845} to $\widehat h = \chi^1_{\frak p_1}\widehat h_1\wedge
\chi^2_{\frak p_2}\widehat h_2$, $Z_{(2)} = Z$ and
$$
Z_{(1)} =
(\psi_{1,\frak p_1}(U_{1,\frak p_1} \cap s_{1,\frak p_1}^{-1}(0)))\,{}_{f_{1}}\times_{f_{2}}
\, (\psi_{2,\frak p_2}(U_{2,\frak p_2} \cap s_{2,\frak p_2}^{-1}(0))).
$$
Note
there exists a good coordinate system
consisting of a single Kuranishi chart  $\mathcal U_{1,\frak p_1} \,{}_{f_{1,\frak p_1}}\times_{f_{2,\frak p_2}}
\, \mathcal U_{2,\frak p_2}$
of $Z_{(1)} \subseteq X_1 \times_M X_2$.
Therefore the left hand side of (\ref{fubinionechart1}) is equal to
$$
\int_{(\mathcal U_{1,\frak p_1} \,{}_{f_{1,\frak p_1}}\times_{f_{2,\frak p_2}}
\, \mathcal U_{2,\frak p_2},\frak S^{\epsilon}_{1,\frak p_1} \,{}_{f_{1,\frak p_1}}\times_{f_{2,\frak p_2}}
\, \frak S^{\epsilon}_{2,\frak p_2})}
\chi^1_{\frak p_1}\widehat h_1\wedge
\chi^2_{\frak p_2}\widehat h_2.
$$
Thus to prove (\ref{fubinionechart1}) it suffices to prove the next lemma.
\begin{lem}\label{lem846}
Let $\mathcal U_i$ be Kuranishi charts of $X_i$,
$f_i : U_i \to M$ smooth maps,
$h_i$  differential forms on $U_i$ and
$\frak S_i$ CF-perturbations of $\mathcal U_i$, for $i=1,2$.
We assume that $f_1$ is strongly submersive
with respect to $\frak S_1$ and
$\frak S_2$ is transversal to $0$.
Then
\begin{equation}\label{formula819}
\int_{(\mathcal U_1 \,{}_{f_1}\times_{f_2}\, \mathcal U_2,\frak S^{\epsilon}_1 \,{}_{f_1}\times_{f_2}\, \frak S_2^{\epsilon})}
h_1\wedge h_2
=
\int_{(\mathcal U_2,\frak S^{\epsilon}_2)} f_2^*
(
f_1!(h_1;\frak S_1^{\epsilon})
)
\wedge h_2.
\end{equation}
\end{lem}
\begin{proof}[Proof of Lemma \ref{lem846}]
Let $\frak S_i = (\{\frak V^i_{\frak r_i}\},\{\mathcal S^i_{\frak r_i}\})$
and $\{\chi_{\frak r_i}^i\}$ be a smooth partition of unities
of orbifolds $U_i$ subordinate to its open cover
$\{U^i_{\frak r_i}\}$.
Then $\{\chi^1_{\frak r_1}\chi^2_{\frak r_2} \mid \frak r_1 \in \frak R_1,
\frak r_2 \in \frak R_2\}$
is a partition of unity subordinate to the covering
$\{U_1\,{}_{f_1}\times_{f_2} U_2 \mid \frak r_1 \in \frak R_1,
\frak r_2 \in \frak R_2\}$.
Therefore, to prove (\ref{formula819}), it suffices to prove
\begin{equation}\label{formula820}
\aligned
&\int_{(\mathcal U^1_{\frak r_1} \,{}_{f_1}\times_{f_2}\, \mathcal U^2_{\frak r_2},\mathcal S^{1\epsilon}_{\frak r_1} \,{}_{f_1}\times_{f_2}\, \mathcal S^{2\epsilon}_{\frak r_2})}
\chi_{\frak r_1}^1h_1\wedge \chi_{\frak r_2}^2 h_2\\
&=
\int_{(\mathcal U^2_{\frak r_2},\mathcal S^{2\epsilon}_{\frak r_2})}
\chi_{\frak r_1}^1\chi_{\frak r_2}^2  f_2^*
(
f_1!(h_1;\mathcal S^{1\epsilon}_{\frak r_1})
)\wedge h_2.
\endaligned
\end{equation}
(\ref{formula820}) is an immediate consequence of the next lemma.
\begin{lem}\label{lem847}
For $i=1,2$ let $N_i$ be smooth manifolds and $f_i : N_i \to M$ smooth maps
and $h_i$ smooth differential forms on $N_i$ of compact support.
Suppose that $f_1$ is a submersion.
Then we have
\begin{equation}\label{form10821}
\int_{N_1\,{}_{f_1}\times_{f_2} N_2} h_1 \wedge h_2
=
\int_{N_2} f_2^* (f_1!(h_1)) \wedge h_2.
\end{equation}
\end{lem}
\begin{proof}
By using a partition of unity, it suffices to prove the lemma
in the following special case:
$N_1 = \R^a \times M$,
$f_1 :  \R^a \times M \to M$ is the projection to the
second factor,
$N_2 =  M \times \R^b$,
$f_2 :  M \times \R^b \to M$ is the projection to the
first factor,  and $M,N_1,N_2$ are open subsets of
Euclidean spaces.
We prove this case below.
\par
In this case $N_1 \,{}_{f_1}\times_{f_2}  N_2 \cong \R^a \times M \times \R^b$.
Let $(x_1,\dots,x_m)$ be a coordinate of $M$,
$(y_1,\dots,y_a)$ a coordinate of $\R^a$ and
$(z_1,\dots,z_b)$ a coordinate of $\R^b$.
Then
$(y_1,\dots,y_a,x_1,\dots,x_m)$ is a coordinate of $N_1$,
$(x_1,\dots,x_m,z_1,\dots,z_b)$ is a coordinate of $N_2$,
and $(y_1,\dots,y_a,x_1,\dots,x_m,z_1,\dots,z_b)$ is a coordinate of
$N_1\,{}_{f_1}\times_{f_2} N_2$.
\par
We may write
$$
\aligned
h_1 &= \sum_{I} g_{1,I}(y_1,\dots,y_a,x_1,\dots,x_m)
dy_1 \wedge \dots \wedge dy_a \wedge  dx_{I}
\\
h_2 &= \sum_{J} g_{2,J}(x_1,\dots,x_m,z_1,\dots,z_b)
dx_J \wedge dz_1 \wedge \dots \wedge dz_b
\endaligned
$$
for certain smooth functions $g_{1,I}, g_{2,J}$.
We write them as $g_{1,I}(y,x), g_{2,J}(x,z)$ for simplicity.
Here $I,J \subset \{ 1, \dots ,m\}$ and $dx_I = dx_{i_1} \wedge \dots \wedge dx_{i_{\vert I \vert}}$ for
$I=\{ i_1, \dots , i_{\vert I \vert} \}$ and $dx_{J}$ is defined in a similar way.
We may assume $I \cup J =\{ 1,\dots , m\}$.
Then the left hand side of (\ref{form10821}) is given by
\begin{equation}
\sum _{I,J}\int_{N_1\,{}_{f_1}\times_{f_2} N_2}
g_{1,I}(y,x) g_{2,J}(x,z)dydx_I dx_J dz.
\end{equation}
On the other hand,  the right hand side of (\ref{form10821}) is given by
\begin{equation}
\sum_{I,J}\int_{N_2}
\left(\int_{y\in \R^a}g_{1,I}(y,f_2(x,z))dy \right)g_{2,J}(x,z) dx_I dx_J dz.
\end{equation}
Therefore Lemma \ref{lem847} is an immediate consequence of
Fubini's theorem.
\end{proof}
This completes the proof of Lemma \ref{lem846}.
\end{proof}
This also completes the proof of Proposition \ref{fubiniprop}.
\end{proof}
Therefore the proof of Theorem \ref{compformulaprof} is now complete.
\end{proof}
We finally remark the following.
\begin{lem}
When CF-perturbations $\tilde{\frak X}_{32}$, $\tilde{\frak X}_{21}$ vary in a
uniform family, we can choose $\epsilon$ in Theorem \ref{compformulaprof} in the  way  independent of the
CF-perturbations in that family.
\end{lem}
The proof goes in the same way as the proof of Proposition \ref{lem761}. So we omit it.

\section{Construction of good coordinate system}
\label{sec:contgoodcoordinate}

In this section we prove Theorem \ref{Them71restate}
together with various addenda and variants.

\subsection{Construction of good coordinate systems: the absolute case}
\label{subsec:constgcsabs}

This subsection will be occupied by the proof of Theorem \ref{Them71restate}.

\begin{proof}[Proof of Theorem \ref{Them71restate}]
\footnote{The proof below is basically the same proof
as written in\cite[Section 7]{foootech}.
(The proof in \cite[Section 7]{foootech}
is a detailed version of one given in \cite[page 957-958]{FO}.)
However we
polish the presentation and reorganize the proof slightly so that
it becomes shorter and easier to read.}
Let $\widehat{\mathcal U}$ be a Kuranishi structure of $Z \subseteq X$.
We use the dimension stratification $\mathcal S_{\frak d}(X,Z;\widehat{\mathcal U})$ as in
Definition \ref{stratadim},
for the inductive construction of a good coordinate
system compatible to the given Kuranishi structure in the sense of Definition \ref{gcsystem}.
The main part of the proof is the proof of Proposition \ref{inductiveprop} below.
We first describe the situation
under which Proposition \ref{inductiveprop} is stated.

\begin{shitu}\label{situation101}
Let $\frak d \in \Z_{\ge 0}$, $Z_0$  a compact subset of
$$
\mathcal S_{\frak d}(X,Z;\widehat{\mathcal U})
\setminus
\bigcup_{\frak d' > \frak d}\mathcal S_{\frak d'}(X,Z;\widehat{\mathcal U}),
$$
and $Z_1$ a compact subset of
$\mathcal S_{\frak d}(X,Z;\widehat{\mathcal U})$. We assume that $Z_1$
contains an open neighborhood of
\begin{equation}\label{biggedZZZ}
\bigcup_{\frak d' > \frak d}\mathcal S_{\frak d'}(X,Z;\widehat{\mathcal U})
\nonumber\end{equation}
in $\mathcal S_{\frak d}(X,Z;\widehat{\mathcal U})$.
\par
We also assume that we have a good coordinate system
${\widetriangle{\mathcal U}}
= (\frak P,\{\mathcal U_{\frak p}\},\{\Phi_{\frak p\frak q}\})$ of
$Z_1^+ \subseteq X$, where $Z^+_1$ is a compact neighborhood
of $Z_1$ in $X$,
and a strict KG-embedding
$$
\widehat{\Phi^1}
= \{\Phi^1_{\frak p p} \mid
p \in {\rm Im}(\psi_{\frak p}) \cap Z^+_1\} : \widehat{\mathcal U}\vert_{Z^+_1} \to {\widetriangle{\mathcal U}},
$$
where $\widehat{\mathcal U}\vert_{Z^+_1}$
is the restriction
of $\widehat{\mathcal U}$,  which is defined in Definition \ref{defn735f} (3).
\par
Let $\mathcal U_{\frak p_0} = (U_{\frak p_0},E_{\frak p_0},s_{\frak p_0},\psi_{\frak p_0})$ be a Kuranishi neighborhood of
$Z^+_0$ such that $\dim U_{\frak p_0} = \frak d$.
Here $Z^+_0$ is a compact neighborhood of $Z_0$ in $X$.
We regard $\mathcal U_{\frak p_0}$  as a good coordinate system
$\widetriangle{\mathcal U_{\frak p_0}}$ that
consists of a single
Kuranishi chart and suppose that we are given a strict KG-embedding
$$
\widehat{\Phi^0}
= \{\Phi^0_{\frak p_0 p} \mid
p \in {\rm Im}(\psi_{\frak p_0}) \cap Z^+_0\} : \widehat{\mathcal U}\vert_{Z^+_0} \to \widetriangle{\mathcal U_{\frak p_0}}.
$$
We put
\begin{equation}\label{form111}
Z_+ = Z_1 \cup Z_0.
\end{equation}
$\blacksquare$
\end{shitu}
\begin{rem}
We remark $Z_+ \subseteq \mathcal S_{\frak d}(X,Z;\widehat{\mathcal U})$.
However in general $Z^+_0$, $Z^+_1$
are not subsets of $\mathcal S_{\frak d}(X,Z;\widehat{\mathcal U})$.
\end{rem}
\begin{prop}\label{inductiveprop}
In Situation \ref{situation101}, there exists a good coordinate system
${\widetriangle{\mathcal U^+}}
= (\frak P^+,\{\mathcal U^+_{\frak p}\},\{\Phi^+_{\frak p\frak q}\})$
of $Z^+_+ \subseteq X$ with the following properties.
Here $Z_+^+$ is a compact neighborhood of $Z_+$ in $X$.
\begin{enumerate}
\item
$\frak P^+ = \frak P \cup \{\frak p_0\}$.
The partial order on $\frak P^+$ is the same as one on $\frak P$
among the elements of $\frak P$ and
$\frak p > \frak p_0$ for any $\frak p \in \frak P$.
\item
If $\frak p \in \frak P$ then $\mathcal U^+_{\frak p}$
is an open subchart $\mathcal U_{\frak p}\vert_{U^+_{\frak p}}$ where
$U^+_{\frak p}$ is an open subset of $U_{\frak p}$.
We have
\begin{equation}
\bigcup_{\frak p \in \frak P}
\psi_{\frak p}^+ ((s_{\frak p}^+)^{-1}(0))
\supset Z_1.
\end{equation}
\item
$\mathcal U^+_{\frak p_0}$ is a restriction of
$\mathcal U_{\frak p_0}$ to an open subset
$U^+_{\frak p_0}$ of $U_{\frak p_0}$.
We have
$$
\psi_{\frak p_0}^+ ((s_{\frak p_0}^+)^{-1}(0))
\supset
Z_0.
$$
\item
The coordinate change $\Phi^+_{\frak p\frak q}$
is the restriction of $\Phi_{\frak p\frak q}$ to
$
U^+_{\frak q} \cap \varphi_{\frak p\frak q}^{-1}(U^+_{\frak p}),
$
if $\frak p,\frak q \in \frak P$.
\item
There exists an open substructure $\widehat{\mathcal U^0}$
of $\widehat{\mathcal U}$, a strict KG-embedding
$\widehat{\Phi^+}
= \{\Phi^+_{\frak p p} \mid
p \in {\rm Im}(\psi^+_{\frak p}) \cap Z^+_+\} : \widehat{\mathcal U^0} \to {\widetriangle{\mathcal U^+}}$
with the following properties.
\par\smallskip
\begin{enumerate}
\item
If $\frak p \in \frak P$ and
$p  \in {\rm Im}(\psi^+_{\frak p}) \cap Z_+^+$,
then we have the commutative diagram:
\begin{equation}\label{diagram103}
\begin{CD}
\mathcal U_{p}^0 @ > {\Phi_{\frak p p}^+} >>
{\mathcal U}_{\frak p}^+  \\
@ VVV @ VVV\\
\mathcal U_{p} @ > {\Phi_{\frak p p}^1} >>{\mathcal U}_{\frak p}
\end{CD}
\end{equation}
where the vertical arrows are embeddings as open subcharts.
(The commutativity of diagram means that the maps coincide when
both sides are defined.)
\item
If
$p\in {\rm Im}(\psi^+_{\frak p_0}) \cap Z^+_+$, then
we have the commutative diagram:
\begin{equation}\label{diagram104}
\begin{CD}
\mathcal U_{p}^0 @ > {\Phi_{\frak p_0 p}^+} >>
{\mathcal U}_{\frak p_0}^+  \\
@ VVV @ VVV\\
\mathcal U_{p} @ > {\Phi_{\frak p_0 p}^0} >>{\mathcal U}_{\frak p_0}
\end{CD}
\end{equation}
where the vertical arrows are embeddings as open subcharts.
(The commutativity of diagram means that the maps coincide when
both sides are defined.)
\end{enumerate}
\end{enumerate}
\end{prop}
\begin{rem}
Note the good coordinate system
${\widetriangle{\mathcal U^+}}$ has
one more Kuranishi chart than ${\widetriangle{\mathcal U}}$.
Moreover ${\widetriangle{\mathcal U^+}}$ is a good coordinate system
of a neighborhood of $Z_{+}$
which contains $Z_1$.
${\widetriangle{\mathcal U}}$
is a good coordinate system
of a neighborhood of $Z_{1}$.
This is the reason why we use the symbol $+$ in
${\widetriangle{\mathcal U^+}}$.
\par
On the other hand, each Kuranishi chart
of ${\widetriangle{\mathcal U^+}}$
is either an open subchart of ${\widetriangle{\mathcal U}}$
or of ${\widetriangle{\mathcal U_{\frak p_0}}}$.
In other words each  Kuranishi chart of
${\widetriangle{\mathcal U^+}}$ is {\it smaller}
than that of ${\widetriangle{\mathcal U}}$
or of ${\widetriangle{\mathcal U_{\frak p_0}}}$.
\end{rem}
\begin{proof}
We take a support system $\mathcal K$ of ${\widetriangle{\mathcal U}}$.
Note by definition we have:
$$
\bigcup_{\frak p} \psi_{\frak p}(\overset{\circ}{\mathcal K_{\frak p}} \cap
s_{\frak p}^{-1}(0))
\supset
Z^+_1.
$$
We take a compact subset $\mathcal K_{\frak p_0}$ of $U_{\frak p_0}$ such that
$$
\psi_{\frak p_0}(\overset{\circ}{\mathcal K_{\frak p_0}} \cap
s_{\frak p_0}^{-1}(0))
\supset
Z^+_0.
$$
Since $
\widehat{\Phi^1}
 : \widehat{\mathcal U}\vert_{Z^+_1} \to {\widetriangle{\mathcal U}}
$
and
$
\widehat{\Phi^0} : \widehat{\mathcal U}\vert_{Z^+_0} \to \widetriangle{\mathcal U_{\frak p_0}}
$
are {\it strict} KG-embeddings, they have the following
properties.
\begin{proper}\label{proper113}
\begin{enumerate}
\item
If $q \in
\psi_{\frak p}(\mathcal K_{\frak p} \cap
s_{\frak p}^{-1}(0))$, then
$\Phi^1_{\frak p q} = ( \varphi^1_{\frak p q},
\widehat{\varphi^1_{\frak p q}})$ is
defined and is an embedding
$\Phi^1_{\frak p q} : \mathcal U_q \to \mathcal U_{\frak p}$.
Note the domain of $\varphi^1_{\frak p q}$ is $U_q$.
We put
$$
 U^1_q = U^1_{\frak p q} = U_q.
$$
\item
If $q \in \psi_{\frak p_0}(\mathcal K_{\frak p_0} \cap
s_{\frak p_0}^{-1}(0))$, then
$\Phi^0_{\frak p_0 q} = (\varphi^0_{\frak p_0 q},
\widehat{\varphi^0_{\frak p_0 q}})$ is
defined and is an embedding
$\Phi^0_{\frak p_0 q} : \mathcal U_q \to \mathcal U_{\frak p_0}$.
Note the domain of $\varphi^0_{\frak p_0 q}$ is $U_q$.
We put
$$
U^0_q = U^0_{\frak p_0 q} = U_q.
$$
\end{enumerate}
\end{proper}
For each $\frak p$ we use
compactness to find a finite  number of points
$q^{\frak p}_1,\dots,q^{\frak p}_{N(\frak p)} \in
 \psi_{\frak p_0}(\mathcal K_{\frak p_0} \cap
s_{\frak p_0}^{-1}(0)) \cap
\mathcal S_{\frak d}(X,Z;\widehat{\mathcal U})
\cap \psi_{\frak p}(\mathcal K_{\frak p} \cap
s_{\frak p}^{-1}(0))$ such that
\begin{equation}
\bigcup_{i=1}^{N(\frak p)}
\psi_{q_i^{\frak p}}(s_{q^{\frak p}_i}^{-1}(0) \cap U^1_{q^{\frak p}_i})
\supset
\psi_{\frak p_0}(\mathcal K_{\frak p_0} \cap
s_{\frak p_0}^{-1}(0))
\cap
\mathcal S_{\frak d}(X,Z;\widehat{\mathcal U})
\cap \psi_{\frak p}(\mathcal K_{\frak p} \cap
s_{\frak p}^{-1}(0)).
\nonumber\end{equation}
We take relatively compact open subsets $U^2_{q_i^{\frak p}}$ of $U^1_{q_i^{\frak p}}$
such that
\begin{equation}\label{qcoversuru}
\bigcup_{i=1}^{N(\frak p)}
\psi_{q^{\frak p}_i}(s_{q^{\frak p}_i}^{-1}(0) \cap U^2_{q^{\frak p}_i})
\supset
\psi_{\frak p_0}(\mathcal K_{\frak p_0} \cap
s_{\frak p_0}^{-1}(0))
\cap
\mathcal S_{\frak d}(X,Z;\widehat{\mathcal U})
\cap
\psi_{\frak p}(\mathcal K_{\frak p} \cap
s_{\frak p}^{-1}(0)).
\end{equation}
We consider the following diagram:
$$
\xymatrix{
\mathcal U^1_{\frak p_0} \ar@{.>}[rr]  && \mathcal U^1_{\frak p}\\
& \mathcal U^2_{{q_i^{\frak p}}}\ar[lu]^{\Phi^0_{\frak p_0 q_i^{\frak p}}}
\ar[ru]_{\Phi^1_{\frak p q_i^{\frak p}}}
}
$$
Here $\mathcal U^2_{{q_i^{\frak p}}}$ is the restriction
of $\mathcal U_{{q_i^{\frak p}}}$ to $U^2_{{q_i^{\frak p}}}$.
\par
Note $\Phi^0_{\frak p_0 q_i^{\frak p}}$ is locally
invertible since $\varphi^0_{\frak p_0 q_i^{\frak p}}$ is
an orbifold embedding between orbifolds of the same dimension.
So we can find an embedding (from an appropriate open subchart)
written in dotted arrow in the diagram.
Those maps for various $q_i^{\frak p}$
however may not coincide on the overlapped part.
We use the next lemma to shrink the domains so that those
maps are glued to define a coordinate change
from an open subchart of $\mathcal U^1_{\frak p_0}$
to an open subchart of $\mathcal U^1_{\frak p}$, which we need to
find to prove Proposition \ref{inductiveprop}.
\par
Note that
Property \ref{proper113} also holds when we replace $q$ by $r$.

\begin{lem}{\rm (Compare \cite[Lemma 7.8]{foootech})}\label{lemma114}
For each $r \in \psi_{\frak p_0}(\mathcal K_{\frak p_0} \cap
s_{\frak p_0}^{-1}(0))
\cap
\mathcal S_{\frak d}(X,Z;\widehat{\mathcal U})
\cap \bigcup_{\frak p} \psi_{\frak p}(\mathcal K_{\frak p} \cap
s_{\frak p}^{-1}(0))$
there exists an open neighborhood $U^3_{r}$ of $o_r$ in
$U_r = U^1_{\frak p r}$ with the following properties.
\begin{enumerate}
\item
If $r \in \psi_{\frak p}(s_{\frak p}^{-1}(0) \cap \mathcal K_{\frak p})$
and
$
\varphi^1_{\frak p r}(U^3_r) \cap \overline{\varphi^1_{\frak p q^{\frak p}_i}(U^2_{q_i^{\frak p}}})
\ne \emptyset
$,
then
\begin{equation}\label{115555}
U^3_{r} \subset U_{q^{\frak p}_i r} \cap \varphi^{-1}_{q^{\frak p}_i r}
(U^1_{\frak p q^{\frak p}_i}).
\end{equation}
\item
If
$
\varphi^0_{\frak p_0 r}(U^3_r) \cap \overline{\varphi^0_{\frak p_0 q^{\frak p}_i}(U^2_{q^{\frak p}_i}})
\ne \emptyset
$,
then
$$
U^3_{r} \subset  U_{q^{\frak p}_i r} \cap \varphi^{-1}_{q^{\frak p}_i r}(U^0_{\frak p_0 q^{\frak p}_i}).
$$
\end{enumerate}
\end{lem}
\begin{proof}
We observe the following 3 points.
\begin{enumerate}
\item[(a)]
(1)(2) above require a finite number of conditions for each $U_r^3$.
\item[(b)]
If a choice of $U_r^3$ satisfies one of those (finitely many) conditions,
then any smaller $U_r^3$ satisfies the same condition.
\item[(c)]
For any one of those (finitely many) conditions,
there exists $U_r^3$ which satisfies that condition.
\end{enumerate}
In fact, (c) is proved in the case of Lemma \ref{lemma114} (1),
for example, as follows.
Suppose $r \in \psi_{\frak p_0}(\mathcal K_{\frak p_0} \cap
s_{\frak p_0}^{-1}(0))
\cap
\mathcal S_{\frak d}(X,Z;\widehat{\mathcal U})
\cap \bigcup_{\frak p} \psi_{\frak p}(\mathcal K_{\frak p} \cap
s_{\frak p}^{-1}(0))$.
Let $\overline{U^2_{q^{\frak p}_i}}$ be the closure of
${U^2_{q^{\frak p}_i}}$ in ${U^1_{q^{\frak p}_i}}$.
This is a compact set.
\par
Suppose $o_r \notin   U_{q^{\frak p}_i r}
\cap \varphi_{q^{\frak p}_i r}^{-1}(\overline{U^2_{q^{\frak p}_i}})$.
Since $o_r \in U_{q^{\frak p}_i r}$, we have
$\varphi_{q^{\frak p}_i r}(o_r) \notin \overline{U^2_{q^{\frak p}_i}}$.
Therefore we can find $U^3_{r}$
such that
$\varphi^1_{\frak p r}(U^3_r) \cap \overline{\varphi^1_{\frak p q^{\frak p}_i}(U^2_{q_i^{\frak p}}}) = \emptyset$.
Thus Lemma \ref{lemma114} (1) is satisfied in this case.
\par
Suppose $o_r \in   U_{q^{\frak p}_i r}
\cap \varphi_{q^{\frak p}_i r}^{-1}(\overline{U^2_{q^{\frak p}_i}})$.
Since $\overline{U^2_{q^{\frak p}_i}} \subset {U^1_{q^{\frak p}_i}}$,
we can find $U^3_{r}$  satisfying (\ref{115555}).
\par
The lemma follows immediately from (a)(b)(c).
\end{proof}

For each $\frak p \in \frak P$ we choose
$J(\frak p)$ points
$$
r^{\frak p}_{1},\dots,r^{\frak p}_{J(\frak p)}
\in
 \psi_{\frak p_0}(s_{\frak p_0}^{-1}(0)  \cap \mathcal K_{\frak p_0})
 \cap \psi_{\frak p}(s_{\frak p}^{-1}(0) \cap \mathcal K_{\frak p})
\cap
\mathcal S_{\frak d}(X,Z;\widehat{\mathcal U})
$$
such that
\begin{equation}\label{rtachicoversuru}
\bigcup_{j=1}^{J(\frak p)}
\psi_{r^{\frak p}_j}(U^3_{r^{\frak p}_j}\cap s^{-1}_{r^{\frak p}_j}(0))
\supset
 \psi_{\frak p_0}(s_{\frak p_0}^{-1}(0) \cap \mathcal K_{\frak p_0}) \cap \psi_{\frak p}(s_{\frak p}^{-1}(0) \cap \mathcal K_{\frak p})
\cap
\mathcal S_{\frak d}(X,Z;\widehat{\mathcal U}).
\end{equation}
We define
\begin{equation}\label{118888}
U^1_{\frak p \frak p_0}
=
\bigcup_{i=1}^{N(\frak p)} \bigcup_{j=1}^{J(\frak p)}
\left(
\varphi^0_{\frak p_0 r_j^{\frak p}}(U^3_{r_j^{\frak p}})
\cap
\varphi^0_{\frak p_0 q^{\frak p}_i}(U^2_{q^{\frak p}_i}))
\right)
\subset
U_{\frak p_0}.
\end{equation}
We put
$$
\aligned
U^0_{q^{\frak p}_i r_j^{\frak p}}
&=
(\varphi^0_{\frak p_0 r_j^{\frak p}})^{-1}
\left(
\varphi^0_{\frak p_0 r_j^{\frak p}}(U^3_{r_j^{\frak p}})
\cap
\varphi^0_{\frak p_0 q^{\frak p}_i}(U^2_{q_i^{\frak p}}))
\right)
\subset U^3_{r_j^{\frak p}},
\\
U^0_{r_j^{\frak p}q^{\frak p}_i}
&=
(\varphi^0_{\frak p_0 q^{\frak p}_i})^{-1}
\left(
\varphi^0_{\frak p_0 r_j^{\frak p}}(U^3_{r_j^{\frak p}})
\cap
\varphi^0_{\frak p_0 q^{\frak p}_i}(U^2_{q_i^{\frak p}}))
\right)
\subset U^2_{q_i^{\frak p}}.
\endaligned
$$
\begin{lem}\label{lem115555}
$$
\psi_{\frak p_0}(s_{\frak p_0}^{-1}(0) \cap U^1_{\frak p \frak p_0} )
\supset
\psi_{\frak p}(\mathcal K_{\frak p} \cap s^{-1}_{\frak p}(0))
\cap \psi_{\frak p_0}(\mathcal K_{\frak p_0} \cap s_{\frak p_0}^{-1}(0))
\cap
\mathcal S_{\frak d}(X,Z;\widehat{\mathcal U}).
$$
\end{lem}
\begin{proof}
This is a consequence of
(\ref{qcoversuru}) and (\ref{rtachicoversuru}).
\end{proof}
\begin{lem}\label{lemma116}
There exists an embedding
$\Phi^+_{\frak p \frak p_0} : \mathcal U\vert_{U^1_{\frak p \frak p_0}}
\to \mathcal U_{\frak p}$ such that:
\begin{equation}\label{formform1177}
\Phi^+_{\frak p \frak p_0}
\circ
\Phi^0_{\frak p_0 q^{\frak p}_i}\vert_{U^2_{q^{\frak p}_i}
\cap (\varphi^1_{\frak p q_i^{\frak p}})^{-1}(U^1_{\frak p\frak p_0})}
=
\Phi^1_{\frak p q^{\frak p}_i}\vert_{U^2_{q^{\frak p}_i}
\cap (\varphi^1_{\frak p q_i^{\frak p}})^{-1}(U^1_{\frak p\frak p_0})}.
\end{equation}
\end{lem}
$$
\xymatrix{
\mathcal U^1_{\frak p_0} \ar@{.>}[rrrr]  &&&& \mathcal U^1_{\frak p}\\
&
\mathcal U^2_{q_i^{\frak p}}\ar[lu]
\ar[rrru]
&&
\mathcal U^2_{q_{i'}^{\frak p}} \ar[ru]\ar[lllu]
\\
&&
\mathcal U^3_{r_j^{\frak p}} \ar[lu] \ar[ru]
}
$$
\begin{proof}
Recalling the definition of $U^1_{\frak p\frak p_0}$ in
(\ref{118888}), we define a map
$\varphi^+_{\frak p \frak p_0} : U^1_{\frak p\frak p_0} \to
U_{\frak p}$ by
\begin{equation}\label{labelfixedqidefcoochange}
\varphi^+_{\frak p \frak p_0}(x)
=
\varphi^1_{\frak p q^{\frak p}_i}(y_i)
\end{equation}
for $x = \varphi^0_{\frak p_0 q^{\frak p}_i}(y_i) \in \varphi^0_{\frak p_0 q^{\frak p}_i}(U^2_{q^{\frak p}_i}) \cap U^1_{\frak p \frak p_0} $.
We will prove that (\ref{labelfixedqidefcoochange}) is well-defined below.
\par
Suppose $x = \varphi^0_{\frak p_0 q^{\frak p}_i}(y_i) = \varphi^0_{\frak p_0 q^{\frak p}_{i'}}(y_{i'})$.
(Here $y_i \in U^2_{q^{\frak p}_i}$, $y_{i'} \in U^2_{q^{\frak p}_{i'}}$.)
There exist $r^{\frak p}_j$ and $z_j \in U^3_{r_j^{\frak p}}$ such that
$x =  \varphi^0_{\frak p_0 r_j^{\frak p}}(z_j)$.
(This is a consequence of (\ref{118888}).)
\par
By Lemma \ref{lemma114} (2),
$z_j \in U_{q^{\frak p}_i r^{\frak p}_j} \cap U_{q^{\frak p}_{i'} r^{\frak p}_j}$
and
$\varphi_{q^{\frak p}_i r^{\frak p}_j}(z_j) \in U^0_{\frak p_0 q^{\frak p}_i}$.
We have
$$
\varphi^0_{\frak p_0 q^{\frak p}_i}(\varphi_{q^{\frak p}_i r^{\frak p}_j}(z_j))
= \varphi^0_{\frak p_0 r_j^{\frak p}}(z_j) = x = \varphi^0_{\frak p_0 q^{\frak p}_i}(y_i).
$$
Here the first equality is a consequence of the fact that
$\widehat{\Phi^0}$ is a strict KG-embedding.
\par
Since $\varphi^0_{\frak p_0 q^{\frak p}_i}$ is injective, we have $y_i = \varphi_{q^{\frak p}_i r^{\frak p}_j}(z_j)$.
In the same way we have
$y_{i'} = \varphi_{q^{\frak p}_{i'} r^{\frak p}_j}(z_j)$.
Therefore
$$
\varphi^1_{\frak p q^{\frak p}_i}(y_i)
=
\varphi^1_{\frak p q^{\frak p}_i}(\varphi_{q^{\frak p}_i r^{\frak p}_j}(z_j))
=
\varphi^1_{\frak p  r^{\frak p}_j}(z_j).
$$
This is a consequence of the fact that
$\widehat{\Phi^1}$ is a strict KG-embedding.
\par
We also have
$$
\varphi^1_{\frak p q^{\frak p}_{i'}}(y_{i'})
=
\varphi^1_{\frak p q^{\frak p}_{i'}}(\varphi_{q^{\frak p}_{i'} r^{\frak p}_j}(z_j))
=
\varphi^1_{\frak p  r^{\frak p}_j}(z_j).
$$
Thus $\varphi^1_{\frak p q^{\frak p}_i}(y_i) = \varphi^1_{\frak p  q^{\frak p}_{i'}}(y_{i'})$
as required.
\par
Thus we have defined a set theoretical map
$\varphi^+_{\frak p \frak p_0} : U^1_{\frak p \frak p_0} \to U_{\frak p}$.
\par
We note that it is locally written as $\varphi^1_{\frak p q^{\frak p}_i} \circ (\varphi_{\frak p_0 q^{\frak p}_i}^0)^{-1}$.
Furthermore $\varphi^0_{\frak p_0 q^{\frak p}_i}$ is an embedding between orbifolds of the
same dimension. Therefore $\varphi^0_{\frak p_0 q^{\frak p}_i}$ is locally a diffeomorphism.
Hence $\varphi_{\frak p \frak p_0}^+$ is locally a composition of
embedding and diffeomorphism and so is an embedding (of orbifolds).
We can define the bundle map  $\widehat{\varphi^+_{\frak p \frak p_0}}$
in the same way.
\par
We note that the condition for
$(\varphi^+_{\frak p \frak p_0},\widehat{\varphi^+_{\frak p \frak p_0}})$
to be an embedding of Kuranishi charts can be checked locally.
Namely it suffices to check it at a neighborhood of each point.
On the other hand, $\Phi^0_{* q^{\frak p}_i}$ is locally an isomorphism and
$(\varphi^+_{\frak p \frak p_0},\widehat{\varphi^+_{\frak p \frak p_0}})$ is locally
$\Phi_{\frak p q^{\frak p}_i} \circ (\Phi^0_{* q^{\frak p}_i})^{-1}$.
Hence $(\varphi^+_{\frak p \frak p_0},\widehat{\varphi^+_{\frak p \frak p_0}})$
is an embedding of Kuranishi charts, as required.
Then (\ref{formform1177}) is immediate from definition.
\end{proof}
We take a support system $\mathcal K^-$ of
$(X,Z^+_1;\widetriangle{\mathcal U})$
such that $(\mathcal K^-,\mathcal K)$ is a support pair of
$(X,Z^+_1;\widetriangle{\mathcal U})$.
In particular, we have
\begin{equation}
\bigcup_{\frak p \in \frak P} \psi_{\frak p}({\rm Int}~ \mathcal K^-_{\frak p}
\cap
s_{\frak p}^{-1}(0))
\supset
Z^+_1.
\end{equation}
We also take a compact subset $\mathcal K_{\frak p_0}^-$  of $U_{\frak p_0}$
such that
$$
{\mathcal K_{\frak p_0}^-}
\subset
{\rm Int}~ \mathcal K_{\frak p_0},
\qquad
Z^+_0 \subset \psi_{\frak p_0}
(s_{\frak p_0}^{-1}(0) \cap {\rm Int}~ \mathcal K_{\frak p_0}^-).
$$
Put $U'_{\frak p} = {\rm Int}~ \mathcal K^-_{\frak p}$
and $U'_{\frak p_0} = {\rm Int}~ \mathcal K^-_{\frak p_0}$.
\begin{lem}\label{lem119}
There exist open subsets $U''_{\frak p} \subseteq U'_{\frak p}$,
$U''_{\frak p_0} \subseteq U'_{\frak p_0}$ with the
following properties:
\begin{enumerate}
\item
We put
$U''_{\frak p \frak p_0}
=
U^1_{\frak p\frak p_0}
\cap U''_{\frak p_0} \cap
(\varphi_{\frak p\frak p_0}^+)^{-1}(U''_{\frak p})$.
Then
$$
\psi_{\frak p}(U''_{\frak p} \cap s_{\frak p}^{-1}(0))
\cap \psi_{\frak p_0}(U''_{\frak p_0} \cap s_{\frak p_0}^{-1}(0))
=
\psi_{\frak p_0}(U''_{\frak p\frak p_0} \cap s_{\frak p_0}^{-1}(0)).
$$
\item
$
\bigcup_{\frak p} \psi_{\frak p}(U''_{\frak p} \cap s_{\frak p}^{-1}(0))
\supseteq Z^+_1.
$
\item
$
\psi_{\frak p_0}(U''_{\frak p_0} \cap s_{\frak p_0}^{-1}(0))
\supseteq Z^+_0.
$
\end{enumerate}
\end{lem}
This lemma implies that $\Phi_{\frak p \frak p_0}^+\vert_{U''_{\frak p\frak p_0}}$
is a coordinate change from $\mathcal U_{\frak p_0}\vert_{U''_{\frak p_0}}$
to $\mathcal U_{\frak p}\vert_{U''_{\frak p}}$ in the strong sense.
In particular, Definition \ref{coordinatechangedef} (3) is satisfied.
\begin{proof}
We take compact subsets $Z_{\frak p} \subset X$
such that
\begin{enumerate}
\item[(i)]
$
Z_{\frak p}\cap Z^+_{0}
\subseteq
\psi_{\frak p_0}(U^1_{\frak p\frak p_0} \cap s_{\frak p_0}^{-1}(0)).
$
\item[(ii)]
$\psi_{\frak p}(U'_{\frak p}\cap s_{\frak p}^{-1}(0)) \supset Z_{\frak p}$.
\item[(iii)]
$\bigcup Z_{\frak p} \supset Z^+_1$.
\end{enumerate}
Existence of such $Z_{\frak p}$ is a consequence of
Lemma \ref{lem115555}.
\par
We next take open subsets $W_{\frak p}, W_{\frak p_0} \subset X$
such that
\begin{enumerate}
\item[(a)]
$
W_{\frak p}\cap W_{\frak p_0}
\subseteq
\psi_{\frak p_0}(U^1_{\frak p\frak p_0} \cap s_{\frak p_0}^{-1}(0)).
$
\item[(b)]
$\psi_{\frak p}(U'_{\frak p}\cap s_{\frak p}^{-1}(0)) \supset W_{\frak p} \supset Z_{\frak p}$.
\item[(c)]
$\psi_{\frak p_0}(U'_{\frak p_0}\cap s_{\frak p_0}^{-1}(0)) \supset W_{\frak p_0} \supset Z^+_{0}$.
\end{enumerate}
Existence of such $W_{\frak p}, W_{\frak p_0}$ is a consequence of
(i)(ii)(iii) above.
\par
We then take open subsets $U''_{\frak p} \subset U'_{\frak p}$,
$U''_{\frak p_0} \subset U'_{\frak p_0}$
such that
\begin{enumerate}
\item[(A)]
$W_{\frak p}
\supset \psi_{\frak p}(\overline{U''_{\frak p}}\cap s_{\frak p}^{-1}(0))
\supset \psi_{\frak p}({U''_{\frak p}}\cap s_{\frak p}^{-1}(0))
\supset Z_{\frak p}$.
\item[(B)]
$W_{\frak p_0}
\supset \psi_{\frak p_0}(\overline{U''_{\frak p_0}}\cap s_{\frak p_0}^{-1}(0))
\supset
\psi_{\frak p_0}({U''_{\frak p_0}}\cap s_{\frak p_0}^{-1}(0))
\supset Z^+_{0}$.
\item[(C)]
$\overline{U''_{\frak p}} \subset U'_{\frak p}$,
$\overline{U''_{\frak p_0}} \subset U'_{\frak p_0}$
are compact.
\end{enumerate}
Existence of such $U''_{\frak p}$, $U''_{\frak p_0}$
is a consequence of (b)(c).
\par
Lemma \ref{lem119} (3) is nothing but the last inclusion of (B).
Lemma \ref{lem119} (2) follows from (iii)(A).
\par
We finally prove Lemma \ref{lem119} (1).
\begin{sublem}\label{sublem11111}
$$
\psi_{\frak p}(U''_{\frak p} \cap s_{\frak p}^{-1}(0))
\cap
\psi_{\frak p_0}(U''_{\frak p_0} \cap s_{\frak p_0}^{-1}(0))
\subseteq
\psi_{\frak p_0}
(U^1_{\frak p\frak p_0} \cap s_{\frak p_0}^{-1}(0)).
$$
\end{sublem}
\begin{proof}
This is a consequence of (a)(A)(B).
\end{proof}
Note
$U''_{\frak p \frak p_0}
=
U^1_{\frak p\frak p_0}
\cap U''_{\frak p_0} \cap
(\varphi^+_{\frak p\frak p_0})^{-1}(U''_{\frak p})$
by definition.
\begin{sublem}
$$
\psi_{\frak p}(U''_{\frak p} \cap s_{\frak p}^{-1}(0))
\cap
\psi_{\frak p_0}(U''_{\frak p_0} \cap s_{\frak p_0}^{-1}(0))
\subseteq
\psi_{\frak p_0}
(U''_{\frak p\frak p_0} \cap s_{\frak p_0}^{-1}(0)).
$$
\end{sublem}
\begin{proof}
Let $\psi_{\frak p}(x) = \psi_{\frak p_0}(y)$ be in the left hand side.
Since $\psi_{\frak p_0}$ is injective, Sublemma \ref{sublem11111}
implies
$y \in U^1_{\frak p\frak p_0}  \cap U''_{\frak p_0} \cap s_{\frak p_0}^{-1}(0)$.
Since $\psi_{\frak p}(\varphi_{\frak p\frak p_0}^+(y)) =
\psi_{\frak p}(x)$, the injectivity of $\psi_{\frak p}$ implies
$\varphi_{\frak p\frak p_0}^+(y) \in U''_{\frak p}$.
Therefore $y \in U''_{\frak p \frak p_0}$ as required.
\end{proof}
The opposite inclusion
$$
\psi_{\frak p}(U''_{\frak p} \cap s_{\frak p}^{-1}(0))
\cap
\psi_{\frak p_0}(U''_{\frak p_0} \cap s_{\frak p_0}^{-1}(0))
\supseteq
\psi_{\frak p_0}
(U''_{\frak p\frak p_0} \cap s_{\frak p_0}^{-1}(0))
$$
is obvious.
We have thus proved Lemma \ref{lem119} (1).
The proof of Lemma  \ref{lem119} is complete.
\end{proof}
Thus $U''_{\frak p}$, $U''_{\frak p_0}$ together with
the restriction of
$\Phi^1_{\frak p \frak q}$, $\Phi_{\frak p \frak p_0}^+\vert_{U''_{\frak p
\frak p_0}}$ satisfy the conditions of
Definition \ref{gcsystem} except (7) and (8).
Then we use  Shrinking Lemma \cite[Theorem 2.7]{foooshrink}
to shrink $U''_{\frak p}$, $U''_{\frak p_0}$ and $U''_{\frak p \frak p_0}$ so that
Definition \ref{gcsystem} (7) and (8)
also hold.
\par
We now change the notations $U''_{\frak p}$, $U''_{\frak p_0}$
to $U_{\frak p}$, $U_{\frak p_0}$
so that they denote the Kuranishi neighborhoods
of the good
coordinate system we have obtained above.
\par\medskip
To complete the proof of Proposition \ref{inductiveprop}
it remains to prove Proposition \ref{inductiveprop} (5).
For this purpose, we will shrink $U_{\frak p}$, $U_{\frak p_0}$ further
as follows.
\par
We first take open subsets $U^{(1)}_{\frak p} \subset U_{\frak p}$
and $U^{(1)}_{\frak p_0} \subset U_{\frak p_0}$
with the following properties:
\begin{enumerate}
\item
$\bigcup_{\frak p \in \frak P} \psi_{\frak p}(U^{(1)}_{\frak p}
\cap s_{\frak p}^{-1}(0)) \supset Z_1$.
\item
$\bigcup_{\frak p \in \frak P} \psi_{\frak p}(U^{(1)}_{\frak p}
\cap s_{\frak p}^{-1}(0))
\subset {\rm Int}\, Z^+_1$.
\item
$\psi_{\frak p_0}(U^{(1)}_{\frak p_0}
\cap s_{\frak p_0}^{-1}(0)) \supset Z_0$.
\item
$\psi_{\frak p_0}(U^{(1)}_{\frak p_0}
\cap s_{\frak p_0}^{-1}(0))
\subset {\rm Int}\, Z^+_0$.
\end{enumerate}
By Lemma \ref{lem320} we have a good coordinate system of
$Z_+^{++} \subset X$ whose Kuranishi charts are
$\mathcal U_{\frak p}\vert_{U^{(1)}_{\frak p}}$
and $\mathcal U_{\frak p_0}\vert_{U^{(1)}_{\frak p_0}}$.
(Here $Z_+^{++}$ is a compact neighborhood of $Z_+$.
The compact neighborhood $Z^+_+$ we will obtain is a subset of
$Z_+^{++}$.)
We denote this good coordinate system by
$\widetriangle{\mathcal U^{(1)}}$.
\par
By the assumptions put in Situation \ref{situation101} the following holds
after replacing $\widehat{\mathcal U}$ by its open
substructure if necessary.
\begin{enumerate}
\item[(A)]
If $p \in \psi_{\frak p_0}(U_{\frak p_0}^{(1)} \cap s_{\frak p_0}^{-1}(0)) \cap Z_+^{++}$,
then there exists an embedding
of Kuranishi charts $
\Phi^{(1)}_{\frak p_0 p} : \mathcal U_p \to \mathcal U_{\frak p_0}^{(1)}$.
Here $\Phi^{(1)}_{\frak p_0 p}$ is an open restriction
of $\Phi^{0}_{\frak p_0 p}$
\item[(B)]
If $p \in\psi_{\frak p}(U_{\frak p}^{(1)} \cap s_{\frak p}^{-1}(0)) \cap Z_+^{++}$and
$\frak p \in \frak P$, then there exists an embedding
of Kuranishi charts $\Phi^{(1)}_{\frak p p} :\mathcal U_p \to \mathcal U_{\frak p}$.
Here $\Phi^{(1)}_{\frak p p}$ is an open restriction of $\Phi^{1}_{\frak p p}$
\item[(C)]
If $p \in {\rm Im}(\psi_{\frak p}) \cap  Z_+^{++}$,
$q \in {\rm Im}(\psi_{\frak q}) \cap Z_+^{++}$,
$q \in  {\rm Im}(\psi_{ p})$
and $\frak q \le \frak p$ then
the following diagram commutes.
\begin{equation}\label{diag33repeat}
\begin{CD}
\mathcal U_{q}\vert_{U_{pq}
\cap \varphi_{\frak q q}^{-1}(U^{(1)}_{\frak p \frak q})} @ > {\Phi^{(1)}_{\frak q q}} >>
{\mathcal U}^{(1)}_{\frak q}\vert_{{U}^{(1)}_{\frak p\frak q}}  \\
@ V{\Phi_{pq}}VV @ VV{\Phi^{(1)}_{\frak p \frak q}}V\\
\mathcal U_{p} @ > {\Phi^{(1)}_{\frak p p}} >>{\mathcal U}_{\frak p}
\end{CD}
\end{equation}
\end{enumerate}
Note we use
Proposition \ref{inductiveprop} (2)(4),
when we take $U^{(1)}_{\frak p_0}$, $U^{(1)}_{\frak p}$
to show the properties (A)(B) above.
\par
We also note that, to define a KG-embedding $\widehat{\mathcal U} \to {\widetriangle{\mathcal U^+}}$, we will use $\Phi^{(1)}_{\frak p p}$, $\Phi^{(1)}_{\frak p}$ or their open
restrictions. Therefore the only remaining nontrivial part of the
proof of Proposition \ref{inductiveprop} (5) is the commutativity of the
next diagram
\begin{equation}\label{diag33new}
\begin{CD}
\mathcal U^0_{q}\vert_{U^0_{pq}
\cap \varphi_{\frak p_0 q}^{-1}(U_{\frak p \frak p_0})} @ > {\Phi_{\frak p_0 q}} >>
{\mathcal U}^+_{\frak p_0}\vert_{{U}_{\frak p\frak p_0}
}  \\
@ V{\Phi_{pq}\vert_{U^0_{pq}\cap \varphi_{\frak p_0 q}^{-1}(U_{\frak p \frak p_0})}}VV @ VV{\Phi^+_{\frak p \frak p_0}}V\\
\mathcal U^0_{p} @ > {\Phi_{\frak p p}} >>{\mathcal U}^+_{\frak p}
\end{CD}
\end{equation}
which we need to prove under the assumption
$p \in {\rm Im}(\psi_{\frak p}) \cap Z_+^+$,
$q \in {\rm Im}(\psi_{\frak p_0}) \cap Z_+^+$,
$q \in  {\rm Im}(\psi_{p})$.
Here
$\widehat{\mathcal U^0}$ is an  open substructure of
 $\widehat{\mathcal U}$ and $\mathcal U^0_{p}$ is its Kuranishi chart.
$\Phi_{\frak p p}$
(resp. $\Phi_{\frak p_0 q}$) is an open restriction of  $\Phi^{(1)}_{\frak p p}$
(resp. $\Phi^{(1)}_{\frak p_0 q}$). ${\mathcal U}^+_{\frak p}$
(resp. ${\mathcal U}^+_{\frak p_0}$) is an
open subchart of $\mathcal U^{(1)}_{\frak p}$ (resp. $\mathcal U^{(1)}_{\frak p_0}$)
and $\Phi^+_{\frak p \frak p_0}$
is an open restriction of $\Phi^{(1)}_{\frak p \frak p_0}$.
\par
We will use the next lemma to show the commutativity of
Diagram (\ref{diag33new}).
We take a support system $\{\mathcal K_{\frak p}^+
\mid \frak p \in \frak P\} \cup \{\mathcal K^+_{\frak p_0}\}$ of the good coordinate system on $Z_+$ that we have already obtained.
Note that
$\mathcal K^+_{\frak p_0} \subset \overset{\circ}{\mathcal K_{\frak p_0}}$,
$\mathcal K^+_{\frak p} \subset \overset{\circ}{\mathcal K_{\frak p}}$.
\begin{lem}\label{119lem}
There exist open subsets $U_{\frak p}^+ \subset U_{\frak p}$,
$U_{\frak p_0}^+ \subset U_{\frak p_0}$ with the following properties.
\begin{enumerate}
\item
$
\psi_{\frak p_0}(\mathcal K_{\frak p_0}^+ \cap
s_{\frak p_0}^{-1}(0) \cap U_{\frak p_0}^+)
\cup
\bigcup_{\frak p} \psi_{\frak p}(\mathcal K^+_{\frak p} \cap
s_{\frak p}^{-1}(0) \cap U_{\frak p}^+)
\supset
Z_+.
$
\item
If
$
p \in \psi_{\frak p_0}(\mathcal K_{\frak p_0}^+ \cap
s_{\frak p_0}^{-1}(0) \cap U_{\frak p_0}^+)
\cap \psi_{\frak p}(\mathcal K^+_{\frak p} \cap
s_{\frak p}^{-1}(0) \cap U_{\frak p}^+),
$
then there exist $q^{\frak p}_i$ ($i=1,\dots,N(\frak p)$) such that
$$
p  \in \psi_{q^{\frak p}_i}(s_{q^{\frak p}_i}^{-1}(0) \cap U_{q^{\frak p}_i}^2).
$$
\end{enumerate}
\end{lem}
Recall $Z_+ \subset \mathcal S_{\frak d}(X,Z;\widehat{\mathcal U})$.
(See (\ref{form111}).)
We also note that we do {\it  not} assume
$p \in \mathcal S_{\frak d}(X,Z;\widehat{\mathcal U})$
in Item (2) above.
\begin{proof}
We choose a decreasing sequence of relatively compact open subsets
$U_{\frak p_0}^{[a]} \subset U$, $U_{\frak p}^{[a]}\subset U_{\frak p}$,
($a=1,2,\dots$) such that
$$
\overline{U_{\frak p}^{[a+1]}} \subset U_{\frak p}^{[a]},
\quad
\overline{U_{\frak p_0}^{[a+1]}} \subset U_{\frak p_0}^{[a]}
$$
and
\begin{equation}\label{new1118}
\bigcap_{a=1}^{\infty} U_{\frak p_0}^{[a]}
= \psi_{\frak p_0}^{-1}(Z_0) \cap \mathcal K_{\frak p_0}^+,
\quad
\bigcap_{a=1}^{\infty} U_{\frak p}^{[a]}
= \psi_{\frak p}^{-1}(Z_1) \cap \mathcal K^+_{\frak p}.
\end{equation}
Then, for each $a$, the choice $U_{\frak p_0}^+ = U_{\frak p_0}^{[a]}$
and $U_{\frak p}^+ = U_{\frak p}^{[a]}$ satisfies (1) above.
\par
On the other hand, we can use (\ref{qcoversuru})
to prove that (2) is satisfied for sufficiently large $a$ as follows.
Suppose that (2) is not satisfied for $a_n \to \infty$.
Then we have
\begin{equation}\label{phakokoniiru}
p_n \in \psi_{\frak p_0}(\mathcal K_{\frak p_0}^+ \cap
s^{-1}(0) \cap \overline{U_{\frak p_0}^{[a_n]}})
\cap \psi_{\frak p}(\mathcal K^+_{\frak p} \cap
s_{\frak p}^{-1}(0) \cap \overline{U_{\frak p}^{[a_n]}})
\end{equation}
 and
\begin{equation}\label{phaittenai}
p_n  \notin
\bigcup_{i=1}^{N(\frak p)}
\psi_{q^{\frak p}_i}(s_{q^{\frak p}_i}^{-1}(0) \cap U_{q^{\frak p}_i}^{2}).
\end{equation}
Note that we may assume $\frak p$ is independent of $n$
by taking a subsequence if necessary,
since  $\frak p \in \frak P$ and $\frak P$ is a finite set.
\par
Since the right hand side of (\ref{phakokoniiru}) is contained in
$\psi_{\frak p_0}(\mathcal K_{\frak p_0}^+ \cap
s^{-1}(0) \cap \overline{U_{\frak p_0}^{[a_1]}})
\cap \psi_{\frak p}(\mathcal K^+_{\frak p} \cap
s_{\frak p}^{-1}(0) \cap \overline{U_{\frak p}^{[a_1]}})$
which is
compact and independent of $n$,
we may take a subsequence and may assume that
$p_n$ converges.
Then (\ref{new1118}) and (\ref{phakokoniiru}) imply
$$
\lim_{n\to \infty} p_n
\in Z_0 \cap \psi_{\frak p_0}(\mathcal K_{\frak p_0}^+
\cap s_{\frak p_0}^{-1}(0))
\cap
Z_1 \cap \psi_{\frak p}(\mathcal K^+_{\frak p}
\cap s_{\frak p}^{-1}(0)).
$$
This contradicts to (\ref{phaittenai}) since
the right hand side of  (\ref{phaittenai})
is a neighborhood of $Z_0 \cap \psi_{\frak p_0}(\mathcal K_{\frak p_0}^+
\cap s_{\frak p_0}^{-1}(0))
\cap
Z_1 \cap \psi_{\frak p}(\mathcal K^+_{\frak p}
\cap s_{\frak p}^{-1}(0))$.
\end{proof}
We consider ${\mathcal U}^+_{\frak p} =
{\mathcal U}_{\frak p}\vert_{U_{\frak p}^+}$,
${\mathcal U}^+_{\frak p_0} =
{\mathcal U}\vert_{U_{\frak p_0}^+}$.
Put
$$
\frak s_{\frak p}^+ = \frak s_{\frak p}\vert_{{U}^+_{\frak p}},
\quad
\frak s_{\frak p_0}^+ = \frak s_{\frak p_0}\vert_{{U}^+_{\frak p_0}}
$$
and denote by $\psi^+_{\frak p}$ and $\psi^+_{\frak p_0}$
the restrictions of the parametrization
$\psi_{\frak p}$ and $\psi_{\frak p_0}$
to $(\frak s_{\frak p}^+)^{-1}(0)$ and
$(\frak s_{\frak p_0}^+)^{-1}(0)$,
respectively.
\par
Next we note the following lemmata.

\begin{lem}\label{lem111441}
There exists an open neighborhood
$U^0_p$ of $o_p$ in $U_p$ for each $p$ with the following properties.
\begin{enumerate}
\item If
$p  \in \psi_{q^{\frak p}_i}(s_{q^{\frak p}_i}^{-1}(0) \cap U_{q^{\frak p}_i}^2)$, then
\begin{equation}\label{eq1115}
U^0_p \subset U_{q^{\frak p}_i p}
\cap
\varphi_{q^{\frak p}_i p}^{-1}( U^2_{q^{\frak p}_i}).
\end{equation}
\item
If $p \in {\rm Im}(\psi_{q_i^{\frak p}})$ and $q \in \psi_{p}(U_{p}^0 \cap s_p^{-1}(0))$,
then
 $q \in {\rm Im}(\psi_{q_i^{\frak p}})$.
\end{enumerate}
\end{lem}
\begin{proof}
We remark that, for each $p$, there exists only a finite number of $q^{\frak p}_i$
satisfying the assumptions of (1) or (2).
The lemma immediately follows from this remark.
\end{proof}
\begin{lem}\label{111515}
For each $p,q$ with $q \in \psi_{p}(U_{p}^0 \cap s_p^{-1}(0))$
there exists an open neighborhood $U_{pq}^0$ of $o_q$ in $U_{pq} \cap U^0_q$ with
the following properties.
\begin{enumerate}
\item
$\varphi_{pq}(U_{pq}^0) \subset U_p^0$.
\item
If $p \in {\rm Im}(\psi_{q_i^{\frak p}})$ then
$\varphi_{pq}(U_{pq}^0) \subset U_{q_i^{\frak p}p}$.
\end{enumerate}
\end{lem}
\begin{proof}
Since $\varphi_{pq}(o_q) \in U^0_p$ for
$q \in \psi_{p}(U_{p}^0 \cap s_p^{-1}(0))$,
(1) holds for a sufficiently small neighborhood $U^0_{pq}$ of $o_q$.
Lemma \ref{lem111441}(2)  implies
(2) for a sufficiently small neighborhood $U^0_{pq}$ of $o_q$.
\end{proof}
The pair $\mathcal U_p\vert_{U_p^0}$, $\Phi_{pq}\vert_{U_{pq}^0}$ define
an open substructure of $\widehat{\mathcal U}$
by Lemma \ref{111515} (1).
\par
Now we will prove the commutativity of Diagram (\ref{diag33new}).
We take $Z_+^+$ so that it is contained in the
left hand side of Lemma \ref{119lem} (1).
Suppose $p \in {\rm Im}(\psi^+_{\frak p}) \cap Z_+^+$,
$q \in {\rm Im}(\psi^+_{\frak p_0}) \cap Z_+^+$,
$q \in  \psi_{p}(U_p^0\cap s_p^{-1}(0))$.
Let
$$
y \in U_{pq}^0 \cap \varphi_{\frak p_0 q}^{-1}(U_{\frak p\frak p_0}).
$$
There exists $i$ such that
$
p  \in \psi_{q^{\frak p}_i}(s_{q^{\frak p}_i}^{-1}(0) \cap U_{q^{\frak p}_i}^2)
$
by Lemma \ref{119lem} (2).
Then $
q  \in \psi_{q^{\frak p}_i}(s_{q^{\frak p}_i}^{-1}(0) \cap U_{q^{\frak p}_i}^2)
$
by Lemma \ref{lem111441} (2).
\par
Now we have
$$
\aligned
(\varphi^+_{\frak p \frak p_0}\circ \varphi^0_{\frak p_0 q})(y)
&=
\varphi^+_{\frak p \frak p_0}(\varphi^0_{\frak p_0 q^{\frak p}_i}(\varphi_{q^{\frak p}_i q}(y)))
\\
&= \varphi^1_{\frak p  q^{\frak p}_i}(\varphi_{q^{\frak p}_i q}(y))
\\
&= \varphi^1_{\frak p  q^{\frak p}_i}(\varphi_{q^{\frak p}_i p}(\varphi_{p q}(y)))
\\
&=
\varphi^1_{\frak p   p}(\varphi_{p q}(y)).
\endaligned$$
Here the equality in the first line follows from (\ref{eq1115})
(with $p$ replaced by $q$) and the fact that $\widehat{\Phi^0}$
is a KG-embedding.
The equality of the second line is the definition of $\varphi^+_{\frak p \frak p_0}$, that
is, (\ref{labelfixedqidefcoochange}).
The equality of the third line is the consequence of
the compatibility of coordinate change of Kuranishi
structure and Lemma \ref{111515} (2).
The equality of the fourth line is the consequence of the
fact that $\widehat{\Phi^1}$ is a strict KG-embedding.
\par
We have proved the commutativity of Diagram (\ref{diag33new})
for the maps between base orbifolds.
The proof of the commutativity of  the maps between bundles
is the same.
The proof of Proposition \ref{inductiveprop}
is now complete.
\end{proof}
We use Proposition \ref{inductiveprop} to prove Theorem
\ref{Them71restate} as follows.
\par
We will construct a good coordinate system of
a closed subset
$\mathcal S_{\frak d}(X,Z;\widehat{\mathcal U})$ of $Z$
by a downward induction on $\frak d$.
\par
Let ${\widetriangle{\mathcal U_{\frak d+1}}}
= (\frak P,\{\mathcal U_{\frak p}\},\{\Phi_{\frak p\frak q}\})$
be a good coordinate system of a compact neighborhood of
$\mathcal S_{\frak d+1}(X,Z;\widehat{\mathcal U})$.
We will construct one of
$\mathcal S_{\frak d}(X,Z;\widehat{\mathcal U})$.
\par
We put
$$
\frak P(\frak d) =
\{
\frak p \in \frak P \mid \dim U_{\frak p} > \frak d
\}.
$$
We take a support system $\mathcal K_{\frak p}$ of
$\widetriangle{\mathcal U_{\frak d+1}}$
and put
$$
B
=
\mathcal S_{\frak d}(X,Z;\widehat{\mathcal U})
\setminus
\bigcup_{\frak p \in \frak P(\frak d)}
\psi_{\frak p}(s_{\frak p}^{-1}(0) \cap {\rm Int}\, \mathcal K_{\frak p}).
$$
Then $B$ is a compact subset of
$
\mathcal S_{\frak d}(X,Z;\widehat{\mathcal U})
\setminus
\bigcup_{\frak d' > \frak d}
\mathcal S_{\frak d'}(X,Z;\widehat{\mathcal U})
$.
We take a finite number of points $p_1,\dots,p_N
\in B$
such that
$$
\bigcup_{i=1}^N \psi_{p_i}(s_{p_i}^{-1}(0))
\supset B.
$$
We take compact subsets $K_{p_i} \subset U_{p_i}$ such that
$$
\bigcup_{i=1}^N \psi_{p_i}(s_{p_i}^{-1}(0) \cap K_{p_i})
\supset B.
$$
We put $Z_i = \psi_{p_i}(s_{p_i}^{-1}(0))$
and
$$
Z(0) =
\bigcup_{\frak p \in \frak P(\frak d)}
\psi_{\frak p}(s_{\frak p}^{-1}(0) \cap \mathcal K_{\frak p}) \cap Z.
$$
\par
Now using  Proposition \ref{inductiveprop}
we can construct a good coordinate system on
a compact neighborhood of
$$
Z(0) \cup \bigcup_{i=1}^n Z_i
$$
by an induction over $n$.
We thus obtain a good coordinate system
on $\mathcal S_{\frak d}(X,Z;\widehat{\mathcal U})$.
\par
Now the downward induction over $\frak d$ is complete and we obtain
a good coordinate system
on $\mathcal S_{0}(X,Z;\widehat{\mathcal U}) = Z$.
The proof of Theorem
\ref{Them71restate} is complete.
\end{proof}
\begin{rem}
\begin{enumerate}
\item
The heart of the proof of Theorem
\ref{Them71restate}  is Lemma \ref{lemma116}.
There (and in several other places) we crucially use the
fact that equality between embeddings of orbifolds
is a local property. (It implies that the equality between
embeddings of Kuranishi charts is also a local property).
Thanks to this fact, we can check  various equalities
by looking finer and finer charts.
For this main idea of the proof to work,
we assumed our orbifold to be effective and
the maps between them to be embeddings.
Without this restriction the argument
will become more cumbersome and lengthy.
(We however believe that one can prove somewhat similar
results without assuming effectivity of orbifold.
Joyce may have done it in \cite{joyce2}.)
\item
The proof of Theorem  \ref{Them71restate} in \cite{foootech}
is basically the same as one we gave above.
However the argument of the proof of Lemma \ref{lemma116}
appeared twice in \cite{foootech} once during the proof of \cite[Lemma 7.5]{foootech}
and once during the proof of  \cite[Lemma 7.24]{foootech}.
We reorganize the proof so that we need to use it only once.
Also the argument
used to prove the main theorem of \cite{foooshrink}
appeared twice in \cite{foootech}.
Besides these simplifications,
we defined the notion of embedding of Kuranishi charts
and use it systematically in this document.
For example the statement of \cite[Lemma 7.24]{foootech}
is nothing but the existence of an
appropriate embedding of Kuranishi charts.
\par
These changes make the presentation of the proof simpler and
make the proof shorter, in
this document, although its mathematical contents are the same as those in \cite{foootech}.
\par
We also prove compatibility of the good coordinate system
to the Kuranishi structure (that is the existence of the
embedding of the latter to the former)
explicitly. (Namely we proved Proposition \ref{inductiveprop} (5),
explicitly.)
In \cite{foootech}, the proof of this part was said similar and was omitted.
The proof we gave in this article is indeed similar to the other part of the proof.
We repeated the argument here for completeness' sake.
\end{enumerate}
\end{rem}
In Subsections \ref{subsec:constgcsrel}-\ref{subsec:moreversionegcs2},
we will prove several variants of Theorem \ref{Them71restate}.

\subsection{Construction of good coordinate systems: the relative case 1}
\label{subsec:constgcsrel}

This section will be occupied by the proofs of Propositions \ref{prop518} and \ref{prop519}.

\begin{proof}[Proof of Propositions \ref{prop518}
and \ref{prop519}]
We will prove Proposition \ref{prop519} in this subsection in detail.
Proposition \ref{prop518} then follows by putting
$\widehat{\mathcal U^+}
= \widehat{\mathcal U_1^+} = \widehat{\mathcal U_2^+}$.
\par
\par
We will use the next proposition in addition to Proposition \ref{inductiveprop}
in our inductive construction of good coordinate system.

\begin{prop}\label{prop1111}
Under Situation \ref{situation101},
we assume in addition that {for each $a=1,2$
there exists a Kuranishi structure
$\widehat{\mathcal U^{+}_a}$ of $Z \subseteq X$} with the following properties.
\begin{enumerate}
\item[(a)]
There exists a strict KK-embedding $\widehat{\Phi^2_a} : \widehat{\mathcal U}
\to \widehat{\mathcal U^{+}_a}$.
\item[(b)]
There exists a GK-embedding $\widehat{\Phi^3_a} : \widetriangle{\mathcal U}
\to \widehat{\mathcal U^{+}_a}\vert_{Z_1^+}$ such that the
composition
$\widehat{\Phi^3_a}\circ\widehat{\Phi^1} :
\widehat{\mathcal U}\vert_{Z_1^+} \to {\widetriangle{\mathcal U}}
\to {\widehat{\mathcal U^{+}_a}}\vert_{Z_1^+}$ is an open restriction of $\widehat{\Phi^2_a}\vert_{Z_1^+}$.
\item[(c)]
There exists a GK-embedding $\widehat{\Phi^4_a} : \widetriangle{{\mathcal U}_{\frak p_0}}
\to {\widehat{\mathcal U^{+}_a}}\vert_{Z_0^+}$ such that
the composition $\widehat{\Phi^4_a}\circ\widehat{\Phi^0} :
\widehat{\mathcal U}\vert_{Z_{0}^+}
\to \widetriangle{\mathcal U_{\frak p_0}} \to {\widehat{\mathcal U^{+}_a}}
\vert_{Z_0^+}$ is an open restriction of $\widehat{\Phi^2_a}\vert_{Z_0^+}$.
\end{enumerate}
Then there exists ${\widetriangle{\mathcal U^+}}$ such
as in the conclusion of Proposition \ref{inductiveprop} that the following holds in addition.
\begin{enumerate}
\item
There exists a GK-embedding
$\widehat{\Phi^5_a} : {\widetriangle{\mathcal U^{+}}}
\to {\widehat{\mathcal U^{+}_a}}$.
\item The composition
$\widehat{\Phi^5_a} \circ \widehat{\Phi^{+}} :
\widehat{\mathcal U}\vert_{Z_+^+} \to {\widetriangle{\mathcal U^+}}
\to {\widehat{\mathcal U^{+}_a}}$ is an open restriction of
$\widehat{\Phi^2_a}\vert_{Z_+^+}$.
\item
The restriction of $\widehat{\Phi^5_a}\vert_{Z_1}$ to a strictly open substructure
of ${\widetriangle{\mathcal U^+}}\vert_{Z_1}$ coincides with $\widehat{\Phi^3_a}\vert_{Z_1}$.
\item
The restriction of  $\widehat{\Phi^5_a}\vert_{Z_0}$ to a strictly open substructure
of ${\widetriangle{\mathcal U^+}}\vert_{Z_0}$ coincides with
 $\widehat{\Phi^4_a}\vert_{Z_0}$.
\end{enumerate}
\end{prop}
\begin{proof}
Let $\widetriangle{\mathcal U^+}$ be the good coordinate
system constructed in the proof of Proposition
\ref{inductiveprop}.
During the proof of Proposition \ref{inductiveprop}
we took a support system $\mathcal K$ of $\widetriangle{\mathcal U}$
and a compact subset $\mathcal K_{\frak p_0}$ of
$U_{\frak p_0}$.
We put $\mathcal K^+ = \mathcal K\cup \{\mathcal K_{\frak p_0}\}
= \{\mathcal K^+_{\frak p} \mid \frak p \in \frak P^+\}$.
Then $\mathcal K^+$ is a support system of  $\widetriangle{\mathcal U^+}$.
We denote by $U^3_{a;\frak p}(q_i^{\frak p})$ and
$U^4_{a;\frak p_0}(q_i^{\frak p})$ the domains of
$\varphi^3_{a;q_i^{\frak p} \frak p}$
and $\varphi^4_{a;q_i^{\frak p}\frak p_0}$,
which are part of data carried by $\widehat{\Phi^3_a}$ and $\widehat{\Phi^4_a}$.
During this proof we choose $U^1_{q_i^{\frak p}}$
so that the following properties are also satisfied:
\begin{equation}\label{formula1125}
\varphi^1_{\frak p q^{\frak p}_i}(U^1_{q^{\frak p}_i})
\subset U^3_{a;\frak p}(q_i^{\frak p}), \quad
\varphi^0_{\frak p_0 q^{\frak p}_i}(U^1_{q^{\frak p}_i})
\subset U^4_{a;\frak p_0}(q_i^{\frak p})
\qquad \text{for $a=1,2$}.
\end{equation}
Here $U^3_{a;\frak p}(q_i^{\frak p})$ and
$U^4_{a;\frak p_0}(q_i^{\frak p})$ are domains of
$\varphi^3_{a;q_i^{\frak p} \frak p}$
and $\varphi^4_{a;q_i^{\frak p}\frak p_0}$,
which are parts of $\widehat{\Phi^3_a}$ and
$\widehat{\Phi^4_a}$.
Note that we required Property \ref{proper113} to define $U^1_{q}$.
We can certainly require
$$\varphi^1_{\frak p q}(U^1_{q}) \subset U^3_{a;\frak p}(q),\quad
\varphi^0_{\frak p_0 q}(U^1_{q}) \subset U^4_{a;\frak p_0}(q)
\qquad \text{for $a=1,2$}
$$
in addition,
by replacing $\widehat{\mathcal U}$ by its open substructure
if necessary.
Then (\ref{formula1125}) is satisfied.
\par
We will further shrink $U^+_{\frak p}$ to $U^{+\prime}_{\frak p}$
so that it satisfies the conclusion of Proposition \ref{prop1111}.
(Here $U_{\frak p}^+$ is a Kuranishi neighborhood
which is a part of $\widetriangle{\mathcal U^+}$.)
\par
Let $p \in \psi^+_{\frak p}(U_{\frak p}^+ \cap (s^+_{\frak p})^{-1}(0))$. We will construct
a neighborhood
$U^5_{\frak p}(p) \subset U_{\frak p}^+$ of $p$ and
an embedding of Kuranishi charts
$\Phi^5_{a;p\frak p} : \mathcal U_{\frak p}^+\vert_{U^5_{\frak p}(p)}
\to \mathcal U^{+}_{a,p}$ for $a=1,2$.
Here $\mathcal U^{+}_{a,p}$ is the Kuranishi neighborhood which
${\widehat{\mathcal U^{+}_a}}$ assigns to $p\in Z$.
\par\smallskip
\noindent(Case 1)
$\frak p \in \frak P$.
\par
By assumption (b) there exists a neighborhood
$U^3_{a;\frak p}(p) \subset U_{\frak p}$ of $o_{\frak p}(p)$
and an embedding
$\Phi^3_{a;p\frak p} : \mathcal U_{\frak p}\vert_{U^3_{a;\frak p}(p)}
\to \mathcal U^{+}_{a,p}$.
We take
$$
U^{5\prime}_{\frak p}(p) =  U^+_{\frak p} \cap  \bigcap_{a=1}^2U^3_{a;\frak p}(p)
$$
and $\Phi^{5\prime}_{a;p\frak p}: \mathcal U_{\frak p}\vert_{U^{5\prime}_{\frak p}(p)}
\to \mathcal U^{+}_{a,p}$ to be the restriction of $\Phi^3_{a;p\frak p}$.
\par\smallskip
\noindent(Case 2)
$\frak p = \frak p_0$.
\par
By assumption (c) there exists a neighborhood $U^4_{a;\frak p_0}(p) \subset U_{\frak p_0}$ of $o_{\frak p_0}(p)$
and an embedding
$\Phi^4_{a;p\frak p_0} : \mathcal U_{\frak p_0}\vert_{U^4_{a;\frak p_0}(p)}
\to \mathcal U^{+}_{a;p}$.
We take
$$
U^{5\prime}_{\frak p_0}(p) = U^+_{\frak p_0}\cap  \bigcap_{a=1}^2U^4_{a;\frak p_0}(p)
$$
and $\Phi^{5\prime}_{a;p\frak p_0}: \mathcal U_{\frak p_0}\vert_{U^{5\prime}_{\frak p_0}(p)}
\to \mathcal U^{+}_{a;p}$ to be the restriction of $\Phi^4_{a;p\frak p_0}$.
\par\smallskip
Most of the properties required for $(\{\Phi^{5\prime}_{a;p\frak p}\},\{U^{5\prime}_{\frak p}(p)\})$
to be a GK-embedding is a direct consequence of Assumptions (a), (b), (c).
The nontrivial part to check is the commutativity of
Diagram (\ref{diagram58}) in the following case:
Here we will further shrink $U^+_{\frak p}$, $U^+_{\frak p_0}$, $U^{5\prime}_{\frak p}(p)$,
$U^{5\prime}_{\frak p_0}(q)$ to
$U^{+-}_{\frak p}$, $U^{+-}_{\frak p_0}$, $U^{5}_{\frak p}(p)$,
$U^{5}_{\frak p_0}(q)$, respectively,
by Lemmas \ref{lem11166} and \ref{1112label} below.
After these shrinking,
we will prove the commutativity of Diagram
(\ref{diagram58sec11}), where
$$
p \in \psi^+_{\frak p}((s^+_{\frak p})^{-1}(0) \cap U^{+ -}_{\frak p}), \quad
q \in \psi^+_{\frak p_0}((s^+_{\frak p_0})^{-1}(0) \cap U^{+ -}_{\frak p_0})
\cap \psi^+_{\frak p}((s^+_{\frak p})^{-1}(0) \cap U^{5}_{\frak p}(p)).
$$
\begin{equation}\label{diagram58sec11}
\begin{CD}
\mathcal U_{\frak p_0}\vert_{
(U^+_{\frak p\frak p_0} \cap (\varphi_{\frak p\frak p_0}^+)^{-1}(U^{5}_{\frak p}(p)))
\cap ( U^{5}_{\frak p_0}(q)
\cap(\varphi_{q\frak p_0}^4)^{-1}({U}_{a;pq}^{+}))} @ > {\Phi_{a;q\frak p_0}^4} >>
{\mathcal U}^+_{a;q}\vert_{{U}_{a;pq}^+}  \\
@ V{\Phi^+_{\frak p\frak p_0}}VV @ VV{\Phi_{a;pq}^{+}}V\\
\mathcal U_{\frak p}\vert_{U^{5}_{\frak p}(p)} @ > {\Phi_{a;p\frak p}^3} >>{\mathcal U}_{a;p}^{+}
\end{CD}
\end{equation}
Here
$\Phi_{a;pq}^{+}$ is the coordinate change of the Kuranishi structure
$\widehat{\mathcal U^+_a}$ whose domain is ${U}_{a;pq}^+$, and
$\Phi^+_{\frak p\frak p_0}$ is the coordinate change
of the good coordinate system $\widetriangle{\mathcal U^+}$ whose
domain is $U^+_{\frak p\frak p_0}$.
Thus Diagram  (\ref{diagram58sec11}) is nothing but Diagram
(\ref{diagram58}) in the current context.
During the proof below, we use the notations used in the proof of Proposition \ref{inductiveprop}.
\par
Now we begin with describing the shrinkings.
For each $i$, let $U^{2 \prime}_{q_i^{\frak p}}$ be a relatively compact open neighborhood
of $q_i^{\frak p}$ in $U^{2}_{q_i^{\frak p}}$ such that
\begin{equation}\label{qcoversuru22}
\aligned
&\bigcup_{i=1}^{N(\frak p)}
\psi^+_{q^{\frak p}_i}((s_{q^{\frak p}_i}^+)^{-1}(0)
\cap U^{2\prime}_{q^{\frak p}_i})
\\
&\supset
\psi^+_{\frak p_0}(\mathcal K_{\frak p_0} \cap
(s^+_{\frak p_0})^{-1}(0))
\cap
\mathcal S_{\frak d}(X,Z;\widehat{\mathcal U})
\cap
\psi^+_{\frak p}(\mathcal K_{\frak p} \cap
(s^+_{\frak p})^{-1}(0)).
\endaligned
\end{equation}
Such a choice is possible because of (\ref{qcoversuru}).
(Recall $U^2_{q_i^{\frak p}}$ is defined so that
(\ref{qcoversuru}) holds.)

\begin{lem}\label{lem11166}
We can take relatively compact opens subsets
$U^{+-}_{\frak p} \subset U_{\frak p}^+$
and $U^{+-}_{\frak p_0} \subset U_{\frak p_0}^+$
with the following properties.
\begin{enumerate}
\item
$
Z_1 \subset
\bigcup_{\frak p}\psi^+_{\frak p}(U^{+-}_{\frak p}
\cap (s^+_{\frak p})^{-1}(0))
$.
\item
$
Z_0 \subset
\psi^+_{\frak p_0}(U^{+-}_{\frak p_0}
\cap (s^+_{\frak p_0})^{-1}(0))
$.
\item
We put
$U^{+-}_{\frak p\frak p_0} = U_{\frak p\frak p_0} \cap U^{+-}_{\frak p_0}
\cap (\varphi^+_{\frak p\frak p_0})^{-1}(U^{+-}_{\frak p})$.
Then
$$
\overline{U^{+-}_{\frak p\frak p_0}}
\subset
\bigcup_{i=1}^{N(\frak p)}
\varphi^0_{\frak p_0 q^{\frak p}_i}(U^{2\prime}_{q^{\frak p}_i}).
$$
\end{enumerate}
\end{lem}
\begin{proof}
We note that the good coordinate system
$\widetriangle{\mathcal U^+}$ defines a Hausdorff space $\vert\widetriangle{\mathcal U^+}\vert$
and
$U_{\frak p}^+$, $U_{\frak p_0}^+$, $Z_0$, $Z_1$ can be identified
as subspaces of $\vert\widetriangle{\mathcal U^+}\vert$.
Therefore the existence of  such $U^{+-}_{\frak p} \subset U_{\frak p}^+$
and $U^{+-}_{\frak p_0} \subset U_{\frak p_0}^+$
is an easy consequence of (\ref{qcoversuru22}).
\end{proof}
By Lemma \ref{lem11166} (1)(2), we can apply Lemma \ref{lem320} to $\widetriangle{\mathcal U^+}$
and $U^{+-}_{\frak p}$, $U^{+-}_{\frak p_0}$
and obtain a good coordinate system whose
Kuranishi neighborhoods are  $U^{+-}_{\frak p}$, $U^{+-}_{\frak p_0}$.
We denote this good coordinate system by  $\widetriangle{\mathcal U^{+-}}$.
It still satisfies the conclusion of  Proposition
\ref{inductiveprop}.

\begin{lem}\label{1112label}
We can take relatively compact
neighborhoods $U^{5}_{\frak p}(p) \subseteq U^{5\prime}_{\frak p}(p)$ and $U^{5}_{\frak p_0}(q) \subseteq U^{5\prime}_{\frak p_0}(q)$
of $o_{\frak p}(p)$  and  $o_{\frak p_0}(q)$ in
$U^+_{\frak p}$ and $U^+_{\frak p_0}$, respectively,
so that if
$$
\aligned
p &\in \psi^+_{\frak p}((s^+_{\frak p})^{-1}(0) \cap U^{+-}_{\frak p})\\
q &\in \psi^+_{\frak p_0}(U^{+-}_{\frak p_0} \cap (s^+_{\frak p_0})^{-1}(0))
\cap
\psi^+_{\frak p}(U^{5}_{\frak p}(p) \cap (s^+_{\frak p})^{-1}(0))
\endaligned$$
then there exists $i \in \{1,\dots, N(\frak p)\}$
such that:
\begin{enumerate}
\item
$
(\varphi^1_{\frak p q^{\frak p}_i})^{-1}(U^{5}_{\frak p}(p)) \subset
U^{2\prime}_{q^{\frak p}_i}
$.
\item
$
U^{5}_{\frak p_0}(q) \subset
\varphi^0_{\frak p_0 q^{\frak p}_i}(U^{2\prime}_{q^{\frak p}_i}).
$
\item
$
U^{5}_{\frak p}(p) \subset
(\varphi^3_{a;p \frak p})^{-1}
(U^+_{a;q_i^{\frak p} p})
\cap U^3_{a;\frak p}(q_i^{\frak p}).
$
\item
$U^{5}_{\frak p_0}(q)
\subset U^4_{a;\frak p_0}(q_i^{\frak p})$.
\item
$\varphi^4_{a;q\frak p_0}(U^{5}_{\frak p_0}(q))
\subset
U^+_{{a;q_i^{\frak p}q}}
\cap (U^+_{a;pq} \cap
(\varphi^+_{a;pq})^{-1}(U^+_{a;q_i^{\frak p}p}))$.
\end{enumerate}
\end{lem}

\begin{proof}
We first observe that
Lemma \ref{lem11166} (3) and Definition \ref{coordinatechangedef} (3)
 imply
\begin{equation}\label{112626}
\aligned
&{\rm Close}\left(
\psi^+_{\frak p}((s^+_{\frak p})^{-1}(0) \cap U^{+-}_{\frak p})
\cap
\psi^+_{\frak p_0}(U^{+-}_{\frak p_0} \cap (s^+_{\frak p_0})^{-1}(0))\right) \\
&\subset
\bigcup_{i=1}^{N(\frak p)}
\psi^+_{\frak p_0}(\varphi^0_{\frak p_0  q^{\frak p}_i}(U^{2\prime}_{q^{\frak p}_i})
\cap (s^+_{\frak p_0})^{-1}(0)).
\endaligned
\end{equation}
Here  ${\rm Close}$ denotes the closure. Since
$$
\psi_{\frak p_0}\circ \varphi^0_{\frak p_0  q^{\frak p}_i}
=
\psi_{q^{\frak p}_i}
=
\psi_{\frak p} \circ \varphi^1_{\frak p  q^{\frak p}_i}
$$
on
$U^{2\prime}_{q^{\frak p}_i}
\cap s_{q^{\frak p}_i}^{-1}(0)$,
(\ref{112626}) implies
\begin{equation}\label{11262622}
\aligned
&{\rm Close}\left(\psi^+_{\frak p}((s^+_{\frak p})^{-1}(0) \cap U^{+-}_{\frak p})
\cap
\psi^+_{\frak p_0}(U^{+-}_{\frak p_0} \cap (s^+_{\frak p_0})^{-1}(0))
\right)\\
&\subset
\bigcup_{i=1}^{N(\frak p)}
\psi^+_{\frak p}(\varphi^1_{\frak p  q^{\frak p}_i}(U^{2\prime}_{q^{\frak p}_i})
\cap (s^+_{\frak p})^{-1}(0)).
\endaligned
\end{equation}
Moreover
\begin{equation}\label{11262622+}
\psi^+_{\frak p_0}(\varphi^0_{\frak p_0  q^{\frak p}_i}(U^{2\prime}_{q^{\frak p}_i})
\cap (s^+_{\frak p_0})^{-1}(0))
=
\psi_{\frak p}^+(\varphi^1_{\frak p  q^{\frak p}_i}(U^{2\prime}_{q^{\frak p}_i})
\cap (s^+_{\frak p})^{-1}(0)).
\end{equation}
By (\ref{11262622}) we can find compact subsets
$C_{q^{\frak p}_i}$ of
$U^{2\prime}_{q^{\frak p}_i}$ such that
\begin{equation}\label{112626222}
\aligned
&\psi^+_{\frak p}((s^+_{\frak p})^{-1}(0) \cap U^{+-}_{\frak p})
\cap
\psi^+_{\frak p_0}(U^{+-}_{\frak p_0} \cap (s^+_{\frak p_0})^{-1}(0)) \\
&\subset
\bigcup_{i=1}^{N(\frak p)}
\psi^+_{\frak p}(\varphi^1_{\frak p  q^{\frak p}_i}(C_{q^{\frak p}_i})
\cap (s^+_{\frak p})^{-1}(0)).
\endaligned
\end{equation}
\begin{sublem}
We may choose $U^{5}_{\frak p}(p)$
to be a sufficiently small neighborhood of $o_{\frak p}(p)$
so that the following holds.
\par
Let
$p \in \psi^+_{\frak p}((s^+_{\frak p})^{-1}(0) \cap U^{+-}_{\frak p})$.
Then there exists $i \in \{1,\dots, N(\frak p)\}$ such that
Lemma \ref{1112label} (1) and (3) hold and
\begin{equation}\label{11262622333}
(\psi^+_{\frak p_0})^{-1}\left(
\psi^+_{\frak p_0}(U^{+-}_{\frak p_0} \cap (s^+_{\frak p_0})^{-1}(0))
\cap
\psi^+_{\frak p}(U^{5}_{\frak p}(p) \cap (s^+_{\frak p})^{-1}(0))
\right)
\subset
 \varphi^0_{\frak p_0 q^{\frak p}_i}(U^{2\prime}_{q^{\frak p}_i}).
\end{equation}
\end{sublem}
\begin{proof}
Let
$p \in \psi^+_{\frak p}((s^+_{\frak p})^{-1}(0) \cap U^{+-}_{\frak p})$.
We may choose $U^{5}_{\frak p}(p)$ sufficiently small so that
the following holds.
\begin{enumerate}
\item[(*)]
If
$$
\psi^+_{\frak p}(U^{5}_{\frak p}(p) \cap (s^+_{\frak p})^{-1}(0))
\cap
\psi^+_{\frak p}(\varphi^1_{\frak p  q^{\frak p}_i}(C_{q^{\frak p}_i})
\cap (s^+_{\frak p})^{-1}(0))\ne
\emptyset,
$$
then
$$
p \in
\psi^+_{\frak p}(\varphi^1_{\frak p  q^{\frak p}_i}(U^{2\prime}_{q^{\frak p}_i})
\cap (s^+_{\frak p})^{-1}(0)).
$$
\end{enumerate}
In fact we can prove the contraposition of (*) by using
the compactness of $C_{q^{\frak p}_i}$.
\par
Let
$q \in \psi^+_{\frak p_0}(U^{+-}_{\frak p_0} \cap (s^+_{\frak p_0})^{-1}(0))
\cap
\psi^+_{\frak p}(U^{5}_{\frak p}(p) \cap (s^+_{\frak p})^{-1}(0))$.
By (\ref{112626222}), there exists $i \in \{1,\dots, N(\frak p)\}$
such that
$q \in \psi^+_{\frak p}(\varphi^1_{\frak p  q^{\frak p}_i}(C_{q^{\frak p}_i})
\cap (s^+_{\frak p})^{-1}(0))$.
For any such $i$, (*) implies
$
p \in
\psi_{\frak p}^+(\varphi^1_{\frak p  q^{\frak p}_i}(U^{2\prime}_{q^{\frak p}_i})
\cap (s^+_{\frak p})^{-1}(0)).
$
We can further shrink $U^{5}_{\frak p}(p)$, if necessary,
so that Lemma \ref{1112label} (1)
is satisfied.
\par
Therefore  (\ref{formula1125}) and
$U^{2\prime}_{q^{\frak p}_i} \subset U^1_{q^{\frak p}_i}$
imply
$
o_{\frak p}(p) \in  U^3_{a;\frak p}(q_i^{\frak p}).
$
We remark that
$
p \in
\psi^+_{\frak p}(\varphi^1_{\frak p  q^{\frak p}_i}(U^{2\prime}_{q^{\frak p}_i})
\cap (s^+_{\frak p})^{-1}(0))
$
implies $o_{\frak p}(p) \in U^3_{a;p\frak p}$.
Moreover
$$
\varphi^3_{a;p\frak p}(o_{\frak p}(p)) = o_{p}(p)
\in U^+_{a;q^{\frak p}_i p}.
$$
(In fact $U^+_{a;q^{\frak p}_i p}$ is nonempty since
$
p \in
\psi_{q^{\frak p}_i}(s_{q^{\frak p}_i}^{-1}(0) \cap
U_{q^{\frak p}_i})
$ follows from Lemma \ref{1112label} (1).)
Hence
$
o_{\frak p}(p) \in
(\varphi^3_{a;p \frak p})^{-1}
(U^+_{q_i^{\frak p} p})
\cap U^3_{a;\frak p}(q_i^{\frak p})
$.
Therefore we can shrink $U^{5}_{\frak p}(p)$  so that
Lemma \ref{1112label} (3) is satisfied.
\par
To prove (\ref{11262622333}) it suffices to show
\begin{equation}\label{form1125}
\psi^+_{\frak p}(U^{5}_{\frak p}(p) \cap (s^+_{\frak p})^{-1}(0))
\subset
\psi^+_{\frak p_0}( \varphi^0_{\frak p_0 q^{\frak p}_i}(U^{2\prime}_{q^{\frak p}_i})).
\end{equation}
The right hand side is
$\psi^+_{\frak p}(\varphi^1_{\frak p  q^{\frak p}_i}(U^{2\prime}_{q^{\frak p}_i})
\cap (s^+_{\frak p})^{-1}(0))$
by (\ref{11262622+}).
Therefore (\ref{form1125}) follows from Lemma \ref{1112label} (1).
\end{proof}
By (\ref{11262622333}) we find that
$o_{\frak p_0}(q) \in \varphi^0_{\frak p_0 q^{\frak p}_i}(U^{2\prime}_{q^{\frak p}_i})$.
Since $\varphi^0_{\frak p_0 q^{\frak p}_i}$ an
open mapping, we can choose $U^{5}_{\frak p_0}(q)$
so that Lemma \ref{1112label} (2) holds.
\par
We will next shrink $U^{5}_{\frak p_0}(q)$ so that
Lemma \ref{1112label} (4)(5) are satisfied.
Lemma \ref{1112label} (2) implies
$o_{\frak p_0}(q) \in \varphi^0_{\frak p_0 q^{\frak p}_i}(U^{2\prime}_{q^{\frak p}_i})
\subseteq \varphi^0_{\frak p_0 q^{\frak p}_i}(U^{1}_{q^{\frak p}_i})$.
Therefore by (\ref{formula1125})
$
o_{\frak p_0}(q) \in U^4_{a;\frak p_0}(q_i^{\frak p}).
$
Hence we can find a neighborhood $U^{5}_{\frak p_0}(q)$
of $o_{\frak p_0}(q)$ such that Lemma \ref{1112label} (4) is satisfied.
\par
Note
Lemma \ref{1112label} (1)(2) implies
\begin{equation}\label{1126262233343}
q \in \psi_{q_i^{\frak p}}(U_{q_i^{\frak p}}
\cap s_{q_i^{\frak p}}^{-1}(0)),
\qquad
p \in \psi_{q_i^{\frak p}}(U_{q_i^{\frak p}}
\cap s_{q_i^{\frak p}}^{-1}(0)).
\end{equation}
In fact, the first half is a consequence of Lemma \ref{1112label}
(2) and
$
q \in \psi^+_{\frak p_0}(U^{5}_{\frak p_0}(q)
\cap (s^+_{\frak p_0})^{-1}(0))
$.
The second half follows from
Lemma \ref{1112label}
(1) as we mentioned already.
\par
Using the existence of a {\it strict} KK-embedding
$\widehat{\Phi^2_a} : \widehat{\mathcal U} \to \widehat{\mathcal U^+_a}$
the second formula of (\ref{1126262233343}) implies
$p \in \psi_{a;q_i^{\frak p}}(U^+_{a;q_i^{\frak p}}
\cap s_{q_i^{\frak p}}^{-1}(0))$.
Therefore
the coordinate change $\varphi^+_{a;q_i^{\frak p}p}$
is defined by Definition \ref{kstructuredefn}.
\par
The first formula of (\ref{1126262233343}) and the existence of
a strict KK-embedding
$\widehat{\Phi^2_a} : \widehat{\mathcal U} \to \widehat{\mathcal U^+_a}$
imply
$$
q \in \psi_{a;q_i^{\frak p}}(U^+_{a;q_i^{\frak p}}
\cap s_{a;q_i^{\frak p}}^{-1}(0)).
$$
Since Lemma \ref{1112label} (1) implies
$q \in \psi^+_{\frak p}((s^+_{\frak p})^{-1}(0) \cap U^{5}_{\frak p}(p))
\subseteq \psi_{a;p}(s_p^{-1}(0) \cap U^+_{a;p})$,
the coordinate change $\varphi_{a;pq}^+$ is defined.
\par
We may choose $U^5_{\frak p}(p)$ such that
$\varphi_{a;p \frak p}^3(U^5_{\frak p}(p)) \subset U^+_{a;q_i^{\frak p} p}$
and so
$$
o_q \in U^+_{{a;q_i^{\frak p}q}}
\cap (U^+_{a;pq} \cap
(\varphi^+_{a;pq})^{-1}(U^+_{a;q_i^{\frak p}p})).
$$
Therefore we can find a neighborhood $U^{5}_{\frak p_0}(q)$
of $o_{\frak p_0}(q)$ so that Lemma \ref{1112label} (5) is satisfied.
The proof of Lemma \ref{1112label} is now complete.
\end{proof}
Now we are ready to prove the commutativity of
Diagram (\ref{diagram58sec11}).
Our assumption is
$p \in \psi^+_{\frak p}(U^{+-}_{\frak p} \cap (s^+_{\frak p})^{-1}(0))$ and
$q \in \psi^+_{\frak p_0}(U^{+-}_{\frak p_0} \cap (s^+_{\frak p_0})^{-1}(0))
\cap
\psi^+_{\frak p}(U^{5}_{\frak p}(p) \cap (s^+_{\frak p})^{-1}(0)).
$
We choose $i$ as in Lemma \ref{1112label}.
Let
$$
y \in (U^+_{\frak p\frak p_0} \cap (\varphi_{\frak p\frak p_0}^+)^{-1}(U^{5}_{\frak p}(p)))
\cap ( U^{5}_{\frak p_0}(q)
\cap(\varphi_{a;q\frak p_0}^4)^{-1}({U}_{a;pq}^{+})).
$$
Then, since $y\in U^{5}_{\frak p_0}(q)$, Lemma \ref{1112label} (2)  implies
$y = \varphi^0_{\frak p_0  q^{\frak p}_i}(\tilde y)$
with $\tilde y \in U^{2 \prime}_{q^{\frak p}_i}$.
\par
Now we calculate
$$
\aligned
(\varphi^+_{a;q_i^{\frak p}p}
\circ \varphi^{+}_{a;pq}\circ\varphi^4_{a;q \frak p_0})(y)
&=
\varphi^{+}_{a;q_i^{\frak p}q}(\varphi^4_{a;q \frak p_0}(y)) \\
&=
\varphi^4_{a;q^{\frak p}_i \frak p_0}(y)\\
&=
\varphi^4_{a;q^{\frak p}_i \frak p_0}( \varphi^0_{\frak p_0  q^{\frak p}_i}(\tilde y))
\\
&=
\varphi^2_{a; q^{\frak p}_i}(\tilde y).
\endaligned
$$
The first equality is the compatibility condition
of the coordinate change of $\widehat{\mathcal U^+_a}$.
(Here we use Lemma \ref{1112label} (5)
to apply the compatibility condition.)
The second equality is the consequence of the fact that
$\varphi^4_{a;q \frak p_0}$ is a part of the object consisting GK-embedding.
(We use Lemma \ref{1112label} (4)(5) to apply the compatibility condition.)
The third line is the definition of $\tilde y$.
The fourth equality follows from (\ref{formula1125}) and Assumption (c).
\par
On the other hand, we have
$$
\aligned
(\varphi^+_{a;q_i^{\frak p}p}
\circ\varphi^3_{a;p \frak p}\circ \varphi^+_{\frak p\frak p_0})(y)
&=
\varphi^+_{a;q_i^{\frak p}p}
(\varphi^3_{a;p \frak p}(\varphi^+_{\frak p\frak p_0}(y))) \\
&=
\varphi^3_{a;q_i^{\frak p} \frak p}(\varphi^+_{\frak p\frak p_0}(y)) \\
&=
\varphi^3_{a;q_i^{\frak p} \frak p}(\varphi^+_{\frak p\frak p_0}( \varphi^0_{\frak p_0  q^{\frak p}_i}(\tilde y))) \\
&=
\varphi^3_{a;q_i^{\frak p} \frak p}(\varphi^1_{\frak p  q^{\frak p}_i}(\tilde y)) \\
&=
\varphi^2_{a; q^{\frak p}_i}(\tilde y).
\endaligned
$$
The first equality is obvious.
The second equality  follows from the fact that
$\widehat{\Phi^3}$ is a GK-embedding
and Lemma \ref{1112label} (3).
(Note
$\varphi^+_{\frak p\frak p_0}(y) \in U^{5}_{\frak p}(p) \subset U^{3}_{\frak p}(p)$.)
The third equality is the definition of $\tilde y$.
The fourth equality is the definition of $\varphi^+_{\frak p\frak p_0}$,
that is (\ref{formform1177}).
The fifth equality follows from (\ref{formula1125}) and Assumption (b).
\par
Therefore we obtain
$$
(\varphi^+_{a;q_i^{\frak p}p}
\circ \varphi^{+}_{a;pq}\circ\varphi^4_{a;q \frak p_0})(y)
=
(\varphi^+_{a;q_i^{\frak p}p}
\circ\varphi^3_{a;p \frak p}\circ \varphi^+_{\frak p\frak p_0})(y).
$$
Since $\varphi^+_{a;q_i^{\frak p}p}$ is injective, we have
proved the commutativity of the
Diagram (\ref{diagram58sec11}).
\par
The proof of Proposition \ref{prop1111} is complete.
\end{proof}
The rest of the proof of  Proposition \ref{prop519}
is mostly the same as the proof of
Theorem
\ref{Them71restate}, using Proposition \ref{prop1111}
in addition to Proposition \ref{inductiveprop}.
We use the same notation as in the last part of the
proof of Theorem
\ref{Them71restate}.
We  construct
a good coordinate system of a compact neighborhood of
$$
Z \cup \bigcup_{i=1}^n Z_i
$$
together with its GK-embedding to $\widehat{\mathcal U^{+}_a}$
by induction.
Suppose we have a
good coordinate system
of a compact neighborhood of
$
Z \cup \bigcup_{i=1}^{n-1} Z_i
$
together with its embedding to $\widehat{\mathcal U^{+}_a}$.
To apply Propositions \ref{inductiveprop}, \ref{prop1111},
we need to find a good coordinate system on $Z_n$
together with the GK-embedding to $\widehat{\mathcal U^{+}_a}$.
We use the following lemma for this purpose.
\begin{lem}\label{charthuyasilemma}
Let $\widehat{\mathcal U}$ be a Kuranishi structure
and $\widehat{\mathcal U^{+}_a}$ its
thickenings for $a=1,2$. Then for each $p
\in \mathcal S_{\frak d}(X,Z;\widehat{\mathcal U})$
there exists a good coordinate system
$\widetriangle{\mathcal U_{\frak p_0}}$
of $\{p\} \subseteq X$ with the following
properties.
\begin{enumerate}
\item
$\widetriangle{\mathcal U_{\frak p_0}}$ consists of a single Kuranishi chart $\mathcal U_{\frak p_0}$
such that $\dim U_{\frak p_0} = \frak d$
and $\mathcal U_{\frak p_0}$ is a restriction to an open set of a
Kuranishi chart of a point $p$ of $\widehat{\mathcal U}$.
\item
For each $a=1,2$,
there exists a GK-embedding
$\widetriangle{\mathcal U_{\frak p_0}}
\to \widehat{\mathcal U_a^{+}}\vert_{Z_0}$, where $Z_0$ is a
compact neighborhood of $p$ in $Z$.
\end{enumerate}
\end{lem}
\begin{proof}
Let $\mathcal U_p$
(resp. $\mathcal U^{+}_{a;p}$) be the Kuranishi chart of $p$
induced by the Kuranishi structure $\widehat{\mathcal U}$
(resp. $\widehat{\mathcal U^{+}_{a}}$).
Let $O_{a;p}$ be the neighborhood of $p$ in $X$ as in
Definition \ref{thickening}.
(Here we use thickness of $\widehat{\mathcal U^{+}_a}$.)
We put $O_p = O_{1;p} \cap O_{2;p}$.
We take an open neighborhood $U_{\frak p_0}$ of $o_p$ in $U_p$ such that
\begin{equation}
\psi_{p}(U_{\frak p_0} \cap s_p^{-1}(0)) \subset O_p,
\end{equation}
and set $\mathcal U_{\frak p_0} = \mathcal U_p\vert_{U_{\frak p_0}}$.
It is obvious that $\mathcal U_{\frak p_0}$
satisfies (1).
Let us prove (2). The subset $Z_0$ is any compact
neighborhood of $p$ contained in $O_{p}$.
\par
Let $q \in \psi(U_{\frak p_0} \cap s_{\frak p_0}^{-1}(0)) \cap Z_0$.
We take $W_{a;p}(q) \subset U_p$ as in
Definition \ref{thickening} and put
$U(q) = W_{1;p}(q) \cap W_{2;p}(q)$.
Then by Definition \ref{thickening} (2) (a) we have
$$
\varphi_{a;p}(U(q)) \subset \varphi^{+}_{a;pq}(U^{+}_{a;pq}).
$$
Since $\varphi^{+}_{a;pq}$ is injective,
we have a set theoretical map
$\varphi_{a;q\frak p_0} : U(q) \to U^{+}_{a;pq}$
such that $\varphi^{+}_{a;pq} \circ \varphi_{a;q\frak p_0}
= \varphi_{a;p}$.
Since $\varphi^{+}_{a;pq}$ and $\varphi_{a;p}$ are embeddings between orbifolds,
the map  $\varphi_{a;q \frak p_0}$ is an embedding of orbifolds.
We can use Definition \ref{thickening} (2) (b)
in the same way to obtain
an embedding of vector bundles $\widehat\varphi_{a;q\frak p_0} :
E_{\frak p_0} \to E^+_{a;q}$ such that it covers $\varphi_{a;q\frak p_0}$
and satisfies $\widehat\varphi^{+}_{a;pq}\circ \widehat\varphi_{a;q\frak p_0}
= \widehat\varphi_{a;p}$.
We thus obtain $\Phi_{a;q\frak p_0} = (U(q),\widehat\varphi_{a;q\frak p_0})$ for each
$q \in \psi(U_{\frak p_0} \cap s_{\frak p_0}^{-1}(0))$.
It is easy to see that $\Phi_{a;q\frak p_0}$ defines an embedding of
a good coordinate system to a Kuranishi structure.
The proof of Lemma \ref{charthuyasilemma} is complete.
\end{proof}
\begin{rem}
The proof of this lemma is the place where we use the assumption that
$\widehat{\mathcal U^{+}_a}$ is
a thickening of $\widehat{\mathcal U}$.
\end{rem}
Since  ${\mathcal U}_{\frak p_0}$ is an open subchart of a
Kuranishi chart of $\widehat{\mathcal U}$
there exists a KG-embedding $\widehat{\mathcal U}
\to  \widetriangle{\mathcal U_{\frak p_0}}$.
Moreover the composition
$\widehat{\mathcal U}
\to\widetriangle{\mathcal U_{\frak p_0}} \to \widehat{\mathcal U^{+}_a}$
coincides with the given KK-embedding on a neighborhood of $p$.
In other words Assumption (c) of Proposition \ref{prop1111} is
satisfied.
Therefore we apply Proposition \ref{prop1111}
inductively to complete the proof of
Proposition \ref{prop519}.
\end{proof}

\subsection{KG-embeddings and compatible perturbations}
\label{subsec:movingmulsectionetc}

The present section will be occupied by the proofs of Propositions \ref{le614}, \ref{pro616}
and Lemmata \ref{le714}, \ref{le7155}, \ref{lem92929}.

\begin{proof}[Proof of Proposition \ref{le614} and Lemma \ref{le714}]
We use the same induction scheme as in the proof of
Theorem
\ref{Them71restate}.
\begin{lem}\label{lem11115}
We consider the situation of Proposition \ref{inductiveprop}
and Situation \ref{situation101}.
\begin{enumerate}
\item
Suppose there exists a system of multivalued perturbations
$\widetriangle{\frak s} = \{\frak s^{n}_{\frak p}\}$ of ${\widetriangle{\mathcal U}}$,
$\widetriangle{\frak s_{\frak p_0}} = \{\frak s^{n}_{\frak p_0}\}$ of $\widetriangle{\mathcal U_{\frak p_0}}$,
and $\widehat{\frak s} = \{\frak s^{n}_{p}\}$ of ${\widehat{\mathcal U}}$.
We assume that they are compatible with the KG-embeddings
$\widehat{\Phi^0} : {\widehat{\mathcal U}}\vert_{Z_0^+} \to \widetriangle{\mathcal U_{\frak p_0}}$
and
$\widehat{\Phi^1} : {\widehat{\mathcal U}}\vert_{Z_1^+} \to {\widetriangle{\mathcal U}}$.
\par
Then there exists a system of multivalued perturbations
$\widetriangle{\frak s^+} =  \{\frak s^{n +}_{\frak p}\}$ of ${\widetriangle{\mathcal U^+}}$
with the following properties.
\begin{enumerate}
\item
$\widetriangle{\frak s^+}$, $\widehat{\frak s}$
are compatible with the strict KG-embedding $\widehat{\Phi^+} :
{\widehat{\mathcal U_0}}\vert_{Z_+^+} \to {\widetriangle{\mathcal U^+}}$,
where ${\widehat{\mathcal U_0}}$ is an open substructure of
$\widehat{\mathcal U}$.
\item
If $\frak p \in \frak P$,
then $\frak s^{n +}_{\frak p}$ is the restriction of
$\frak s^{n}_{\frak p}$ to $\mathcal U_{\frak p}^+$.
\item
In case of $\frak p_0$,
$\frak s^{n +}_{\frak p_0}$ is the restriction of
$\frak s^{n}_{\frak p_0}$ to $\mathcal U_{\frak p_0}^+$.
\item If $\widetriangle{\frak s}$, $\widehat{\frak s}$
are transversal to $0$, then so is  $\widetriangle{\frak s^+}$.
\end{enumerate}
\item
Suppose there exists a CF-perturbation
$\widetriangle{
\frak S} = \{\frak S^{\epsilon}_{\frak p}\}$ of ${\widetriangle{\mathcal U}}$,
$\widetriangle{\frak S_{\frak p_0}} = \{\frak S_{\frak p_0}^{\epsilon}\}$ of $\widetriangle{\mathcal U_{\frak p_0}}$, and
$\widehat{
\frak S} = \{\frak S^{\epsilon}_{p}\}$ of $\widehat{{\mathcal U}}$.
We assume that they are compatible with the KG-embeddings
$\widehat{\Phi^0} : {\widehat{\mathcal U}}\vert_{Z_0^+} \to \widetriangle{\mathcal U_{\frak p_0}}$
and
$\widehat{\Phi^1} : {\widehat{\mathcal U}}\vert_{Z_1^+} \to {\widetriangle{\mathcal U}}$.
\par
Then there exists a CF-perturbation
$\widetriangle{
\frak S^+} =\{\frak S^{\epsilon +}_{\frak p}\}$ of ${\widetriangle{\mathcal U^+}}$
with the following properties.
\begin{enumerate}
\item
$\widehat{
\frak S}$, $\widetriangle{
\frak S^+}$
are compatible with the strict KG-embedding $\widehat{\Phi^+} :
{\widehat{\mathcal U}}_0 \to {\widetriangle{\mathcal U^+}}$,
where ${\widehat{\mathcal U_0}}$ is an open substructure of
$\widehat{\mathcal U}$.
\item
If $\frak p \in \frak P$,
then $\frak S^{\epsilon +}_{\frak p}$ is the restriction of
$\frak S^{\epsilon}_{\frak p}$ to $\mathcal U_{\frak p}^+$.
\item
In case of $\frak p_0$,
$\frak S^{\epsilon +}_{\frak p_0}$ is the restriction of
$\frak S^{\epsilon}$ to $\mathcal U_{\frak p_0}^+$.
\item If $\widetriangle{\frak S}$, $\widehat{\frak S}$
are transversal to $0$, then so is  $\widetriangle{\frak S^+}$.
\end{enumerate}
\item
Suppose there exists a differential form (resp. strongly continuous map
to a manifold $M$)
$\widetriangle h = \{h_{\frak p}\}$ (resp. $
\widetriangle f = \{f_{\frak p}\}$) on ${\widetriangle{\mathcal U}}$,
$\widetriangle{h_{\frak p_0}}$ (resp. $\widetriangle{f_{\frak p_0}}$) on $\widetriangle{\mathcal U_{\frak p_0}}$, and
$\widehat{h} = \{h_{p}\}$ (resp. $\widehat{f} = \{f_{p}\}$) of ${\widehat{\mathcal U}}$.
We assume $(\widehat{\Phi^{0}})^*(\widetriangle{h_{\frak p_0}}) =
\widehat{h}\vert_{Z_0^+}$
(resp. $\widetriangle{f_{\frak p_0}} \circ \widehat{\Phi^{0}}=
\widehat{f}\vert_{Z_0^+}$)
and
$(\widehat{\Phi^{1}})^*(\widetriangle h) = \widehat{h}\vert_{Z_1^+}$
(resp. $\widetriangle f \circ \widehat{\Phi^{1}} = \widehat{f}\vert_{Z_1^+}$).
\par
Then there exists a differential form (resp. strongly continuous map to $M$)
$\widetriangle{h^+} =  \{h^{+}_{\frak p}\}$ (resp.$\widetriangle{f^+}
= \{f^{+}_{\frak p}\}$ ) of ${\widetriangle{\mathcal U^+}}$
with the following properties.
\begin{enumerate}
\item
$(\widehat{\Phi^{+}})^*(\widetriangle{h^+}) = \widehat{h}\vert_{\widehat{\mathcal U_0}}$
(resp. $\widetriangle{f^+} \circ \widehat{\Phi^{+}} = \widehat{f}\vert_{\widehat{\mathcal U_0}}$) holds.
Here $\widehat{\mathcal U_0}$ is an open substructure  of
$\widehat{\mathcal U}$ and
$\widehat{\Phi^+} : {\widehat{\mathcal U}}_0 \to {\widetriangle{\mathcal U^+}}$ is a
strict KG-embedding.
\item
If $\frak p \in \frak P$,
then $h^+_{\frak p}$ (resp. $f^+_{\frak p}$) is a restriction of
$h_{\frak p}$ (resp.$f_{\frak p}$ ) to $\mathcal U_{\frak p}^+$.
\item
In case of $\frak p_0$,
$h^+_{\frak p_0}$ (resp. $f^+_{\frak p_0}$) is a restriction of
$h_{\frak p_0}$ (resp. $f_{\frak p_0}$) to $\mathcal U_{\frak p_0}^+$.
\end{enumerate}
\item Suppose we are in the situation of (1). Let $\widehat f$,
$\widetriangle f$ be as in (3).
\begin{enumerate}
\item
If $\widehat f$,
$\widetriangle f$ are strongly submersive with respect to $\widehat{\frak s}$
and $\widetriangle{\frak s}$, respectively, then
$\widetriangle{f^{+}}$ is
strongly submersive with respect to $\widetriangle{\frak s^+}$.
\item
Let $g : N\to M$ be a smooth map such that
$\widehat f$,
$\widetriangle f$ are strongly transversal to $g$ with respect to
$\widehat{\frak s}$
and $\widetriangle{\frak s}$, respectively.
Then $\widetriangle{f^{+}}$ is
strongly transversal to $g$ with respect to $\widetriangle{\frak s^+}$.
\end{enumerate}
\item Suppose we are in the situation of (2). Let $\widehat f$,
$\widetriangle f$ be as in (3).
\begin{enumerate}
\item
If  $\widehat f$,
$\widetriangle f$ are strongly submersive with respect to $\widehat{\frak S}$
and $\widetriangle{\frak S}$, respectively, then
$\widetriangle{f^{+}}$ is
strongly submersive with respect to $\widetriangle{\frak S^+}$.
\item
Let $g : N\to M$ be a smooth map such that
$\widehat f$,
$\widetriangle f$ are strongly transversal to $g$ with respect to
$\widehat{\frak S}$
and $\widetriangle{\frak S}$, respectively.
Then $\widetriangle{f^{+}}$ is
strongly transversal to $g$ with respect to $\widetriangle{\frak S^+}$.
\end{enumerate}
\end{enumerate}
\end{lem}
\begin{proof}
We will prove (1),(4). The proofs of (2),(3) and (5) are entirely similar.
\par
We note that (1) (b) and (c) uniquely determine
$\{\frak s^{n +}_{\frak p}\}$.
Its compatibility with the coordinate change
$\Phi^+_{\frak p\frak q}$ for $\frak p,\frak q\in \frak P$ follows
from the assumption, that is, the compatibility of
$\{\frak s^{n}_{\frak p}\}$ with
$\Phi_{\frak p\frak q}$.
Therefore to complete the proof, it suffices to check the following
three points.
\begin{enumerate}
\item[(I)]
$\frak s^{n +}_{\frak p}$, $\frak s^{n +}_{\frak p_0}$
are compatible with $\Phi^+_{\frak p \frak p_0}$.
\item[(II)]
Statement (1)(a) holds.
\item[(III)]
Statements (1)(d) and (4) hold.
\end{enumerate}
\begin{proof}[Proof of (I)]
Let $y \in U^+_{\frak p\frak p_0}$.
Since $U^+_{\frak p\frak p_0} \subset U^1_{\frak p\frak p_0}$,
(\ref{118888}) implies that there exist $q_i^{\frak p}$ and
$\tilde y \in U^2_{q_i^{\frak p}}$ such that
$\varphi^0_{\frak p_0 q_i^{\frak p}}(\tilde y) = y$.
We can take a representative
$(\frak s^{n}_{\frak p_0,k})_{k=1,\dots,\ell}$
of $\frak s^{n}_{\frak p_0}$
(resp. $({\frak s}^{n}_{q_i^{\frak p},k})_{k=1,\dots,\ell}$
of $\frak s^{n}_{q_i^{\frak p}}$)
on a neighborhood $U_y$ of $y$ (resp. $U_{\tilde y}$ of $\tilde y$)
such that
\begin{equation}\label{112555}
\frak s^{n}_{\frak p_0,k}(\varphi^0_{\frak p_0 q_i^{\frak p}}(\tilde z))
=
\widehat{\varphi^0_{\frak p_0 q_i^{\frak p}}}(\frak s^{n}_{q_i^{\frak p},k}(\tilde z))
\end{equation}
holds for any $\tilde z \in U_{\tilde y}$.
This is a consequence of the compatibility of
$\{\frak s_{\frak p_0}^{n}\}$ and $\{\frak s_{\frak p}^{n}\}$
with $\widehat{\Phi^0}$.
\par
We put $\overline y = \varphi^+_{\frak p\frak p_0}(y) \in U^+_{\frak p}$.
Then we have $\overline y = \varphi^1_{\frak pq_i^{\frak p}}(\tilde y)$.
(This is the definition of $\varphi^+_{\frak p\frak p_0}$.)
We can take a representative
$(\frak s^{n}_{\frak p,k})_{k=1,\dots,\ell}$
of $\frak s^{n}_{\frak p}$ such that
\begin{equation}\label{112666}
\frak s^{n}_{\frak p,k}(\varphi^1_{\frak p q_i^{\frak p}}(\tilde z))
=
\widehat{\varphi^1_{\frak p q_i^{\frak p}}}(\frak s^{n}_{q_i^{\frak p},k}(\tilde z))
=
(\widehat{\varphi^+_{\frak p \frak p_0}}\circ
\widehat{\varphi^0_{\frak p_0 q_i^{\frak p}}})
(\frak s^{n}_{q_i^{\frak p},k}(\tilde z))
\end{equation}
holds for any $\tilde z \in U_{\tilde y}$.
Here the first equality is a consequence of the compatibility of
$\{\frak s^{n}_{\frak p}\}$ and $\{\frak s^{n}_{p}\}$
with $\widehat{\Phi^1}$
and the second equality is the definition of
$\widehat{\varphi^+_{\frak p \frak p_0}}$.
\par
We also have
\begin{equation}\label{112777}
\varphi^1_{\frak p q_i^{\frak p}}(\tilde z)
=
\varphi^+_{\frak p \frak p_0}(\varphi^0_{\frak p_0 q_i^{\frak p}}(\tilde z))
\end{equation}
by definition of $\varphi^+_{\frak p \frak p_0}$,
which is (\ref{formform1177}).
Then (\ref{112555}), (\ref{112666}), (\ref{112777}) imply
$$
\frak s^{n}_{\frak p,k}(\varphi^+_{\frak p \frak p_0}(z))
=
\widehat{\varphi^+_{\frak p \frak p_0}}(\frak s^{n}_{{\frak p}_0,k}(z))
$$
for all $z \in \varphi^0_{\frak p_0 q_i^{\frak p}}(U^2_{q_i^{\frak p}})$.
Since $\varphi^0_{\frak p_0 q_i^{\frak p}}(U^2_{q_i^{\frak p}})$ is a
neighborhood of $y$, this implies (I).
\end{proof}
\begin{proof}[Proof of (II)]
Suppose $p \in {\rm Im}(\psi^+_{\frak p}) \cap Z_1^+$,
$\frak p \in \frak P$.
We need to prove $(\Phi^{+}_{\frak p p})^*(\frak s^{n +}_{\frak p})
= \frak s^{n}_p$,
where $\Phi_{\frak p p}^+ : \mathcal U_{0,p} \to \mathcal U_{\frak p}$.
This is a consequence of the fact that
$\{\frak s^{n}_{\frak p}\}$, $\{ \frak s^{n}_p\}$
are compatible with the embedding $
\widehat{\Phi^1} : \widehat{\mathcal U}\vert_{Z_1^+}
\to {\widetriangle{\mathcal U}}$.
We can prove the case of $\frak p = \frak p_0$ in the
same way.
\end{proof}
\begin{proof}[Proof of (III)]
This is an immediate consequence of the fact that
transversality to $0$, strong submersivity, and
transversality to a map, is preserved under the restriction to
open subsets.
\end{proof}
The proof of Lemma \ref{lem11115} is complete.
\end{proof}
Using Lemma \ref{lem11115},  we can prove
Proposition \ref{le614} and  Lemma \ref{le714},
by inspecting the proof of Theorem
\ref{Them71restate}.
\end{proof}
\begin{proof}
[Proof of Proposition \ref{pro616}  and Lemmata \ref{le7155}, \ref{lem92929}]
We use the next lemma for the proof.
\begin{lem}\label{lemma1116}
Suppose we are in the situation of Proposition
\ref{prop1111}.
\begin{enumerate}
\item
We suppose that the assumptions of Lemma \ref{lem11115} (1)
are satisfied, in addition. Moreover we assume that
there exists a multivalued perturbation
$\widehat{\frak s^+_a} = \{\frak s^{n +}_{a;p}\}$ of $\widehat{\mathcal U^{+}_a}$ such that:
\begin{enumerate}
\item[(i)]
$\widehat{\frak s^+_a}$, $\widehat{\frak s}$ are compatible with the
KK-embedding
$\widehat{\Phi^2_a} : \widehat{\mathcal U} \to \widehat{\mathcal U^{+}_a}\vert_{Z_1^+}$.
\item[(ii)]
$\widehat{\frak s^+_a}$, $\widetriangle{\frak s}$ are compatible with the GK-embedding
$\widehat{\Phi^3_a} : {\widetriangle{\mathcal U}} \to \widehat{\mathcal U^{+}_a}\vert_{Z_0^+}$.
\item[(iii)]
$\widehat{\frak s^+_a}$, $\widetriangle{\frak s_{\frak p_0}}$ are compatible with the GK-embedding
$\widehat{\Phi^4_a} :  \widetriangle{\mathcal U_{\frak p_0}} \to \widehat{\mathcal U^{+}_a}$.
\end{enumerate}
Then we can take the multivalued perturbation
$\widetriangle{\frak s^+} = \{\frak s^{n +}_p\}$
of $\widetriangle{\mathcal U^+}$ as in Lemma \ref{lem11115} (1)
such that
$\widehat{\frak s^+_a}$, $\widetriangle{\frak s^+}$ are compatible with
the GK-embedding $\widehat{\Phi^5_a} : \widetriangle{\mathcal U^+} \to \widehat{\mathcal U^{+}_a}$.
\item
A  statement similar to (1) for the CF-perturbations
holds.
\item
A  statements similar to (1) for the differential forms and for
strongly continuous maps hold.
\end{enumerate}
\end{lem}
We omit the precise statement for (2)(3) above. We believe
that it is not difficult to find it for the reader.
\begin{proof}
We prove (1). The proofs of (2) (3) are entirely similar.
Let $p \in \psi^+_{\frak p}((s^+_{\frak p})^{-1}(0)) \cap  Z$.
It suffices to show that $\widehat{\frak s^+_a}$, $\widetriangle{\frak s^+}$
are compatible with the embedding
$\Phi^5_{a;p\frak p} : \mathcal U^+_{\frak p}\vert_{U^+_{\frak p}(p)}
\to \mathcal U_{a;p}^{+}$.
In case $\frak p \in \frak P$ this is a consequence of
Lemma \ref{lem11115} (1)(b) and (ii).
In case $\frak p = \frak p_0$ this is a consequence of
Lemma \ref{lem11115} (1)(c) and (iii).
\end{proof}
Using Lemma \ref{lemma1116}, we can prove Proposition \ref{pro616}  and Lemma \ref{le7155} by inspecting the
proof of Proposition \ref{prop519}.
\par
Lemma \ref{lem92929} is immediate from construction.
\end{proof}
\subsection{Extension of good coordinate systems: the relative case 2}
\label{subsec:moreversionegcs2}

The present section will be occupied by the proofs of Proposition \ref{prop7582752} and Lemma \ref{lem753753}.

\begin{proof}[Proof of Proposition \ref{prop7582752}]
In the statement of Proposition \ref{prop7582752} we used the symbols
$Z_1, Z_2$ for compact subsets of $X$.
In the proof below we use
the symbols $\mathcal Z_{(1)}, \mathcal Z_{(2)}$
for the compact subsets $Z_1, Z_2$ in  Proposition \ref{prop7582752}
to distinguish them
from $Z_0, Z_1$ that appear in Proposition \ref{inductiveprop}.
\par
To prove Proposition \ref{prop7582752} we use the same
induction scheme as the proof of Theorem \ref{Them71restate}.
We will modify the statement of Proposition \ref{inductiveprop}
to Lemma \ref{lem1118} below.
We begin with modifying Situation \ref{situation101}.

In Proposition \ref{prop7582752}
we considered ${\widetriangle{\mathcal U^{1}}}$.
We write  ${\widetriangle{\mathcal U^{(1)}}}$ hereafter in this subsection
in place of ${\widetriangle{\mathcal U^{1}}}$.
(We also write ${\widetriangle{\mathcal U^{(2)}}}$ hereafter in this subsection
in place of ${\widetriangle{\mathcal U^{2}}}$.)
Let  ${\widetriangle{\mathcal U^{(1)}}}
= (\frak P(\mathcal Z_{(1)}),\{\mathcal U^{(1)}_{\frak p'}\},
\{\Phi^{(1)}_{\frak p'\frak q'}\})$.
(We denote elements of $\frak P(\mathcal Z_{(1)})$ by $\frak p'$, $\frak q'$,
that is, by small German characters with prime.)
\par
Let $(\mathcal K^{(1)},\mathcal K^{(1) +})$ be a support pair of
${\widetriangle{\mathcal U^{(1)}}}$.
We put
$$
Z_{\frak p'} = \psi^{(1)}_{\frak p'}((s_{\frak p'}^{(1)})^{-1}(0) \cap \mathcal K^{(1)}_{\frak p'}).
$$
Let
$
\frak U(Z_{\frak p'})
$
be a relatively compact open neighborhood of $Z_{\frak p'}$
in
$
\psi_{\frak p'}^{(1)}((s^{(1)}_{\frak p'})^{-1}(0) \cap {\rm Int}\,
\mathcal K^{(1) +}_{\frak p'})
$.
In Proposition \ref{prop7582752} a Kuranishi structure
$\widehat{\mathcal U^2}$ is given
as  a part of the assumption.
During the proof, we will write $\widehat{\mathcal U}$
in place of $\widehat{\mathcal U^2}$.

\begin{shitu}\label{situation117}
Let
$\frak d \in \Z_{\ge 0}$, and let
$Z_0$ be a compact subset of
$$
\mathcal S_{\frak d}(X,\mathcal Z_{(2)};\widehat{\mathcal U})
\setminus
\bigcup_{\frak d' > \frak d}\mathcal S_{\frak d'}(X,\mathcal Z_{(2)};\widehat{\mathcal U})
\setminus
\bigcup_{\frak p' \in \frak P(\mathcal Z_{(1)})} Z_{\frak p'},
$$
and  $Z_1$  a compact subset of
$$
\mathcal S_{\frak d}(X,\mathcal Z_{(2)};\widehat{\mathcal U})
\cup
\bigcup_{\frak p' \in \frak P(\mathcal Z_{(1)})} \frak U(Z_{\frak p'}).
$$
We assume that $Z_1$ contains an open neighborhood
of $\mathcal Z_{(1)} \cup \bigcup_{\frak d' > \frak d}\mathcal S_{\frak d'}(X,
\mathcal Z_{(2)};\widehat{\mathcal U})$ in
$\mathcal S_{\frak d}(X,\mathcal Z_{(2)};\widehat{\mathcal U})$.
We put
$$
Z_+ = Z_0 \cup Z_1.
$$
\par
Let
${\widetriangle{\mathcal U}}
= (\frak P,\{\mathcal U_{\frak p}\},\{\Phi_{\frak p\frak q}\})$  be
a good coordinate system on a compact neighborhood $Z_1^+$ of
$Z_1$.
We assume that $\frak P$ is written as
$$
\frak P = \frak P({\mathcal Z}_{(1)}) \cup \frak P_0
$$
(disjoint union) and the inclusions
$\frak P({\mathcal Z}_{(1)})
\to \frak P$, $\frak P_0
\to \frak P$  preserve the partial order.
Moreover we assume
that, for $\frak p' \in \frak P({\mathcal Z}_{(1)})$,
the Kuranishi chart
$\mathcal U_{\frak p'}$ of ${\widetriangle{\mathcal U}}$
is an open
subchart of the Kuranishi chart $\mathcal U^{(1)}_{\frak p'}$
of ${\widetriangle{\mathcal U^{(1)}}}$
and
$$
U_{\frak p'} \cap \mathcal Z_{(1)} = U_{\frak p'}^{(1)}
\cap \mathcal Z_{(1)} .
\footnote{Compare Definition \ref{defn735f} (1)(b).}
$$
Furthermore we assume $\dim U_{\frak p} \ge \frak d$
for $\frak p \in \frak P_0$.
\par
Let
$\widehat{\Phi^1}
= \{\Phi^1_{\frak p p} \mid
p \in {\rm Im}(\psi_{\frak p}) \cap Z^+_1\} : \widehat{\mathcal U}\vert_{Z^+_1} \to {\widetriangle{\mathcal U}}$
be a strict KG-embedding
such that, for $\frak p' \in \frak P({\mathcal Z}_{(1)})$, the embedding
$\Phi^1_{\frak p' p}$ is an open restriction of one that is a part of
the given KG-embedding $\widehat{\mathcal U}\vert_{\mathcal Z_{(1)}}
\to {\widetriangle{\mathcal U
^{(1)}}}$.
Let $Z_0^+$ be a compact neighborhood of $Z_0$ in
$\mathcal S_{\frak d}(X,
\mathcal Z_{(2)};\widehat{\mathcal U})$ and
$\mathcal U_{\frak p_0} = (U_{\frak p_0},E_{\frak p_0},s_{\frak p_0},\psi_{\frak p_0})$  a Kuranishi neighborhood of
$Z_0^+$ such that $\dim U_{\frak p_0} = \frak d$.
We regard $\mathcal U_{\frak p_0}$  as a good coordinate system
$\widetriangle{\mathcal U_{\frak p_0}}$ that
consists of a single
Kuranishi chart and suppose that we are given a strict KG-embedding
$\Phi^0
= \{\Phi^0_{\frak p_0} \mid
p \in {\rm Im}(\psi) \cap Z_0\} : \widehat{\mathcal U}\vert_{Z^+_0} \to {\widetriangle{\mathcal U_{\frak p_0}}}$.
\par
We put $\frak P^+ = \frak P \cup \{\frak p_0\}$,
so $\frak P^+ \supset \frak P(\mathcal Z_{(1)})$.
$\blacksquare$
\end{shitu}
\begin{lem}\label{lem1118}
In Situation \ref{situation117}
there exists a good coordinate system
${\widetriangle{\mathcal U
^{+}}}$ of $Z_+^+$ satisfying conclusions
(1)-(5) of  Proposition \ref{inductiveprop}.
Here $Z_+^+$ is a compact neighborhood of $Z_+$ in $X$.
\par
Moreover the following holds.
\begin{enumerate}
\item[(6)]
If $\frak p' \in \frak P(\mathcal Z_{(1)}) \subset \frak P^+$,
then the Kuranishi chart
$\mathcal U^+_{\frak p'}$ of ${\widetriangle{\mathcal U^+}}$
is an open
subchart of the Kuranishi chart $\mathcal U^{(1)}_{\frak p'}$
of ${\widetriangle{\mathcal U^{(1)}}}$
and
$$
U^+_{\frak p'} \cap \mathcal Z_{(1)} = U^{(1)}_{\frak p'}
\cap \mathcal Z_{(1)}.
$$
\end{enumerate}
\end{lem}
\begin{proof}
We put
$$
\frak P_{\ge \frak d}
=
\frak P \setminus \{\frak p' \in \frak P(\mathcal Z_{(1)})
\mid \dim U_{\frak p'} < \frak d\}.
$$
Then
${\widetriangle{\mathcal U_{\ge \frak d}}}
= (\frak P_{\ge \frak d},\{\mathcal U_{\frak p}
\mid \frak p \in \frak P_{\ge \frak d}\},\{\Phi_{\frak p\frak q}
\mid \frak p,\frak q \in \frak P_{\ge \frak d}, \frak p \ge \frak q\})$
is a good coordinate system of any compact subset of
$\mathcal Z_{(2)} \cap \bigcup_{\frak p' \in \frak P_{\ge \frak d}}
{\rm Im}(\psi_{\frak p'})$.
\par
We take $\frak U'(Z_{\frak p'})$ which
is an open neighborhood of $Z_{\frak p'}$ and is relatively compact in
$\frak U(Z_{\frak p'})
$.
We put
$$
Z'_1 = Z_1 \setminus \bigcup_{\frak p' \in \frak P(\mathcal Z_{(1)}),
\dim U_{\frak p'} < \frak d} \frak U'(Z_{\frak p'}).
$$
We observe that we are then in Situation \ref{situation101},
where ${\widetriangle{\mathcal U_{\ge \frak d}}}$
(resp. $Z'_1$) plays the role of
${\widetriangle{\mathcal U}}$
(resp. $Z_1$) in Situation \ref{situation101}.
We apply Proposition \ref{inductiveprop} to our situation
and obtain ${\widetriangle{\mathcal U
^{+ \prime}}}$.
Note that the union of the sets of Kuranishi charts of ${\widetriangle{\mathcal U
^{+ \prime}}}$
and $\{\mathcal U^{(1)}_{\frak p'} \mid \frak p' \in
\frak P(\mathcal Z_{(1)}),
\dim U_{\frak p'} < \frak d\}$
has most of the properties we need to prove.
The only point to take care of is that,
for $\frak p' \in
\frak P(\mathcal Z_{(1)})$ with $\dim U_{\frak p'} < \frak d$, neither the coordinate change $\Phi^+_{\frak p_0\frak p'}$
nor  $\Phi^+_{\frak p'\frak p_0}$ is  defined.

Let $\frak P^{+\prime}$ be the partial ordered set
appearing in $\widetriangle{\mathcal U^{+\prime}}$.
Then $\frak p_0 \in \frak P^{+\prime}$ and $\mathcal U^{+\prime}
_{\frak p_0}$ is a Kuranishi chart  that is an open subchart of
$\mathcal U_{\frak p_0}$.
To take care of the point mentioned above we shrink $\mathcal U^{+\prime}
_{\frak p_0}$ to $\mathcal U^{+}
_{\frak p_0}$
so that these two coordinates will not intersect,
as follows.
\begin{sublem}\label{sublem119}
There exists an open subset $U^{+}
_{\frak p_0}$ of $U^{+ \prime}
_{\frak p_0}$ such that the following holds.
\begin{enumerate}
\item
If
$\frak p' \in
\frak P(\mathcal Z_{(1)}),
\dim U_{\frak p'} < \frak d$ then
$$
\psi^+_{\frak p_0}(s_{\frak p_0}^{-1}(0) \cap U^{+}
_{\frak p_0})
\cap \frak U'(Z_{\frak p'})
= \emptyset
.
$$
\item
$
\psi_{\frak p_0}(s_{\frak p_0}^{-1}(0) \cap U^{+ \prime}
_{\frak p_0})
\cap \mathcal S_{\frak d}(X,\mathcal Z_{(1)}; \widehat{\mathcal U})
=
\psi_{\frak p_0}(s_{\frak p_0}^{-1}(0) \cap U^{+}
_{\frak p_0})
\cap \mathcal S_{\frak d}(X,\mathcal Z_{(1)};\widehat{\mathcal U}).
$
\end{enumerate}
\end{sublem}
\begin{proof}
By definition, we have
$Z_{\frak p'}
\cap
\mathcal S_{\frak d}(X,\mathcal Z_{(1)};\widehat{\mathcal U}) = \emptyset
$
for $\frak p' \in
\frak P(\mathcal Z_{(1)}),
\dim U_{\frak p'} < \frak d$.
Therefore we may choose $\frak U(Z_{\frak p'})$ so that
$$
\frak U(Z_{\frak p'})
\cap
\mathcal S_{\frak d}(X,\mathcal Z_{(1)};\widehat{\mathcal U}) = \emptyset
$$
for such $\frak p'$.
In fact,
$\mathcal S_{\frak d}(X,\mathcal Z_{(1)};\widehat{\mathcal U})$ is a
closed set.
Since $\frak U'(Z_{\frak p'})$ is relatively compact in
$\frak U(Z_{\frak p'})$, we have
$$
\overline{\frak U'(Z_{\frak p'})}
\cap
\mathcal S_{\frak d}(X,\mathcal Z_{(1)};\widehat{\mathcal U})
= \emptyset.
$$
Sublemma \ref{sublem119} is an immediate consequence of this fact.
\end{proof}
We now put
$
\mathcal U^{+}_{\frak p_0}
=
\mathcal U^{+\prime}_{\frak p_0}\vert_{U^{+}
_{\frak p_0}}
$.
For $\frak p' \in
\frak P(\mathcal Z_{(1)})$ with
$\dim U_{\frak p'} < \frak d$, we take an open subset
$
U^{+}_{\frak p'} \subset U^{(1)}_{\frak p'}
$
such that
$$
\frak U'(Z_{\frak p'})
=
\psi_{\frak p'}^{(1)}((s_{\frak p'}^{(1)})^{-1}(0)
\cap U^{+}
_{\frak p'}).
$$
Then we put $
\mathcal U^{+}_{\frak p'}
=
\mathcal U^{+\prime}_{\frak p'}\vert_{U^{+}
_{\frak p'}}
$.
For $\frak p \in \frak P^{+ \prime} \setminus \{\frak p_0\}$
we put
$
\mathcal U^{+}_{\frak p}
=
\mathcal U^{+\prime}_{\frak p}$.
We define a partial order $\le$ on $\frak P_+
= \frak P^{+\prime} \cup \frak P(\mathcal Z_{(1)})$
such that  $\le$ coincides with the
partial orders on $\frak P^{+\prime}$ and on
$\frak P(\mathcal Z_{(1)})$.
Moreover we define $\le$ so that
for $\frak p' \in \frak P^{+\prime}$ with
$\dim U_{\frak p'} < \frak d$,
neither $\frak p'\le \frak p_0$ nor $\frak p'\ge \frak p_0$.

\par
We can define coordinate change among them by
restricting of the coordinate change of either
$\widetriangle{{\mathcal U
^{+ \prime}}}$
or of ${\widetriangle{\mathcal U}}$.
Sublemma \ref{sublem119} (1)
implies that these two cases exhaust the
cases we need to define coordinate change.
The proof of Lemma \ref{lem1118} is complete.
\end{proof}
Using Lemma \ref{lem1118}
we discuss in the same way as the last step of the proof of
Theorem \ref{Them71restate}
to complete the proof of  Proposition \ref{prop7582752}.
\end{proof}
\begin{proof}[Proof of Lemma \ref{lem753753}]
Using Lemma \ref{lem1118}, we can prove it in the same way as
in Subsection \ref{subsec:movingmulsectionetc}.
\end{proof}

\section{Construction of CF-perturbations}
\label{sec:contfamilyconstr}

In this section, we give a thorough detail of the
proof of  existence of CF-perturbations
with respect to which a given weakly submersive map becomes
strongly submersive.
We also prove its relative version.

\subsection{Construction of CF-perturbations
on a single chart}
\label{subsec:confapersingle}

We first study the case of a single Kuranishi chart.
\begin{shitu}\label{situ121}
$\mathcal U = (U,E,s,\psi)$ is a Kuranishi chart of $X$
and $f : U \to M$ is a smooth map.
$\blacksquare$
\end{shitu}

\begin{shitu}\label{situ122}
In Situation \ref{situ121}, we assume
 that $g : N \to M$ is a smooth map between manifolds and
that $f$ is  transversal to $g$.
$\blacksquare$
\end{shitu}
The main result of Subsection \ref{situ121} is
Proposition \ref{prop123123} below.
We recall the following well-known definition.
\begin{defn}
A sheaf (of sets) $\mathscr F$ on a topological space $V$ is said to be
{\it soft} \index{sheaf $\mathscr{S}$ ! soft} if the restriction map
$$
\mathscr F(V) \to \mathscr F(K)
$$
is surjective for any closed subset $K$ of $V$.
(We note
$\mathscr F(K) = \varinjlim_{W \supset K, \text{open}} \mathscr F(W)$.)
\end{defn}
\begin{prop}\label{prop123123}
Suppose we are in Situation \ref{situ121}.
\begin{enumerate}
\item The sheaf $\mathscr S$ in Proposition \ref{prop721} is soft.
\item
The sheaf $\mathscr S_{\pitchfork 0}$ in
Lemma-Definition \ref{strosubsemiloc} is soft.
\item
Suppose $f$ is a submersion. Then,
the sheaf $\mathscr S_{f \pitchfork}$ in
Lemma-Definition \ref{strosubsemiloc} is soft.
\item
In Situation \ref{situ122}
the sheaf $\mathscr S_{f \pitchfork g}$ in
Lemma-Definition \ref{strosubsemiloc} is soft.
\end{enumerate}
\end{prop}
The rest of this subsection will be occupied by the proof of this proposition.
\begin{proof}
We first prove (1).
We use partition of
unity to glue sections of $\mathscr S$.
Note our sheaf $\mathscr S$ is a sheaf of sets. Nevertheless we can
apply partition of
unity, as we will discuss below.
\begin{shitu}\label{situ12999}
Let $A$ be a subset of $U$ and
let $\{U_{\frak r} \mid \frak r \in \frak R\}$ be a locally finite open cover
of a subset $A$ in $U$ and $\{\chi_{\frak r}\}$ a smooth partition
of unity subordinate to this covering.
In other words, $\chi_{\frak r}  : U \to [0,1]$ is a smooth function of $U$ which has
compact support in $U_{\frak r}$, and
$$
\sum_{\frak r \in \frak R} \chi_{\frak r}(x) =1
$$
for $x \in A$.
\par
We assume that an element $\frak S_{\frak r} \in \mathscr S(A)$ is given for each $\frak r\in \frak R$.
$\blacksquare$
\end{shitu}
Below we will define the sum
\begin{equation}
\sum_{\frak r} \chi_{\frak r} \frak S_{\frak r} \in \mathscr S
(A).
\end{equation}
For $x \in A$, let $\frak V_x = (V_x,\Gamma_x,E_x,\phi_x,\widehat\phi_x)$
be an orbifold chart of $(U,\mathcal E)$
at $x$.
We may assume that, for each $\frak r$ with $x \in U_{\frak r}$,
we are given a representative $\mathcal S_{\frak r}$ of  $\frak S_{\frak r}$
on a neighborhood of $x$.
It consists of $\frak V_{\frak r} = (V_{\frak r},\Gamma_{\frak r},E_{\frak r},\psi_{\frak r},
\widehat{\psi}_{\frak r})$ and
$(W_{\frak r},\omega_{\frak r},\{{\frak s}_{\frak r}^{\epsilon} \mid \epsilon\})$
where
$\frak V_{\frak r} = (V_{\frak r},\Gamma_{\frak r},E_{\frak r},\psi_{\frak r},
\widehat{\psi}_{\frak r})$ is an orbifold chart of $(U,\mathcal E)$ at $x$
and $(W_{\frak r},\omega_{\frak r},\{{\frak s}_{\frak r}^{\epsilon} \mid \epsilon\})$ is
as in Definition \ref{defn73ss}.
\par
We put
\begin{equation}
\frak R(x) = \{\frak r \in \frak R \mid x \in {\rm Supp}(\chi_{\frak r})\}.
\end{equation}
\par
By shrinking $V_x$ if necessary we may assume
$\psi_x(V_x) \subset U_{\frak r}$ for each $\frak r \in \frak R(x)$ and
$\chi_{\frak r} \equiv 0$ on $\psi_x(V_x)$ for each $\frak r \notin \frak R(x)$.
\par
Furthermore we may choose $U_{x}$   so that there
exist
\begin{equation}\label{coorchange124}
\aligned
&h_{\frak r x} : \Gamma_{x} \to \Gamma_{\frak r}, \\
&\widetilde{\varphi}_{\frak r x} : V_{x} \to V_{\frak r}, \\
&\breve\varphi_{\frak r x} : V_{x} \times E_{x} \to  E_{\frak r}
\endaligned
\end{equation}
as in Situation \ref{opensuborbifoldchart},
for each $\frak r \in \frak R(x)$.
(See Lemma \ref{lem2622}.)

\begin{defn}\label{defn123}
We put
$$
W_{x} = \prod_{\frak r \in \frak R(x)} W_{\frak r},
\quad
\omega_{x} = \prod_{\frak r \in \frak R(x)} \omega_{\frak r}.
$$
We define $\frak s^{\epsilon}_{x}
: V_{x} \times W_{x} \to E_{x}$
by the following formula:
\begin{equation}\label{formula12555}
\frak s^{\epsilon}_{x}(y,(\xi_{\frak r})_{\frak r \in \frak R(x)})
=
s_{x}(y) + \sum_{\frak r \in  \frak R(x)}
\chi_{\frak r}(\psi_x(y))
g_{\frak r,y}^{-1}(\frak s_{\frak r}^{\epsilon}(\widetilde{\varphi}_{\frak r x}(y),\xi_{\frak r})
- s_{\frak r}(\widetilde{\varphi}_{\frak r x}(y)).
\end{equation}
Here $s_{x} : V_{x} \to E_{x}$ and $s_{\frak r} : V_{\frak r} \to E_{\frak r}$
are the local expressions of the Kuranishi map
(Definition \ref{defnlocex}.)
and
$g_{\frak r,y} : E_{x} \to E_{\frak r}$
is defined by
$
\breve{\varphi}_{\frak r x}(y,\xi) = g_{\frak r,y}(\xi).
$
\par
We put
$\mathcal S_x = (W_{x},\omega_{x},\{\frak s^{\epsilon}_{x}\})$.
\end{defn}
\begin{lem}\label{lem1241}
\begin{enumerate}
\item
$\mathcal S_{x}$ is a CF-perturbation
of $\mathcal U$ on $\frak V_{x}$.
\item
The germ $[\mathcal S_{x}] \in \mathscr S_x$ represented by $\mathcal S_{x}$
depends only on $\{\chi_{\frak r}\}$, $\{\mathcal S_{\frak r}\}$,
$x$ and is independent of the choices of $\frak V_x$,
the coordinate changes (\ref{coorchange124}),
and the representatives of $\{\mathcal S_{\frak r}\}$.
\item
$x \mapsto [\mathcal S_{x}] \in \mathscr S_x$ defines a
(global) section of the
sheaf $\mathscr S$.
\end{enumerate}
\end{lem}
\begin{proof}
Statement (1) is an immediate consequence of the construction.
\par
We prove Statement (2).
We first prove independence of the coordinate changes (\ref{coorchange124}).
Let $(h'_{\frak r x},\widetilde{\varphi}'_{\frak r x},\breve\varphi'_{\frak  r x})$ be an alternative choice.
Then there exists $\gamma_{\frak r} \in \Gamma_{\frak r}$
such that
$h'_{\frak r x}
=\gamma_{\frak r}h_{\frak r x}\gamma_{\frak r}^{-1}
$,
$
\widetilde{\varphi}'_{\frak r x}
= \gamma_{\frak r}\widetilde{\varphi}_{\frak r x}
$,
$
\breve\varphi'_{\frak r x}
= \gamma_{\frak r}\breve\varphi_{\frak r x}.
$
The third equality implies $g'_{\frak r,y} = \gamma_{\frak r} g_{\frak r,y}$
by Lemma \ref{lem2715}.
Let $\frak s^{\epsilon \prime}_{x}$ be obtained from this
alternative choice. Then we have
\begin{equation}\label{form126126}
\aligned
&\frak s^{\epsilon \prime}_{x}(y,\xi) \\
&=
s_{x}(y) + \sum_{\frak r \in  \frak R(x)}
\chi_{\frak r}(\psi_x(y)) (g'_{\frak r,y})^{-1}
(\frak s_{\frak r}^{\epsilon}(\widetilde{\varphi}'_{\frak r x}(y),\xi_{\frak r})
- s_{\frak r}(\widetilde{\varphi}'_{\frak r x}(y))) \\
&=
s_{x}(y) + \sum_{\frak r \in  \frak R(x)}
\chi_{\frak r}(\psi_x(y)) g_{\frak r,y}^{-1}
(\frak s_{\frak r}^{\epsilon}(\widetilde{\varphi}_{\frak r x}(y),\gamma_{\frak r}^{-1}\xi_{\frak r})
- s_{\frak r}(\widetilde{\varphi}_{\frak r x}(y))).
\endaligned
\end{equation}
Here we use $\Gamma_{\frak r}$ equivariance of
$\frak s_{\frak r}^{\epsilon}$ and of
$s_{\frak r}$.
\par
We define a $\Gamma_{x}$ action on $W_{\frak r}$ by $\mu\cdot\xi
= h_{\frak r x}(\mu)\xi$. We write $W_{\frak r}$ with this action by $W_{\frak r}^{h_{\frak r x}}$.
The notation $W_{\frak r}^{h'_{\frak r x}}$ is defined in a similar way.
Its product  in Definition \ref{defn123} is denoted by $W_x^h$ and $W_x^{h'}$, respectively.
\par
Then $\xi_{\frak r} \mapsto \gamma_{\frak r}^{-1}\xi_{\frak r}$
(resp. $(\xi_{\frak r}) \mapsto (\gamma_{\frak r}^{-1}\xi_{\frak r})$)
is a $\Gamma_{x}$ equivariant linear map $: W_{\frak r}^{h_{\frak r x}}
\to W_{\frak r}^{h'_{\frak r x}}$
(resp. $W_x^h \to W_x^{h'}$).
Therefore (\ref{form126126}) implies that the equivalence
class $[\mathcal S_{x}]$ is independent of the choices of
the coordinate changes (\ref{coorchange124}).
\par
Secondly we prove independence of the representative of $\mathcal S_{\frak r}$.
We consider one of $\frak r_0 \in \frak R(x)$ and
take an alternative choice $\mathcal S'_{\frak r_0}$ of $\mathcal S_{\frak r_0}$.
It suffices to  consider the case when other $\mathcal S_{\frak r}$'s
for $\frak r \ne \frak r_0$ are the same for both.
We may also assume that $\mathcal S'_{\frak r_0}$ is also
a projection
of  $\mathcal S_{\frak r_0}$.
Then it is immediate from definition that
$\mathcal S'_x$  obtained by using $\mathcal S'_{\frak r_0}$
is a projection of  $\mathcal S_x$ which is obtained by using $\mathcal S_{\frak r_0}$.
We have thus proved the
independence of the representative of $\mathcal S_{\frak r}$.
\par
Thirdly we prove independence of the orbifold chart
$\frak V_x = (V_x,\Gamma_x,E_x,\phi_x,\widehat\phi_x)$.
Let
$\frak V'_x = (V'_x,\Gamma'_x,E'_x,\phi'_x,\widehat\phi'_x)$
and suppose we obtain  $\mathcal S'_x$ when we use
$\frak V'_x$.
\par
By shrinking $V'_x$ if necessary we may assume that there exists
a coordinate change $(h_x,\widetilde{\varphi}_x,\breve\varphi_x)$ from the
chart $\frak V'_x$ to $\frak V_x$.
Let
$(h_{\frak r x},\widetilde{\varphi}_{\frak r x},\breve\varphi_{\frak r x})$
be the coordinate change as in (\ref{coorchange124}).
Then by putting
$$
h'_{\frak r x} = h_{\frak r x}\circ h_x,
\quad
\widetilde{\varphi}'_{\frak r x} = \widetilde{\varphi}_{\frak r x}\circ \widetilde{\varphi}_x,
\quad
\breve\varphi'_{\frak r x} = \breve\varphi_{\frak r x}\circ
\breve\varphi_x,
$$
$(h'_{\frak r x},\widetilde{\varphi}'_{\frak r x},\breve\varphi'_{\frak r x})$
becomes a coordinate change from $\frak V'_x$ to
$\frak V_{\frak r}$ as in (\ref{coorchange124}).
\par
Then
\begin{equation}\label{form12612622}
\aligned
&\frak s^{\epsilon \prime}_{x}(y,\xi)\\
&=
s_{x}(y) + \sum_{\frak r \in  \frak R(x)}
\chi_{\frak r}(\psi'_x(y)) (g'_{\frak r,y})^{-1}
(\frak s_{\frak r}^{\epsilon}(\widetilde{\varphi}'_{\frak r x}(y),\xi_{\frak r})
- s_{\frak r}(\widetilde{\varphi}'_{\frak r x}(y))) \\
&=
s_{x}(y) + \sum_{\frak r \in  \frak R(x)}
\chi_{\frak r}(\psi_x(\widetilde{\varphi}_x(y))) g_{\frak r,y}^{-1}
(\frak s_{\frak r}^{\epsilon}(\widetilde{\varphi}_{\frak r x}(\widetilde{\varphi}_x(y)),\gamma_{\frak r}^{-1}\xi_{\frak r})
- s_{\frak r}(\widetilde{\varphi}_{\frak r x}(\widetilde{\varphi}_x(y))))
\\
&= \frak s^{\epsilon}_{x}(\widetilde{\varphi}_x(y),\xi).
\endaligned
\nonumber
\end{equation}
This implies the required independence of the coordinate $\frak V_x$.
The proof of Statement (2) is complete.
\par
We now prove Statement (3).
Let $\mathcal S_x = (W_{x},\omega_{x},\{\frak s^{\epsilon}_{x}\})$
as above.
Suppose $y \in s_{x}^{-1}(0) \cap U_x$
and $y = \phi_x(\tilde y)$.
We denote $\Gamma_{\tilde y} = \{\gamma \in \Gamma_x \mid
\gamma \tilde y  = \tilde y\}$ and take a $\Gamma_{\tilde y}$ invariant neighborhood $V_y$ of $\tilde y$.
Then
$\frak V_y = (V_y,\Gamma_y,E_x,\phi_x\vert_{V_y},\widehat\phi_x\vert_{V_y})$
is an orbifold chart of $(U,E)$ at $y$.
It is easy to see that $\mathcal S_x\vert_{V_y} = (W_{x},\omega_{x},\{\frak s^{\epsilon}_{x}\vert_{V_y \times E_x}\})$ is a CF-perturbation on $\frak V_y$.
\begin{sublem}\label{sublem128}
$\mathcal S_x\vert_{V_y}$ is equivalent to $\mathcal S_y$  in the sense of Definition \ref{conmultiequiv11}.
\end{sublem}
\begin{proof}
We consider $\frak R(y) = \{\frak r \in \frak R \mid
y \in {\rm Supp}(\chi_{\frak r})\}$.
Since we chose $U_x \subseteq U_{\frak r}$ such that
$U_x \cap \text{\rm Supp} (\chi_{\frak r}) = \emptyset$ for
${\frak r} \notin {\frak R}(x)$, we have $\frak R(y)
\subseteq \frak R(x)$.
Therefore there exists an obvious projection
$$
\pi : \widehat W_x = \prod_{\frak r\in \frak R(x)} W_{\frak r}
\to \widehat W_y = \prod_{\frak r\in \frak R(y)} W_{\frak r}.
$$
It is easy to see that
$\pi!(\omega_x) = \omega_y$.
\par
We may choose $V_y$ so small  that for $z \in V_y$ and $\frak r \in \frak R(x)
\setminus \frak R(y)$ we have $\chi_{\frak r}(z) = 0$.
Therefore by definition
$$
\frak s^{\epsilon}_{x}(\tilde\varphi_x(z),\xi)
=
\frak s^{\epsilon}_{y}(\tilde\varphi_y(z),\pi(\xi))
$$
for $z \in V_y$.
Thus $\mathcal S_y$ is a projection of $\mathcal S_x\vert_{V_y}$.
\end{proof}
Statement (3) follows from Sublemma \ref{sublem128} and Lemma \ref{lem723}.
\end{proof}
\begin{defn}\label{defn1213}
We denote by
$$
\sum_{\frak r} \chi_{\frak r} \frak S_{\frak r}
$$
the element $x \mapsto [\mathcal S_{x}] \in \mathscr S_x$
of $\mathscr S(A)$ obtained by Lemma \ref{lem1241}.
\end{defn}
\begin{rem}
Suppose $\frak S \in \mathscr S(U)$ and
$\{U_{\frak r}\mid \frak r \in \frak R\}$ is a locally
finite cover of $U$. We can define an element of $\mathscr S(U)$ by
$$
\sum_{\frak r} \chi_{\frak r} \frak S\vert_{U_{\frak r}}
$$
as above.
(Here $\frak S\vert_{U_{\frak r}} \in \mathscr S(U_{\frak r})$
is the restriction of $\frak S$.) However this section is in general {\it different} from the originally given
$\frak S \in \mathscr S(U)$.
\end{rem}
Proposition \ref{prop123123} (1)  follows easily from Definition \ref{defn1213} and the results
we proved above.
(See also the end of this subsection where the proof of Proposition \ref{prop123123} (2)(3)(4) are completed.)
\par
We next prove Proposition \ref{prop123123} (2)(3)(4).
We begin with the next definition.
\begin{defn}\label{strongtransvers}
Suppose we are in Situation \ref{situ121}.
Let $\frak S_x \in \mathscr S_x$ be a germ of the sheaf $\mathscr S$ at $x \in U$.
\par
We say $\frak S_x$ is {\it strongly transversal} \index{sheaf $\mathscr{S}$ ! strongly transversal germ}\index{CF-perturbation !strongly transversal germ}
if its representative $(W_x,\omega_x,\{\frak s^{\epsilon}_x\})$
(which is defined on the orbifold chart $\frak V_x = (V_x,\Gamma_x,E_x,\phi_x,\widehat\phi_x)$ (Definition \ref{defn61}(3)))
has the following properties.
\begin{enumerate}
\item
For all sufficiently small $\epsilon > 0$, the map $\frak s^{\epsilon}_x : V_x \times W_x \to E_x$
is transversal to $c \in E_x$ on a neighborhood of $\{o_x\} \times {\rm Supp}(\omega_x)$ for any $c \in E_x$.
(Here $o_x \in V_x$ is the point such that $\phi_x(o_x) = x$.)
\item
For $\xi \in {\rm Supp}(\omega_x)$ and $c= \frak s^{\epsilon}_x(o_x,\xi)$ the projection
$$
T_{(o_x,\xi)} (\frak s^{\epsilon}_x)^{-1}(c) \to T_{o_x} V_x
$$
is surjective.
\end{enumerate}
We write  $(\mathscr S_{\pitchfork\pitchfork 0})_x$ the set of all germs of the sheaf $\mathscr S$ at $x \in U$
that is strongly transversal.
\par
It is easy to see that there exists a subsheaf of $\mathscr S$ whose stalk at $x$ is $(\mathscr S_{\pitchfork\pitchfork 0})_x$.
We denote this sheaf by $\mathscr S_{\pitchfork\pitchfork 0}$.
\end{defn}
\begin{rem}
It is easy to see that the above properties (1)(2) are independent of the choice
of the representative $(W_x,\omega_x,\{\frak s^{\epsilon}_x\})$ and of the
orbifold chart $(V_x,\Gamma_x,E_x,\phi_x,\widehat\phi_x)$
but depend only on $\frak S_x$.
\end{rem}
\begin{lem}\label{lemma12770}
In Situation \ref{situ121}, the set $(\mathscr S_{\pitchfork\pitchfork 0})_x$ is nonempty.
\end{lem}
\begin{proof}
Let
$\frak V_x = (V_x,\Gamma_x,E_x,\phi_x,\widehat\phi_x)$
be an orbifold chart of $(U,\mathcal E)$
at $x$.
We put $\widehat W_x = E_x$ and $W_x$ is a sufficiently
small $\Gamma_x$ invariant neighborhood of $0$ in  $\widehat W_x$ and
$\omega_x$ is a $\Gamma_x$ invariant differential form of compact
support  on $W_x$ of degree $\dim W_x$ with $\int \omega_x =1$.
We define
\begin{equation}
\frak s^{\epsilon}(x,\xi) = s(x) + \epsilon \xi.
\end{equation}
It is easy to see that
$(W_x,\omega_x,\{{\frak s}_x^{\epsilon} \mid \epsilon\})$ is a
CF-perturbation on $\frak V_x$.
Moreover it is easy to show that the projection
$
(\frak s^{\epsilon})^{-1}(c) \to V_x
$
is a submersion for any $c$.
 Lemma \ref{lemma12770} follows.
\end{proof}
\begin{lem}\label{lem1244440}
Suppose we are in Situation \ref{situ121}.
\begin{enumerate}
\item
$(\mathscr S_{\pitchfork\pitchfork 0})_x \subseteq (\mathscr S_{\pitchfork})_x$.
\item
If $f$ is a submersion at $x$, then,
$(\mathscr S_{\pitchfork\pitchfork 0})_x \subseteq (\mathscr S_{f \pitchfork})_x$.
\item
In Situation \ref{situ122}, we have $(\mathscr S_{\pitchfork\pitchfork 0})_x \subseteq (\mathscr S_{f \pitchfork g})_x$.
\end{enumerate}
\end{lem}
This is immediate from the definition.
\begin{cor}\label{lem124444}
Suppose we are in Situation \ref{situ121}.
\begin{enumerate}
\item
The stalk $(\mathscr S_{\pitchfork 0})_x$ is
nonempty for any $x \in U$.
\item
Suppose $f$ is a submersion at $x$. Then,   the
stalk $(\mathscr S_{f \pitchfork})_x$ is
nonempty.
\item
In Situation \ref{situ122}, the stalk $(\mathscr S_{f \pitchfork g})_x$ is
nonempty.
\end{enumerate}
\end{cor}
This is an immediate consequence of Lemmata \ref{lemma12770} and \ref{lem1244440}.
\par
To prove Proposition \ref{prop123123} (2)(3)(4), we need one more result
(Proposition \ref{lem1211515} below.)
\begin{shitu}\label{situ121515}
Suppose we are in Situation \ref{situ12999}.
We put $\frak R = \{\frak r_0\} \cup \frak R'$ and assume that
for $\frak r \in \frak R'$ the section $\frak S_{\frak r} \in \mathscr S(U_{\frak r})$
is strongly transversal. $\blacksquare$
\end{shitu}
\begin{prop}\label{lem1211515}
In Situation \ref{situ121515},
the following holds for $i=1,2,3,4$.
\par
If $\frak S_{\frak r_0}$ has Property
$(i)$ below
and $\frak S_{\frak r} \in \mathscr S_{\pitchfork\pitchfork 0}(U_{\frak r})$
for $\frak r \in \frak R'$,
then the sum
$$
\frak S = \sum_{\frak r\in \frak R} \chi_{\frak r} \frak S_{\frak r}
$$
has the same property $(i)$.
\begin{enumerate}
\item
$\frak S_{\frak r_0} \in \mathscr S_{\pitchfork\pitchfork 0}(U_{\frak r_0})$.
\item
$\frak S_{\frak r_0} \in \mathscr S_{\pitchfork 0}(U_{\frak r_0})$.
\item
$\frak S_{\frak r_0} \in \mathscr S_{f \pitchfork}(U_{\frak r_0})$.
\item
We are in Situation \ref{situ122} and
$\frak S_{\frak r_0} \in \mathscr S_{f \pitchfork g}(U_{\frak r_0})$.
\end{enumerate}
\end{prop}
\begin{proof}
To prove Proposition \ref{lem1211515} we rewrite the strong transversality
as follows.
Let $\mathcal S_x = (W_x,\omega_x,\{\frak s^{\epsilon}_x\})$ be a representative
of a germ $\mathscr S_x$ which is defined on an orbifold chart $\frak V_x = (V_x,\Gamma_x,E_x,\phi_x,\widehat\phi_x)$
of $(U,\mathcal E)$.
\begin{lem}\label{lem1217}
$\mathcal S_x$ is strongly transversal if and only if
the derivative
$$
\nabla^W_{(x,\xi)} \frak s^{\epsilon}_x : T_{\xi}W_x \to T_cE_x
$$
in $W_x$ direction is surjective for all $\xi $ in the support of $\omega_x$.
Here $c = \frak s^{\epsilon}_x(o_x,\xi)$.
\end{lem}
\begin{proof}
We consider the following commutative diagram where all the
horizontal and vertical lines are exact.
$$
\begin{CD}
&&
&&
0
&&
0
\\
&&
&& @VVV  @VVV
\\
&&
&&
T_{o_x,\xi} (\frak s_x^{\epsilon})^{-1}(c)
@>>>
T_{o_x}V_x
\\
&&
&& @ VVV @VVV
\\
0 @>>>T_{\xi}W_x
@>>>T_{(o_x,\xi)}(V_x \times W_x)
@>>>
T_{o_x}V_x
@>>> 0
\\
&& @ VVV @VVV @VVV\\
0 @>>>
T_c E_x @>>>T_c E_x
@>>>0
\\
&& && @VVV
\\
&&&& 0
\end{CD}
$$
The required strong transversality is nothing but the surjectivity of the
second horizontal map
$: T_{o_x,\xi} (\frak s_x^{\epsilon})^{-1}(c) \to T_{o_x}V_x$
and the map $
\nabla^W_{(x,\xi)} \frak s^{\epsilon}_x : T_{\xi}W_x \to T_cE_x
$
is the second vertical map.
The equivalence of the surjectivities of them
is a consequence of simple diagram chase.
\end{proof}
The next lemma is a half of the proof of Proposition  \ref{lem1211515}.
\begin{lem}\label{lem121818}
Suppose we are in Situation \ref{situ121515} and
$\chi_{\frak r_1}(x) \ne 0$ for some $\frak r_1 \in \frak R'$.
Then the germ $\frak S_x$ of $\frak S$ at $x$ is strongly transversal.
\end{lem}
\begin{proof}
A representative of $\frak S_x$ is $(W_x,\omega_x,\frak s_x^{\epsilon})$
where
\begin{equation}
\frak s^{\epsilon}_{x}(y,(\xi_{\frak r})_{\frak r \in \frak R(x)})
=
s_{x}(y) + \sum_{\frak r \in  \frak R(x)}
\chi_{\frak r}(\psi_x(y))
g_{\frak r,y}^{-1}(\frak s_{\frak r}^{\epsilon}(\tilde{\varphi}_{\frak r x}(y),\xi_{\frak r})
- s_{\frak r}(\tilde{\varphi}_{\frak r x}(y)).
\end{equation}
Here $\frak R(x) \subseteq \frak R$ and $\frak r_1 \in \frak R(x)$.
The derivative of $\frak s^{\epsilon}_{x}$ in $W_{\frak r_1}$ direction is
\begin{equation}\label{form12118}
\chi_{\frak r_1}(\psi_x(y))
g_{\frak r_1,y}^{-1}(\nabla^{W_{\frak r_1}}\frak s_{\frak r_1}^{\epsilon}\vert_{(\tilde{\varphi}_{\frak r x}(y),\xi_{\frak r_1})}).
\end{equation}
By Lemma \ref{lem1217} the derivative $\nabla^{W_{\frak r_1}}\frak s_{\frak r_1}^{\epsilon}$ is surjective
(to $T_cE_{x}$).
Therefore (\ref{form12118}) is surjective.
Therefore by Lemma \ref{lem1217} $\frak S_x$ is strongly transversal.
\end{proof}
Now we are ready to complete  the proof of Proposition \ref{lem1211515}.
By Lemma \ref{lem121818} the  germ $\frak S_x$ has the property claimed in Proposition \ref{lem1211515}
unless $\chi_{\frak r}(x) = 0$ for all $\frak r \in \frak R'$.
We may also assume that
$\frak s^{\epsilon}(x,(\xi_{\frak r})_{\frak r \in \{\frak r_0\} \cup \frak R(x)}) = 0$.
We consider such a point $x$. A
representative of $\frak S_x$ is $(W_x,\omega_x,\frak s_x^{\epsilon})$
where
\begin{equation}
\aligned
\frak s^{\epsilon}_{x}(y,&(\xi_{\frak r})_{\frak r \in \{\frak r_0\} \cup \frak R(x)}) \\
=
s_{x}(y) & +
\chi_{\frak r_0}(\psi_x(y))
g_{\frak r_0,y}^{-1}(\frak s_{\frak r_0}^{\epsilon}(\tilde{\varphi}_{\frak r_0 x}(y),\xi_{\frak r})
- s_{\frak r}(\tilde{\varphi}_{\frak r x}(y))\\
&+ \sum_{\frak r \in  \frak R(x)}
\chi_{\frak r}(\psi_x(y))
g_{\frak r,y}^{-1}(\frak s_{\frak r}^{\epsilon}(\tilde{\varphi}_{\frak r x}(y),\xi_{\frak r})
- s_{\frak r}(\tilde{\varphi}_{\frak r x}(y)).
\endaligned
\end{equation}
Here $\frak R(x) \subseteq \frak R'$.
We remark that $\chi_{\frak r_0}(x) = 1$ and takes maximum there.
(Note $\chi_{\frak r}$ is a smooth map to $[0,1]$.)
Therefore the first derivative at $x$ of $\chi_{\frak r_0}$ is zero.
In a similar way  we can show that the first derivatives at $x$ of
$\chi_{\frak r}$ are all zero.
\par
We also remark that $\frak s^{\epsilon}_{\frak r_0}(x,\xi_{\frak r_0}) = 0$.
Therefore
\begin{equation}\label{form121010}
T_{(x,\xi_{\frak r_0})}(\frak s^{\epsilon}_{\frak r_0})^{-1}(0)
\times \prod_{\frak r \in \frak  R(x)} T_{\xi_{\frak r}}W_{\frak r}
\subseteq
T_{(x,(\xi_{\frak r})_{\frak r \in \{\frak r_0\} \cup \frak R(x)})}
(\frak s^{\epsilon}_{x})^{-1}(0).
\end{equation}
(\ref{form121010}) implies that
if $\frak S_{\frak r_0}$ has Property
$(i)$ at $x$ then $
\frak S = \sum_{\frak r} \chi_{\frak r} \frak S_{\frak r}
$
has the same property $(i)$ at $x$,
where $\chi_{\frak r_0}(x) = 1$.
\par
This fact together with Lemmata \ref{lem124444} and \ref{lem121818}
imply Proposition  \ref{lem1211515} .
\end{proof}
We are now in the position to complete the proof of Proposition  \ref{prop123123} (2)(3)(4).
Let $K \subset U$ be a closed subset and $\frak S_K \in \mathscr S(K)$.
By definition there exists an open neighborhood $U_{\frak r_0}$ of $K$
such that  $\frak S_K$ is a restriction of $\frak S_{\frak r_0}
\in \mathscr S(U_{\frak r_0})$.
We take an index set $\frak R'$ and an open covering
\begin{equation}\label{form121111}
U = U_{\frak r_0} \cup \bigcup_{\frak r \in \frak R'} U_{\frak r}
\end{equation}
with the following properties.
\begin{enumerate}
\item[(a)]
The covering (\ref{form121111}) is locally finite.
\item[(b)]
$\mathscr S_{\pitchfork\pitchfork 0}(U_{\frak r}) \ne \emptyset$ for $\frak r \in \frak R'$.
\item[(c)]
$K \cap U_{\frak r} = \emptyset$ for $\frak r \in \frak R'$.
\end{enumerate}
Existence of such covering is a consequence of paracompactness of $U$ and
Lemma \ref{lemma12770}.
\par
Let $\frak S_{\frak r} \in \mathscr S_{\pitchfork\pitchfork 0}(U_{\frak r})$
and $\chi_{\frak r}$ a partition of unity subordinate to the
covering (\ref{form121111}).
We put
$$
\frak S = \sum_{\frak r \in \{\frak r_0\} \cup \frak R'} \chi_{\frak r}\frak S_{\frak r}.
$$
Property (c) implies that $\frak S$ restricts to $\frak S_K$.
We have thus proved the softness of $\mathscr S$.
\par
To prove the softness of $\mathscr S_{\pitchfork 0}$,
we may assume $\frak S_{\frak r_0} \in
\mathscr S_{\pitchfork 0}(U_{\frak r_0})$.
Then by Proposition \ref{lem1211515},
$\frak S \in  \mathscr S_{\pitchfork 0}(U)$.
\par
The proof of softness of $\mathscr S_{f \pitchfork }$ and of $\mathscr S_{f \pitchfork g}$
is the same.
The proof of Proposition  \ref{prop123123} is complete.
\end{proof}
\begin{rem}
The same argument also proves softness of $\mathscr S_{\pitchfork\pitchfork 0}$.
\end{rem}

\subsection{Embedding of Kuranishi charts and
extension of CF-perturbations}
\label{subsec:extembandcfp}
In this subsection, we study the case of good coordinate system
and in the next we prove Theorem \ref{existperturbcont},
and its relative version Proposition \ref{existperturbcontrel}.
\par
For that purpose we need to
study  other kinds of extension.
Namely we will study extension of a CF-perturbation
defined on an embedded orbifold to its neighborhood.
(Then, by Proposition  \ref{prop123123}, we can extend a CF-perturbation defined on this neighborhood.)
We consider the following Situation \ref{sit121111}.
\begin{shitu}\label{sit121111}
Let $\mathcal U_i = (U_i,\mathcal E_i,s_i,\psi_i)$ $(i=1,2)$ be
Kuranishi charts of $X$,
$\Phi_{21} = (\varphi_{21},\widehat{\varphi}_{21}) : \mathcal U_1 \to \mathcal U_2$ an embedding
of Kuranishi charts
and $K$ a closed subset of $X$ such that
$K \subset \psi_1(s_1^{-1}(0)) \subset \psi_2(s_2^{-1}(0))$.
\par
Hereafter in this subsection, we regard $K$ also as a subset of
$U_1$ or $U_2$ via the parameterizations $\psi_1$, $\psi_2$ respectively.
\par
Let $f_2 : U_2 \to M$ be a strongly smooth map and we put
$f_1 = f_2 \circ \varphi_{21}$.
$\blacksquare$
\end{shitu}
We now consider the following commutative diagram.
\begin{equation}
\begin{CD}
\varphi_{21}^{\star}\mathscr S_{\sharp}^{\mathcal U^1\triangleright \mathcal U^2}(U_1)
@>{\frak i_{KU_1}}>>
\varphi_{21}^{\star}\mathscr S_{\sharp}^{\mathcal U^1\triangleright \mathcal U^2}(K)
\\
@V{\Phi^*_{21}}VV @V{\Phi^*_{21}}VV
\\
\mathscr S_{\sharp}^{\mathcal U^1}(U_1)
@>{\frak i_{KU_1}}>>
\mathscr S_{\sharp}^{\mathcal U^1}(K)
\end{CD}
\end{equation}
Here $\sharp$ stands for any of $\pitchfork 0$, $f \pitchfork$,
$f \pitchfork g$, or $\pitchfork\pitchfork 0$.
Recall
$\Phi^*_{21}$ denotes the restriction map
and $\varphi_{21}^{\star}$ stands for pullback sheaf.
(See Definition \ref{deflem743} (4).)
\begin{prop}\label{prop1221}
The following holds for ${\sharp}=\pitchfork 0$, $f \pitchfork$,
$f \pitchfork g$, or $\pitchfork\pitchfork 0$.
\par
Let $\frak S^1 \in \mathscr S_{\sharp}^{\mathcal U^1}(U_1)$,
$\frak S^{2,K} \in \varphi_{21}^{\star}\mathscr S_{\sharp}^{\mathcal U^1\triangleright \mathcal U^2}(K)$
such that
\begin{equation}
\frak i_{KU_1}(\frak S^1) = \Phi^*_{21}(\frak S^{2,K}).
\end{equation}
Then for any compact subset $Z \subset U_1$ containing $K$,
there exists
$\frak S^{2} \in \varphi_{21}^{\star}\mathscr S_{\sharp}^{\mathcal U^1\triangleright \mathcal U^2}(Z)$ such that
\begin{equation}\label{form1214}
\frak i_{KZ}(\frak S^2) = \frak S^{2,K},
\qquad
\Phi^*_{21}(\frak S^{2}) = \frak S^{1}.
\end{equation}
\end{prop}
The proof occupies the rest of this subsection.
\begin{proof}
For the proof, we will use the notion of bundle extension data.
To give its definition, we first introduce the definition of a tubular neighborhood of
an orbifold embedding.
\begin{defn}\label{lem123000}
Let $X \to Y$ be an embedding of orbifolds, $Z\subset X$  a
compact subset and $U$ be an open neighborhood of $K$ in $Y$. We say that $\pi : U \to X$ is
{\it diffeomorphic to the
projection of normal bundle}
\index{diffeomorphic to the
projection of normal bundle} if the following holds.
\par
Let ${\rm pr} : N_XY \to X$ be the normal bundle. Then there exists
a neighborhood $U'$  of $Z \subset X \subset N_XY$ and a diffeomorphism
of orbifolds $h : U' \to U$ such that $\pi\circ h = {\rm pr}$.
We also require that $h(x) = x$ for $x$ in a neighborhood of $Z$ in  $X$.
\end{defn}
Using this, we give the definition of a bundle extension datum.
\begin{defn}\label{defn1230}
Let $\mathcal U^i$ $(i=1,2)$ be Kuranishi charts
and $\Phi_{21} : \mathcal U^1 \to \mathcal U^2$ be an embedding
of Kuranishi charts.
Let $Z \subset U_1$ be a compact subset.
\par
A quadruple $(\pi_{12},\tilde\varphi_{21},\Omega_{12}, \Omega_1)$
is called a {\it bundle extension data}
\index{bundle extension data}
of  $(\Phi_{21},Z)$ if
\begin{enumerate}
\item
$\pi_{12} : \Omega_{12} \to \Omega_1$ is a continuous map,
where $\Omega_{12}$ is a neighborhood of $\varphi_{21}(Z)$
in $U_2$ and $\Omega_{1}$ is a neighborhood of $Z$
in $U_1$.
\item
$\pi_{12}$ is diffeomorphic to the projection
of the normal bundle in a neighborhood of $Z$ in the sense of Definition \ref{lem123000}.
\item $\tilde\varphi_{21} : \pi_{12}^*\mathcal E_1 \to \mathcal E_2$
is an embedding of vector bundle.
(See Definition \ref{def26222}.)
\item
The map
$$
\varphi_{21}^{*} \pi^*_{12} \mathcal E_{1} \to \varphi_{21}^*\mathcal E_{1}
$$
that is induced from $\tilde\varphi_{21}$
and $\varphi_{21}$ coincides with the bundle
map
$$
\widehat{\varphi}_{21} : \mathcal E_{1} \to \mathcal E_{2}
$$
which covers $\varphi_{21}$.
\end{enumerate}
Here $\pi^*_{12}\mathcal E_1$ and $\varphi_{12}^*\mathcal E_1$ denote the pullback bundles.
\end{defn}

\begin{rem}
Note here $\pi_{12}$ is a map between orbifolds but is not an embedding.
So it violates our thesis that we consider only an embedding as
a morphism of orbifolds.
This map however is identified with a restriction of the projection of
vector bundle so is also a map discussed in Section \ref{sec:ofd}.
We use this fact to define the pullback $\pi_{12}^*\mathcal E$ of the vector bundle
$\mathcal E$. (See Definition-Lemma \ref{pullbackbyproj}.)
\end{rem}
\begin{lem}\label{lem1224}
Suppose we are in Situation \ref{sit121111}.
Then for any  compact subset $Z$ of $U_1$ there exists a bundle extension datum
$(\pi_{12},\tilde\varphi_{21},\Omega_{12}, \Omega_1)$.
\end{lem}
\begin{proof}
Let $\Omega_1$ be a relatively compact open neighborhood
of $Z$.
Let $\pi : N_{U_1}U_2 \to U_1$ be the normal bundle.
Then by taking a Riemannian metric and exponential map,
we can find an open neighborhood $\tilde\Omega_{12}$ of the zero
section of $N_{U_1}U_2\vert_{\Omega_1}$ and a diffeomorphism $\tilde\Omega_{12} \to \Omega_{12}$
onto an open neighborhood $\Omega_{12}$ of $\varphi_{21}(\Omega_1)$.
(See for example \cite[Lemma 6.5]{fooooverZ}.)
Therefore we find $\pi_{12} : \Omega_{12} \to U_1$ as the composition of
the diffeomorphism $I _{12}: \Omega_{12} \cong \tilde\Omega_{12}$ and the projection
$\pi : N_{U_1}U_2 \to U_1$.
\par
The diffeomorphism $I _{12}$ is homotopic to the
composition of $\pi_{12} : \Omega_{12} \to \Omega_1$ and
$\varphi'_{21} : \Omega_1 \to \tilde\Omega_{12}$.
(Here $\varphi'_{21}$ is the embedding as the zero section of vector
bundle.)
Therefore
$\mathcal E_{2}\vert_{\Omega_{12}}$ is
isomorphic to
the pullback of $\mathcal E_2$ by $I _{12}^{-1}\circ \varphi'_{21}\circ \pi_{12}$.
(See Proposition \ref{homotopicpulback}.)  Note $I _{12}^{-1}\circ \varphi'_{21}\circ \pi_{12} = \varphi_{21}\circ \pi_{12}$.
Thus $({\varphi}_{21}\circ \pi_{12})^* \mathcal E_2 \cong \mathcal E_2$.
\par
$\widehat{\varphi}_{21} : \mathcal E_1 \to \mathcal E_2$
induces a bundle inclusion
$
\pi_{12}^*\mathcal E_1 \to ({\varphi}_{21}\circ \pi_{12})^* \mathcal E_2.
$
Therefore by using isomorphism
$({\varphi}_{21}\circ \pi_{12})^* \mathcal E_2 \cong \mathcal E_2$
we obtain the required embedding:
$
\tilde\varphi_{21} : (\pi_{12}^*\mathcal E_1)\vert_{\Omega_{12}} \to (\mathcal E_2)\vert_{\Omega_{12}}.
$
\par
We can take the homotopy $I _{12} \sim \varphi'_{21}\circ \pi_{12}$ so that its composition with
the inclusion $\varphi_{21}\vert_{\Omega_1} : \Omega_1 \to \Omega_{12} \subset U_2$  is the trivial homotopy between
${\varphi}_{21} : \Omega_1 \to \tilde\Omega_{12}$
 and
${\varphi}'_{21}\circ \pi \circ {\varphi}_{21} = {\varphi}_{21}$.
We can use this fact to check Definition \ref{defn1230} (4).
It is easy to check Definition \ref{defn1230} (1)(2)(3).
\end{proof}

Let $(\pi_{12},\tilde\varphi_{21},\Omega_{12}, \Omega_1)$
be a bundle extension datum as in Lemma \ref{lem1224}
and $\frak S^1 \in \mathscr S(\Omega_1)$.
We will define $\frak S^2 \in \mathscr S(\Omega_{12})$
which is compatible with $\frak S^1$ below.
\par
Let $x_2 \in \Omega_{12}$. We put $x_1 = \pi_{12}(x_2)$.
Since $\pi_{12}$ is diffeomorphic to a
restriction of a projection of a vector bundle
we can find orbifold charts
$(V_i,\Gamma_i,\phi_i)$ of $x_i$ in $U_i$
such that $(V_i,\Gamma_i,\phi_i)$, $\pi_{12}$ have the following properties:
(See Definition \ref{defn2613} (1).)
\begin{proper}\label{proper1227}
\begin{enumerate}
\item
$\Gamma_1 = \Gamma_2$.
\item
$V_2$ is identified with an open neighborhood of $V_1 \times \{0\}$ in $V_1 \times F$
where $F$ is a vector space, which is the fiber of the normal bundle $N_{U_1}U_2$.
\item
$\Gamma_1 = \Gamma_2$ acts on $V_1$ and has a linear action on $F$.
\item
The diagram
\begin{equation}\label{diagram1215}
\xymatrix{
V_1 \times F \ar@{{<}-^{)}}[r]\ar[d]  & V_2\ar[r]^{\phi_{2}}\ar[d]  & \Omega_{12}\ar^{\pi_{12}}[d] \\
V_1    \ar@{=}[r]   & V_1\ar^{\phi_1}[r] & \Omega_1
}
\end{equation}
commutes, where the first vertical arrow is the projection to the first factor.
\end{enumerate}
\end{proper}
For each given representative
$(W_{1},\omega_{1},\{\frak s_1^{\epsilon}\})$ of $\frak S^1$,
we  define
$\frak s_2^{\epsilon} : V_2 \times W_1 \to E_2$ by
\begin{equation}\label{defextSSSSS}
\frak s_2^{\epsilon}(y,\xi)
=
s_2(y)
+
\tilde\varphi_{21}
\left(
\frak s_1^{\epsilon}(\pi(y),\xi)
- s_1(\pi(y))
\right).
\end{equation}
\begin{lem}
Define
$(W_2,\omega_2,\{\frak s_2^{\epsilon}\}) = (W_1,\omega_1,
\{\frak s_2^{\epsilon}\})$.
Then it is a CF-perturbation of $\mathcal U_2$.
\end{lem}
\begin{proof}
The conditions (1),(2) and (4) of
Definition \ref{defn73ss} obviously hold by definition.
\par
The $C^1$ convergence,
$\lim_{\epsilon \to 0}\frak s_2^{\epsilon} = s_2$,
is a consequence of that of $\lim_{\epsilon \to 0}\frak s_1^{\epsilon} = s_1$.
If $y \in V_1 \times \{0\}$ then $s_2(y) = s_1(y)$,
the second term of the right hand side of (\ref{defextSSSSS})
is $\frak s_1^{\epsilon}(y,\xi) -s_1(y)$.
Therefore the right hand side of (\ref{defextSSSSS})
coincides with $\frak s_1^{\epsilon}(y,\xi)$
on $V_1 \times \{0\}$.
Definition \ref{defn73ss} (3) holds.
\end{proof}
\begin{lem}\label{lem1225new}
The equivalence class of $(W_2,\omega_2,\{\frak s_2^{\epsilon}\})$ is
independent of the choice of the orbifold charts
$(V_i,\Gamma_i,\phi_i)$ of $x_i$ in $U_i$
satisfying (1)(2)(3)(4) above.
\end{lem}
\begin{proof}
The proof is easy from definition and Lemma \ref{lem2622}.
\end{proof}
We remark that the equivalence class of $(W_2,\omega_2,
\{\frak s^{\epsilon}_2\})$
depends on the choice of bundle extension datum.
\begin{lem}\label{lemma1217}
Suppose $(W_1,\omega_1,\{\frak s^{\epsilon}_1\})$ is equivalent to
$(W'_1,\omega'_1,\{\frak s^{\epsilon \prime}_1\})$ in the sense of
Definition \ref{conmultiequiv11}.
We take and fix a bundle extension datum
$(\pi_{12},\tilde\varphi_{21},\Omega_{12}, \Omega_1)$.
We use it to define  $\frak s^{\epsilon}_2$ (resp. $\frak s^{\epsilon \prime}_2$)
by formula (\ref{defextSSSSS}) for $(W_1,\omega_1,\{\frak s^{\epsilon}_1\})$ (resp. for
$(W'_1,\omega'_1,\{\frak s^{\epsilon \prime}_1\})$)
and obtain $(W_1,\omega_1,\{\frak s^{\epsilon}_2\})$
(resp. $(W'_1,\omega'_1,\{\frak s^{\epsilon \prime}_2\})$).
\par
Then $(W_1,\omega_1,\{\frak s^{\epsilon}_2\})$ is equivalent to
$(W'_1,\omega'_1,\{\frak s^{\epsilon \prime}_2\})$.
\end{lem}
\begin{proof}
It suffices to show the lemma in the case when
$(W_1,\omega_1,\{\frak s^{\epsilon}_1\})$ is a projection of
$(W'_1,\omega'_1,\{\frak s^{\epsilon \prime}_1\})$.
Let $\Pi : W'_1 \to W_1$ be the projection such that
$\Pi!(\omega'_1) = \omega_1$ and
$\frak s_1^{\epsilon}(y,\Pi(\xi)) = \frak s_1^{\epsilon \prime}(y,\xi)$.
Then (\ref{defextSSSSS}) implies
$\frak s_2^{\epsilon}(y,\Pi(\xi)) = \frak s_2^{\epsilon \prime}(y,\xi)$.
This implies the lemma.
\end{proof}
We use these lemmata to prove Proposition \ref{prop1221} as follows.
\par
Take a bundle extension datum
$(\pi_{12},\tilde\varphi_{21},\Omega_{12}, \Omega_1)$
as in Lemma \ref{lem1224}.
\par
We put $\mathcal O^1 = \Omega_1$, and denote $\mathcal O^2 = \Omega_{12}$ which is an open neighborhood of $K$ in $U_2$ where $\frak S^{2,K}$ is defined on.
We take open sets $O^i_1,O^i_2 \subset U_i$ ($i=1,2$) such that
\begin{equation}
K \subset O^i_1 \subset \overline{O^i_1} \subset O^i_2 \subset \overline{O^i_2}  \subset \mathcal O^i
\end{equation}
and $\varphi_{21}^{-1}(O^2_j) = O^1_j$.
Put $C = \overline{O^1_2} \setminus O^1_1$.
\par
Then we take a smooth function $\chi : U_1 \to [0,1]$ such that
\begin{equation}
\chi(p) =
\begin{cases}
1 &\text{on a neighborhood of $\overline{O^1_1}$}, \\
0  &\text{on a neighborhood of $U_1 \setminus {O^1_2}$}.
\end{cases}
\end{equation}
\par
Now we are ready to define an extension $\mathfrak S^2$ of $\mathfrak S^1$.
Let $\frak z \in \mathcal O^1$. We will define $(\mathfrak S^2)_{\frak z}
\in \mathscr S_{\varphi_{21}(\frak z)}$ in the following three cases separately.
\par\smallskip
\noindent (Case 1): $\frak z \in \overline{O^1_1}$.
\par
In this case, we set
$(\mathfrak S^2)_{\frak z}$ to be a germ of  $\frak S^{2,K}$.
\par
\noindent (Case 2):  $\frak z \in \mathcal O^1 \setminus {O^1_2}$.
\par
In this case, we use $(\pi_{12},\tilde\varphi_{21},\Omega_{12}, \Omega_1)$ to extend $\frak S^1$ to
$\mathcal O^2$ by Formula (\ref{defextSSSSS}).
We then obtain $(\mathfrak S^2)_{\frak z}$.
\par
\noindent (Case 3):  $\frak z \in C$.
\par
We take a
representative of the germ of $\frak S^{1}$ at ${\frak z}$, which we
denote by  $(\frak V^1_{{\frak z}},\mathcal S^{1}_{{\frak z}})$.
We put
$\mathcal S^{1}_{{\frak z}}
= (W,\omega,\{\frak s^{\epsilon}_1\})$
and
$\frak V^1_{{\frak z}}
= (V_1,E_1,\Gamma_1,\phi_1,\widehat{\phi}_1)$.
Then, by shrinking $O^2$ if necessary, we may assume that the germ
$\frak S^{2,K}_{{\frak z}}
\in (\mathscr S^{\mathcal U^1\triangleright \mathcal U^2})_{\varphi_{21}(\frak z)}$ is represented by
$(\frak V^2_{{\frak z}},\mathcal S^{2,K}_{{\frak z}})$
such that
\begin{enumerate}
\item
$\frak V^2_{{\frak z}} = (V_2,E_2,\Gamma_2,\phi_2,\widehat{\phi}_2)$,
where
$(V_i,\Gamma_i,\phi_i)$, $\pi_{12}$
have Property \ref{proper1227}.
\item
$\mathcal S^{2,K}_{{\frak z}}
= (W,\omega,\{\frak s^{\epsilon}_{2,K}\})$
where
$W$ and $\omega$ are the same as those appearing in $\frak S^{1}$.
\end{enumerate}
The map $\tilde \varphi_{21}$, which is a part of bundle extension data,
rise to a map
$$
\breve\varphi_{21,\frak z} :
V_2 \times E_1 \to E_2
$$
which is $\Gamma_2$ equivariant.
\par
Now we define
$(\mathfrak S^2)_{\frak z}
\in \mathscr S_{\varphi_{21}(\frak z)}$
as $(W,\omega,\{\frak s^{\epsilon}_{2}\})$
where
\begin{equation}\label{defext3030}
\aligned
\frak s^{\epsilon}_{2}(y,\xi)
=
s_{2}(y)
&+
\chi([y]) (\frak s^{\epsilon}_{2,K}(y,\xi) - s_{2}(y)) \\
&+
(1-\chi([y]))
\breve \varphi_{21,\frak z}
\left(y,
\frak s^{\epsilon}_1(\pi(y),\xi))
- s_{1}(\pi(y))
\right).
\endaligned
\end{equation}
Here $\pi : V_2 \to V_1$ is a restriction of the projection $V_1 \times F \to V_1$,
which represents the map $\pi_{12}$, which is a part of bundle extension data.
We denote by $[y] \in V_2/\Gamma_2$ the equivalence class of $y$ which we identify with
an element of $U_2$ by an abuse of notation.
\begin{lem}
\begin{enumerate}
\item
$(\mathfrak S^2)_{\frak z}$ defined above is
independent of the choices made in the definition and depends only
on $(\pi,\tilde \varphi_{21})$, $\chi$ and  $\frak S^{2,K}$,
 $\frak S^{1}$.
\item
Moreover $(\mathfrak S^2)_{\frak z}$ for various $\frak z$ defines a section of
$\mathfrak S^2$ of $\varphi_{21}^{\star}\mathscr S$.
\item
$\mathfrak S^2$ is a section of
$\varphi_{21}^{\star}\mathscr S^{\mathcal U^1\triangleright \mathcal U^2}$
and
$\Phi_{21}^*(\mathfrak S^2) = \mathfrak S^1$.
\end{enumerate}
\end{lem}
\begin{proof}
In Case 1, the well-defined-ness is obvious.
In Case 2,  the well-defined-ness follows from Lemmata \ref{lem1225new} and \ref{lemma1217}.
To prove the well-defined-ness in Case 3, it suffices to consider the case of projection, which follows immediately
from (\ref{defext3030}).
We have thus proved Statement (1).
\par
To prove Statement (2), it suffices to show the next two facts (a)(b).
\begin{enumerate}
\item[(a)]
If $\frak z \in \overline{O^1_1} \cap C$ then $(\mathfrak S^2)_{\frak z}$
obtained by applying  Case 1 is equivalent to $(\mathfrak S^2)_{\frak z}$
obtained by applying Case 3.
\item[(b)]
If $\frak z \in  (U'_1 \setminus {O^1_2}) \cap C$ then $(\mathfrak S^2)_{\frak z}$
obtained by applying  Case 2 is equivalent to $(\mathfrak S^2)_{\frak z}$
obtained by applying  Case 3.
\end{enumerate}
To prove (a) we remark that $\chi = 1$ on a neighborhood of $\frak z$.
Therefore (\ref{defext3030}) becomes
$$
\frak s^{\epsilon}_{2}(y,\xi)
=
s_{2}(y)
+
(\frak s^{\epsilon}_{2,K}(y,\xi) - s_{2}(y))
=
\frak s^{\epsilon}_{2,K}(y,\xi),
$$
as required.
\par
To prove (b) we remark that $\chi = 0$ on a neighborhood of $\frak z$.
Therefore (\ref{defext3030}) becomes
$$
\frak s^{2,\epsilon}_{\frak z}(y,(\xi^1,\xi^2))
=
s_2(y) +
\breve \varphi_{21,\frak z}
\left(y,
\frak s^{\epsilon}_1(\pi(y),\xi))
- s_{1}(\pi(y))
\right).
$$
The right hand side coincides with (\ref{defextSSSSS}),
as required.
\end{proof}
By the way how we defined $(\mathfrak S^2)_{\frak z}$ in Case 1,
it is easy to see that the restriction of $\frak S^2\vert_{U_2''}$ to $U_2'' \cap O^2_1$
is equivalent to the restriction of $\frak S^{2,K}$ on $U_2'' \cap O^2_1$.
This implies the first formula of (\ref{form1214})
\par
The second formula of (\ref{form1214}) follows from the way
how we defined $(\mathfrak S^2)_{\frak z}$ in Cases 2 and 3.
In fact if $y = \varphi_{21}(x)$ then (\ref{defext3030})
becomes
\begin{equation}
\aligned
\frak s^{\epsilon}_{2}(y,\xi)
=
s_{2}(y)
&+
\chi([x]) (\frak s^{\epsilon}_{2,K}(\varphi_{21}(x),\xi) - s_{2}(y)) \\
&+
(1-\chi([x]))
\breve \varphi_{21,\frak z}
\left(\varphi_{21}(x),
\frak s^{\epsilon}_1(x,\xi))
- s_{1}(x)
\right) \\
=
s_{1}(x)&+
\chi([x]) (\frak s^{\epsilon}_{1}(x,\xi) - s_{1}(x)) \\
&+
(1-\chi([x]))
(\breve \varphi_{21,\frak z}
\left(
\frak s^{\epsilon}_1(x,\xi))
- s_{1}(x)
\right) \\
&\!\!\!\!\!\!\!\!\!\!\!\!\!\!\!\!\!\!\!\!\!\!= \frak s^{\epsilon}_{1}(x,\xi).
\endaligned
\end{equation}
\par
To complete the proof of Proposition \ref{prop1221},
it remains to prove that if
$\frak S^1 \in \mathscr S_{\star}^{\mathcal U^1}(U_1)$,
$\frak S^{2,K} \in (\varphi_{21}^{\star}\mathscr S_{\sharp}^{\mathcal U^1\triangleright \mathcal U^2})(K)$
then $\frak S^2 \in (\varphi_{21}^{\star}\mathscr S_{\sharp}^{\mathcal U^1\triangleright \mathcal U^2})(Z)$.
This is actually an immediate consequence
of the fact that the condition for a section of $\mathscr S^{\mathcal U^2}$
to be a section of $\mathscr S_{\sharp}^{\mathcal U^2}$
is an open condition.
\end{proof}

\subsection{Construction of CF-perturbations
on good coordinate system}
\label{subsec:cfpgoodcsys}

In this subsection we complete the proof of Theorem  \ref{existperturbcont}.

\begin{proof}[Proof of Theorem  \ref{existperturbcont}]
We first discuss the absolute case.
We will construct a CF-perturbation on the Kuranishi charts $\mathcal U_{\frak p}$ by the
{\it upward} induction on $\frak p \in \frak P$ with respect
to the partial order of $\frak P$.
(We remark that during the construction of good coordinate
system we used downward induction on the partial order of
$\frak P$. So our construction here goes the opposite direction from
that of Section \ref{sec:contgoodcoordinate}.)
\par
We say a subset $\frak F \subseteq \frak P$  a
{\it filter} \index{filter} if
$\frak p,\frak q \in \frak P$, $\frak p \ge \frak q$,
$\frak p \in \frak F$ imply $\frak q \in \frak F$.
In this article we regard $\emptyset$ as a filter.
(This may be different from the usual convention.)
Let $\frak F \subseteq \frak P$ be a filter.
The main ingredient of the proof of Theorem  \ref{existperturbcont}
is Proposition \ref{existontiiindc1t} below, which we prove by an
upward induction on $\#\frak F$.
To state Proposition \ref{existontiiindc1t} we need several notations.
For an arbitrary subset $\frak F \subseteq \frak P$
we put
\begin{equation}\label{form1210}
\frak T(\frak F,\mathcal K)
=
\bigcup_{\frak p \in \frak F}
\psi_{\frak p}(s_{\frak p}^{-1}(0) \cap \mathcal K_{\frak p}).
\end{equation}
This is a compact subset of $X$.
We next define a good coordinate system
on a neighborhood of $\frak T(\frak F,\mathcal K)$.
We take $\mathcal K^+$ so that
$(\mathcal K,\mathcal K^+)$ is a support pair of
${\widetriangle{\mathcal U}}$ and put
\begin{equation}\label{form1211}
\frak Y(\frak F)
=
\bigcup_{\frak p \in \frak F}\psi_{\frak p}(s_{\frak p}^{-1}(0) \cap
\mathcal {\rm Int}\,K^+_{\frak p}).
\end{equation}
We consider the set of Kuranishi charts
$\{ \mathcal U_{\frak p}\vert_{ {\rm Int}\,K^+_{\frak p}}
\mid \frak p \in \frak F\}$.
Together with the restrictions of the coordinate changes of
${\widetriangle{\mathcal U}}$
it defines a good coordinate system of $\frak Y(\frak F)$ on
$\frak T(\frak F,\mathcal K)$.
We denote this good coordinate system by
${\widetriangle{\mathcal U}}(\frak F,\mathcal K^+)$.
Note
$\{\mathcal K_{\frak p} \mid \frak p \in \frak F\}$ is a
support system of ${\widetriangle{\mathcal U}}(\frak F,\mathcal K^+)$.
We denote it by $\mathcal K(\frak F)$.
\begin{rem}
The main difference between ${\widetriangle{\mathcal U}}(\frak F,\mathcal K^+)$
and ${\widetriangle{\mathcal U}}$ lies on the fact that
for ${\widetriangle{\mathcal U}}(\frak F,\mathcal K^+)$
we take the Kuranishi charts $\mathcal U_{\frak p}$
with $\frak p \in \frak F$ only.
\end{rem}
We now prove the following proposition by an induction.
This inductive proof is the same as one we had written
in \cite[page 955 line 17-24]{FO} in a similar situation of
multivalued perturbation.  Here we provide much more detail.
\begin{prop}\label{existontiiindc1t}
For any filter $\frak F \subseteq \frak P$,
there exists a CF-perturbation
$\widetriangle{\frak S^{\frak F}}$ of
$({\widetriangle{\mathcal U}}(\frak F,\mathcal K^+),\mathcal K(\frak F))$
on $\frak T(\frak F,\mathcal K)$.
It satisfies the following properties.
\begin{enumerate}
\item
$\widetriangle{\frak S^{\frak F}}$ is transversal to $0$.
\item
If $\widetriangle{f} : (X,Z;\widetriangle{\mathcal U}) \to M$ is a weakly submersive strongly
smooth map, then $\widetriangle{\frak S^{\frak F}}$ can be taken so that $\widetriangle{f}$
is strongly submersive with respect to $\widetriangle{\frak S^{\frak F}}$.
\item
If $\widetriangle{f} : (X,Z;\widetriangle{\mathcal U}) \to M$ is a strongly
smooth map which is weakly transversal to $g : N \to M$,
then $\widetriangle{\frak S^{\frak F}}$ can be taken so that $\widetriangle{f}$
is strongly transversal to $g$ with respect to $\widetriangle{\frak S^{\frak F}}$.
\end{enumerate}
\end{prop}
\begin{proof}
The proof is by induction of $\# \frak F$. If $\frak F = \emptyset$,
there is nothing to show.
\par
Suppose Proposition \ref{existontiiindc1t} is proved
for all $\frak F'$ with $\# \frak F' < \# \frak F$.
We will prove the case of $\frak F$.
\par
Let $\frak p_0$ be a maximal element of $\frak F$.
Then $\frak F_- =\frak F \setminus \{ \frak p_0\}$ is a filter.
By the standing induction hypothesis, there exists a
CF-perturbation
$\widetriangle{\frak S^{\frak F_-}}$ of
$({\widetriangle{\mathcal U}}(\frak F_-,\mathcal K^+),\mathcal K(\frak F_-))$
on $\frak T(\frak F_-,\mathcal K)$.
\par
We will extend this CF-perturbation to
a CF-perturbation $\widetriangle{\frak S^{\frak F}}$
in two steps.
In the first  step we use Proposition \ref{prop1221}
to extend it to a
CF-perturbation on
$({\widetriangle{\mathcal U}}(\frak F_,\mathcal K^+),\mathcal K(\frak F_-))$
of $\frak T(\frak F_-,\mathcal K)$.
We then obtain $\widetriangle{\frak S^{\frak F}}$ in the second step.
\par
Note the difference between
${\widetriangle{\mathcal U}}(\frak F_-,\mathcal K^+)$
and
${\widetriangle{\mathcal U}}(\frak F,\mathcal K^+)$
is that an open subchart of $\mathcal U_{\frak p_0}$
is included in the latter.
So, the main integrant of the first step is defining
a CF-perturbation on $\mathcal U_{\frak p_0}$
in a neighborhood of $\frak T(\frak F_-,\mathcal K)$.
The second step uses Proposition \ref{prop123123}
to extend it to a Kuranishi neighborhood of $\frak T(\frak F,\mathcal K)$.
\par
We remark that by definition there exists a GG-embedding
$$
\widetriangle{\Phi_{\frak F\frak F_-}} :
{\widetriangle{\mathcal U}}(\frak F_-,\mathcal K^+)
\to {\widetriangle{\mathcal U}}(\frak F,\mathcal K^+).
$$
(We recall $\frak F = \frak F_- \cup \{\frak p_0\}$.)
By induction hypothesis,
we have  a CF-perturbation
$\widetriangle{{\frak S}^{\frak F_-}}$ of
$({\widetriangle{\mathcal U}}(\frak F_-,\mathcal K^+),\mathcal K(\frak F_-))$
on $\frak T(\frak F_-,\mathcal K)$.
\par
We consider an ideal $\frak I \subseteq \frak F_-$.
Namely $\frak I$ is a subset such that $\frak p \in \frak I$,
$\frak q \in \frak F_-$ and $
\frak p \le \frak q$ imply $\frak q \in \frak I$.
\par
In the next lemma, we take a neighborhood
$U_{\frak p_0}(\frak I)$ of
$\psi_{\frak p_0}^{-1}(\frak T(\frak I,\mathcal K)) \cap {\mathcal K}_{\frak p_0}$
in $\mathcal K_{\frak p_0}^+$.
Then we replace the Kuranishi neighborhood $\mathcal U_{\frak p_0}$,
which is a member of the good coordinate system $\widetriangle{\mathcal U^{\frak F}}$,
by $\mathcal U_{\frak p_0}\vert_{U_{\frak p_0}(\frak I)}$.
We denote the resulting good coordinate system by
$\widetriangle{\mathcal U}(\frak F,\frak I;\mathcal K^+)$.
It is a good coordinate system of
$\frak T(\frak F_-,\mathcal K)$.
(The difference between $\widetriangle{\mathcal U}(\frak F,\frak I;\mathcal K^+)$ and
$\widetriangle{\mathcal U^{\frak F_-}}$ is that
$\widetriangle{\mathcal U}(\frak F,\frak I;\mathcal K^+)$ has one
more Kuranishi chart $\mathcal U_{\frak p_0}\vert_{U_{\frak p_0}(\frak I)}$
than $\widetriangle{\mathcal U^{\frak F_-}}$.)
\par
We can define a GG embedding $\widetriangle{\mathcal U^{\frak F_-}}
\to \widetriangle{\mathcal U}(\frak F,\frak I;\mathcal K^+)$
so that the map induced on the index set of the Kuranishi charts is an
obvious embedding $\frak F_- \to \frak F$, which we
denote by $\widehat\Phi_{\frak F_-\frak F;\frak I}$.
\begin{lem}\label{122222}
For any ideal $\frak I \subseteq \frak F_-$, there exist
$U_{\frak p_0}(\frak I)$ and a CF-perturbation
$\widetriangle{\frak S^{\frak F}}(\frak I)$ of $(\widetriangle{\mathcal U}(\frak F,\frak I;\mathcal K^+),
\frak T(\frak F_-,\mathcal K))$
with the following properties.
\begin{enumerate}
\item
$\widetriangle{\frak S^{\frak F}}(\frak I)$, $\widetriangle{\frak S^{\frak F_-}}$
are compatible with the embedding
$\widehat\Phi_{\frak F_-\frak F;\frak I}$.
\item
\begin{enumerate}
\item
If $\widetriangle{{\frak S}^{\frak F_-}}$  is transversal to $0$
so is $\widetriangle{\frak S^{\frak F}}(\frak I)$.
\item
If  $\widetriangle f$ is strongly submersive with respect to
$\widetriangle{{\frak S}^{\frak F_-}}$
and $\widetriangle f$ is weakly submersive then
$\widetriangle f$ is strongly submersive with respect to $\widetriangle{\frak S^{\frak F}}(\frak I)$.
\item
If $\widetriangle f$ is strongly transversal to $g : N \to M$ with respect to
$\widetriangle{{\frak S}^{\frak F_-}}$
and $\widetriangle f$ is weakly transversal to $g : N \to M$  then
$\widetriangle f$ is strongly transversal to $g : N \to M$ with respect to
$\widetriangle{\frak S^{\frak F}}(\frak I)$.
\end{enumerate}
\end{enumerate}
\end{lem}
\begin{proof}
Recall
$\frak I \subseteq \frak F_- = \frak F \setminus \{\frak p_0\}$.
The proof is by an upward induction over $\#\frak I$.
For $\frak I = \emptyset$, we put $U_{\frak p_0}(\emptyset) = \emptyset$.
Then $\widetriangle{\mathcal U}(\frak F,\emptyset;\mathcal K^+)
= \widetriangle{\mathcal U}(\frak F_-;\mathcal K^+)
$.
We put
$\widetriangle{\frak S^{\frak F}}(\emptyset)
= \widetriangle{\frak S^{\frak F_-}}$. It is easy to see that
it has the required properties.
\par
We assume that the lemma is proved for $\frak I'$ with
$\#\frak I' <\#\frak I$ and will prove the case of $\frak I$.
\par
Let $\frak p_1$ be a minimal element of $\frak I$ and
denote $\frak I_- = \frak I \setminus \{\frak p_1\}$.
Then $\frak I_-$ is an ideal. We use the induction hypothesis
to obtain $\widetriangle{\frak S^{\frak F}}(\frak I_-)$ where $\widetriangle{\frak S^{\frak F}}(\frak I_-)$ is a
CF-perturbation of
$(\widetriangle{\mathcal U}(\frak F,\frak I_-;\mathcal K^+),
\frak T(\frak F_-,\mathcal K))$.
\par
\par
We apply Proposition \ref{prop1221} by putting:
\begin{equation}\label{form121212}
\aligned
&\frak S^{2,K} = (\widetriangle{\frak S^{\frak F}}(\frak J_-))_{\frak p_0},
\qquad
\frak S^1 = (\widetriangle{\frak S^{\frak F_-}})_{\frak p_1}, \\
&
\mathcal U_1 = \mathcal U_{\frak p_1},\quad
\mathcal U_2 = \mathcal U_{\frak p_0},
\qquad
K = \frak T(\frak I_-,\mathcal K).
\endaligned
\end{equation}
We denote by $(\widetriangle{\frak S^{\frak F}}(\frak J_-))_{\frak p_0}$
a CF-perturbation induced by
$\widetriangle{\frak S^{\frak F}}(\frak J_-)$
on an open subchart $\mathcal U_{\frak p_0}\vert_{U_{\frak p_0}(\frak I_-)}$
of $\mathcal U_{\frak p_0}$.
The definition of $(\widetriangle{\frak S^{\frak F_-}})_{\frak p_1}$ is similar.
\par
We now apply Proposition \ref{prop1221} to obtain $U_{\frak p_0}(\frak I)$ and
a CF-perturbation
$(\widetriangle{\frak S^{\frak F}}(\frak J))_{\frak p_0}$
of an open subchart $\mathcal U_{\frak p_0}\vert_{U_{\frak p_0}(\frak I)}$
of $\mathcal U_{\frak p_0}$.
\par
Among the CF-perturbations on Kuranishi charts consisting $\widetriangle{\frak S^{\frak F}}(\frak J)$, we replace
$(\widetriangle{\frak S^{\frak F}}(\frak J_-))_{\frak p_0}$
by $(\widetriangle{\frak S^{\frak F}}(\frak J))_{\frak p_0}$
and obtain a required CF-perturbation $\widetriangle{\frak S^{\frak F}}(\frak I)$.
\end{proof}
We have thus completed Step 1 of the proof of
Proposition \ref{existontiiindc1t}.
\par
We take the case $\frak I = \frak F_-$ of Lemma \ref{122222}.
Then we have $\widetriangle{\frak S^{\frak F}}(\frak F_-)$
which include $(\widetriangle{\frak S^{\frak F}}(\frak F_-))_{\frak p_0}$
such that:
\begin{enumerate}
\item
$(\widetriangle{\frak S^{\frak F}}(\frak F_-))_{\frak p_0}$ is a CF-perturbation of  $\mathcal U_{\frak p_0}$ on
$\frak T(\frak F_-,\mathcal K)$.
\item
$\widetriangle{\frak S^{\frak F}}(\frak F_-)$, $\widetriangle{\frak S^{\frak F_-}}$
are compatible with the embedding
$\widehat{\Phi_{\frak F\frak F_-}}$ on a neighborhood of $\frak T(\frak F_-,\mathcal K)$.
\end{enumerate}
We now apply Proposition \ref{prop123123} to
$\mathcal U_{\frak p_0}$ and
find that
\begin{equation}\label{1213formula}
\mathscr F_{\sharp}(V) \to \mathscr F_{\sharp}(K)
\end{equation}
is surjective.
Here $\sharp$ is one of ${\pitchfork 0}$, ${f \pitchfork}$,
${f \pitchfork g}$.
We put
\begin{equation}\label{12132formula}
K = \frak T(\frak F_-,\mathcal K),
\qquad
V = \text{an open neighborhood of $\mathcal K_{\frak p_0}$}.
\end{equation}
Therefore we can extend
$(\widetriangle{\frak S^{\frak F}}(\frak F_-))_{\frak p_0}$
to a CF perturbation
$(\widetriangle{\frak S^{\frak F}}(\frak F))_{\frak p_0}$
on a Kuranishi neighborhood of $\frak T(\frak F,\mathcal K)$
that is an open subchart of $\mathcal U_{\frak p_0}$.
\par
Replacing $(\widetriangle{\frak S^{\frak F}}(\frak F_-))_{\frak p_0}$
in $\widetriangle{\frak S^{\frak F}}(\frak F_-)$
by this extension we obtain $\widetriangle{\frak S^{\frak F}}$.
\par
The proof of Proposition \ref{existontiiindc1t} is complete.
\end{proof}
Theorem  \ref{existperturbcont} is the case
of Proposition \ref{existontiiindc1t} when
$\frak F = \frak P$.
\end{proof}
\begin{proof}[Proof of Proposition \ref{existperturbcontrel}]
Proposition \ref{existperturbcontrel}
is a relative version of  Theorem \ref{existperturbcont}
and the proof is  mostly the same.
We use the symbol $\mathcal Z_{(i)}$ in place of $Z_i$ during the proof of
Proposition \ref{existperturbcontrel}
(since we used $Z$ for other objects in this section already).
\par
We replace (\ref{form1210}) by
$\frak T(\frak F,\mathcal K;\mathcal Z_{(1)}) = \frak T(\frak F,\mathcal K)
\cup \mathcal  Z_{(1)}$ and
(\ref{form1211}) by
$\frak Y(\frak F) \cup \mathcal Z_{(1)} =\frak Y(\frak F;\mathcal Z_{(1)})$.
We have a good coordinate system ${\widetriangle{\mathcal U}}(\frak F,\mathcal K^+) \cup {\widetriangle{\mathcal U^{\mathcal Z_{(1)}}}}$
on $\frak Y(\frak F;\mathcal Z_{(1)})$, which we denote by
${\widetriangle{\mathcal U}}(\frak F,\mathcal K^+;\mathcal Z_{(1)})$.
We take a support system $\mathcal K^{\mathcal Z_{(1)}}$ of
${\widetriangle{\mathcal U^{\mathcal Z_{(1)}}}}$.
Together with $\mathcal K(\frak F)$ it gives a support system of
 $\frak T(\frak F,\mathcal K;\mathcal Z_{(1)})$, which we denote by
$\mathcal K(\frak F;\mathcal Z_{(1)})$.
\par
Proposition \ref{existontiiindc1t}
is replaced by the following:
\begin{enumerate}
\item[(*)]
There exists a CF-perturbation
$\widetriangle{\frak S^{\frak F,\mathcal Z_{(1)}}}$ of
$({\widetriangle{\mathcal U}}(\frak F,\mathcal K^+;\mathcal Z_{(1)}),
\mathcal K(\frak F;\mathcal Z_{(1)}))$
on $\frak T(\frak F,\mathcal K;\mathcal Z_{(1)})$.
The same properties as (1)(2)(3) of Proposition \ref{existontiiindc1t}
are satisfied.
\par
Moreover the CF-perturbation $\widetriangle{\frak S^{\frak F,\mathcal Z_{(1)}}}$ coincides with a restriction of
$\widetriangle{\frak S^1}$ (given as a part of assumption) on the charts of
$\widetriangle{\mathcal U^{(1)}}$.
\end{enumerate}
We prove (*) by the same induction as the proof of
Proposition \ref{existontiiindc1t}.
Namely we prove:
\begin{enumerate}
\item[(**)]
There exists a CF-perturbation
$\widetriangle{\frak S^{\frak F,\mathcal Z_{(1)}}}(\frak I)$ of
$({\widetriangle{\mathcal U}}(\frak F,\frak I;\mathcal K^+;\mathcal Z_{(1)}),\mathcal K(\frak F;\mathcal Z_{(1)}))$
on $\frak T(\frak I,\mathcal K;\mathcal Z_{(1)})$
with the following properties.
\begin{enumerate}
\item
$\widetriangle{\frak S^{\frak F,\mathcal Z_{(1)}}}(\frak I)$,
$\widetriangle{{\frak S}^{\frak F_-,\mathcal Z_{(1)}}}$
are compatible with the embedding
$\widehat{\Phi_{\frak F\frak F_-;\frak I;\mathcal Z_{(1)}}}$ on a neighborhood of $\frak T(\frak I,\mathcal K;\mathcal Z_{(1)})$.
\item
The same properties as (1)(2)(3) of Proposition \ref{existontiiindc1t}
are satisfied.
\end{enumerate}
\end{enumerate}
\par\smallskip
Here we take $U_{\frak p_0}(\frak I)$ and define ${\widetriangle{\mathcal U}}(\frak F,\frak I;\mathcal K^+;\mathcal Z_{(1)})$
in the same way as ${\widetriangle{\mathcal U}}(\frak F,\frak I;\mathcal K^+)$.
The embedding $\widehat{\Phi_{\frak F\frak F_-;\frak I;\mathcal Z_{(1)}}}$ is obtained from
$\widehat{\Phi_{\frak F\frak F_-}}$ by using the identity embedding for the charts
in ${\widetriangle{\mathcal U^{\mathcal Z_{(1)}}}}$.

The proof of (**) is by the same induction as the proof of Lemma \ref{122222},
where we replace (\ref{form121212}) by
\begin{equation}
\aligned
&
\frak S^{2,K} = (\widetriangle{\frak S^{\frak F,\mathcal Z_{(1)}}}(\frak J_-))_{\frak p_0},
\qquad
\frak S_1 = (\widetriangle{\frak S^{\frak F_-,\mathcal Z_{(1)}}})_{\frak p_1},
\\
&\mathcal U_1 = \mathcal U_{\frak p_1},
\qquad\qquad\quad\,\,
\mathcal U_2 = \mathcal U_{\frak p_0},
\quad\qquad\qquad
K = \frak T(\frak I_-,\mathcal K;\mathcal Z_{(1)}).
\endaligned
\nonumber\end{equation}
Using (**) we can prove (*)  in the same way as the last step
in the proof of Theorem \ref{existperturbcont}.
The proof of Proposition   \ref{existperturbcontrel}  is complete.
\end{proof}\par

\section{Construction of multisections}
\label{sec:constrsec}

In Sections \ref{sec:constrsec} we discuss the
multivalued perturbation,
especially its existence result, Theorem \ref{prop621}.
This result will be used in Section \ref{sec:onezerodim}.
One of the advantages using multivalued perturbations is it enables us to work with $\Q$ coefficients.
In the construction based on de Rham theory
we can work only with $\R$ or $\C$.
For many applications, it is enough to work with coefficients $\R$ or $\C$.
For these cases, we do not need to use the results of
Sections \ref{sec:constrsec} and  \ref{sec:onezerodim}.

\subsection{Construction of multisection on a single chart}
\label{subsec:musecexsingle}

The proof of Theorem \ref{prop621}
is similar to the proof of Theorem \ref{existperturbcont}.
\par
We begin with proving a version of Proposition  \ref{prop123123} (2)(4).
(We remark that the multisection version of Proposition  \ref{prop123123} (3) does
not seem to exist.)
\begin{prop}\label{prop127777ver}
In Situation \ref{situ121}, let $K \subset U$ be a
compact subset and $\widetriangle{\frak s_K} = \{\frak s_{K}^n\}$  a multivalued
perturbation of $\mathcal U$ on a neighborhood of $K$ transversal to $0$.
\par
Let $\Omega \subset U$ be an relatively compact open subset such that
$K \subseteq \Omega$.
\begin{enumerate}
\item
There exists a multivalued
perturbation
$\widetriangle{\frak s} = \{\widetriangle{\frak s^n}\}$ of $\mathcal U$ on $\Omega$ such that $\widetriangle{\frak s}$  is
transversal to $0$ and is
equal to $\widetriangle{\frak s_K}$ on a neighborhood of $K$.
\item
Suppose we are in Situation \ref{situ122}.
Then we may choose $\widetriangle{\frak s}$ so that $f$ is strongly transversal to $g$ with respect to $\widetriangle{\frak s^n}$
for sufficiently large $n$. (See Definition \ref{transofdvect}.)
\end{enumerate}
\end{prop}
\begin{proof}
We use Lemma \ref{lemma12770} to obtain
$\{\frak V_{\frak r} \vert \frak r \in \frak R\}$ and $\{\mathcal  S_{\frak r} \vert
\frak r \in \frak R \}$ with the following properties.
Let $U_0$ be an open neighborhood of $K$
which we will fix later.
\begin{proper}
\begin{enumerate}
\item
$\frak V_{\frak r} = (V_{\frak r},\Gamma_{\frak r},E_{\frak r},\phi_{\frak r},\widehat\phi_{\frak r})$
is an  orbifold chart of $U$. We put $U_{\frak r} = {\rm Im}(\phi_{\frak r})$ and assume
$U_{\frak r} \cap K = \emptyset$.
\item
\begin{equation}
U_0 \cup \bigcup_{\frak r \in \frak R} U_{\frak r}
\end{equation}
is an open covering of $\overline{\Omega}$.
\item
$\mathcal  S_{\frak r} = (W_{\frak r},\omega_{\frak r},\{\frak s_{\frak r}^{\epsilon}\})$
is a CF-perturbation of $\mathcal U$ on $\frak V_{\frak r}$.
(Definition \ref{contipertlocalrest}.)
\item
$\mathcal  S_{\frak r}$ is strongly transversal. (Definition \ref{strongtransvers}.)
\end{enumerate}
 \end{proper}
We may assume that the given $\Gamma_{\frak r}$ action on $W_{\frak r}$
is effective, by replacing $W_{\frak r}$ with the product $W_{\frak r} \times W'$
if necessary,
where $W'$ is a faithful representation of the finite group $\Gamma_{\frak r}$.
Then we define $\frak s_{\frak r}^{\epsilon}$ to be the pull-back of
a multisection defined on $W_{\frak r}$ by the projection $W_{\frak r} \times W' \to W_{\frak r}$
and define our $\omega_{\frak r}$ on $W_{\frak r} \times W'$ to be the product form
of a smooth $\Gamma_{\frak r}$ invariant top degree form of compact support
on $W'$ and the given smooth top degree form on $W_{\frak r}$.
\par
For each $\frak r$ we take an open subset
$W_{\frak r}^0$ of $W_{\frak r}$ such that
\begin{equation}\label{form132222}
\gamma W_{\frak r}^0 \cap W_{\frak r}^0 = \emptyset,
\end{equation}
if $\gamma \in \Gamma_{\frak r} \setminus \{1\}$
(we use effectivity of $\Gamma_{\frak r}$ action here) and  put
\begin{equation}\label{form1217777}
W_0 = \prod_{\frak r \in \frak R}W_{\frak r}^0.
\end{equation}
\par
For each $x \in K$ we take a representative of $\frak s^n_K$ in an
orbifold chart $\frak V_x$
at $x$ and denote it by
$(\frak s^n_{x,1},\dots,\frak s^n_{x,\ell(x)})$.
We take a finite number of points $x_1,\dots,x_k$ of $K$ so that
$\bigcup_{a=1}^k U_{{x_a}}  \supset K$
and the sections
$\frak s^n_{x,k}$ are transversal to $0$ on $\bigcup_{a=1}^k U_{{x_a}}$ which
is possible by the openness of transversality condition.
We put $\frak V_a = \frak V_{x_a}$,
$U_a = U_{x_a}$, $U^0_a = U^0_{x_a}$,
and
$(\frak s_{1}^{a,n},\dots,\frak s_{\ell_a}^{a,n})
= (\frak s^n_{x_a,1},\dots,\frak s^n_{x_a,\ell({x_a})})$.
Then we define
\begin{equation}
\frak U(K) = \bigcup_{a=1}^k U_{a}.
\end{equation}
We also take a relatively compact subset $U_{a}^0$
of $U_{a}$ for each $a$
such that
\begin{equation}
\frak U_0(K) = \bigcup_{a=1}^k U^0_{a} \supset K.
\end{equation}
We then fix an open neighborhood $U_0$ of $K$ to be
$U_0 = \frak U_0(K)$ .
\par
Let $\{\chi_0\} \cup \{\chi_{\frak r} \mid \frak r \in \frak R\}$
be a set of strongly smooth functions $\chi_* : \Omega^+ \to [0,1]$
satisfying the following properties.
\begin{proper}
\begin{enumerate}
\item
The support of $\chi_0$ is contained in $U_0$.
The support of $\chi_{\frak r}$
is contained in $U_{\frak r}$.
\item
$
\chi_0 + \sum_{\frak r \in \frak R} \chi_{\frak r} \equiv 1
$
on $\overline{\Omega}$.
\end{enumerate}
\end{proper}
\par
For each $x \in \overline{\Omega}$ we take $U_x$ with the following properties.

\begin{proper}\label{prop12888rev}
\begin{enumerate}
\item
If $x \in U_0 = \frak U_0(K)$  then $U_x \subset \frak U_0(K)$.
\item
If $x \in U_{\frak r}$ for $\frak r\in \frak R$, then
$U_x \subset U_{\frak r}$. Moreover
there exists $(h_{\frak r x},\tilde\varphi_{\frak r x},\breve\varphi_{\frak r x})$ as in Property \ref{proper728}.
\item
If $x \notin U_{\frak r}$ then $\overline U_x \cap {\rm Supp}(\chi_{\frak r}) =
\emptyset$.
\item
If  $x \in \frak U(K)$ then $U_x \subset U_{a}$
for some $a \in \{1,\dots,k\}$. Moreover
there exists $(h_{ax},\tilde\varphi_{ax},\breve\varphi_{ax}) : \frak V_x \to \frak V_a$ as in Property \ref{proper728}.
\end{enumerate}
\end{proper}
For each
$
\vec\xi = (\xi_{\frak r})_{\frak r\in \frak R} \in W_0
$
and $n   \in \Z_{\ge 0}$ we will define a multivalued perturbation
$\{\frak s^{n,\vec\xi}_{x}\}$ on $\frak V_x$ as follows.
\par\smallskip
\noindent (Case 1)
$x \in U_0 = \frak U_0(K)$  then $U_x \subset \frak U_0(K)$.
\par
Then $U_x \subset \frak U_0(K)$.
We take $a \in \{1,\dots,k\}$ such that $x \in U^{0}_{a}$.
By Property \ref{prop12888rev} (4) we can pullback
$(\frak s_{1}^{a,n},\dots,\frak s_{\ell_a}^{a,n})$ to $\frak V_{x}$.
This pullback is our $\frak s_{x}^{n,\vec\xi}$. (This is independent of $\vec\xi$.)
\par\smallskip
\par\smallskip
\noindent (Case 2) $x \in \frak U(K) \setminus \frak U_0(K)$.
\par
We take $a \in \{1,\dots k\}$ such that $U_x \subset U_{a}$.
By Property \ref{prop12888rev} (4)  we have
$(h_{ax},\tilde\varphi_{ax},\breve\varphi_{ax})$.
We put
\begin{equation}\label{formula123rev}
\frak R(x) = \{\frak r \in \frak R
\mid x \in {\rm Supp}(\chi_{\frak r})\}.
\end{equation}
We take $(h_{\frak rx},\tilde\varphi_{\frak rx},\breve\varphi_{\frak r x})$
as in Property \ref{prop12888rev}  (3) for each $\frak r \in \frak R(x)$.
We put
$$
I = \{1,\dots,\ell_a\}\times \prod_{\frak r \in \frak R(x)}\Gamma_{\frak r}.
$$
We define $\frak s^{n,\vec\xi}_{x,i}
: V_{x}  \to E_{x}$
for each $i \in I$
by  Formula (\ref{furmula124rev}) below.
We take a sequence $\epsilon_n > 0$ with $\lim_{n\to\infty}\epsilon_n = 0$
and fix it throughout the proof. (For example we may take $\epsilon_n = 1/n$.)
We put $i = (j,(\gamma_{\frak r}))$.
\begin{equation}\label{furmula124rev}
\aligned
\frak s^{n,\vec\xi}_{x,i}(y)
=
s_{x}(y) &+
\chi_{0}([y]) g_{a,y}^{-1}(\frak s^{n}_{a,j}(\tilde\varphi_{ax}(y)) - s_{a}(\tilde\varphi_{ax}(y)))
\\
&+\sum_{\frak r \in  \frak R(x)}
\chi_{\frak r}([y])
 g_{\frak r,y}^{-1}(
 \frak s_{\frak r}^{\epsilon_n}(\tilde\varphi_{\frak r x}(y),\gamma_{\frak r}^{-1}\xi_{\frak r})
 -  s_{\frak r}(\tilde\varphi_{\frak r x}(y))
\endaligned
\end{equation}
Explanation of the notations in Formula (\ref{furmula124rev}) is in order.
$s_{a} : V_{a} \to E_{a}$ is the representative of the Kuranishi map,
$g_{a,y} : E_{a} \to E_{x}$
is defined by
$
\breve{\varphi}_{ax}(y,\eta) = g_{a,y}(\eta),
$
$s_{\frak r} : V_{\frak r} \to E_{\frak r}$ is the representative of the Kuranishi map,
$g_{\frak r,y} : E_{\frak r} \to E_{x}$
is defined by
$
\breve{\varphi}_{\frak r x}(y,\eta) = g_{\frak r,y}(\eta).
$
$[y] \in V_a/\Gamma_a$ is the equivalence class of $y$, which we regard
as an element of $U$.
\par\smallskip
\noindent (Case 3) $x \in U \setminus \frak U(K)$.
\par
We define $\frak R(x)$ by (\ref{formula123rev}).
We take $(h_{\frak rx},\tilde\varphi_{\frak rx},\breve\varphi_{\frak rx})$
as in Property \ref{prop12888rev} (3) for each $x\in \frak R(x)$.
We put
$$
I = \prod_{\frak r \in \frak R(z)}\Gamma_{\frak r}.
$$
We define $\frak s^{n,\vec\xi}_{x,i}
: V_{x}  \to E_{x}$
for each $i = (\gamma_{\tau})_{\tau \in \frak R(x)}\in I$
by the following formula.
\begin{equation}\label{furmula124rev22}
\frak s^{n,\vec\xi}_{x,i}(y)
=
s_{x}(y) +
\sum_{\frak r \in  \frak R(x)}
\chi_{\frak r}([y])
 g_{\frak r,y}^{-1}(
 \frak s_{\frak r}^{\epsilon_n}(\tilde\varphi_{\frak r x}(y),\gamma_{\frak r}^{-1}\xi_{\frak r})
 -  s_{\frak r}(\tilde\varphi_{\frak r x}(y)).
 \end{equation}
Here the notations are the same as (\ref{furmula124rev}).
\begin{lem}\label{lem1227}
$(\frak s^{n,\vec\xi}_{x,i})_{i\in I}$ defines a multisection on $\frak V_{x}$.
Moreover it satisfies the following properties.
\begin{enumerate}
\item
In Case 1, it is independent of the choice of $a$ and $(h_{ax},\tilde\varphi_{ax},\breve\varphi_{ax})$.
\item
In Case 2, it is independent of the choice of 
$(h_{\frak r\frak z},\tilde\varphi_{\frak rx},\breve\varphi_{\frak rx})$,
and $a$, $(h_{ax},\tilde\varphi_{ax},\breve\varphi_{ax})$.
\item
In Case 3, it is independent of the choice of
$(h_{\frak rx},\tilde\varphi_{\frak rx},\breve\varphi_{\frak rx})$.
\item
The restriction of $(\frak s^{n,\vec\xi}_{x,i})_{i\in I}$ to $U_{x} \cap U_{x'}$
is equivalent to
restriction of $(\frak s^{n,\vec\xi}_{x',i})_{i\in I}$ to $U_{x} \cap U_{x'}$.
\end{enumerate}
\end{lem}
\begin{proof}
In Cases 2 and 3 we will show that $(\frak s^{\epsilon,\vec\xi}_{x,i}( y))_{i\in I}$
is a permutation of $(\gamma\frak s^{n,\vec\xi}_{x,i}(\gamma^{-1}y))_{i\in I}$ for $\gamma \in \Gamma_{x}$.
We calculate
\begin{equation}\label{form1222new}
\gamma g_{\frak r,\gamma^{-1}y}^{-1}(\xi_{\frak r})
=
\gamma \tilde\varphi_{\frak r x}(\gamma^{-1}y,\xi)
=
\tilde\varphi_{\frak r x}(y,\gamma\xi)
= g_{\frak r,y}^{-1}(\gamma\xi_{\frak r}).
\end{equation}
This implies that the third term of (\ref{furmula124rev}) and
the second term of (\ref{furmula124rev22})
is invariant under $\gamma$ action modulo permutation of the
indices in $I$.
The second term of (\ref{furmula124rev})
is invariant under $\gamma$ action modulo permutation of $1,\dots,a$
since $(\frak s_{a,1}^n,\dots,\frak s_{a,\ell_a}^n)$ is a multisection.
We have thus proved that $(\frak s^{n,\vec\xi}_{x,i})_{i\in I}$ is a
multisection.
(In Case 1 this fact is obvious.)
\par
Statement (1) is a consequence of the definition of multivalued perturbation, that is
the well-defined-ness of $\frak s_K$.
\par
To prove Statement (2), we observe that different choices of $(h_{\frak rx},\tilde\varphi_{\frak rx},\breve\varphi_{\frak rx})$ are related
one another by the action of $\gamma\in \Gamma_{\frak r}$.
(Lemma \ref{lem2715}.)
Then using (\ref{form1222new}) we can show that
the third  term of (\ref{furmula124rev}) changes by the permutation of $\gamma_{\frak r}$.
The first and the second terms of (\ref{furmula124rev}) do not change.
\par
By changing $(h_{ax},\tilde\varphi_{ax},\breve\varphi_{ax})$
the second term of (\ref{furmula124rev}) changes by the permutation of $j$
and the third  term of (\ref{furmula124rev}) does not change.
\par
The proof of (3) is easier.
\par
To prove (4) it suffices to consider the case $U_{x'} \subset U_{x}$.
If Case 1 is applied to both $x$ and $x'$ it is a consequence of the well-defined-ness
of the restriction of multisection.
\par
If Case 3 is applied to both, then $\frak R(x') \subset \frak R(x)$.
In this case the right hand side of (\ref{furmula124rev22})
is independent of $\Gamma_{\frak r}$ factor for $\frak r \notin \frak R(x')$
on $U_{\frak z'}$ .
Therefore the restriction of $(\frak s^{n,\vec\xi}_{x,i})_{i\in I}$ to $U_{x} \cap U_{x'}$ is a
permutation of the
$\prod_{\frak r \in \frak R(\frak z') \setminus \frak R(\frak z)} \# \Gamma_{\frak r}$ iteration of
the restriction of $(\frak s^{\epsilon,\vec\xi}_{x',i})_{i\in I}$ to $U_{\frak z} \cap U_{x'}$.
\par
When Case 2 is applied to both we can prove (4) by combining the above those two cases.
\par
What remains to prove is the case where the Case 2 is applied to $x$ and Case 1 or Case 3 is applied to
$x'$.
If Case 1 is applied to
$x'$ then $\chi_0$ becomes 1 on  $U_{x'}$. Therefore
the required equivalence follows from the well-defined-ness of the
restriction of the multisection.
If Case 3 is applied to $x'$ then $\chi_0$ becomes 0 on  $U_{x'}$.
Therefore the second term of (\ref{furmula124rev}) vanishes.
So the restriction of $(\frak s^{n,\vec\xi}_{x,i})_{i\in I}$ to $U_{x'}$ is a
permutation of the
$\ell_a\prod_{\frak r \in \frak R(x') \setminus \frak R(x)} \# \Gamma_{\frak r}$ iteration of
the restriction of $(\frak s^{n,\vec\xi}_{x',i})_{i\in I}$ to $U_{x'}$.
\end{proof}
\begin{lem}\label{lem1228}
For each sufficiently large $n$, the set of $\vec\xi$ such that
$(\frak s^{n,\vec\xi}_{x,i})_{i\in I}$ is transversal to $0$
for all $x \in \overline{\Omega}$ is
open and dense in $W_0$.
\par
If we are in Situation \ref{situ122} in addition, then the set of $\vec\xi$
such that $f$ is strongly transversal to $g$ with respect to $(\frak s^{n,\vec\xi}_{x,i})_{i\in I}$
is dense in $W_0$.
\end{lem}
\begin{proof}
It suffices to show the conclusion on a neighborhood  $U_x$ of each fixed $x \in \overline{\Omega}$, since
we can cover $\overline{\Omega}$ by countably many such $U_{x}$'s.
\par
In Case 1 this follows from the assumption that $\frak s_K$ is transversal to $0$.
\par
In Case 2
the set
$$
\{(y,\vec \xi) \mid \frak s^{n,\vec\xi}_{x,i}(y) = 0\}
$$
is a smooth submanifold of $V_{x} \times W_{0}$.
Therefore we can prove the lemma by applying
Sard's theorem to its projection to $W_{0}$.
(We use (\ref{form132222}) here.)
\par
We consider Case 3.
Suppose $\chi_0([x]) \ne 1$.
Then we can shrink the domain $U_{x}$ if necessary
and assume $\chi_{\frak r}([y]) \ne 0$
on $U_x$. Then by the same argument as Case 2
we can show the required transversality for the dense of
of $\xi_{\frak r}$.
\par
Suppose $\chi_0([x]) = 1$.
Then all the functions $\chi_{\frak r}$ together with its first derivative
is small in a neighborhood of $x$.
Moreover the first derivative of $\chi_0$
is small in a neighborhood of $x$.
We also remark that $\frak s_K$ is transversal to $0$ at $x$.
Therefore by using the openness of transversality,
we can shrink the neighborhood $U_x$ so that
$(\frak s^{n,\vec\xi}_{x,i})_{\in I}$ reminds to be transversal to $0$
on $U_x$ for any $\xi$ contained in a compact subset of $W_0$.
\par
The proof of the second half is similar.
In Case 2 above we consider
the set
$$
\{(y,\vec \xi,z)  \in V_x \times W \times N \mid \frak s^{n,\vec\xi}_{x,i}(y) = 0, \,\, f(y) = g(z)\}.
$$
This is a smooth submanifold of $V_x \times W \times N$.
Therefore we have the required transversality result in this case by applying Sard's theorem to the projection to $W$
from this manifold.
The rest of the proof is entirely the same.
\end{proof}
The proof of Proposition \ref{prop127777ver} is complete.
\end{proof}

\subsection{Compatible system of bundle extension data}
\label{subsection:bdlextcompa}
There is one nasty point in proving a multisection version of Proposition \ref{prop127777ver}.
(See Subsection \ref{subsec:nastyreason}.)
To go around it, we
will take the bundle extension data for each
coordinate change of the good coordinate system
so that they are compatible to one another.
This is closely related to the idea of compatible system
of tubular neighborhoods by Mather.
However in our case its construction is rather easy.

In this section, when we consider a bundle extension datum
$(\pi_{12},\tilde{\varphi}_{12},\Omega_{12},\Omega_1)$
of $(\Phi,\mathcal K)$,
we sometimes shrink $\Omega_{12}, \Omega_1$
and restrict $(\pi_{12},\tilde\varphi_{21})$ thereto respectively.
It will be a bundle extension datum if $\Omega_{12}$, $\Omega_1$
still remain to be neighborhoods of  $\varphi_{21}(K)$, $K$ respectively.
So we sometimes say $(\pi_{12},\tilde\varphi_{21})$ is a bundle extension data
without specifying $\Omega_{12}, \Omega_1$.

\begin{defn}\label{defn1232}
Let ${\widetriangle{\mathcal U}}$ be a
good coordinate system and $\mathcal K$  a
support system.
\par
We call  $(
\{\pi_{\frak q\frak p}\},\{\tilde\varphi_{\frak p\frak q}\},\{\Omega_{\frak q\frak p}\},\{\Omega_{\frak p}\})$
a {\it  system of bundle extension data}
\index{system of bundle extension data} of $({\widetriangle{\mathcal U}},\mathcal K)$
if they have the following properties.
(Here $\frak r \le \frak q \le \frak p$ are elements of $\frak P$.)
\begin{enumerate}
\item
$(\pi_{\frak q\frak p},\tilde\varphi_{\frak p\frak q},
\Omega_{\frak q\frak p},\Omega_{\frak p})$
is a bundle extension datum of
$(\Phi_{\frak p\frak q},\mathcal K_{\frak p\frak q})$,
where $\mathcal K_{\frak p\frak q}
= \varphi_{\frak p\frak q}^{-1}(\mathcal K_{\frak p}) \cap \mathcal K_{\frak q}$.
\item
$\pi_{\frak r\frak q} \circ \pi_{\frak q\frak p} = \pi_{\frak r\frak p}$
on a neighborhood of
$\varphi^{-1}_{\frak q\frak r}(\varphi^{-1}_{\frak p\frak q}(\mathcal K_{\frak p}))
\cap \varphi^{-1}_{\frak q\frak r}(\mathcal K_{\frak q})
\cap \varphi^{-1}_{\frak p\frak r}(\mathcal K_{\frak p})
\cap \mathcal K_{\frak r}$.
\item
$\tilde\varphi_{\frak p\frak q} \circ \tilde\varphi_{\frak q\frak r} = \tilde\varphi_{\frak p\frak r}$
on a neighborhood of
$\varphi^{-1}_{\frak q\frak r}(\varphi^{-1}_{\frak p\frak q}(\mathcal K_{\frak p}))
\cap \varphi^{-1}_{\frak q\frak r}(\mathcal K_{\frak q})
\cap \varphi^{-1}_{\frak p\frak r}(\mathcal K_{\frak p})
\cap \mathcal K_{\frak r}$.
\end{enumerate}
\end{defn}
The precise meaning of equality in Condition (3) is as follows.
We pull back $\tilde\varphi_{\frak q\frak r} : \pi_{\frak r\frak q}^*\mathcal E_{\frak r} \to \mathcal E_{\frak q}$
by $\pi_{\frak q\frak p}$ and obtain
$\pi_{\frak q\frak p}^*\tilde\varphi_{\frak q\frak r} : \pi_{\frak r\frak p}^*\mathcal E_{\frak r} \to \pi_{\frak q\frak p}^*\mathcal E_{\frak q}$.
We compose it with $\tilde\varphi_{\frak p\frak q}: \pi_{\frak q\frak p}^*\mathcal E_{\frak q} \to \mathcal E_{\frak p}$
and obtain a map $: \pi_{\frak r\frak p}^*\mathcal E_{\frak r} \to \mathcal E_{\frak p}$.
We denote it by $\tilde\varphi_{\frak p\frak q} \circ \tilde\varphi_{\frak q\frak r}$.
Condition (3)  requires that this map coincides with $\tilde\varphi_{\frak p\frak r}$.
\par
During the discussion of this section, we shrink $,\{\Omega_{\frak q\frak p}\},\{\Omega_{\frak p}\}$ several times
and restrict $\pi_{\frak q\frak p},\tilde\varphi_{\frak p\frak q}$ thereto.
We call
$(
\{\pi_{\frak q\frak p}\},\{\tilde\varphi_{\frak p\frak q}\})$ a compatible system of bundle extension data
sometimes, in case we do not need to specify the domain.
\begin{lem}
If $\pi_{\frak r\frak q}$, $\pi_{\frak q\frak p}$ are diffeomorphic
to the projections of normal bundles, then
the composition
$\pi_{\frak r\frak q}\circ \pi_{\frak q\frak p}$ is diffeomorphic
to the projections of normal bundle.
\end{lem}
The proof is easy and is omitted.
\begin{prop}\label{prop12333}
For any pair $({\widetriangle{\mathcal U}},\mathcal K)$
there exists a system of bundle extension data associated thereto.
\end{prop}
\begin{proof}
We first construct $\{\pi_{\frak q \frak p}\}$.
\begin{lem}\label{lem1226}
Let $\mathcal K$ be a support system of ${\widetriangle{\mathcal U}}$.
Then there exists $\Omega_{\frak q\frak p}$ for $\frak p > \frak q$,
$\Omega_{\frak p}$ and $\pi_{\frak q\frak p}$ with the following properties.
\begin{enumerate}
\item
$\Omega_{\frak q\frak p}$ is a neighborhood of $\varphi_{\frak p\frak q}
(\mathcal K_{\frak q} \cap U_{\frak p\frak q}) \cap \mathcal K_{\frak p}$
in $U_{\frak p}$.
\item
$\Omega_{\frak p}$ is a neighborhood of $\mathcal K_{\frak p}$ in
$U_{\frak p}$.
\item
$\pi_{\frak q\frak p} : \Omega_{\frak q\frak p} \to \Omega_{\frak q}$
is a continuous map which is diffeomorphic to the restriction of a
projection of a vector bundle to a neighborhood of $0$ section.
(See Definition \ref{lem123000}.)
\item
$\pi_{\frak q\frak p} \circ \varphi_{\frak p\frak q} = {\rm id}$ on a
neighborhood of $\varphi^{-1}_{\frak p\frak q}(\mathcal K_{\frak p})
\cap \mathcal K_{\frak q}$.
\item
$\pi_{\frak r\frak q} \circ \pi_{\frak q\frak p} = \pi_{\frak r\frak p}$
on a neighborhood of
$\varphi^{-1}_{\frak q\frak r}(\varphi^{-1}_{\frak p\frak q}(\mathcal K_{\frak p}))
\cap \varphi^{-1}_{\frak q\frak r}(\mathcal K_{\frak q})
\cap \varphi^{-1}_{\frak p\frak r}(\mathcal K_{\frak p})
\cap \mathcal K_{\frak r}$.
\end{enumerate}
\end{lem}
\begin{proof}
The proof is by induction on the number ${\rm differ}(\frak p,\frak q)$, which
we define now.
\begin{defn}
For $\frak p, \frak q \in \frak P$ with $\frak q < \frak p$
we put
$$
{\rm differ}(\frak p,\frak q)
=
\max \{ n \mid \exists \frak r_1,\dots,\frak r_n \in \frak P,
\,\, \frak q = \frak r_1 < \dots < \frak r_n = \frak p \}.
$$
\end{defn}
We will prove the following statement by induction on $n$.
\begin{enumerate}
\item[(*)]
The conclusion of Lemma \ref{lem1226} holds for
$\frak p,\frak q$ with ${\rm differ}(\frak p,\frak q) \le n$.
\end{enumerate}
Note we will shrink $\Omega_{\frak q\frak p}$, $\Omega_{\frak p}$
several times during the proof
and restrict $\pi_{\frak q\frak p}$ to the shrinked domain.
We however use the same symbol for the shrinked open sets and
retraction on it.
\par
If $n=0$ there is nothing to prove.
Next assuming (*) holds for all
$\frak p, \, \frak q$ with ${\rm differ}(\frak p,\frak q) < n$, we will prove the case of
pair $(\frak p,\frak q)$ with ${\rm differ}(\frak p,\frak q)=n$.
Let $\frak p,\frak q \in \frak P$ with $\frak p > \frak q$ and
${\rm differ}(\frak p,\frak q) = n$.
\par
For $\frak r \in \frak P$ with $\frak p > \frak r > \frak q$
we take a neighborhood
$
\Omega_{\frak q\frak r\frak p}
$
of
$$
\varphi^{-1}_{\frak r\frak q}(\varphi^{-1}_{\frak p\frak r}(\mathcal K_{\frak p}))
\cap \varphi^{-1}_{\frak r\frak q}(\mathcal K_{\frak r})
\cap \varphi^{-1}_{\frak p\frak q}(\mathcal K_{\frak p})
\cap \mathcal K_{\frak q}.
$$
The composition
$\pi_{\frak q\frak r} \circ \pi_{\frak r\frak p}$ is defined there.
\begin{sublem}
We may take
$\Omega_{\frak q\frak r\frak p}$ such that
$$
\pi_{\frak q\frak r} \circ \pi_{\frak r\frak p}
=
\pi_{\frak q\frak r'} \circ \pi_{\frak r'\frak p}
$$
holds on $\Omega_{\frak q\frak r\frak p} \cap \Omega_{\frak q\frak r'\frak p}$.
\end{sublem}
\begin{proof}
By Definition \ref{gcsystem} (5), we may shrink $\Omega_{\frak q\frak r\frak p}$
so that
$\Omega_{\frak q\frak r\frak p} \cap \Omega_{\frak q\frak r'\frak p}
\ne \emptyset$ implies either $\frak r < \frak r'$ or
$\frak r' > \frak r$.
We may assume $\frak r < \frak r'$ without loss of generality.
Then, by the induction hypothesis, we have
$$
\pi_{\frak q\frak r} \circ \pi_{\frak r\frak p}
=
\pi_{\frak q\frak r'} \circ \pi_{\frak r'\frak r} \circ \pi_{\frak r\frak p}
=
\pi_{\frak q\frak r'} \circ \pi_{\frak r'\frak p}
$$
on a neighborhood of
$$
\aligned
&\varphi^{-1}_{\frak r\frak q}(\varphi^{-1}_{\frak p\frak r}(\mathcal K_{\frak p}))
\cap \varphi^{-1}_{\frak r\frak q}(\mathcal K_{\frak r})
\cap \varphi^{-1}_{\frak p\frak q}(\mathcal K_{\frak p})
\cap \mathcal K_{\frak q}\\
&\cap
\varphi^{-1}_{\frak r'\frak q}(\varphi^{-1}_{\frak p\frak r'}(\mathcal K_{\frak p}))
\cap \varphi^{-1}_{\frak r'\frak q}(\mathcal K_{\frak r'})
\cap \varphi^{-1}_{\frak p\frak q}(\mathcal K_{\frak p})
\cap \mathcal K_{\frak q}.
\endaligned$$
because
$$
{\rm differ}(\frak r,\frak q), \, {\rm differ}(\frak p,\frak r), \, {\rm differ}(\frak p,\frak r'),\,
{\rm differ}(\frak r',\frak r), \, {\rm differ}(\frak r',\frak q),\, {\rm differ}(\frak p,\frak r') < n.
$$
The sublemma follows easily.
\end{proof}
Thus we can define $\pi_{\frak q \frak p}$ to be $\pi_{\frak q\frak r} \circ \pi_{\frak r\frak p}$
on a neighborhood of
\begin{equation}\label{form1212}
\varphi^{-1}_{\frak p\frak q}(\mathcal K_{\frak p})
\cap \mathcal K_{\frak q}
\cap
\bigcup_{\frak r}
\varphi^{-1}_{\frak r\frak q}(\varphi^{-1}_{\frak p\frak r}(\mathcal K_{\frak p}))
\cap \varphi^{-1}_{\frak r\frak q}(\mathcal K_{\frak r})
\end{equation}
which will then satisfy the required properties.
Then by shrinking the neighborhood of (\ref{form1212}) a bit,
we can use Proposition \ref{prop2949} to extend it to a neighborhood of
$\varphi^{-1}_{\frak p\frak q}(\mathcal K_{\frak p})
\cap \mathcal K_{\frak q}$.
Thus Lemma \ref{lem1226} is proved by induction on $n$.
\end{proof}
\begin{rem}
The proof of this lemma is easier than the proof of existence of
system of normal bundles of Mather since in our case stratification is
locally linear ordered by Definition \ref{gcsystem} (5).
\end{rem}
\begin{lem}\label{lem12343}
Let $(\{\Omega_{\frak q\frak p}\},\{\Omega_{\frak p}\},\{\pi_{\frak q\frak p}\})$
be as in Lemma \ref{lem1226}.
Then by shrinking $\Omega_{\frak q\frak p}$ and $\Omega_{\frak p}$ if necessary,
there exists $\tilde\varphi_{\frak p\frak q}$ such that
 $(
\{\pi_{\frak q\frak p}\},\{\tilde\varphi_{\frak p\frak q}\},\{\Omega_{\frak q\frak p}\},\{\Omega_{\frak p}\})$
becomes a
 system of bundle extension data of $({\widetriangle{\mathcal U}},\mathcal K)$.
\end{lem}
\begin{proof}
We can prove Lemma \ref{lem12343} by the same induction as Lemma \ref{lem1226}.
\end{proof}
The proof of Proposition \ref{prop12333} is complete.
\end{proof}
\begin{defn}\label{1219}
In the situation of Definition \ref{defn1230},
let $\frak s_1$ (resp. $\frak s_2$) be mutlisections
of $\mathcal U_1$ (resp. $\mathcal U_2$)
defined on a neighborhood of $K$ (resp. $\varphi_{21}(K)$).
\par
We say that $\frak s_1$ and $\frak s_2$ are {\it compatible with the bundle extension data }
$(\pi_{12},\tilde\varphi_{21})$
if
$\frak s_1$ and $\frak s_2$ are represented both by $\ell$ multisection and
there exists a permutation $\sigma : \{1,\dots,\ell\} \to \{1,\dots,\ell\}$
such that
\begin{equation}\label{form1226}
(y,\frak s_{2,i}(y))
=
\tilde\varphi_{21}(y,\frak s_{1,\sigma(i)}(\pi_{12}(y))),
\end{equation}
holds if $y$ is in a neighborhood of $\varphi_{21}(K)$ in $U_2$.
Here $\sigma$ depends on $\pi_{12}(y)$.
\end{defn}
\begin{rem}
The equality (\ref{form1226}) implies that  $\frak s_1$ and $\frak s_2$ are compatible with
the embedding $\Phi_{21}$ automatically.
\end{rem}
\begin{defn}
In the situation of Definition
\ref{defn1232}, let $\frak  s =\{\frak s_{\frak p}^{n}\}$ be a multivalued
perturbation of $({\widetriangle{\mathcal U}},\mathcal K)$.
\par
We say that $\frak  s$ is compatible with the
compatible system of bundle extension data $(
\{\pi_{\frak q\frak p}\},\{\tilde\varphi_{\frak p\frak q}\})$
if for each $\frak p > \frak q$ and $n \in \Z_{\ge 0}$, the pair
$\frak s_{\frak p}^{n}$, $\frak s_{\frak q}^{n}$
is compatible with $(\pi_{\frak q\frak p},\tilde\varphi_{\frak p\frak q})$
in the sense of Definition \ref{1219}.
\end{defn}
\subsection{Embedding of Kuranishi charts and extension of multisections}
\label{subsec:extmultisec}

We now consider the following situation.
\begin{shitu}\label{situ1219rev}
In the situation of Definition  \ref{1219},
let $K \subset X$ be a compact subset
contained in a relatively compact open subset $W \subset X$,
$\{\frak s_{K,2}^{n}\}$  a mutivalued perturbations of $\mathcal U_2$ on a neighborhood of $K$,
and $\{\frak s_1^{n}\}$  a mutivalued perturbations of $\mathcal U_1$ on a neighborhood of
$\overline W$.
\par
We assume that $\{\frak s_{K,2}^n\}$ and a restriction of $\{\frak s^{n}_1\}$ to a neighborhood of $K$ are
compatible with $\Phi_{21}$ and also compatible with the
restriction of the bundle extension data
$(\pi_{12},\tilde\varphi_{21})$ to a neighborhood of $K$.
$\blacksquare$
\end{shitu}
\begin{prop}\label{prop1220rev}
In Situation \ref{situ1219rev} there exists a multivalued perturbation $\{\frak s_2^{n}\}$  of
$\mathcal U_2$
on a neighborhood of $\overline W$
such that
\begin{enumerate}
\item The restriction of $\{\frak s_2^{n}\}$ to a neighborhood of $K$
coincides with $\{\frak s_{K,2}^n\}$ .
\item  $\{\frak s_{K,2}^n\}$ , $\{\frak s_1^{n}\}$ are compatible with $\Phi_{21}$
 on a neighborhood of $K$.
Moreover they are compatible with the bundle extension data
$(\pi_{12},\tilde\varphi_{21})$ in a neighborhood of $\overline W$.
\item  If  $\{\frak s_n^{K,2}\}$, $\{\frak s_1^{n}\}$ are transversal to $0$
for sufficiently large $n$ in addition then
$\frak s_2^{n}$ can be chosen to be
transversal to $0$ for sufficiently large $n$.
\item
Suppose that $g : N \to M$ is a smooth map between manifolds and
$f_1$ is strongly transversal to $g$ with respect to $\frak s_1^{n}$ for sufficiently
large $n$. Then we may choose $\frak s_2^{n}$ such that
$f_2$ is strongly transversal to $g$ with respect to $\frak s_2^{n}$ for sufficiently
large $n$.
\end{enumerate}
\end{prop}
\begin{proof}
We define $\frak s_{2,i}(\pi_{12}(y))$ by
\begin{equation}\label{form1226rev}
(y,\frak s_{2,i}(y))
=
\tilde\varphi_{21}(y,\frak s_{1,i}(\pi_{12}(y))),
\end{equation}
(Note this is Formula (\ref{form1226}) except we take $\sigma$ to be the identity.)
Statements (1) and (2) are obvious. We can shrink the neighborhood of $\overline W$ so that statements (3) and (4) hold.
\end{proof}
\begin{rem}
The proof of Proposition \ref{prop1220rev} is much simpler than that of
Proposition \ref{prop1221}. In fact in the proof of Proposition \ref{prop1221}
we use a bump function to glue two CF-perturbations.
Here Property (2) is automatic without using bump function since
we assumed $\{\frak s_{K,2}^n\}$ is compatible with bundle extension data.
\end{rem}
We need a few more definitions:
\begin{defn}\label{bundleextembedding}
Let $\widehat{\Phi} = (\{\Phi_{\frak p}\},\frak i) : {\widetriangle{\mathcal U}} \to {\widetriangle{\mathcal U^+}}$ be an
embedding of good coordinate systems $\mathcal K$,  $\mathcal K^+$ be their
support systems such that $\varphi_{\frak p}(\mathcal K_{\frak p}) \subset \mathcal K^+_{\frak i(\frak p)}$.
Let $\Xi =
(\{\pi_{\frak q\frak p}\},\{\tilde\varphi_{\frak p\frak q}\},\{\Omega_{\frak q\frak p}\},\{\Omega_{\frak p}\})$ and $\Xi^+ = (
\{\pi^+_{\frak q\frak p}\},\{\tilde\varphi^+_{\frak p\frak q}\},\{\Omega^+_{\frak q\frak p}\},\{\Omega^+_{\frak p}\})$  be systems of bundle extension data of
$({\widetriangle{\mathcal U}},\mathcal K)$
and $({\widetriangle{\mathcal U^+}},\mathcal K^+)$
respectively.
\par
A bundle extension datum of $(\widehat{\Phi},\mathcal K,\Xi,\Xi^+)$ consists of the objects
$\{(\tilde\varphi_{\frak p},\pi_{\frak p})\}$ that satisfy the following properties:
\begin{enumerate}
\item For each $\frak p \in \frak P$ we have a bundle extension datum
$(\tilde\varphi_{\frak p},\pi_{\frak p})$ of the embedding
$\Phi_{\frak p}$ on $\mathcal K_{\frak p}$.
\item
If $\frak q < \frak p$, then
\begin{enumerate}
\item
$\pi_{\frak q\frak p} \circ \pi_{\frak i(\frak p)}
=
\pi_{\frak q} \circ \pi^+_{\frak i(\frak q)\frak i(\frak p)}$
on a neighborhood of
$(\varphi^+_{\frak i(\frak p)\frak i(\frak q)}\circ \varphi_{\frak q})(\varphi^{-1}_{\frak p\frak q}(\mathcal K_{\frak p})
\cap \mathcal K_{\frak q})$.
\item
$\tilde\varphi_{\frak p} \circ \tilde\varphi_{\frak p\frak q} = \tilde\varphi^+_{\frak i(\frak p)\frak i(\frak q)} \circ \tilde\varphi_{\frak q}$
on a neighborhood of
$\varphi^{-1}_{\frak p\frak q}(\mathcal K_{\frak p})
\cap \mathcal K_{\frak q}$.
\end{enumerate}
Here $(\pi_{\frak q\frak p},\tilde\varphi_{\frak p\frak q})$
is a part of the bundle extension datum of $({\widetriangle{\mathcal U}},\mathcal K)$
and  $(\pi^+_{\frak i(\frak q)\frak i(\frak p)}, \tilde\varphi^+_{\frak i(\frak p)\frak i(\frak q)})$
is a part of the bundle extension datum of $({\widetriangle{\mathcal U^+}},\mathcal K^+)$.
(See Diagram \ref{diag33--}.)
\end{enumerate}
\end{defn}
\begin{defn}\label{def1325}
In the situation of Definition \ref{bundleextembedding}, let $\widetriangle{\frak s}$ (resp. $\widetriangle{\frak s^+}$) be a multivalued
perturbation of $({\widetriangle{\mathcal U}},\mathcal K)$ (resp. $({\widetriangle{\mathcal U^+}},\mathcal K^+)$),
such that $\widetriangle{\frak s}$ (resp. $\widetriangle{\frak s^+}$) is compatible with the
bundle extension data $\Xi$ (resp. $\Xi^+$).
\par
We say that $\widetriangle{\frak s}$, $\widetriangle{\frak s^+}$ are compatible with $\widehat{\Phi},
\{(\tilde\varphi_{\frak p},\pi_{\frak p})\}$
if for each $\frak p$, the pair $\frak s_{\frak p}$, $\frak s^+_{\frak i(p)}$
are compatible with
$(\tilde\varphi_{\frak p},\pi_{\frak p})$
in the sense of Definition \ref{1219}.
\end{defn}
\begin{proof}[Proof of Propositions \ref{prop621}]
We will work out the induction scheme of the proof of Proposition  \ref{existperturbcont}
for a multivalued perturbation compatible with the system of bundle extension data
as produced in Proposition \ref{prop12333}.
The detail is now in order.
\par
We write the bundle extension data we use by $\Xi$.
We use the notations in the proof of Proposition  \ref{existperturbcont}.
We replace Proposition \ref{existontiiindc1t} by the following.
\begin{prop}\label{existontiiindc1trev}
There exists a multivalued of perturbations
$\widetriangle{{\frak s}^{\frak F}}$ of
$({\widetriangle{\mathcal U}}(\frak F,\mathcal K^+),\mathcal K(\frak F))$
on $\frak T(\frak F,\mathcal K)$ with the following properties.
\begin{enumerate}
\item $\widetriangle{{\frak s}^{\frak F}}$ is compatible with the
bundle extension data which is a restriction of $\Xi$ to $({\widetriangle{\mathcal U}}(\frak F,\mathcal K^+),\mathcal K(\frak F))$.
\item $\widetriangle{{\frak s}^{\frak F}}$ is transversal to $0$.
\item
If $\widetriangle{f} : (X,Z;\widetriangle{\mathcal U}) \to M$ is weakly transversal to $g : N \to M$, then we may take
$\widetriangle{{\frak s}^{\frak F}}$ so that the restriction of $\widetriangle f$ is strongly transversal to $g$
with respect to $\widetriangle{{\frak s}^{\frak F}}$.
\end{enumerate}
\end{prop}
\begin{proof}
To prove Proposition \ref{existontiiindc1trev} we replace Lemma \ref{122222} by the following.
(We use the notation of Lemma \ref{122222}.)
\begin{lem}\label{122222rev}
For any ideal $\frak I \subseteq \frak F_-$, there exist an open neighborhood $U_{\frak p_0}(\frak I)$
of $\psi_{\frak p_0}^{-1}(\frak T(\frak I,\mathcal K)) \cap {\mathcal K}_{\frak p_0}$
 in $U_{\frak p_0}$ and a multivalued perturbation
$\widetriangle{\frak s^{\frak F}}(\frak I)$ of
$({\widetriangle{\mathcal U}}(\frak F,\mathcal K^+),\mathcal K(\frak F))$
on $\frak T(\frak I,\mathcal K)$
with the following properties.
\begin{enumerate}
\item $\widetriangle{\frak s^{\frak F}}(\frak I)$ is compatible  with the system of the bundle
extension data obtained by restricting $\Xi$.
\item
$\widetriangle{\frak s^{\frak F}}(\frak I)$, $\widetriangle{{\frak s}^{\frak F_-}}$
are compatible with the embedding
$\widehat\Phi_{\frak F\frak F_-;\frak I}$
and its bundle extension data in the sense of Definition \ref{def1325}.
\item
$\widetriangle{\frak s^{\frak F}}(\frak I)$ is transversal to $0$.
\item
If $\widetriangle f : (X,\widetriangle{\mathcal U}) \to M$ is weakly transversal to $g : N\to M$ then
we can choose  $\widetriangle{\frak s^{\frak F}}(\frak I)$ such that the restriction of $\widetriangle f$
is strongly transversal to $g$ with respect to it.
\end{enumerate}
\end{lem}
We remark that the embedding $\widehat\Phi_{\frak F\frak F_-;\frak I}$
is obtained by restricting the coordinate change
of ${\widetriangle{\mathcal U}}$. Therefore $\Xi$ induces a
bundle extension datum of each embedding of the Kuranishi charts
which makes up  $\widehat\Phi_{\frak F\frak F_-;\frak I}$.
The compatibility condition for $\widehat\Phi_{\frak F\frak F_-;\frak I}$, Definition \ref{bundleextembedding} (2),
is a consequence of the compatibility condition
for $\Xi$, Definition \ref{defn1232}.
We thus obtain a bundle extension data of embedding $\widehat\Phi_{\frak F\frak F_-;\frak I}$
that we mentioned in Item (2).
\begin{proof}[Proof of Lemma \ref{122222rev}]
The proof is the same as the proof of Proposition \ref{122222}.
Namely we replace Proposition \ref{prop1221} by Proposition \ref{prop1220rev}.
\end{proof}
We are now ready to complete the proof of Proposition \ref{existontiiindc1trev}.
In the proof of Proposition \ref{existontiiindc1t} we replace Lemma \ref{122222}
by Lemma \ref{122222rev}. We also replace Proposition  \ref{prop123123}
by Proposition \ref{prop127777ver}.
This proves Proposition \ref{existontiiindc1trev}.
\end{proof}
Theorem \ref{prop621} follows from Proposition \ref{existontiiindc1trev}.
\end{proof}
\subsection{Relative version of the existence of multisection}
\label{subsec:relexmulti}
We next prove a relative version of Theorem \ref{prop621}.
We need a relative version of Proposition \ref{prop12333}.
\begin{shitu}\label{sotu1248}
\begin{enumerate}
\item
Let ${\widetriangle{\mathcal U^{\mathcal Z_{(1)}}}}$ be a good coordinate system
of $X$ and $\mathcal Z_{(1)} \subset X$ is a compact subset.
Let $\mathcal K^{\mathcal Z_{(1)}}$ be a support system of
${\widetriangle{\mathcal U^{\mathcal Z_{(1)}}}}$ in the sense of Definition \ref{defn738}
(1).
\item
Let ${\widetriangle{\mathcal U}}$ be a good coordinate system of
$\mathcal Z_{(2)} \subset X$ with $\mathcal Z_{(1)} \subset {\rm Int}\,\mathcal Z_{(2)}$ and
suppose ${\widetriangle{\mathcal U^{\mathcal Z_{(1)}}}}$ strictly
extends to ${\widetriangle{\mathcal U}}$ in the sense of
Definition \ref{defn735f} (4).
Let $\mathcal K$ be a support system of ${\widetriangle{\mathcal U}}$
which extends $\mathcal K^{\mathcal Z_{(1)}}$ in the sense of
Definition \ref{defn738} (2).$\blacksquare$
\end{enumerate}
\end{shitu}
\begin{defn}
In Situation \ref{sotu1248} (1), we call  $(
\{\pi^{\mathcal Z_{(1)}}_{\frak q\frak p}\},\{\tilde\varphi^{\mathcal Z_{(1)}}_{\frak p\frak q}\},\{\Omega^{\mathcal Z_{(1)}}_{\frak q\frak p}\},\{\Omega^{\mathcal Z_{(1)}}_{\frak p}\})$
a {\it  system of bundle extension data} of $({\widetriangle{\mathcal U^{\mathcal Z_{(1)}}}},\mathcal K^{\mathcal Z_{(1)}})$
if they have the following properties.
\begin{enumerate}
\item
$(\pi^{\mathcal Z_{(1)}}_{\frak q\frak p},\tilde\varphi^{\mathcal Z_{(1)}}_{\frak p\frak q},
\Omega^{\mathcal Z_{(1)}}_{\frak q\frak p},\Omega^{\mathcal Z_{(1)}}_{\frak p})$
is a bundle extension data of
$(\Phi^{\mathcal Z_{(1)}}_{\frak p\frak q},\mathcal K^{\mathcal Z_{(1)}}_{\frak p\frak q})$,
where $\mathcal K^{\mathcal Z}_{\frak p\frak q}
= (\varphi^{\mathcal Z_{(1)}}_{21})^{-1}(\mathcal K^{\mathcal Z_{(1)}}_{\frak p}) \cap \mathcal K^{\mathcal Z_{(1)}}_{\frak q}$.
\item
$\pi^{\mathcal Z_{(1)}}_{\frak r\frak q} \circ \pi^{\mathcal Z_{(1)}}_{\frak q\frak p} = \pi^{\mathcal Z_{(1)}}_{\frak r\frak p}$
on a neighborhood of
$(\varphi^{\mathcal Z_{(1)}}_{\frak q\frak r})^{-1}(\varphi^{\mathcal Z_{(1)}}_{\frak p\frak q})^{-1}(\mathcal K^{\mathcal Z_{(1)}}_{\frak p}))
\cap (\varphi^{\mathcal Z_{(1)}}_{\frak q\frak r})^{-1}(\mathcal K^{\mathcal Z_{(1)}}_{\frak q})
\cap (\varphi^{\mathcal Z_{(1)}}_{\frak p\frak r})^{-1}(\mathcal K^{\mathcal Z_{(1)}}_{\frak p})
\cap \mathcal K^{\mathcal Z_{(1)}}_{\frak r}$.
\item
$\tilde\varphi^{\mathcal Z_{(1)}}_{\frak p\frak q} \circ \tilde\varphi^{\mathcal Z_{(1)}}_{\frak q\frak r} = \tilde\varphi^{\mathcal Z_{(1)}}_{\frak p\frak r}$
on a neighborhood of
$(\varphi^{\mathcal Z_{(1)}}_{\frak q\frak r})^{-1}((\varphi^{\mathcal Z_{(1)}}_{\frak p\frak q})^{-1}(\mathcal K^{\mathcal Z_{(1)}}_{\frak p}))
\cap (\varphi^{\mathcal Z_{(1)}}_{\frak q\frak r})^{-1}(\mathcal K^{\mathcal Z_{(1)}}_{\frak q})
\cap (\varphi^{\mathcal Z_{(1)}}_{\frak p\frak r})^{-1}(\mathcal K^{\mathcal Z_{(1)}}_{\frak p})
\cap \mathcal K^{\mathcal Z_{(1)}}_{\frak r}$.
\end{enumerate}
\end{defn}
\begin{prop}\label{prop1250}
 In Situation \ref{sotu1248} (1)+(2), let $\Xi^{\mathcal Z_{(1)}}$ be a bundle extension data of $({\widetriangle{\mathcal U^{\mathcal Z_{(1)}}}},\mathcal K^{\mathcal Z_{(1)}})$. Then there exists a bundle extension data $\Xi$ of $({\widetriangle{\mathcal U}},\mathcal Z_{(2)};\mathcal K)$
which coincides with $\Xi^{\mathcal Z_{(1)}}$ in a neighborhood of $\mathcal Z_{(1)}$.
\end{prop}
The proof is the same as the proof of
Proposition \ref{existontiiindc1trev}.
\begin{shitu}\label{1252situ}
In Situation \ref{sotu1248} (1) + (2), suppose we are given systems of bundle extension data $\Xi^{\mathcal Z_{(1)}}$ and
$\Xi$ respectively as in Proposition \ref{prop1250}.
\par
Let
$\widetriangle{\frak s^{\mathcal Z_{(1)}}}$ is a multivalued perturbation of $({\widetriangle{\mathcal U^{\mathcal Z_{(1)}}}},\mathcal K^{\mathcal Z_{(1)}})$. We assume that $\widetriangle{\frak s^{\mathcal Z_{(1)}}}$ is transversal to $0$. $\blacksquare$
\end{shitu}
\begin{prop}\label{existperturbmultires}
\par
In Situation \ref{1252situ}, there exists multivalued perturbation $\widetriangle{\frak s}$
of $({\widetriangle{\mathcal U}},\mathcal K)$ such that
\begin{enumerate}
\item $\widetriangle{\frak s}$ is transversal to $0$.
\item $\widetriangle{\frak s}$ coincides with $\frak s^{\mathcal Z_{(1)}}$ in a neighborhood of $\mathcal Z_{(1)}$.
\item  $\widetriangle{\frak s}$ is compatible with $\Xi$.
\item
If $\widetriangle{f} : (X,\mathcal Z_{(2)};\widetriangle{\mathcal U}) \to M$ is weakly transversal to $g : N \to M$
and the restriction of $\widetriangle{f}$ is strongly transversal to $g$ with respect to $\widetriangle{\frak s}$ then
we may choose $\widetriangle{\frak s}$ such that $\widetriangle{f}$
is strongly transversal to $g$ with respect to $\widetriangle{\frak s}$.
\end{enumerate}
\end{prop}
\begin{proof}
We can modify the proof of Proposition \ref{existperturbcontrel}  in exactly the same way
as we modified the proof of Proposition  \ref{existperturbcont}  to the proof of Theorem \ref{prop621}.
It proves Proposition \ref{existperturbmultires}.
\end{proof}

\subsection{Remark on the number of branches of extension of
multisection}
\label{subsec:nastyreason}
We now explain a certain delicate point we encounter when we try  to extend a
multisection or multivalued perturbation. (This point does not appear while
we extend a CF-perturbation.)
\par
Note that during the proof of Proposition \ref{prop1221} we used the same parameter space
while we extend our CF-perturbation on a subset of $U_1$ to its
neighborhood in $U_2$.
To extend a multivalued perturbation it is important
that we do not change the number of branches.
This point is mentioned in \cite[page 955 line 20-24]{FO} and
\cite[Remark 6.4]{foootech}.
We elaborate it below.
\par
Let $U_1 \subset U_2$ be an embedded submanifold.
Suppose $U_1$ is expressed as the union $U_1 = U_{1,1} \cup U_{1,2}$ of two
open subsets $U_{1,1}, \, U_{1,2}$ and
we are given multisection $\frak s$ on $U_1$.
We also assume that $\frak s\vert_{U_{1,j}}$ is extended
to its neighborhood in $U_2$ and
we denote  this extension by $\tilde{\frak s}_{j}$.
(We assume that the number of branches of $\tilde{\frak s}_{j}$
is the same as one of $\tilde{\frak s}$ around each point of $U_{1,j}$.)
We try to glue $\tilde{\frak s}_{1}$ and $\tilde{\frak s}_{2}$
to obtain an extension of $\frak s$ to a neighborhood of $U_1$ in $U_2$.
\par
Let $p \in U_{1,1} \cap U_{1,2}$. We represent $\frak s$ in a neighborhood of $p$
as $(s_1,\dots,s_{\ell})$, where $s_i$ are branches of $\frak s$ in a neighborhood of $p$.
We might say that extension $\tilde{\frak s}_{j}$ gives
$(s_{j,1},\dots,s_{j,\ell})$ and
we might try to  glue them as
\begin{equation}\label{1226form}
s_i(p) =
\chi(\pi(p)) s_{1,i}(p) + (1-\chi(\pi(p))) s_{2,i}(p)
\end{equation}
where $\chi$ is a function on $U_1$ such that
$\chi$ has a  support in $U_{1,1}$ and $1-\chi$ has
a support in $U_{1,2}$. $\pi$ is the projection of the
tubular neighborhood of $U_1$ in $U_2$.
\par
However, we can not define $s_i$ by (\ref{1226form}), because of the following problem.
As we explained in Remark \ref{rem614}, the way we take a representative  $(s_1,\dots,s_{\ell})$ of $\frak s$ is not unique.
(The uniqueness modulo permutation also fails.)
Therefore, when we extend $\frak s$ to $\tilde{\frak s}_{j}$, we may take
a representative different from $(s_1,\dots,s_{\ell})$ in a neighborhood of $p$.
So to add $s_{1,i}(p)$ and $s_{2,i}(p)$ may not make sense.
\par
Note the representative $(s_{j,1}(p),\dots,s_{j,\ell}(p))$ does
make sense at $p$ (modulo permutation). Namely it makes sense point-wise.
\par
We remark that  Definition \ref{1219} works point-wise.
Namely the right hand side of (\ref{form1226}) makes sense
point-wise and we do not need to take the representative of $\frak s_i^1$
{\it locally} but we only need to take its representative {\it at} $y = \pi(p)$.
This is because we fix a bundle extension datum which is well-defined
globally (and is compatible with coordinate changes etc. in the sense of
Definition \ref{defn1232}).
In fact this is {\it the} reason why we use the system of bundle extension data
(Definition \ref{defn1232}) to extend multivalued perturbations.
\par
We note that, in Formula (\ref{defext3030}), we do not use
system of bundle extension data to extend
 {\it CF-perturbations}.
\par
The definition of CF-perturbation is similar
to that of multivalued perturbation where
we use an open set of vector space in the former and
a finite set in the latter.
The  definition of equivalence of CF-perturbation in Definition \ref{conmultiequiv11}
is different from that of multivalued perturbation in Definition \ref{defn62}.
In fact the former is a local condition and the latter is a
point-wise condition.
\par
We can slightly modify the definition of  multivalued perturbation
by imitating the way taken in the case of CF-perturbation.
Then we do not need to use system of bundle extension data.
\par
In this document we use system of bundle extension data
since we want the definition of the multisection to be
exactly the same as one in \cite{FO}.
We want to do so since the definition of
\cite{FO} had been used by various people including ourselves.

\section{Zero and one dimensional cases via multisection}
\label{sec:onezerodim}

In Sections \ref{sec:contfamily} - \ref{sec:composition}, we discussed smooth correspondence and defined
virtual fundamental chain based on de Rham theory and
CF-perturbations.
In this section, we discuss another method based on multivalued perturbation.
Here we restrict ourselves to the case when the dimension of K-spaces of our interest is $1$, $0$ or negative,
and define virtual fundamental chain over $\Q$ in the $0$ dimensional case.
In spite of this restriction, the argument of this section is enough for the purpose,
for example, to prove all the results stated in \cite{FO}.
We recall that
in \cite{FO} we originally used a triangulation of the zero set of multisection
to define a virtual fundamental chain.
In this section we present a different way from \cite{FO}.
Namely, we use Morse theory in place of triangulation.
This change will make the relevant argument simpler and shorter for this restricted case.
We will explain the thorough detail about
the triangulation of the  zero set of multisection elsewhere.
\subsection{Statements of the results}
\label{subsec:zero1state}
We start with the following :
\begin{lem}\label{lem1311111}
Let ${\widetriangle{\mathcal U}}$ be a good coordinate system
(which may or may not have boundary or corner), $\mathcal K$
its support system,
and $\widetriangle{\frak s} = \{\frak s^{n}_{\frak p} \mid \frak p \in \frak P\}$
a multivalued perturbation of $({\widetriangle{\mathcal U}},\mathcal K)$
(Definition \ref{defn612}.)
We assume that
$\widetriangle{\frak s}$ is transversal to $0$.
Then there exists a natural number $n_0$ with the following properties.
\begin{enumerate}
\item
If
the dimension of $(X,{\widetriangle{\mathcal U}})$ is negative,
then
$
(\frak s_{\frak p}^{n})^{-1}(0) \cap \vert \mathcal K\vert= \emptyset
$
for $n \ge n_0$.
\item
If the dimension of $(X,{\widetriangle{\mathcal U}})$ is $0$,
then
$
(\frak s_{\frak p}^{n})^{-1}(0) \cap
\vert\partial{\widetriangle{\mathcal U}}\vert\cap \vert \mathcal K\vert
= \emptyset
$
for $n \ge n_0$.
Moreover
there exists a neighborhood $\frak U(X)$ of $X$ in
$\vert{\widetriangle{\mathcal U}}\vert \cap \vert \mathcal K\vert$
such that the intersection
$(\frak s_{\frak p}^{n})^{-1}(0) \cap \frak U(X)$ is a finite set for
any $n \ge n_0$.
\end{enumerate}
\end{lem}
\begin{proof}
(1) is obvious. Using the fact that the dimension of the boundary of
$(X,{\widetriangle{\mathcal U}})$ is negative, we have
$
(\frak s_{\frak p}^{\epsilon})^{-1}(0) \cap \vert\partial{\widetriangle{\mathcal U}}\vert
\cap \vert \mathcal K\vert
= \emptyset
$
in (2) also.
The finiteness of the order of the set $(\frak s_{\frak p}^{\epsilon})^{-1}(0) \cap \frak U(X)
\cap \vert \mathcal K\vert$
is a consequence of its compactness,
Corollary \ref{cor69}.
\end{proof}
We now consider:
\begin{shitu}\label{Situation132}
Let ${\widetriangle{\mathcal U}}$ be a good coordinate system
(which may or may not have boundary or corner),
and assume a support system
$\mathcal K$ thereof is given.
Let $\widetriangle{\frak s} = \{\frak s^{\epsilon}_{\frak p} \mid \frak p \in \frak P\}$ be
a multivalued perturbation of $({\widetriangle{\mathcal U}},\mathcal K)$.
We assume:
\begin{enumerate}
\item
$(X,{\widetriangle{\mathcal U}})$ is oriented,
\item
$\widetriangle{\frak s}$
is transversal to $0$.
\end{enumerate}
Consider another support system $\mathcal K'$
with
$\mathcal K' < \mathcal K$.
We take a neighborhood $\frak U(X)$ of $X$
as in Corollary \ref{cor69} for
$\mathcal K_2 = \mathcal K'$ and $\mathcal K_3 = \mathcal K$.
$\blacksquare$
\end{shitu}
\begin{defn}\label{gcswithperturb}
In Situation \ref{Situation132} we call $(\widetriangle{\mathcal U}, \widetriangle{\frak s})$
a {\it good coordinate system with multivalued perturbation}
\index{multivalued perturbation ! good coordinate system with multivalued perturbation}
\index{good coordinate system ! with multivalued perturbation} of $X$.
\end{defn}
\begin{defn}\label{defn13223}
In Situation \ref{Situation132}, we assume $\dim (X,{\widetriangle{\mathcal U}}) = 0$.
We consider
$p \in U_{\frak p} \cap \frak U(X) \cap  \vert \mathcal K\vert$
such that $\frak s_{\frak p}^{n}(p) = 0$.
(This means that there is a branch of
$\frak s_{\frak p}^{n}$ that vanishes at $p$.)
Let
$\frak V = (V,\Gamma,E,\psi,\hat\psi)$ be an orbifold chart of $(U_{\frak p},E_{\frak p})$
at $p$.
We take a representative $(\frak s_{\frak p,1}^{n},\dots,
\frak s_{\frak p,\ell}^{n})$ of $\frak s_{\frak p}^{n}$
on $\frak V$.
Let $\tilde p \in V$ such that $[\tilde p] = p$.
\begin{enumerate}
\item
For $i=1,\dots,\ell$ we put:
$$
\epsilon_{p,i} =
\begin{cases}
0  &  \text{if $\frak s_{\frak p,i}^{n}(\tilde p) \ne 0$.} \\
+1 & \text{if $\frak s_{\frak p,i}^{n}(\tilde p) = 0$
and (\ref{iso1311111}) below is orientation preserving.}\\
-1 & \text{if $\frak s_{\frak p,i}^{n}(\tilde p) = 0$
and (\ref{iso1311111}) below is orientation reversing.}
\end{cases}
$$
In the current case of virtual dimension 0, the transversality hypothesis
implies that the derivative
\begin{equation}\label{iso1311111}
D_{\tilde p} \frak s_{\frak p,i}^{n} :
T_{\tilde p}V \to E_{\tilde p}
\end{equation}
becomes an isomorphism at every point $\tilde p$ satisfying $\frak s_{\frak p,i}^{n}(\tilde p) = 0$.
\item
The {\it multiplicity} $m_p$ of
$(\frak s_{\frak p}^{n})^{-1}(0)$ at $p$ is a rational number
and is defined by
$$
m_p = \frac{1}{\ell\#\Gamma}\sum_{i=1}^{\ell} \epsilon_{p,i}.
$$
\end{enumerate}
\end{defn}
\begin{lem}\label{lem134}
\begin{enumerate}
\item
The multiplicity $m_p$ in Definition \ref{defn13223}
is independent of the choice of representative
$(\frak s_{\frak p,1}^{n},\dots,
\frak s_{\frak p,\ell}^{n})$.
\item
If $q \in U_{\frak p\frak q}$ and $p = \varphi_{\frak p\frak q}(q)$,
then the multiplicity at $p$ is equal to the multiplicity at $q$.
\end{enumerate}
\end{lem}
\begin{proof}
This is immediate from the definition.
\end{proof}
\begin{defn}\label{defn1355}
In Situation \ref{Situation132} we assume $\dim (X,{\widetriangle{\mathcal U}}) = 0$.
We define the {\it virtual fundamental chain}
$[(X,{\widetriangle{\mathcal U}},\mathcal K',\widetriangle{\frak s^{n}})]$ of $(X,{\widetriangle{\mathcal U}},\mathcal K',\widetriangle{\frak s^{n}})$
by
\begin{equation}\label{defvfcdim0}
[(X,{\widetriangle{\mathcal U}},\mathcal K',\widetriangle{\frak s^{n}})]
=
\sum_{p \in \frak U(X) \cap \vert\mathcal K'\vert
\cap \bigcup_{\frak p \in \frak P}(\frak s_{\frak p}^{n})^{-1}(0)} m_p.
\end{equation}
This is a rational number.
Here the sum in the right hand side of (\ref{defvfcdim0}) is defined as follows.
Let us consider the {\it disjoint} union
$$
\bigcup_{\frak p \in \frak P}
(\frak U(X) \cap \mathcal K'_{\frak p}\cap (\frak s_{\frak p}^{n})^{-1}(0))
\times \{\frak p\}.
$$
We define a relation $\sim$ on it
by $(p,\frak p) \sim (q,\frak q)$
if $\frak p \le \frak q$, $q  = \varphi_{\frak q\frak p}(p)$
or $\frak q \le \frak p$, $p  = \varphi_{\frak p\frak q}(q)$.
This is an equivalence relation by
Definition \ref{gcsystem} (7).
The set of the equivalence classes
is denoted by
$\frak U(X) \cap \vert\mathcal K'\vert\cap \bigcup_{\frak p \in \frak P}(\frak s_{\frak p}^{n})^{-1}(0)$.
By Lemma \ref{lem134} (2) the multiplicity
$m_p$ is a well-defined function on this set.
\end{defn}

\begin{rem}
We note that in case $(X,{\widetriangle{\mathcal U}})$
has a boundary, the number
$[(X,{\widetriangle{\mathcal U}},\mathcal K',\widetriangle{\frak s^{n}})]$
depends on the choice of the multivalued perturbation $\widetriangle{\frak s} =\{\frak s^{n}_{\frak p} \}$.
It also depends on $n$.
\end{rem}
The next result is a multivalued perturbation version of Proposition \ref{indepofukuracont}.
\begin{prop}\label{prop14777}
Let $\mathcal K^1,\mathcal K^2, \mathcal K^3$ be support systems
with
$\mathcal K^1 < \mathcal K^2< \mathcal K^3 = \mathcal K$.
Let $X, \widetriangle{\mathcal U}, \widetriangle{\frak s} = \{\widetriangle{\frak s^{n}}\}$ and $\mathcal K$
be as in Definition \ref{defn1355}.
\begin{enumerate}
\item
The number $[(X,{\widetriangle{\mathcal U}},\mathcal K',\widetriangle{\frak s^{n}})]$ in
(\ref{defvfcdim0}) is independent of $\frak U(X)$ for all sufficiently large $n$.
Here $\mathcal K'$ is either $\mathcal K_1$ or $\mathcal K_2$.
\item
We have
$$
[(X,{\widetriangle{\mathcal U}},\mathcal K_1,\widetriangle{\frak s^{n}})]
=
[(X,{\widetriangle{\mathcal U}},\mathcal K_2,\widetriangle{\frak s^{n}})]
$$
for all sufficiently large $n$.
\end{enumerate}
\end{prop}
\begin{proof}
(1)
Let $\frak U'(X)$ be an alternative choice.
By Corollary \ref{cor69} we have
$$
\frak U(X) \cap \vert\mathcal K'\vert
\cap \bigcup_{\frak p \in \frak P}(\frak s_{\frak p}^{n})^{-1}(0)
=
\frak U'(X) \cap \vert\mathcal K'\vert
\cap \bigcup_{\frak p \in \frak P}(\frak s_{\frak p}^{n})^{-1}(0).
$$
Independence of the multiplicity $m_p$ can be proved
in the same way as the argument in Definition \ref{defn1355}.
(1) follows.
\par
(2) By Proposition \ref{lem715} we have
$$
\frak U(X) \cap \vert\mathcal K_1\vert
\cap \bigcup_{\frak p \in \frak P}(\frak s_{\frak p}^{n})^{-1}(0)
=
\frak U(X) \cap \vert\mathcal K_2\vert
\cap \bigcup_{\frak p \in \frak P}(\frak s_{\frak p}^{n})^{-1}(0).
$$
Independence of the multiplicity $m_p$ can be proved
in the same way as the argument in Definition \ref{defn1355}.
(2) follows.
\end{proof}
\begin{conven}
Since $[(X,{\widetriangle{\mathcal U}},\mathcal K',\widetriangle{\frak s^{n}})]$
is independent of $\mathcal K'$ by Proposition \ref{prop14777}, we will write it as
$[(X,{\widetriangle{\mathcal U}},\widetriangle{\frak s^{n}})]$
hereafter.
\end{conven}

The main result of this section is the following.

\begin{thm}\label{prop13777}
Let $(X,{\widetriangle{\mathcal U}},\widetriangle{\frak s})$ be as in
Situation \ref{Situation132}. We assume
$\dim (X,{\widetriangle{\mathcal U}}) = 1$.
We consider its normalized boundary
$\partial(X,{\widetriangle{\mathcal U}}) =
(\partial X,\partial{\widetriangle{\mathcal U}})$
where $\widetriangle{\frak s}$ induces a multivalued perturbation
$\widetriangle{\frak s_{\partial}}$ thereof and
$(\partial X,\partial{\widetriangle{\mathcal U}},\widetriangle{\frak s_{\partial}^{n}})$
is as in Situation \ref{Situation132} with
$\dim (\partial X,\partial{\widetriangle{\mathcal U}}) = 0$.
Then the following formula holds.
$$
[(\partial X,\partial{\widetriangle{\mathcal U}},\widetriangle{\frak s_{\partial}^{n}})]
=0.
$$
\end{thm}
\begin{rem}\label{rem1820moved}
Here we remark a slightly delicate point about the definition of
transversality of the multivalued perturbation
$\widetriangle{\frak s} = \{\widetriangle{\frak s^{n}}\}$.
\par
We remark that for a CF-perturbation we studied
a family parameterized by $\epsilon$ which is a positive real  number
close to $0$.
For a multivalued perturbation we considered
a sequence of multisections $\widetriangle{\frak s^{n}}$,
where $n$ is an integer.
In other words in the case of CF-perturbation
the parameter space is uncountable while
in the case of multivalued perturbation
the parameter space is countable.
\par
In many parts of the story of multivalued perturbation
we can consider $\widetriangle{\frak s} = \{\widetriangle{\frak s^{\epsilon}}\}$
 in place of $\{\widetriangle{\frak s^{n}}\}$.
However we need a countable family of objects to apply
Baire's category theorem (See the end of the proof of Lemma \ref{lem1313}).
\par
When we discuss transversality of $\widetriangle{\frak s} = \{\widetriangle{\frak s^{\epsilon}}\}$
we may consider one of the following two versions:
\begin{enumerate}
\item
Fix sufficiently small $\epsilon > 0$ and define
the transversality of  $\widetriangle{\frak s^{\epsilon}}$ as a multisection.
\item
We consider the whole family $\widetriangle{\frak s}$ as a multisection on
$(X,{\widetriangle{\mathcal U}}) \times (0,\epsilon_0)$.
\end{enumerate}
Sard's theorem implies that if $\widetriangle{\frak s} = \{\widetriangle{\frak s^{\epsilon}}\}$ is transversal to $0$
in the sense of (2) then {\it for generic $\epsilon$} the multisection
$\widetriangle{\frak s^{\epsilon}}$ is transversal to $0$
in the sense of (1).
\par
The transversality we  need to define virtual fundamental chain
is one in the sense of (1).
\end{rem}
\begin{rem}
The point we elaborate in Remark \ref{rem1820moved} is related
to the $n$-dependence of the virtual fundamental chain
$[(X,{\widetriangle{\mathcal U}},
\widetriangle{\frak s^{n}})]$
of a $0$ dimensional good coordinate system as follows.
\par
This phenomenon occurs only in case when $(X,{\widetriangle{\mathcal U}})$
has a boundary.
\par
Suppose $\widetriangle{\frak s}$ is transversal to $0$ in the sense of
(2) above.
Then for sufficiently small generic $\epsilon$,
$\widetriangle{\frak s^{\epsilon}}$ is transversal to $0$ at $\epsilon$.
So we can define
the rational number
$[(X,{\widetriangle{\mathcal U}},\widetriangle{\frak s^{\epsilon}})]$.
\par
On the other hand, there is a discrete subset
$S \subset (0,\epsilon_0)$ such that if $\epsilon_1 \in S$
then $\widetriangle{\frak s^{\epsilon_1}}$ may not be transversal to $0$.
In particular the zero set of $\widetriangle{\frak s^{\epsilon_1}}$ may
intersect with the boundary
$\partial(X,{\widetriangle{\mathcal U}})$.
Therefore
$$
\lim_{\epsilon \uparrow \epsilon_1} [(X,{\widetriangle{\mathcal U}},
\widetriangle{\frak s^{\epsilon}})]
\ne
\lim_{\epsilon \downarrow \epsilon_1}
[(X,{\widetriangle{\mathcal U}},\widetriangle{\frak s^{\epsilon}})],
$$
in general.
In other words, wall crossing may occur at $\epsilon_1$.
By this reason the virtual fundamental
chain $ [(X,{\widetriangle{\mathcal U}},\widetriangle{\frak s^{\epsilon}})]$ depends on $\epsilon$.
\par
We remark that in the case when we use CF-perturbations
a similar phenomenon happens and the integration along the
fiber depends on $\epsilon$.
However this phenomenon appears there in a slightly different way.
Suppose we consider the zero dimensional case.
Then the integration along the fiber (of the function $1$
and the map $X \to $ point)
$
[(X,{\widetriangle{\mathcal U}},\widetriangle{\frak S^{\epsilon}})]
$
defines a real number for each $\epsilon$.
This is defined for all sufficiently small $\epsilon$ but is
not a constant function of $\epsilon$.
\par
The number we obtain using the CF-perturbation
is a real number and so changes continuously as $\epsilon$ varies.
(It is easy to see from definition that  $[(X,{\widetriangle{\mathcal U}},
\widetriangle{\frak s^{\epsilon}})]$
is a smooth function of $\epsilon$.)
The number we obtain using multivalued perturbation
is a rational number. So it jumps. In other words
it is impossible to obtain virtual fundamental
chain $[(X,{\widetriangle{\mathcal U}},\widetriangle{\frak s^{\epsilon}})],$
for all $\epsilon$ in case $\widetriangle{\frak s}$ is a
multivalued perturbations.
\end{rem}

\subsection{A simple Morse theory on K-space}
\label{subsec:kuramorse}

\begin{proof}[Proof of Theorem \ref{prop13777}]
We can prove this theorem by taking an appropriate perturbation
$\widetriangle{\frak s^{n}}$
so that its zero set has a triangulation.
This is the proof given in \cite[Theorem 6.2]{FO}.
(Theorem \ref{prop13777} is a special case of \cite[Theorem 6.2]{FO}
where $Y$ is a point.)
Here we give an alternative proof without using triangulation.
\par
Let $(\mathcal K',\mathcal K)$ be a support pair of ${\widetriangle{\mathcal U}}$.
Note $\widetriangle{\frak s}$ is a multivalued perturbation of $({\widetriangle{\mathcal U}},\mathcal K)$.
We consider a strongly smooth function
$f : U(\vert\mathcal K\vert) \to [0,\infty)$
of a neighborhood $U(\vert\mathcal K\vert)$ of
$\vert\mathcal K\vert$ in $\vert{\widetriangle{\mathcal U}}\vert$
such that
\begin{equation}\label{133formula}
f^{-1}(0) = \vert\partial{\widetriangle{\mathcal U}}\vert
\cap U(\vert\mathcal K\vert).
\end{equation}
\begin{defn}
We say that $f$ is {\it normally positive at the boundary}
\index{normally positive function at the boundary}
if the following holds.
\par
Let $p \in \partial U_{\frak p} \cap U(\vert\mathcal K\vert)$ and
we identify its neighborhood with
$(W \times [0,1)^k)/\Gamma$ where $p$ corresponds to
$(p_0,(0,\dots,0))$ and $p_0 \in W$ is an interior point.
Let $\vec v \in T_{(p_0,0,\dots,0)}(W \times [0,1)^k)$
such that $\vec v = (\vec v_0,v_1,\dots,v_k)$ with
$v_i > 0$ for all $i=1,\dots,k$.
Then
\begin{equation}\label{normallypositive}
\vec v(f) > 0.
\end{equation}
\end{defn}
\begin{rem} The above definition can be rephrased as follows.
Consider the conormal space
$$
N^*_{(p_0,0,\dots,0)} (W \times \{0\}) \subset T^*_{(p_0,0,\dots,0)} (W \times [0,1]^k).
$$
Then $f$ is normally positive at $(p_0,0,\dots,0)$ if and only if
$df(p_0,0,\dots,0)$ is contained in the cone
$$
\aligned
C_+(N^*_{(p_0,0,\dots,0)} (W \times \{0\})) & := \{\alpha \in N^*_{(p_0,0,\dots,0)} (W \times \{0\})
\mid \alpha(\vec v) > 0, \, \vec v = (\vec v_0,v_1,\dots,v_k) \\
& \text{ with },
v_i > 0 \text{ for all }\, i=1,\dots,k\}.
\endaligned
$$
\end{rem}
\begin{lem}\label{139lem}
There exists a strongly smooth function
$f$ as above
which is normally positive at the boundary and
satisfies (\ref{133formula}).
\end{lem}
\begin{proof}
We take $\mathcal K^+$ such that
$(\mathcal K,\mathcal K^+)$ is a support pair.
Let $x \in \vert\mathcal K\vert$.
We then take a maximal $\frak p$ such that
$x \in \mathcal K_{\frak p}$.
We take a sufficiently small neighborhood $\Omega_x$
of $x$ in $\vert{\widetriangle{\mathcal U}}\vert$
such that $\Omega_x \cap \mathcal K_{\frak q} \ne \emptyset$
implies $\frak q \le \frak p$.
Then we may slightly shrink  $\Omega_x$ so that
$\Omega_x \cap \vert\mathcal K\vert$
is contained in $\Omega_x \cap \mathcal K^+_{\frak p}$.
Note $\Omega_x \cap \mathcal K^+_{\frak p}$ is an orbifold
with corner.
Therefore a neighborhood of $x$ in it is identified with
a point $(x_0,(0,\dots,0))$ in
$(V_x \times [0,1)^{k_x})/\Gamma_x$.
Here $x \in \overset{\circ}{S}_{k_x}(\mathcal K^+_{\frak p})$.
We can choose this coordinate so that the $\Gamma_x$
action on $V_x \times [0,1)^{k_x}$ is given in the form as
$$
(y,(t_1,\dots,t_{k_x})) \mapsto
(\gamma(y,(t_1,\dots,t_{k_x})),(t_{\sigma(1)},\dots,t_{\sigma(k_x)})),
$$
that is, the action on $[0,1)^{k_x}$ factor is by permuting its components.
We define a function $f_x$ on $\Omega_x \cap \mathcal K^+_{\frak p}$
by
\begin{equation}\label{cornertrans}
f_x(y,(t_1,\dots,t_{k_x}))
=
\begin{cases}
 t_1  t_2  \cdots  t_{k_x} &\text{if $k_x >0$,} \\
1 &  \text{if $k_x =0$.}
\end{cases}
\end{equation}
We may regard it as a strongly smooth function on a neighborhood of
$W_x = \Omega_x \cap \vert\mathcal K\vert$.
We may assume that if
$x \in \overset{\circ}{S}_{k_x}(\mathcal K^+_{\frak p})$
then
$W_x \cap \overset{\circ}{S}_{k_x+1}(\mathcal K_{\frak p})
= \emptyset$ for any $\frak p$.
Let $W_{0,x}$ be a relatively compact neighborhood of $x$ in
$W_x$.
\par
We take finitely many points $x_i$, $i=1,\dots,N$
of $\vert\mathcal K\vert$ such that
$$
\bigcup_{i=1}^N W_{0,x_i}
\supseteq \vert\mathcal K\vert.
$$
Then there exists a strongly smooth functions $\chi_i$ on a neighborhood of
$\vert\mathcal K\vert$ to $[0,1]$ such that
\begin{enumerate}
\item
The support of $\chi_i$ is in $W_{0,x_i}$.
\item
$
\sum_{i=1}^N \chi_i \equiv 1.
$
\end{enumerate}
We can prove the existence of such $\chi_i$ in the same way as in the proof of
Proposition \ref{pounitexi}
by using Lemma \ref{bumpfunctionlemma}.
We put
$$
f = \sum_{i=1}^N \chi_i f_{x_i}.
$$
By (\ref{cornertrans}) this function has the required properties.
\end{proof}
Let $\mathcal K''$ be a support system with
$\mathcal K' < \mathcal K'' < \mathcal K$.
\begin{defn}
Let $f$ be a strongly smooth function defined on a neighborhood of
$\vert\mathcal K\vert$.
We say that
a point $p \in \vert\mathcal K\vert$ is a {\it critical point}  of $f$ if
there exists $\frak p$ such that $p \in \mathcal K_{\frak p}$
and $p$ is a critical point
of the restriction of $f$ to a neighborhood of $p$ in
$U_{\frak p}$.
We denote by ${\rm Crit}(f)$ the set of all
critical points of $f$.
\end{defn}
\begin{defn}
Let $f$ be a strongly smooth function defined on a neighborhood of
$\vert\mathcal K\vert$.
We say that $f$ is a {\it Morse function} \index{Morse function}
if the restriction of $f$ to each $U_{\frak p}$
is a Morse function. (Namely its Hessian at all the critical points
are nondegenerate.)
\end{defn}
\begin{rem}
Suppose ${f} : \vert\mathcal K\vert \to \R$ is a
strongly smooth function and $\mathcal K'$ is a support system with $\mathcal K' < \mathcal K$,
then the restriction ${f}\vert_{\vert\mathcal K'\vert} : \vert\mathcal K'\vert \to \R$ is
strongly smooth. If ${f}$ is Morse so is ${f}\vert_{\vert\mathcal K'\vert}$.
\end{rem}
We use the following:
\begin{lem}\label{lem131111}
There exist finitely many orbifold charts
$(V_i,\Gamma_i,\psi_i)$ of $U_{\frak p_i}$
($\frak p_i \in \frak P$)
and smooth embeddings
$
h_i : [0,1] \to V_{i}
$
such that
\begin{equation}
\frak U(X) \cap
\bigcup_{\frak p\in \frak P}
(\frak s_{\frak p}^{n})^{-1}(0) \cap \mathcal K''_{\frak p}
\subset
\bigcup_{i=1}^N [(\psi_i\circ h_i)((0,1))].
\end{equation}
Here $\frak U(X)$ is as in Corollary \ref{cor69}.
\end{lem}
\begin{proof}
Since all the branches of $\frak s^n_{\frak p}$ are transversal to $0$
and the (virtual) dimension is $1$, locally the zero set of $\frak s^n_{\frak p}$
is a one dimensional manifold. The lemma then follows from compactness of the left hand side.
(Corollary \ref{cor69}.)
\end{proof}
\begin{prop}\label{morseufnction}
Suppose that $\dim U_{\frak p}\ge 2$ for each $\frak p$.
Then there exists a strongly smooth function
$f$ on a neighborhood $\frak U(X)$ of $X$ in $\vert\mathcal K\vert$
such that
\begin{enumerate}
\item
$f$ is normally positive at the boundary.
\item
$f$ satisfies (\ref{133formula}).
\item
$f$ is a Morse function.
\item
The composition $f\circ \psi_i \circ h_i : (0,1) \to \R$ is a Morse
function for each $
h_i : [0,1] \to V_{i}
$ in Lemma \ref{lem131111}
\end{enumerate}
\end{prop}
\begin{proof}
The proof of the lemma is a minor modification of a standard argument. We give a proof below for
completeness' sake.
We will use certain results concerning denseness of the set of Morse functions,
which we will prove in Subsection \ref{subsec:moresedense}.
\par
We make a choice of Fr\'echet space we work with.
Take a support system $\mathcal K^+$ such that
$(\mathcal K,\mathcal K^+)$ is a support pair and define the set
\begin{equation}
C^{\infty}(\mathcal K^+)
=
\{
(f_{\frak p})_{\frak p \in \frak P}
\in \prod_{\frak p \in \frak P}C^{\infty}{(\mathcal K^+_{\frak p})}
\mid
f_{\frak p} \circ \varphi_{\frak p\frak q}
=f_{\frak q},
\text{on $\varphi_{\frak p\frak q}^{-1}(\mathcal K^+_{\frak p})
\cap \mathcal K^+_{\frak q}$}
\}.
\end{equation}
Here $C^{\infty}(\mathcal K^+_{\frak p})$ is the space of
$C^{\infty}$ functions of an orbifold $U_{\frak p}$
defined on its compact subset $\mathcal K^+_{\frak p}$
and so is a Fr\'echet space with respect to the
$C^{\infty}$ topology. Then the set
$C^{\infty}(\mathcal K^+)$ is a
closed subspace of a finite product of the Fr\'echet spaces
and so is a Fr\'echet space.
\par
Let $f_0 \in C^{\infty}(\mathcal K^+_{\frak p})$ be a function satisfying
(1)(2) above. (Existence of such $f_0$ follows from
Lemma \ref{139lem}.)
We can take another neighborhood $\frak U_0$ of $\partial X$
in $\vert\mathcal K\vert$
such that $\text{\rm Crit}f_0 \cap \frak U_0 = \emptyset$.
We take a neighborhood $\frak U_1$ of $\partial X$
in $\vert\mathcal K\vert$ such that
$\overline{\frak U_1} \subset \frak U_0$.
Let $C^{\infty}(\mathcal K^+)_0$
be the set of all $f \in C^{\infty}(\mathcal K^+)$
that vanish on $\overline{\frak U_1}$.
Then $C^{\infty}(\mathcal K^+)_0$ itself is a Fr\'echet space.
\begin{lem}\label{lem1313}
For each $p \in X \setminus \frak U_0$
there exists a compact neighborhood
$\frak U_p$ of $p$ such that the set of $g\in C^{\infty}(\mathcal K^+)_0$
satisfying conditions (1)(2) below is a dense subset of
$C^{\infty}(\mathcal K^+)_0$.
\begin{enumerate}
\item
$f_0+g$
is a Morse function on $ \frak U_p$.
\item
The composition $(f_0+g)\circ \psi_i \circ h_i : (0,1) \to \R$ is a Morse
function on $h_i^{-1}(\psi_i^{-1}(\frak U_p))$
for each $
h_i : [0,1] \to V_{i}
$ in Lemma \ref{lem131111}
\end{enumerate}
\end{lem}
\begin{proof}
Let $\frak p \in \frak P$ be the maximal
element of $\{\frak p \mid p \in \vert\mathcal K_{\frak p}\vert\}$.
(Maximal element exists because of Definition \ref{gcsystem} (5).)
Let $\Omega_p$ be a neighborhood of $p$ in
${\rm Int}\,\mathcal K^+_{\frak p}$.
We may choose $\Omega_p$ sufficiently small
so that the following $(*)$ holds.
\begin{enumerate}
\item[(*)]
If $\mathcal K_{\frak q} \cap \Omega_p \ne \emptyset$,
then $\frak q \le \frak p$
and
$K^+_{\frak q} \cap \Omega_p$ is an open subset
of $U_{\frak q}$.
\end{enumerate}
\begin{sublem}\label{sublem1314}
For each $\frak q \le \frak p$ the set of
$g\in C^{\infty}(\mathcal K^+)_0$
satisfying conditions (a)(b) below is an open dense subset of
$C^{\infty}(\mathcal K^+)_0$ for each $n \ge n_0$.
\begin{enumerate}
\item[(a)]
The restriction of
$
f_0+g
$
to $\mathcal K_{\frak q} \cap \Omega_p$ is a
Morse function on
$\mathcal K_{\frak q} \subset U_{\frak q}$.
\item[(b)]
The composition $(f_0+g)\circ h_i : (0,1) \to \R$ is a Morse
function on $h_i^{-1}(\psi_i^{-1}(\Omega_p))$
for each $
h_i : [0,1] \to V_{i}
$ in Lemma \ref{lem131111}
\end{enumerate}
\end{sublem}
\begin{proof}
This is an immediate consequence of Propositions \ref{prop1427} and \ref{lem1430},
which we will prove in Subsection \ref{subsec:moresedense}.
\end{proof}
Lemma \ref{lem1313} follows from Sublemma \ref{sublem1314}
and Baire's category theorem.
\par
We remark that we have used the fact that we consider only countably many $n$ here.
\end{proof}
Now Proposition \ref{morseufnction} easily follows from Lemma \ref{lem1313}.
\end{proof}
Now we go back to the proof of
Theorem \ref{prop13777}.
We fix $n \ge n_0$.
(We choose $n_0$ so that the compactness  in
Corollary \ref{cor69}
holds.)
We take a strongly smooth function $f$ that satisfies (1)-(4) of
Proposition \ref{morseufnction}.
We take the maps $h_i:[0,1] \to U_{\frak p_i}$ as in Lemma \ref{lem131111}.
We take $0 < a_i < b_i < 1$ such that
\begin{equation}
\frak U(X) \cap
\bigcup_{\frak p\in \frak P}
(\frak s_{\frak p}^{n})^{-1}(0) \cap \mathcal K''_{\frak p}
\subset
\bigcup_{i=1}^N [(\psi_i\circ h_i)([a_i,b_i])].
\end{equation}
and consider the union
$$
S = f({\rm Crit}(f)\cap \mathcal K''_{\frak p})
\cup \bigcup_{i=1}^N (f\circ\psi_i\circ h_i)({\rm Crit}(f \circ\psi_i\circ h_i) \cap [a_i,b_i]).
$$
Here ${\rm Crit}(f \circ\psi_i\circ h_i)$ is the critical point set of the
function $f \circ\psi_i\circ h_i : [0,1] \to [0,\infty)$.
Proposition \ref{morseufnction} implies that $S$ is a finite subset of $[0,\infty)$.
($S$ depends on $n$.)
\begin{lemdef}
Suppose $s \notin S$. We consider
$X^s = X \cap f^{-1}(s)$
and $U^s_{\frak p} =  f^{-1}(s) \cap {\rm Int} \,\mathcal K_{\frak p}$.
\par
Then restricting $\mathcal U_{\frak p}$ to $U^s_{\frak p}$
and restricting the coordinate changes thereto, we obtain
a good coordinate system on $X^s$ of dimension $0$.
We write it as ${\widetriangle{\mathcal U}}\vert_{\{U^s_{\frak p}\}}$.
\par
The restriction of $\frak s^{\epsilon}_{\frak p}$ defines a
multivalued perturbation on  ${\widetriangle{\mathcal U}}\vert_{\{U^s_{\frak p}\}}$.
We write it as $\widetriangle{\frak s^{n,s}}
= \{ \frak s^{n,s}_{\frak p}\}$. The multisections $\widetriangle{\frak s^{n,s}}$ are transversal to $0$.
\end{lemdef}
The proof is obvious from definition.
For $s \in [0,\infty) \setminus S$ we have
$$
[(X^s,{\widetriangle{\mathcal U}}\vert_{ \{U^s_{\frak p}\}},
\widetriangle{\frak s^{n,s}})]
\in \Q
$$
by Definition \ref{defn1355}.

\begin{lem}\label{1314}
For each $s_0 \in (0,\infty)$ there exists a positive number $\delta$
such that
$[(X^s,{\widetriangle{\mathcal U}}\vert_{ \{U^s_{\frak p}\}},
\widetriangle{\frak s^{n,s}})]$
is independent of $s \in (s_0-\delta,s_0+\delta) \setminus S$.
\end{lem}
\begin{lem}\label{lem1315}
There exists $\delta > 0$ such that
$S \cap [0,\delta) = \emptyset$ and
\begin{equation}\label{ddd137}
[(X^s,{\widetriangle{\mathcal U}}\vert_{\{U^s_{\frak p}\}},\widetriangle{\frak s^{n,s}})]
=[(\partial X,\partial{\widetriangle{\mathcal U}},\widetriangle{\frak s^{n}_{\partial}})].
\end{equation}
\end{lem}
\begin{proof}[Proof of Lemma \ref{1314}]
Let
$$
\{p_i \mid i=1,\dots,I\}
= \frak U(X) \cap \vert{\widetriangle{\mathcal U}} \vert_{ \{U^{s_0}_{\frak p}\}}\vert
\cap \bigcup_{\frak p}(\frak s_{\frak p}^{n,{s_0}})^{-1}(0).
$$
The right hand side is a finite set by Proposition \ref{morseufnction} (4).
\par
For each $p_i$ we take $\frak p_i$
with $p_i \in {\rm Int} \,\mathcal K_{\frak p_i}$, and a representative
$(\frak s^{n}_{\frak p_i,j})_{j=1,\dots,\ell_i}$ of $\frak s^{n}_{\frak p_i}$
on an orbifold chart $\frak V_{p_i} = (V_{p_i},\Gamma_{p_i},E_{p_i},\psi_{p_i},\hat\psi_{p_i})$ of $(U_{\frak p_i},\mathcal E_{\frak p_i})$
 at $p_i$.
\par
Then for all sufficiently small $\delta>0$ and $s \in (s_0-\delta,s_0+\delta) \setminus S$, we have the following.
When we put $\hat f_i = f \circ \psi_{p_i} : V_{p_i} \to \R$,
\begin{enumerate}
\item
$s$ is a regular value of
$\hat f_i$,
\item
$\hat{f}_i^{-1}(s)$ intersects transversally to
$(\frak s^{n}_{\frak p_i,j})^{-1}(0)$.
\end{enumerate}
Moreover, we can orient $(\frak s^{n}_{\frak p_i,j})^{-1}(0)$
for each $i,j$ so that
$$
[(X^s,{\widetriangle{\mathcal U}}\vert_{\{U^s_{\frak p}\}},\widetriangle{\frak s^{n,s}})]
=\sum_{i=1}^I
\sum_{j=1}^{\ell_i}
\frac{1}{\ell_i\#\Gamma_{\Gamma_{p_i}}} \hat{f}_i^{-1}(s)
\cdot (\frak s^{n}_{\frak p_i,j})^{-1}(0).
$$
Here $\cdot$ in the right hand side is the intersection number, that is,
the order of the intersection counted with sign.
We use compactness of
$\frak U(X) \cap  \bigcup_{\frak p}(\frak s^{n}_{\frak p})^{-1}(0)$
to show that the intersection number
$\hat{f}_i^{-1}(s)
\cdot (\frak s^{\epsilon}_{\frak p_i,j})^{-1}(0)$
is independent of  $s \in (s_0-\delta,s_0+\delta) \setminus S$.
Thus
Lemma \ref{1314} follows.
Note we use Proposition \ref{prop14777} during this argument.
\end{proof}
\begin{proof}[Proof of Lemma \ref{lem1315}]
The existence of $\delta$ with $S \cap [0,\delta) = \emptyset$
is an immediate consequence of (\ref{normallypositive}).
The formula (\ref{ddd137}) can be proved in the same way as
the proof of Lemma \ref{1314}.
\end{proof}
Now we are ready to complete the proof of
Theorem \ref{prop13777}.
Lemma \ref{1314} implies that
$[(X^s,{\widetriangle{\mathcal U}}\vert_{ \{U^s_{\frak p}\} },
\widetriangle{\frak s^{n,s}})]
$ is independent of $s \in (0,\infty) \setminus S$.
(Note we did not assume $s_0 \notin S$ in Lemma \ref{1314}.)
Therefore Lemma \ref{lem1315} implies that it is equal to
$[(\partial X,\partial{\widetriangle{\mathcal U}},
\widetriangle{\frak s^{n}_{\partial}})]$.
On the other hand, by compactness of $X$ we find that
$X^s$ is empty for sufficiently large $s$.
Hence $[(\partial X,\partial{\widetriangle{\mathcal U}},
\widetriangle{\frak s^{n}_{\partial}})] = 0$ as required.
\end{proof}
We have the following corollary of Theorem \ref{prop13777}.
In particular, we can use it to prove that
Gromov-Witten invariants are
independent of the choice of multivalued perturbations.
\begin{cor}
\begin{enumerate}
\item
If ${\widetriangle{\mathcal U}}$ is an oriented and $0$ dimensional
good coordinate system without boundary of $X$,
then the rational number
$[(X,{\widetriangle{\mathcal U}},\widetriangle{\frak s^{n}})]$
is independent of the multivalued perturbation $\{
\widetriangle{\frak s^{n}}\}$
for all sufficiently large $n$.
\item
Let $(X,{\widehat{\mathcal U}})$ be an oriented Kuranishi space
without boundary. Let ${\widetriangle{\mathcal U}}$ be
an oriented good coordinate system and a KG-embedding
${\widehat{\mathcal U}} \to {\widetriangle{\mathcal U}}$.
Then the rational number $[(X,{\widetriangle{\mathcal U}},\widetriangle{\frak s^{n}})]
$
is independent of the choice of ${\widetriangle{\mathcal U}}$ and $\widetriangle{\frak s^{n}}$
but depends only on the Kuranishi structure $(X,{\widehat{\mathcal U}})$ itself.
\end{enumerate}
\end{cor}
\begin{proof}
(1)
Let $\{\widetriangle{\frak s^{(k) n}}\}$, $k=0,1$ be the two choices of multivalued perturbations.
We consider the direct product of the good
coordinate system
${\widetriangle{\mathcal U}} \times [0,1]$
on $X \times [0,1]$.  The pair
thereon so such that its restriction to
${\widetriangle{\mathcal U}} \times \{0\}$
(resp. ${\widetriangle{\mathcal U}} \times \{1\}$)
is $\widetriangle{\frak s^{(0) n}}$ (resp.
$\widetriangle{\frak s^{(1) n}}$)
and such that it is transversal to $0$.
In fact, we can take it so that it is constant in $[0,1]$ direction
in a neighborhood of ${\widetriangle{\mathcal U}} \times \partial [0,1]$.
We apply Proposition \ref{prop13777} to
$(X \times [0,1],{\widetriangle{\mathcal U}} \times [0,1],
\widetriangle{\frak s^{n}})$ to obtain
$
[X,{\widetriangle{\mathcal U}},\widetriangle{\frak s^{(0) n}}]
=
[X,{\widetriangle{\mathcal U}},\widetriangle{\frak s^{(1) n}}].
$
\par
Independence of $n$ then follows by observing that
$\{\widetriangle{\frak s^{n+n_0}}\}$ is  a multivalued perturbation if
$\{\widetriangle{\frak s^{n}}\}$ is.
\par
(2)
Let ${\widetriangle{\mathcal U^{j}}}$, $j=1,2$ be
two good coordinate systems which are compatible with
${\widehat{\mathcal U}}$.
We consider the direct product Kuranishi structure
$(X \times [0,1],{\widehat{\mathcal U}} \times [0,1])$
on $X \times [0,1]$.
In the same way as the proof of Proposition \ref{relextendgood}, we may assume that
${\widetriangle{\mathcal U^{1}}}$  we put on
${\widehat{\mathcal U}} \times \{0\}$ can be extended
to a neighborhood of it in ${\widehat{\mathcal U}} \times [0,1]$.
Similarly
${\widetriangle{\mathcal U^{2}}}$  we put on
${\widetriangle{\mathcal U}} \times \{1\}$ can be extended
to a neighborhood of it in ${\widehat{\mathcal U}} \times [0,1]$.
We then use Proposition \ref{prop7582752} to
show that there exists a good coordinate system
${\widetriangle{\mathcal U}}$ on $X \times [0,1]$
so that its restriction to $X \times \{0\}$
(resp. $X \times \{1\}$) becomes
${\widetriangle{\mathcal U^{1}}}$
(resp. ${\widetriangle{\mathcal U^{2}}}$.)
Now we use Proposition \ref{prop1220rev}   to find a multivalued perturbation
$\widetriangle{\frak s^{n}}$ on ${\widetriangle{\mathcal U}}$.
(2) now follows from Theorem \ref{prop13777}.
\end{proof}
\begin{defn}
In the situation of (2) we call the rational number
$[(X,{\widetriangle{\mathcal U}},\widetriangle{\frak s^n})]
$
the {\it virtual fundamental class }of $(X,{\widehat{\mathcal U}})$
and write $[(X,{\widehat{\mathcal U}})] \in \Q$.
\end{defn}
We can also prove the following analogue of Proposition \ref{cobordisminvsmoothcor}
\begin{cor}\label{cobordisminvsmoothcormulti}
Let $\frak X_i = (X_i,\widehat{\mathcal U^i})$ be $K$ spaces without boundary
of dimension $0$.
Suppose that there exists a $K$ space
$\frak Y = (Y,\widehat{\mathcal U})$  (but without corner) such that
$$
\partial \frak Y = \frak X_1 \cup -\frak X_2.
$$
Here $-\frak X_2$ is the smooth correspondence $\frak X_2$ with
opposite orientation.
Then we have
\begin{equation}\label{chomotopyrelation2}
[(X_1,\widehat{\mathcal U^1})]
=
[(X_2,\widehat{\mathcal U^2})].
\end{equation}
\end{cor}
\begin{proof}
Using Proposition \ref{prop13777} in place of Stokes' formula
the proof of Corollary \ref{cobordisminvsmoothcormulti}  goes in the same way as the proof of
Proposition \ref{cobordisminvsmoothcor}.
\end{proof}
\subsection{Denseness of the set of Morse functions on orbifold}
\label{subsec:moresedense}

In this subsection we review the proof of the denseness of the set of Morse functions on orbifolds.
We consider the case of one orbifold chart.
The case when we have several orbifold charts is the same but we do not need it.
All the results of this subsection should be well-known. We include it here
only for completeness' sake.
\begin{shitu}\label{shitu1426}
Let $V$ be a manifold on which a finite group $\Gamma$ acts effectively.
We denote by $C^{k}_{\Gamma}(V)$ the set of all $\Gamma$ invariant
$C^k$ functions on $V$.
$\blacksquare$
\end{shitu}
We take and fix a $\Gamma$ invariant Riemannian metric on $V$, which we use
in the proof of some of the lemmata below.
\begin{prop}\label{prop1427}
Suppose we are in Situation \ref{shitu1426}.
The set of all  $\Gamma$ invariant smooth Morse functions on $V$ is a countable intersection of  open dense subsets in $C^{\infty}_{\Gamma}(V)$.
\end{prop}
\begin{proof}
Let $K$ be a compact subset of $V$.
It suffices to prove that the set of all the functions in $C^{2}_{\Gamma}(V)$
which are Morse on $K$ is open and dense.
The openness is obvious.
We will prove that it is also dense.
\par
For $p \in X$ we put
$$
\Gamma_p = \{ \gamma \in \Gamma \mid \gamma p = p\},
$$
and define
\begin{equation}
\aligned
\overset{\circ}{X}(n)  &= \{p \in X \mid \#\Gamma_p = n\}, \\
{X}(n)  &= \{p \in X \mid \#\Gamma_p \ge n\}.
\endaligned
\end{equation}
Note $\overset{\circ}{X}(n)/\Gamma$ is a smooth manifold.
\begin{lem}
Let $p \in \overset{\circ}{X}(n)$. Then $p$ is a critical point of $f$
if and only if $p$ is a critical point of $f\vert_{\overset{\circ}{X}(n)}$.
\end{lem}
\begin{proof}
This is a consequence of the fact that the directional derivative $X[f]$ is zero
if $X \in T_{p}X$ is perpendicular to $\overset{\circ}{X}(n)$.
This fact follows from the $\Gamma$ invariance of $f$.
\end{proof}
We define the following sets.
$$
\aligned
A(n) &= \{ f \in C^{\infty}_{\Gamma}(V) \mid
\text{all the critical points of $f$ on ${X}(n) \cap K$ is Morse.}\} \\
B(n) & = A(n+1)
\cap
\{ f \in C^{\infty}_{\Gamma}(V) \mid
\text{the restriction of $f$ to $\overset{\circ}{X}(n) \cap K$ is Morse.}\}
\endaligned
$$\
They are open sets.
\begin{lem}\label{lemlem1429}
If $A(n+1)$ is dense then $B(n)$ is dense.
\end{lem}
\begin{proof}
Let $W$ be a relatively compact open subset of $K \cap \overset{\circ}{X}(n)$.
We define a $C^1$-map
$
F : W \times C^{2}_{\Gamma}(\overline W) \to T^*\overset{\circ}{X}(n)
$
by
\begin{equation}
F(x,f) = D_xf \in T^*\overset{\circ}{X}(n).
\end{equation}
(Since $f$ is $\Gamma$ invariant and $x \in \overset{\circ}{X}(n)$
it follows $D_xf \in T^*\overset{\circ}{X}(n)$.)
It is easy to see that $F$ is transversal to the submanifold $\overset{\circ}{X}(n) \subset T^*\overset{\circ}{X}(n)$.
Here we identify $\overset{\circ}{X}(n)$ with the zero section of $T^*\overset{\circ}{X}(n)$.
We put
$$
\frak W = \{(x,f) \in W \times C^{2}_{\Gamma}(\overline W)
\mid F(x,f) \in \overset{\circ}{X}(n) \subset T^*\overset{\circ}{X}(n)\}.
$$
$\frak W$ is a sub-Banach manifold of the Banach manifold
$W \times C^{2}_{\Gamma}(\overline W)$.
Moreover the restriction of the projection
$$
{\rm pr} : \frak W \to  C^{2}_{\Gamma}(\overline W)
$$
is a Fredholm map.
Therefore by Sard-Smale theorem the regular value of
${\rm pr}$ is dense.
\begin{sublem}\label{lem1428}
If $f$ is a regular value of ${\rm pr}$ then $f\vert_{\overset{\circ}{X}(n)}$ is Morse on $W$.
\end{sublem}
\begin{proof}
Let $x \in W$ be a critical point of $f$. Then
$(x,f) \in \frak W$.
We consider the following commutative diagram where all
the vertical and horizontal lines are exact.
$$
\begin{CD}
&&
&&
0
&&
0
\\
&&
&& @VVV  @VVV
\\
&&
&&
T_x\overset{\circ}{X}(n)
@>>>
T^*_x\overset{\circ}{X}(n)
\\
&&
&& @ VVV @VVV
\\
0 @>>>T_{(x,f)}\frak W
@>>>T_x\overset{\circ}{X}(n) \oplus T_fC^{2}_{\Gamma}(\overline W)
@>{\overline{D_{(x,f)}F}}>>
\frac{T_{(x,o)}T^*\overset{\circ}{X}(n)}{T_x\overset{\circ}{X}(n)} = T^*_x\overset{\circ}{X}(n)
@>>> 0
\\
&& @ VVV @VVV @VVV\\
0 @>>>
T_fC^{2}_{\Gamma}(\overline W) @>>>T_fC^{2}_{\Gamma}(\overline W)
@>>>0
\\
&& && @VVV
\\
&&&& 0
\end{CD}
$$
Here $\overline{D_{(x,f)}F}$ is the composition of
$D_{(x,f)}F : T_x\overset{\circ}{X}(n) \oplus T_fC^{2}_{\Gamma}(\overline W)
\to T_{(x,o)}T^*\overset{\circ}{X}(n)$ and the projection.
Since $f$ is a regular value the first vertical line is surjective.
We can use it to show that the first horizontal line $T_x\overset{\circ}{X}(n) \to T_x^*\overset{\circ}{X}(n)$ is
surjective by a simple diagram chase.
This map is identified with the Hessian at $x$ of $f\vert_{\overset{\circ}{X}(n)}$.
The sublemma follows.
\end{proof}
We observe that if $f \in A(n+1)$ then the set of critical points in $\overset{\circ}{X}(n) \cap K$
is compact. (This is because it does not have accumulation points on ${X}(n+1) \cap K$.)
Therefore Lemma \ref{lemlem1429} follows from Sublemma \ref{lem1428} and
Sard-Smale theorem.
\end{proof}
\begin{lem}\label{lem1431}
If $B(n)$ is dense then $A(n)$ is dense.
\end{lem}
\begin{proof}
Let $f \in B(n)$. We remark that the set of critical points of $f\vert_{\overset{\circ}{X}(n)}$
on $\overset{\circ}{X}(n) \cap K$ is a finite set.
This is because $f\vert_{\overset{\circ}{X}(n)}$ is a Morse function on
$\overset{\circ}{X}(n) \cap K$ and $f\vert_{\overset{\circ}{X}(n)}$
does not have accumulation points on ${X}(n+1) \cap K$.
Let $p_1,\dots,p_m$ be the critical points of
$f\vert_{\overset{\circ}{X}(n)}$ on $\overset{\circ}{X}(n) \cap K$.
Note the Hessian of $f$ at those points are non-degenerate on
$T_{p_i}\overset{\circ}{X}(n)$ but may  be degenerate
in the normal direction to $\overset{\circ}{X}(n)$.
We choose functions $\chi_i$ and $V_i$
with the following properties.
\begin{enumerate}
\item
$V_i$ is a neighborhood of $p_i$.
\item
The support of $\chi_i$ is in $V_i$
\item
$\chi_i \equiv 1$ in a neighborhood of $p_i$.
\item
$\overline V_i$ ($i=1,\dots,m$) are disjoint.
\item $\overline V_i \cap X(n+1) = \emptyset$.
\item
If $\gamma p_i = p_j$, $\gamma \in \Gamma$,  then $\gamma V_i = V_j$ and
$\chi_j \circ \gamma = \chi_i$.
\end{enumerate}
We use our $\Gamma$ invariant Riemannian metric and put
$$
f_n(x) = d(x,\overset{\circ}{X}(n))^2
$$
We choose $V_i$ small so that $\chi_i f_n$ is a smooth function.
(We use Item (5) above here.)
Now
\begin{equation}\label{form1413}
f_{\epsilon} = f + \epsilon \sum_{i=1}^m \chi_i f_n
\end{equation}
is a Morse function for sufficiently small positive $\epsilon$.
Item (6) implies that this function is $\Gamma$ invariant.
Therefore $f_{\epsilon} \in A(n)$. Moreover $f_{\epsilon}$ converges to $f$ as
$\epsilon \to 0$.
\end{proof}
Lemmata \ref{lemlem1429} and \ref{lem1431} imply that
$A(1)$ is dense.
The proof of Proposition \ref{prop1427} is complete.
\end{proof}
\begin{rem}
The set of functions whose gradient flow is Morse-Smale is {\it not} dense
in general
in the case of orbifold.
This is an important point where Morse theory of orbifold is different from one of manifold.
In fact we can use virtual fundamental chain
technique to work out the theory of Morse homology for orbifold.
\end{rem}
\begin{prop}\label{lem1430}
Suppose we are in Situation \ref{shitu1426}.
Let $h : [0,1] \to V$ be a smooth embedding.
Then the set of all the  functions $f$ in $C^{\infty}_{\Gamma}(V)$
such that $f \circ h$ is a Morse function on $(0,1)$
is a countable intersection of  open dense subsets.
\end{prop}
\begin{proof}
The proof is similar to Proposition \ref{prop1427}.
Let $0<c<1/2$.
We put
$$
\aligned
\overset{\circ}T(n,c) &= \{t \in [c,1-c] \mid h(t) \in \overset{\circ}{X}(n)\}, \\
T(n,c) &= \{t \in [c,1-c] \mid h(t) \in {X}(n)\}, \\
S(n,c) &= \{t_0 \in \overset{\circ}T(n,c) \mid (dh/dt)(t_0) \in T_{h(t_0)}\overset{\circ}{X}(n) \}.
\endaligned
$$
We also put
$$
\aligned
C(n,c) &=
 \{ f \in C^{\infty}_{\Gamma}(V) \mid
\text{all the critical points of $f \circ h$ on $T(n,c)$ is Morse.}\} \\
D(n,c) &=
C(n+1,c) \cap
 \{ f \in C^{\infty}_{\Gamma}(V) \mid
\text{all the critical points of $f \circ h$ on $S(n,c)$ is Morse.}\}
\endaligned
$$
\begin{lem}\label{lem1434}
If $C(n+1,c)$ is dense then $D(n,c)$ is dense.
\end{lem}
\begin{proof}
We fix $c$ and will prove $D(n,c)$ is dense.
We take $\epsilon >0$ such that  $c - \epsilon > 0$.
Let $U_1$ be a sufficiently small neighborhood of $X(n+1)$
which we choose later.
We take $U_2$ a neighborhood $X(n) \setminus U_1$
and $\pi : U_2 \to \overset{\circ}X(n)$ be the projection of normal
bundle. We may assume $U_1$, $U_2$ and $\pi$ are $\Gamma$ equivariant.
We take a compact subset $Z$ of $U_2 \cap  \overset{\circ}X(n)$.
We take a set $S \subset [c-\epsilon,1-c+\epsilon]$ with the following properties.
\begin{enumerate}
\item
${\rm Int}S \supset S(n,c) \cap h^{-1}(Z)$.
\item
$S$ is a finite disjoint union of closed intervals.
\item
The composition $\pi \circ h$ is an embedding on $S(n,c)$.
\end{enumerate}
Note the differential of $\pi \circ h$ is injective on $S(n,c) \cap h^{-1}(Z)$.
Therefore we can take such $S$.
We define a map
$
G : {\rm Int} S \times C^{2}_{\Gamma}(K) \to \R
$
by
\begin{equation}
G(t_0,f) = \frac{d( f\circ \pi \circ h)}{dt}(t_0).
\end{equation}
Using the fact that $\pi \circ h$ is an embedding
and $\overset{\circ}X(n)/\Gamma$ is a smooth
manifold we can easily show that $G$ is
transversal to $0$. We put
$\frak N = G^{-1}(0) \subset {\rm Int} S \times C^{2}_{\Gamma}(K)$.
It is a Banach submanifold.
The restriction of the projection
$\frak N \to C^{2}_{\Gamma}(K)$ is a Fredholm map.
In the same way as Subemma \ref{lem1428} we can show that if
$f \in C^{2}_{\Gamma}(K)$ is a regular value of $\frak N \to C^{2}_{\Gamma}(K)$
then $f\circ \pi \circ h$ is a Morse function on ${\rm Int}\,S$.
We remark that at the point of $S(n,c)$
the Hessian of $f\circ \pi \circ h$ is the same as
the Hessian of $f\circ h$.
Thus by Sard-Smale theorem we find that the
set of $f \in C^{2}_{\Gamma}(K)$
such that all the critical points on $f\circ h$ on $S(n,c)  \cap h^{-1}(Z)$
is Morse, is dense.
\begin{sublem}\label{sublem1435}
Let $f \in C(n+1,c)$. Then there exists a
compact set $P \subset \overset{\circ}T(n,c)$ such that
there is no critical point of $f\circ h$ on $\overset{\circ}T(n,c) \setminus P$.
\end{sublem}
\begin{proof}
If not there is a sequence $t_i \in \overset{\circ}T(n,c)$ such that
$t_i$ is a critical point of $f\circ h$ and no
subsequence of $t_i$  converges to an element of $\overset{\circ}T(n,c)$.
We may assume that $t_i$ converges. Then the limit $t$
should be an element of $T(n+1,c)$.
Since $f \in C(n+1,c)$ the function
$f\circ h$ is Morse at $t$. Therefore there is no critical point of $f\circ h$
other than $p$ in a neighborhood of $t$.
This is a contradiction.
\end{proof}
Now let $f \in C(n+1,c)$.
We choose $U_1$, $U_2$ and $Z$ such that $Z
\supset h(P)$ and $P$ is as in Sublemma \ref{sublem1435}.
Then there exists a sequence of functions $f_i$ such that $f_i$ converges to $f$
and all the  critical points of $f_i \circ h$  on  $S(n,c)  \cap h^{-1}(Z)$ are Morse.
Therefore $f_i \in D(n,c)$ for all sufficiently large $i$.
The proof of Lemma \ref{lem1434}
is complete.
\end{proof}
\begin{lem}\label{Dndense}
If $D(n,c)$ is dense then $C(n,c)$ is dense.
\end{lem}
\begin{proof}
Let $f \in D(n,c)$.
\begin{sublem}
The set
$$
Q = \{ t \in \overset{\circ}T(n,c) \mid  \text{$t$ is a critical point of $f\circ h$} \}
$$
is a finite set.
\end{sublem}
\begin{proof}
Suppose $Q$ is an infinite set and $t_i \in Q$ is
its infinitely many points. We may also assume that $t_i$
converges to $t$. If $t \in  \overset{\circ}T(n,c) \setminus S(n,c)$ then
$h$ is transversal to $\overset{\circ}{X}(n)$ at $t$.
Therefore $h(t_i) \notin \overset{\circ}{X}(n)$ for large $i$.
This contradicts to $t_i \in Q$.
Therefore $t \in S(n,c) \cup T(n+1,c)$.
Since $f \in D(n,c)$ the composition $f \circ h$ is Morse at $t$.
Therefore there is no critical point of $f \circ h$ other than $t$
in a neighborhood of $t$. This contradicts $\lim t_i = t$.
\end{proof}
Let $Q \setminus S(n,c) = \{t_1,\dots,t_n\}$ and we put
$p_i = h(t_i) \in \overset{\circ}T(n,c)$.
We define $V_i$ and $\chi_i$ in the same way as the proof of
Lemma \ref{lem1431}
and define $f_{\epsilon}$ in the same way as (\ref{form1413}).
It is easy to see that $f_{\epsilon} \in C(n,c)$
and $f_{\epsilon}$ converges to $f$ as $\epsilon \to 0$.
The proof of Lemma \ref{Dndense} is complete.
\end{proof}
By Lemmata \ref{lem1434} and \ref{Dndense} we have proved
that $C(1,c)$ is dense for all $c \in (0,1/2)$. The proof of Proposition \ref{lem1430} is complete.
\end{proof}
\begin{rem}
We remark that we do not assume $h$ to be $\Gamma$ invariant in Proposition \ref{lem1430}.
\end{rem}


\section{Appendices : Orbifold and orbibundle by local coordinate}
\label{sec:ofd}

In this section we describe the story of orbifold as far as we
need in this document.
We restrict ourselves to effective orbifolds
and  regard only embeddings as morphisms.
The category $\mathscr{OB}_{\rm ef,em}$ where objects are effective
orbifolds and morphisms are embeddings among them
is naturally a $1$ category. Moreover
it has the following property.
We consider the forgetful map
$$
\frak{forget} : \mathscr{OB}_{\rm ef,em} \to \mathscr{TOP}
$$
where $\mathscr{TOP}$ is the category of topological spaces.
Then
$$
\frak{forget} : \mathscr{OB}_{\rm ef,em}(c,c')
\to \mathscr{TOP}(\frak{forget}(c),\frak{forget}(c'))
$$
is injective.
In other words, we can check the equality between morphisms
set-theoretically.
This is a nice property, which we use extensively in the main
body of this document.
If we go beyond this category then we need to distinguish
carefully the two notions, two morphisms are equal,
two morphisms are isomorphic.
It will then makes the argument much more complicated.
\footnote{We need to use several maps between underlying topological
spaces of orbifolds, such as projection of bundles or covering maps.
In case we include those maps, we need to see carefully whether
set-theoretical equality is enough to show various properties of them
are preserved.}
\par
We emphasize that there is nothing new in this section.
The story of orbifold is  classical and is well-established.
It has been used in various branches of mathematics
since its invention by Satake \cite{satake} in more than 50 years ago.
Especially, if we restrict ourselves to
effective orbifolds, the story of orbifolds is nothing more than
straightforward generalization of the theory of
smooth manifolds. The only important issue is the
observation that for effective orbifolds almost everything work
in the same way as manifolds.

\subsection{Orbifolds and embedding among them}
\label{subsec;ofds}
\begin{defn}\label{2661}
Let $X$ be a paracompact Hausdorff space.
\begin{enumerate}
\item
An {\it orbifold chart of $X$ (as a topological space)}
\index{orbifold ! orbifold chart}
is a triple $(V,\Gamma,\phi)$ such that
$V$ is a manifold, $\Gamma$ is a finite group acting smoothly
and effectively on $V$ and $\phi : V \to X$
is a $\Gamma$ equivariant continuous map\footnote{The $\Gamma$ action
on $X$ is trivial.} which induces a
homeomorphism $\overline\phi : V/\Gamma \to X$ onto an open subset of $X$.
We assume that there exists $o \in V$ such that $\gamma o = o$
for all $\gamma \in \Gamma$.
We call $o$ the  {\it base point}.
\index{orbifold ! base point of an orbifold chart}
We say $(V,\Gamma,\phi)$ is an orbifold chart {\it at $x$} if
$x = \phi(o)$.
We call $\Gamma$ the {\it isotropy group},
\index{orbifold ! isotropy group of an orbifold chart}
$\phi$ the {\it local uniformization map}
\index{orbifold ! local uniformization map}  and
$\overline\phi$ the
{\it parametrization}
\index{orbifold ! parametrization}.
\item
Let $(V,\Gamma,\phi)$ be an orbifold chart and $p \in V$.
We put $\Gamma_p = \{ \gamma \in \Gamma \mid \gamma p = p\}$.
Let $V_p$ be a $\Gamma_p$ invariant open neighborhood of $p$
in $V$. We assume the map
$\overline\phi : V_p/\Gamma_p \to X$ is injective.
(In other words, we assume that
$\gamma V_{p} \cap V_p \ne \emptyset$ implies $\gamma \in \Gamma_p$.)
We say such triple $(V_{p},\Gamma_p,\phi\vert_{V_p})$
a {\it subchart} \index{orbifold ! subchart} of $(V,\Gamma,\phi)$.
\item
Let $(V_i,\Gamma_i,\phi_i)$  $(i=1,2)$ be orbifold charts of $X$.
We say that they are {\it compatible}
\index{orbifold ! compatible (orbifold chart)} if the following holds
for each $p_1 \in V_1$ and $p_2 \in V_2$ with
$\phi_1(p_1) = \phi_2(p_2)$.
\begin{enumerate}
\item
There exists a group isomorphism $h : (\Gamma_1)_{p_1}
\to (\Gamma_2)_{p_2}$.
\item
There exists an $h$ equivariant diffeomorphism $\tilde\varphi : V_{1,p_1}
\to V_{2,p_2}$. Here $V_{i,p_i}$ is a $(\Gamma_i)_{p_i}$
equivariant subset of $V_i$ such that
 $(V_{i,p_i},(\Gamma_i)_{p_i},\phi\vert_{V_{i,p_i}})$
 is a subchart.
 \item
 $\phi_2 \circ \tilde\varphi = \phi_1$ on $V_{1,p_1}$.
\end{enumerate}
\item
A {\it representative of an orbifold structure}
\index{orbifold ! representative of an orbifold structure} on $X$
is a set of orbifold charts $\{(V_i,\Gamma_i,\phi_i) \mid i \in I\}$
such that each two of the charts are compatible in the sense of (3)
above and
$
\bigcup_{i\in I} \phi_i(V_i) = X,
$
is a locally finite open cover of $X$.
\end{enumerate}
\end{defn}
\begin{defn}\label{def262220}
Suppose that $X$, $Y$ have representatives
of orbifold structures $\{(V^X_i,\Gamma^X_i,\phi^X_i) \mid i \in I\}$
and $\{(V^Y_j,\Gamma^Y_j,\phi^Y_j) \mid j \in J\}$,
respectively.
A continuous map $f : X \to Y$ is said to be an
{\it embedding}
\index{orbifold ! embedding of orbifolds}\index{embedding ! of orbifolds}
if the following holds.
\begin{enumerate}
\item
$f$ is an embedding of  topological spaces.
\item
Let $p \in V^X_i$, $q \in V^Y_j$ with
$f(\phi_i(p)) = \phi_j(q)$. Then we have the following.
\begin{enumerate}
\item
There exists an isomorphism of groups
$h_{p;ji} : (\Gamma_i^X)_p \to (\Gamma_j^Y)_q$.
\item
There exist  $V^X_{i,p}$ and $V^Y_{j,q}$ such that
$(V^X_{i,p},(\Gamma^X_i)_p,\phi_i\vert_{V^X_{i,p}})$
is a subchart for $i=1,2$.
There exists an $h_{p;ji}$ equivariant embedding of manifolds
$\tilde f_{p;ji}: V^X_{i,p} \to V^Y_{j,q}$.
\item
The diagram below commutes.
\begin{equation}\label{diag2633}
\begin{CD}
V^X_{i,p} @ >{\tilde f_{p;ji}}>>
V^Y_{j,q}  \\
@ V{\phi_{i}}VV @ VV{\phi_{j}}V\\
X @ > {f} >> Y
\end{CD}
\end{equation}
\end{enumerate}
\end{enumerate}
Two orbifold embeddings are said to be {\it equal} if they coincide set-theoretically.
\end{defn}
\begin{rem}
Note that an embedding of effective orbifolds
is a continuous map $f : X \to Y$ of
underlying topological spaces, which has
the properties (2) above.
\par
When we study morphisms among not-necessary effective
orbifolds or morphisms between effective orbifolds which is not
necessarily an embedding, then such a morphism is a
continuous map $f : X \to Y$ of underlying topological spaces
{\it plus} certain additional data.
For example, we consider noneffective orbifold
that is a point with an action of a nontrivial finite group $\Gamma$.
Then the morphism from this noneffective orbifold to itself
contains a datum which is an automorphism of the group $\Gamma$.
(Two such morphisms $h_1,h_2$ are equivalent if there exists an
inner automorphism $h$ such that $h_1 = h \circ h_2$.)
\end{rem}

\begin{lem}\label{lem26444}
\begin{enumerate}
\item
The composition of embeddings is an embedding.
\item
The identity map is an embedding.
\item
If an embedding is a homeomorphism,
then its inverse is also an embedding.
\end{enumerate}
\end{lem}
The proof is easy and is left to the reader.

\begin{defn}\label{defn285555}
\begin{enumerate}
\item We call an embedding of orbifolds a {\it diffeomorphism}
\index{orbifold ! diffeomorphism} if it is a
homeomorphism in addition.
\item
We say that two representatives of orbifold structures on $X$ are
{\it equivalent} if the identity map regarded as a map between
$X$ equipped with those two representatives of orbifold structures
is a diffeomorphism. This is an equivalence relation by Lemma \ref{lem26444}.
\item
An equivalence class of the equivalence relation (2) is called
an {\it orbifold structure} of $X$.
\index{orbifold ! orbifold structure}
An {\it orbifold} is a pair of topological space and its
orbifold structure. \index{orbifold ! orbifold}
\item
The condition for $X \to Y$ to be an embedding
does not change if we replace representatives of orbifold structures
to equivalent ones. So we can define the notion of
an {\it embedding of orbifolds}.
\item
If $U$ is an open subset of an orbifold $X$, then
there exists a unique orbifold structure on $U$ such that
the inclusion $U \to X$ is an embedding.
We call $U$ with this orbifold structure an {\it open suborbifold}.
\index{orbifold ! open suborbifold}
\end{enumerate}
\end{defn}
\begin{defn}\label{defn26550}
\begin{enumerate}
\item
Let $X$ be an orbifold. An orbifold chart $(V,\Gamma,\phi)$
of underlying topological space $X$
in the sense of Definition \ref{2661} (1)
is called an {\it orbifold chart of an orbifold}  $X$ if
the map $\overline\phi : V/\Gamma \to X$ induced by $\phi$ is an
embedding of orbifolds.
\index{orbifold ! orbifold chart}
\item
Hereafter when $X$ is an orbifold,
an `orbifold chart' always means an orbifold chart of an orbifold in the sense of
(1).
\item
In case when an orbifold structure on $X$ is given,
a representative of its orbifold structure is
called an {\it orbifold atlas}.
\index{orbifold ! orbifold atlas}
\item
Two orbifold charts  $(V_i,\Gamma_i,\phi_i)$
are said to be isomorphic if there exist a
group isomorphism $h : \Gamma_1 \to \Gamma_2$
and an $h$-equivariant diffeomorphism $\tilde\varphi : V_1 \to V_2$
such that $\phi_2 \circ \tilde\varphi = \phi_1$.
The pair $(h,\tilde\varphi)$ is said to be the {\it isomorphism} between
two orbifold charts.
\end{enumerate}
\end{defn}
\begin{prop}\label{prop266}
In the situation of Definition \ref{defn26550} (4), suppose
$(h,\tilde\varphi)$ and $(h',\tilde\varphi')$ are both
isomorphisms between
two orbifold charts $(V_1,\Gamma_1,\phi_1)$ and $(V_2,\Gamma_2,\phi_2)$.
Then there exists $\mu \in \Gamma_2$
such that
\begin{equation}\label{262}
h'(\gamma) = \mu h(\gamma) \mu^{-1},
\qquad
\tilde\varphi'(x) = \mu \tilde\varphi(x).
\end{equation}
\par
On the contrary, if $(h,\tilde\varphi)$ is an isomorphism between
orbifold charts then $(h',\tilde\varphi')$ defined by (\ref{262})
is also an isomorphism between
orbifold charts.
In particular, any automorphism of orbifolds charts
$(h,\tilde\varphi)$ is given
by $h(\gamma) = \mu\gamma \mu^{-1}$, $\tilde\varphi(x) = \mu x$ where $\mu$ is an
element of $\Gamma$.
\end{prop}
\begin{proof}
The proposition follows immediately from the next lemma.
\begin{lem}\label{lem21tenhatena}
Let $V_1$, $V_2$ be manifolds on which finite groups $\Gamma_1$,
$\Gamma_2$ act effectively and smoothly.
We assume $V_1$ is connected.
Let $(h_i,\tilde\varphi_i)$ $(i=1,2)$ be pairs such that
$h_i : \Gamma_1 \to \Gamma_2$ are injective group homomorphisms
and $\tilde\varphi_i : V_1 \to V_2$ are $h_i$-equivariant embeddings of manifolds.
Moreover, we assume that the induced maps ${\varphi}_i :V_1/\Gamma_1
\to V_2/\Gamma_2$
are embeddings of orbifolds.
Furthermore we assume that the induced map ${\varphi}_1 : V_1/\Gamma_1
\to V_2/\Gamma_2$ coincides with ${\varphi}_2$ set-theoretically.
\par
Then there exists $\mu \in \Gamma_2$ such that
$$
\tilde\varphi_2(x) = \mu\tilde\varphi_1(x),
\qquad
h_2(\gamma) = \mu h_1(\gamma)\mu^{-1}.
$$
\end{lem}
\begin{proof}
For the sake of simplicity we prove only
the case when Condition \ref{convinv} below is satisfied.
Let $X$ be an orbifold.
For a point $x \in X$ we take its orbifold chart $(V_x,\Gamma_x,\psi_x)$.
We say $x \in {\rm Reg}(X)$ if
$\Gamma_x = \{1\}$,
and put ${\rm Sing}(X) = X \setminus {\rm Reg}(X)$.
\begin{conds}\label{convinv}
We assume that $\dim {\rm Sing}(X) \le \dim X -2$.
\end{conds}
This condition is satisfied if $X$ is oriented.  (In fact, Condition \ref{convinv} fails only
when there exists an element of $\Gamma_x$ (an isotropy group of some orbifold chart)
whose action is $(x_1,x_2,\dots,x_n) \mapsto (-x_1,x_2,\dots,x_n)$ for some coordinate $(x_1,\dots,x_n)$.
Therefore we can always assume Condition \ref{convinv}
in the study of Kuranishi structure, by adding a trivial factor which is acted by
the induced representation of $t \mapsto -t$ to both the obstruction bundle and to the Kuranishi neighborhood.)
\par
Let $x_0 \in V_1^0$.
By assumption there exists uniquely $\mu \in \Gamma_2$
such that $\tilde\varphi_2(x_0) = \mu\tilde\varphi_1(x_0)$.
By Condition \ref{convinv} the subset $V^0_1$ is connected.
Therefore the above element $\mu$ is independent of $x_0 \in V_1^0$
by uniqueness.
Since $V_1^0$ is dense, we conclude $\tilde\varphi_2(x) = \mu\tilde\varphi_1(x)$ for any $x \in V_1$.
Now, for $\gamma \in \Gamma_1$, we calculate
$$
h_1(\gamma)\tilde\varphi_1(x_0)
=
\tilde\varphi_1(\gamma x_0)
=
\mu^{-1}\tilde\varphi_2( \gamma x_0)
=
\mu^{-1}h_2(\gamma)\tilde\varphi_2(x_0)
=
\mu^{-1}h_1(\gamma)\mu \tilde\varphi_1(x_0).
$$
Since the induced map is an embedding of orbifold,
it follows that the isotropy group of $\tilde\varphi_1(x_0)$
is trivial. Therefore
$h_1(\gamma) = \mu^{-1}h_2(\gamma)\mu$
as required.
\end{proof}
The proof of Proposition \ref{prop266} is complete.
\end{proof}
\begin{defn}\label{defn281010}
Let $X$ be an orbifold.
\begin{enumerate}
\item
A function $f : X \to \R$ is said to be a {\it smooth function} if
for any orbifold chart $(V,\Gamma,\phi)$ the composition
$f\circ \phi : V \to \R$ is smooth.
\index{orbifold ! smooth function on orbifold}
\item
A {\it differential form}
\index{orbifold ! differential form on orbifold}
\index{differential form ! on oribifold}
on an orbifold $X$ assigns a
$\Gamma$ invariant differential form $h_{\frak V}$ on $V$ to
each orbifold chart $\frak V = (V,\Gamma,\phi)$
such that the following holds.
\begin{enumerate}
\item
If $(V_1,\Gamma_1,\phi_1)$ is isomorphic to
$(V_2,\Gamma_2,\phi_2)$ and $(h,\tilde\varphi)$ is an isomorphism,
then $\tilde\varphi^*h_{\frak V_2} = h_{\frak V_1}$.
\item
If $\frak V_p = (V_p,\Gamma_p,\phi_p)$ is a subchart of $\frak V =(V,\Gamma,\phi)$,
then $h_{\frak V}\vert_{V_p} = h_{\frak V_p}$.
\end{enumerate}
\item
An $n$ dimensional orbifold $X$ is said to be {\it orientable} if there exists
a differential $n$-form $\omega$ such that $\omega_{\frak V}$
never vanishes.
\index{orbifold ! orientation}
\item
Let $\omega$ be an $n$-form as in (3) and
$\frak V = (V,\Gamma,\phi)$  an orbifold chart. Then
we give $V$ an orientation so that it is compatible with $\omega_{\frak V}$.
The $\Gamma$ action preserves the orientation.
We call such $(V,\Gamma,\phi)$ equipped with an orientation of $V$,
an {\it oriented orbifold chart}.
\index{orbifold ! oriented orbifold chart}
\item
Let $\bigcup_{i\in I}U_i = X$ be an open covering of an orbifold $X$.
A {\it smooth partition of unity subordinate to the covering} $\{U_i\}$
\index{orbifold ! partition of unity} is
a set of functions $\{\chi_i\mid i\in I\}$ such that:
\begin{enumerate}
\item
$\chi_i$ are smooth functions.
\item
The support of $\chi_i$ is contained in $U_i$.
\item
$\sum_{i\in I}\chi_i = 1$.
\end{enumerate}
\end{enumerate}
\end{defn}
\begin{lem}
For any locally finite open covering of an orbifold $X$ there exists a smooth
partition of unity subordinate to it.
\end{lem}
We omit the proof, which is an obvious analogue of the standard proof
for the case of manifolds.

\begin{defn}\label{orbifolddefn}
An {\it orbifold with corner} is defined in the same way.
\index{orbifold ! orbifold with corner}
We require the following.
\begin{enumerate}
\item In Definition \ref{2661} (1) we assume that $V$ is a manifold with
corners.
\item
Let $S_k(V)$ be the set of points which lie on the codimension $k$ corner
and $\overset{\circ}S_k(V) = S_k(V) \setminus \bigcup_{k' > k}S_{k'}(V)$.
We require that $\Gamma$ action on each connected component of
$\overset{\circ}S_k(V)$ is effective.
(Compare Condition \ref{effectivitycorner}.)
\item
For an embedding of orbifolds with corners
\index{orbifold ! embedding of orbifolds with corners}
\index{embedding ! of orbifolds with corners}
we require that the map
$\tilde f$ in
Definition \ref{def262220} (c) satisfies
$\tilde f(\overset{\circ}S_k(V_1)) \subset \overset{\circ}S_k(V_2)$.
\end{enumerate}
\end{defn}

\begin{lem}\label{lem26999}
Let $X_i$ $(i=1,2)$ be orbifolds and $\varphi_{21} : X_1 \to X_2$  an embedding.
Then we can find an orbifold atlas
$\{\frak V^i_{\frak r} =
\{ (V_{\frak r}^i, \Gamma_{\frak r}^i, \phi_{\frak r}^i )\} \mid \frak r \in \frak R_i\}$ with the following properties.
\begin{enumerate}
\item $\frak R_1 \subseteq \frak R_2$.
\item  $V^2_{\frak r} \cap \varphi_{21}(X_1) \ne \emptyset$
if and only if $\frak r \in \frak R_1$.
\item
If $\frak r \in \frak R_1$ then
$\varphi_{21}^{-1}(\phi_{\frak r}^2(V^2_{\frak r})) =
\phi_{\frak r}^1(V^1_{\frak r})$
and there exists
$(h_{\frak r,21},\tilde\varphi_{\frak r,21})$
such that:
\begin{enumerate}
\item
$h_{\frak r,21} : \Gamma^1_{\frak r} \to \Gamma^2_{\frak r}$
is a group isomorphism.
\item
$\tilde\varphi_{\frak r,21} : V^1_{\frak r} \to V^2_{\frak r}$
is an $h_{\frak r,21}$-equivariant embedding of smooth manifolds.
\item
The next diagram commutes.
\begin{equation}\label{diagin2611}
\begin{CD}
V^1_{\frak r} @ >{\tilde\varphi_{\frak r,21}}>>
V^2_{\frak r}  \\
@ V{\phi^{\frak r}_{1}}VV @ VV{\phi^{\frak r}_{2}}V\\
X_1 @ > {\varphi_{21}} >> X_2
\end{CD}
\end{equation}
\end{enumerate}
\item
In case $X_i$ has a boundary or corners we may choose
our charts so that the following is satisfied.
\begin{enumerate}
\item
$V^i_{\frak r}$ is an open subset of $\overline V^i_{\frak r}
\times [0,1)^{d({\frak r})}$, where $d(\frak r)$ is independent of $i$
and $\overline V^i_{\frak r}$ is a manifold without boundary.
\item
There exists a point $o^i(\frak r)$ which is fixed by all
$\gamma \in \Gamma^i_{\frak r}$ such that
$[0,1)^{d({\frak r})}$ components of $o^i(\frak r)$ are all $0$.
\item We put
$$
\varphi_{\frak r,21}(\overline y',(t'_1,\dots,t'_{d({\frak r})}))
=
(\overline y,(t_1,\dots,t_{d({\frak r})})).
$$
Then $t_i = 0$ if and only if $t'_i = 0$.
\end{enumerate}
\end{enumerate}
We may take our charts finer than given coverings of $X_1$ and $X_2$.
\end{lem}
\begin{proof}
For each $x \in X_1$ we can find  orbifold charts $\frak V^i_{x}$
for $i=1,2$, such that
$\varphi_{21}^{-1}(U^2_{x}) = U^1_{x}$,
$x \in U^1_{x}$ and that
 there exists a representative
$(h_{x,21},\tilde\varphi_{x,21})$
of embedding $U^1_{x} \to U^2_{x}$ that is a restriction
of $\varphi_{21}$.
In case $X_i$ has a boundary or corners,
we choose them so that (4) is also satisfied.
\par
We cover $X_1$ by finitely many of such $U_{x_j}^{1}$.
This is our
$\{\frak V^1_{\frak r} \mid \frak r \in \frak R_1\}$.
Then we have $\{\frak V^2_{\frak r} \mid \frak r \in \frak R_1\}$,
satisfying (3)(4) and that covers
$\varphi_{21}(X_1)$. We can extend it to
$\{\frak V^2_{\frak r} \mid \frak r \in \frak R_2\}$
so that (1)(2) are also satisfied.
\end{proof}
\begin{defn}
We call $(h_{\frak r,21},\tilde\varphi_{\frak r,21})$
a {\it local representative of embedding}  $\varphi_{\frak r,21}$ on the charts
$\frak V^1_{\frak r}$, $\frak V^2_{\frak r}$.
\index{orbifold ! local representative of embedding}
\index{embedding ! local representative of embedding of orbifolds}
\end{defn}
\begin{lem}
If $(h_{\frak r,21},\tilde\varphi_{\frak r,21})$, $(h'_{\frak r,21},\tilde\varphi'_{\frak r,21})$
are local representatives of an embedding of the same
charts
$\frak V^1_{\frak r}$, $\frak V^2_{\frak r}$,
then there exists $\mu \in \Gamma_2$
such that
$$
\tilde\varphi'_{\frak r,21}(x) = \mu\tilde\varphi_{\frak r,21}(x),
\qquad
h_{\frak r,21}^{\prime}(\gamma) = \mu h_{\frak r,21}(\gamma)\mu^{-1}.
$$
\end{lem}
This is a consequence of Lemma \ref{lem21tenhatena}.
\begin{lem}\label{smoothstruemb}
Let $X$ be a topological space, $Y$  an orbifold,
and $f : X \to Y$  an embedding of topological spaces.
Then the orbifold structure of $X$ by which $f$ becomes an
embedding of orbifolds is unique if exists.
\end{lem}
\begin{proof}
Let $X_1$, $X_2$ be orbifolds whose underlying topological spaces
are both $X$ and such that $f_i : X_i \to Y$  are embeddings of orbifolds for
$i=1,2$. We will prove that the identity map
$\rm{id} : X_1 \to X_2$ is a diffeomorphism
of orbifolds.
Since the condition for a homeomorphism
to be a diffeomorphism of orbifolds is a local condition,
it suffices to check it on a neighborhood of each point.
Let $p \in X$ and $q = f(p)$.
We take a representative $(h_i,\tilde\varphi_i)$ of the orbifold embeddings
$f_i : X_i \to Y$ using the orbifold charts
$\frak V^i_p = (V^i_p,\Gamma^i_p,\phi^i_p)$ of $X$ and
$\frak V_q = (V_q,\Gamma_q,\phi_q)$ of $Y$.
The maps $h_i : \Gamma^i_p \to \Gamma^i_q$ are group isomorphisms.
So we have a group isomorphism
$
h = h_2^{-1}\circ h_1 : \Gamma^1_p \to \Gamma^2_p.
$
Since $\tilde\varphi_1(V^1)/\Gamma_p = \tilde\varphi_2(V^2)/\Gamma_p$
set-theoretically, we have
$\tilde\varphi_1(V^1_p) = \tilde\varphi_2(V^2_p) \subset V_q$.
They are smooth submanifolds since $f_i$ are  embeddings of
orbifolds. Therefore
$
\varphi = \tilde\varphi_2^{-1}\circ \tilde\varphi_1
$
is defined in a neighborhood of the base point $o^i_p$
and is a diffeomorphism.
$(h,\tilde\varphi)$ is a local representative of ${\rm id}$.
\end{proof}

\subsection{Vector bundle on orbifold}
\label{subsec:vectorbundle}

\begin{defn}\label{defn2613}
Let $(X,\mathcal E,\pi )$ be a pair of an orbifold $X$ and $\pi : \mathcal E \to X$  a continuous map
between their underlying topological spaces.
Hereafter we write $(X,\mathcal E)$ in place of $(X,\mathcal E,\pi)$.
\begin{enumerate}
\item
An {\it orbifold chart} of $(X,\mathcal E)$ is $(V,E,\Gamma,\phi,\widehat\phi)$
with the following properties.
\index{orbifold ! orbifold chart of a vector bundle}
\begin{enumerate}
\item
$\frak V = (V,\Gamma,\phi)$ is an orbifold chart of the orbifold $X$.
\item $E$ is a finite dimensional vector space equipped with
a linear $\Gamma$
action.
\item
$(V \times E,\Gamma,\widehat\phi)$ is an orbifold chart of the
orbifold $\mathcal E$.
\item
The diagram below commutes set-theoretically,
\begin{equation}\label{diag26399}
\begin{CD}
V \times E @ >{\widehat\phi}>>
\mathcal E  \\
@ V{}VV @ VV{{\pi}}V\\
V @ > {\phi} >> X
\end{CD}
\end{equation}
where the left vertical arrow is the projection to the
first factor.
\end{enumerate}
\item
In the situation of (1), let  $( V_{p},\Gamma_p,\phi\vert_{V_p})$
be a subchart of $(V,\Gamma,\phi)$ in the sense of
Definition \ref{2661} (2).
Then  $(V_{p},E,\Gamma_p,\phi\vert_{V_p},\widehat\phi\vert_{V_p \times E})$
is an orbifold chart of $(X,\mathcal E)$.
We call it a {\it subchart}.
\index{orbifold ! subchart of an orbifold chart of vector bundle}
\item
Let $(V^i,E^i,\Gamma^i,\phi^i,\widehat{\phi^i})$  $(i=1,2)$ be orbifold charts of $(X,\mathcal E)$.
We say that they are {\it compatible} if the following holds
for each $p_1 \in V^1$ and $p_2 \in V^2$ with
$\phi^1(p_1) = \phi^2(p_2)$:
There exist open neighborhoods $V^i_{p_i}$ of $p_i \in V^i$ such that:
\par
\begin{enumerate}
\item
There exists an isomorphism
$(h,\tilde\varphi) : (V^1,\Gamma^1,\phi^1)\vert_{V^1_{p_1}} \to (V^2,\Gamma^2,\phi^2)\vert_{V^2_{p_2}}$
between orbifold charts of $X$, which are subcharts.
\item
There exists an isomorphism
$(h,\tilde{\hat{\varphi}}) : (V^1\times E^1,\Gamma^1,\phi^1)\vert_{V^1_{p_1}\times E^1} \to
(V^2\times E^2,\Gamma^2,\phi^2)\vert_{V^2_{p_2} \times E^2}$
between orbifold charts of $\mathcal E$, which are subcharts.
\item
For each $y \in V^1_{p_1}$ the map: $ E^1 \to E^2$,
$\xi \to \pi_{E^2}\tilde{\hat{\varphi}}(y,\xi)$ is a linear isomorphism.
Here $\pi_{E^2} : V^2 \times E^2 \to E^2$ is the projection.
\end{enumerate}
\item
A {\it representative of a vector bundle structure} of $(X,\mathcal E)$
is a set of orbifold charts $\{(V_i,E_i,\Gamma_i,\phi_i,\widehat\phi_i) \mid i \in I\}$
such that any two of the charts are compatible in the sense of (3)
above and
$$
\bigcup_{i\in I} \phi_i(V_i) = X,
\quad
\bigcup_{i\in I} \widehat\phi_i(V_i \times E_i) = \mathcal E,
$$
are locally finite open covers.
\index{orbifold ! representative of a vector bundle structure on orbifold}
\end{enumerate}
\end{defn}
\begin{defn}\label{def26222}
Suppose $(X^*,\mathcal E^*)$ $(* = a,b)$ have representatives
of vector bundle structures $\{(V^*_i,E^*_i,\Gamma^*_i,\phi^*_i,\widehat\phi^*_i) \mid i \in I^*\}$,
respectively.
A pair of orbifold embeddings $(f,\widehat f)$,
$f : X^a \to X^b$, $\widehat f : \mathcal E^a \to \mathcal E^b$ is said to be an
{\it embedding of vector bundles} if the following holds.
\index{orbifold ! embedding of vector bundles on orbifolds}
\index{embedding ! of vector bundles on orbifolds}
\begin{enumerate}
\item
Let
$p \in V^a_i$, $q \in V^b_j$ with
$f(\phi^a_i(p)) = \phi^b_j(q)$.
Then there exist open subcharts
$(V^a_{i,p}\times E^a_{i,p},\Gamma^a_{i,p},\widehat\phi^a_{i,p})$
and
$(V^b_{j,q}\times E^b_{j,q},\Gamma^b_{j,q}\widehat\phi^b_{j,q})$
and a local representative
$(h_{p;i,j},f_{p;i,j},\widehat f_{p;i,j})$ of the embeddings $f$ and $\widehat f$
such that
for each $y \in V^a_i$ the map
$\xi \mapsto \pi_{E^b}(\widehat f_{p;i,j}(y,\xi))$,
$E^a_{i,p} \to E^b_{j,q}$ is a linear embedding. Here $\pi_{E^b} : V^b \times E^b \to E^b$ is the projection.
\item
The diagram below commutes set-theoretically.
\begin{equation}\label{diag2633}
\begin{CD}
\mathcal E^a @ >{\widehat f}>>
\mathcal E^b  \\
@ V{\pi_{E^a}}VV @ VV{\pi_{E^b}}V\\
X^a @ > {f} >> X^b
\end{CD}
\end{equation}
\end{enumerate}
Two orbifold embeddings are said to be {\it equal}
if they coincide set-theoretically
as pairs of maps.
\end{defn}
\begin{lem}\label{lem26444AA1}
\begin{enumerate}
\item
A composition of embeddings of vector bundles is an embedding.
\item
A pair of identity maps is an embedding.
\item
If an embedding of vector bundles is a pair of homeomorphisms,
then the pair of their inverses is also an embedding.
\end{enumerate}
\end{lem}
The proof is easy and is omitted.
\begin{defn}\label{defn2820}
Let
$(X,\mathcal E,\pi)$ be as in Definition \ref{defn2613}.
\begin{enumerate}
\item
An embedding of vector bundles is said to be an {\it isomorphism}
\index{orbifold ! isomorphism of vector bundles on orbifolds}
if it is a pair of diffeomorphisms of orbifolds.
\item
We say that two representatives of a vector bundle structure of $(X,\mathcal E)$
are
{\it equivalent} if the pair of identity maps regarded as a map between
$(X,\mathcal E)$ equipped with those two representatives of  vector bundle
structures
is an isomorphism. This is an equivalence relation by Lemma \ref{lem26444AA1}.
\item
An equivalence class of the equivalence relation (2) is called
a {\it vector bundle structure} of $(X,\mathcal E)$.
\index{orbifold ! vector bundle structure}
\item
A pair $(X,\mathcal E)$ together with its vector bundle
structure is called a {\it vector bundle} on $X$.
\index{orbifold ! vector bundle on orbifold}
\index{orbifold ! total space of vector bundle}
\index{orbifold ! base space of vector bundle}
\index{orbifold ! projection of vector bundle}
We call $\mathcal E$ the {\it total space}, $X$ the
{\it base space}, and $\pi : \mathcal E \to X$ the {\it projection}.
\item
The condition for $(f,\widehat f) : (X^a,\mathcal E^a) \to (X^b,\mathcal E^b)$
to be an embedding
does not change if we replace representatives of vector bundle
structures
to equivalent ones. So we can define the notion of
an {\it embedding of vector bundles}.
\item
We say $(f,\widehat f)$ is an embedding {\it over the orbifold embedding $f$.}
\end{enumerate}
\end{defn}
\begin{rem}
\begin{enumerate}
\item
We may use the terminology `orbibundle' in place of vector bundle.
We use this terminology in case we emphasize that it is different
from the vector bundle over the underlying topological space.
\index{orbifold ! orbibundle}
\item
We sometimes simply say $\mathcal E$ is a vector bundle on an
orbifold $X$.
\end{enumerate}
\end{rem}
\begin{defn}\label{defn2655}
\begin{enumerate}
\item
Let $(X,\mathcal E)$ be a vector bundle. We call an orbifold chart $(V,E,\Gamma,\phi,\widehat\phi)$
in the sense of Definition \ref{defn2613} (1)
of underlying pair of topological spaces $(X,\mathcal E)$
an {\it orbifold chart of a vector bundle}  if
the pair of maps $(\overline\phi,\overline{\widehat\phi}) : (V/\Gamma,(V\times E)/\Gamma) \to (X,\mathcal E)$ induced from $(\phi,\widehat\phi)$ is an
embedding of vector bundles.
\item
If $(V,E,\Gamma,\phi,\widehat\phi)$
is an orbifold chart of a vector bundle $(X,\mathcal E)$ we say a pair $(E,\widehat\phi)$
a {\it trivialization} \index{orbifold ! trivialization of vector bundle} of our vector bundle on $V/\Gamma$.
\item
Hereafter when $(X,\mathcal E)$ is a vector bundle,
its `orbifold chart' always means an orbifold chart of a vector bundles in the
sense of (1).
\item
In case when a vector bundle structure on $(X,\mathcal E)$ is given,
a representative of this vector bundle structure is
called an {\it orbifold atlas} of $(X,\mathcal E)$.
\index{orbifold ! orbifold atlas}
\item
Two orbifold charts  $(V_i,E_i,\Gamma_i,\phi_i,\widehat \phi_i)$
are said to be {\it isomorphic} if there exist an
isomorphism $(h,\tilde\varphi)$ of orbifold charts
$(V_1,\Gamma_1,\phi_1) \to (V_2,\Gamma_2,\phi_2)$
and an isomorphism
$(h,\tilde{\hat{\varphi}})$ of orbifold charts
$(V_1\times E_1,\Gamma_1,\widehat\phi_1) \to (V_2\times E_2,\Gamma_2,\widehat\phi_2)$
such that they induce an embedding of vector bundles
$(\varphi,\hat\varphi) : (V_1/\Gamma_1,(V_1\times E_1)/\Gamma_1) \to (V_2/\Gamma_2,(V_2\times E_2)/\Gamma_2)$.
The triple
$(h,\tilde\varphi,\tilde{\hat{\varphi}})$ is called an {\it isomorphism} of
orbifold charts. \index{orbifold ! isomorphism of orbifold charts of vector bundles}
\end{enumerate}
\end{defn}
\begin{lem}\label{lem2619}
Let $(X^b,\mathcal E^b)$ be a vector bundle over an orbifold $X^b$
and $f : X^a \to X^b$  an embedding of orbifolds.
Let $\mathcal E^a = X^a \times_{X^b} \mathcal E^b$ be the fiber product
in the category of topological space.
By the definition of the fiber product we have maps
$\mathcal E^a \to X^a$ and $\mathcal E^a \to \mathcal E^b$.
We write them $\pi$ and $\widehat f$ respectively.
Then the exists a unique structure of vector bundle on $(X^a,\mathcal E^a)$
such that the projection is the above map $\pi$  and that
$(f,\widehat f)$ is an embedding of vector bundles.
\end{lem}
\begin{proof}
Let
$\{\frak V^*_{\frak r} \mid \frak r \in \frak R_*\}$, $*=a,b$
be orbifold atlases where
$\frak V^*_{\frak r} = (V^*_{\frak r},\Gamma^*_{\frak r},\phi^*_{\frak r})$.
Let
$(V^b_{\frak r},E^b_{\frak r},\Gamma^b_{\frak r},\phi^b_{\frak r},\widehat\phi^b_{\frak r})$
be orbifold atlas of the vector bundle $(X^b,\mathcal E^b)$.
Let
$(h_{\frak r,ba},\tilde\varphi_{\frak r,ba})$ be
a local representative of the embedding $f$ on the charts
$\frak V^a_{\frak r}$, $\frak V^b_{\frak r}$.
We put $E^a_{\frak r} = E^b_{\frak r}$, on which $\Gamma^a_{\frak r}$ acts
by the isomorphism $h_{\frak r,ba}$.
By definition of fiber product, there exists uniquely a map
$\widehat\phi^a_{\frak r} : V_{\frak r}^b \times E_{\frak r}^b \to \mathcal E^a$
such that the next diagram commutes.
\begin{equation}\label{diag2619diag}
\begin{CD}
V_{\frak r}^a @ <{\pi}<<
V_{\frak r}^a \times E_{\frak r}^b   @>{\tilde\varphi_{\frak r,ba}\times id}>> V_{\frak r}^b \times E_{\frak r}^b\\
@ V{\phi^a_{\frak r}}VV @ VV{\widehat\phi^a_{\frak r}}V
@VV{\widehat\phi^b_{\frak r}}V\\
X^a @ < {\pi} <<\mathcal E^a @>{\hat f}>>\mathcal E^b
\end{CD}
\end{equation}
In fact,
$$
f\circ \phi^a_{\frak r}\circ \pi
= \phi^b_{\frak r}   \circ \varphi_{\frak r,ba} \circ \pi
= \phi^b_{\frak r} \circ \pi  \circ (\tilde\varphi_{\frak r,ba}\times id)
= \pi \circ \widehat\phi^b_{\frak r} \circ (\tilde\varphi_{\frak r,ba}\times id).
$$
Thus $\{(V^a_{\frak r},E^a_{\frak r},\Gamma^a_{\frak r},\phi^a_{\frak r},\widehat\phi^a_{\frak r}) \mid \frak r \in \frak R\}$
is an atlas of the vector bundle $(X^a,\mathcal E^a)$.
\end{proof}
\begin{defn}
We call the vector bundle in Lemma \ref{lem2619} the {\it pullback}
\index{orbifold ! pullback of vector bundle on orbifold}
and write $f^*(X^b,\mathcal E^b)$.
(Sometimes we write $f^*\mathcal E^b$ by an abuse of notation.)
\par
In case $X^a$ is an open subset of $X^b$ equipped with
open substructure we call {\it restriction}
\index{orbifold ! restriction of vector bundle on orbifold}
in place of pullback of $\mathcal E^b$
and write $\mathcal E^b\vert_{X^a}$ in place of $f^*\mathcal E^b$.
\end{defn}
\begin{lem}\label{lem2622}
In the situation of Lemma \ref{lem26999}
suppose in addition that $\mathcal E^i$ is a vector bundle over $X^i$
and $\widehat\varphi_{21} : \mathcal E^1 \to \mathcal E^2$
is an embedding of vector bundles over $\varphi_{21}$.
Then in addition to the conclusion of Lemma \ref{lem26999},
there exists $\tilde{\hat{\varphi}}_{\frak r;21} :
V^{1}_{\frak r} \times E_{\frak r}^1
\to V^{2}_{\frak r} \times E_{\frak r}^2$
that is an $h_{\frak r;21}$ equivariant embedding of manifolds
with the following properties.
\begin{enumerate}
\item
The next diagram commutes.
\begin{equation}\label{diagin26777}
\begin{CD}
V^1_{\frak r}\times E^1_{\frak r} @ >{\tilde{\hat{\varphi}}_{\frak r,21}}>>
V^2_{\frak r} \times E^2_{\frak r} \\
@ V{\widehat\phi_{\frak r}^{1}}VV @ VV{\widehat\phi_{\frak r}^{2}}V\\
\mathcal E^1 @ > {\widehat\varphi_{21}} >> \mathcal E^2
\end{CD}
\end{equation}
\item
For each $y \in V^1_{\frak r}$ the map
$\xi \mapsto \pi_2(\tilde{\hat{\varphi}}_{\frak r,21}(y,\xi))$
$: E^1_{\frak r} \to E^2_{\frak r}$ is a linear embedding.
\end{enumerate}
\end{lem}
The proof is similar to the proof of Lemma \ref{lem26999}
and is omitted.
\begin{defn}\label{def28262826}
We call $(h_{\frak r,21},\tilde\varphi_{\frak r,21},\tilde{\hat{\varphi}}_{\frak r,21})$
a {\it local representative of embedding}  $(\varphi_{21},\widehat\varphi_{21})$
on the charts
$(V^1\times E^1,\Gamma^1,\widehat\phi^1)$, $ (V^2\times E^2,\Gamma^2,\widehat\phi^2)$.
\index{orbifold ! local representative of embedding of vector bundle on orbifold}
\index{embedding ! local representative of embedding of vector bundle on orbifold}
\end{defn}
\begin{lem}\label{lem2715}
If $(h_{\frak r,21},\tilde\varphi_{\frak r,21},\tilde{\hat{\varphi}}_{\frak r,21})$,
$(h'_{\frak r,21},\tilde\varphi'_{\frak r,21},\tilde{\hat{\varphi}}'_{\frak r,21})$
are local representatives of an embedding of vector bundles of the same
charts
$(V^1\times E^1,\Gamma^1,\widehat\phi^1)$, $ (V^2\times E^2,\Gamma^2,\widehat\phi^2)$,
then there exists $\mu \in \Gamma^2$
such that
$$
\tilde\varphi_{\frak r,21}'(x) = \mu\tilde\varphi_{\frak r,21}(x),
\quad
\tilde{\hat{\varphi}}_{\frak r,21}'(x,\xi) = \mu\tilde{\hat{\varphi}}_{\frak r,21}(x,\xi)
\quad
h_{\frak r,21}'(\gamma) = \mu h_{\frak r,21}(\gamma)\mu^{-1}.
$$
\end{lem}
\begin{proof}
This is a consequence of Lemma \ref{lem21tenhatena}.
\end{proof}
\begin{rem}
In Situation \ref{opensuborbifoldchart} we introduced the
notation
$(h_{\frak r,21},\tilde\varphi_{\frak r,21},\breve{{\varphi}}_{\frak r,21})$
where $\breve{{\varphi}}_{\frak r,21}$ is related to
$\tilde{\hat{\varphi}}_{\frak r,21}$ by the formula
$$
\tilde{\hat{\varphi}}_{\frak r,21}(y,\xi)
=(\tilde{{\varphi}}_{\frak r,21}(y),\breve{{\varphi}}_{\frak r,21}(y,\xi)).
$$
\end{rem}
We use the pullback of vector bundles in a different situation.
Let $\mathcal E^i$, $i=1,2$, are vector bundles over an orbifold $X$.
We take the Whitney sum  bundle $\mathcal E^1 \oplus \mathcal E^2$.
Let $\vert\mathcal E^1 \oplus \mathcal E^2\vert$ be its total
space. There exists a projection
\begin{equation}\label{projfromWhe}
\vert\mathcal E^1 \oplus \mathcal E^2\vert
\to \vert\mathcal E^2\vert.
\end{equation}

\begin{defnlem}\label{pullbackbyproj}
$\vert\mathcal E^1 \oplus \mathcal E^2\vert$
has a structure of vector bundle over  $\vert\mathcal E^2\vert$
such that (\ref{projfromWhe}) is the projection.
We write it as
$\pi_{\mathcal E^2}^*\mathcal E^1$
and call the {\it pullback} of $\mathcal E^1$ by the projection
$\pi_{\mathcal E^2} : \vert\mathcal E^2\vert \to X$.
\index{orbifold ! pullback of vector bundle by a projection of
another vector bundle}
\par
When $U$ is an open subset of $\vert\mathcal E^2\vert$ and
$\pi : U \to X$ is the restriction of $\pi_{\mathcal E^2}$ to $U$,
the pullback $\pi^*_{\mathcal E^2}\mathcal E^1$ is by definition
the restriction of $\pi_{\mathcal E^2}^*\mathcal E^1$ to $U$.
\end{defnlem}
The proof is immediate from definition.
\begin{rem}
We note that the total space $\vert\mathcal E^1 \oplus \mathcal E^2\vert$
is {\it not} a fiber product
$\vert\mathcal E^1 \vert \times_X \vert\mathcal E^2\vert$.
In fact, if $X$ is a point and $\mathcal E^1
= \mathcal E^2= \R^n/\Gamma$ with linear $\Gamma$ action,
then the fiber of $\vert\pi_{\mathcal E^2}^*\mathcal E^1 \vert \to
\vert\mathcal E^2\vert \to X$ at $[0]$ is
$(E^1 \times E^2)/\Gamma$.
The fiber of the map
$\vert\mathcal E^1 \vert \times_X \vert\mathcal E^2\vert \to X$
at $[0]$ is $(E^1/\Gamma) \times (E^2/\Gamma)$.
\end{rem}
\begin{defn}
Let $(X,\mathcal E)$ be a vector bundle. A {\it section}
of $(X,\mathcal E)$ is an embedding of orbifolds $s : X \to \mathcal E$
such that the composition of $s$ and the projection
is the identity map (set-theoretically).
\index{orbifold ! section of a vector bundle on orbifold}
\end{defn}

\begin{lem}\label{lem2626}
Let $\{(V_{\frak r},E_{\frak r},\Gamma_{\frak r},
\psi_{\frak r},\widehat\psi_{\frak r})\mid \frak r \in \frak R\}$
be an atlas of $(X,\mathcal E)$.
Then a section of  $(X,\mathcal E)$
corresponds one to one to the following object.
\begin{enumerate}
\item
For each $\frak r$ we have a $\Gamma_r$ equivariant
map $s_{\frak r} : V_{\frak r} \to E_{\frak r}$,
which is compatible in the sense of (2) below.
\item
Suppose $\phi_{\frak r_1}(x_1) = \phi_{\frak r_2}(x_2)$.
Then the definition of orbifold atlas implies that
there exist subcharts
$(V_{\frak r_i,x_i},E_{\frak r_i,x_i},\Gamma_{\frak r_i,x_i},
\phi_{\frak r_i,x_i},\widehat\phi_{\frak r_i})$
of the orbifold charts
$(V_{\frak r_i},E_{\frak r_i},\Gamma_{\frak r_i},
\phi_{\frak r_i},\widehat\phi_{\frak r_i})$
at $x_i \in V_{\frak r_i}$ for $i = 1,2$ and
an isomorphism of charts
$$
\aligned
(h^{\frak r,p}_{12},\tilde\varphi^{\frak r,p}_{12},
\tilde{\hat{\varphi}}^{\frak r,p}_{12}) ~:~
&(V_{\frak r_2,x_2},E_{\frak r_2},\Gamma_{\frak r_2,x_2},
\phi_{\frak r_2,x_2},\widehat\phi_{\frak r_2})
\\
&\to (V_{\frak r_1,x_1},E_{\frak r_1,x_1},\Gamma_{\frak r_1,x_1},
\phi_{\frak r_1,x_1},\widehat\phi_{\frak r_1}).
\endaligned
$$
\par
Now we require  the next equality:
\begin{equation}\label{sectioncompati}
\tilde{\hat{\varphi}}^{\frak r,p}_{12}(s_{\frak r_1}(y,\xi)) =
s_{\frak r_2}(\tilde\varphi^{\frak r,p}_{12}(y),\xi).
\end{equation}
\end{enumerate}
\end{lem}
\begin{proof}
The proof is mostly the same as the corresponding
standard result
in the case of vector bundle on a manifold or on a topological space.
Let $s : X \to \mathcal E$ be a section, which is an orbifold embedding.
Let $p \in \phi_{\frak r}(V_{\frak r})$.
Then there exist a subchart $(V_{\frak r,p},\Gamma_{\frak r,p},\phi_{\frak r,p})$
of $\frak V_{\frak r}$ and
a subchart $(\widehat V_{\frak r,\tilde p},
\Gamma_{\frak r,\tilde p},\phi_{\frak r,\tilde p})$
of $(V_{\frak r}\times E_{\frak r},\Gamma_{\frak r},\phi_{\frak r,\tilde p})$
such that  a representative  $(h',\tilde\varphi')$
of $s$ exists on this subcharts.
Since $\pi \circ s = $identity set-theoretically,
it follows that $\pi_1(\tilde\varphi(y)) \equiv y \mod \Gamma_p$
for any $y \in V_{\frak r,p}$.
We take $y$ such that $\Gamma_{y} = \{1\}$.
Then, there exists {\it uniquely} $\mu \in \Gamma_p$
such that $\pi_1(\tilde\varphi'(y)) \equiv \mu y$.
By continuity this $\mu$ is independent of $y$.
(We use Condition \ref{convinv} here.)
\par
We replace $\tilde p$ by $\mu^{-1}\tilde p$
and $(\widehat V_{\frak r,\tilde p},\Gamma_{\frak r,\tilde p},\phi_{\frak r,\tilde p})$
by $(\mu^{-1}\widehat V_{\frak r,\tilde p},\mu^{-1}\Gamma_{\tau,\tilde p} \mu,\phi_{\frak r,\tilde p}\circ \mu)$
and $(h',\tilde\varphi')$ by $(h'\circ {\rm conj}_{\mu},\tilde\varphi' \circ \mu^{-1})$.
(Here ${\rm conj}_{\mu}(\gamma) = \mu \gamma \mu^{-1}$.)
Therefore we may assume $\pi_1(\tilde\varphi'(y)) = y$.
Note that $\tilde\varphi'$ is $h'$-equivariant and $\pi_1$ is ${\rm id}$-equivariant.
Here ${\rm id}$ is the identity map $\Gamma_{\frak r,y} \to \Gamma_{\frak r,y}$.
Therefore the identity map $V_{\frak r,p} \to V_{\frak r,p}$ is
$h'$ equivariant.
Hence $h' = {\rm id}$.
\par
In sum we have the following.
(We put $s_{\frak r,p} = \tilde\varphi'$.)
For a sufficiently small $\frak V_{\frak r,p}$
there exists uniquely a map $
s_{\frak r,p} : V_{\frak r,p} \to V_{\frak r,p} \times E_{\frak r}$
such that
\begin{enumerate}
\item[(a)]
$
\pi_1(s_{\frak r,p}(x)) =x
$
\item[(b)]
$s_{\frak r,p}$ is equivariant with respect to the embedding
$\Gamma_{\frak r,p} \to \Gamma_{\frak r}$.
(Recall $\Gamma_{\frak r,p} =
\{\gamma \in \Gamma_{\frak r} \mid \gamma p = p\}$.)
\item[(c)]
$({\rm id},s_{\frak r,p})$ is a local representative of $s$.
\end{enumerate}
We can use uniqueness of such $s_{\frak r,p}$ to glue them
to obtain a map $V_{\frak r} \to V_{\frak r} \times E_{\frak r}$.
By (a)  this map is of the form
$x \mapsto (x,\frak s_{\frak r}(x))$.
This is the map $\frak s_{\frak r}$ in (1).
Since $x \mapsto \gamma^{-1}\frak s_{\frak r}(\gamma x)$
also has the same property, the uniqueness implies that
$\frak s_{\frak r}$ is $\Gamma_{\frak r}$
equivariant.
(\ref{sectioncompati}) is also a consequence of
the uniqueness.
\par
We thus find a map from the set of sections to the
set of $(s_{\frak r})_{\frak r \in \frak R}$ satisfying (1)(2).
The construction of the converse map is obvious.
\end{proof}
The next lemma is proved during the proof of Lemma \ref{lem2626}.

\begin{lem}\label{lem2627}
Let $(V_{\frak r},E_{\frak r},\Gamma_{\frak r},\phi_{\frak r},\widehat\phi_{\frak r})$ be an orbifold chart of
$(X,\mathcal E)$ and $s$ a section of $(X,\mathcal E)$.
Then there exists uniquely a $\Gamma$ equivariant map
$s_{\frak r} : V_{\frak r} \to E_{\frak r}$ such that the following diagram commutes.
\begin{equation}\label{diagin26277XXrev}
\begin{CD}
V_{\frak r}\times E_{\frak r} @ >{\widehat\phi_{\frak r}}>>
\mathcal E_{\frak r} \\
@ A{{\rm id} \times s_{\frak r}}AA @ AA{s}A\\
V_{\frak r} @ > {\phi_{\frak r}} >> X
\end{CD}
\end{equation}
\end{lem}
\begin{defn}\label{defnlocex}
We call the system of maps $s_{\frak r}$ the {\it local expression} of $s$ in the orbifold
chart $(V_{\frak r},E_{\frak r},\Gamma_{\frak r},\phi_{\frak r},\widehat\phi_{\frak r})$.
\index{orbifold ! local expression of a section on orbifold chart}
\end{defn}
We next review the proof of a few well-known facts on
pullback bundle etc.. Those proofs are straightforward generalization of the
corresponding results in manifold theory. We include them for
completeness' sake only.
\begin{prop}\label{homotopicpulback}
Let $\mathcal E$ is a vector bundle on $X \times [0,1]$, where
$X$ is an orbifold. We identify $X \times \{0\}$, $X \times \{1\}$
with $X$ in an obvious way. Then
there exists an isomorphism of vector bundles
$$
I : \mathcal E\vert_{X \times \{0\}} \cong \mathcal E\vert_{X \times \{1\}}.
$$
Suppose in addition that there exists a compact set $K \subset X$
an its neighborhood $V$ and isomorphism
$$
I_0 : \mathcal E\vert_{V \times [0,1]} \cong \mathcal E\vert_{V \times \{0\}} \times [0,1]
$$
Then we may choose $I$ so that it coincides with the isomorphism induced by
$I_0$ on $K$.
\par
In the case $K$ is a submanifold we may take $K=V$.
\end{prop}
\begin{proof}
To prove the proposition we use the notion of connection of vector
bundle on orbifolds.
Note a vector field on an orbifold is a section of the tangent bundle.
\begin{defn}
A {\it connection}\index{connection} of a vector bundle $(X,\mathcal E)$ is
an $\R$
linear map
$$
\nabla : C^{\infty}(TX) \otimes_{\R} C^{\infty}(\mathcal E)
\to C^{\infty}(\mathcal E)
$$
such that $\nabla_X(V) = \nabla(X,V)$ satisfies
$$
\nabla_{fX}(V) = f\nabla_X(V), \qquad
\nabla_{X}(fV) = f\nabla_X(V) + X(f)V.
$$
\end{defn}
Here $C^{\infty}(\mathcal E)$ is the vector space consisting of
all the smooth sections of $\mathcal E$.
\par
For any connection $\nabla$ and piecewise smooth map $
\ell : [a,b] \to X$ we obtain parallel transport
$$
{\rm Pal}^{\nabla} : \mathcal E_{\ell(a)}  \to \mathcal E_{\ell(b)}
$$
in the same way as the case of manifold.
\begin{rem}
Here $\mathcal E_{\ell(a)}$ is the fiber of $\mathcal E$ at $\ell(a) \in  X$
and is defined as follows.
We take a chart $(V_{\frak r},E_{\frak r},\Gamma_{\frak r},\psi_{\frak r},\widehat\psi_{\frak r})$ of $(\mathcal E,X)$
at $\ell(a)$. Then $\mathcal E_{\ell(a)} = E_{\frak r}$.
If $(V_{\frak r'},E_{\frak r'},\Gamma_{\frak r'},\psi_{\frak r'},\widehat\psi_{\frak r'})$ is another chart
we can identify $E_{\frak r}$ and $E_{\frak r'}$ by
$\xi \mapsto \breve{\varphi}_{\frak r'\frak r}(\xi,y)$ where $\psi_{\frak r}(y) = \ell(a)$ and
$ \breve{\varphi}_{\frak r'\frak r} : V_{\frak r} \times E_{\frak r} \to E_{\frak r'}$ is a part of the
coordinate change.
(Situation \ref{opensuborbifoldchart}.)
\par
Note the identification $\xi \mapsto \breve{\varphi}_{\frak r'\frak r}(\xi,y)$ is well-defined up to the $\Gamma_{\ell(a)}
= \{ \gamma \in \Gamma_{\frak r} \mid \gamma(y) = y\}$ action.
\par
The parallel transport $
{\rm Pal}^{\nabla} : \mathcal E_{\ell(a)}  \to \mathcal E_{\ell(b)}
$ is well-defined up to $\Gamma_{\ell(a)} \times \Gamma_{\ell(b)}$ action.
\end{rem}
\begin{lem}
Any vector bundle $\mathcal E$ over orbifold $X$ has a connection.
Moreover if a connection is given for $\mathcal E\vert_{V}$
where $V$ is an open neighborhood of a compact subset $K$ of $X$,
then we can extend it without changing it on $K$.
In the case $K$ is a submanifold we may take $K=V$.
\end{lem}
The proof is an obvious analogue of the proof of the existence of
connection of a vector bundle on a manifold, which uses partition of unity.
\par
We start the proof of Proposition \ref{homotopicpulback}.
We take a connection of $\mathcal E\vert_V$.
We then take direct product connection on
$\mathcal E\vert_{V \times \{0\}} \times [0,1]$,
and use $I_0$ to obtain a connection on
$ \mathcal E\vert_{V \times [0,1]}$.
We extend it to a connection on $\mathcal E$
without changing it on $K \times [0,1]$.
Let $x \in X$ then we can use parallel transportation
along the path $t \mapsto (x,t)$ to obtain an isomorphism
$\mathcal E_{(x,0)} \cong \mathcal E_{(x,1)}$.
We thus obtain set theoretical map
$$
\vert\mathcal E\vert_{X \times \{0\}}\vert \cong \vert\mathcal E\vert_{X \times \{1\}}\vert.
$$
It is easy to see that it induces an isomorphism of vector bundles.
\par
Using the fact that our connection is direct product on $K \times [0,1]$,
we can check the second half of the statement.
\end{proof}
\begin{defn}
We say two embeddings of orbifold $f_i : X \to Y$
($i=1,2$) to be {\it isotopic}
\index{isotopic} each other
if there exists an embedding of orbifolds
$H : X\times [0,1] \to Y \times [0,1]$ such that
the second factor of $H(x,t)$ is $t$ and that
$$
H(x,0) = (f_1(x),0) \qquad
H(x,1) = (f_2(x),1).
$$
Suppose $V \subset X$ and $f_1 = f_2$ on a neighborhood $V$ of $K$. We say
$f_1$ is {\it isotopic to $f_2$ relative to $K$}
\index{isotopic relative to a compact subset} if we may take $H$ such that
\begin{equation}\label{homoisidentity}
H(x,t) = (f_1(x),t) = (f_2(x),t)
\end{equation}
for $x$ in a neighborhood of $K$.
\par
In the case $K$ is a submanifold we may take $K=V$
and then (\ref{homoisidentity}) holds for $x \in K$.
\end{defn}
\begin{cor}\label{cor2939}
Let $f_i : X \to Y$  be two embeddings which are isotopic and
$\mathcal E$ is a vector bundle on $Y$. Then
the pullback bundle $f_1^*\mathcal E$ is isomorphic
to $f_2^*\mathcal E$.
\par
If $f_1 = f_2$ on a neighborhood of $K \subset X$ and
$f_1$ is isotopic to $f_2$ relative to $K$
then we may choose the isomorphism $f_1^*\mathcal E \cong f_2^*\mathcal E$
so that its restriction to $K$ is the identity map.
\end{cor}
\begin{proof}
This is an immediate consequence of
Proposition \ref{homotopicpulback} and the definition.
\end{proof}
We next recall Definition \ref{lem123000} which we copy below.
\begin{defn}\label{lem1230002}
Let $f : X \to Y$ be an embedding of orbifolds and $K\subset X$  a
compact subset and $U$ be an open neighborhood of $K$ in $Y$. We say that
a continuous map $\pi : U \to X$ is diffeomorphic to the
projection of normal bundle if the following holds.
\par
Let ${\rm pr} : N_XY \to X$ be the normal bundle. Then there exists
a neighborhood $U'$  of $K$ in $ N_XY$, (Note $K\subset X \subset N_XY$.) and a diffeomorphism
of orbifolds $h : U' \to U$ such that $\pi\circ h = {\rm pr}$.
We also require that $h(x) = x$ for $x$ in a neighborhood of $K$ in  $X$.
\end{defn}
\begin{defnlem}\label{defpullbackbundlenbd}
Suppose $\pi : U \to X$ is diffeomorphic to the
projection of normal bundle as in Definition \ref{lem1230002}
and $\mathcal E$ be a vector bundle on $X$. We define
$\pi^*\mathcal E$, the pullback bundle as follows.
\par
Let $h$, $U'$ be as in Definition \ref{lem1230002}.
We defined a pull back bundle ${\rm pr}^*\mathcal E$ on $N_XY$
in Definition \ref{pullbackbyproj}.
We put
$$
\pi^*\mathcal E = (h^{-1})^*{\rm pr}^*\mathcal E\vert_{U'}.
$$
This is independent of the choice of $(U',h)$ in the following sense.
Let $U'_i$, $h_i$ ($i=1,2$) be two choices.
Then we can shrink $U$ and $U'_i$ so that the restriction of $h_i$
becomes an isomorphism between them.
Then
\begin{equation}\label{294949}
(h_1^{-1})^*{\rm pr}^*\mathcal E\vert_{U'_1}
\cong
(h_2^{-1})^*{\rm pr}^*\mathcal E\vert_{U'_2}.
\end{equation}
Moreover the isomorphism (\ref{294949})
can be taken so that the following holds in addition.
We regard $K \subset U$.
Then by definition it is easy to see that
the restriction of both sides of (\ref{294949})
are canonically identified with
the restriction of $\mathcal E$ to $K \subset X$.
The isomorphism  (\ref{294949}) becomes the identity map
on $K$ by this isomorphism.
\end{defnlem}
\begin{proof}
We can replace $U$ by a smaller open neighborhood so that
$h_1^{-1} : U \to N_XY$ is isotopic to $h_2^{-1} : U \to N_XY$.
(See the proof of Proposition \ref{prop2942} below.)
Then (\ref{294949}) follows from Corollary \ref{cor2939}.
\par
The second half of the claim follows also from the second half of
Corollary \ref{cor2939}.
\end{proof}
Actually the pullback bundle is independent of $\pi$ but depend only
on $U$ in the situation of Definition \ref{lem1230002}.
In fact we have
\begin{prop}\label{prop2942}
Let $\pi_i : U \to X$ be as in Definition \ref{lem1230002} for $i=1,2$.
Then there exists a neighborhood $U_0$ of $X$ in $Y$ and
$f : U_0 \to U$ such that
\begin{enumerate}
\item
$
\pi_2 \circ f = \pi_1
$
\item
$f : U_0 \to U$ is isotopic to identity relative to $X$.
\end{enumerate}
\end{prop}
\begin{proof}
Let $h_i : U'_i \to U_i$ be as in  Definition \ref{lem1230002}.
We put $f = h_2 \circ h^{-1}_1$ which is defined for sufficiently
small  $U_0$. If suffices to show that $f$ is isotopic to identity.
We first prove it in the case when the following additional assumption
is satisfied. (We will  remove this assumption later.)
\begin{assump}\label{assym2929}
For any $x \in K \subset N_XY$ the first derivative at $x$, $D_xf : T_x(N_XY) \to T_x(N_XY)$ is the
identity map.
\end{assump}
\par
In the case of manifold we can prove Proposition \ref{prop2942} in this case,  by observing $f$ is $C^1$-close to the identity map.
Then for example by using minimal geodesic we can show that $f$ is isotopic to identity.
\par
In the case of orbifold we need to work out this last step a bit more
carefully since the number
$$
\inf \{ r \mid \text{if $d(x,y) < r$ the minimal geodesic joining $x$ and $y$ is unique}\}
$$
can be $0$ in general.
\par
We will prove certain lemmata to clarify this point.
We need certain digression to state the lemmata.
We can define the notion of Riemannian metric of orbifold $X$ in an obvious way.
For $p\in X$ we have a geodesic coordinate
$(TB_p(c_p),\Gamma_p,\psi_p)$ where
$$
TB_p(c_p) = \{ \xi \in T_cX \mid \Vert \xi\Vert < c\}.
$$
The group $\Gamma_p$ is the isotropy group of the orbifold chart of $X$ at $p$.
The uniformization map $\psi : TB_p(c) \to X$ is defined by using
minimal geodesic in the same way as the usual Riemannian geometry.
We remark that this map is well defined up to the action of $\Gamma_p$.
\par
We need to take the number $c$ small  so that $\psi$ induces homeomorphism
$TB_p(c)/\Gamma_p \to X$. We can not choose $c$ uniformly away from $0$
even in the compact set. (This is because $d(p,q) < c_p$ must imply
$\#\Gamma_q \le \# \Gamma_p$.)
However we can prove the following.
\par
Let $X$ be an orbifold and $Z$ be a compact set. Suppose
$B_{c_0}(Z) = \{x \mid d(x,Z) \le c_0\}$ is complete with respect to the metric induced by
the Riemannian metric.
\begin{lem}
Let $Z \subset X$ be a compact subset.
Then there exists a finite set $\{p_i \mid j \in J\} \subset Z$ and $0 < c_j < c_0$
such that
\begin{enumerate}
\item The geodesic coordinate $(TB_{p_j}(c_j),\Gamma_{p_j},\psi_{p_j})$
exists.
\item
$$
Z \subset \bigcup_{j} \psi_{p_j}(TB_{p_j}(c_j/2)).
$$
\end{enumerate}
\end{lem}
The proof is immediate from the compactness of $Z$.
We call such $\{(TB_{p_j}(c_j),\Gamma_{p_j},\psi_{p_j}) \mid j\}$
a {\it geodesic coordinate system}\index{geodesic coordinate system} of $(X,Z)$.
We put $P = \{p_j \mid j =1,\dots, J\}$.
\begin{defn}\label{lem2945}
We fix a geodesic coordinate system of $(X,K)$.
Let $Z_0 \subset Z$ be a compact subset containing $P$. Suppose $F : U \to X$ be an
embedding of orbifold where $U \supset Z$ is an open neighborhood of $Z$.
We say $F$ is {\it $C^1$-$\epsilon$ close to identity} \index{$C^1$-$\epsilon$ close to identity} on $Z_0$ if the following holds.
\begin{enumerate}
\item
$F(B_{p_j}(c_j/2)) \subset B_{p_j}(c_j)$.
\item
There exists
$\tilde F_{j} : B_{p_j}(c_j/2) \to B_{p_j}(c_j)$
such that:
\begin{enumerate}
\item
$
\psi_{p_j} \circ  \tilde F_{j} = F \circ \psi_{p_j}$.
\item
$d(x,\tilde F_{j}(x)) < \epsilon$ for $x \in TB_{p_j}(c_j/2) \cap \psi_j^{-1}(Z_0)$.
\item
$d(D_x\tilde F_{j},id) < \epsilon$ for $x \in TB_{p_j}(c_j/2) \cap \psi_j^{-1}(Z_0)$.
\end{enumerate}
\smallskip
Here $d$ in Item (b) is the standard metric on Euclidean space $T_{p_j}X$
(together with metric induce by our Riemannian metric),
$d$ in Item (c) is a distance in the space of $n\times n$ matrices.
($n = \dim X$. We use our Riemannian metric to define a metric
on this space of matrices, which is a vector space of dimension $n^2$ with metric.)
\end{enumerate}
\end{defn}
\begin{lem}\label{lem2946}
For each $Z$ and a geodesic coordinate system of $(X,Z)$
there exists $\epsilon$ such that the following
holds for any $Z_0$ and $F : U \to X$.
\par
If $F$ is $C^1$-$\epsilon$ close to identity on $Z_0$
then $F$ is isotopic to the identity on $Z_0$.
\par
Moreover for any $\delta$ there exists $\epsilon(\delta)$ that
that if $F$ is $C^1$-$\epsilon(\delta)$ close to identity on $Z_0$
then the isotopy from $F$ to the identity map is taken
to be  $C^1$-$\delta$ close to identity on $Z_0$.
\end{lem}
\begin{proof}
We first observe that if $\epsilon$ is sufficiently small then the maps
$\tilde F_j$ satisfying Definition \ref{lem2945} (2) (a),(b) and (c) is unique.
In fact such $\tilde F_j$ is unique up to the action of
$\Gamma_{p_j}$.
Since the $\Gamma_{p_j}$ action is effective,
$\Gamma_{p_j}$ is a finite group and $p_j \in Z_0$,
we find that at most one such $\tilde F_j$ can satisfy (c).
(Note the map $\Gamma_{p_j} \to O(n)$ taking the linear part at $p_j$ of the action
is injective since $\Gamma_{p_j}$ action has a fixed point and is effective.)
\par
We define for $t \in [0,1]$ a map
$$
\tilde F_{t,j} : V_j \to TB_{p_j}(c_j/2)
$$
as follows. Here $V_j$ is a sufficiently small neighborhood of
$TB_{p_j}(c_j/2) \cap \psi_j^{-1}(Z_0)$.
We take a Riemannian metric on $V_j$ which is a pullback
of our Riemannian metric on $X$ by $\psi_{p_j}$.
By choosing $\epsilon$ sufficiently small and using (b),
there exists a unique minimal geodesic $\ell_{x,j} : [0,1] \to TB_{p_j}(c_j/2+2\epsilon)$ joining
$x$ to $\tilde F_{j}(x)$. We put
$$
\tilde F_{t,j}(x) = \ell_{x,j}(t).
$$
In the same way as the proof of the uniqueness of $\tilde F_j$ we can show that there exists
$F_t$ such that
$
\psi_{p_j} \circ  \tilde F_{t,j} = F_t \circ \psi_{p_j}$.
Using Definition \ref{lem2945} (2) (b) and (c) we can show that $F_t$ together with its first derivative is
close to the identity map. We can use it to show that $F_t$ is a  diffeomorphism to
its image.
Thus
$F_t$ is the required isotopy from $F$ to the identity map.
\end{proof}
We now use Lemma \ref{lem2946} to prove Proposition \ref{prop2942}
under Assumption \ref{assym2929}.
\par
Note we put $f = h_2 \circ h^{-1}_1$ and we want to show that $f$ is isotopic to
the identity map in a neighborhood of $K \subset X \subset Y$.
\par
We take a finite cover $
K \subset \bigcup_{j} \psi_{p_j}(TB_{p_j}(c_j/2))
$
by geodesic coordinate, where $p_j \in K$.
We take a compact neighborhood $Z \subset K$ such that
$Z \subset \bigcup_{j} \psi_{p_j}(TB_{p_j}(c_j/2))$.
We apply Lemma \ref{lem2945} to obtain $\epsilon$.
Note $f$ is the identity map on $K$. Moreover its first derivative
is identity at $K$ by Assumption \ref{assym2929}.
Therefore we can find a compact neighborhood $Z_0$ of $K$ sufficiently small
so that $f$ is $C^1$-$\epsilon$ close to identity on $Z_0$.
Thus Lemma \ref{lem2945}  implies
that $f$ is isotopic to the identity map.
The proof of Proposition \ref{prop2942} under the additional Assumption \ref{assym2929}
is complete.
\par
To remove Assumption \ref{assym2929},
we use the next lemma.
\begin{lem}\label{lem292333}
Let $U$ be an open neighborhood of $K$ in $N_XY$
and $F : U \to N_XY$ be an open embedding of orbifold.
\par
Assume $F =$ identity on a neighborhood of $K$ in $X$
and $D_xF(V) \equiv V \mod T_xV$ for $x$ in a neighborhood of $K$.
\par
Then there exists a smaller neighborhood $U'$ of $K$ such that
the restriction of $F$ to $U'$ is isotopic to the embedding satisfying Assumption \ref{assym2929}.
\end{lem}
\begin{proof}\label{lem2947}
We take the first derivative of $F$ at points $x$ in a neighborhood of  $K$ in $X$ and obtain
$$
D_x F : T_x N_XY \to N_XY
$$
Note $T_xX \subset N_XY$ is preserved by this map
and $T_x N_XY = T_xX \oplus (N_XY)_x$.
Therefore there exists linear bundle map
$$
H : N_XY \to TX
$$
on a neighborhood of $K$ such that
$$
D_xF (V_1,V_2) = (V_1 + H_x(V_2), V_2),
$$
where $V_1 \in  T_xX$ and $V_2 \in (N_XY)_x$.
Now we define $G_t : U' \to N_XY$ as follows.
($U'$ is a small neighborhood of $K$ in $N_XY$.)
Let $(x,V) \in U'$. Here $x$ is in a neighborhood of $K$ in $X$ and $V \in (N_XY)_x$.
We take a geodesic $\ell : [0,1] \to X$ of constant speed with $\ell(0) = x$ and $D\ell/dt(0) = H_x(V)$.
Let $\ell_{\le t}$ be its restriction to $[0,t]$.
Then $G_t(x,V) = (\ell(t),{\rm Pal}_{\ell_{\le t}}(V))$,
here ${\rm Pal}_{\ell_{\le t}}(V) \in (N_XY)_{\ell(t)}$ is the
parallel transport along $\ell_{\le t}$.
By construction the first derivative of $G_t$ at a point in $K$ is given by
$$(V_1,V_2)
\mapsto
(V_1 + tH_x(V_2), V_2),
$$
which is invertible.
Therefore,
if $V'$ is sufficiently small neighborhood of $K$, then the restriction of $G_t$ is an embedding $V' \to N_XY$.
\par
Note $F\circ G_1^{-1}$ satisfies Assumption \ref{assym2929}.
The proof of Lemma \ref{lem292333} is complete.
\end{proof}
Using Lemma \ref{lem292333} we can reduce the general case of
Proposition \ref{prop2942} to the case when Assumption \ref{assym2929}
is satisfied.
The proof of Proposition \ref{prop2942} is complete.
\end{proof}

We used the next result in Subsection \ref{subsection:bdlextcompa}.
\begin{prop}\label{prop2949}
Let $f : X \to Y$ be an embedding of orbifolds and $K_i$  compact subsets of $X$
for $i=1,2$ such that $K_1 \subset {\rm Int}K_2$.
Suppose $U_1$ is an open neighborhood $K_1$ in $Y$ and $\pi_1 : U_1 \to X$
such that it is diffeomorphic to the normal bundle.
\par
Then there exists an open neighborhood $U_2$ of $K_2$ in $Y$ and $\pi_2 : U_2 \to X$
such that it is diffeomorphic to the normal bundle
and $\pi_1 = \pi_2$ on a open neighborhood of $K_1$.
\end{prop}
\begin{proof}
In the case of manifold this is a standard result and can be proved using
isotopy extension lemma. By applying Lemma \ref{lem2946} we can
prove it in the same way in orbifold case.
For completeness' sake we give detail of the proof below.
\par
We first apply \cite[Lemma 6.5]{fooooverZ} to obtain
$U'_2$ and $\pi'_2 : U'_2 \to X$ such that
$U'_2 \supset K_2$ and $(U'_2,\pi'_2)$ is diffeomorphic to the normal bundle.
\par
We modify it to obtain $\pi_2$ so that $\pi_1 = \pi_2$ on an open neighborhood of $K_1$.
The detail follows.
\par
Let $W_1^{(i)}$ be a neighborhood of $K_1$ in $X$ such that
$$
\overline{W_1^{(1)}} \subset W_1^{(2)} \subset \overline{W_1^{(2)}}
\subset U_1 \cap X,
$$
Let $\Omega$ be an open subset of $U_1$ with
$$
\overline{W_1^{(2)}} \subset \Omega
\subset \overline\Omega \subset U_1.
$$
Later on we will choose $\Omega$ so that it is in a sufficiently small
neighborhood of $X$.
We put
$$
V_1^{(i)} = \pi_1^{-1}(W_1^{(i)}) \cap  \Omega.
$$
Let $\tilde\chi : {\rm Int}(K_2) \cup \Omega \to [0,1]$ be  a smooth function such that
$$
\tilde\chi =
\begin{cases}
1  &\text{on $W_1^{(1)}$} \\
0  &\text{on the complement of $W_1^{(2)}$}.
\end{cases}
$$
and put $\chi = \tilde\chi \circ \pi'_2$.
We may choose $\Omega$ small so that the following holds.
$$
\chi =
\begin{cases}
1  &\text{on $V_1^{(1)}$} \\
0  &\text{on the complement of $V_1^{(2)}$}.
\end{cases}
$$
\par
Let $Z =  \overline{W_1^{(2)}} \setminus  {W_1^{(1)}}$.
We take a neighborhood $U'$ of $Z$ and restrict
$\pi_1$ and $\pi'_2$ there.
Then we can apply  Proposition \ref{prop2942}
to prove that there exists an isotopy $F_t : U' \to X$ such that
$\pi'_2 \circ F_1 = \pi_1$ and $F_0$ is the identity map.
Now we put
$$
\pi_2(x)
=
\begin{cases}
(\pi_2 \circ F_{\chi(x)})(x) &\text{on $U'$} \\
\pi_1 &\text{ on $V_1^{(1)}$} \\
\pi'_2(x) &\text{elsewhere on $U'_2$}.
\end{cases}
$$
It is easy to see that they are glued to define a map.
To complete the proof it suffices to show that
$x \mapsto F_{\chi(x)}(x)$ is an embedding $: U' \to X$.
We will prove it below.
\par
We first assume that  Assumption \ref{assym2929} is satisfied for the map $f : U' \to X$ with
$\pi'_2 \circ f = \pi_1$.
In this case we may choose the isotopy $F_t$ to be arbitrary close to the identity
map in $C^1$ sense by taking $\Omega$ small.
(This is the consequence of the second half of Lemma \ref{lem2946}.)
Therefore the first derivative of $x \mapsto F_{\chi(x)}(x)$ is close to
identity. It follows that this map is an embedding.
\par
We finally show that we can choose $(U'_2,\pi'_2)$ so that
Assumption \ref{assym2929} is satisfied for the map $f : U' \to X$ with
$\pi'_2 \circ f = \pi_1$.
\par
We consider the fiber $\pi_1^{-1}(x)$ of $\pi_1$.
We may choose the Riemannian metric of $X$ in a neighborhood of $K_1$
so that this fiber is perpendicular to $X$ for any $x$  in a neighborhood of $K_1$.
We now extend this Riemannian metric to the whole $X$.
We use this Riemannian metric and exponential map in the
normal direction to identify a neighborhood of $K_2$ with a
normal bundle and to obtain $U'_2$ and $\pi'_2$.  Then Assumption \ref{assym2929} is satisfied.
\par
The proof of Proposition \ref{prop2949} is now complete.
\end{proof}

\newpage

\bibliographystyle{amsalpha}

\include{index}
\printindex
\end{document}